\documentclass[12pt,a4paper,twoside,reqno]{amsart}   

\usepackage{latexsym}
\usepackage{amsmath}
\usepackage{amsfonts}
\usepackage{amssymb}
\usepackage{lscape}
\usepackage{boxedminipage}
\usepackage[matrix,arrow]{xy}
\usepackage{longtable}

\numberwithin{equation}{section}

\theoremstyle{plain}
\newtheorem{thm}{Theorem}[section]
\newtheorem{cor}[thm]{Corollary}
\newtheorem{lem}[thm]{Lemma}
\newtheorem{prop}[thm]{Proposition}

\theoremstyle{definition}
\newtheorem{Def}[thm]{Definition}
\newtheorem{rem}[thm]{Remark}

\newcommand{\R}{\mathrm{I\!R}}
\newcommand{\N}{\mathrm{I\!N}}

\newcommand{\PP}{\mathrm{I\!P}}

\newcommand{\Z}{\mathchoice {\hbox{$\sf\textstyle Z\kern-0.4em
Z$}}{\hbox{$\sf\textstyle Z\kern-0.4em Z$}}{\hbox{$\sf\scriptstyle
Z\kern-0.3em Z$}}{\hbox{$\sf\scriptscriptstyle Z\kern-0.2em Z$}}}

\newcommand{\Q}{\mathchoice {\setbox0=\hbox{$\displaystyle\rm
Q$}\hbox{\raise0.15\ht0\hbox to0pt{\kern0.4\wd0\vrule
height0.8\ht0\hss}\box0}}{\setbox0=\hbox{$\textstyle\rm
Q$}\hbox{\raise0.15\ht0\hbox to0pt{\kern0.4\wd0\vrule
height0.8\ht0\hss}\box0}}{\setbox0=\hbox{$\scriptstyle\rm
Q$}\hbox{\raise0.15\ht0\hbox to0pt{\kern0.4\wd0\vrule
height0.7\ht0\hss}\box0}}{\setbox0=\hbox{$\scriptscriptstyle\rm
Q$}\hbox{\raise0.15\ht0\hbox to0pt{\kern0.4\wd0\vrule
height0.7\ht0\hss}\box0}}}

\newcommand{\OO}{\mathchoice {\setbox0=\hbox{$\displaystyle\rm
O$}\hbox{\hbox to0pt{\kern0.4\wd0\vrule
height0.9\ht0\hss}\box0}}{\setbox0=\hbox{$\textstyle\rm O$}\hbox{\hbox
to0pt{\kern0.4\wd0\vrule
height0.9\ht0\hss}\box0}}{\setbox0=\hbox{$\scriptstyle\rm O$}\hbox{\hbox
to0pt{\kern0.4\wd0\vrule
height0.9\ht0\hss}\box0}}{\setbox0=\hbox{$\scriptscriptstyle\rm
O$}\hbox{\hbox to0pt{\kern0.4\wd0\vrule height0.9\ht0\hss}\box0}}}

\newcommand{\GL}{\mathrm{GL}}
\newcommand{\SL}{\mathrm{SL}}
\newcommand{\SO}{\mathrm{SO}}

\newcommand{\eps}{\varepsilon}
\newcommand{\vi}{\varphi}
\newcommand{\vkap}{\varkappa}

\newcommand{\qmq}[1]{\quad\mbox{#1}\quad}
\newcommand{\Menge}[2]{\{\,#1\,|\,#2\,\}}

\newcommand{\arcosh}{\mathrm{arcosh}}

\newcommand{\vol}{\mathop{\mathrm{vol}}\nolimits}
\newcommand{\sign}{\mathop{\mathrm{sign}}\nolimits}

\newcommand{\SU}{\mathrm{SU}}

\newcommand{\RE}{\mathop{\mathrm{Re}}\nolimits}
\newcommand{\IM}{\mathop{\mathrm{Im}}\nolimits}

\newcommand{\tr}{\mathop{\mathrm{tr}}\nolimits}

\newcommand{\ord}{\mathop{\mathrm{ord}}\nolimits}
\newcommand{\polord}{\mathop{\mathrm{polord}}\nolimits}

\newcommand{\Res}{\mathrm{Res}}
\newcommand{\As}{\mathrm{As}}

\newcommand{\Mengegr}[2]{\{\,#1\,{\bigr |}\,#2\,\}}

\newcommand{\wt}{\widetilde}
\newcommand{\wh}{\widehat}
\newcommand{\calD}{\mathcal{D}}
\newcommand{\calM}{\mathcal{M}}
\newcommand{\calO}{\mathcal{O}}
\newcommand{\calR}{\mathcal{R}}

\newcommand{\C}{\mathchoice {\setbox0=\hbox{$\displaystyle\rm
C$}\hbox{\hbox to0pt{\kern0.4\wd0\vrule
height0.95\ht0\hss}\box0}}{\setbox0=\hbox{$\textstyle\rm C$}\hbox{\hbox
to0pt{\kern0.4\wd0\vrule
height0.95\ht0\hss}\box0}}{\setbox0=\hbox{$\scriptstyle\rm C$}\hbox{\hbox
to0pt{\kern0.4\wd0\vrule
height0.95\ht0\hss}\box0}}{\setbox0=\hbox{$\scriptscriptstyle\rm
C$}\hbox{\hbox to0pt{\kern0.4\wd0\vrule height0.95\ht0\hss}\box0}}}

\newcommand{\unity}{{1\!\!\!\:\mathrm{l}}}
\newcommand{\Pot}{\mathsf{Pot}}
\newcommand{\Mon}{\mathsf{Mon}}
\newcommand{\Div}{\mathsf{Div}}
\newcommand{\Jac}{\mathrm{Jac}}

\begin{document}
\begin{titlepage}
{
\renewcommand{\baselinestretch}{2.0}
\LARGE
\begin{center}
\textbf{A spectral theory  \\ for simply periodic solutions \\ of the sinh-Gordon equation}
\vspace{2cm}

\Large
Habilitationsschrift \\
zur Erlangung der venia legendi \\
im Fach Mathematik \\
an der Universit\"at Mannheim 

\vspace{2cm}
vorgelegt von \\
\textbf{\LARGE Sebastian Klein} \\
aus K\"oln
\vspace{2cm}

{\normalsize 2015} \\
\end{center}
}
\vfill
\clearpage
\thispagestyle{empty}
\strut\vfill

\end{titlepage}
\addtocounter{page}{-1}

\title[Spectral theory for simply periodic solutions of sinh-Gordon]{A spectral theory \\ for simply periodic solutions \\ of the sinh-Gordon equation}

\author[S.~Klein]{Sebastian Klein}

\address{\newline Universit\"at Mannheim, 68131 Mannheim, Germany}
\email{s.klein@math.uni-mannheim.de}




\abstract
{In this work a spectral theory for 2-dimensional, periodic, complex-valued solutions \,$u$\, of the sinh-Gordon equation is developed. Spectral data
for such solutions are defined (following \textsc{Hitchin} and \textsc{Bobenko}) 
and the space of spectral data is described by an asymptotic characterization. Using methods of asymptotic estimates,
the inverse problem for the spectral data is solved along a line, i.e.~the solution \,$u$\, is reconstructed on a line from the spectral data. Finally 
a Jacobi variety and Abel map for the spectral curve is constructed; they are used to describe the change of the spectral data under translation of the solution \,$u$\,.}
\endabstract

\maketitle

\bigskip

\tableofcontents

\part{Introduction}

\section{Introduction}
\label{Se:intro}

The objective of the present work is to develop a spectral theory for periodic, complex-valued solutions \,$u: X \to \C$\, of the 2-dimensional (i.e.~\,$X\subset \C$\,) sinh-Gordon equation
$$ \Delta u + \sinh(u)=0 \;. $$
We call such solutions \emph{simply periodic} when we wish to emphasize the difference to \emph{doubly periodic} solutions (which have two linear independent
periods).

We first discuss the importance of the sinh-Gordon equation.
One of the most significant reasons why the sinh-Gordon equation is interesting (apart from its relation to soliton theory)
is that solutions of the sinh-Gordon equation arise from minimal surfaces resp.~constant mean curvature surfaces (CMC surfaces) without umbilical points
(=zeros of the Hopf differential) in
the 3-dimensional real space forms: 
the Euclidean space \,$\R^3$\,, the 3-sphere \,$S^3$\, and the hyperbolic 3-space \,$H^3$\, (in \,$H^3$\, only if the absolute value of the mean curvature is \,$>1$\,).
For example, as we will describe more explicitly in Section~\ref{Se:minimal}, a conformal, minimal immersion \,$f: X \to S^3$\, 
without umbilical points can be reparameterized such that the Hopf differential is equal to \,$E\,\mathrm{d}z^2$\, with a constant value \,$E \in S^1$\,; in this situation
the Gauss equation for the conformal metric \,$g=e^{u/2}\,\mathrm{d}z\,\mathrm{d}\overline{z}$\, with the conformal factor \,$u: X \to \R$\, of \,$f$\, reduces to the sinh-Gordon equation, thus \,$u$\, is a solution of the sinh-Gordon equation. 
Moreover the Codazzi equation reduces to \,$0=0$\,. 
Therefore studying minimal surfaces or CMC surfaces in the 3-dimensional space forms means to study solutions of the sinh-Gordon equation.

Simple examples of CMC surfaces are surfaces with parallel second fundamental form, in \,$\R^3$\, these are the planes, spheres and (circular) cylinders (embedded in the usual way). 
More interesting examples of CMC surfaces in \,$\R^3$\, are the Delaunay surfaces, a class of immersed CMC surfaces which are topologically cylinders and which were first considered by \textsc{Delaunay} in 
\cite{Delaunay:1841}. \textsc{Alexandrov} showed in \cite{Alexandrov:1958} that the only \emph{embedded} compact CMC surfaces in \,$\R^3$\, are the round 2-spheres, and \textsc{Hopf} showed in 
\cite{Hopf:1983} that the only immersed compact CMC surfaces in \,$\R^3$\, of genus \,$g=0$\, (i.e.~homeomorphic to the 2-sphere) are again the round 2-spheres; he in fact conjectured that there are no
immersed compact CMC surfaces in \,$\R^3$\, besides the round 2-spheres at all. Interest in compact immersed CMC surfaces in the 3-dimensional space forms soared, however, after \textsc{Wente} showed in
1986 (see \cite{Wente:1986}) this conjecture to be false by constructing what is now known as the Wente tori: a family of compact, immersed CMC surfaces in \,$\R^3$\, of genus \,$g=1$\,. Later, \textsc{Kapouleas}
found compact immersed CMC surfaces in \,$\R^3$\, at first for every genus \,$g\geq 3$\, (see \cite{Kapouleas:1987}), and then also for \,$g=2$\, (see \cite{Kapouleas:1990}, \cite{Kapouleas:1991}). 
A family of compact minimal surfaces in \,$S^3$\, is the family of Lawson surfaces: they exist for every genus \,$g\geq 1$\,, have a large symmetry group and were obtained by \textsc{Lawson} by applying a
reflection principle to a suitably chosen solution of Plateau's Problem (of finding a minimal surface with prescribed boundary) in \cite{Lawson:1968}. 

A celebrated result was the classification of all immersed CMC tori in \,$\R^3$\, resp.~\,$S^3$\, by \textsc{Pinkall} and \textsc{Sterling} in \cite{Pinkall/Sterling:1989} 
and independently by \textsc{Hitchin} in \cite{Hitchin:1990}. The classification by Pinkall and Sterling is based on proving that CMC tori can be described by polynomial Killing fields; 
a summary of their proof is found in \cite{Helein:2001}, Chapter~9.
While also Pinkall's and Sterling's classification can be interpreted in terms of spectral theory, Hitchin's classification is based on spectral theory in a more explicit way, and we will discuss the latter
classification in some more detail below.



As stated at the beginning, the present work sets out to develop a spectral theory for simply periodic solutions of the sinh-Gordon equation.
The term ``spectral theory'' here refers more precisely to a scheme of studying solutions of a given differential operator by looking at the spectrum of an associated Lax operator.
This scheme was first developed for the Korteweg-de Vries equation (KdV equation).
In this introduction, we will first describe the scheme as it applies to the KdV equation, and then how it can be applied to the sinh-Gordon equation. The study of simply periodic solutions of the sinh-Gordon equation by this method
is the main topic of the present work.

The KdV equation is the partial differential equation \,$u_t = -u_{xxx}+6\,u\,u_x$\, (which was first studied in the 19th century in relation to the modeling of waves on shallow water channels). 
It is an important insight that is due to \textsc{Lax} (see \cite{Lax:1975})  that the KdV equation can be written in the form of what is now known as a \emph{Lax equation}, i.e.~in the form 
\,$\tfrac{\mathrm{d}\ }{\mathrm{d}t} L = [B,L]$\, (where \,$[B,L]$\, is the commutator \,$B\cdot L-L\cdot B$\,) with respect to two other differential operators \,$L$\, and \,$B$\,. 
In this situation, the pair \,$(L,B)$\, is called a \emph{Lax pair}. More specifically, for the KdV equation \,$L=-\tfrac{\mathrm{d}^2\ }{\mathrm{d}x^2} + u$\, is the
1-dimensional Schr\"odinger operator, and \,$B=-4 \tfrac{\mathrm{d}^3\ }{\mathrm{d}x^3} + 3u\,\tfrac{\mathrm{d}\ }{\mathrm{d}x} + 3\,\tfrac{\mathrm{d}\ }{\mathrm{d}x}\,u$\,
is another differential operator of third order. Note that if \,$u$\, is a periodic solution of the KdV equation, then the Lax operators \,$L$\, and \,$B$\, are also periodic. 

The idea then is to study periodic solutions \,$u$\, of the KdV equation via the periodic spectrum of the associated Lax operator \,$L$\,. It is known that this spectrum is a pure point spectrum,
thus the spectrum is defined by the eigenvalue equation \,$L\psi = \lambda\psi$\,, where \,$\lambda\in \C$\, and \,$\psi$\, is a periodic function. 
It is a consequence of the Lax equation that the spectrum of \,$L$\, is invariant under the flows of \,$L$\, and \,$B$\, (and also under an infinite family of ``higher flows'').
The periodic spectrum does not determine \,$u$\, uniquely, however one can define additional data associated to the periodic eigenfunctions \,$\psi$\, (which move under the flows of \,$L$\, and \,$B$\,),
so that the spectrum along with these additional data determine the original solution \,$u$\, of the KdV equation uniquely. The spectrum together with these additional data are called the
\emph{spectral data} of \,$u$\,. 

In relation to the correspondence between the solution \,$u$\, and the spectral data there are two fundamental problems: first the \emph{direct problem} to study the spectral data
corresponding to a solution \,$u$\, and to derive properties of
the spectral data from properties of the solution \,$u$\,; and second, the \emph{inverse problem} to study \,$u$\, by its spectral data, especially
to show that \,$u$\, is determined uniquely by the spectral data and to reconstruct \,$u$\, from its spectral data. 

Note that this approach to the study of periodic solutions \,$u$\, of the KdV equation depends almost exclusively on the study of the spectrum of the 1-dimensional Schr\"odinger operator \,$L$\, with potential \,$u$\, along a period (i.e.~only
with respect to the single variable \,$x$\,).
The other Lax operator \,$B$\, comes into play only when one reconstructs the two-dimensional solution \,$u$\, of the KdV equation from its values along the \,$x$-line by applying the flow of \,$B$\,. 

The Lax equation can be interpreted as a zero curvature condition in the following way:
The second order differential operator \,$L$\, can be rewritten as the 2-dimensional first order differential operator 
\,$-\tfrac{\mathrm{d}\ }{\mathrm{d}x}+\left( \begin{smallmatrix} 0 & 1 \\ u & 0 \end{smallmatrix} \right)$\,
and then the eigenvalue equation \,$L\psi = \lambda\psi$\, takes the form
\begin{equation}
\frac{\mathrm{d}\ }{\mathrm{d}x} \Psi = U_\lambda\Psi
\tag{$*$}
\end{equation}
with the 2-dimensional periodic operator
$$ U_\lambda := \begin{pmatrix} 0 & 1 \\ u-\lambda & 0 \end{pmatrix} \qmq{and} \Psi = \begin{pmatrix} \psi_1 \\ \psi_2 \end{pmatrix} = \begin{pmatrix} \psi \\ \tfrac{\mathrm{d}\ }{\mathrm{d}x}\psi \end{pmatrix} \; . $$
Moreover there exists another 2-dimensional periodic operator \,$V_\lambda$\, also depending on the parameter \,$\lambda$\, so that 
the Lax equation \,$\tfrac{\mathrm{d}\ }{\mathrm{d}t} L = [B,L]$\, is equivalent to the zero curvature condition
\begin{equation}
\frac{\partial U_\lambda}{\partial t} - \frac{\partial V_\lambda}{\partial x} + [U_\lambda,V_\lambda]=0 \; . 
\tag{$\dagger$}
\end{equation}

Results on the study of the KdV equation via the spectral theory for the associated 1-dimensional Schr\"odinger operator in this way are summarized for example in \cite{Kappeler/Poeschel:2003}, Chapter~III. 
The spectral theory for the 1-dimensional Schr\"odinger operator was investigated by 
\textsc{McKean} and \textsc{van Moerbeke} in \cite{McKean/Moerbeke:1975}; this work was later extended
by \textsc{McKean} and \textsc{Trubowitz} to include also the study of the Jacobi variety of the spectral curve in \cite{McKean/Trubowitz:1976}.
Representative for many further results on this topic, I mention only  \cite{Mueller/Schmidt/Schrader:1998}
by \textsc{M\"uller}, \textsc{Schmidt} and \textsc{Schrader}, where certain quasi-periodic solutions are constructed by means of 
complex analysis on hyperelliptic spectral curves of infinite genus.
A very accessible introduction to the spectral theory of the 1-dimensional Schr\"odinger equation is the book \cite{Poeschel-Trubowitz:1987},
which however considers only solutions \,$u$\, of the Schr\"odinger equation with the Dirichlet boundary conditions \,$y(0)=y(1)=0$\,
(choosing what David Hilbert would have called ``that special case which contains all the germs of generality''). The cited papers and books, but especially
\cite{Poeschel-Trubowitz:1987}, have been very inspirational for the treatment of the sinh-Gordon equation in the present work, and several
of the results in the present work correspond to analogous theorems in the theory of the Schr\"odinger equation; the correspondences will be
pointed out.

Other examples of differential equations that have been studied successfully by means of the spectral theory of a corresponding Lax operator \,$L$\, 
include the modified Korteweg-de Vries equation (mKdV equation) and the non-linear Schr\"odinger equation.

In this work we study the sinh-Gordon equation via periodic 2-dimensional operators \,$U_\lambda$\,, \,$V_\lambda$\, associated to a periodic solution \,$u$\, of the sinh-Gordon equation
that satisfy the zero curvature condition $(\dagger)$; from the monodromy of the associated
differential equation $(*)$ we will obtain data that are analogous to the spectral data described above for the 1-dimensional Schr\"odinger operator. 
For a given solution \,$u$\, of the sinh-Gordon equation, the operators \,$U_\lambda$\,, \,$V_\lambda$\, are obtained from the Maurer-Cartan form
\,$\alpha_\lambda=U_\lambda\,\mathrm{d}x + V_\lambda\,\mathrm{d}y$\, of the extended frame 
\,$F_\lambda$\, of the minimal immersion \,$f: X\to S^3$\, corresponding to \,$u$\,. However the dependence of the operators \,$U_\lambda$\, and \,$V_\lambda$\, on \,$\lambda$\, is more complicated
in this setting; consequently the equation $(*)$ can no longer be interpreted as the eigenvalue equation for a differential operator \,$L$\, as was the case for the KdV equation, and thus also
the parameter \,$\lambda$\, can no longer be interpreted as a eigenvalue resp.~a spectral value. Nevertheless, the analogy to the situation for the KdV equation is very close via the interpretation
of the Lax equation as a zero curvature condition, which is the reason
why I take the liberty of applying the term \emph{spectral theory} also to the study of periodic solutions of the sinh-Gordon equation by this method; likewise we will use the adjective ``spectral'' for
the parameter \,$\lambda$\, and similar objects.

This concept of spectral theory has first been applied to the sinh-Gordon equation by \textsc{Hitchin} in \cite{Hitchin:1990} where he studied doubly periodic solutions in relation with his classification of the minimal
tori in \,$S^3$\,. 
These concepts have later been refined by \textsc{Bobenko} and adapted to CMC immersions in all the 3-dimensional space forms. 
\textsc{McIntosh} \cite{McIntosh:2008} has described spectral data from the point of view of considering harmonic maps from a torus into \,$S^2$\,
(such maps correspond to minimal immersions of a torus into \,$S^3$\, by taking the Gauss map). In the present text, we will use this kind of spectral theory to study simply periodic solutions of the sinh-Gordon
equation.

We now describe the construction of spectral data for the sinh-Gordon equation explicitly using the terminology of Bobenko, which we will also use in the present work, and then Hitchin's method for the classification.

Consider a minimal conformal immersion \,$f:X \to S^3$\, without umbilical points into the 3-sphere, corresponding to a solution \,$u:X\to\R$\, of the sinh-Gordon equation.
Then the frame \,$\underline{F}: X \to \SO(4)$\, of \,$f$\,
induces the \,$\mathfrak{so}(4)$-valued flat connection form \,$\underline{\alpha}=\underline{F}^{-1}\cdot \mathrm{d}\underline{F}$\, on \,$X$\, (which is the Maurer-Cartan form of \,$\underline{F}$\,).
In fact, minimal immersions into \,$S^3$\, come in families \,$f_\lambda$\,
with \,$\lambda\in S^1$\, (the \,$f_\lambda$\, differ only in their Hopf differential \,$E\,\mathrm{d}z^2$\,), 
and thus we obtain a corresponding family of frames \,$\underline{F}_\lambda$\, and of flat connection forms \,$\underline{\alpha}_\lambda$\,. 
If we identify \,$\SO(4)$\, with \,$(\SU(2)\times \SU(2))/\Z_2$\, and therefore \,$\mathfrak{so}(4)$\,
with \,$\mathfrak{su}(2) \oplus \mathfrak{su}(2)$\,, then the connection form \,$\underline{\alpha}_\lambda$\, is of the form
\,$\underline{\alpha}_{\lambda}=(\alpha_\lambda,\alpha_{-\lambda})$\, with an \,$\mathfrak{su}(2)$-valued flat connection form \,$\alpha_\lambda$\,.
Moreover, \,$\alpha_\lambda$\, can be written down explicitly in terms of \,$u$\, and \,$\lambda$\,, also for \,$\lambda\in \C^*$\, (instead of only for \,$\lambda\in S^1$\,), and also for complex-valued solutions \,$u$\,
(which do not correspond to immersions, then the Lie groups \,$\SO(4)$\, and \,$\SU(2)$\, are replaced by their complex forms \,$\SO(4,\C)$\, and \,$\SL(2,\C)$\,).
The \,$\mathrm{d}x$-part and the \,$\mathrm{d}y$-part of \,$\alpha_\lambda$\, give the \,$\lambda$-family of Lax pairs for our solution \,$u$\, of the sinh-Gordon equation as explained above.
We integrate \,$\alpha_\lambda$\,
by solving the differential equation \,$\mathrm{d}F_\lambda=\alpha_\lambda\cdot F_\lambda$\, for \,$F_\lambda$\, along a period of \,$u$\,, and consider the monodromy \,$M(\lambda) := F_\lambda(z_1) \cdot F_\lambda(z_0)^{-1}$\,
along a period. Then the hyperelliptic complex curve (possibly with singularities) defined by the eigenvalues of \,$M(\lambda)$\,
$$ \Sigma := \Mengegr{(\lambda,\mu)\in \C^* \times \C}{\det(M(\lambda)-\mu\cdot \unity)=0} $$
will be called the \emph{spectral curve} of \,$u$\, in the present work (but is called the multiplier curve by Hitchin, who reserves the term ``spectral curve'' for a different object, see below), and the
divisor \,$D$\, defining the eigenline bundle of \,$M(\lambda)$\, on \,$\Sigma$\, is the \emph{spectral divisor} of \,$u$\,. \,$(\Sigma,D)$\, is called the \emph{spectral data} for the solution \,$u$\, of the sinh-Gordon equation.
The branch points of \,$\Sigma$\, are analogous to the periodic spectrum of the 1-dimensional Schr\"odinger operator \,$L$\, in the treatment of the KdV equation via spectral theory described above, and 
the spectral data \,$(\Sigma,D)$\, are analogous to the spectral data for \,$L$\,. 

Hitchin's objective is to classify the minimal immersed tori of \,$S^3$\,, corresponding to the case where \,$u$\, is doubly periodic. The two minimal periods of \,$u$\, induce two
different monodromies and therefore Hitchin obtains at first two different sets \,$(\Sigma,D)$\, of spectral data for \,$u$\,. However, because the homotopy group of the torus is abelian,
the two monodromies commute, and therefore the two curves \,$\Sigma$\, are in fact the same. Hitchin constructs a complex curve \,$\wt{\Sigma}$\, on which the eigenlines of the monodromies
actually are holomorphic line bundles, even at the branch points or singularities of \,$\Sigma$\,; he obtains \,$\wt{\Sigma}$\, as a partial desingularization of \,$\Sigma$\,,
namely by desingularizing \,$\Sigma$\, at all those singularities which are also divisor points. It turns out that the resulting curve \,$\wt{\Sigma}$\, has finite geometric genus,
and can therefore be compactified by adding points above \,$\lambda=0$\, and \,$\lambda=\infty$\,. 
This is what permits Hitchin to show that doubly periodic solutions \,$u$\, are of finite type, and therefore ultimately to
classify the doubly periodic solutions of the sinh-Gordon equation, by applying classical results of the theory of compact Riemann surfaces.

We mention that \textsc{Heller} has applied Hitchin's construction of spectral data to compact (closed) immersed surfaces of genus \,$g\geq 2$\, in \,$S^3$\,;
he obtains the most interesting results for surfaces which are ``Lawson symmetric'', i.e.~which have the symmetry group of one of the Lawson surfaces;
he also obtained CMC deformations of such surfaces in \,$S^3$\,. See for example \cite{Heller:2013}, \cite{Heller:2014} and \cite{Heller/Schmitt:2015}.

In the present work, we develop the theory of spectral data for solutions \,$u:X \to \C$\, of the sinh-Gordon equation 
(with \,$X\subset \C$\,)
which have only a single period (are \emph{simply periodic}). Note that we include the case of complex-valued \,$u$\, throughout the work.
By a linear transformation of \,$X$\, we may suppose without loss of generality that \,$0\in X$\, holds and that the period of \,$u$\,
is \,$1$\,. Then \,$X$\, is a horizontal strip in \,$\C$\, and 
$$ u(z+1) = u(z) \qmq{holds for all \,$z\in X$\,.} $$
It is possible to construct the spectral data also for bare Cauchy data \,$(u,u_y)$\, for a solution of the sinh-Gordon equation,
here \,$u$\, and \,$u_y$\, are periodic functions defined only on the interval \,$[0,1]$\,, where \,$u_y$\, gives the derivative
of \,$u$\, in normal direction. This is the viewpoint we will take throughout most of the work. 

In this context we are interested in requiring only as weak regularity conditions for \,$(u,u_y)$\, as possible. 
There are two reasons: First, we are interested in characterizing precisely which divisors on a spectral curve are spectral divisors of some Cauchy data \,$(u,u_y)$\,;
it turns out that every additional differentiability condition imposed on \,$(u,u_y)$\, reduces the space of divisors by an intricate relationship between its divisor points.
By not imposing more regularity than necessary, we obtain a description of the space of divisors that is as simple as possible. 
Second, while any solution \,$u$\, of the sinh-Gordon equation is infinitely differentiable (in fact even real analytic, because the sinh-Gordon equation is elliptic)
on the interior of its domain, we are also interested in the behavior of the solution on the boundary of its domain, where its behavior can be worse. 
For these reasons we only require \,$(u,u_y) \in W^{1,2}([0,1]) \times L^2([0,1])$\,. 

The definition of the spectral data \,$(\Sigma,D)$\, described above carries over to our situation without difficulty. 
As spectral curve, we use the complex curve
\,$\Sigma$\, defined by the eigenvalues of the monodromy (called the multiplier curve by Hitchin), not the partial desingularization \,$\wt{\Sigma}$\, that
Hitchin calls the spectral curve. There are several reasons for this choice: 

The most important reason is that I take the point of view of describing variations of \,$u$\, within its spectral class by variation of the spectral divisor,
not by change of the curve on which the divisor is defined. Because Hitchin's partial desingularization \,$\wt{\Sigma}$\, needs to be chosen in dependence of the spectral
divisor, spectral divisors on the same complex curve \,$\Sigma$\, can induce different partial desingularizations \,$\wt{\Sigma}$\,. Thus using the
partial desingularization would complicate the approach of working on a fixed curve.

Another reason is that in our setting (with \,$u$\, only simply periodic) Hitchin's partial desingularization \,$\wt{\Sigma}$\, would still have infinite geometric genus in general,
so passing from \,$\Sigma$\, to \,$\wt{\Sigma}$\, is not quite as useful in the present situation as it is for Hitchin. Moreover we note that the existence of a partial
desingularization of \,$\Sigma$\, which parameterizes the eigenline bundle of the monodromy is a property of only the sinh-Gordon integrable system, as a consequence of the fact
that the monodromy has values in \,$(2\times 2)$-matrices, which causes \,$\Sigma$\, to be hyperelliptic. 
In the spectral theory for other integrable systems, the monodromies generally have values in larger endomorphism spaces,
hence \,$\Sigma$\, is no longer hyperelliptic (instead it covers the \,$\lambda$-plane with order \,$\geq 3$\,), 
and then there does not exist such a partial desingularization, see Remark~\ref{R:spectrum:locally-free-interpretation}.

However, there is one complication that is caused by using the curve \,$\Sigma$\, as spectral curve, instead of its partial desingularization, and that is that we need to deal
with the fact that it is possible for singularities of \,$\Sigma$\, to occur in the support of the spectral divisor. In this case, the divisor, in the classical sense,
does not convey enough information to determine the monodromy \,$M(\lambda)$\, uniquely, which is of course what we expect of the spectral data. The solution to
this problem is to generalize the concept of a divisor following a concept of \textsc{Hartshorne} \cite{Hartshorne:1986}. In this sense, a \emph{generalized divisor}
is a subsheaf \,$\calD$\, of the sheaf \,$\calM$\, of meromorphic functions on \,$\Sigma$\, which is finitely generated over
the sheaf \,$\calO$\, of holomorphic functions on \,$\Sigma$\,. Such a generalized divisor contains exactly the same information as a classical divisor
in the regular points of \,$\Sigma$\,, but it conveys more information in the singular points. It is the proper concept of a divisor
to ensure that the spectral curve \,$\Sigma$\, together with a generalized divisor \,$\calD$\, uniquely determines the monodromy \,$M(\lambda)$\,,


In our setting, the geometric genus of the spectral curve \,$\Sigma$\, is in general infinite -- the branch points accumulate near \,$\lambda=\infty$\, and near \,$\lambda=0$\, -- so
\,$\Sigma$\, can no longer be compactified. For this reason, the classical results on compact Riemann surfaces which played a significant role
for the study of solutions of finite type, are no longer applicable to \,$\Sigma$\, in the case of infinite type. 

We need to replace these standard results with arguments for non-compact surfaces. To make such arguments feasible, one needs to prescribe
the asymptotic behavior of the spectral curve \,$\Sigma$\, ``at its ends'', i.e.~near the points \,$\lambda=\infty$\, and \,$\lambda=0$\,. 
The information on the asymptotic behavior of \,$\Sigma$\, is the replacement for the fact from the finite type theory that \,$\Sigma$\, can be compactified
at \,$\lambda=\infty$\, and \,$\lambda=0$\,. It is for this reason that asymptotic estimates for the quantities involved in the construction
of the spectral data are a very important, perhaps the most important tool for the conduction of the present study.

Unfortunately there are only very few results on open Riemann surfaces with prescribed asymptotics found in the literature.
One example would be the book \cite{Feldman/Knoerrer/Trubowitz:2003} by \textsc{Feldman/Kn\"orrer/Trubowitz}. From this book, we will indeed use one of
the results from Chapter~1 of the book that apply to general surfaces that are parabolic in the sense of \textsc{Ahlfors} and \textsc{Nevanlinna}. 
However the results in the later part of the book,
which would be very useful to us, depend on very strict geometric hypotheses for the surface under consideration, see Section~5 of \cite{Feldman/Knoerrer/Trubowitz:2003}, and especially
the hypothesis GH2(iv), which requires that the ``size of the handles'' \,$t_j$\, (comparable to the distance between the pairs of branch points) satisfies \,$\sum t_j^\beta < \infty$\, for every
\,$\beta>0$\,. This requirement implies in particular that the \,$t_j$\, are decreasing more rapidly than \,$j^{-n}$\, for every \,$n\geq 1$\,.
Because of our relatively mild requirements on the differentiability of \,$(u,u_y)$\,, the \,$t_j$\, do not decrease sufficiently rapidly by far
for our spectral curve \,$\Sigma$\,, and thus the results from the latter part of \cite{Feldman/Knoerrer/Trubowitz:2003} are not applicable 
in our situation. For this reason we develop some results analogous to classical results on compact Riemann surfaces in this work, as needed.




\smallskip

We now give an overview of the results of the work. In Section~\ref{Se:minimal} we explicitly describe the relationship between solutions \,$u$\,
of the sinh-Gordon equation, and minimal immersions \,$f: X \to S^3$\, without umbilical points. This relationship is of importance because it
provides the \,$\mathfrak{sl}(2,\C)$-valued
flat connection form \,$\alpha_\lambda$\, (depending on the spectral parameter \,$\lambda\in \C^*$\,) and thereby the extended frame
\,$F_\lambda: X \to \SL(2,\C)$\, 
as the solution of the differential equation \,$\mathrm{d}F_\lambda=\alpha_\lambda\,F_\lambda$\,. In the case where \,$u$\, is simply periodic,
the monodromy \,$M(\lambda) \in \SL(2,\C)$\, of the extended frame along the period is then used in Section~\ref{Se:spectrum} to construct the spectral data 
\,$(\Sigma,D)$\, for \,$u$\,. Here \,$\Sigma$\, is the spectral curve, defined by the eigenvalues of \,$M(\lambda)$\,, which turns out to be 
a hyperelliptic complex curve, and \,$D$\, is the divisor on \,$\Sigma$\, defining the eigenline bundle of \,$M(\lambda)$\,. With respect to \,$D$\,,
a difficulty arises in view of the fact that the spectral curve \,$\Sigma$\, can be singular: If the support of \,$D$\, contains singular points,
then the classical divisor does not convey enough information at these points to uniquely determine the monodromy \,$M(\lambda)$\, anymore,
and thus we replace \,$D$\, with a generalized divisor \,$\calD$\, in the sense of \textsc{Hartshorne} as explained above.
We will see in Section~\ref{Se:special} that \,$(\Sigma,\calD)$\, uniquely determines the monodromy \,$M(\lambda)$\, (up to a choice of sign).

The first main part of the work then is to describe the asymptotic behavior of the monodromy \,$M(\lambda)$\, and thereby of the spectral
data \,$(\Sigma,\calD)$\, near \,$\lambda=\infty$\, and near \,$\lambda=0$\,. As explained above, this is the most important instrument for understanding
and controlling the spectral curve in our situation.

We begin in Section~\ref{Se:vacuum} by calculating the spectral data for the ``vacuum solution'' \,$u=0$\, of the sinh-Gordon equation.
This example is of importance because we will describe the asymptotic behavior of the various objects near \,$\lambda\to\infty$\, and \,$\lambda\to0$\,
by comparing the objects for any given simply periodic solution \,$u$\, to the corresponding (explicitly calculated) object for the 
vacuum solution. In this way, one can view the spectral data for the vacuum solution as being ``typical'' for the spectral data
for any simply periodic solution, where the extent of typical-ness is quantified precisely by the asymptotic estimates of the sections to come.

For the investigation of the asymptotic behavior of \,$M(\lambda)$\, and of the spectral data \,$(\Sigma,\calD)$\,, we do not consider
a solution \,$u$\, of the sinh-Gordon equation defined on a horizontal strip in \,$\C$\, of positive height, but rather 
periodic Cauchy data \,$(u,u_y)$\, defined on \,$[0,1]$\,. For the reasons explained above, we only require
\,$(u,u_y) \in W^{1,2}([0,1])\times L^2([0,1])$\,; the condition of periodicity then reduces to \,$u(0)=u(1)$\,. 
Note that for many of the asymptotic estimates we also admit Cauchy data  without the condition of periodicity.
Such ``non-periodic potentials'' of course do not correspond to simply-periodic solutions of the sinh-Gordon equation;
they are considered for the sole purpose of deriving the asymptotic behavior of the extended frame \,$F_\lambda$\, in Proposition~\ref{P:asympfinal:frame}.

The strategy for obtaining asymptotic estimates for the spectral data is to find asymptotic estimates for the monodromy \,$M(\lambda)$\, first
and then to derive asymptotic estimates for the spectral data from them. Similarly as in the case of the 1-dimensional Schr\"odinger equation
(compare \cite{Poeschel-Trubowitz:1987}, Chapter~2), one needs to carry out this process in two stages: First one derives a relatively 
mild asymptotic for the monodromy (called the ``basic asymptotic'' in Section~\ref{Se:asymp}), which then yields first information on the 
spectral data (Section~\ref{Se:excl}): It turns out that the points in the support of the spectral divisor \,$\calD$\, can be enumerated by a sequence 
\,$(\lambda_k,\mu_k)_{k\in \Z}$\, of points on \,$\Sigma$\,, where for \,$|k|$\, large, each \,$(\lambda_k,\mu_k)$\, is in a certain neighborhood of a corresponding
point \,$(\lambda_{k,0},\mu_{k,0})$\, in the divisor of the vacuum; we call these neighborhoods the \emph{excluded domains} \,$\wh{U}_{k,\delta}$\,. 
Note that that the divisor points are (asymptotically) free to move in the excluded domains, they are not restricted to intervals or curves
(as it is the case, for example, for real solutions of the KdV equation, where these intervals are called ``gaps''); the reason is that
we are considering complex, not only real, solutions of the sinh-Gordon equation.

Similarly, \,$\Sigma$\, has two branch points in each \,$\wh{U}_{k,\delta}$\, for \,$|k|$\, large; if they coincide, then a singularity of \,$\Sigma$\,
occurs. This singularity is an ordinary double point; for \,$|k|$\, large this is the only type of singularity that can occur,
but for small \,$|k|$\,, the spectral curve \,$\Sigma$\, can also have singularities of higher order. 

These asymptotic assessments can still be refined: By studying the Fourier transform of the given Cauchy data \,$(u,u_y)$\,,
in Section~\ref{Se:fasymp} it is shown that the difference between the monodromy \,$M(\lambda)$\, and the monodromy \,$M_0(\lambda)$\,
follows for \,$\lambda=\lambda_{k,0}$\, a law of square-summability, which we call the ``Fourier asymptotic'' for the monodromy.
From this fact, it follows that the divisor points and the branch points
also follow a square-summability law with respect to the difference to their counterparts for the vacuum. We show this for the divisor
points in Section~\ref{Se:asympdiv}, whereas the corresponding result for the branch points of the spectral curve is postponed to Section~\ref{Se:asympfinal}
to avoid technical difficulties.

Our next objective is to solve the inverse problem for the monodromy: Suppose that spectral data \,$(\Sigma,\calD)$\, are given, where
the support of \,$\calD$\, satisfies the square-summability law just mentioned, then one would like to reconstruct the monodromy \,$M(\lambda)$\,
from which these spectral data were induced; in particular, we need to show that \,$M(\lambda)$\, is determined uniquely by the spectral data. 
If we write the \,$\SL(2,\C)$-valued monodromy \,$M(\lambda)$\, in the form \,$M(\lambda)=\left( \begin{smallmatrix} a(\lambda) & b(\lambda) \\
c(\lambda) & d(\lambda) \end{smallmatrix} \right)$\, with holomorphic functions \,$a,b,c,d: \C^*\to\C$\,, the spectral divisor directly conveys
information on the zeros of \,$c$\, and on the function values of \,$d$\, (or \,$a$\,) at the zeros of \,$c$\,. The inverse problem for the monodromy
is, in essence, to reconstruct the functions \,$a,b,c,d$\, from these information and from their known asymptotic behavior.

To facilitate the discussion of this problem, we introduce in Section~\ref{Se:As} spaces of holomorphic functions on \,$\C^*$\, or on \,$\Sigma$\,
which have a prescribed asymptotic behavior near \,$\lambda=\infty$\, and/or near \,$\lambda=0$\,. The difference between the functions
comprising \,$M(\lambda)$\, and the corresponding functions of the vacuum will turn out to be contained in certain of these spaces, 
and it is for this reason that these spaces are extremely important in the sequel.

In Section~\ref{Se:interpolate} we address the mentioned problem of reconstructing a holomorphic function on \,$\C^*$\, with a prescribed
asymptotic behavior from the sequence of its zeros, or else from a sequence of its values at certain points. The first of these problems
(reconstruction from the zeros) is reminiscent of Hadamard's Factorization Theorem (see for example~\cite{Conway:1978}, Theorem~XI.3.4, p.~289), but here we have
a holomorphic function defined on \,$\C^*$\, with two accumulation points of the sequence of zeros (\,$\lambda=\infty$\, and \,$\lambda=0$\,),
instead of a function on \,$\C$\, with the zeros accumulating only at \,$\lambda=\infty$\,. For both of the reconstruction problems,
we obtain explicit descriptions of the solution function as an infinite product resp.~sum. For the proof of the corresponding statements,
we require several relatively elementary results on the convergence and estimate of infinite sums and products, these results are 
collected in Appendix~\ref{Ap:inf}.

By applying the explicit product resp.~sum formulas from Section~\ref{Se:interpolate} to the functions comprising the monodromy \,$M(\lambda)$\,,
we can improve the asymptotic description for \,$M(\lambda)$\, one final time: in the description of the Fourier asymptotic in Section~\ref{Se:fasymp},
the points \,$\lambda=\lambda_{k,0}$\, played a special role. In Section~\ref{Se:asympfinal} we are now liberated from this condition
and obtain a square-summability law for \,$M-M_0$\, with respect to arbitrary points e.g.~in excluded domains. We also obtain 
the asymptotic behavior of the branch points of the spectral curve \,$\Sigma$\,, and an asymptotic estimate for the extended frame
\,$F_\lambda$\, itself.

Finally in Section~\ref{Se:special} we characterize those generalized divisors \,$\calD$\, on a spectral curve \,$\Sigma$\,
for which \,$(\Sigma,\calD)$\, is the spectral data of a monodromy \,$M(\lambda)$\,. It turns out that besides the asymptotic behavior
that was described before, one needs a certain compatibility condition (Definition~\ref{D:special:special}(1)) and (as is to be expected)
that the divisor \,$\calD$\, is non-special. Under these circumstances it turns out that \,$(\Sigma,\calD)$\, corresponds to 
a monodromy \,$M(\lambda)$\, and that \,$M(\lambda)$\, is determined uniquely by \,$(\Sigma,\calD)$\, (up to a change of sign of the 
off-diagonal entries). This concludes the solution of the inverse problem for the monodromy.

The next question is the inverse problem for the Cauchy data \,$(u,u_y)$\, themselves, i.e.~to reconstruct \,$(u,u_y)$\, from 
the spectral data \,$(\Sigma,\calD)$\,. The corresponding problem has been solved for the potentials of finite type, and it appears
natural to try to apply that result for our infinite type potentials.

For this purpose we show in Section~\ref{Se:finite} (by a fixed point argument) that the set of finite type divisors is dense in the space of all asymptotic 
divisors. Moreover, in Section~\ref{Se:darboux} we equip the space of potentials with a symplectic form and construct Darboux coordinates
on this symplectic space, building on previous results by \textsc{Knopf} \cite{Knopf:2015}. Using these two results, we are then
able to show relatively easily in Section~\ref{Se:diffeo} that for the open and dense subset \,$\Pot_{tame}$\, of ``tame'' potentials,
the map \,$\Pot_{tame}\to\Div$\, associating to each tame potential (Cauchy data) the corresponding spectral divisor is 
a diffeomorphism onto an open and dense subset of the space of all asymptotic divisors. This essentially solves the inverse problem
for the Cauchy data \,$(u,u_y)$\,.

The final part of the work is concerned with constructing a Jacobi variety and an Abel map for the spectral curve \,$\Sigma$\,.
For the sake of simplicity, we will here restrict ourselves to spectral curves \,$\Sigma$\, which do not have any singularities other
than ordinary double points. Also in this part we admit spectral data corresponding to complex solutions of the sinh-Gordon equation,
it is for this reason that the individual Jacobi coordinates take values in \,$\C$\, (and not only in \,$S^1$\,, which would be familiar
for example from the Jacobi variety for (real) solutions of the 1-dimensional Schr\"odinger operator constructed in 
\cite{McKean/Trubowitz:1976}, Sections~11ff.).

The construction of the Jacobi variety needs to be prepared with more asymptotic estimates: In Section~\ref{Se:jacobiprep} we prove estimates for certain contour
integrals on the spectral curve \,$\Sigma$\,, and in Section~\ref{Se:holo1forms} we study holomorphic 1-forms on \,$\Sigma$\,.
In particular we construct holomorphic 1-forms which have zeros in all the excluded domains, with the exception of a single one, 
and find sufficient conditions for such 1-forms to be square-integrable on \,$\Sigma$\,.

With help of these results we will construct a canonical basis \,$(\omega_n)$\, for the space of holomorphic 1-forms explicitly in Section~\ref{Se:jacobi}.
In the case where the spectral curve \,$\Sigma$\, has no singularities, the existence of such a basis follows from the general theory
of parabolic (in the sense of Ahlfors/Nevanlinna) Riemann surfaces described for example in Chapter~1 of \cite{Feldman/Knoerrer/Trubowitz:2003}. However this is not enough for us,
not only because we also need to handle the case where \,$\Sigma$\, has singularities, but also because we derive asymptotic 
estimates for the \,$\omega_n$\, from their explicit representation which are stronger than those one can obtain abstractly;
we need these stronger estimates in the construction of the Jacobi variety and the Abel map for \,$\Sigma$\,. 

The work is concluded with two applications of this construction: In Section~\ref{Se:jacobitrans} we describe the motions in the
Jacobi variety that correspond to translations of the solution \,$u$\, of the sinh-Gordon equation. Similarly as is known for the
finite type setting, it turns out that these translations correspond to linear motions in the Jacobi variety. And finally,
in Section~\ref{Se:strip} we describe the asymptotic behavior of the spectral data \,$(\Sigma,D)$\, for an actual simply periodic solution
\,$u: X \to \C$\, of the sinh-Gordon equation, given on a horizontal strip of positive height \,$X \subset \C$\,. We see that
in this case the spectral data satisfy an exponential asymptotic law, much steeper than the asymptotic behavior of the spectral
data for Cauchy data \,$(u,u_y)$\,, as is to be expected, solutions of the sinh-Gordon equation being real analytic in the interior
of their domain.

In section~\ref{Se:perspectives} we give a perspective how the present work of research might be extended to study \emph{compact}
immersed minimal surfaces in \,$S^3$\,, especially those of genus \,$g\geq 2$\,.

\smallskip

\paragraph{\textbf{Acknowledgements.}} I would like to express my sincerest gratitude to Professor Martin Schmidt, who has advised me during the
creation of this work. His steady support and help has been invaluable to me. I have learned a lot from him. I would also like to thank 
Professor C.~Hertling for the introductory seminary on algebraic geometry and for further advice,
Dr.~Markus Knopf for discussions on spectral theory and especially on Darboux coordinates, 
Dr.~Alexander Klauer for hints on functional analysis, 
and Tobias Simon and Eva L\"ubcke for general mathematical talk.

\part{Spectral data}

\section{Minimal immersions into the 3-sphere and the sinh-Gordon equation}
\label{Se:minimal}

We begin by describing the relationship between minimal immersions without umbilical points into the 3-sphere
and solutions of the sinh-Gordon equation explicitly, especially to obtain the \,$\mathfrak{sl}(2,\C)$-valued flat connection form \,$\alpha_\lambda$\,;
from the integration of \,$\alpha_\lambda$\, we will obtain spectral data for periodic solutions of the sinh-Gordon equation. 

\label{not:minimal:S3}
Suppose that \,$f: X \to S^3$\, is a conformal immersion of some Riemann surface \,$X$\, into
\label{not:minimal:u}
the 3-sphere \,$S^3(\vkap) \subset \R^4$\, of curvature \,$\vkap>0$\,. Then there exists a function \,$u: X \to \R$\,
such that the pull-back Riemannian metric \,$g= f^* \langle \,\cdot\,,\,\cdot\, \rangle$\,
induced by \,$f$\, on \,$X$\, is given by \,$g = e^{u/2}\,\mathrm{d}z\,\mathrm{d}\overline{z} = e^{u/2}\,(\mathrm{d}x^2+\mathrm{d}y^2) $\,.
\,$u$\, is called the \emph{conformal factor} of \,$f$\,. If we denote the mean curvature function of \,$f$\,
by \,$H$\,, and the Hopf differential of \,$f$\, by \,$E\,\mathrm{d}z^2$\,, the equations
of Gauss and Codazzi read
\begin{align*}
u_{z\overline{z}} + \frac12 (H^2+\vkap)\,e^u-2\,E\,\overline{E}\,e^{-u} & = 0 \\
H_z\,e^u & = 2\,E_{\overline{z}} \; . 
\end{align*}
We now specialize to the situation where \,$f$\, is minimal, i.e.~\,$H=0$\,, and maps into the sphere of radius \,$2$\,, i.e.~\,$\vkap=\tfrac14$\,.
Moreover let us suppose that \,$f$\, has no umbilical points. Then we can choose a holomorphic chart \,$z$\,
of \,$X$\, such that the Hopf differential's function \,$E$\, is constant and of modulus \,$\tfrac14$\,. 
In this situation the equation of Codazzi reduces to \,$0=0$\,, whereas the equation of Gauss
reduces to the sinh-Gordon equation
\begin{equation}
\label{eq:minimal:sinh-gordon}
\Delta u + \sinh(u)=0 \; .
\end{equation}
This shows that minimal immersions \,$f: X \to S^3(\vkap)$\, give rise to solutions of the sinh-Gordon
equation.

We now describe the relationship between minimal immersions into \,$S^3:= S^3(\vkap=\tfrac14)$\, and solutions of
the sinh-Gordon equation in more detail. For this purpose, we suppose \,$X$\, to be simply 
connected.

We denote by \,$f_x$\, and \,$f_y$\, the derivative
of \,$f$\, in the direction \,$x$\, resp.~\,$y$\,, and put \,$e_x := \tfrac{f_x}{\|f_x\|}$\,,
\,$e_y := \tfrac{f_y}{\|f_y\|}$\,. Moreover, we let \,$N$\, be the unit normal field of the immersion \,$f$\,. 
Then \,$\left( e_x, e_y, N \right)$\, is an
orthonormal basis field of \,$TS^3$\, along \,$f$\,, hence the \emph{frame} 
\,$\underline{F}:=  \left( e_x, e_y, N ,\tfrac12f \right): X \to \SO(4)$\, 
is a positively oriented orthonormal basis field of \,$\R^4$\, along \,$f$\,. \,$X$\, being 
simply connected, we can lift \,$\underline{F}$\, to the
universal covering group \,$\SU(2)\times \SU(2)$\, of \,$\SO(4)$\,, obtaining \,$(F^{[1]},F^{[2]}): X \to \SU(2) \times \SU(2)$\,
with \,$F^{[1]},F^{[2]}: X \to \SU(2)$\,.


\label{not:minimal:lambda}
It is clear that the ``integrability condition'' \,$((F^{[1]},F^{[2]})_z)_{\overline z} = ((F^{[1]},F^{[2]})_{\overline z})_z$\, translates
into differential equations for \,$F^{[1]}$\, and \,$F^{[2]}$\,. But
it is an insight due to \textsc{Bobenko} (see \cite{Bobenko:1991b}) that the two components \,$F^{[1]},F^{[2]}$\, of the frame are governed by essentially the \emph{same}
differential equation. Indeed, there exists \,$\lambda\in S^1$\, so that we have
\,$(F^{[1]},F^{[2]})=(F_\lambda,F_{-\lambda})$\,, where \,$F_\lambda$\, satisfies the differential equation \,$\mathrm{d}F_\lambda = \alpha_\lambda\,F_\lambda$\,
with
$$ \alpha_\lambda := \frac14 \begin{pmatrix} -u_z & -\lambda^{-1}\,e^{-u/2} \\ e^{u/2} & u_z \end{pmatrix}\mathrm{d}z
+ \frac14\begin{pmatrix} u_{\overline{z}} & -e^{u/2} \\ \lambda\,e^{-u/2} & -u_{\overline{z}} \end{pmatrix} \mathrm{d}\overline{z} \; . $$
An explicit calculation shows that 
$$ \mathrm{d}\alpha_\lambda + [\alpha_\lambda\wedge\alpha_\lambda] = \frac18 (\Delta u + \sinh(u)) \begin{pmatrix} 1 & 0 \\ 0 & -1 \end{pmatrix} \mathrm{d}z \wedge \mathrm{d}\overline{z} \; , $$
therefore the Maurer-Cartan equation
\,$\mathrm{d}\alpha_{ \lambda} + [\alpha_{ \lambda} \wedge \alpha_{ \lambda}]=0$\, 
for \,$\alpha_{\lambda}$\, is equivalent to the sinh-Gordon equation for \,$u$\,. 
Therefore we have the following equivalence:

\begin{prop}
Suppose \,$u: X \to \C$\, is any differentiable function. Then the following three statements are equivalent:
\begin{enumerate}
\item The function \,$u$\, solves the sinh-Gordon equation \,$\Delta u + \sinh(u)=0$\,.
\item The metric \,$e^{u/2}\,\mathrm{d}z \,\mathrm{d}\overline{z}$\, satisfies the equations of Gauss and Codazzi 
for a minimal immersion into \,$S^3$\, with constant Hopf differential of modulus \,$\tfrac14$\,. 
\item For any, or for all \,$\lambda\in \C^*$\,,
\,$\alpha_\lambda$\, (as defined above) satisfies the Maurer-Cartan equation \,$\mathrm{d}\alpha_\lambda + [\alpha_\lambda\wedge\alpha_\lambda]=0$\,.
\end{enumerate}
\end{prop}

\label{not:minimal:F}
If any of the statements in the Proposition holds, then (3) shows that there is a unique solution \,$F_\lambda: X \to \SU(2)$\, to the differential equation 
\,$\mathrm{d}F_\lambda = \alpha_\lambda \, F_\lambda$\, with \,$F_\lambda(z_0)=\unity$\, (where \,$z_0\in X$\, is fixed) for every \,$\lambda\in\C^*$\,.
It is clear that \,$F_\lambda$\, depends holomorphically on \,$\lambda\in \C^*$\,. In this setting, the family \,$(F_\lambda)_{\lambda\in\C^*}$\, 
is called the \emph{extended frame} of \,$u$\, or of \,$f$\,; the parameter \,$\lambda$\, is called the \emph{spectral parameter}. 

On the other hand, if \,$u: X \to \C$\, is a solution of the sinh-Gordon equation, then 
the metric \,$e^{u/2}\,\mathrm{d}z \,\mathrm{d}\overline{z}$\, together
with any constant Hopf differential \,$E\,\mathrm{d}z^2$\, of modulus \,$\tfrac14$\, satisfies the equations of Gauss and Codazzi
for a minimal immersion into \,$S^3$\,. Thus, if \,$u$\, is real-valued in this setting, then \,$u$\, corresponds
to a minimal conformal immersion \,$f_E: X \to S^3$\, with that metric and that Hopf differential.
Thus we see that minimal immersions into \,$S^3$\, come
in families \,$(f_E)_{E\in (1/4)\,S^1}$\,; each such family is called an \emph{associated family}.

Therefore solutions \,$u:X \to \R$\, of the sinh-Gordon equation are in one-to-one correspondence with
associated families of minimal immersions into \,$S^3$\,.

Note that \,$\alpha_\lambda$\, has certain symmetries. 
In the case where the solution \,$u$\, is simply periodic, the symmetries of \,$\alpha_\lambda$\, also induce symmetries for the extended frame \,$F_\lambda$\,
and for the \emph{monodromy} \,$M(\lambda)=F_\lambda(z_1) \cdot F_\lambda(z_0)^{-1}$\, of \,$F_\lambda$\, along the period. In the following proposition we describe these symmetries,
where we suppose for the sake of simplicity that \,$X\subset \C$\, holds, that we have \,$0 \in X$\,, that the period of \,$u$\, is \,$1$\, (so the real interval \,$[0,1]$\, is contained in \,$X$\,),
and that we consider the base point \,$z_0=0$\,. Then the monodromy of the extended frame is given by \,$M(\lambda) = F_\lambda(z_0+1) \cdot F_\lambda(z_0)^{-1} = F_\lambda(1)$\,. 

One of the applications of these symmetries will be to derive asymptotic estimates for \,$\lambda\to 0$\, from asymptotic estimates for \,$\lambda\to\infty$\,. 

\begin{prop}
\label{P:mono:symmetry}
Let \,$\Phi: z \mapsto \overline{z}$\,. If \,$u: X \to \C$\, is a simply periodic solution of the sinh-Gordon equation as above, then 
\,$\widetilde{u} := -u\circ \Phi$\, and \,$\overline{u}: z \mapsto \overline{u(z)}$\, also are simply periodic solutions of the sinh-Gordon equation.
Moreover we put \,$g := \left( \begin{smallmatrix} 1 & 0 \\ 0 & \lambda \end{smallmatrix} \right)$\,. Then we have
\begin{enumerate}
\item 
\,$\alpha_{\lambda^{-1},u} = g^{-1} \cdot \Phi^*\alpha_{\lambda,\widetilde{u}}\cdot g$\,, \\
\,$F_{\lambda^{-1},u} = g^{-1} \cdot (F_{\lambda,\widetilde{u}} \circ \Phi) \cdot g$\,, \\
\,$M_u(\lambda^{-1}) = g^{-1} \cdot M_{\widetilde{u}}(\lambda) \cdot g$\,.
\item
\,$\alpha_{\overline{\lambda}^{-1},u} = -\overline{\alpha_{\lambda,\overline{u}}^t}$\,\quad (where \,$\alpha^t$\, denotes the transpose of \,$\alpha$\,), \\
\,$F_{\overline{\lambda}^{-1},u} = \overline{(F_{\lambda,\overline{u}}^{-1})^t}$\,, \\
\,$M_u(\overline{\lambda}^{-1}) = \overline{(M_{\overline{u}}(\lambda)^{-1})^t}$\,. 
\end{enumerate}
\end{prop}

\begin{proof}
\emph{For (1).} We have 
$$ \Phi^* \mathrm{d}z = \mathrm{d}\overline{z} \qmq{and} \Phi^* \mathrm{d}\overline{z} = \mathrm{d}z \;, $$
and also
$$ \widetilde{u}_z \circ \Phi = -u_{\overline{z}} \qmq{and} \widetilde{u}_{\overline{z}} \circ \Phi = -u_z \; . $$
Using these equations, we calculate:
\begin{align*}
& g^{-1} \cdot (\Phi^*\alpha_{\lambda,\wt{u}})(z) \cdot g \\
=\; & g^{-1} \cdot \biggr[ \frac14 \begin{pmatrix} -\wt{u}_z \circ \Phi & -\lambda^{-1}\,e^{-(\wt{u} \circ \Phi)/2} \\ e^{(\wt{u} \circ \Phi)/2} & \wt{u}_z  \circ \Phi \end{pmatrix} \Phi^*\mathrm{d}z 
+ \frac14\begin{pmatrix} \wt{u}_{\overline{z}}  \circ \Phi & -e^{(\wt{u} \circ \Phi)/2} \\ \lambda\,e^{-(\wt{u} \circ \Phi)/2} & -\wt{u}_{\overline{z}}  \circ \Phi \end{pmatrix} \Phi^*\mathrm{d}\overline{z} \biggr] \cdot g \\
=\; & g^{-1} \cdot \biggr[ \frac14 \begin{pmatrix} u_{\overline{z}} & -\lambda^{-1}\,e^{u/2} \\ e^{-u/2} & -u_{\overline{z}} \end{pmatrix} \mathrm{d}\overline{z} 
+ \frac14\begin{pmatrix} -u_z & -e^{-u/2} \\ \lambda\,e^{u/2} & u_z \end{pmatrix} \mathrm{d}z \biggr] \cdot g \\
=\; & \frac14 \begin{pmatrix} u_{\overline{z}} & -e^{u/2} \\ \lambda^{-1}\,e^{-u/2} & -u_{\overline{z}} \end{pmatrix} \mathrm{d}\overline{z} 
+ \frac14\begin{pmatrix} -u_z & -\lambda\,e^{-u/2} \\ e^{u/2} & u_z \end{pmatrix} \mathrm{d}z \\
=\; & \alpha_{\lambda^{-1},u}(z) \;.
\end{align*}

Now let \,$\widetilde{F} := g^{-1} \cdot (F_{\lambda,\widetilde{u}} \circ \Phi) \cdot g$\,. Then we have
\begin{align*}
\mathrm{d}\widetilde{F} & = g^{-1} \cdot \mathrm{d}(F_{\lambda,\widetilde{u}} \circ \Phi)\cdot g = g^{-1} \cdot \Phi^*\,\mathrm{d}F_{\lambda,\widetilde{u}}\cdot g = g^{-1} \cdot \Phi^*(\alpha_{\lambda,\widetilde{u}}\,F_{\lambda,\widetilde{u}}) \cdot g \\
& = g^{-1} \cdot (\Phi^*\alpha_{\lambda,\widetilde{u}}) \cdot g \cdot g^{-1}\cdot (F_{\lambda,\widetilde{u}}\circ \Phi) \cdot g = \alpha_{\lambda^{-1},u}\,\widetilde{F} 
\end{align*}
as well as \,$\widetilde{F}(0)=\unity$\,. In other words, \,$\widetilde{F}$\, solves the same initial value problem as \,$F_{\lambda^{-1},u}$\,, and thus we have \,$\widetilde{F}=F_{\lambda^{-1},u}$\,. 

Finally, we have \,$M_u(\lambda^{-1}) = F_{\lambda^{-1},u}(1) = g^{-1} \cdot F_{\lambda,\widetilde{u}}(\Phi(1)) \cdot g = g^{-1} \cdot M_{\widetilde{u}}(\lambda) \cdot g$\,. 

\emph{For (2).} We have
\begin{align*}
\overline{\alpha_{\overline{\lambda}^{-1},\overline{u}}}
& = \overline{\frac14 \begin{pmatrix} -\overline{u}_z & -\overline{\lambda}\,e^{-\overline{u}/2} \\ e^{\overline{u}/2} & \overline{u}_z \end{pmatrix} \mathrm{d}z 
        + \frac14 \begin{pmatrix} \overline{u}_{\overline{z}} & -e^{\overline{u}/2} \\ \overline{\lambda}^{-1}\,e^{-\overline{u}/2} & -\overline{u}_{\overline{z}} \end{pmatrix} \mathrm{d}\overline{z}} \\
& = \frac14 \begin{pmatrix} -u_{\overline{z}} & -\lambda\,e^{-u/2} \\ e^{u/2} & u_{\overline{z}} \end{pmatrix} \mathrm{d}\overline{z}
        + \frac14 \begin{pmatrix} u_{z} & -e^{u/2} \\ \lambda^{-1}\,e^{-u/2} & -u_{z} \end{pmatrix} \mathrm{d}z \\
& = -\alpha_{\lambda,u}^t \; . 
\end{align*}
Now let \,$\wt{F} := \overline{(F_{\lambda,u}^{-1})^t}$\,. We then have
\begin{align*}
\mathrm{d}\wt{F} & = \overline{(\mathrm{d} F_{\lambda,u}^{-1})^t} = -\overline{(F_{\lambda,u}^{-1}\,\mathrm{d}F_{\lambda,u}\,F_{\lambda,u}^{-1})^t} = -\overline{(F_{\lambda,u}^{-1}\,\alpha_{\lambda,u}\,F_{\lambda,u}\,F_{\lambda,u}^{-1})^t} \\
& = -\overline{(F_{\lambda,u}^{-1}\,\alpha_{\lambda,u})^t} = -\overline{\alpha_{\lambda,u}^t}\,\wt{F} = \alpha_{\overline{\lambda}^{-1},\overline{u}}\,\wt{F}
\end{align*}
as well as \,$\wt{F}(0)=\unity$\,. It follows that \,$\wt{F}$\, solves the initial value problem for \,$F_{\overline{\lambda}^{-1},\overline{u}}$\,, and hence \,$F_{\overline{\lambda}^{-1},\overline{u}} = \wt{F}$\, holds. 

Finally, we have \,$(M_u(\lambda)^{-1})^t = (F_{\lambda,u}(1)^{-1})^t = \overline{F_{\overline{\lambda}^{-1},\overline{u}}(1)} = \overline{M_{\overline{u}}(\overline{\lambda}^{-1})}$\,. 
\end{proof}

\section{Spectral data for simply periodic solutions of the sinh-Gordon equation}
\label{Se:spectrum}

We now suppose that a \emph{simply periodic} solution \,$u: X \to \C$\, of the sinh-Gordon equation
\,$\Delta u + \sinh(u)=0$\, on a domain \,$X \subset \C$\, is given. Without loss of generality, we
suppose that \,$0\in X$\, and that \,$1\in \C$\, is the period of \,$u$\,. Then \,$X$\, is a
horizontal strip in \,$\C$\, that contains the real axis \,$\R$\, and we have
$$ u(z+1)=u(z) \qmq{for all \,$z\in X$\,.} $$

\label{not:spectrum:F}
We let \,$F_\lambda$\, be the extended frame corresponding to \,$u$\,, i.e.~for each \,$\lambda\in \C^*$\,,
\,$F_\lambda: X \to \SU(2)$\, satisfies
$$ \mathrm{d}F_\lambda = \alpha_\lambda\,F_\lambda \;,\quad F_\lambda(0)=\unity $$
with
\label{not:spectrum:alpha}
\begin{align}
\label{eq:mono:alphazz}
\alpha_\lambda & := \frac14 \begin{pmatrix} -u_z & -\lambda^{-1}\,e^{-u/2} \\ e^{u/2} & u_z \end{pmatrix}\mathrm{d}z
+ \frac14\begin{pmatrix} u_{\overline{z}} & -e^{u/2} \\ \lambda\,e^{-u/2} & -u_{\overline{z}} \end{pmatrix} \mathrm{d}\overline{z} \\
& = \frac{1}{4} \begin{pmatrix} i\,u_y & -e^{u/2} - \lambda^{-1}\,e^{-u/2} \\ e^{u/2}+\lambda\,e^{-u/2} & -i\,u_y \end{pmatrix} \mathrm{d}x \notag \\
\label{eq:mono:alphaxy}
& \qquad\qquad + \frac{i}{4} \begin{pmatrix} -u_x & e^{u/2}-\lambda^{-1}\,e^{-u/2} \\ e^{u/2}-\lambda\,e^{-u/2} & u_x \end{pmatrix}\mathrm{d}y \; . 
\end{align}

\label{not:spectrum:M}
Although \,$u$\, is periodic, the extended frame \,$F_\lambda$\, generally is not. 
To measure the deviance of \,$F_\lambda$\, from being periodic, we are lead to consider
the \emph{monodromy} \,$M_z(\lambda) := F_\lambda(z+1)\cdot F_\lambda(z)^{-1}$\,. 
It should be noted that the monodromy \,$M$\, itself satisfies a differential equation with respect to differentiation in the base point \,$z$\,. Indeed, we have
\begin{align}
\mathrm{d}_z M_z(\lambda) 
& = \mathrm{d}\bigr(F_\lambda(z+1)\cdot F_\lambda(z)^{-1}\bigr) \notag \\
& = \mathrm{d}F_\lambda(z+1) \cdot F_\lambda(z)^{-1} - F_\lambda(z+1) \cdot F_\lambda(z)^{-1}\cdot \mathrm{d}F_\lambda(z) \cdot F_\lambda(z)^{-1} \notag \\
& = \alpha_\lambda(z+1)\cdot F_\lambda(z+1) \cdot F_\lambda(z)^{-1} - F_\lambda(z+1) \cdot F_\lambda(z)^{-1} \cdot \alpha_\lambda(z) \notag \\
\label{eq:mono:monodromy-dgl}
& = [\alpha_\lambda(z),M_z(\lambda)] \;, 
\end{align}
because \,$\alpha_\lambda(z+1)=\alpha_\lambda(z)$\, holds due to the periodicity of \,$u$\,. Consequently, we have
$$ M_z(\lambda) = F_\lambda(z) \cdot M_{z=0}(\lambda) \cdot F_\lambda(z)^{-1} \; . $$

We use the monodromy \,$M(\lambda):=M_{z=0}(\lambda)$\, at \,$z=0$\, to define \emph{spectral data} for
\,$u$\,. We first define the \emph{spectral curve} by the eigenvalues of \,$M(\lambda)$\,:
\begin{equation}
\label{eq:spectral:Sigma}
\Sigma := \Mengegr{(\lambda,\mu) \in \C^* \times \C}{\det(M(\lambda)-\mu\cdot \unity)=0} \; . 
\end{equation}
\,$\Sigma$\, is a non-compact complex curve (i.e.~a 1-dimensional complex analytic space), possibly with singularities. The natural functions \,$\lambda$\, and \,$\mu$\, on \,$\Sigma$\,
are holomorphic, and they generate the sheaf of holomorphic functions on \,$\Sigma$\,, i.e.~a function \,$f: \Sigma \to \C$\, is holomorphic if
and only if it can locally be extended to a holomorphic function in \,$\lambda$\, and \,$\mu$\, on a neighborhood in \,$\C^* \times \C$\,.

To phrase the definition  \eqref{eq:spectral:Sigma} of \,$\Sigma$\, more explicitly, we write the monodromy \,$M(\lambda)$\, in the form
\begin{equation}
\label{eq:spectral:M-entries}
M(\lambda) = \begin{pmatrix} a(\lambda) & b(\lambda) \\ c(\lambda) & d(\lambda) \end{pmatrix}
\end{equation}
with the holomorphic functions \,$a,b,c,d: \C^* \to \C$\,. Then we have
\begin{equation}
\label{eq:spectral:tr}
\Delta(\lambda) := \tr(M(\lambda)) = a(\lambda) + d(\lambda) \;.
\end{equation}
Moreover, because \,$\alpha_\lambda$\, is trace-free, we have
\begin{equation}
\label{eq:spectral:det}
a(\lambda)\cdot d(\lambda)- b(\lambda) \cdot c(\lambda) = \det(M(\lambda))=1 \; . 
\end{equation}
Equations~\eqref{eq:spectral:tr} and \eqref{eq:spectral:det} show that the characteristic polynomial of \,$M(\lambda)$\, 
is \,$\mu^2-\Delta(\lambda)\cdot \mu+1$\,, and thus we obtain 
\begin{equation}
\label{eq:spectral:Sigma2}
\Sigma = \Mengegr{(\lambda,\mu) \in \C^* \times \C}{\mu^2 - \Delta(\lambda)\cdot \mu+1=0} \; . 
\end{equation}
It follows from this presentation of \,$\Sigma$\, that for any \,$(\lambda,\mu)\in \Sigma$\,, \,$\mu\neq 0$\, holds, that there are 
at least one and at most two points \,$(\lambda,\mu) \in \Sigma$\, that are above some given \,$\lambda\in \C^*$\,, and that 
the holomorphic involution
\begin{equation}
\label{eq:spectral:sigma}
\sigma: \Sigma \to\Sigma,\; (\lambda,\mu) \mapsto (\lambda,\mu^{-1})
\end{equation}
maps \,$\Sigma$\, onto itself. Therefore the complex curve \,$\Sigma$\, is hyperelliptic with the hyperelliptic involution \,$\sigma$\,.

More generally, for a given \,$\lambda \in \C^*$\, we have \,$(\lambda,\mu)\in \Sigma$\, if and only if
\begin{equation}
\label{eq:spectral:mu}
\mu = \frac{1}{2}\left( \Delta(\lambda) \pm \sqrt{\Delta(\lambda)^2 - 4} \right) 
\end{equation}
holds; equivalently \,$(\lambda,\mu)\in \Sigma$\, is characterized by the eigenvalue equation
\begin{equation}
\label{eq:spectral:eigeneq}
(a(\lambda)-\mu)\cdot (d(\lambda)-\mu) = b(\lambda) \cdot c(\lambda) \; .
\end{equation}

Note that \,$\sqrt{\Delta^2-4}=\mu-\mu^{-1}$\, is a holomorphic function on \,$\Sigma$\,, which is
anti-symmetric with respect to \,$\sigma$\,, i.e.~we have \,$\sqrt{\Delta^2-4} \circ \sigma = -\sqrt{\Delta^2-4}$\,. 

It also follows from Equation~\eqref{eq:spectral:mu} that \,$(\lambda_*,\mu_*)\in\Sigma$\, is a fixed point of the hyperelliptic involution
\,$\sigma$\, if and only if \,$\Delta(\lambda_*) = \pm 2$\, (or \,$\Delta(\lambda_*)^2-4=0$\,) holds;
this is the case if and only if \,$\mu_* \in \{ \pm 1\}$\, holds, and in this case, there is only one point in \,$\Sigma$\, 
above \,$\lambda_*$\,. 
Such a point is a branch point of \,$\Sigma$\, if the zero of \,$\Delta\mp 2$\, (or \,$\Delta^2-4$\,) is of odd order, and it is a singularity
of \,$\Sigma$\, if the zero of \,$\Delta\mp 2$\, (or \,$\Delta^2-4$\,) is of order \,$\geq 2$\, (an ordinary double point if the zero is of order \,$2$\,);
in the latter case we call the order of zero of \,$\Delta \mp 2$\, (or \,$\Delta^2-4$\,) the \emph{order of the singularity}.
There are no other branch points or singularities of \,$\Sigma$\,.

In this situation, the zeros of \,$\Delta^2-4$\, (corresponding to the branch points and singularities of \,$\Sigma$\,) 
are analogous to the Dirichlet spectrum in the spectral theory of the 1-dimensional Schr\"odinger equation (see for example, \cite{Poeschel-Trubowitz:1987}, Chapter~2, p.~25f.).

\label{not:spectral:calMcalO}
We will now describe the local complex space structure of \,$\Sigma$\,, in particular near its singularities.
We let \,$\mathcal{O}$\, resp.~\,$\mathcal{M}$\, be the sheaf of holomorphic resp.~of meromorphic functions on \,$\Sigma$\,. 
The sheaf of unitary rings \,$\mathcal{O}$\, is generated by the holomorphic functions \,$\lambda$\, and \,$\mu$\, (or by \,$\lambda$\, and \,$\mu-\mu^{-1}$\,)
in the sense that a function \,$f$\, on \,$\Sigma$\,
is holomorphic if and only if can locally be written in the form
$$ f = f_+(\lambda) + f_-(\lambda)\cdot (\mu-\mu^{-1}) $$
with holomorphic functions \,$f_+$\, and \,$f_-$\, in \,$\lambda$\,. In other words, \,$f$\, is holomorphic at some \,$(\lambda_*,\mu_*)\in \Sigma$\, if and only if it can be extended to a holomorphic function
in \,$\lambda$\, and \,$\mu$\, on a neighborhood of \,$(\lambda_*,\mu_*)$\, in \,$\C^* \times \C^*$\,. 

Near any point \,$(\lambda_*,\mu_*)$\, of \,$\Sigma$\, that is not a fixed point of \,$\sigma$\, (i.e.~that is neither a branch point nor a singularity), the holomorphic
function \,$\lambda-\lambda_*$\, is a local coordinate on \,$\Sigma$\, near \,$(\lambda_*,\mu_*)$\,, 
and near any regular branch point \,$(\lambda_*,\mu_*)$\, of \,$\Sigma$\,, the function \,$\mu-\mu^{-1}$\, is a local coordinate.

\label{not:spectral:whSigma}
To study the structure of \,$\Sigma$\, at singular points, we 
consider the normalization \,$\wh{\Sigma}$\, of \,$\Sigma$\, together with the one-sheeted covering \,$\pi: \wh{\Sigma}\to\Sigma$\, (see for example \cite{deJong:2000}, Section~4.4, p.~161ff.). 
We let \,$\mathcal{O}_{\wh{\Sigma}}$\, resp.~\,$\wh{\mathcal{M}}_{\wh{\Sigma}}$\, be the sheaf of holomorphic functions resp.~of meromorphic functions on \,$\wh{\Sigma}$\,. Then the abelian group \,$H^0(\Sigma,\mathcal{M})$\, of 
meromorphic functions on \,$\Sigma$\, is isomorphic to the abelian group \,$H^0(\Sigma,\pi_*\mathcal{M}_{\wh{\Sigma}})$\, of sections in the direct image sheaf \,$\pi_*\mathcal{M}_{\wh{\Sigma}}$\, via \,$f \mapsto f\circ \pi$\,,
and we identify \,$\pi_*\mathcal{M}_{\wh{\Sigma}}$\, with \,$\mathcal{M}$\, in this way. Under this identification, \,$\wh{\mathcal{O}} := \pi_*\mathcal{O}_{\wh{X}}$\, corresponds to the germs of meromorphic functions on \,$\Sigma$\,
which are locally bounded (\cite{deJong:2000}, Theorem~4.4.15, p.~167). 

If \,$\Sigma$\, has at \,$(\lambda_*,\mu_*)$\, a singularity of odd order \,$2n+1$\, (i.e.~\,$\Delta^2-4$\, has at \,$\lambda_*$\, a zero of order \,$2n+1$\,),
then there is exactly one point in \,$\wh{\Sigma}$\, above \,$(\lambda_*,\mu_*)$\,, this point is a branch point of the normalization,
\,$\sqrt{\lambda-\lambda_*}$\, is a coordinate of \,$\widehat{\Sigma}$\, near that branch point, 
\,$\mathrm{d}\lambda$\, vanishes at \,$(\lambda_*,\mu_*)$\, of first order,
and we have (on \,$\wh{\Sigma}$\,)
\begin{equation}
\label{eq:spectrum:mumu-1-odd}
\mu-\mu^{-1} = \left(\sqrt{\lambda-\lambda_*}\right)^{2n+1} \cdot \vi
\end{equation}
with a non-zero holomorphic function germ \,$\vi \in \wh{\mathcal{O}}$\,. 
It follows that the \,$\delta$-invariant \,$\dim(\wh{\mathcal{O}}_{(\lambda_*,\mu_*)}/\mathcal{O}_{(\lambda_*,\mu_*)})$\, of \,$\Sigma$\, at \,$(\lambda_*,\mu_*)$\, (see \cite{deJong:2000}, Definition~5.2.1(1), p.~186)
equals \,$n$\,. 
A meromorphic function germ \,$f\in \mathcal{M}$\, on \,$\Sigma$\, at \,$(\lambda_*,\mu_*)$\, can be described as a 
meromorphic function germ \,$\wh{f}$\, (i.e.~a Laurent series) in \,$\sqrt{\lambda-\lambda_*}$\,; in this setting
we define the order of a zero (or the order of a pole) of \,$f$\, at \,$(\lambda_*,\mu_*)$\,, denoted by \,$\ord^{\Sigma}(f)$\, (or by \,$\polord^{\Sigma}(f):=-\ord^\Sigma(f)$\,)
as the order of the zero (or pole) of \,$\wh{f}$\, as a function in \,$\sqrt{\lambda-\lambda_*}$\,.
In particular, for a function \,$g$\, that is meromorphic in \,$\lambda$\,, we have \,$\ord^{\Sigma}(g) = 2\cdot \ord^{\C}(g)$\,. 

On the other hand, if \,$\Sigma$\, has at \,$(\lambda_*,\mu_*)$\, a singularity of even order \,$2n$\, (i.e.~\,$\Delta^2-4$\, has at \,$\lambda_*$\, a zero of order \,$2n$\,),
then there are two points in \,$\wh{\Sigma}$\, above \,$(\lambda_*,\mu_*)$\,, these points are not branch points of \,$\wh{\Sigma}$\,,
\,$\lambda-\lambda_*$\, is a coordinate near each of these two points, so \,$\mathrm{d}\lambda$\, does not vanish at either of these points.
A meromorphic function germ \,$f\in \mathcal{M}$\, on \,$\Sigma$\, at \,$(\lambda_*,\mu_*)$\, can be described as a pair
\,$(f_1,f_2)$\, of germs of functions meromorphic in \,$\lambda-\lambda_*$\,, and in this sense we have
\begin{equation}
\label{eq:spectrum:mumu-1-even}
\mu-\mu^{-1} = \bigr((\lambda-\lambda_*)^n \cdot \vi \;,\; -(\lambda-\lambda_*)^n \cdot \vi \bigr)
\end{equation}
with a non-zero holomorphic function germ \,$\vi$\, in \,$\lambda$\,. 
It follows that the \,$\delta$-invariant \,$\dim(\wh{\mathcal{O}}_{(\lambda_*,\mu_*)}/\mathcal{O}_{(\lambda_*,\mu_*)})$\, of \,$\Sigma$\, at \,$(\lambda_*,\mu_*)$\, again equals \,$n$\,. 
We define the order of a zero (or the order of a pole) of a meromorphic function germ \,$f$\,
at \,$(\lambda_*,\mu_*)$\,, denoted by \,$\ord^{\Sigma}(f)$\, (or by \,$\polord^{\Sigma}(f) := -\ord^\Sigma(f)$\,), as
the sum of the two zero (or pole) orders 
\begin{align*}
\ord^\Sigma(f) := \ord^{\C}(f_1)+\ord^{\C}(f_2) \qmq{or}
& \polord^\Sigma(f) := \polord^{\C}(f_1)+ \polord^{\C}(f_2)\; ,
\end{align*}
where \,$f$\, is represented by \,$(f_1,f_2)$\, as above.
In particular, for a function \,$g$\, that is meromorphic in \,$\lambda$\,, we again have \,$\ord^{\Sigma}(g) = 2\cdot \ord^{\C}(g)$\,. 

\bigskip

\label{not:spectral:Lambda}
The spectral curve \,$\Sigma$\, alone does not fully characterize the monodromy \,$M(\lambda)$\,. It describes
the eigenvalues of \,$M(\lambda)$\,, but not the corresponding eigenvectors. 
Therefore we define a second spectral object, namely the divisor associated to the bundle \,$\Lambda$\, of eigenvectors of \,$M(\lambda)$\,,
which is a holomorphic line bundle at least on \,$\Sigma'$\,, the Riemann surface of regular points of \,$\Sigma$\,. 

To obtain an explicit description of \,$\Lambda$\,, we note that \,$(v_1,v_2)\in \C^2$\, is an eigenvector of \,$M(\lambda)$\, corresponding
to the eigenvalue \,$\mu \in \C^*$\, if and only if either of the two equivalent (by the eigenvalue equation~\eqref{eq:spectral:eigeneq}) equations
$$ (a(\lambda)-\mu)\cdot v_1 + b(\lambda)\cdot v_2 = 0 \qmq{or} c(\lambda)\cdot v_1 + (d(\lambda)-\mu)\cdot v_2 =0 \; $$
holds.
Therefore \,$(\tfrac{\mu-d(\lambda)}{c(\lambda)},1)$\,
is a global meromorphic section of \,$\Lambda$\,, and thus the divisor on \,$\Sigma$\, characterizing the line bundle \,$\Lambda$\,
is the pole divisor \,$D$\, of the meromorphic function \,$\tfrac{\mu-d}{c}$\, on \,$\Sigma$\,.

The divisor \,$D$\, is analogous in the spectral theory for the 1-dimensional Schr\"odinger operator \,$L$\, to the constants \,$\kappa_n$\,
associated to the periodic eigenfunctions of \,$L$\,, that uniquely determine the potential of \,$L$\, together with the periodic spectrum of \,$L$\,
(see for example, \cite{Poeschel-Trubowitz:1987}, p.~59 and Theorem~3.5, p.~62).
We (preliminarily) call \,$D$\, the \emph{spectral divisor} of the monodromy \,$M(\lambda)$\,.


\begin{prop}
\label{P:spectral:D-regular}
If \,$(\lambda_*,\mu_*)$\, is a regular point of \,$\Sigma$\, that is a member of the support of \,$D$\,, say with multiplicity \,$m\geq 1$\,,
then \,$\ord_{\lambda_*}^{\C}(c)=m$\, and \,$\ord_{(\lambda_*,\mu_*)}^{\Sigma}(\mu-a)\geq m$\, holds. 
In particular \,$c(\lambda_*)=0$\, and \,$a(\lambda_*)=\mu_*$\, holds. Moreover, if \,$(\lambda_*,\mu_*)$\, is a regular branch point of \,$\Sigma$\,,
then \,$m=1$\, holds.
\end{prop}

\begin{proof}
Let \,$(\lambda_*,\mu_*) \in \Sigma$\, be a regular point of \,$\Sigma$\, that is in the support of \,$D$\, with multiplicity \,$m$\,,
and therefore a pole of \,$\tfrac{\mu-d}{c}$\, of order \,$m$\,. 


Let us first suppose that \,$(\lambda_*,\mu_*)$\, is not a branch point of \,$\Sigma$\,, i.e.~that \,$\Delta^2-4$\, is non-zero at \,$\lambda_*$\,. Then
both \,$\lambda$\, and \,$\mu$\, are local coordinates of \,$\Sigma$\, near \,$(\lambda_*,\mu_*)$\,, and therefore the order of zeros or poles of holomorphic functions
at \,$(\lambda_*,\mu_*)$\, do not depend on whether we view the functions on \,$\Sigma$\, or on \,$\C^*$\,. Moreover, it is not possible for both of the functions
\,$\mu-a$\, and \,$\mu-d$\, to be zero at \,$(\lambda_*,\mu_*)$\,; it follows by Equation~\eqref{eq:spectral:eigeneq} that either \,$\mu-a$\, or \,$\mu-d$\,
has a zero of order at least \,$m$\,. If \,$\mu-d$\, had a zero of order at least \,$m$\,, then \,$\tfrac{\mu-d}{c}$\,
would not have a pole at \,$(\lambda_*,\mu_*)$\,, and therefore \,$(\lambda_*,\mu_*)$\, would not be in the support of \,$D$\,. 
Thus we see that \,$\mu-d$\, is not zero at \,$(\lambda_*,\mu_*)$\,, and therefore the pole order \,$m$\, of \,$\tfrac{\mu-d}{c}$\, equals the exact order of the zero of \,$c$\,.
Moreover Equation~\eqref{eq:spectral:eigeneq} shows that \,$\mu-a$\, has a zero of order at least \,$m$\,. 

Now consider the case that \,$(\lambda_*,\mu_*)$\, is a regular branch point of \,$\Sigma$\,, i.e.~that \,$\Delta^2-4$\, has a simple zero at \,$\lambda_*$\,. 
Then \,$\mu_*\in \{\pm 1\}$\, and therefore \,$a(\lambda_*)=d(\lambda_*)=\mu_*$\, holds. Thus we see that as function on \,$\C^*$\,, \,$(a-d)^2$\, has a zero of order at least \,$2$\,;
the equation (which follows from Equations~\eqref{eq:spectral:tr} and \eqref{eq:spectral:det})
$$ \Delta^2 - 4 = (a+d)^2 -4(ad-bc) = (a-d)^2 + 4bc $$
therefore shows that \,$bc$\, can have a zero of order at most \,$1$\,. Thus \,$\ord^{\C}_{\lambda_*}(c)=1$\, and \,$b(\lambda_*)\neq 0$\,
holds. Because of \,$\ord^{\Sigma}_{(\lambda_*,\mu_*)}(\lambda)=2$\,, we thus obtain \,$\ord^{\Sigma}_{(\lambda_*,\mu_*)}(c)=2$\, and \,$\ord^{\Sigma}_{(\lambda_*,\mu_*)}(\mu-d)\geq 1$\,. Because
\,$\tfrac{\mu-d}{c}$\, has a pole at \,$(\lambda_*,\mu_*)$\,, in fact \,$\ord^{\Sigma}_{(\lambda_*,\mu_*)}(\mu-d)= 1$\, is the only possibility. Thus the pole order \,$m$\, of \,$\tfrac{\mu-d}{c}$\, equals
\,$2-1=1$\,, and hence \,$\ord^{\C}_{\lambda_*}(c)=m$\, holds.
\end{proof}



The preceding Proposition shows in particular that 
if the support of \,$D$\, is contained in the set of regular points of \,$\Sigma$\,, then we have
\begin{equation}
\label{eq:spectral:D-classical}
D = \Menge{(\lambda,\mu)^m \in \C^* \times \C}{c(\lambda)=0, \mu=a(\lambda), m=\ord_\lambda^{\C}(c)} \; . 
\end{equation}

However, in general it is possible for poles of \,$\tfrac{\mu-d}{c}$\, to lie in singularities of \,$\Sigma$\,, and in these points,
the spectral divisor \,$D$\, as defined above does not contain enough information to completely characterize the behavior of the
monodromy \,$M(\lambda)$\,. To handle this case, we need to generalize the concept of a divisor in an appropriate way such that
the necessary additional information at the singularities of \,$\Sigma$\, at which \,$\tfrac{\mu-d}{c}$\, is not holomorphic is included.
It turns out that the most suitable generalization of the concept of a divisor for the present problem has been introduced by \textsc{Hartshorne}
in \cite{Hartshorne:1986}, \S 1. We now describe Hartshorne's concept of a generalized divisor, which he introduced for general 
Gorenstein curves in \cite{Hartshorne:1986} (and later, in even more generality), in our setting, i.e.~on a hyperelliptic complex
plane curve \,$\Sigma$\,. 

As before, we denote by \,$\mathcal{O}$\, resp.~by \,$\mathcal{M}$\, the sheaf of holomorphic functions resp.~of meromorphic functions
on \,$\Sigma$\,. 
A \emph{generalized divisor} is a subsheaf \,$\mathcal{D}$\, of \,$\mathcal{M}$\, that is finitely generated
over \,$\mathcal{O}$\,, and we say that \,$\mathcal{D}$\, is \emph{positive} if \,$\mathcal{O}\subset \mathcal{D}$\, holds.
For a positive generalized divisor \,$\mathcal{D}$\,, the \emph{support} of \,$\mathcal{D}$\, is the set of points \,$(\lambda,\mu)\in \Sigma$\, 
for which \,$\mathcal{D}_{(\lambda,\mu)} \neq \mathcal{O}_{(\lambda,\mu)}$\, holds. The map
\,$\Sigma \to \Z,\; (\lambda,\mu) \mapsto \dim\bigr(\mathcal{D}_{(\lambda,\mu)}/\mathcal{O}_{(\lambda,\mu)}\bigr)$\,
defines a divisor on \,$\Sigma$\, in the usual sense, which we call the \emph{underlying classical divisor} of \,$\mathcal{D}$\,. 

\begin{Def}
\label{D:spectral:D-generalized}
The \emph{(generalized) spectral divisor} of the monodromy 
$$ M(\lambda) = \left( \begin{smallmatrix} a(\lambda) & b(\lambda) \\ c(\lambda) & d(\lambda) \end{smallmatrix} \right) $$
is the subsheaf \,$\calD$\, of \,$\calM$\, on \,$\Sigma$\, generated by the meromorphic functions \,$1$\, and \,$\tfrac{\mu-d}{c}$\, over \,$\mathcal{O}$\,.
\end{Def}

\label{not:spectral:D-classical}
The generalized spectral divisor \,$\calD$\, is the proper replacement for the pole divisor of \,$\tfrac{\mu-d}{c}$\, in the classical sense, which we regarded above as spectral divisor. It is obviously positive, 
and at regular points of \,$\Sigma$\,, \,$\mathcal{D}$\, is equivalent to that pole divisor.
It should be noted that it is possible for a singular point \,$(\lambda_*,\mu_*)$\, to be in the support of \,$\calD$\,
even if \,$\tfrac{\mu-d}{c}$\, does not have a pole at this point; this happens if \,$\tfrac{\mu-d}{c}$\, is locally bounded (and therefore holomorphic on the normalization \,$\wh{\Sigma}$\,) at \,$(\lambda_*,\mu_*)$\,,
but not holomorphic on \,$\Sigma$\, at this point. Therefore \,$\tfrac{\mu-d}{c}$\, is not necessarily a function of maximal pole order in \,$\calD_{(\lambda_*,\mu_*)}$\,. 
From now on, the \emph{classical spectral divisor} \,$D$\, is always
the underlying classical divisor of \,$\calD$\, (rather than the pole divisor of \,$\tfrac{\mu-d}{c}$\, in the classical sense).
We will see below that \,$D$\, is given by Equation~\eqref{eq:spectral:D-classical}, now without the restriction to regular points.

We introduced the generalized spectral divisor because it conveys more information than a classical divisor at the singular points of \,$\Sigma$\,. 
We will now study in more detail exactly which additional information is conveyed.
The following proposition gives the key for understanding the 
structure of generalized divisors on the hyperelliptic complex curve \,$\Sigma$\, in general: 

\begin{prop}
\label{P:spectrum:locally-free}
Any generalized divisor \,$\mathcal{D}$\, on \,$\Sigma$\, is locally free over a unique subsheaf of rings \,$\mathcal{R}$\, of \,$\wh{\calO}$\,. 

More specifically, if \,$(\lambda_*,\mu_*)$\, is a singular point of \,$\Sigma$\,, and 
\,$f \in \mathcal{D}_{(\lambda_*,\mu_*)}$\, is of maximal pole order, let \,$j_0\geq 0$\, be the largest integer so that 
\,$\tfrac{\mu-\mu^{-1}}{(\lambda-\lambda_*)^{j_0}} \cdot f \in \mathcal{D}_{(\lambda_*,\mu_*)}$\, holds, and let \,$\mathcal{R}_{(\lambda_*,\mu_*)}$\, be the subring of \,$\wh{\mathcal{O}}_{(\lambda_*,\mu_*)}$\,
generated by \,$\tfrac{\mu-\mu^{-1}}{(\lambda-\lambda_*)^{j_0}}$\, and \,$\mathcal{O}_{(\lambda_*,\mu_*)}$\,. Then \,$\mathcal{D}_{(\lambda_*,\mu_*)}$\, is generated by \,$f$\,
over \,$\mathcal{R}_{(\lambda_*,\mu_*)}$\,, or by \,$f$\, and \,$\tfrac{\mu-\mu^{-1}}{(\lambda-\lambda_*)^{j_0}}\cdot f$\, over \,$\mathcal{O}_{(\lambda_*,\mu_*)}$\,. 
Here we have \,$j_0 \in \{0,\dotsc,n\}$\,, where the order of the singularity \,$(\lambda_*,\mu_*)$\, is either \,$2n$\, or \,$2n+1$\,. 
\end{prop}

Note that in this proposition, the case \,$j_0=0$\, is equivalent to \,$\calR_{(\lambda_*,\mu_*)}=\calO_{(\lambda_*,\mu_*)}$\,; in this case \,$\calD$\, is at \,$(\lambda_*,\mu_*)$\, a locally free
divisor on \,$\Sigma$\,. The case \,$j_0=n$\, is equivalent to \,$\calR_{(\lambda_*,\mu_*)} = \wh{\calO}_{(\lambda_*,\mu_*)}$\,; in this case \,$\calD$\, corresponds at \,$(\lambda_*,\mu_*)$\, to a locally
free divisor on the normalization \,$\wh{\Sigma}$\, of \,$\Sigma$\,. Moreover, if the maximal pole order of function germs in \,$\calD_{(\lambda_*,\mu_*)}$\, is \,$0$\, (then the germs in 
\,$\calD_{(\lambda_*,\mu_*)}$\, are locally bounded), then \,$1$\, is a function of maximal pole order in \,$\calD_{(\lambda_*,\mu_*)}$\, and therefore \,$\calR_{(\lambda_*,\mu_*)} = \calD_{(\lambda_*,\mu_*)}$\, holds. 

\begin{rem}
\label{R:spectrum:locally-free-interpretation}
As noted in the Introduction, \textsc{Hitchin} uses in his classification of the minimal tori in \,$S^3$\, in \cite{Hitchin:1990} a certain partial desingularization of \,$\Sigma$\, as spectral curve;
his spectral curve is precisely the complex curve \,$\wt{\Sigma}$\, ``between'' \,$\Sigma$\, and its normalization \,$\wh{\Sigma}$\, whose sheaf of holomorphic functions is the sheaf of rings \,$\calR$\,
from Proposition~\ref{P:spectrum:locally-free}.

We also note that for the validity of Proposition~\ref{P:spectrum:locally-free}, the fact that \,$\Sigma$\, is hyperelliptic is
of crucial importance. For every \,$m\geq 3$\,, there are examples of complex curves which are \,$m$-fold (branched) coverings
above \,$\PP^1$\, and on which there exist generalized divisors that are not locally free. An example for \,$m=3$\,
is described in \cite{Schmidt:1996}, Example~9.3, p.~93f..
\end{rem}

\begin{proof}[Proof of Proposition~\ref{P:spectrum:locally-free}.]
We need to consider the situation only at singularities of \,$\Sigma$\,. 
Let \,$(\lambda_*,\mu_*)\in \Sigma$\, be a singular point of \,$\Sigma$\,, say of order \,$2n+1$\, or \,$2n$\, with \,$n\geq 1$\,.
In what follows, we consider the stalks at \,$(\lambda_*,\mu_*)$\, of the sheaves under consideration, and omit the subscript \,${}_{(\lambda_*,\mu_*)}$\, for them.
We suppose that \,$(\lambda_*,\mu_*)$\, is in the support of \,$\mathcal{D}$\,, i.e.~that \,$\mathcal{D}/\mathcal{O} \neq \{0\}$\, holds.



We now choose \,$f\in \mathcal{D}$\, such that its pole order (as defined above) is maximal among the elements of \,$\mathcal{D}$\,,
and let \,$\widehat{\calD}$\, resp.~\,$\widetilde{\calD}$\, be the generalized divisor generated by \,$f$\, over \,$\wh{\mathcal{O}}$\, resp.~over \,$\mathcal{O}$\,.
Because \,$\mathcal{D}$\, is an \,$\mathcal{O}$-module, we have \,$\wt{\mathcal{D}}\subset \mathcal{D}$\,. We also have 
\,$\mathcal{D} \subset \wh{\mathcal{D}}$\, (in the case where \,$(\lambda_*,\mu_*)$\, is a singularity of even order, both components \,$f_\nu$\, of the
representant \,$(f_1,f_2)$\, of \,$f$\, have maximal pole order); from this inclusion it also follows that 
\,$\wh{\calD}$\, is the generalized divisor generated by \,$\calD$\, over \,$\wh{\calO}$\,. 


Because of \,$\dim(\wh{\mathcal{O}}/\mathcal{O})=n$\,, we have \,$\dim(\wh{\mathcal{D}}/\wt{\mathcal{D}})=n$\,. 
Moreover, \,$(\vi_j)_{j=1,\dotsc,n}$\, with \,$\vi_j := \tfrac{\mu-\mu^{-1}}{(\lambda-\lambda_*)^j}$\, is a basis of \,$\wh{\mathcal{O}}/\mathcal{O}$\,, 
and therefore \,$(f_j)_{j=1,\dotsc,n}$\, with \,$f_j := \vi_j\cdot f$\, is a basis of \,$\wh{\mathcal{D}}/\wt{\mathcal{D}}$\,.


For \,$g\in \wh{\mathcal{D}}/\wt{\mathcal{D}}$\, written as \,$g = \sum_{j=1}^n t_j\,f_j$\, with \,$t_j\in \C$\,, we let \,$j_0(g) \in \{0,\dotsc,n\}$\,
be the largest \,$j_0(g)\geq 1$\, so that \,$t_{j_0(g)}\neq 0$\,, or \,$j_0(g)=0$\, for \,$g=0$\,. Further we let \,$j_0 \in \{0,\dotsc,n\}$\, be the largest value of \,$j_0(g)$\,
that occurs for any \,$g\in \mathcal{D}/\wt{\mathcal{D}}$\,. Then we claim that \,$(f_j)_{j=1,\dotsc,j_0}$\, is a basis of
\,$\mathcal{D}/\wt{\mathcal{D}}$\, (with \,$\calD/\wt{\calD}=\{0\}$\, if \,$j_0=0$\,). It only needs to be shown that \,$f_j\in \mathcal{D}/\mathcal{\wt{D}}$\, holds for
\,$j\in \{1,\dotsc,j_0\}$\,, and to show this, we let \,$g = \sum_{j=1}^n t_j\,f_j \in \mathcal{D}/\wt{\mathcal{D}}$\,
be with \,$t_{j_0}\neq 0$\, and \,$t_j=0$\, for all \,$j>j_0$\,. Then \,$g\cdot (\lambda-\lambda_0)^{j_0-1} = t_{j_0}\,f_1 \in
\mathcal{D}/\wt{\mathcal{D}}$\, and therefore \,$f_1 \in \mathcal{D}/\wt{\mathcal{D}}$\, holds. Now suppose that we have already
shown that \,$f_1,\dotsc,f_{j_1-1} \in \mathcal{D}/\wt{\mathcal{D}}$\, holds for some \,$j_1 \in \{2,\dotsc,j_0\}$\,. Then we also
have \,$(g - \sum_{j=1}^{j_1-1} t_j\,f_j) \cdot (\lambda-\lambda_*)^{j_0-j_1} = t_{j_0}\,f_{j_1} \in \mathcal{D}/\wt{\mathcal{D}}$\, and
thus \,$f_{j_1} \in \mathcal{D}/\wt{\mathcal{D}}$\,. 

Therefore \,$\mathcal{D}$\, is locally free at \,$(\lambda_*,\mu_*)$\,; it is generated by \,$f$\, over the subring \,$\mathcal{R}$\, of
\,$\wh{\calO}$\, generated by \,$\tfrac{\mu-\mu^{-1}}{(\lambda-\lambda_k)^{j_0}} \in \wh{\calO}$\, and
\,$\mathcal{O}$\,. 
%
\end{proof}


\begin{prop}
\label{P:spectrum:msj0}
Let \,$\calD$\, be a positive generalized divisor on \,$\Sigma$\, and \,$(\lambda_*,\mu_*)$\, be a singular point of \,$\Sigma$\, that is in the support of \,$\calD$\,. We let \,$m := \dim(\calD_{(\lambda_*,\mu_*)}/\calO_{(\lambda_*,\mu_*)})$\,
be the degree of \,$\calD$\, at \,$(\lambda_*,\mu_*)$\,, \,$s\geq 0$\, be the maximal pole order that occurs for a germ in \,$\calD_{(\lambda_*,\mu_*)}$\,, and \,$j_0$\, be the number from Proposition~\ref{P:spectrum:locally-free}.
Then we have \,$m=s+j_0$\,. 
\end{prop}

\begin{proof}
As before, we work in the stalks of sheaves at \,$(\lambda_*,\mu_*)$\, and omit the subscript \,${}_{(\lambda_*,\mu_*)}$\,. 

Like in the proof of Proposition~\ref{P:spectrum:locally-free} we let \,$\wh{\calD}$\, be the 
generalized divisor generated by \,$\mathcal{D}$\, over \,$\wh{\mathcal{O}}$\,. 
Because \,$\calD$\, and \,$\wh{\calD}$\, are positive, we then have the inclusions \,$\calO \subset \calD \subset \wh{\calD}$\, and \,$\calO \subset \wh{\calO} \subset \wh{\calD}$\,, which imply the 
short exact sequences
$$ 0 \longrightarrow \calD/\calO \longrightarrow \wh{\calD}/\calO \longrightarrow \wh{\calD}/\calD \longrightarrow 0 
\qmq{and} 0 \longrightarrow \wh{\calO}/\calO \longrightarrow \wh{\calD}/\calO \longrightarrow \wh{\calD}/\wh{\calO} \longrightarrow 0 $$
and thereby the equations
\begin{align*}
\dim(\calD/\calO) - \dim(\wh{\calD}/\calO) + \dim(\wh{\calD}/\calD)& = 0 \\
\qmq{and} \dim(\wh{\calO}/\calO)-\dim(\wh{\calD}/\calO)+\dim(\wh{\calD}/\wh{\calO})& =0 \;, 
\end{align*}
whence
$$ m= \dim(\mathcal{D}/\mathcal{O}) = \dim(\wh{\mathcal{O}}/\mathcal{O}) + \dim(\wh{\mathcal{D}}/\wh{\mathcal{O}})- \dim(\wh{\mathcal{D}}/\mathcal{D})  $$
follows.  We have \,$\dim(\wh{\mathcal{O}}/\mathcal{O})=n$\, (this is the \,$\delta$-invariant of the singularity \,$(\lambda_*,\mu_*)$\, of \,$\Sigma$\,), \,$\dim(\wh{\calD}/\wh{\calO})=s$\, and 
\,$\dim(\wh{\mathcal{D}}/\mathcal{D})=n-j_0$\, (compare the proof of Proposition~\ref{P:spectrum:locally-free}). Thus we obtain
$$ m = n + s - (n-j_0) = s+j_0 \; . $$
\end{proof}

The following two propositions give a ``geometric'' interpretation of the data \,$\mathcal{R}$\, and \,$j_0$\, introduced by
Proposition~\ref{P:spectrum:locally-free} in the case where \,$\mathcal{D}$\, is the spectral divisor of a monodromy \,$M(\lambda)$\,. 

\label{not:spectral:eigensheaf}
We let \,$\mathcal{O}_M$\, be the sheaf of \,$(2\times 2)$-matrices of holomorphic functions in \,$\lambda$\,
which commute with the monodromy \,$M(\lambda)$\, for every \,$\lambda \in \C^*$\,. We say that a meromorphic function \,$\vi$\, on \,$\Sigma$\, 
is an \emph{eigenvalue}
of \,$\mathcal{O}_M$\,, if there exists locally a section \,$N(\lambda)$\, of \,$\mathcal{O}_M$\, so that \,$\vi$\,
is an eigenvalue of the matrix \,$N(\lambda)$\,, and we denote by \,$\mathcal{R}_M$\, the sheaf of eigenvalues of \,$\mathcal{O}_M$\,. 
It is clear that \,$\lambda \cdot \unity$\, and \,$M(\lambda)$\,
are global sections in \,$\mathcal{O}_M$\,, and therefore the functions \,$\lambda$\, and \,$\mu$\, are eigenvalues of
\,$\mathcal{O}_M$\,. Hence we have \,$\mathcal{O} \subset \mathcal{R}_M$\,. 

Conversely, if \,$N(\lambda)$\, is a \,$(2\times 2)$-matrix of holomorphic functions in \,$\lambda$\, that commutes with \,$M(\lambda)$\,,
then \,$N(\lambda)$\, acts on the eigenvector \,$\left( \begin{smallmatrix} \mu-d \\ c \end{smallmatrix} \right)$\, of \,$M(\lambda)$\, 
by
$$ N\cdot \begin{pmatrix} \mu-d \\ c \end{pmatrix} = \vi\cdot \begin{pmatrix} \mu-d \\ c \end{pmatrix} $$
with a meromorphic function \,$\vi$\, in \,$\lambda$\,, and in this situation \,$\vi$\, is a section in \,$\mathcal{R}_M$\,. 

\begin{prop}
\label{P:spectrum:specdiv-Rj0}
Let \,$\mathcal{D}$\, be the generalized spectral divisor of the monodromy \,$M(\lambda)$\, and \,$(\lambda_*,\mu_*)\in \Sigma$\,,
and let \,$\mathcal{R}_{(\lambda_*,\mu_*)}$\, and \,$j_0$\, be the data associated to \,$\mathcal{D}_{(\lambda_*,\mu_*)}$\, in Proposition~\ref{P:spectrum:locally-free}.
Then we have \,$\mathcal{R}_{(\lambda_*,\mu_*)}=\mathcal{R}_{M,(\lambda_*,\mu_*)}$\,, and \,$j_0$\, is the order of the zero of \,$M(\lambda)-\tfrac12\,\Delta(\lambda)\,\unity$\, at \,$\lambda_*$\,
(defined as the minimum of the order of the zeros of the entries of this \,$(2\times 2)$-matrix). 
\end{prop}

\begin{proof}
As before, we omit the subscript \,${}_{(\lambda_*,\mu_*)}$\,. 

By Proposition~\ref{P:spectrum:locally-free}, \,$\mathcal{R}$\, is uniquely characterized as the 
subring of \,$\wh{\calO}$\, over which \,$\mathcal{D}$\, is locally free. Because of this uniqueness property,
\,$\mathcal{R}$\, is the largest subring of \,$\wh{\calO}$\, 
that acts on \,$\mathcal{D}$\,. To show that \,$\mathcal{R}=\mathcal{R}_M$\, holds, it therefore suffices to show
that \,$\mathcal{R}_M$\, acts on \,$\mathcal{D}$\,, and that for any meromorphic, locally bounded function \,$\vi$\,
with \,$\vi\cdot \mathcal{D} \subset \mathcal{D}$\,, we have \,$\vi \in \mathcal{R}_M$\,. 


We first show that \,$\mathcal{R}_M$\, acts on \,$\mathcal{D}$\,. Let \,$\vi \in \mathcal{R}_M$\, be given. By definition,
there exists \,$N\in \mathcal{O}_M$\, so that \,$\vi \, \left( \begin{smallmatrix} f_1 \\ f_2 \end{smallmatrix} \right)
= N \cdot \left( \begin{smallmatrix} f_1 \\ f_2 \end{smallmatrix} \right)$\, holds, where \,$f_1 := \tfrac{\mu-d}{c}$\, and \,$f_2:=1$\,
are the two generators of \,$\mathcal{D}$\, over \,$\mathcal{O}$\,. Because \,$\calD$\, is a \,$\calO$-module, and the entries
of \,$N$\, are in \,$\calO$\,, it follows that \,$\vi\, f_1, \vi\,f_2 \in \mathcal{D}$\, holds. 
Because \,$f_1,f_2$\, generate \,$\mathcal{D}$\,, we in fact have \,$\vi\cdot \mathcal{D} \subset \mathcal{D}$\,.

We now show that \,$\mathcal{R}_M$\, is the largest subring of \,$\wh{\calO}$\, that acts on \,$\mathcal{D}$\,. For this purpose,
let \,$\vi$\, be a meromorphic, locally bounded function with \,$\vi\cdot \mathcal{D} \subset \mathcal{D}$\,. We will show that
\,$\vi \in \mathcal{R}_M$\, holds. We decompose \,$\vi$\, into its symmetric and its anti-symmetric part, i.e.~we write 
$$ \vi(\lambda,\mu) = \vi_+(\lambda) + \vi_-(\lambda)\cdot (\mu-\mu^{-1}) $$
with meromorphic functions \,$\vi_+$\, and \,$\vi_-$\, in \,$\lambda$\,. 
The symmetric part \,$\vi_+(\lambda)$\, is meromorphic in \,$\lambda$\, and locally bounded, and therefore holomorphic in \,$\lambda$\,.
Because of \,$\mathcal{O} \subset \mathcal{R}_M$\,, we thus have \,$\vi_+ \in \mathcal{R}_M$\,. It therefore remains to show that 
\,$\vi_-(\lambda)\cdot (\mu-\mu^{-1})\in\mathcal{R}_M$\, holds, and for this, it suffices to show that
the entries of the \,$(2\times 2)$-matrix \,$N(\lambda) := \vi_-(\lambda) \cdot (M(\lambda)-M(\lambda)^{-1})$\, are
holomorphic in \,$\lambda$\,, i.e.~that \,$N(\lambda)$\, does not have a pole.

Because the vectors \,$\left( \begin{smallmatrix} f_1(\lambda,\mu) \\ f_2(\lambda,\mu) \end{smallmatrix} \right)
= \left( \begin{smallmatrix} \tfrac{\mu-d(\lambda)}{c(\lambda)} \\ 1 \end{smallmatrix} \right)$\, 
and \,$\left( \begin{smallmatrix} f_1(\lambda,\mu^{-1}) \\ f_2(\lambda,\mu^{-1}) \end{smallmatrix} \right)
= \left( \begin{smallmatrix} \tfrac{\mu^{-1}-d(\lambda)}{c(\lambda)} \\ 1 \end{smallmatrix} \right)$\, are
linear independent in \,$\C^2$\, for any \,$(\lambda,\mu)\in \Sigma$\, with \,$\mu\neq \mu^{-1}$\, and \,$c(\lambda)\neq 0$\,,
any meromorphic function \,$h$\, can uniquely be represented in the form
\begin{equation}
\label{eq:spectrum:hg1g2}
h(\lambda,\mu) = h_1(\lambda) \cdot f_1(\lambda,\mu) + h_2(\lambda)\cdot f_2(\lambda,\mu)
\end{equation}
with meromorphic functions \,$h_1,h_2$\,. 
In this situation we have \,$h\in \mathcal{D}$\, if and only if \,$h_1$\, and \,$h_2$\, are holomorphic.
(Indeed, if \,$h_1,h_2$\, are holomorphic, then we have \,$h\in \mathcal{D}$\, simply because \,$\mathcal{D}$\,
is an \,$\mathcal{O}$-module generated by \,$f_1$\, and \,$f_2$\,. Conversely, if \,$h\in \mathcal{D}$\, holds,
we have a representation \,$h = \wt{h}_1\,f_1 + \wt{h}_2\,f_2$\, with holomorphic \,$\wt{h}_1, \wt{h}_2$\,
because \,$\mathcal{D}$\, is generated by \,$f_1$\, and \,$f_2$\,; and because of the uniqueness of the representation
\eqref{eq:spectrum:hg1g2}, we see that \,$h_1 = \wt{h}_1$\,, \,$h_2=\wt{h}_2$\, are holomorphic.)

Now we have 
\begin{equation}
\label{eq:spectrum:Nacts}
N(\lambda) \cdot \begin{pmatrix} f_1 \\ f_2 \end{pmatrix} = \vi_-(\lambda)\cdot (\mu-\mu^{-1})\cdot \begin{pmatrix} f_1 \\ f_2 \end{pmatrix} \;,
\end{equation}
where \,$(\mu-\mu^{-1})\,\left( \begin{smallmatrix} f_1 \\ f_2 \end{smallmatrix} \right)\in \mathcal{D}$\, holds. Because of the hypothesis
\,$\vi_-\cdot\mathcal{D} \subset \mathcal{D}$\,, it follows that the components of the vector in \eqref{eq:spectrum:Nacts} are in \,$\mathcal{D}$\,,
whence it follows by the preceding statement that the entries of \,$N(\lambda)$\, are holomorphic. 

This completes the proof that \,$\mathcal{R}=\mathcal{R}_M$\, holds. 

\,$j_0\geq 0$\, is by definition the largest number so that
\,$\tfrac{\mu-\mu^{-1}}{(\lambda-\lambda_*)^{j_0}}\in \mathcal{R}_M$\, holds. For any \,$j\geq 0$\,, we have \,$\tfrac{\mu-\mu^{-1}}{(\lambda-\lambda_*)^j}\in \mathcal{R}_M$\, if and only if
\,$\tfrac{1}{(\lambda-\lambda_*)^j}\cdot (M(\lambda)-M(\lambda)^{-1})$\, is holomorphic and therefore a section in \,$\mathcal{O}_M$\,. 
This is the case if and only if the entries of \,$M(\lambda)-M(\lambda)^{-1}$\, all have zeros of multiplicity at least \,$j$\,
at \,$\lambda_*$\,. Thus \,$j_0$\, is the multiplicity of the zero of \,$M(\lambda)-M(\lambda)^{-1}$\, at \,$\lambda_*$\,. 

Because of \,$\det(M(\lambda))=1$\,, we have
$$ M(\lambda)-M(\lambda)^{-1} = \begin{pmatrix} a(\lambda)-d(\lambda) & 2b(\lambda) \\ 2c(\lambda) & d(\lambda)-a(\lambda) \end{pmatrix} 
= 2\cdot \left( M(\lambda)-\tfrac12\,\Delta(\lambda)\,\unity \right)\;, $$
and therefore \,$j_0$\, is the multiplicity of the zero of \,$M(\lambda)-\tfrac12\,\Delta(\lambda)\,\unity$\, at \,$\lambda_*$\,. 
\end{proof}


The following proposition describes further properties of the spectral divisors belonging to monodromies of the kind studied here. 
Proposition~\ref{P:spectrum:dims}(1) shows that the underlying classical divisor \,$D$\, of the generalized spectral divisor \,$\calD$\, of a monodromy \,$M(\lambda)$\,
is given by Equation~\eqref{eq:spectral:D-classical}, now without the restriction that the support of \,$D$\, be contained in the regular points of \,$\Sigma$\,. 
We will use the property of Proposition~\ref{P:spectrum:dims}(2)
in Section~\ref{Se:special} as one of several properties that characterize the spectral divisors among all the generalized divisors on \,$\Sigma$\,.


\begin{prop}
\label{P:spectrum:dims}
Let \,$\mathcal{D}$\, be the generalized spectral divisor of the monodromy \,$M(\lambda)=\left( \begin{smallmatrix} a(\lambda) & b(\lambda) \\ c(\lambda) & d(\lambda) \end{smallmatrix} \right)$\, and \,$(\lambda_*,\mu_*)\in \Sigma$\, 
a singular point in the support of \,$\calD$\,. 
Let \,$m := \dim(\mathcal{D}_{(\lambda_*,\mu_*)}/\calO_{(\lambda_*,\mu_*)})$\, be the degree of \,$(\lambda_*,\mu_*)$\, in \,$\calD$\,. 
\begin{enumerate}
\item \,$m=\ord^{\C}(c)$\,.
\item There exists \,$g \in \calD_{(\lambda_*,\mu_*)}$\, so that \,$\eta := g-\tfrac{\mu-\mu^{-1}}{(\lambda-\lambda_*)^m}$\, is a meromorphic function germ solely in \,$\lambda$\,, with \,$\polord^{\Sigma}(\eta)\leq \polord^{\Sigma}(g)$\,. 
\item If \,$m=1$\,, then we have \,$j_0=1$\, for the number \,$j_0$\, from Proposition~\ref{P:spectrum:locally-free}, and every germ in \,$\calD_{(\lambda_*,\mu_*)}$\, is locally bounded 
(i.e.~we have \,$s=0$\, for the maximal pole order occurring for germs in \,$\calD_{(\lambda_*,\mu_*)}$\,). 
\end{enumerate}
\end{prop}

\begin{proof}
As before, we omit the subscript \,${}_{(\lambda_*,\mu_*)}$\,. 

\emph{For (1).}
Let us abbreviate \,$\ell:=\ord^{\C}(c)$\,, then we have \,$c = \gamma\cdot (\lambda-\lambda_*)^\ell$\, with some invertible \,$\gamma\in \mathcal{O}$\,. Because \,$\mathcal{D}$\, is generated by 
\,$1$\, and \,$\tfrac{\mu-d}{c}$\, over \,$\mathcal{O}$\,, it is therefore also generated by \,$1$\, and \,$\tfrac{\mu-d}{(\lambda-\lambda_*)^\ell}$\, over \,$\mathcal{O}$\,. Thus \,$\mathcal{D}/\mathcal{O}$\, is spanned
as a linear space by
\,$\tfrac{\mu-d}{(\lambda-\lambda_*)^j}$\, with \,$j=1,\dotsc,\ell$\,. In fact, \,$\left( \tfrac{\mu-d}{(\lambda-\lambda_*)^j}\right)_{j=1,\dotsc,\ell}$\, is a basis of \,$\mathcal{D}/\mathcal{O}$\,:
Otherwise there would exist a non-trivial linear combination of \,$\left( \tfrac{\mu-d}{(\lambda-\lambda_*)^j}\right)_{j=1,\dotsc,\ell}$\, that is a member of \,$\calO$\,; by multiplying this
linear combination with an appropriate power of \,$(\lambda-\lambda_*)$\,, we would obtain \,$\tfrac{\mu-d}{\lambda-\lambda_*} \in \calO$\,. Hence the anti-symmetric part of \,$\tfrac{\mu-d}{c}$\,,
which equals \,$\tfrac{\mu-\mu^{-1}}{2\,(\lambda-\lambda_*)}$\,, would also be a member of \,$\calO$\,. But this is a contradiction to the hypothesis that \,$(\lambda_*,\mu_*)$\, is a singular point of \,$\Sigma$\,. 
Therefore we have \,$m=\dim(\mathcal{D}/\mathcal{O})=\ell$\,. 

\emph{For (2).}
We again write \,$c=\gamma\cdot (\lambda-\lambda_*)^m$\,, then we have \,$\frac{1}{\gamma}\cdot \tfrac{\mu-d}{(\lambda-\lambda_*)^m}=\tfrac{\mu-d}{c}\in\calD$\, and therefore also \,$g := 2\cdot \tfrac{\mu-d}{(\lambda-\lambda_*)^m}
\in \calD$\,. With this choice of \,$g$\,,
$$ \eta := g - \frac{\mu-\mu^{-1}}{(\lambda-\lambda_*)^m} = \frac{2\cdot (\mu-d)-(\mu-\mu^{-1})}{(\lambda-\lambda_*)^m} = \frac{\mu+\mu^{-1}-2d}{(\lambda-\lambda_*)^m} = \frac{\Delta-2d}{(\lambda-\lambda_*)^m} $$
is a meromorphic function in \,$\lambda$\,. 

The decomposition of \,$g$\, into its even and odd parts is
$$ g = \frac{\Delta-2d}{(\lambda-\lambda_*)^m} + \frac{\mu-\mu^{-1}}{(\lambda-\lambda_*)^m} \; , $$
and therefore we have
\begin{align*}
\polord^{\Sigma}(g) & = \max\left\{ \polord^{\Sigma}\left( \frac{\Delta-2d}{(\lambda-\lambda_*)^m} \right), \polord^{\Sigma}\left( \frac{\mu-\mu^{-1}}{(\lambda-\lambda_*)^m} \right) \right\} \\
& = \max\left\{ \polord^{\Sigma}(\eta), 2m-\wh{n} \right\} \geq \polord^{\Sigma}(\eta) \; ,
\end{align*}
where \,$\wh{n}$\, denotes the order of the singularity at \,$(\lambda_*,\mu_*)$\,, i.e.~the multiplicity of the zero of \,$\Delta^2-4$\, at \,$\lambda_*$\,. 
%
%
%

\emph{For (3).} By Proposition~\ref{P:spectrum:msj0} we have \,$s+j_0=m=1$\, and therefore either \,$s=0, j_0=1$\, or \,$s=1, j_0=0$\,. Assume \,$s=1, j_0=0$\,. 
Because of \,$s>0$\,, the function \,$\tfrac{\mu-d}{c}$\, is then of maximal pole order in \,$\calD$\, and thus has a pole of order \,$s=1$\, at \,$(\lambda_*,\mu_*)$\,.
On the other hand we have
$$ \frac{\mu-d}{c} = \frac{\Delta-2d}{2c} + \frac{\mu-\mu^{-1}}{2c} \; . $$
Here it follows from (1) that \,$\ord^{\C}(c)=m=1$\, and therefore \,$\ord^{\Sigma}(c) \leq 2$\, holds. On the other hand, \,$\ord^{\Sigma}(\mu-\mu^{-1})\geq 2$\, holds because \,$(\lambda_*,\mu_*)$\, is a singular
point of \,$\Sigma$\,, and therefore \,$\tfrac{\mu-\mu^{-1}}{2c}$\, is locally bounded near \,$(\lambda_*,\mu_*)$\,. Moreover we have \,$\Delta(\lambda_*)=\pm 2 = 2d(\lambda_*)$\,,
hence \,$\Delta-2d$\, has a zero at \,$\lambda_*$\,, and therefore \,$\tfrac{\Delta-2d}{2c}$\, is also locally bounded near \,$(\lambda_*,\mu_*)$\,. It follows that \,$\tfrac{\mu-d}{c}$\, is locally bounded
near \,$(\lambda_*,\mu_*)$\,, which is a contradiction. Thus we have \,$s=0,j_0=1$\,. 
\end{proof}

\part{The asymptotic behavior of the spectral data}

\section{The vacuum solution}
\label{Se:vacuum}

The most obvious solution of the sinh-Gordon equation \,$\Delta u + \sinh(u)=0$\, is the ``vacuum'' \,$u=0$\,, corresponding to a minimal surface of zero sectional curvature, i.e.~to a flat minimal cylinder. 
In the present section, we calculate the monodromy and the spectral data of the vacuum. This example is of particular importance to us because in the coming sections, we will describe the asymptotic
behavior of the spectral data in the general situation by comparing these data to the corresponding spectral data of the vacuum.

By plugging \,$u=0$\, into Equation~\eqref{eq:mono:alphaxy} we see that the connection form \,$\alpha_0 := \alpha_{u=0}$\, corresponding to the vacuum is given by
\begin{equation}
\label{eq:vacuum:alpha0}
\alpha_0 = \frac{1}{4} \begin{pmatrix} 0 & -(1+\lambda^{-1}) \\ 1+\lambda & 0 \end{pmatrix}\mathrm{d}x + \frac{i}{4} \begin{pmatrix} 0 & 1-\lambda^{-1} \\ 1-\lambda & 0 \end{pmatrix} \mathrm{d}y \; .
\end{equation}
Because \,$\alpha_0$\, does not depend on the point \,$z=x+iy$\,, and its \,$\mathrm{d}x$\,-component and its \,$\mathrm{d}y$-component
$$ A := \frac{1}{4} \begin{pmatrix} 0 & -(1+\lambda^{-1}) \\ 1+\lambda & 0 \end{pmatrix} \qmq{resp.} \widetilde{A} := \frac{i}{4} \begin{pmatrix} 0 & 1-\lambda^{-1} \\ 1-\lambda & 0 \end{pmatrix} $$
commute, we can compute the extended frame \,$F_0=F_0(z,\lambda)$\, corresponding to the vacuum simply by
$$ F_0(x+iy,\lambda) = \exp(xA) \cdot \exp(y\widetilde{A}) \; . $$
We carry out this computation explicitly. We have
\begin{equation}
\label{eq:vacuum:zeta}
A^2 = -\frac{1}{16}\,(1+\lambda)\,(1+\lambda^{-1})\,\unity = -\zeta(\lambda)^2\cdot\unity \qmq{with} \zeta(\lambda) := \frac14\,\left({\lambda}^{1/2} + {\lambda}^{-1/2}\right) 
\end{equation}
and
\begin{equation}
\label{eq:vacuum:wtzeta}
\wt{A}^2 = -\frac{1}{16}\,(1-\lambda)(1-\lambda^{-1})\,\unity = \wt{\zeta}(\lambda)^2\cdot\unity \qmq{with} \wt{\zeta}(\lambda) := \frac14\,\left({\lambda}^{1/2} - {\lambda}^{-1/2}\right) \; ,
\end{equation}
and therefore for every \,$n\in\N_0$\,
$$ A^{2n} = (-1)^n\,\zeta(\lambda)^{2n}\cdot \unity \qmq{and} \wt{A}^{2n} = \wt{\zeta}(\lambda)^{2n}\cdot \unity $$
and
\begin{align*}
A^{2n+1} & = (-1)^n\,\zeta(\lambda)^{2n}\,A = (-1)^n\,\zeta(\lambda)^{2n+1} \left( \begin{smallmatrix} 0 & -{\lambda}^{-1/2} \\ {\lambda}^{1/2} & 0 \end{smallmatrix} \right) \\
\qmq{and} \wt{A}^{2n+1} & = \wt{\zeta}(\lambda)^{2n}\,\wt{A} = -i\,\wt{\zeta}(\lambda)^{2n+1} \left( \begin{smallmatrix} 0 & {\lambda}^{-1/2} \\ {\lambda}^{1/2} & 0 \end{smallmatrix} \right) \; .
\end{align*}
Thus we obtain
\begin{align*}
\exp(xA) & = \sum_{n=0}^\infty \frac{x^n}{n!}A^n = \sum_{n=0}^\infty \frac{x^{2n}}{(2n)!} (-1)^n \zeta(\lambda)^{2n}\unity + \sum_{n=0}^\infty \frac{x^{2n+1}}{(2n+1)!} \,(-1)^n\,\zeta(\lambda)^{2n+1} \left( \begin{smallmatrix} 0 & -{\lambda}^{-1/2} \\ {\lambda}^{1/2} & 0 \end{smallmatrix} \right) \\
& = \cos(x\,\zeta(\lambda))\cdot \unity + \sin(x\,\zeta(\lambda))\cdot \left( \begin{smallmatrix} 0 & -{\lambda}^{-1/2} \\ {\lambda}^{1/2} & 0 \end{smallmatrix} \right) 
\end{align*}
and
\begin{align*}
\exp(y\wt{A}) & = \sum_{n=0}^\infty \frac{y^n}{n!}\wt{A}^n = \sum_{n=0}^\infty \frac{y^{2n}}{(2n)!} \wt{\zeta}(\lambda)^{2n}\unity - i \sum_{n=0}^\infty \frac{y^{2n+1}}{(2n+1)!} \wt{\zeta}(\lambda)^{2n+1} \left( \begin{smallmatrix} 0 & {\lambda}^{-1/2} \\ {\lambda}^{1/2} & 0 \end{smallmatrix} \right) \\
& = \cosh(y\,\wt{\zeta}(\lambda))\cdot \unity - i\,\sinh(y\,\wt{\zeta}(\lambda))\cdot\left( \begin{smallmatrix} 0 & {\lambda}^{-1/2} \\ {\lambda}^{1/2} & 0 \end{smallmatrix} \right) \; . 
\end{align*}
Finally we obtain
\begin{align}
& F_0(x+iy,\lambda) \notag \\
= & \exp(xA) \cdot \exp(y\widetilde{A}) \notag \\
= & \left( \cos(x\,\zeta(\lambda))\cdot \unity + \sin(x\,\zeta(\lambda))\,\left( \begin{smallmatrix} 0 & -{\lambda}^{-1/2} \\ {\lambda}^{1/2} & 0 \end{smallmatrix} \right) \right) \notag \\
& \qquad\qquad \cdot \left( \cosh(y\,\wt{\zeta}(\lambda))\cdot \unity - i\,\sinh(y\,\wt{\zeta}(\lambda))\,\left( \begin{smallmatrix} 0 & {\lambda}^{-1/2} \\ {\lambda}^{1/2} & 0 \end{smallmatrix} \right) \right)\notag \\
= & \left( \begin{smallmatrix} 
    \cos(x\zeta(\lambda))\cosh(y\wt{\zeta}(\lambda)) + i\,\sin(x\zeta(\lambda))\sinh(y\wt{\zeta}(\lambda))
    & -{\lambda}^{-1/2}\,\left( i\,\cos(x\zeta(\lambda))\sinh(y\wt{\zeta}(\lambda)) + \sin(x\zeta(\lambda))\cosh(y\zeta(\lambda)) \right) 
    \\ {\lambda}^{1/2}\,\left( -i\,\cos(x\zeta(\lambda))\sinh(y\wt{\zeta}(\lambda)) + \sin(x\zeta(\lambda))\cosh(y\zeta(\lambda)) \right) 
    & \cos(x\zeta(\lambda))\cosh(y\zeta(\lambda)) - i\,\sin(x\zeta(\lambda))\sinh(y\wt{\zeta}(\lambda))
\end{smallmatrix} \right) \notag \\
\label{eq:vacuum:F0}
=& \left( \begin{matrix} \cos\left( x\zeta(\lambda)-iy\wt{\zeta}(\lambda)\right) & -{\lambda}^{-1/2} \sin\left( x\zeta(\lambda) + iy\wt{\zeta}(\lambda) \right) \\
    {\lambda}^{1/2}\cdot \sin\left( x\zeta(\lambda) - iy\wt{\zeta}(\lambda) \right) & \cos\left( x\zeta(\lambda)+iy\wt{\zeta}(\lambda)\right)
\end{matrix} \right) \;. 
\end{align}
Note that all the entries of \,$F_0$\, are even in \,${\lambda}^{1/2}$\,, and therefore indeed define holomorphic functions in \,$\lambda\in\C^*$\,.

In particular, we have for the monodromy of the vacuum with respect to the base point \,$z_0=0$\,
\begin{equation}
\label{eq:vacuum:M0}
M_0(\lambda)  = F_{0}(1,\lambda) = \begin{pmatrix} \cos(\zeta(\lambda)) & -{\lambda}^{-1/2}\,\sin(\zeta(\lambda)) \\ {\lambda}^{1/2}\,\sin(\zeta(\lambda)) & \cos(\zeta(\lambda)) \end{pmatrix} 
=: \begin{pmatrix} a_0(\lambda) & b_0(\lambda) \\ c_0(\lambda) & d_0(\lambda) \end{pmatrix}\; . 
\end{equation}
We will use the names \,$a_0,\dotsc,d_0$\, for the component functions of the monodromy of the vacuum throughout the entire book without any further reference,
and likewise
\begin{equation}
\label{eq:vacuum:Delta0}
\Delta_0(\lambda) := \tr(M_0(\lambda)) = 2\cos(\zeta(\lambda)) \; .
\end{equation}

As consequence of Equation~\eqref{eq:spectral:Sigma2}, the spectral curve of the vacuum is given by
\begin{align}
\Sigma_0 & = \Mengegr{(\lambda,\mu)\in \C^* \times \C}{\mu = \tfrac12\left( \Delta_0(\lambda) \pm \sqrt{\Delta_0(\lambda)^2-4} \right)} \notag \\
& = \Mengegr{(\lambda,\mu)\in \C^* \times \C}{\mu = \cos(\zeta(\lambda)) \pm \sqrt{\cos(\zeta(\lambda))^2-1}} \notag \\
\label{eq:vacuum:Sigma0}
& = \Mengegr{(\lambda,\mu)\in \C^* \times \C}{\mu = e^{\pm i\,\zeta(\lambda)} } \; . 
\end{align}
This curve has no branch points above \,$\C^*$\,. It has double points at all those \,$\lambda\in \C^*$\, for which
\,$\zeta(\lambda)$\, is an integer multiple of \,$\pi$\,; these values of \,$\lambda$\, are exactly the following:
\begin{equation}
\label{eq:vacuum:lambdak0-def}
\lambda_{k,0} := 8\pi^2k^2 + 4\pi k \sqrt{4\pi^2k^2-1} - 1 \qmq{with \,$k\in \Z$\,.}
\end{equation}

\begin{proof}
Let us choose for \,$\sqrt{\cdots}$\, the standard branch of the square root function (so that we have
\,$\RE(\sqrt{\lambda})\geq 0$\, for all \,$\lambda\in \C$\,). Then the equation \,$\zeta(\lambda)=k\pi$\, with \,$k\in \Z$\, has a solution
only if \,$k\geq 0$\,. In this case we have
\,$\zeta(\lambda)=\tfrac14(\lambda^{1/2}+\lambda^{-1/2})=k\pi$\,
if and only if \,$\lambda - 4k\pi\sqrt{\lambda}+1=0$\, holds. The two possible values of \,$\sqrt{\lambda}$\,
that correspond to this equation are \,$\sqrt{\lambda}= 2k\pi \pm \sqrt{4k^2\pi^2-1}$\,, and this yields
$$ \lambda = \left( 2k\pi \pm \sqrt{4k^2\pi^2-1} \right)^2 = 8\pi^2k^2 \pm 4\pi k \sqrt{4\pi^2k^2-1} - 1 \;. $$
Therefore the entirety of solutions of \,$\zeta(\lambda) \in \Z\,\pi$\, is given by \eqref{eq:vacuum:lambdak0-def},
where now \,$k$\, again runs through all of \,$\Z$\,. 
\end{proof}


We have \,$\lambda_{k,0}\in \R$\, for all \,$k\in \Z$\,, and \,$|\lambda_{k,0}| > 1$\, resp.~\,$|\lambda_{k,0}|<1$\, if \,$k>0$\, resp.~\,$k<0$\,. Moreover, \,$\lambda_{k,0}$\, tends to \,$\infty$\, resp.~to \,$0$\,
for \,$k\to\infty$\, resp.~\,$k\to-\infty$\,. By using the Taylor expansion \,$\sqrt{1+x} = 1 + \tfrac12 x - \tfrac18 x^2 + \tfrac{1}{16}x^3 + O(x^4)$\,, we obtain the more specific asymptotic assessments:
\begin{align}
\lambda_{k,0} & = 16\pi^2 k^2 - 2 + O(k^{-2}) \text{ for \,$k\to\infty$\,} \notag \\
\label{eq:vacuum:lambdak0-asymp}
\qmq{and} \lambda_{k,0} & = \frac{1}{16\pi^2}k^{-2} + \frac{1}{128\pi^4}k^{-4} + O(k^{-6}) \text{ for \,$k\to -\infty$\,.} 
\end{align}

We finally calculate the spectral divisor of the vacuum. The points
\,$(\lambda,\mu)\in \Sigma_0$\, of the classical divisor are determined by the equations
\,$c_0(\lambda) := \sqrt{\lambda}\cdot \sin(\zeta(\lambda))=0$\, and \,$\mu = a_0(\lambda) := \cos(\zeta(\lambda))$\,.
Again \,$c_0(\lambda)=0$\, holds if and only if \,$\lambda=\lambda_{k,0}$\, for some \,$k\in \Z$\,, 
all these zeros of \,$c_0$\, are of simple multiplicity, and we have 
\begin{equation}
\label{eq:vacuum:muk0-def}
\mu_{k,0}:=a_0(\lambda_{k,0}) = \cos(k\pi)=(-1)^k \; . 
\end{equation}
Thus, the classical spectral divisor of the vacuum is given by the divisor on \,$\Sigma_0$\,
\begin{equation}
\label{eq:vacuum:D0}
D_0 := \Menge{(\lambda_{k,0},\mu_{k,0})}{k\in \Z} = \Menge{(\lambda_{k,0},(-1)^k)}{k\in\Z} \; .
\end{equation}
Notice that for the vacuum, the support of the spectral divisor coincides with the set of double points of the spectral curve.

\label{not:vacuum:gen-D0}
To determine the generalized spectral divisor \,$\calD_0$\, of the vacuum, we look at the associated data \,$\calR$\, and \,$j_0$\, from Proposition~\ref{P:spectrum:locally-free}. 
Because every point \,$(\lambda_{k,0},\mu_{k,0})$\, in the support of \,$\calD_0$\, is of degree \,$m=1$\,, Proposition~\ref{P:spectrum:dims}(3)
shows that we have \,$j_0=1$\,, and thus \,$j_0$\, equals the \,$\delta$-invariant of that point in \,$\Sigma$\,. Hence we have \,$\calR_{(\lambda_{k,0},\mu_{k,0})} = \wh{\calO}_{(\lambda_{k,0},\mu_{k,0})}$\,.
By Proposition~\ref{P:spectrum:dims}(3) the maximal pole order \,$s$\, occurring in \,$(\calD_0)_{(\lambda_{k,0},\mu_{k,0})}$\, is \,$s=0$\,, and therefore
\,$(\calD_0)_{(\lambda_{k,0},\mu_{k,0})}$\, is generated by \,$1$\, over \,$\calR_{(\lambda_{k,0},\mu_{k,0})} = \wh{\calO}_{(\lambda_{k,0},\mu_{k,0})}$\,. Thus
\,$(\calD_0)_{(\lambda_{k,0},\mu_{k,0})}  = (\wh{\calO}_0)_{(\lambda_{k,0},\mu_{k,0})}$\, holds, where \,$\wh{\calO}_0$\, is the direct image in \,$\Sigma_0$\, of the sheaf of holomorphic functions on the normalization \,$\wh{\Sigma}_0$\,
of \,$\Sigma_0$\,. Therefore we have \,$\calD_0 = \wh{\calO}_0$\,.

\section{The basic asymptotic of the monodromy}
\label{Se:asymp}

In the present section, we will prove the basic asymptotic estimates for the monodromy. We will refine our asymptotic assessment
in Section~\ref{Se:fasymp} and again in Section~\ref{Se:asympfinal}. The asymptotic estimate of Theorem~\ref{T:asymp:basic} 
(resp.~Proposition~\ref{P:asymp:more}, where we require one more degree of differentiability) is
analogous to the basic estimates of \cite{Poeschel-Trubowitz:1987}, Theorem~1.3, p.~13 in the treatment of the 1-dimensional Schr\"odinger equation.

Beginning with this section, we suppose that for the potential (simply periodic solution of the sinh-Gordon equation) \,$u$\,, only 
periodic \emph{Cauchy data} \,$(u,u_y)$\, on the real line are given. As explained in the Introduction, we want the requirements on the
differentiability of \,$u$\, and \,$u_y$\, to be as relaxed as possible. Specifically, we only require that \,$u$\, is in the Sobolev-space 
of weakly once-differentiable functions with square-integrable derivative, i.e.~\,$u \in W^{1,2}([0,1])$\, and that \,$u_y$\, is
square-integrable, i.e.~\,$u_y \in L^2([0,1])$\,. Note that \,$u$\, is in particular continuous, so individual function values \,$u(x)$\, of \,$u$\, are well-defined. 

We regard \,$u$\, and \,$u_y$\, as being extended periodically to the real line, and we
define ``mixed derivatives'' of \,$u$\, in the natural way by using both \,$u$\, and \,$u_y$\,, e.g.~\,$u_z := \tfrac12
(u_x-i\,u_y)$\,, where \,$u_x$\, is the Sobolev-derivative of \,$u$\,, and \,$u_y$\, is the given function. 

\label{not:asymp:Potnp}
We denote the space of such \emph{non-periodic potentials} by \,$\Pot_{np} := \Menge{(u,u_y)}{u \in W^{1,2}([0,1]), u_y \in L^2([0,1])}$\,.
Via the norm
$$ \|(u,u_y)\|_{\Pot} := \sqrt{\|u\|_{W^{1,2}}^2 + \|u_y\|_{L^2}^2} \;, $$
\,$\Pot_{np}$\, becomes a Banach space. 

Because our ultimate interest is with periodic solutions of the sinh-Gordon equation, we are of course most interested in \emph{periodic potentials}. 
\,$u$\, and \,$u_y$\, are at first only defined on \,$[0,1]$\,, so the condition of periodicity for \,$u$\, is simply \,$u(0)=u(1)$\,; this condition
is well-defined because \,$u$\, is continuous. \,$u_y$\, is only square-integrable, so we cannot access individual function values of \,$u_y$\,;
for this reason we do not impose a similar condition of periodicity on \,$u_y$\,. We then regard \,$u$\, and \,$u_y$\, also as extended periodically to the real line.
The space of \emph{periodic potentials} is thus given by
\begin{equation}
\label{eq:asymp:Pot}
\Pot := \Mengegr{ (u,u_y) \in \Pot_{np}}{u(0)=u(1)} \,,
\end{equation}
it is a complex hyperplane in \,$\Pot_{np}$\,. 

Also for such Cauchy data \,$(u,u_y)\in \Pot_{np}$\, (in the place of a solution \,$u$\, of the sinh-Gordon equation defined on an entire horizontal strip in \,$\C$\,),
we can define the \,$\mathrm{d}x$-part of the flat connection \,$1$-form \,$\alpha$\,
\begin{equation}
\label{eq:asymp:alpha}
\alpha = \frac{1}{4} \begin{pmatrix} i\,u_y & -e^{u/2} - \lambda^{-1}\,e^{-u/2} \\ e^{u/2}+\lambda\,e^{-u/2} & -i\,u_y \end{pmatrix} \mathrm{d}x \;.
\end{equation}
Although \,$\alpha$\, is in general only square-integrable and not locally Lipschitz continuous (as a consequence of our conditions of differentiability on \,$(u,u_y)$\,),
the following Lemma~\ref{L:asymp:existence} shows that we still have an
extended frame \,$F_\lambda: \R \to \SL(2,\C)$\, along the real line associated to \,$(u,u_y)$\,, i.e.~a solution to the initial value problem
\begin{equation}
\label{eq:asymp:Fode}
F_\lambda'(x) = \alpha_\lambda(x)\cdot F_\lambda(x) \qmq{and} F_\lambda(0)=\unity \;. 
\end{equation}
\label{not:asymp:M}
From \,$F_\lambda$\, we obtain the monodromy at \,$x=0$\, as \,$M(\lambda) := F_\lambda(1)$\,. If \,$(u,u_y)$\, is periodic (\,$(u,u_y)\in\Pot$\,),
we define spectral data \,$(\Sigma,D)$\, resp.~\,$(\Sigma,\calD)$\, for \,$(u,u_y)$\, exactly as in Section~\ref{Se:spectrum}. If \,$(u,u_y)$\, is not periodic,
we do not define a spectral curve \,$\Sigma$\, for \,$(u,u_y)$\, (if one were to define \,$\Sigma$\, by Equation~\ref{eq:spectral:Sigma}, it would turn out
that the branch points of \,$\Sigma$\, do not satisfy any reasonable asymptotic law, so such a definition would be of no interest), but we still
define the classical spectral divisor \,$D$\, of \,$(u,u_y)$\, as a point multi-set in \,$\C^* \times \C^*$\, by Equation~\eqref{eq:spectral:D-classical}. 

\begin{rem}
Lemma~\ref{L:asymp:existence} shows that for the existence of the extended frame \,$F_\lambda$\, as a solution of \eqref{eq:asymp:Fode}, it in fact suffices for \,$\alpha$\, to be integrable
(instead of square-integrable). 
Most of the asymptotic estimates of the present section, and of Section~\ref{Se:fasymp} could also be carried out in an analogous way, 
if the space of potentials were extended from \,$W^{1,2}([0,1])\times L^2([0,1])$\, to \,$W^{1,p}([0,1]) \times L^p([0,1])$\, with a fixed \,$p>1$\,. 
We will phrase several of the following lemmas (but not the full proof of the asymptotic estimates) in such a way that this possibility becomes apparent.
\end{rem}

In what follows, we denote by \,$|\dotsc|$\, also an arbitrary sub-multiplicative matrix norm on the space \,$\C^{2\times 2}$\, of complex \,$(2\times 2)$-matrices.

\begin{lem}
\label{L:asymp:existence}
Let \,$\alpha \in L^p([0,1],\C^{2\times 2})$\, with \,$p\geq 1$\,. Then the sum
$$ F: [0,1] \to \C^{2\times 2},\; x\mapsto \unity + \sum_{n=1}^\infty \int_{t_1=0}^x \int_{t_2=0}^{t_1} \cdots
\int_{t_n=0}^{t_{n-1}} \alpha(t_1)\,\alpha(t_2)\,\dotsc \alpha(t_n)\,\mathrm{d}^n t $$
converges in \,$W^{1,p}([0,1],\C^{2\times 2})$\, to the solution of the initial value problem
$$ F' = \alpha\,F \qmq{and} F(0)=\unity \; , $$
and we have \,$|F(x)| \leq \exp(\|\alpha\|_{L^1})$\, for every \,$x\in[0,1]$\,. 
\end{lem}

\begin{proof}
For every \,$n \geq 1$\,, we put
$$ F_n(x) := \int_{t_1=0}^x \int_{t_2=0}^{t_1} \cdots \int_{t_n=0}^{t_{n-1}} \alpha(t_1)\,\alpha(t_2)\,\dotsc \alpha(t_n)\,\mathrm{d}^n t \; . $$
We then have \,$F_n \in W^{1,p}([0,1],\C^{2\times 2})$\, with
$$ F_n'(x) = \alpha(x) \cdot F_{n-1}(x) $$
(where we put \,$F_0(x) := \unity$\,). Moreover, we have 
\begin{align*}
|F_n(x)| & = \left| \int_{t_1=0}^x \int_{t_2=0}^{t_1} \cdots
\int_{t_n=0}^{t_{n-1}} \alpha(t_1)\,\alpha(t_2)\,\dotsc \alpha(t_n)\,\mathrm{d}^n t \right| \\
& \leq \int_{t_1=0}^x \int_{t_2=0}^{t_1} \cdots \int_{t_n=0}^{t_{n-1}} |\alpha(t_1)| \cdot |\alpha(t_2)|
\cdot \dotsc \cdot |\alpha(t_n)| \mathrm{d}^n t \\
& = \frac{1}{n!}\,\int_{t_1=0}^x \int_{t_2=0}^{x} \cdots \int_{t_n=0}^{x} |\alpha(t_1)| \cdot |\alpha(t_2)|
\cdot \dotsc \cdot |\alpha(t_n)| \mathrm{d}^n t \\
& = \frac{1}{n!}\, \left( \int_0^x |\alpha(t)|\,\mathrm{d}t \right)^n \leq \frac{1}{n!}\, \|\alpha\|_{L^1}^n \; . 
\end{align*}
It follows that the series \,$\sum_{n\geq 1} F_n(x)$\, converges absolutely and uniformly for \,$x\in [0,1]$\,,
whereas the series \,$\sum_{n\geq 1} F_n' = \alpha \cdot \sum_{n\geq 0} F_n$\, converges in \,$L^p([0,1],\C^{2\times 2})$\,. 
Therefore the series defining \,$F$\, converges in \,$W^{1,p}([0,1])$\,, and we have \,$F' = \alpha\cdot F$\,
and \,$F(0)=\unity$\,. Moreover, it follows from the preceding estimate that 
\,$|F(x)| \leq \exp(\|\alpha\|_{L^1})$\, holds for every \,$x\in[0,1]$\,. 
\end{proof}


We now turn to the description of the asymptotic behavior of the monodromy \,$M(\lambda)$\,. For
estimating the error of the asymptotic approximation of \,$M(\lambda)$\,, the function
\begin{equation}
\label{eq:asymp:w}
w(\lambda) := |\cos(\zeta(\lambda))| + |\sin(\zeta(\lambda))|
\end{equation}
will be important. We jot down a few facts about this function:

\begin{prop}
\label{P:asymp:w}
\begin{enumerate}
\item We have \,$\tfrac12\,e^{|\IM(\zeta(\lambda))|} \leq w(\lambda) \leq 2\,e^{|\IM(\zeta(\lambda))|}$\, for all \,$\lambda\in \C^*$\,. 
\item \,$w$\, is bounded on every horizontal strip in the \,$\zeta$-plane.
\end{enumerate}
\end{prop}

\begin{proof}
\emph{For (1).}
For any \,$\lambda\in\C^*$\,, we have
\begin{align*}
|\cos(\zeta(\lambda))| 
& = \left| \frac12\left(e^{i\zeta(\lambda)}+e^{-i\zeta(\lambda)}\right)\right| \leq \max \left\{ \left| e^{i\zeta(\lambda)} \right|,  \left| e^{-i\zeta(\lambda)} \right| \right\} \\
& = \max \left\{ e^{\IM \zeta(\lambda)}, e^{-\IM \zeta(\lambda)} \right\} = e^{|\IM \zeta(\lambda)|} 
\end{align*}
and likewise
$$ |\sin(\zeta(\lambda))| \leq e^{|\IM \zeta(\lambda)|} \;, $$
hence
$$ w(\lambda) = |\cos(\zeta(\lambda))| + |\sin(\zeta(\lambda))| \leq 2\,e^{|\IM(\zeta(\lambda))|} \; . $$
Concerning the lower bound for \,$w(\lambda)$\,, we have
\begin{align}
w(\lambda) & = |\cos(\zeta(\lambda))| + |\sin(\zeta(\lambda))| = \left| \frac12\left(e^{i\zeta(\lambda)}+e^{-i\zeta(\lambda)}\right) \right| + \left| \frac{1}{2i}\left(e^{i\zeta(\lambda)}-e^{-i\zeta(\lambda)} \right) \right| \notag \\
& = \frac12\,|e^{i\zeta(\lambda)}|\cdot \left( \left|1+e^{-2i\zeta(\lambda)}\right| + |1-e^{-2i\zeta(\lambda)}| \right) \notag \\
\label{eq:vacuum:zeta:M0lowerpre}
& = \frac12\,e^{|\IM(\zeta(\lambda))|}\cdot \left( \left|1+e^{x+iy}\right| + |1-e^{x+iy}| \right) 
\end{align}
with \,$x := \RE(-2i\zeta(\lambda)),\;y := \IM(-2i\zeta(\lambda))$\,. We further have 
\begin{align*}
\left|1\pm e^{x+iy}\right|^2 & = \left| (1\pm e^x\,\cos(y)) \pm i\cdot e^x\,\sin(y) \right|^2 = (1\pm e^x\,\cos(y))^2 + (e^x\,\sin(y))^2 \\
& = 1 \pm 2\,e^x\,\cos(y) + e^{2x} \geq 1 \pm 2\,e^x\,\cos(y) \; . 
\end{align*}
Depending on the sign of \,$\cos(y)$\,, for at least one choice of the sign \,$\pm$\,, we have \,$1 \pm 2\,e^x\,\cos(y) \geq 1$\,, and therefore it follows from \eqref{eq:vacuum:zeta:M0lowerpre} that
\,$w(\lambda) \geq \tfrac12\,e^{|\IM(\zeta(\lambda))|}$\, holds.

\emph{For (2).} This is an immediate consequence of (1).
\end{proof}

\begin{thm}
\label{T:asymp:basic}
Let \,$(u,u_y)\in \Pot_{np}$\, be given. We put \,$\tau := e^{-(u(0)+u(1))/4}$\, and \,$\upsilon := e^{(u(1)-u(0))/4}$\,. (Note that we have \,$\upsilon=1$\, for \,$(u,u_y)\in \Pot$\,.) 

We compare the monodromy  \,$M(\lambda) = \begin{pmatrix} a(\lambda) & b(\lambda) \\ c(\lambda) & d(\lambda) \end{pmatrix}$\, of \,$(u,u_y)$\, to the monodromy
 \,$M_0(\lambda) = \begin{pmatrix} a_0(\lambda) & b_0(\lambda) \\ c_0(\lambda) & d_0(\lambda) \end{pmatrix}$\, of the vacuum (see Equation~\eqref{eq:vacuum:M0}).
For every \,$\eps>0$\, there exists \,$R>0$\,, such that:
\begin{enumerate}
\item
For all \,$\lambda\in\C$\, with \,$|\lambda|\geq R$\,, we have
\begin{align*}
|a(\lambda)-\upsilon\,a_0(\lambda)| & \leq \eps \, w(\lambda) \\
|b(\lambda)-\tau^{-1}\,b_0(\lambda)| & \leq \eps \, |\lambda|^{-1/2}\,w(\lambda) \\
|c(\lambda)-\tau\,c_0(\lambda)| & \leq \eps \, |\lambda|^{1/2}\,w(\lambda) \\
|d(\lambda)-\upsilon^{-1}\,d_0(\lambda)| & \leq \eps \, w(\lambda) \; . 
\end{align*}
\item
For all \,$\lambda\in\C^*$\, with \,$|\lambda|\leq \tfrac{1}{R}$\,, we have
\begin{align*}
|a(\lambda)-\upsilon^{-1}\,a_0(\lambda)| & \leq \eps \, w(\lambda) \\
|b(\lambda)-\tau\,b_0(\lambda)| & \leq \eps \, |\lambda|^{-1/2}\,w(\lambda) \\
|c(\lambda)-\tau^{-1}\,c_0(\lambda)| & \leq \eps \, |\lambda|^{1/2}\,w(\lambda) \\
|d(\lambda)-\upsilon\,d_0(\lambda)| & \leq \eps \, w(\lambda) \; . 
\end{align*}
\end{enumerate}
Moreover if \,$P$\, is a relatively compact subset of \,$\Pot_{np}$\,, then \,$R$\, can in (1) and (2) be chosen uniformly (in dependence of \,$\eps$\,) for \,$(u,u_y)\in P$\,. 
\end{thm}


Most of the remainder of the section is concerned with the proof of this theorem. One important instrument for the proof
will be the following lemma, concerning the comparison of solutions of the differential equations \,$\mathrm{d}F=\alpha F$\,
and \,$\mathrm{d}F=(\alpha+\beta)F$\, with suitable \,$\alpha,\beta$\,.

\begin{lem}
\label{L:asymp:compare}
Suppose that \,$p\geq 1$\, is fixed, and \,$\alpha,\beta \in L^p([0,1],\C^{2\times 2})$\, are given. We put \,$\wt{\alpha} := \alpha+\beta$\,,
and let \,$F, \wt{F} \in W^{1,p}([0,1],\C^{2\times 2})$\, be the solutions of
\,$F'=\alpha F$\, and \,$\wt{F}'=\wt{\alpha}\wt{F}$\, with \,$F(0)=\wt{F}(0)=\unity$\,,
see Lemma~\ref{L:asymp:existence}.

\begin{enumerate}
\item
We have
\begin{equation}
\label{eq:asymp:compare:wtF-series}
\wt{F}(x) = F(x) + \sum_{n=1}^\infty \wt{F}_n(x) \;, 
\end{equation}
where we put for \,$n\geq 1$\,
$$ \wt{F}_n(x) = F(x)\cdot \int_{t_1=0}^x \int_{t_2=0}^{t_1} \cdots \int_{t_n=0}^{t_{n-1}} \prod_{j=1}^n F(t_{j})^{-1}\,\beta(t_j)\,F(t_j)\,\mathrm{d}^nt \;, $$
and the infinite sum in Equation~\eqref{eq:asymp:compare:wtF-series} converges in \,$W^{1,p}([0,1],\C^{2\times 2})$\,. 
Moreover, there exists a constant \,$C_{\|\alpha\|}>0$\, depending only on \,$\|\alpha\|_{L^1}$\,, such that for all \,$x\in [0,1]$\,
we have
\begin{equation}
\label{eq:asymp:compare:wtF-estimate}
|\wt{F}(x)| \leq \exp(C_{\|\alpha\|} \cdot \|\beta\|_{L^1}) \cdot |F(x)| \; .
\end{equation}
If we fix \,$R>0$\,, then 
the convergence of the sum in \eqref{eq:asymp:compare:wtF-series}
is  uniform for all \,$\alpha,\beta\in L^p([0,1],\C^{2\times 2})$\, with \,$\|\alpha\|_{L^1},\|\beta\|_{L^1} \leq R$\,.
\item
In the previous setting, we now suppose that \,$\alpha$\, is constant with respect to \,$x$\,, and that \,$\alpha^2 = a\cdot \unity$\,
holds with some \,$a\in \C^*$\,. Then we have:
\begin{enumerate}
\item
\,$F(x)=\exp(x\,\alpha)$\, for all \,$x\in [0,1]$\,. Therefore \,$F$\, extends to the group homomorphism \,$\R\to\GL(2,\C),\;x\mapsto F(x):=\exp(x\,\alpha)$\,, i.e.~we have
besides \,$F(0)=\unity$\, for all \,$x,y\in \R$\,
\begin{equation}
\label{eq:asymp:compare:F-homo}
F(x+y) = F(x)\cdot F(y) \qmq{and} F(x)^{-1} = F(-x) \; . 
\end{equation}
Moreover, \,$|F(x)| \leq \exp(x\,|\alpha|)$\, holds.
\item
There exist unique \,$\beta_+,\beta_-\in L^p([0,1],\C^{2\times 2})$\, with \,$\beta = \beta_+ + \beta_-$\, and such that
\,$\beta_+$\, commutes with \,$\alpha$\, and \,$\beta_-$\, anti-commutes with \,$\alpha$\,. We also have the following
rule of commutation:
\begin{equation}
\label{eq:asymp:compare:F-beta}
\beta_\pm(x) \cdot F(y) = F(\pm y) \cdot \beta_\pm(x) \; . 
\end{equation}

Then we have for every \,$n\geq 1$\,
$$ \wt{F}_n(x) = \sum_{\eps \in \{\pm 1\}^n} \int_{t_1=0}^x \int_{t_2=0}^{t_1} \cdots \int_{t_n=0}^{t_{n-1}} F(\xi_\eps(t))\,\beta_{\eps_1}(t_1)
\, \cdots \, \beta_{\eps_n}(t_n)\,\mathrm{d}^nt \;, $$
where for every \,$\eps := (\eps_1,\dotsc,\eps_n)\in \{\pm 1\}^n$\, and \,$t=(t_1,\dotsc,t_n)\in [0,x]^n$\,
with \,$0\leq t_n \leq t_{n-1} \leq \dotsc \leq t_1 \leq x$\,, we put
$$ \xi_\eps(t) := x-2t_{\nu_1} + 2t_{\nu_2} -+\dotsc + 2(-1)^jt_{\nu_j} \;, $$
where \,$1 \leq \nu_1 <\dotsc<\nu_j \leq n$\, are those indices \,$\nu\in \{1,\dotsc,n\}$\, for which we have \,$\eps_\nu=-1$\,. 
We have \,$\xi_{\eps}(t) \in [-x,x]$\,. 
\item
Now suppose that \,$\beta^{[1]},\beta^{[2]} \in L^p([0,1],\C^{2\times 2})$\, are given. We apply the situation of (b)
to \,$\beta^{[\nu]}$\, (for \,$\nu\in \{1,2\}$\,), denoting the quantities corresponding to \,$\beta^{[\nu]}$\, by the superscript \,${}^{[\nu]}$\,. 
Then there exists a constant \,$C>0$\,, depending only on an upper bound for \,$\|\max\{|\beta^{[\nu]}_\pm|\}\|_{L^1}$\,, such that
$$ |\wt{F}^{[1]}(x)-\wt{F}^{[2]}(x)| \leq C \cdot |F(x)| \cdot \bigr\|\max\{|\beta^{[1]}_+-\beta^{[2]}_+|, |\beta^{[1]}_--\beta^{[2]}_-|\}\bigr\|_{L^1} $$
holds.
\end{enumerate}
\end{enumerate}
\end{lem}

\begin{proof}
\emph{For (1).}
We begin by showing the claims on the convergence of the infinite sum in \eqref{eq:asymp:compare:wtF-series}. We note that we have for any \,$n\geq 1$\,:
\,$\wt{F}_n\in W^{1,p}([0,1],\C^{2\times 2})$\, and 
\begin{align}
\wt{F}_n'(x) & = F'(x) \cdot \int_{t_1=0}^x \int_{t_2=0}^{t_1} \cdots \int_{t_n=0}^{t_{n-1}} \prod_{j=1}^n F(t_{j})^{-1}\,\beta(t_j)\,F(t_j)\,\mathrm{d}^nt \notag \\
& \qquad\qquad + F(x) \cdot F(x)^{-1}\,\beta(x)\,F(x)\,\int_{t_2=0}^{x} \cdots \int_{t_n=0}^{t_{n-1}} \prod_{j=2}^n F(t_{j})^{-1}\,\beta(t_j)\,F(t_j)\,\mathrm{d}^{n-1}t \notag \\
\label{eq:asymp:compare:wtF'}
& = \alpha(x)\,\wt{F}_{n}(x) + \beta(x)\,\wt{F}_{n-1}(x)
\end{align}
(with \,$\wt{F}_0 := F$\,). Moreover, for \,$x\in [0,1]$\,, we have
\begin{align*}
|\wt{F}_n(x)| & \leq |F(x)|\cdot \int_{t_1=0}^x \int_{t_2=0}^{t_1} \cdots \int_{t_n=0}^{t_{n-1}} \prod_{j=1}^n |F(t_{j})^{-1}| \cdot |\beta(t_j)|\cdot |F(t_j)|\, \mathrm{d}^nt \\
& \leq |F(x)|\cdot \frac{1}{n!} \, \int_{t_1=0}^x \int_{t_2=0}^{x} \cdots \int_{t_n=0}^{x} \prod_{j=1}^n |F(t_{j})^{-1}| \cdot |\beta(t_j)|\cdot |F(t_j)|\, \mathrm{d}^nt \\
& = |F(x)|\cdot \frac{1}{n!} \, \left( \int_0^x |F(t)^{-1}| \cdot |\beta(t_j)|\cdot |F(t)|\, \mathrm{d}t \right)^n \\
& \leq |F(x)|\cdot \frac{1}{n!}\,\left(\exp(\|\alpha\|_{L^1})^2 \, \|\beta\|_{L^1}\right)^n \;,
\end{align*}
where the last inequality follows from Lemma~\ref{L:asymp:existence}. It follows that the series \,$\sum_{n=1}^\infty \wt{F}_n(x)$\,
converges absolutely and uniformly in \,$x$\,, whereas the series
$$ \sum_{n=1}^\infty \wt{F}_n'(x) = \alpha(x) \cdot \sum_{n=1}^\infty \wt{F}_n(x) + \beta(x) \cdot \sum_{n=0}^\infty \wt{F}_n(x) $$
converges in \,$L^p([0,1],\C^{2\times 2})$\,. Therefore \,$\sum_{n=1}^\infty \wt{F}_n$\, converges in the Sobolev space \,$W^{1,p}([0,1],\C^{2\times 2})$\,. 
Moreover, these convergences are also uniform when 
 \,$\alpha$\, and \,$\beta$\, vary with \,$\|\alpha\|_{L^1}, \|\beta\|_{L^1} \leq R$\, (for a fixed \,$R>0$\,). Moreover, we have 
\,$1+\sum_{n=1}^\infty |\wt{F}_n(x)| \leq |F(x)|\cdot \exp(C_{\|\alpha\|}\cdot \|\beta\|_{L^1})$\, with \,$C_{\|\alpha\|} := \exp(\|\alpha\|_{L^1})^2$\,.
 
Therefore Equation~\eqref{eq:asymp:compare:wtF-series} defines an element \,$\wt{F} \in W^{1,p}([0,1],\C^{2\times 2})$\,. 
We have \,$\wt{F}_n(0)=0$\, for all \,$n\geq 1$\,, and therefore
\,$\wt{F}(0)=F(0)=\unity$\,. It remains to show that \,$\wt{F}$\, solves the differential equation 
\,$\wt{F}' = \wt{\alpha}\wt{F}$\,. Indeed, because the series 
\,$\sum_{n=1}^\infty \wt{F}_n$\, converges in \,$W^{1,p}([0,1],\C^{2\times 2})$\,, we have
\begin{align*}
\wt{F}'(x)
& = F'(x) + \sum_{n=1}^\infty \wt{F}_n'(x) \overset{\eqref{eq:asymp:compare:wtF'}}{=} \alpha(x)\,F(x) + \sum_{n=1}^\infty \bigr( \alpha(x)\,\wt{F}_n(x) + \beta(x)\,\wt{F}_{n-1}(x) \bigr) \\
& = (\alpha(x)+\beta(x))\cdot \left( F(x) + \sum_{n=1}^\infty \wt{F}_n(x)\right) = \wt{\alpha}(x)\cdot \wt{F}(x) \;. 
\end{align*}

\emph{For (2).} (a) is obvious. For (b), we put 
$$ \beta_{\pm} := \frac12\,(\beta \pm \alpha\,\beta\,\alpha^{-1}) \; . $$
Then we clearly have \,$\beta_++\beta_-=\beta$\,. Because of the hypothesis \,$\alpha^2=a\cdot \unity$\,, we also have
\,$\alpha^{-1}=\tfrac{1}{a}\cdot \alpha$\,, and therefore
$$ \alpha \cdot \beta_\pm = \frac12\,(\alpha\,\beta \pm \alpha^2\,\beta\,\alpha^{-1}) = \frac12\,(\alpha\,\beta \pm \beta\,\alpha)
= \pm \frac12\,(\beta \pm \alpha\,\beta\,\alpha^{-1})\,\alpha = \pm \beta_\pm\cdot \alpha \;, $$
i.e.~\,$\alpha$\, commutes with \,$\beta_+$\, and anti-commutes with \,$\beta_-$\,; herefrom also Equation~\eqref{eq:asymp:compare:F-beta} follows.
It is clear that this decomposition of \,$\beta$\, is unique.

We now calculate \,$\wt{F}_n$\, in the present situation, using the properties of \,$F$\, given by 
Equations~\eqref{eq:asymp:compare:F-homo} and \eqref{eq:asymp:compare:F-beta}:

\begin{align}
\wt{F}_n(x)
& = F(x)\cdot \int_{t_1=0}^x \int_{t_2=0}^{t_1} \cdots \int_{t_n=0}^{t_{n-1}} \prod_{j=1}^n F(t_{j})^{-1}\,\beta(t_j)\,F(t_j)\,\mathrm{d}^nt \notag \\
& = \int_0^x \int_0^{t_1} \cdots \int_0^{t_{n-1}} \!\!\!\!\!\!\!\! F(x-t_1) \, \beta(t_1) \, F(t_1-t_2) \, \beta(t_2) \dotsc \beta(t_n) \, F(t_n) \mathrm{d}^nt \notag \\
& = \sum_{\eps \in \{\pm 1\}^n} \int_0^x \int_0^{t_1} \cdots \int_0^{t_{n-1}} \!\!\!\!\!\!\!\! F(x-t_1) \, \beta_{\eps_1}(t_1) \, F(t_1-t_2) \, \beta_{\eps_2}(t_2) \dotsc \beta_{\eps_n}(t_n) \, F(t_n) \mathrm{d}^nt \notag \\
& = \sum_{\eps \in \{\pm 1\}^n} \int_0^x \int_0^{t_1} \cdots \int_0^{t_{n-1}} \!\!\!\!\!\!\!\! F(\xi_\eps(t))\,\beta_{\eps_1}(t_1)\,\beta_{\eps_2}(t_2)\dotsc \beta_{\eps_n}(t_n)\,\mathrm{d}^n t \;, 
\end{align}
where we define for any \,$\eps=(\eps_1,\dotsc,\eps_n)\in\{\pm 1\}^n$\, and \,$t\in \R^n$\,
\begin{align*}
\xi_\eps(t) & := (x-t_1) + \eps_1\,(t_1-t_2) + \eps_1\,\eps_2\,(t_2-t_3) + \dotsc + \eps_1\cdots\eps_{n-1}\,(t_{n-1}-t_n) + \eps_1\cdots\eps_n\,t_n \\
& = x-2t_{\nu_1} + 2t_{\nu_2} -+\dotsc + 2(-1)^jt_{\nu_j} 
\end{align*}
with the \,$\nu_j$\, defined as in the statement of part (2)(b) of the lemma. 

For (c), we abbreviate \,$\beta^* := \max\{|\beta^{[1]}_\pm|, |\beta^{[2]}_\pm|\} \in L^p([0,1])$\, and \,$\delta := \max\{|\beta^{[1]}_+-\beta^{[2]}_+|,
|\beta^{[1]}_--\beta^{[2]}_-|\} \in L^p([0,1])$\,. Then we 
have for any \,$n\geq 1$\,, any \,$\eps=(\eps_1,\dotsc,\eps_n)\in \{\pm 1\}^n$\,,
and any \,$0\leq t_n \leq \dotsc \leq t_1 \leq x$\,
\begin{align*}
& \beta^{[1]}_{\eps_1}(t_1)\dotsc\beta^{[1]}_{\eps_n}(t_n) - \beta^{[2]}_{\eps_1}(t_1)\dotsc\beta^{[2]}_{\eps_n}(t_n) \\
= & \sum_{j=1}^n \beta^{[2]}_{\eps_1}(t_1)\cdots \beta^{[2]}_{\eps_{j-1}}(t_{j-1}) \cdot (\beta^{[1]}_{\eps_j}(t_j)-\beta^{[2]}_{\eps_j}(t_j))
\cdot \beta^{[1]}_{\eps_{j+1}}(t_{j+1}) \cdots \beta^{[1]}_{\eps_n}(t_n)
\end{align*}
and therefore by (b)
\begin{align*}
& |\wt{F}_n^{[1]}(x)-\wt{F}_n^{[2]}(x)| \\
\leq\; & \sum_{\eps \in \{\pm 1\}^n} \int_{t_1=0}^x \int_{t_2=0}^{t_1} \cdots \int_{t_n=0}^{t_{n-1}} \underbrace{|F(\xi_\eps(t))|}_{\leq |F(x)|} \cdot |\beta^{[1]}_{\eps_1}(t_1) \, \cdots \, \beta^{[1]}_{\eps_n}(t_n)-\beta^{[2]}_{\eps_1}(t_1) \, \cdots \, \beta^{[2]}_{\eps_n}(t_n)|\,\mathrm{d}^nt \\
\leq\; & |F(x)|\cdot \sum_{\eps \in \{\pm 1\}^n} \sum_{j=1}^n \int_{t_1=0}^x \int_{t_2=0}^{t_1} \cdots \int_{t_n=0}^{t_{n-1}} 
|\beta^{[2]}_{\eps_1}(t_1)|\cdots |\beta^{[2]}_{\eps_{j-1}}(t_{j-1})| \cdot |\beta^{[1]}_{\eps_j}(t_j)-\beta^{[2]}_{\eps_j}(t_j)| \\
& \hspace{8cm} \cdot |\beta^{[1]}_{\eps_{j+1}}(t_{j+1})| \cdots |\beta^{[1]}_{\eps_n}(t_n)| \, \mathrm{d}^nt \\
\leq\; & |F(x)|\cdot \sum_{\eps \in \{\pm 1\}^n} \sum_{j=1}^n \int_{t_1=0}^x \int_{t_2=0}^{t_1} \cdots \int_{t_n=0}^{t_{n-1}} 
\beta^*(t_1)\cdots \beta^*(t_{j-1}) \cdot \delta(t_j) 
\cdot \beta^*(t_{j+1}) \cdots \beta^{*}(t_n) \, \mathrm{d}^nt \\
\leq\; & |F(x)| \cdot 2^n \cdot n \cdot \frac{1}{(n-1)!} \cdot \|\beta^*\|_{L^1}^{n-1} \cdot \|\delta\|_{L^1} \;,
\end{align*}
whence
\begin{align*}
|\wt{F}^{[1]}(x) - \wt{F}^{[2]}(x)| 
& \leq \sum_{n=1}^\infty |\wt{F}_n^{[1]}(x) - \wt{F}_n^{[2]}(x)| \leq |F(x)| \cdot \|\delta\|_{L^1} \sum_{n=1}^\infty \frac{n}{(n-1)!}\cdot 2^n \cdot \|\beta^*\|_{L^1}^{n-1} \\
& \leq 8\,|F(x)|\cdot \|\delta\|_{L^1} \cdot \|\beta^*\|_{L^1} \cdot (1+\exp(2\,\|\beta^*\|_{L^1})) = C\cdot |F(x)| \cdot \|\delta\|_{L^1}
\end{align*}
follows with \,$C := 8\|\beta^*\|_{L^1} \cdot (1+\exp(2\|\beta^*\|_{L^1}))$\,. 
\end{proof}

\begin{proof}[Proof of Theorem~\ref{T:asymp:basic}.]
For the purposes of the proof, we fix an arbitrary branch of the square root function on some slitted complex plane, denoted by \,$\lambda^{1/2}$\, (or \,$\sqrt{\lambda}$\,).
Below, we will show the claimed asymptotic assessment for \,$\lambda$\, in this slitted plane; it then extends to all of \,$\C^*$\, by an argument of continuity.

We will primarily prove the asymptotic behavior for \,$\lambda\to\infty$\,; only in the end
will we use Proposition~\ref{P:mono:symmetry} to transfer that result to the case \,$\lambda\to 0$\,. 

The first important step in the proof is to \emph{regauge} \,$\alpha$\,,
that is to say, we choose a function \,$g=g_\lambda: \R\to \GL(2,\C)$\, depending holomorphically on the spectral parameter 
\,$\lambda\in\C^*$\, such that \,$g_\lambda\in W^{1,2}([0,1],\C^{2\times 2})$\, (\,$g$\, should be periodic if \,$u$\, is periodic), and then pass  from \,$\alpha$\, and \,$F$\, to
\begin{equation}
\label{eq:asymp:regauge}
\widetilde{\alpha} = g^{-1}\,\alpha\,g - g^{-1}\,g' \qmq{and} \widetilde{F}(x) = g(x)^{-1} \cdot F(x) \cdot g(0) \; .
\end{equation}
We then again have \,$\wt{\alpha} \in L^2([0,1],\C^{2\times 2})$\, and \,$\wt{F} \in W^{1,2}([0,1],\C^{2\times 2})$\,, and \,$\wt{F}$\, satisfies the following ordinary initial value problem, 
analogous to Equation~\eqref{eq:asymp:Fode}
\begin{equation}
\label{eq:asymp:wtFode}
\wt{F}'(x) = \wt{\alpha}(x)\cdot \wt{F}(x) \qmq{and} \wt{F}(0)=\unity  \; .
\end{equation}
We also put \,$\wt{M}(\lambda) = \wt{F}_\lambda(1)$\,.

The \,$g$\, we are going to use for the regauging of \,$\alpha$\, is
\begin{equation}
\label{eq:asymp:expand:g}
g := \begin{pmatrix} 1 & 0 \\ 0 & \lambda^{1/2} \end{pmatrix} \cdot \begin{pmatrix} e^{u/4} & 0 \\ 0 & e^{-u/4} \end{pmatrix}
= \begin{pmatrix} e^{u/4} & 0 \\ 0 & \lambda^{1/2}\,e^{-u/4} \end{pmatrix} \; . 
\end{equation}
In the definition of \,$g$\,, the first factor represents the transition from an ``untwisted'' to a ``twisted''
representation of the minimal surface; it is also needed to ensure that the leading term (with respect to \,$\lambda$\,)
of  \,$\wt{\alpha}$\, is no longer nilpotent (which would imply that its adjunct operator is not invertible, and would
cause problems later on). The second factor serves to make the leading term of \,$\wt{\alpha}$\,
independent of \,$u$\,, thereby making asymptotic assessments feasible.

Via Equation~\eqref{eq:mono:alphaxy} and the equations
\begin{equation}
\label{eq:asymp:expand:g-1}
g^{-1} = \begin{pmatrix} e^{-u/4} & 0 \\ 0 & \lambda^{-1/2}\,e^{u/4} \end{pmatrix} \qmq{and}
g' = \frac14\,\begin{pmatrix} u_x\,e^{u/4} & 0 \\ 0 & -\lambda^{1/2}\,u_x\,e^{-u/4} \end{pmatrix} \;,
\end{equation}
we obtain 
\begin{align}
\wt{\alpha} & = g^{-1}\,\alpha\,g - g^{-1}\,g' 
= \begin{pmatrix} - \tfrac12\,u_z & -\tfrac{1}{4}\,\lambda^{1/2}- \tfrac14\,\lambda^{-1/2}\,e^u \\ \tfrac{1}{4}\,\lambda^{1/2} + \frac14\,\lambda^{-1/2}\,e^{-u} & \frac12\,u_z\end{pmatrix} \notag\\
\label{eq:asymp:basic:wtalpha}
& = \wt{\alpha}_0 + \beta + \gamma 
\end{align}
with
$$
\wt{\alpha}_0 := \zeta(\lambda)\,\begin{pmatrix} 0 & -1 \\ 1 & 0 \end{pmatrix} ,\;
\beta :=  - \frac12\,u_z\,\begin{pmatrix} 1 & 0 \\ 0 & -1 \end{pmatrix} ,\; 
\gamma := \frac14\,\lambda^{-1/2}\,\begin{pmatrix} 0 & -e^u+1 \\ e^{-u}-1 & 0 \end{pmatrix} \; .
$$

We will now asymptotically compare the monodromy \,$\wt{M}(\lambda)$\, with the monodromy \,$\wt{M}_0(\lambda)
:= \wt{F}_{0,\lambda}(1)$\, corresponding to \,$\wt{F}_0' = \wt{\alpha}_0\,\wt{F}_0$\,, \,$\wt{F}_0(0) = \unity$\,. 
More specifically, if we denote by \,$|\dotsc|$\, also the maximum absolute row sum norm%
\footnote{Note that the maximum absolute row sum norm is the operator norm associated to the maximum norm on \,$\C^2$\,.
Therefore it is sub-multiplicative.}
for \,$(2\times 2)$-matrices,
we will show that for given \,$\eps>0$\, there exists \,$R>0$\, so that for \,$\lambda\in \C$\, with \,$|\lambda|\geq R$\,,
we have
\begin{equation}
\label{eq:asymp:claim}
|\wt{M}(\lambda)-\wt{M}_0(\lambda)| \leq \eps \cdot w(\lambda) \; .
\end{equation}
Note that
\begin{equation}
\label{eq:asymp:wtM0}
\wt{M}_0(\lambda) = \begin{pmatrix} \cos(\zeta(\lambda)) & -\sin(\zeta(\lambda)) \\ \sin(\zeta(\lambda)) & \cos(\zeta(\lambda)) \end{pmatrix}
\end{equation}
and therefore \,$w(\lambda)=|\wt{M}_0(\lambda)|$\, holds.

To see that the claimed estimate in Theorem~\ref{T:asymp:basic}(1) follows from~\eqref{eq:asymp:claim}, we write the non-regauged monodromy as
\,$M(\lambda) = \left( \begin{smallmatrix} a(\lambda) & b(\lambda) \\ c(\lambda) & d(\lambda) \end{smallmatrix} \right)$\, as in the statement of the theorem.
Then we have
\begin{align*}
\wt{M}(\lambda) & = g(1)^{-1} \, M(\lambda)\,g(0) \\ 
& = \begin{pmatrix} e^{-(u(1)-u(0))/4}\,a(\lambda) & \lambda^{1/2}\,e^{-(u(0)+u(1))/4}\,b(\lambda) \\ \lambda^{-1/2}\,e^{(u(0)+u(1))/4}\,c(\lambda) & e^{(u(1)-u(0))/4}\,d(\lambda) \end{pmatrix} \;.
\end{align*}
By plugging this equation, as well as Equation~\eqref{eq:asymp:wtM0} into the estimate~\eqref{eq:asymp:claim}, we obtain Theorem~\ref{T:asymp:basic}(1).

It thus suffices to prove the estimate \eqref{eq:asymp:claim}.

The expression \,$\gamma$\, can be neglected in our asymptotic consideration. The reason for this is, roughly speaking, that \,$\gamma$\, is of order \,$|\lambda|^{-1/2}$\, for \,$|\lambda|\to\infty$\,.
To prove more precisely that \,$\gamma$\, can be neglected, 
we apply Lemma~\ref{L:asymp:compare}(2)(c) with \,$\alpha=\alpha_0$\,, \,$\beta^{[1]} = \beta$\,
and \,$\beta^{[2]} = \beta+\gamma$\, to compare the resulting solutions \,$\wt{F}^{[1]}$\, and \,$\wt{F}^{[2]}$\,. Notice that
\,$\beta$\, anti-commutes with \,$\alpha_0$\,, thus we have \,$\beta^{[1]}_+=0$\, and \,$\beta^{[1]}_-=\beta$\,. Moreover, we
have
$$ \gamma_+ = \frac{1}{4}\,\lambda^{-1/2}\,(\cosh(u)-1)\,\begin{pmatrix} 0 & -1 \\ 1 & 0 \end{pmatrix}
\qmq{and} \gamma_- = -\frac{1}{4}\,\lambda^{-1/2}\,\sinh(u)\,\begin{pmatrix} 0 & 1 \\ 1 & 0 \end{pmatrix} \; . $$
It follows that we have \,$\beta_\pm^{[1]}-\beta_\pm^{[2]} = \gamma_\pm = O(|\lambda|^{-1/2})$\,, and
that \,$\|\max\{|\beta_\pm^{[\nu]}|\}\|_{L^1}$\, is bounded for \,$\lambda\to \infty$\,. It now follows from
Lemma~\ref{L:asymp:compare}(2)(c) that there exists \,$C>0$\, such that
\,$|\wt{F}^{[1]}-\wt{F}^{[2]}| \leq C \cdot |\lambda|^{-1/2}\cdot w(\lambda)$\, holds. By choosing \,$|\lambda|$\, large,
we can therefore ensure \,$|\wt{F}^{[1]}-\wt{F}^{[2]}| \leq \tfrac12 \, \eps \, w(\lambda)$\,. 

In other words, it follows from this argument that it suffices to show that for given \,$\eps>0$\, there exists \,$R>0$\, so
that for \,$\lambda\in \C$\, with \,$|\lambda|>R$\, and \,$x\in [0,1]$\, we have
\begin{equation}
\label{eq:asymp:claim2}
|E_\lambda(1)-E_{0,\lambda}(1)| \leq \eps \cdot w(\lambda) 
\end{equation}
for the solution \,$E=E_\lambda$\, of the initial value problem
\begin{equation}
\label{eq:asymp:Edgl}
E'(x) = (\wt{\alpha}_0 + \beta(x))\,E(x) \qmq{with} E(0)=\unity 
\end{equation}
and for the solution \,$E_0(x)=E_{0,\lambda}(x)=\exp(x\,\wt{\alpha}_0)$\, of
\begin{equation*}
E_0'(x) = \wt{\alpha}_0\,E_0(x) \qmq{with} E_0(0)=\unity \; ;
\end{equation*}
we have
\begin{equation}
\label{eq:asymp:E0}
E_0(x) = \wt{F}_0(x) = \begin{pmatrix} \cos(\zeta(\lambda)\,x) & -\sin(\zeta(\lambda)\,x) \\ \sin(\zeta(\lambda)\,x) & \cos(\zeta(\lambda)\,x) \end{pmatrix} \;.
\end{equation}



Because \,$\beta$\, anti-commutes with \,$\wt{\alpha}_0$\,, we have by Lemma~\ref{L:asymp:compare}(2)(b) 
\begin{equation}
\label{eq:asymp:E-power-series}
E = E_0 + \sum_{n=1}^\infty E_n 
\end{equation}
with
\begin{align}
E_n(x) & = \int_0^x \int_0^{t_1} \cdots \int_0^{t_{n-1}} \!\!\!\!\!\!\!\! E_0(x-t_1)\,\beta(t_1) \,E_0(t_1-t_2)\,\beta(t_2) E_0(t_2-t_3) \cdots \beta(t_n)\,E_0(t_n)\,\mathrm{d}^n t  \notag \\
\label{eq:asymp:En-2}
& = \int_0^x \int_0^{t_1} \cdots \int_0^{t_{n-1}} \!\!\!\!\!\!\!\! E_0(\xi(t))\,\beta(t_1) \,\beta(t_2) \cdots \beta(t_n) \,\mathrm{d}^n t  \;,
\end{align}
where
$$ \xi(t) := x-2t_1+2t_2-2t_3+-\dotsc+2(-1)^n\,t_n \; . $$


In the integral of Equation~\eqref{eq:asymp:En-2}, we now carry out the substitution \,$(t_1,\dotsc,t_n) \mapsto (s_1,\dotsc,s_n)$\, with \,$s_j = t_j-s_{j+1}$\, for \,$j\leq n-1$\, and \,$s_n = t_n$\,, i.e.~
\begin{align*}
s_n & = t_n \\
s_{n-1} & = t_{n-1}-s_n = t_{n-1}-t_n \\
s_{n-2} & = t_{n-2}-s_{n-1} = t_{n-2}-t_{n-1}+t_n \\
s_{n-3} & = t_{n-3}-s_{n-2} = t_{n-3}-t_{n-2}+t_{n-1}-t_n \\
& \dotsc \\
s_2 & = t_2-s_3 = t_2-t_3+t_4-+ \dotsc + (-1)^n\,t_n \\
s_1 & = t_1-s_2 = t_1-t_2+t_3-+ \dotsc + (-1)^{n+1}\,t_n \; . 
\end{align*}
Then we have \,$t_j = s_j+s_{j+1}$\, for \,$j\leq n-1$\, and \,$t_n=s_n$\,. 
Thus, the corresponding mapping \,$\Phi: (s_j) \mapsto (t_j)$\, is a diffeomorphism with \,$\det \Phi'=1$\, from 
$$ U := \Menge{(s_1,\dotsc,s_n)\in \R^n}{\forall j=1,\dotsc,n: s_j\geq 0, \; s_1+s_2 \leq x,\; \forall j=3,\dotsc,n: s_j \leq s_{j-2}} $$
onto the simplex
$$ O := \Menge{(t_1,\dotsc,t_n)\in \R^n}{0 \leq t_n \leq t_{n-1} \leq \cdots \leq t_1 \leq x} $$
which is the domain of integration in \eqref{eq:asymp:En-2}.
We also see
$$ \xi(t) = x-2s_1 \; . $$
By carrying out the substitution in \eqref{eq:asymp:En-2}, we thus obtain
\begin{align}
E_n(x) & = \int_U E_0(x-2s_1)\,\beta(s_1+s_2)\,\beta(s_2+s_3)\,\dotsc\,\beta(s_{n-1}+s_n)\,\beta(s_n)\,\mathrm{d}^ns \notag \\
& = \int_{s_1=0}^x E_0(x-2s_1)\int_{s_2=0}^{x-s_1} \beta(s_1+s_2) \int_{s_3=0}^{s_1} \beta(s_2+s_3) \int_{s_4=0}^{s_2} \beta(s_3+s_4) \notag\\
& \qquad\qquad \cdots \int_{s_n=0}^{s_{n-2}} \beta(s_{n-1}+s_n)\,\beta(s_n)\,\mathrm{d}^ns \notag \\
\label{eq:asymp:EnGn}
& = \int_{s_1=0}^x E_0(x-2s_1)\,G_n(s_1)\,\mathrm{d}s_1
\end{align}
with 
\begin{equation}
\label{eq:asymp:Gn-def}
G_n(s_1) := H_n(x-s_1,s_1) \;,
\end{equation}
where \,$H_n$\, is defined by
\begin{equation}
\label{eq:asymp:H1-def}
H_1(s_0,s_1) = \beta(s_1)
\end{equation}
and 
$$ H_n(s_0,s_1) = \int_{s_2=0}^{s_0} \beta(s_1+s_2)\,H_{n-1}(s_1,s_2)\,\mathrm{d}s_2 \qmq{for \,$n\geq 2$\,.} $$
We thus have for \,$n\geq 2$\,
\begin{equation}
\label{eq:asymp:Hn-def}
H_n(s_0,s_1) := \int_{s_2=0}^{s_0} \beta(s_1+s_2) \int_{s_3=0}^{s_1} \beta(s_2+s_3) \cdots \int_{s_n=0}^{s_{n-2}} \beta(s_{n-1}+s_n)\,\beta(s_n)\,\mathrm{d}^{n-1}s \; .
\end{equation}

We will now show%
\footnote{The proof given below shows that we in fact have \,$G_n \in L^\infty([0,1],\C^{2\times 2})$\, for \,$n\geq 2$\,.}
for all \,$n\geq 1$\,
\begin{equation}
\label{eq:asymp:Gn}
G_n \in L^2([0,1],\C^{2\times 2}) \qmq{and} \|G_n\|_2 \leq \frac{1}{(\lfloor n/2 \rfloor-1)!}\cdot \|\beta\|_2^n \; . 
\end{equation}
Indeed, for \,$n=1$\, we have \,$G_1 = \beta$\,, and thus \eqref{eq:asymp:Gn} holds in this case. For \,$n=2$\,, we have
$$ G_2(s_1) = H_2(x-s_1,s_1) = \int_{0}^{x-s_1} \beta(s_1+s_2)\,\beta(s_2)\,\mathrm{d}s_2 $$
and therefore \,$G_2 \in W^{1,1}([0,1],\C^{2\times 2}) \subset L^2([0,1],\C^{2\times 2})$\, and \,$\|G_2\|_2 \leq \|\beta\|_2^2$\, by the Cauchy-Schwarz inequality. For even \,$n\geq 4$\,,
we fix \,$s_0,s_1\in [0,1]$\, and define the simplices
\begin{align*}
S_+ & := \Menge{(s_2,s_4,\dotsc,s_n)\in\R^{n/2}}{0 \leq s_n \leq \dotsc \leq s_4 \leq s_2 \leq s_0} \;, \\
S_- & := \Menge{(s_3,s_5,\dotsc,s_{n-1})\in \R^{n/2-1}}{0 \leq s_{n-1} \leq \dotsc \leq s_5 \leq s_3 \leq s_1} \; . 
\end{align*}
Then we have
\begin{align*}
& H_n(s_0,s_1) \\
= &  \int_{S_-} \int_{S_+} \!\!\!\! \beta(s_2+s_1)\,\beta(s_2+s_3)\,\beta(s_4+s_3)\,\beta(s_4+s_5)\,\dotsc\,\beta(s_n+s_{n-1})\,\beta(s_n)\,\mathrm{d}^{n/2}s\,\mathrm{d}^{n/2-1}s \; .
\end{align*}
For the inner integral, we have by the Cauchy-Schwarz inequality
$$ \left | \int_{S_+} \beta(s_2+s_1)\,\beta(s_2+s_3)\,\beta(s_4+s_3)\,\beta(s_4+s_5)\,\dotsc\,\beta(s_n+s_{n-1})\,\beta(s_n)\,\mathrm{d}^{n/2}s \right| \leq \|\beta\|_2^n \;, $$
and thus we have
$$ |H_n(s_0,s_1)| \leq \vol(S_-) \cdot \|\beta\|_2^n = \frac{s_1^{n/2-1}}{(n/2-1)!}\cdot \|\beta\|_2^n \leq \frac{1}{(n/2-1)!}\cdot \|\beta\|_2^n \; . $$
Thus we see \,$G_n(s) = H_n(x-s,s) \in L^{\infty}([0,1],\C^{2\times 2}) \subset L^2([0,1],\C^{2\times 2})$\, and \,$\|G_n\|_2 \leq \|H_n\|_\infty \leq \tfrac{1}{(n/2-1)!}\cdot \|\beta\|_2^n$\,. 
Therefore we obtain \eqref{eq:asymp:Gn} for even \,$n\geq 4$\,. 
For odd \,$n\geq 3$\,, we can argue similarly, obtaining
$$ |H_n(s_0,s_1)| \leq \frac{1}{((n-1)/2)!}\cdot \|\beta\|_2^n \; , $$
and therefore also in this case \eqref{eq:asymp:Gn}, completing its proof.

As a consequence of \eqref{eq:asymp:Gn}, we have
$$ \sum_{n=1}^\infty \|G_n\|_2 \leq \sum_{n=1}^\infty \frac{1}{(\lfloor n/2 \rfloor-1)!}\cdot \|\beta\|_2^n < \infty $$
(because the power series \,$\sum_{n=1}^\infty \frac{1}{(\lfloor n/2 \rfloor-1)!}\cdot r^n$\, is convergent for every real \,$r$\,), and thus we have
$$ G := \sum_{n=1}^\infty G_n \;\in\; L^2([0,1],\C^{2\times 2}) \; . $$
It follows via Equations~\eqref{eq:asymp:E-power-series} and \eqref{eq:asymp:EnGn} that the series \,$\sum_{n=0}^\infty E_n$\, converges in \,$W^{1,2}([0,1],\C^{2\times 2})$\, and that we have
$$ E(x)-E_0(x) = \int_0^x E_0(x-2s)\,G(s)\,\mathrm{d}s \; . $$

We now use this representation of \,$E-E_0$\, to show the estimate~\eqref{eq:asymp:claim}.
For this purpose, we set \,$x=1$\,, and fix \,$\delta>0$\, at first arbitrarily. Then we have for all \,$s\in [0,1]$\,: \,$|E_0(1-2s)| \leq |E_0(1)|=w(\lambda)$\,, and for all
\,$s\in [\delta,1-\delta]$\,: \,$|E_0(1-2s)| \leq |E_0(1-2\delta)| \leq e^{-2\delta\,|\IM(\zeta(\lambda))|}\cdot w(\lambda)$\,. Thus we obtain
\begin{align*}
|E(1)-E_0(1)| & \leq \int_0^1 |E_0(1-2s)|\cdot |G(s)|\,\mathrm{d}s \\
& = \int_0^\delta |E_0(1-2s)|\cdot |G(s)|\,\mathrm{d}s + \int_\delta^{1-\delta} |E_0(1-2s)|\cdot |G(s)|\,\mathrm{d}s \\
& \qquad\qquad + \int_{1-\delta}^1 |E_0(1-2s)|\cdot |G(s)|\,\mathrm{d}s \\
& \leq w(\lambda) \cdot \|G|[0,\delta]\|_1 + e^{-2\delta\,|\IM(\zeta(\lambda))|}\cdot w(\lambda) \cdot \|G\|_1 + w(\lambda) \cdot \|G|[1-\delta,1]\|_1 \; . 
\end{align*}
Because of \,$G\in L^1([0,1],\C^{2\times 2})$\,, we can now choose \,$\delta>0$\, such that we have \,$\|G|[0,\delta]\|_1 \leq \tfrac{\eps}{4}$\, and \,$\|G|[1-\delta,1]\|_1 \leq \tfrac{\eps}{4}$\,. 
Dependent on this \,$\delta$\,, there further exists \,$C>0$\, such that \,$e^{-2\delta\,C}\cdot \|G\|_1 \leq \tfrac{\eps}{2}$\,. For all
\,$\lambda\in\C^*$\, with \,$|\IM(\zeta(\lambda))| \geq C$\, we then obtain 
\begin{equation}
\label{eq:asymp:EE0-zetaC}
|E(1)-E_0(1)| \leq \eps \cdot w(\lambda)
\end{equation}
and therefore \eqref{eq:asymp:claim2}.

It remains to show that \eqref{eq:asymp:claim2} also holds within \,$\Mengegr{\lambda\in\C^*}{|\IM(\zeta(\lambda))|\leq C}$\, for \,$\lambda$\, of 
sufficiently large absolute value. We have
$$ E(1)-E_0(1) = \int_0^1 E_0(1-2s)\,G(s)\,\mathrm{d}s = \frac12 \int_{-1}^1 E_0(s)\,G\left( \tfrac{1-s}{2}\right) \,\mathrm{d}s \;, $$
where the entries of \,$E_0(s)$\, are 
\begin{equation}
\label{eq:asymp:cossinexp}
\cos(\zeta(\lambda)s) = \frac12(e^{is\zeta(\lambda)} + e^{-is\zeta(\lambda)}) \qmq{and} \pm\sin(\zeta(\lambda)s) = \pm\frac{1}{2i}(e^{is\zeta(\lambda)} - e^{-is\zeta(\lambda)}) \;,
\end{equation}
therefore it follows from the variant of Riemann-Lebesgue's Lemma given in the following lemma (Lemma~\ref{L:asymp:riemannlebesgue}), 
applied with the compact set \,$N\subset L^1([-1,1])$\, comprised of the four component functions of \,$G(\tfrac{1-s}{2})$\,,
that there exists \,$R>0$\, so that  \eqref{eq:asymp:claim2} holds for \,$\lambda\in \C^*$\, with \,$|\RE(\zeta(\lambda))|\geq R$\, and \,$|\IM(\zeta(\lambda))|\leq C$\,. 
This completes the proof of (1).

To prove that \,$R$\, can be chosen uniformly in dependence on \,$\eps$\, for relatively compact subsets \,$P$\, of \,$\Pot_{np}$\,, we note that for \,$(u,u_y)\in P$\,, \,$\|(u,u_y)\|_{\Pot}$\, is bounded, and therefore
\,$\|\beta\|_{L^1}$\, and \,$\|\gamma\|_{L^1}$\, are also bounded for \,$(u,u_y)\in P$\,. For this reason, the estimates leading up to \eqref{eq:asymp:EE0-zetaC} are uniform for \,$(u,u_y)\in P$\,,
and thus we obtain \eqref{eq:asymp:claim2} for \,$\lambda\in \C^*$\, with \,$|\IM(\zeta(\lambda))|\geq C$\, for all \,$(u,u_y)\in P$\,, where \,$C>0$\, where \,$C>0$\, is a constant depending on \,$\eps$\,.
Moreover, we let \,$N$\, be the topological closure of the set of the component functions of \,$G(\tfrac{1-s}{2})$\, where \,$(u,u_y)$\, now runs through all of \,$P$\,; \,$N$\, is compact, because \,$P$\,
is relatively compact. By applying Lemma~\ref{L:asymp:riemannlebesgue} with this \,$N$\,, we obtain \eqref{eq:asymp:claim2} for \,$|\RE(\zeta(\lambda))|\geq R$\, and \,$|\IM(\zeta(\lambda))|\leq C$\,
for all \,$(u,u_y)\in P$\,. 

To prove part (2) of the theorem (i.e.~the case \,$\lambda\to 0$\,), we let
\,$\wt{u}(x) := -u(x)$\, and \,$\wt{u}_y(x) := u_y(x)$\,; 
then we have \,$(\wt{u},\wt{u}_y)\in \Pot$\,. If we denote the monodromy corresponding to \,$u$\,
resp.~to \,$\wt{u}$\, by 
$$ M_u(\lambda)= \left( \begin{smallmatrix} {a}(\lambda) & {b}(\lambda) \\ {c}(\lambda) & {d}(\lambda) \end{smallmatrix} \right) \qmq{resp.~by} M_{\wt{u}}(\lambda) = \left( \begin{smallmatrix} \wt{a}(\lambda) & \wt{b}(\lambda) \\ \wt{c}(\lambda) & \wt{d}(\lambda) \end{smallmatrix} \right) \;, $$ 
then we have by Proposition~\ref{P:mono:symmetry}(1)
$$ M_u(\lambda) = \left( \begin{smallmatrix} 1 & 0 \\ 0 & \lambda \end{smallmatrix} \right) \cdot M_{\wt{u}}(\lambda^{-1}) \cdot \left( \begin{smallmatrix} 1 & 0 \\ 0 & \lambda^{-1} \end{smallmatrix} \right) \;, $$
i.e.
\begin{equation}
\label{eq:asymp:symmetry}
\begin{pmatrix} a(\lambda) & b(\lambda) \\ c(\lambda) & d(\lambda) \end{pmatrix} = \begin{pmatrix} \wt{a}(\lambda^{-1}) & \lambda^{-1}\cdot \wt{b}(\lambda^{-1}) \\ \lambda\cdot \wt{c}(\lambda^{-1}) & \wt{d}(\lambda^{-1}) \end{pmatrix} \; .
\end{equation}
By the result for \,$\lambda\to\infty$\,, applied to \,$\wt{u}$\,, we have for \,$|\lambda^{-1}|\geq R$\,
\begin{align*}
|\wt{a}(\lambda^{-1})-e^{(\wt{u}(1)-\wt{u}(0))/4}\,\cos(\zeta(\lambda^{-1}))| & \leq \eps \, w(\lambda^{-1}) \\
|\wt{b}(\lambda^{-1})-(-\lambda^{1/2}\,e^{(\wt{u}(0)+\wt{u}(1))/4}\,\sin(\zeta(\lambda^{-1}))| & \leq \eps \, |\lambda|^{1/2}\,w(\lambda^{-1}) \\
|\wt{c}(\lambda^{-1})-\lambda^{-1/2}\,e^{-(\wt{u}(0)+\wt{u}(1))/4}\,\sin(\zeta(\lambda^{-1}))| & \leq \eps \, |\lambda|^{-1/2}\,w(\lambda^{-1}) \\
|\wt{d}(\lambda^{-1})-e^{-(\wt{u}(1)-\wt{u}(0))/4}\,\cos(\zeta(\lambda^{-1}))| & \leq \eps \, w(\lambda^{-1}) \; . 
\end{align*}
By Equation~\eqref{eq:asymp:symmetry}, and \,$\zeta(\lambda^{-1})=\zeta(\lambda)$\,, \,$w(\lambda^{-1})=w(\lambda)$\,, \,$\wt{u}(0)=-u(0)$\,, we now obtain
for \,$|\lambda|\leq 1/R$\,
\begin{align*}
|a(\lambda)-e^{-(u(1)-u(0))/4}\,\cos(\zeta(\lambda))| & \leq \eps \, w(\lambda) \\
|\lambda\cdot b(\lambda)-(-\lambda^{1/2}\,e^{-(u(0)+u(1))/4}\,\sin(\zeta(\lambda))| & \leq \eps \, |\lambda|^{1/2}\,w(\lambda) \\
|\lambda^{-1}\cdot c(\lambda)-\lambda^{-1/2}\,e^{(u(0)+u(1))/4}\,\sin(\zeta(\lambda))| & \leq \eps \, |\lambda|^{-1/2}\,w(\lambda) \\
|d(\lambda)-e^{(u(1)-u(0))/4}\,\cos(\zeta(\lambda))| & \leq \eps \, w(\lambda) 
\end{align*}
and therefore the statement claimed for \,$\lambda\to0$\,. It is clear that also the statement on the uniformness of the estimate for \,$(u,u_y)\in P$\, transfers to the case \,$\lambda\to 0$\,. 
\end{proof}

The following Lemma, which is a variant of the Riemann-Lebesgue Lemma, has been used in the preceding proof:

\begin{lem}
\label{L:asymp:riemannlebesgue}
Let \,$N$\, be a compact subset of \,$L^1([a,b])$\,. 
For any given \,$\eps,C>0$\, there then exists \,$R>0$\, such that for every \,$\zeta\in \C$\, with \,$|\IM(\zeta)|\leq C$\,, \,$|\RE(\zeta)|\geq R$\, and every \,$g\in N$\, 
we have
$$ \left| \int_a^b e^{-2\pi i\zeta t}\,g(t)\,\mathrm{d}t \right| \leq \eps \; . $$
\end{lem}

\begin{proof}
We extend the functions in \,$L^1([a,b])$\, to 
\,$\R$\, by zero, then we have \,$N \subset L^1(\R)$\,. The map
$$ \Phi: L^1(\R)\times [-C,C] \to L^1(\R),\; (g,y) \mapsto (t \mapsto e^{2\pi yt}\cdot g(t)) $$
is continuous, hence the image \,$\wt{N} := \Phi(N\times [-C,C])$\, is a compact set in \,$L^1(\R)$\,. 

Therefore there exist finitely many \,$f_1,\dotsc,f_n\in L^1(\R)$\, so that
\,$\wt{N} \subset \bigcup_{k=1}^n B(f_k,\tfrac{\eps}{2})_{L^1}$\,. By the classical 
Riemann-Lebesgue Lemma (\cite{Grafakos:2008}, Proposition~2.2.17,
p.~105), there exists \,$R>0$\,, so that we have \,$|\widehat{f}_k(x)|\leq\tfrac{\eps}{2}$\, for every
\,$x\in \R$\, with \,$|x|\geq R$\, and every \,$k\in \{1,\dotsc,n\}$\,; here we denote for any 
\,$f\in L^1(\R)$\, by 
$$ \widehat{f}(x) := \int_{\R} e^{-2\pi ixt}\,f(t)\,\mathrm{d}t $$
the Fourier transform of \,$f$\,. 



Now let \,$g\in N$\, and \,$\zeta=x+iy \in \C$\, be given with \,$|x| \geq R$\,, \,$|y| \leq C$\,. By construction,
there exists some \,$k\in \{1,\dotsc,n\}$\, with \,$\|\Phi(g,y)-f_k\|_1 \leq \tfrac{\eps}{2}$\,,
and with this \,$k$\, we have
$$ \|\widehat{\Phi(g,y)} - \widehat{f}_k\|_\infty \leq \|\Phi(g,y)-f_k\|_1 \leq \frac{\eps}{2} $$
by the Hausdorff-Young inequality for the case \,$p'=\infty$\,, \,$p=1$\, (\cite{Grafakos:2008}, Proposition~2.2.16,
p.~104). We now have
\begin{gather*}
\left| \int_a^b e^{-2\pi i\zeta t}\,g(t)\,\mathrm{d}t \right| 
= \left| \int_a^b e^{-2\pi ix t}\,e^{2\pi y t}\,g(t)\,\mathrm{d}t \right| 
= \left| \int_a^b e^{-2\pi ix t}\,\Phi(g,y)(t)\,\mathrm{d}t \right| = \left| \widehat{\Phi(g,y)}(x) \right| \\
\leq \left|\widehat{\Phi(g,y)}(x) - \widehat{f}_k(x)\right| +  \left| \widehat{f}_k(x) \right| 
\leq \frac{\eps}{2} + \frac{\eps}{2} = \eps \; . 
\end{gather*}
\end{proof}

\begin{cor}
\label{C:asymp:basic-deltabc}
In the setting of Theorem~\ref{T:asymp:basic} we have the following additional facts:
\begin{enumerate}
\item There exists \,$C>0$\, (dependent on \,$(u,u_y)$\,) so that we have
\begin{align*}
|a(\lambda)| & \leq C\, w(\lambda) \\
|b(\lambda)| & \leq C\,|\lambda|^{-1/2}\,w(\lambda) \\
|c(\lambda)| & \leq C\,|\lambda|^{1/2}\,w(\lambda) \\
|d(\lambda)| & \leq C\, w(\lambda) \;.
\end{align*}
\item
There exists for every \,$\eps>0$\, some \,$R>0$\,, so that for every
\,$\lambda\in\C^*$\, with \,$|\lambda|\geq R$\, or \,$|\lambda| \leq \tfrac{1}{R}$\,, we have
\begin{enumerate}
\item \,$|b(\lambda)\,c(\lambda)-b_0(\lambda)\,c_0(\lambda)| \leq \eps\,w(\lambda)^2$\,.
\end{enumerate}
We now suppose \,$(u,u_y)\in \Pot$\,. Then we define \,$\Delta(\lambda):=\tr(M(\lambda))=a(\lambda)+d(\lambda)$\,,
and for every \,$\eps>0$\, there exists some \,$R>0$\,, so that for every
\,$\lambda\in\C^*$\, with \,$|\lambda|\geq R$\, or \,$|\lambda| \leq \tfrac{1}{R}$\,, we have
\begin{enumerate}
\addtocounter{enumii}{1}
\item \,$|\Delta(\lambda)-\Delta_0(\lambda)| \leq \eps\,w(\lambda)$\,. 
\item \,$|(\Delta(\lambda)^2-4)-(\Delta_0(\lambda)^2-4)| \leq \eps\,w(\lambda)^2$\,.
\end{enumerate}
\end{enumerate}
The constants \,$C$\, resp.~\,$R$\, can be chosen uniformly for \,$(u,u_y) \in P$\,, where \,$P$\, is a relatively 
compact subset of \,$\Pot_{np}$\,. 
\end{cor}

\begin{proof}
\emph{For (1).}
By Theorem~\ref{T:asymp:basic} there exists \,$R>0$\, (corresponding to \,$\eps=1$\,) so that we have for \,$\lambda\in \C^*$\, with \,$|\lambda|\geq R$\, or \,$|\lambda|\leq \tfrac{1}{R}$\,
$$ |a(\lambda)-\upsilon\,a_0(\lambda)| \leq w(\lambda) $$
and therefore 
$$ |a(\lambda)| \leq |\upsilon|\cdot |a_0(\lambda)| + |a(\lambda)-\upsilon\,a_0(\lambda)| \leq |\upsilon|\cdot w(\lambda)+w(\lambda) = (|\upsilon|+1)\cdot w(\lambda) \; . $$
Because \,$\{\tfrac{1}{R}\leq |\lambda|\leq R\}$\, is compact, \,$a$\, is bounded on this set by some constant \,$C_1>0$\,. Because \,$w(\lambda)\geq 1$\, holds for all \,$\lambda\in\C^*$\,, we then have
\,$|a(\lambda)| \leq C\,w(\lambda)$\, for all \,$\lambda\in\C^*$\, with \,$C := \max\{|\upsilon|+1,C_1\}$\,. The claims on \,$b$\,, \,$c$\, and \,$d$\, are shown similarly. 

\emph{For (2).} For any \,$\wt{\eps}>0$\, we have by Theorem~\ref{T:asymp:basic} some \,$R>0$\, so that the estimates from that theorem hold for \,$|\lambda|\geq R$\, resp. for \,$|\lambda|\leq\tfrac{1}{R}$\,. 

We then have for \,$|\lambda|\geq R$\,
\begin{align*}
|b(\lambda)\,c(\lambda) - b_0(\lambda)\,c_0(\lambda)|
& = |b(\lambda)\,c(\lambda) - \tau^{-1}\,b_0(\lambda)\,\tau\,c_0(\lambda)| \\
& \leq |b(\lambda)-\tau^{-1}\,b_0(\lambda)|\cdot|c(\lambda)| + |\tau|^{-1}\,|b_0(\lambda)| \cdot |c(\lambda)-\tau\,c_0(\lambda)| \\
& \leq \wt{\eps}\,|\lambda|^{-1/2}\,w(\lambda)\cdot C\,|\lambda|^{1/2}\,w(\lambda) + |\tau|^{-1}\,|\lambda|^{-1/2}\,w(\lambda)\cdot \wt{\eps}\,|\lambda|^{1/2}\,w(\lambda) \\
& = (C+|\tau|^{-1})\,\wt{\eps}\,w(\lambda)^2 \; . 
\end{align*}
A similar calculation gives for \,$|\lambda|\leq \tfrac{1}{R}$\, 
$$ |b(\lambda)\,c(\lambda) - b_0(\lambda)\,c_0(\lambda)| \leq (C+|\tau|)\,\wt{\eps}\,w(\lambda)^2 \; . $$
By choosing \,$\wt{\eps} := (C+\max\{|\tau|,|\tau|^{-1}\})^{-1}\cdot \eps$\,, we obtain (2)(a).

We now suppose that we have \,$(u,u_y)\in \Pot$\, and therefore \,$\upsilon=1$\,. 
We then have
$$ |\Delta(\lambda)-\Delta_0(\lambda)| \leq |a(\lambda)-a_0(\lambda)| + |d(\lambda)-d_0(\lambda)| \leq 2\,\wt{\eps}\,w(\lambda) \;, $$
therefore we obtain (2)(b) by choosing \,$\wt{\eps} := \eps/2$\,. Finally, we have
\begin{align*}
|(\Delta(\lambda)^2-4)-(\Delta_0(\lambda)^2-4)| 
& = |\Delta(\lambda)+\Delta_0(\lambda)| \cdot |\Delta(\lambda)-\Delta_0(\lambda)| \\
& \leq (2|\Delta_0(\lambda)| + |\Delta(\lambda)-\Delta_0(\lambda)|) \cdot |\Delta(\lambda)-\Delta_0(\lambda)| \\
& \overset{(b)}{\leq} (4\,w(\lambda)+\wt{\eps}\,w(\lambda))\cdot \wt{\eps}\,w(\lambda) = (4+\wt{\eps})\,\wt{\eps}\,w(\lambda)^2 \; . 
\end{align*}
By choosing \,$\wt{\eps}>0$\, such that \,$(4+\wt{\eps})\, \wt{\eps} = \eps$\,, we obtain (2)(c).
\end{proof}

The following corollary shows that the derivatives (with respect to \,$\lambda$\,) of the functions comprising the monodromy
can also be estimated via Theorem~\ref{T:asymp:basic}.
This is basically a consequence of the fact that these functions depend holomorphically on \,$\lambda$\,, whence it follows that 
Cauchy's inequality is applicable.
For reference, we note
\begin{align}
\label{eq:asymp:a0'}
a_0'(\lambda) = d_0'(\lambda) & = \frac{1-\lambda}{8\,\lambda^2}\,c_0(\lambda) \\ 
\label{eq:asymp:b0'}
b_0'(\lambda) & = \frac{1-\lambda}{8\,\lambda}\,a_0(\lambda)-\frac{1}{2\,\lambda}\,b_0(\lambda) \\
\label{eq:asymp:c0'}
c_0'(\lambda) & = \frac{\lambda-1}{8\,\lambda}\,a_0(\lambda) + \frac{1}{2\,\lambda}\,c_0(\lambda) \; . 
\end{align}

\begin{cor}
\label{C:asymp:Mprime}
In the situation of Theorem~\ref{T:asymp:basic}, for every \,$\eps>0$\, there exists \,$R'>0$\, such that:
\begin{enumerate}
\item
For all \,$\lambda\in\C$\, with \,$|\lambda|\geq R'$\,, we have
\begin{align*}
|a'(\lambda)-\upsilon\,a_0'(\lambda)| & \leq \eps \, |\lambda|^{-1/2}\,w(\lambda) \\
|b'(\lambda)-\tau^{-1}\,b_0'(\lambda)| & \leq \eps \, |\lambda|^{-1}\,w(\lambda) \\
|c'(\lambda)-\tau\,c_0'(\lambda)| & \leq \eps \,w(\lambda) \\
|d'(\lambda)-\upsilon^{-1}\,d_0'(\lambda)| & \leq \eps \, |\lambda|^{-1/2}\,w(\lambda) \; . 
\end{align*}
\item
For all \,$\lambda\in\C^*$\, with \,$|\lambda|\leq \tfrac{1}{R'}$\,, we have
\begin{align*}
|a'(\lambda)-\upsilon^{-1}\,a_0'(\lambda)| & \leq \eps \,|\lambda|^{-3/2}\, w(\lambda) \\
|b'(\lambda)-\tau\,b_0'(\lambda)| & \leq \eps \, |\lambda|^{-2}\,w(\lambda) \\
|c'(\lambda)-\tau^{-1}\,c_0'(\lambda)| & \leq \eps \, |\lambda|^{-1}\,w(\lambda) \\
|d'(\lambda)-\upsilon\,d_0'(\lambda)| & \leq \eps \, |\lambda|^{-3/2}\,w(\lambda) \; . 
\end{align*}
\end{enumerate}
The constant \,$R'$\, can be chosen uniformly for \,$(u,u_y) \in P$\,, where \,$P$\, is a relatively 
compact subset of \,$\Pot_{np}$\,. 
\end{cor}

\begin{proof}
We show the estimates for the function \,$a'$\,, the other functions are handled analogously.
Let \,$\eps>0$\, be given, and 
fix \,$\delta>0$\,. By Theorem~\ref{T:asymp:basic} there exists \,$R>0$\, such that for \,$\lambda\in \C$\, with \,$|\lambda|\geq R$\, we have
$$ |a(\lambda)-\upsilon\,a_0(\lambda)| \leq \frac{\eps\,\delta}{2\,e^\delta}\,w(\lambda) \; . $$
There exists \,$R'>R$\, such that for all \,$\lambda\in \C$\, with \,$|\lambda|\geq R'$\,, the closure \,$\overline{U}$\, of 
\,$U := \Mengegr{\lambda'\in \C}{|\lambda'-\lambda|<\delta\,|\lambda|^{1/2}}$\, is entirely contained in \,$\{|\lambda'|\geq R\}$\,. 
For \,$\lambda'\in \overline{U}$\,, we have \,$w(\lambda') \leq 2\,e^\delta\,w(\lambda)$\, by Proposition~\ref{P:asymp:w}(1),
and therefore by Cauchy's inequality, applied to the holomorphic function \,$a-a_0$\,:
\begin{align*}
|a'(\lambda)-\upsilon\,a_0'(\lambda)| & \leq \frac{1}{\delta\,|\lambda|^{1/2}} \max_{\lambda'\in\partial U} |a(\lambda')-\upsilon\,a_0(\lambda')| \\
& \leq \frac{1}{\delta\,|\lambda|^{1/2}} \cdot \frac{\eps\,\delta}{2\,e^\delta}\cdot 2\,e^\delta\,w(\lambda) = \eps\,|\lambda|^{-1/2}\,w(\lambda) \; . 
\end{align*}
For the case \,$\lambda\to0$\,, we take \,$U:=\Mengegr{\lambda'\in \C^*}{|\lambda'-\lambda|<\delta\,|\lambda|^{3/2}}$\,, and again choose \,$R'>R$\,, such that for \,$|\lambda|\leq \tfrac{1}{R'}$\,,
we have \,$\overline{U} \subset \{|\lambda'|\leq \tfrac{1}{R'}\}$\,. We again have \,$w(\lambda') \leq 2\,e^\delta\,w(\lambda)$\, for \,$\lambda'\in \overline{U}$\,, and therefore
by Cauchy's inequality
\begin{align*}
|a'(\lambda)-\upsilon^{-1}\,a_0'(\lambda)| & \leq \frac{1}{\delta\,|\lambda|^{3/2}} \max_{\lambda'\in\partial U} |a(\lambda')-\upsilon^{-1}\,a_0(\lambda')| \\
& \leq \frac{1}{\delta\,|\lambda|^{3/2}} \cdot \frac{\eps\,\delta}{2\,e^\delta}\cdot 2\,e^\delta\,w(\lambda) = \eps\,|\lambda|^{-3/2}\,w(\lambda) \; . 
\end{align*}
\end{proof}

For Cauchy data \,$(u,u_y)$\, that are once more differentiable, we obtain an even better asymptotic estimate.
For stating this asymptotic estimate, we introduce the space of (non-periodic) once more differentiable potentials
\begin{equation}
\label{eq:asymp:Potnp1}
\Pot_{np}^1 := \Menge{(u,u_y)}{u \in W^{2,2}([0,1]), u_y \in W^{1,2}([0,1])} \;,
\end{equation}
and phrase the statement in terms of the regauged monodromy \,$E$\, from the proof of Theorem~\ref{T:asymp:basic}:

\enlargethispage{3em}

\begin{prop}
\label{P:asymp:more}
Suppose that we have \,$(u,u_y) \in \Pot^1_{np}$\,. Then there exist constants
\,$R,C>0$\,, such that for every \,$\lambda\in \C$\, with \,$|\lambda|\geq R$\, and \,$x\in [0,1]$\,, we have
$$ |E_\lambda(x)-E_{0,\lambda}(x)| \leq \frac{C}{\sqrt{|\lambda|}}\,w(\lambda) \; . $$
If \,$P$\, is a relatively compact subset of \,$\Pot^1_{np}$\,, then 
the constants \,$R,C$\, can be chosen uniformly for \,$(u,u_y) \in P$\,.
\end{prop}

\begin{proof}
We use the objects and notations from the proof of Theorem~\ref{T:asymp:basic}. In particular
we use the regauging described at the beginning at that proof. But now we can use the higher
degree of differentiability of \,$(u,u_y)$\, 
to eliminate the \,$\lambda^0$-component from \,$\wt{\alpha}$\, 
by regauging \,$\wt{\alpha}$\, resp.~\,$\wt{F}$\, once more, with 
\begin{equation}
\label{eq:asymp:more:h}
h := \begin{pmatrix} \lambda^{-1/2}\,u_z & 1 \\ -1 & -\lambda^{-1/2}\,u_z \end{pmatrix} \; . 
\end{equation}
We have
\begin{equation}
\label{eq:asymp:more:h-1}
h^{-1} = \frac{1}{1-\lambda^{-1}\,u_z^2}\begin{pmatrix} -\lambda^{-1/2}\,u_z & -1 \\ 1 & \lambda^{-1/2}\,u_z \end{pmatrix}
\qmq{and} \mathrm{d}h = \begin{pmatrix} \lambda^{-1/2}\,u_{zx} & 0 \\ 0 & -\lambda^{-1/2}\,u_{zx} \end{pmatrix} \;, 
\end{equation}
and from this, we obtain for the again regauged \,$\wt{\alpha}$\, by an explicit calculation
\begin{align*}
\overline{\alpha} & := h^{-1}\,\wt{\alpha}\,h - h^{-1}\,\mathrm{d}h 
= \wt{\alpha}_0 + \lambda^{-1/2}\cdot (\overline{\beta}_++\overline{\beta}_-) + O(|\lambda|^{-1}) \\
& \qquad \text{ with } \overline{\beta}_+ := \begin{pmatrix} 0 & \tfrac12\,u_z^2-\tfrac14\,\cosh(u)+1 \\ -\tfrac12\,u_z^2 + \tfrac14\,\cosh(u)-1 & 0 \end{pmatrix} \\
& \qquad \text{ and } \overline{\beta}_- := \begin{pmatrix} 0 & -u_{zx}-\tfrac14\,\sinh(u) \\ -u_{zx}-\tfrac14\,\sinh(u) & 0 \end{pmatrix} \; . 
\end{align*}
We also consider \,$\overline{E}(x) := h^{-1}(x)\,E(x)\,h(0)$\,. From Equations~\eqref{eq:asymp:E0}, \eqref{eq:asymp:more:h} and \eqref{eq:asymp:more:h-1}, we obtain by an explicit calculation
\begin{equation}
\label{eq:asymp:more:h-1E0h}
h^{-1}\,E_0\,h = E_0 + O(|\lambda|^{-1/2}) \; . 
\end{equation}

By Lemma~\ref{L:asymp:compare}(2)(b) we have \,$\overline{E} = E_0 + \sum_{n=1}^\infty \overline{E}_n$\, with
$$ \overline{E}_n(x) = \lambda^{-n/2}\,\sum_{\eps\in\{\pm 1\}^n} \int_0^x \int_0^{t_1} \cdots \int_0^{t_{n-1}} \!\!\!\!\!\!\!\!\!\! E_0(\xi_\eps(t))\,\overline{\beta}_{\eps_1}(t_1)\,\overline{\beta}_{\eps_2}(t_2)\dotsc\,\overline{\beta}_{\eps_n}(t_n)\,\mathrm{d}^nt + O(|\lambda|^{-1}) \;, $$
where \,$\xi_\eps$\, is defined as in Lemma~\ref{L:asymp:compare}(2)(b). We have \,$\xi_\eps(t) \in [-x,x]$\, and therefore
$$ |E_0(\xi_\eps(t))| \leq |E_0(x)| \leq |E_0(1)| = w(\lambda) \; . $$
We therefore obtain 
$$ |\overline{E}_n(x)| \leq |\lambda|^{-n/2} \cdot 2^n \cdot \frac{1}{n!} \cdot w(\lambda) \cdot (\|\overline{\beta}_+\|_\infty + \|\overline{\beta}_-\|_\infty )^n \; .$$
It follows that the sum \,$\sum_{n=0}^\infty \overline{E}_n$\, converges uniformly to some \,$\overline{E}$\,, which is a solution of \,$\overline{E}'(x)=\overline{\alpha}(x)\,\overline{E}(x)$\,, \,$\overline{E}(0)=\unity$\,,
and we have
\begin{equation}
\label{eq:asymp:more:ovEasymp}
|\overline{E}(x)-\overline{E}_0(x)| \leq \left( \exp(2\,|\lambda|^{-1/2}\,(\|\overline{\beta}_+\|_\infty + \|\overline{\beta}_-\|_\infty )) - 1 \right)\cdot w(\lambda) \leq C_1\,|\lambda|^{-1/2}\,w(\lambda)
\end{equation}
with a suitable constant \,$C_1>0$\,. 

Because we have
$$ E(x)-E_0(x) = h(x)\,\bigr(\overline{E}(x)-\overline{E}_0(x) + E_0(x) - h(x)^{-1}\,E_0(x)\,h(0)\bigr)\,h(0)^{-1} \;; $$
it follows from the estimates \eqref{eq:asymp:more:ovEasymp} and \eqref{eq:asymp:more:h-1E0h}, along with the fact that \,$h$\, and \,$h^{-1}$\, are  bounded with respect to \,$\lambda\in\C^*$\,,
that the estimate given in the Proposition holds.
\end{proof}

\section{Basic behavior of the spectral data}
\label{Se:excl}

We would like to show that the spectral data \,$(\Sigma,\mathcal{D})$\, corresponding to a potential behave like the spectral data \,$(\Sigma_0,\mathcal{D}_0)$\, of the vacuum (as described in Section~\ref{Se:vacuum})
``asymptotically'' for \,$\lambda$\, near \,$\infty$\,, and for \,$\lambda$\, near \,$0$\,. In particular, we would like to show that the classical spectral divisor \,$D$\, corresponding to \,$\mathcal{D}$\,, 
like the classical spectral divisor \,$D_0$\, of the vacuum, is composed of a \,$\Z$-sequence of points
\,$(\lambda_k,\mu_k)_{k\in\Z}$\,, and that for \,$|k|$\, large, \,$(\lambda_k,\mu_k) \in D$\, is near \,$(\lambda_{k,0},\mu_{k,0})\in D_0$\,. 
Similarly, we will show that the set of zeros of \,$\Delta^2-4$\, with multiplicities (corresponding to the branch points resp.~singularities
of the spectral curve \,$\Sigma$\,) is enumerated by two \,$\Z$-sequences \,$(\vkap_{k,1})$\, and \,$(\vkap_{k,2})$\, such that for \,$|k|$\, large,
\,$\vkap_{k,1}$\, and \,$\vkap_{k,2}$\, are near \,$\lambda_{k,0}$\,. 

For a first quantified description of these asymptotic claims, 
we introduce the concept of so-called \emph{excluded domains} around the double points of \,$\Sigma_0$\,. 

For this, let \,$\delta>0$\, be given. If \,$\delta$\, is sufficiently small, certainly for \,$\delta < \pi-\tfrac12$\,,
the open set \,$\Mengegr{\lambda\in\C^*}{|\zeta(\lambda)-\zeta(\lambda_{k,0})|<\delta}$\, 
(where \,$\zeta(\lambda)$\, is here and henceforth again defined as in Equation~\eqref{eq:vacuum:zeta})
has two connected components for every \,$k\in \Z \setminus \{0\}$\,,
one of them contained in \,$\{|\lambda|>1\}$\,, the other in \,$\{|\lambda|<1\}$\,; whereas it is connected
for \,$k=0$\,. 

\begin{proof}
\,$O := \Mengegr{\lambda\in\C^*}{|\zeta(\lambda)-\zeta(\lambda_{k,0})|<\delta}$\, is invariant under
\,$\lambda\mapsto \lambda^{-1}$\,, thus it can be connected only if \,$O$\, contains at least one \,$\lambda\in \C^*$\,
with \,$|\lambda|=1$\,. Because of \,$\zeta(e^{it})=\tfrac14(e^{it/2}+e^{-it/2})=\tfrac12\,\cos(t/2)$\,,
we see that \,$|\lambda|=1$\, implies \,$|\zeta(\lambda)|\in [0,\tfrac12]$\,. Therefore \,$O$\, cannot contain
any \,$\lambda$\, with \,$|\lambda|=1$\, if \,$\delta < \pi-\tfrac12$\,.
\end{proof}

\begin{Def}
\label{D:vac2:excldom}
Let \,$0<\delta<\pi-\tfrac12$\, be given. Then we put for \,$k\in \Z$\,
$$ U_{k,\delta} := \begin{cases}
\Mengegr{\lambda \in \C^*}{\bigr|\zeta(\lambda)-\zeta(\lambda_{k,0}) \bigr|<\delta,\;|\lambda|>1} & \text{ for \,$k>0$\,} \\
\Mengegr{\lambda \in \C^*}{\bigr|\zeta(\lambda) \bigr|<\delta} & \text{ for \,$k=0$\,} \\
\Mengegr{\lambda \in \C^*}{\bigr|\zeta(\lambda)-\zeta(\lambda_{k,0}) \bigr|<\delta,\;|\lambda|<1} & \text{ for \,$k<0$\,} 
\end{cases} \; . $$
We call \,$U_{k,\delta}$\, the \emph{excluded domain} near \,$\lambda_{k,0}$\, of radius \,$\delta$\,. We also
consider the \emph{union \,$U_\delta := \bigcup_{k\in \Z} U_{k,\delta}$\, of all excluded domains} of radius \,$\delta$\, 
and the \emph{area outside the excluded domains} \,$V_\delta := \C^* \setminus U_\delta$\,. 
\end{Def}

In the sequel, we will use the domains \,$U_{k,\delta}$\,, \,$U_\delta$\, and \,$V_\delta$\, permanently
without further explicit reference.

\begin{prop}
\label{P:vac2:excldom-new}
There exist constants \,$C_1,C_2>0$\,, such that for every \,$0<\delta<\pi-\tfrac12$\,, every \,$k\in \Z$\, and every \,$\lambda_1,\lambda_2\in U_{k,\delta}$\, we have
\begin{enumerate}
\item 
For \,$k>0$\,: \,$C_1 \cdot |k|^{-1} \cdot |\lambda_1-\lambda_2| \leq |\zeta(\lambda_1)-\zeta(\lambda_2)| \leq C_2 \cdot |k|^{-1}\cdot |\lambda_1-\lambda_2|$\, \\
For \,$k<0$\,: \,$C_1 \cdot |k|^{3} \cdot |\lambda_1-\lambda_2| \leq |\zeta(\lambda_1)-\zeta(\lambda_2)| \leq C_2 \cdot |k|^{3}\cdot |\lambda_1-\lambda_2|$\, 
\item
For \,$k>0$\,: \,$B(\lambda_{k,0},\delta C_1 |k|) \;\subset\; U_{k,\delta} \;\subset\; B(\lambda_{k,0},\delta C_2 |k|)$\, \\
For \,$k<0$\,: \,$B(\lambda_{k,0},\delta C_1 |k|^{-3}) \;\subset\; U_{k,\delta} \;\subset\; B(\lambda_{k,0},\delta C_2 |k|^{-3})$\, \\
Here \,$B(z_0,r) := \Menge{z\in \C}{|z-z_0|<r}$\, is the open Euclidean ball of radius \,$r>0$\, around \,$z_0\in \C$\, in \,$\C$\,. 
\end{enumerate}
\end{prop}

\begin{proof}
\emph{For (1).}
For \,$k>0$\, we have
\begin{align*}
\zeta(\lambda_1)-\zeta(\lambda_2)
& = \frac{1}{4}\cdot \left(\sqrt{\lambda_1}-\sqrt{\lambda_2}+\frac{1}{\sqrt{\lambda_1}}-\frac{1}{\sqrt{\lambda_2}} \right) \\
& = \frac{1}{4}\cdot \left( 1- \frac{1}{\sqrt{\lambda_1\,\lambda_2}} \right) \cdot \frac{1}{\sqrt{\lambda_1}+\sqrt{\lambda_2}} \cdot (\lambda_1-\lambda_2) \; .
\end{align*}
Because of \,$\lambda_\nu \in U_{k,\delta}$\, for \,$\nu\in\{1,2\}$\, and \,$k>0$\,, we have \,$|\lambda_\nu|>2$\,, and therefore \,$\tfrac12 \leq \left| 1-\tfrac{1}{\sqrt{\lambda_1\,\lambda_2}}\right| \leq \tfrac32$\,, and moreover,
\,$\left| \frac{1}{\sqrt{\lambda_1}+\sqrt{\lambda_2}} \right|$\, is of the order of \,$k$\,. Herefrom, the claimed statement follows.

For \,$k<0$\,, we use the fact that \,$\zeta(\lambda_\nu)=\zeta(\lambda_\nu^{-1})$\, holds, and that we have \,$\lambda_\nu^{-1} \in U_{-k,\delta}$\, where \,$-k>0$\,. By application of the first part of the proof,
we thus obtain
$$ C_1 \cdot |k|^{-1} \cdot |\lambda_1^{-1}-\lambda_2^{-1}| \leq |\zeta(\lambda_1)-\zeta(\lambda_2)| \leq C_2 \cdot |k|^{-1} \cdot |\lambda_1^{-1}-\lambda_2^{-1}| \; . $$
Because of \,$|\lambda_1^{-1}-\lambda_2^{-1}| = \tfrac{|\lambda_1-\lambda_2|}{|\lambda_1\,\lambda_2|}$\, and because \,$\tfrac{1}{|\lambda_\nu|}$\, is of order \,$k^2$\,, we obtain the claimed statement also in this case.

\emph{For (2).}
This is an immediate consequence of (1).
\end{proof}

\begin{prop}
\label{P:vac2:excldom-M0}
Let \,$0<\delta<\pi-\tfrac12$\, be given. 
\begin{enumerate}
\item \,$w(\lambda)$\, is bounded on \,$U_\delta$\,, more precisely we have for every \,$\lambda\in U_\delta$\,:
$$ \frac12 \leq w(\lambda) \leq 2\,e^\delta \; . $$
\item \,$\cot\circ \zeta$\, is bounded on \,$V_\delta$\,.
\item There exists a constant \,$C>1$\, (dependent on \,$\delta$\,) so that \,$|\sin(\zeta(\lambda))| \leq w(\lambda) \leq C\,|\sin(\zeta(\lambda))|$\, holds for all \,$\lambda\in V_\delta$\,.
\end{enumerate}
\end{prop}

\begin{proof}
\emph{For (1).}
Let \,$\lambda\in U_\delta$\, be given. Then there exists \,$k\in \Z$\, with \,$\lambda\in U_{k,\delta}$\, and therefore \,$|\zeta(\lambda)-\zeta(\lambda_{k,0})|<\delta$\,. We have \,$\zeta(\lambda_{k,0})=|k|\pi\in \R$\,, and therefore
it follows that \,$|\IM(\zeta(\lambda))| < \delta$\,. Thus we obtain from Proposition~\ref{P:asymp:w}(1):
$$ \frac12 \leq \frac12\,e^{|\IM(\zeta(\lambda))|} \leq w(\lambda) \leq 2\,e^{|\IM(\zeta(\lambda))|} \leq 2\,e^\delta \; . $$

\emph{For (2).} 
It suffices to show that \,$\cot$\, is bounded on \,$\C \setminus \bigcup_{k\in \Z} B(k\pi,\delta)$\,, and for this,
it in turn suffices to show that \,$\cot$\, is bounded on the sets
$$ M_1 := \Mengegr{x+iy \in \C}{|x-k\pi|\geq \tfrac{\delta}{2} \text{ for all \,$k\in \Z$\, and } |y| \leq \tfrac{\delta}{2}}  $$
and 
$$ M_2 := \Mengegr{x+iy \in \C}{|y| \geq \tfrac{\delta}{2}} \;. $$

\,$M_1$\, is the union of translations by integer multiples of \,$\pi$\, of the rectangle
\,$Q := [\tfrac{\delta}{2},\pi-\tfrac{\delta}{2}]\times[-\tfrac{\delta}{2},\tfrac{\delta}{2}]$\,. \,$\cot$\, is bounded
on the compact set \,$Q$\,; because it is \,$\pi$-periodic, it follows that \,$\cot$\, is also bounded on all 
of \,$M_1$\,. 

To show that \,$\cot$\, is bounded on \,$M_2$\,, we use the formula
$$ \cot(x+iy) = \frac{\sin(2x)}{\cosh(2y) - \cos(2x)} - i\, \frac{\sinh(2y)}{\cosh(2y)-\cos(2x)} \;. $$
If \,$x+iy\in M_1$\, is given, we have
$$ \left| \frac{\sin(2x)}{\cosh(2y)-\cos(2x)} \right| \leq \frac{1}{\cosh(\delta)-1} $$
and
$$ \left| \frac{\sinh(2y)}{\cosh(2y)-\cos(2x)} \right| \leq \frac{|\sinh(2y)|}{\cosh(2y)-1} \;; $$
note that \,$y\mapsto \tfrac{|\sinh(2y)|}{\cosh(2y)-1}$\, is real for \,$y\neq 0$\,, and tends to \,$1$\, for
\,$y\to \pm \infty$\,, hence is bounded on \,$\{|y|\geq \tfrac{\delta}{2}\}$\,. Thus it follows that 
\,$\cot$\, is bounded on \,$M_2$\,.

\emph{For (3).} Because of 
\begin{equation}
\label{eq:vac2:excldom-M0:M0}
w(\lambda) = |\cos(\zeta(\lambda))| + |\sin(\zeta(\lambda))| \;,
\end{equation}
the inequality \,$|\sin(\zeta(\lambda))| \leq w(\lambda)$\, is obvious. On the other hand, \,$\cot \circ \zeta$\, is bounded on \,$V_\delta$\, by (2), thus there exists \,$C>1$\, with \,$|\cot(\zeta(\lambda))| \leq C-1$\, for \,$\lambda\in V_\delta$\,.
It follows from Equation~\eqref{eq:vac2:excldom-M0:M0} that we have
$$ w(\lambda) = |\sin(\zeta(\lambda))| \cdot \left( |\cot(\zeta(\lambda))| + 1 \right) \leq |\sin(\zeta(\lambda))|\cdot C \; . $$

\end{proof}

We now suppose that ``a monodromy'' that has the asymptotic behavior of Theorem~\ref{T:asymp:basic} is given, 
that is to say: Suppose that a \,$(2\times 2)$-matrix 
$$ M(\lambda) = \left( \begin{matrix} a(\lambda) & b(\lambda) \\ c(\lambda) & d(\lambda)\end{matrix} \right) $$
of holomorphic functions \,$a,b,c,d: \C^* \to \C$\, is given, and denote by
$$ M_0(\lambda)  = \begin{pmatrix} a_0(\lambda) & b_0(\lambda) \\ c_0(\lambda) & d_0(\lambda) \end{pmatrix} 
= \begin{pmatrix} \cos(\zeta(\lambda)) & -{\lambda}^{-1/2}\,\sin(\zeta(\lambda)) \\ {\lambda}^{1/2}\,\sin(\zeta(\lambda)) & \cos(\zeta(\lambda)) \end{pmatrix} \; , $$
the monodromy of the vacuum, compare Equation~\eqref{eq:vacuum:M0}. Then we require of \,$M(\lambda)$\, that \,$\det(M(\lambda))=1$\, holds for
all \,$\lambda\in \C^*$\, and that there exist 
constants \,$\tau,\upsilon \in \C^*$\,, such that for every \,$\eps>0$\, there exists \,$R>0$\, so that for every \,$\lambda\in \C^*$\, with
\,$|\lambda| \geq R$\, resp.~with \,$|\lambda|\leq \tfrac{1}{R}$\, the estimates of Theorem~\ref{T:asymp:basic}(1) resp.~(2) hold.
We note that in this setting, the Corollaries~\ref{C:asymp:basic-deltabc} and \ref{C:asymp:Mprime} also hold.

In this setting, we also use the trace functions \,$\Delta(\lambda) := a(\lambda)+d(\lambda)$\, 
and \,$\Delta_0(\lambda) = a_0(\lambda)+d_0(\lambda)=2\cos(\zeta(\lambda))$\, again. In the case \,$\upsilon=1$\, we define the spectral curve \,$\Sigma$\,, the
generalized spectral divisor \,$\mathcal{D}$\, associated to \,$M(\lambda)$\, and the underlying classical divisor \,$D$\, 
in the way described in Section~\ref{Se:spectrum}. For \,$\upsilon\neq 1$\, we avoid defining the spectral curve (as it would not be of any interest, because it would not
be asymptotically close to any ``reasonable'' curve), but we still define the classical spectral divisor \,$D$\, as a multi-set of points in \,$\C^*\times \C^*$\, by Equation~\eqref{eq:spectral:D-classical}. 

Theorem~\ref{T:asymp:basic} shows that this situation is at hand in particular when a potential \,$(u,u_y) \in \Pot_{np}$\, is given
and \,$M(\lambda)$\, is the monodromy associated to it; for \,$(u,u_y) \in \Pot$\, we have \,$\upsilon=1$\,. 

\newpage

\label{not:excl:excldom-Sigma}
Whenever we have a spectral curve \,$\Sigma$\,, we consider also the pre-images \,$\wh{U}_{k,\delta}$\,, \,$\wh{U}_\delta$\, and \,$\wh{V}_\delta$\, 
in \,$\Sigma$\, of the 
excluded domains \,$U_{k,\delta}$\,, their union \,$U_\delta := \bigcup_{k\in \Z} U_{k,\delta}$\, and the complement \,$V_\delta := \C^* \setminus U_\delta$\,.
More explicitly, we define for \,$\delta>0$\, 
and \,$k\in \Z$\, 
\begin{align*}
\wh{U}_{k,\delta} & := \Menge{(\lambda,\mu)\in \Sigma}{\lambda\in U_{k,\delta}} \;, \\
\wh{U}_{\delta} & := \Menge{(\lambda,\mu)\in \Sigma}{\lambda\in U_{\delta}} = \bigcup_{k\in \Z} \wh{U}_{k,\delta} \;, \\
\wh{V}_{\delta} & := \Menge{(\lambda,\mu)\in \Sigma}{\lambda\in V_{\delta}} = \Sigma \setminus \wh{U}_\delta \; . 
\end{align*}
We call also the \,$\wh{U}_{k,\delta}$\, \emph{excluded domains} (in \,$\Sigma$\,). 

\begin{Def}
Let \,$X\subset \C^*$\, resp.~\,$\wh{X}\subset \Sigma$\, be any subset, \,$\delta>0$\, 
and \,$n\in \N_0$\,. 
\begin{enumerate}
\item We say that the excluded domains \,$U_{k,\delta}$\, resp.~\,$\wh{U}_{k,\delta}$\, \emph{asymptotically} contain (exactly) \,$n$\, elements of \,$X$\, resp.~of \,$\wh{X}$\,, if there exists \,$k_0\in \N$\, such that for all \,$k\in \Z$\, with \,$|k| \geq k_0$\,
we have \,$\#(X\cap U_{k,\delta})=n$\, resp.~\,$\#(\wh{X}\cap \wh{U}_{k,\delta})=n$\,.
\item We say that the excluded domains \,$U_{k,\delta}$\, resp.~\,$\wh{U}_{k,\delta}$\, \emph{asymptotically and totally} contain (exactly) \,$n$\, elements of \,$X$\, resp.~of \,$\wh{X}$\,, if there exists \,$k_0\in \N$\, such that for all \,$k\in \Z$\, with \,$|k| \geq k_0$\,
we have
$$ \#(X\cap U_{k,\delta})=n \qmq{and} \#(X \setminus \bigcup_{|k|\geq k_0} U_{k,\delta}) = n\cdot (2k_0-1) $$
resp.~
$$ \#(\wh{X}\cap \wh{U}_{k,\delta})=n \qmq{and} \#(\wh{X} \setminus \bigcup_{|k|\geq k_0} \wh{U}_{k,\delta}) = n\cdot (2k_0-1) \;. $$
\end{enumerate}
\end{Def}

The following statement, which establishes the fundamental structure of the spectral data, is analogous to the ``Counting Lemma'' 
in the treatment of the 1-dimensional Schr\"odinger equation in \cite{Poeschel-Trubowitz:1987} (Lemma~2.2, p.~27).

\begin{prop}
\label{P:excl:basic}
\begin{enumerate}
\item
In the case \,$\upsilon=1$\,, the excluded domains \,$U_{k,\delta}$\, asymptotically and totally contain exactly two zeros (with multiplicity) of \,$\Delta^2-4$\,.

Therefore the zeros of \,$\Delta^2-4$\, can then be enumerated by two sequences \,$(\vkap_{k,1})_{k\in \Z}$\, and \,$(\vkap_{k,2})_{k\in \Z}$\, in \,$\C^*$\, 
in such a way that for every
\,$0<\delta<\pi-\tfrac12$\, there exists \,$N\in \N$\, so that \,$\vkap_{k,1},\vkap_{k,2} \in U_{k,\delta}$\, for all \,$k\in \Z$\, with \,$|k|\geq N$\,.
\item The excluded domains \,$U_{k,\delta}$\, asymptotically and totally contain exactly one zero (with multiplicity) of \,$c$\,. 

Therefore the zeros of \,$c$\, can be enumerated by a sequence \,$(\lambda_k)_{k\in \Z}$\, in \,$\C^*$\, in such a way 
that for every \,$0<\delta<\pi-\tfrac12$\, there exists \,$N\in \N$\, so that
\,$\lambda_k \in U_{k,\delta}$\, holds for all \,$k\in \Z$\, with \,$|k|\geq N$\,. 
\item We have \,$D=\Menge{(\lambda_k,\mu_k)}{k\in \Z}$\,, where we put \,$\mu_k := a(\lambda_k)$\,. 
Then we have
$$ \lim_{k\to\infty} (\mu_k-\upsilon\cdot \mu_{k,0}) = \lim_{k\to-\infty} (\mu_k-\upsilon^{-1}\cdot \mu_{k,0}) = 0 \; . $$
\end{enumerate}
These estimates hold uniformly for the set of monodromies \,$M(\lambda)$\, corresponding to a relatively compact set \,$P$\, 
of potentials in \,$\Pot_{np}$\,. 
\end{prop}

\begin{proof}
\emph{For (1).} 
Let \,$0<\delta<\pi-\tfrac12$\, be given, and let \,$C_1>1$\, be the constant from Proposition~\ref{P:vac2:excldom-M0}(3) so that
\begin{equation}
\label{eq:excl:basic:w-sin}
\forall \lambda\in V_\delta \; : \; |\sin(\zeta(\lambda))| \leq w(\lambda) \leq C_1 \cdot |\sin(\zeta(\lambda))|
\end{equation}
holds.
Because of \,$\upsilon=1$\,, by Corollary~\ref{C:asymp:basic-deltabc}(2)(c) there then exists \,$R>0$\,, 
such that for every \,$\lambda\in\C^*$\, with \,$|\lambda|\geq R$\, or \,$|\lambda|\leq\tfrac{1}{R}$\,,
we have
$$ |(\Delta(\lambda)^2-4)-(\Delta_0(\lambda)^2-4)| \leq \frac{2}{C_1^2}\, w(\lambda)^2 \; . $$
If we additionally have \,$\lambda \in V_\delta$\,, then it follows by the estimate \eqref{eq:excl:basic:w-sin} that
\begin{equation}
\label{eq:excl:basic:estimate}
|(\Delta(\lambda)^2-4)-(\Delta_0(\lambda)^2-4)| \leq \frac{2}{C_1^2} \cdot (C_1\,|\sin(\zeta(\lambda))|)^2 = \frac12 |\Delta_0(\lambda)^2-4| < |\Delta_0(\lambda)^2-4|
\end{equation}
holds. 

There exists \,$N\in \N$\, such that the boundaries of the excluded domains \,$U_{k,\delta}$\, with \,$|k|> N$\, and moreover the boundaries of the topological annuli
$$ T_k := \Mengegr{\lambda\in\C^*}{|\zeta(\lambda)| \leq \pi \, \left(k+\tfrac12 \right)} $$
with \,$k\geq N$\, are contained in \,$V_\delta \cap \{|\lambda|\geq R \text{ or } |\lambda|\leq \tfrac{1}{R}\}$\,. It then follows from \eqref{eq:excl:basic:estimate} by Rouch\'e's Theorem that \,$\Delta^2-4$\, and \,$\Delta_0^2-4$\, 
have the same number of zeros on each \,$U_{k,\delta}$\, and on each \,$T_k$\, with \,$|k|\geq N$\,. 
Therefore \,$\Delta^2-4$\, has two zeros on each \,$U_{k,\delta}$\, with \,$|k|\geq N$\,, and \,$2\cdot(2k+1)$\, zeros on each \,$T_k$\, with \,$|k|\geq N$\,.

Herefrom (1) follows: 
For \,$|k|\geq N$\,, we let \,$\vkap_{k,1}$\, and \,$\vkap_{k,2}$\, be the two zeros of \,$\Delta^2-4$\, on \,$U_{k,\delta}$\,, and then we name the remaining \,$2\cdot(2N+1)$\, 
zeros as \,$\vkap_{k,1}, \vkap_{k,2}$\,, where \,$k\in \{-N,\dotsc,N\}$\,.

\emph{For (2).}
We would like to apply a similar argument to \,$c$\, as we used for \,$\Delta^2-4$\, in the proof of (1). However, in the basic asymptotics from Theorem~\ref{T:asymp:basic}, 
\,$c$\, in general approximates different functions for \,$\lambda\to\infty$\,
and for \,$\lambda\to0$\,, namely \,$\tau\,c_0$\, resp.~$\tau^{-1}\,c_0$\,. 
But to be able to apply Rouch\'e's Theorem especially for the topological annuli \,$T_k$\,, we need
a comparison function that approximates \,$c$\, both near \,$\lambda\to \infty$\, and \,$\lambda\to0$\,.

To obtain such a comparison function, we consider the constant potential \,$(\wt{u},\wt{u}_y) \in \Pot$\, with \,$\wt{u} = -2\,\ln(\tau)$\,, 
\,$\wt{u}_y=0$\,, and the associated monodromy \,$\wt{M}(\lambda) = \left( \begin{smallmatrix} \wt{a}(\lambda) & \wt{b}(\lambda) \\ \wt{c}(\lambda) & \wt{d}(\lambda) \end{smallmatrix} \right)$\,.
Then we also have \,$e^{-(\wt{u}(0)+\wt{u}(1))/4}=\tau$\,, and therefore for given \,$\eps>0$\,
there exists \,$R>0$\, such that for \,$|\lambda|\geq R$\, we have
$$ \qmq{both} |c(\lambda)-\tau\,c_0(\lambda)| \leq \eps\,|\lambda|^{1/2}\,w(\lambda) \qmq{and} |\wt{c}(\lambda)-\tau\,c_0(\lambda)| \leq \eps\,|\lambda|^{1/2}\,w(\lambda) \; $$
(the first inequality by our hypothesis on \,$M(\lambda)$\,, and the second inequality from Theorem~\ref{T:asymp:basic}),
wherefrom we obtain
\begin{equation}
\label{eq:excl:basic:cwtc-pre}
|c(\lambda)-\wt{c}(\lambda)| \leq 2\eps\,|\lambda|^{1/2}\,w(\lambda) \; .
\end{equation}
If we additionally have \,$\lambda\in V_\delta$\,, then we have \,$|\tau\,c_0(\lambda)| \leq |\wt{c}(\lambda)| + |\wt{c}(\lambda)-\tau\,c_0(\lambda)| \leq |\wt{c}(\lambda)|+\eps\,C_1\,|c_0(\lambda)|$\,
and therefore
$$ |\lambda|^{1/2}\,w(\lambda) \leq C_1\,|c_0(\lambda)| \leq \frac{C_1}{|\tau|-\eps\,C_1} \, |\wt{c}(\lambda)| \; . $$
By combining this estimate with \eqref{eq:excl:basic:cwtc-pre}, we obtain
$$ |c(\lambda)-\wt{c}(\lambda)| \leq \frac{2\,C_1\,\eps}{|\tau|-\eps\,C_1} \, |\wt{c}(\lambda)| \; . $$
By choosing \,$\eps>0$\, sufficiently small, we can thus achieve for all \,$\lambda\in V_\delta$\,, \,$|\lambda|\geq R$\,
$$ |c(\lambda)-\wt{c}(\lambda)| \leq \frac{1}{2} \, |\wt{c}(\lambda)| < |\wt{c}(\lambda)| \; . $$
By a similar argument, we obtain the same estimate also for \,$\lambda\in V_\delta$\, with \,$|\lambda|\leq \tfrac{1}{R}$\,. By an analogous application of Rouch\'e's Theorem as in the proof of (1),
it follows that there exists \,$N\in \N$\,, so that \,$c$\, and \,$\wt{c}$\, have the same number of zeros on the excluded domains \,$U_{k,\delta}$\, with \,$|k|\geq N$\, and on the topological annuli
\,$T_{k}$\, with \,$k\geq N$\,. To complete the proof of (2), it therefore suffices to show that the statement of (2) holds for \,$\wt{c}$\, in the place of \,$c$\,. 

For this purpose, we calculate \,$\wt{c}$\, and its zeros explicitly. Proceeding similarly as in the investigation of the vacuum in Section~\ref{Se:vacuum}, we note that the \,$1$-form \,$\wt{\alpha}_\lambda$\, associated to \,$(\wt{u},\wt{u}_y)$\, (Equation~\eqref{eq:asymp:alpha}) is independent of \,$x$\,,
and therefore \,$\wt{M}(\lambda) = \exp(\wt{\alpha}_{\lambda})$\, holds. By an explicit calculation similar to the one in Section~\ref{Se:vacuum}, we obtain
$$ \wt{c}(\lambda) = \sqrt{\frac{\lambda\,\tau+\tau^{-1}}{\lambda^{-1}\,\tau+\tau^{-1}}}\cdot \sin(\xi(\lambda)) \qmq{with} \xi(\lambda) = \frac14\,\sqrt{(\lambda\,\tau+\tau^{-1})\cdot(\lambda^{-1}\,\tau+\tau^{-1})} \; . $$
\,$\lambda\in \C^*$\, is a zero of \,$\wt{c}$\, if and only if 
$$ \qmq{either} \xi(\lambda)^2 = (k\pi)^2 \text{ with \,$k\in \N$\,} \qmq{or} \sqrt{\frac{\lambda\,\tau+\tau^{-1}}{\lambda^{-1}\,\tau+\tau^{-1}}} \cdot \xi(\lambda)=0 $$
holds. For \,$k\in \N$\,, the equation \,$\xi(\lambda)^2=(k\pi)^2$\, yields the quadratic equation 
$$ \lambda^2 - (16\,\pi^2\,k^2-\tau^2-\tau^{-2})\cdot \lambda + 1 = 0 $$
and thereby the two zeros of \,$\wt{c}$\,
$$ \wt{\lambda}_{\pm k} := \frac{1}{2}\cdot \left( \, 16\,\pi^2\,k^2-\tau^2-\tau^{-2} \pm \sqrt{(16\,\pi^2\,k^2-\tau^2-\tau^{-2})^2-4} \,\right) $$
By comparing \,$\wt{\lambda}_{\pm k}$\, with the asymptotic assessments of \,$\lambda_{k,0}$\, in \eqref{eq:vacuum:lambdak0-asymp}, we see that for \,$k\to\infty$\,
\,$\wt{\lambda}_k - \lambda_{k,0} = O(1)$\, and \,$\wt{\lambda}_{-k}-\lambda_{-k,0} = O(k^{-4})$\, holds. Therefore we have \,$\wt{\lambda}_k \in U_{k,\delta}$\, for \,$|k|$\, sufficiently large.

Finally, the equation \,$\sqrt{\tfrac{\lambda\,\tau+\tau^{-1}}{\lambda^{-1}\,\tau+\tau^{-1}}} \cdot \xi(\lambda)=0$\, is equivalent to \,$\lambda\,\tau+\tau^{-1}=0$\,, which yields one further zero:
\,$\wt{\lambda}_0 = -\tau^{-2}$\,. This shows that \,$\wt{c}$\, satisfies the property of (2), and therefore the proof of (2) is completed.

\emph{For (3).}
Let \,$\eps>0$\, be given. For \,$k>0$\,, we have
\begin{align*}
\mu_k - \upsilon\,\mu_{k,0} & = a(\lambda_k)-\upsilon\,a_0(\lambda_{k,0}) = a(\lambda_k)-\upsilon\,a_0(\lambda_k) + \upsilon\cdot \left( a_0(\lambda_k)-a_0(\lambda_{k,0}) \right) \\
& = a(\lambda_k) - \upsilon\,a_0(\lambda_k) + \upsilon\cdot \int_{\lambda_{k,0}}^{\lambda_k} a_0'(\lambda)\,\mathrm{d}\lambda \; . 
\end{align*}
We now fix \,$0<\delta<\pi-\tfrac12$\, (to be chosen later).
By the asymptotic behavior of \,$a$\, and the fact that \,$w(\lambda)$\, is bounded on \,$U_\delta$\, (Proposition~\ref{P:vac2:excldom-M0}(1)),
there exists \,$N\in \N$\, such that for all \,$k\geq N$\, we have
$$ |a(\lambda_k)-\upsilon\,a_0(\lambda_k)| \leq \frac{\eps}{2} \; . $$
We now estimate \,$\int_{\lambda_{k,0}}^{\lambda_k} a_0'(\lambda)\,\mathrm{d}\lambda$\,:
$$ \left| \int_{\lambda_{k,0}}^{\lambda_k} a_0'(\lambda)\,\mathrm{d}\lambda \right| \leq |\lambda_k-\lambda_{k,0}| \cdot \max_{\lambda\in U_{k,\delta}} |a_0'(\lambda)| \; . $$
We have \,$a_0'(\lambda) = -\frac18\,(\lambda^{-1/2}-\lambda^{-3/2})\,\sin(\zeta(\lambda))$\, by Equation~\eqref{eq:asymp:a0'}; 
because \,$\sin$\, is bounded on any horizontal strip in the \,$\zeta$-plane, there exists \,$C>0$\, (independent of the choice 
of \,$\delta<\pi-\tfrac12$\,), such that \,$\max_{\lambda\in U_{k,\delta}} |a_0'(\lambda)|\leq \tfrac{C}{k}$\, holds
for all \,$k\geq N$\,. We now set \,$\delta:=\tfrac{\eps}{2C}$\,, then it follows 
from (2) and Proposition~\ref{P:vac2:excldom-new}(2) that by increasing \,$N$\, if needed, we can obtain 
\,$|\lambda_k-\lambda_{k,0}|\leq \tfrac{\eps}{2C\,|\upsilon|}\,k$\, for \,$k\geq N$\,. For such \,$k$\, we then have \,$\left| \int_{\lambda_{k,0}}^{\lambda_k} a_0'(\lambda)\,\mathrm{d}\lambda \right| \leq \tfrac{\eps}{2\,|\upsilon|}$\,
and thus 
$$ |\mu_k-\upsilon\cdot \mu_{k,0}| \leq |a(\lambda_k)-\upsilon\,a_0(\lambda_k)| + |\upsilon|\cdot \left| \int_{\lambda_{k,0}}^{\lambda_k} a_0'(\lambda)\,\mathrm{d}\lambda \right| \leq \frac{\eps}{2} + |\upsilon|\cdot \frac{\eps}{2\,|\upsilon|} = \eps \;. $$
This shows that \,$\lim_{k\to\infty} (\mu_k-\upsilon\cdot \mu_{k,0})=0$\, holds. The limit \,$\lim_{k\to-\infty} (\mu_k-\upsilon^{-1}\cdot \mu_{k,0})=0$\, 
is shown similarly.

Finally, we note that the fact that the estimates (1)--(3) hold uniformly for a relatively compact set \,$P$\, of potentials follows from the corresponding statements in the propositions
in Section~\ref{Se:asymp}.
\end{proof}

\begin{rem}
\label{R:excl:constpot}
The trace function \,$\wt{\Delta}$\, corresponding to the constant potential \,$(\wt{u},\wt{u}_y)$\, studied in the proof of Proposition~\ref{P:excl:basic}(2) equals
\,$\wt{\Delta}(\lambda)=2\,\cos(\xi(\lambda))$\,, and thus we have \,$\wt{\Delta}^2(\lambda)-4=-4\,\sin(\xi(\lambda))^2$\,. 
It follows that for the zeros \,$\wt{\vkap}_{k,\nu}$\, 
of \,$\wt{\Delta}^2-4$\, we have for \,$k\neq 0$\,: \,$\wt{\vkap}_{k,1}=\wt{\vkap}_{k,2}=\wt{\lambda}_k$\,. 
However, for \,$k=0$\,, we have \,$\wt{\vkap}_{k,1}=-\tau^{-2}\neq -\tau^2=\wt{\vkap}_{k,2}$\, (for \,$\tau\not\in\{\pm 1,\pm i\}$\,). 

This shows that the spectral curve \,$\wt{\Sigma}$\, corresponding to \,$\wt{u}$\, has geometric genus \,$1$\,: 
Compared to the spectral curve of the vacuum, which has only double points, no branch points, and therefore geometric genus \,$0$\,,
one double point (namely the one for \,$k=0$\,) has been ``opened'' into a pair of branch points.

The minimal surfaces in \,$S^3$\, of spectral genus \,$1$\,, and how they arise from surfaces 
of spectral genus \,$0$\, by means of deformation,
have been studied extensively by \textsc{Kilian}, \textsc{Schmidt} and \textsc{Schmitt}
in \cite{Kilian-Schmidt-Schmitt:2013}. Such minimal surfaces are ``bubbletons'' (i.e.~solitons of the sinh-Gordon equation);
they are constant along all lines parallel to the \,$x$-axis.
Moreover, in \cite{Kilian-Schmidt-Schmitt:2013}, Section~3.3, 
an explicit presentation of the solution of the
sinh-Gordon equation corresponding to our potential \,$(\wt{u},\wt{u}_y)$\, in terms of elliptic functions
is given.
\end{rem}

\begin{prop}[Uniqueness statements.]
\label{P:excl:unique}
\begin{enumerate}
\item Suppose \,$c,\wt{c}:\C^*\to \C$\, are two holomorphic functions with the following asymptotic behavior: There are constants \,$\tau,\wt{\tau} \in \C^*$\, such that for every \,$\eps>0$\, there exists \,$R>0$\,,
such that we have for all \,$\lambda \in \C$\, with \,$|\lambda|\geq R$\,
$$ |c(\lambda)-\tau\,c_0(\lambda)| \leq \eps\,|\lambda|^{1/2}\,w(\lambda) \qmq{resp.} |\wt{c}(\lambda)-\wt{\tau}\,c_0(\lambda)| \leq \eps\,|\lambda|^{1/2}\,w(\lambda) $$
and for all \,$\lambda\in\C^*$\, with \,$|\lambda|\leq \tfrac{1}{R}$\, 
$$ |c(\lambda)-\tau^{-1}\,c_0(\lambda)| \leq \eps\,|\lambda|^{1/2}\,w(\lambda) \qmq{resp.} |\wt{c}(\lambda)-\wt{\tau}^{-1}\,c_0(\lambda)| \leq \eps\,|\lambda|^{1/2}\,w(\lambda) \; . $$
If \,$c$\, and \,$\wt{c}$\, have the same zeros (counted with multiplicity), then either \,$c=\wt{c}$\,, \,$\tau=\wt{\tau}$\, or \,$c=-\wt{c}$\,, \,$\tau=-\wt{\tau}$\, holds.
\item Suppose further that \,$a,\wt{a}:\C^*\to \C$\, are two holomorphic functions with the following asymptotic behavior: 
There are constants \,$\upsilon,\wt{\upsilon} \in \C^*$\, such that for every \,$\eps>0$\, there exists \,$R>0$\,,
such that we have for all \,$\lambda \in \C$\, with \,$|\lambda|\geq R$\,
$$ |a(\lambda)-\upsilon\,a_0(\lambda)| \leq \eps\,w(\lambda) \qmq{resp.} |\wt{a}(\lambda)-\wt{\upsilon}\,a_0(\lambda)| \leq \eps\,w(\lambda) $$
and for all \,$\lambda\in\C^*$\, with \,$|\lambda|\leq \tfrac{1}{R}$\, 
$$ |a(\lambda)-\upsilon^{-1}\,a_0(\lambda)| \leq \eps\,w(\lambda) \qmq{resp.} |\wt{a}(\lambda)-\wt{\upsilon}^{-1}\,a_0(\lambda)| \leq \eps\,w(\lambda) \; . $$
If we have \,$a(\lambda)=\wt{a}(\lambda)$\, and \,$\ord_\lambda(a-\wt{a}) \geq \ord_\lambda(c)$\, for every \,$\lambda\in\C^*$\, with \,$c(\lambda)=0$\,, then \,$a=\wt{a}$\, and \,$\upsilon=\wt{\upsilon}$\, holds.
\item Suppose that \,$f,\wt{f}:\C^*\to \C$\, are two holomorphic functions with the following asymptotic behavior: There are constants \,$\tau_\pm,\wt{\tau}_\pm \in \C^*$\, such that for every \,$\eps>0$\, there exists
\,$R>0$\, such that we have for all \,$\lambda\in \C$\, with \,$|\lambda|\geq R$\,
$$ |f(\lambda)-\tau_+\,(\Delta_0(\lambda)^2-4)| \leq \eps\,w(\lambda)^2 \qmq{resp.} |\wt{f}(\lambda)-\wt{\tau}_+\,(\Delta_0(\lambda)^2-4)| \leq \eps\,w(\lambda)^2 $$
and for all \,$\lambda\in \C^*$\, with \,$|\lambda|\leq \tfrac{1}{R}$\, 
$$ |f(\lambda)-\tau_-\,(\Delta_0(\lambda)^2-4)| \leq \eps\,w(\lambda)^2 \qmq{resp.} |\wt{f}(\lambda)-\wt{\tau}_-\,(\Delta_0(\lambda)^2-4)| \leq \eps\,w(\lambda)^2 \;. $$
If \,$f$\, and \,$\wt{f}$\, have the same zeros (counted with multiplicity), then there is a constant \,$A\in \C^*$\, such that \,${f}=A\cdot \wt{f}$\, and \,${\tau}_\pm = A\cdot \wt{\tau}_\pm$\, holds.
\end{enumerate}
\end{prop}

\begin{proof}
\emph{For (1).}
Because \,$c$\, and \,$\widetilde{c}$\, have the same zeros, \,$\vi := c/\widetilde{c}$\, is holomorphic and non-zero on all of \,$\C^*$\,. 
For fixed \,$\eps,\delta>0$\,, there exists \,$R>0$\, so that for 
\,$\lambda\in V_\delta$\, with \,$|\lambda|\geq R$\, we have
$$ |c(\lambda)-\tau\,c_0(\lambda)| \leq \eps\,|\lambda|^{1/2}\,w(\lambda) \leq 
\eps\,C_1\,|c_0(\lambda)| \;, $$
where \,$C_1>1$\, is the constant from Proposition~\ref{P:vac2:excldom-M0}(3), and hence
\begin{equation}
\label{eq:f-zero:uniq:f}
|c(\lambda)| \leq \left( |\tau| + \eps\,C_1 \right) \cdot |c_0(\lambda)| \; . 
\end{equation}
Likewise, from
$$ |\widetilde{c}(\lambda)-\widetilde{\tau}\,c_0(\lambda)| \leq \eps\,|\lambda|^{1/2}\,w(\lambda) \leq 
\eps\,C_1\,|c_0(\lambda)| $$
we obtain
\begin{equation}
\label{eq:f-zero:uniq:wtf}
|\widetilde{c}(\lambda)| \geq \left( |\widetilde{\tau}| - \eps\,C_1 \right) \cdot |c_0(\lambda)| \; . 
\end{equation}
From the estimates \eqref{eq:f-zero:uniq:f} and \eqref{eq:f-zero:uniq:wtf} we see that for \,$\lambda\in V_\delta$\, with \,$|\lambda| \geq R$\,, we have
$$ |\vi(\lambda)| = \frac{|c(\lambda)|}{|\widetilde{c}(\lambda)|} \leq \frac{|\tau| + \eps\,C_1}{|\widetilde{\tau}| - \eps\,C_1}=:C_2 \; . $$
If we now choose \,$\eps>0$\, with \,$|\widetilde{\tau}|-\eps\,C_1>0$\,, it follows that \,$\vi$\, is bounded on \,$V_\delta$\, near \,$\lambda\to\infty$\,. 
If \,$k\in \N$\, is so large that \,$|\lambda|\geq R$\, holds for all \,$\lambda\in U_{k,\delta}$\,, then we have by the maximum principle for holomorphic functions for every \,$\lambda \in U_{k,\delta}$\,
$$ |\vi(\lambda)| \leq \max_{\lambda' \in \partial U_{k,\delta}} |\vi(\lambda')| \leq C_2 \; . $$
Hence \,$\vi$\, is bounded on a neighborhood of \,$\lambda=\infty$\,, and thus can be extended holomorphically in \,$\lambda=\infty$\,. 

An analogous argument shows that \,$\vi$\, can be extended holomorphically also in \,$\lambda=0$\,. Thus \,$\vi$\, is a holomorphic function on the Riemann sphere \,$\PP^1(\C)$\,, and hence constant.
In other words, there exists \,$A\in\C^*$\, with \,$\wt{c}=A\cdot c$\,. 

Let us now consider the sequence of real numbers \,$(\xi_k)_{k\in\N}$\, with \,$|\xi_k|>1$\,, \,$\zeta(\xi_k)=(2k+\tfrac12)\pi$\,. Then we have \,$\lim \xi_k = \infty$\,, \,$\sin(\zeta(\xi_k))=1$\,,
\,$\cos(\zeta(\xi_k))=0$\,, and thus \,$c_0(\xi_k)=\xi_k^{1/2}$\, and \,$w(\xi_k)=1$\,. Therefore there exists for any \,$\eps>0$\, an \,$R>0$\, such that 
$$ |c(\xi_k)-\tau\,\xi_k^{1/2}| \leq \eps\,|\xi_k|^{1/2} \qmq{and} |A\,c(\xi_k)-\wt{\tau}\,\xi_k^{1/2}| \leq \eps\,|\xi_k|^{1/2} $$
holds. From the second inequality, we get \,$|c(\xi_k)-\tfrac{\wt{\tau}}{A}\,\xi_k^{1/2}| \leq \tfrac{\eps}{|A|}\,|\xi_k|^{1/2}$\,, hence we obtain from the first inequality
$$ \left| \left( \tau - \frac{\wt{\tau}}{A} \right) \xi_k^{1/2} \right| \leq \eps\,\left(1+\frac{1}{|A|}\right)\,|\xi_k|^{1/2} \;, $$
and therefore \,$\left| \tau - \tfrac{\wt{\tau}}{A} \right| \leq \eps\,(1+\tfrac{1}{|A|})$\, follows. Because the latter inequality holds true for every \,$\eps>0$\,, we have
\,$\left| \tau - \tfrac{\wt{\tau}}{A} \right|=0$\, and thus \,$\wt{\tau} = A\cdot \tau$\,. The analogous argument for a sequence \,$(\xi_k)$\, with \,$\xi_k\to 0$\, gives \,$\wt{\tau}^{-1}=A\cdot \tau^{-1}$\,. 
Therefrom we conclude \,$A\in \{\pm 1\}$\,, and the claim follows.

\emph{For (2).}
We first note that if we denote the sequence of the zeros of \,$c$\, by \,$(\lambda_k)_{k\in \Z}$\, as in Proposition~\ref{P:excl:basic}(2)
and put \,$\mu_k := a(\lambda_k) = \wt{a}(\lambda_k)$\,, then we have by Proposition~\ref{P:excl:basic}(3) (applied to both the function \,$a$\,
and the function \,$\wt{a}$\,)
$$ \lim_{k\to\infty} (\mu_k-\upsilon\,\mu_{k,0}) = \lim_{k\to\infty} (\mu_k-\wt{\upsilon}\,\mu_{k,0}) =0 $$
and therefore 
$$ \lim_{k\to\infty} \bigr( (\wt{\upsilon}-\upsilon)\cdot \mu_{k,0} \bigr) = 0 \;. $$
Because of \,$\mu_{k,0}=(-1)^k$\,, \,$\upsilon=\wt{\upsilon}$\, follows.

Now consider \,$\psi(\lambda) := \lambda\cdot \tfrac{a(\lambda)-\wt{a}(\lambda)}{c(\lambda)}$\,. 
Because of the hypothesis \,$\ord_\lambda(a-\wt{a})\geq \ord_\lambda(c)$\, for all \,$\lambda\in \C^*$\, with \,$c(\lambda)=0$\,, \,$\psi$\, is holomorphic on all of \,$\C^*$\,, even at the zeros of \,$c$\,.

We again fix \,$\eps,\delta>0$\,. 
Then there exists \,$R>0$\, so that for every \,$\lambda\in\C^*$\, with \,$|\lambda|\geq R$\, or \,$|\lambda|\leq\tfrac{1}{R}$\, we have
$$ |a(\lambda)-\upsilon\,a_0(\lambda)|\leq \eps\,w(\lambda) \qmq{and} |\wt{a}(\lambda)-\upsilon\,a_0(\lambda)|\leq \eps\,w(\lambda) $$
and therefore
$$ |a(\lambda)-\wt{a}(\lambda)|\leq 2\eps\,w(\lambda) \; . $$
If we additionally have \,$\lambda\in V_\delta$\, and let \,$C_1>0$\, be the constant from Proposition~\ref{P:vac2:excldom-M0}(3), then we deduce
$$ |a(\lambda)-\wt{a}(\lambda)|\leq 2\eps\,C_1\,|\lambda|^{-1/2}\,|c_0(\lambda)| \; . $$
We then also have by the asymptotic for \,$c$\, (after suitably enlarging \,$R$\, if necessary)
$$ |c(\lambda)-\tau\,c_0(\lambda)| \leq \eps\,|\lambda|^{1/2}\,w(\lambda) \leq \eps\,C_1\,|c_0(\lambda)| $$
and therefore 
$$ |c(\lambda)| \geq (|\tau|-\eps\,C_1)\,|c_0(\lambda)| \; . $$
Thus we obtain on \,$V_\delta\cap \{|\lambda|\geq R \text{ or } |\lambda|\leq\tfrac{1}{R}\}$\,:
\begin{align*}
|\psi(\lambda)| & = \left| \lambda\cdot \frac{a(\lambda)-\wt{a}(\lambda)}{c(\lambda)} \right| 
\leq |\lambda|\cdot \frac{2\eps\,C_1\,|\lambda|^{-1/2}\,|c_0(\lambda)|}{(|\tau|-\eps\,C_1)\,|c_0(\lambda)|} = \frac{2\eps\, C_1}{|\tau|-\eps\,C_1}\,|\lambda|^{1/2} \; . 
\end{align*}
By applying the maximum principle to \,$\psi$\, on the excluded domains contained in \,$\{|\lambda|\geq R \text{ or } |\lambda|\leq\tfrac{1}{R}\}$\, and then on \,$\{\tfrac{1}{R}\leq|\lambda|\leq R\}$\,, it
follows that there exists a constant \,$C_2>0$\, such that 
$$ |\psi(\lambda)| \leq C_2\cdot |\lambda|^{1/2} \qmq{holds for all \,$\lambda\in\C^*$\,.} $$
From this estimate, it follows first that \,$\psi$\, can be extended holomorphically in \,$\lambda=0$\, by setting \,$\psi(0)=0$\,, and then that the holomorphic function \,$\psi$\, on the complex plane \,$\C$\, is constant. 
We therefore have
\,$\psi=0$\,, whence \,$a=\wt{a}$\, follows.

\emph{For (3).}
The proof is analogous to that of (1). We abbreviate \,$f_0(\lambda) := \Delta_0(\lambda)^2-4 = -4\,\sin(\zeta(\lambda))^2$\,. 

Because \,$f$\, and \,$\widetilde{f}$\, have the same zeros, \,$\vi := f/\wt{f}$\, is holomorphic and non-zero on all of \,$\C^*$\,. 
For fixed \,$\eps,\delta>0$\,, there exists \,$R>0$\, so that for 
\,$\lambda\in V_\delta$\, with \,$|\lambda|\geq R$\, we have
$$ |f(\lambda)-\tau_+\,f_0(\lambda)| \leq \eps\,w(\lambda)^2 \leq \eps\,C_1^2\,|\sin(\zeta(\lambda))|^2 = \tfrac{1}{4}\,\eps\,C_1^2\,|f_0(\lambda)| \;, $$
where \,$C_1>1$\, again is the constant from Proposition~\ref{P:vac2:excldom-M0}(3),
and hence
\begin{equation}
\label{eq:f-zero:uniq:delta-f}
|f(\lambda)| \leq \left( |\tau_+| + \tfrac{1}{4}\,\eps\,C_1^2 \right) \cdot |f_0(\lambda)| \; . 
\end{equation}
From the analogous inequality for \,$\wt{f}$\,,
$$ |\widetilde{f}(\lambda)-\widetilde{\tau}_+\,f_0(\lambda)| \leq 
\tfrac{1}{4}\,\eps\,C_1^2\,|f_0(\lambda)| \,, $$
we obtain
\begin{equation}
\label{eq:f-zero:uniq:delta-wtf}
|\widetilde{f}(\lambda)| \geq \left( |\widetilde{\tau}_+| - \frac{1}{4}\,\eps\,C_1^2 \right) \cdot |f_0(\lambda)| \; . 
\end{equation}
From the estimates \eqref{eq:f-zero:uniq:delta-f} and \eqref{eq:f-zero:uniq:delta-wtf} we see that for \,$\lambda\in V_\delta$\, with \,$|\lambda| \geq R$\,, we have
$$ |\vi(\lambda)| = \frac{|f(\lambda)|}{|\widetilde{f}(\lambda)|} \leq \frac{|\tau_+| + \tfrac14\,\eps\,C_1^2}{|\widetilde{\tau}_+| - \tfrac14\,\eps\,C_1^2}=:C_3 \; . $$
If we now choose \,$\eps>0$\, with \,$|\widetilde{\tau}_+|-\tfrac14\,\eps\,C_1^2>0$\,, it follows that \,$\vi$\, is bounded on \,$V_\delta$\, near \,$\lambda\to\infty$\,. 
If \,$k\in \N$\, is so large that \,$|\lambda|\geq R$\, holds for all \,$\lambda\in U_{k,\delta}$\,, then we have by the maximum principle for holomorphic functions for every \,$\lambda \in U_{k,\delta}$\,
$$ |\vi(\lambda)| \leq \max_{\lambda' \in \partial U_{k,\delta}} |\vi(\lambda')| \leq C_3 \; . $$
Hence \,$\vi$\, is bounded on \,$\C^*$\, near \,$\lambda=\infty$\,, and thus can be extended holomorphically in \,$\lambda=\infty$\,. 

An analogous argument shows that \,$\vi$\, can be extended holomorphically also in \,$\lambda=0$\,. Thus \,$\vi$\, is a holomorphic function on the Riemann sphere \,$\PP^1(\C)$\,, and hence constant.
That is, there exists \,$A\in\C^*$\, with \,$f=A\cdot \wt{f}$\,. 

Let us now consider the sequence of real numbers \,$(\xi_k)_{k\in\N}$\, with \,$|\xi_k|>1$\,, \,$\zeta(\xi_k)=(2k+\tfrac12)\pi$\,. Then we have \,$\lim \xi_k = \infty$\,, \,$\sin(\zeta(\xi_k))=1$\,,
\,$\cos(\zeta(\xi_k))=0$\, and \,$w(\xi_k)=1$\,. For any \,$\eps>0$\, there thus exists \,$R>0$\, so that 
$$ |A\cdot\wt{f}(\xi_k)+4\,\tau_+| \leq \eps \qmq{and} |\wt{f}(\xi_k)+4\,\wt{\tau}_+| \leq \eps $$
holds. From the second inequality, we get \,$|A\cdot \wt{f}(\xi_k)+4\,A\,\wt{\tau}_+| \leq \eps\cdot|A|$\,, hence
$$ \left| \tau_+ - A\cdot \wt{\tau}_+ \right| \leq (1+|A|) \cdot \eps \;, $$
whence \,$\tau_+ = A\cdot \wt{\tau}_+$\, follows. The analogous argument also gives \,$\tau_- = A\cdot \wt{\tau}_-$\,. 
\end{proof}

\section{The Fourier asymptotic of the monodromy}
\label{Se:fasymp}

Beyond the ``basic'' asymptotic of the monodromy described in the Section~\ref{Se:asymp}, which applies to all \,$\lambda$\, for which
\,$|\lambda|$\, is sufficiently large resp.~small, we will also need another type of asymptotic estimate that specifically relates
\,$(M(\lambda_{k,0}))_{k\in\Z}$\, to certain Fourier coefficients. In particular, that series is \,$\ell^2$-summable. Because 
of the relation to Fourier coefficients, we will call this type of asymptotic ``Fourier asymptotic''. Via the Fourier asymptotic, 
we will prove a refinement of the basic asymptotic description of the spectral data from Section~\ref{Se:excl}.

In the treatment of the 1-dimensional Schr\"odinger equation, our Fourier asymptotic as described in the following Theorem, is analogous to 
the refined estimate of \cite{Poeschel-Trubowitz:1987}, Theorem~1.4, p.~16.

\enlargethispage{3em}

\begin{thm}
\label{T:fasymp:fourier}
Let \,$(u,u_y)\in\Pot_{np}$\, be given and put \,$\tau := e^{-(u(0)+u(1))/4}$\, and \,$\upsilon := e^{(u(1)-u(0))/4}$\,. 
We denote by \,$a_k$\, resp.~\,$b_k$\, the cosine resp.~the sine Fourier coefficients of \,$u_z$\,, and by \,$\wt{a}_k,\wt{b}_k$\, the Fourier
coefficients of \,$-u_{\overline z}$\,:
\begin{align}
a_k & := \int_0^1 u_z(t)\,\cos(2\pi k t)\,\mathrm{d}t & \wt{a}_k := -\int_0^1 u_{\overline z}(t)\,\cos(2\pi k t)\,\mathrm{d}t \notag \\
\label{eq:fasymp:fourier:akbk}
b_k & := \int_0^1 u_z(t)\,\sin(2\pi k t)\,\mathrm{d}t & \wt{b}_k := -\int_0^1 u_{\overline z}(t)\,\sin(2\pi k t)\,\mathrm{d}t \; .
\end{align}
Let \,$(r_{ij,k})_{k\in \Z}$\, for \,$i,j\in\{1,2\}$\, be the sequences defined by the following equations for \,$k\geq 0$\,:
\begin{align*}
(-1)^k\, M(\lambda_{k,0}) & = \begin{pmatrix} \upsilon & 0 \\ 0 & \upsilon^{-1} \end{pmatrix} + \frac12 \begin{pmatrix} -\upsilon\,a_k & -\lambda_{k,0}^{-1/2}\,\tau^{-1}\,b_k \\ \lambda_{k,0}^{1/2}\,\tau\,b_k & \upsilon^{-1}\,a_k \end{pmatrix} \!\!+\!\! \begin{pmatrix} r_{11,k} & \lambda_{k,0}^{-1/2}\,r_{12,k} \\ \lambda_{k,0}^{1/2}\,r_{21,k} & r_{22,k} \end{pmatrix} \\
(-1)^k\,M(\lambda_{-k,0}) & = \begin{pmatrix} \upsilon^{-1} & 0 \\ 0 & \upsilon \end{pmatrix} + \frac12 \begin{pmatrix} -\upsilon^{-1}\,\wt{a}_k & -\lambda_{-k,0}^{-1/2}\,\tau\,\wt{b}_k \\ \lambda_{-k,0}^{1/2}\,\tau^{-1}\,\wt{b}_k & \upsilon\,\wt{a}_k \end{pmatrix}\\
& \hspace{6cm}  + \begin{pmatrix} r_{11,-k} & \lambda_{-k,0}^{-1/2}\,r_{12,-k} \\ \lambda_{-k,0}^{1/2}\,r_{21,-k} & r_{22,-k} \end{pmatrix} \;.
\end{align*}
For every \,$p>1$\, we then have \,$(r_{ij,k})_{k\in \Z} \in \ell^p(\Z)$\,, and there exists a constant \,$C>0$\, (dependent on \,$p$\,)
with  \,$\|r_{ij}\|_{\ell^p} \leq C\cdot \|(u,u_y)\|_{\Pot}^2$\,.
\end{thm}

The remainder of this section is concerned with the proof of the above theorem. The strategy for the proof is to decompose the given potential \,$(u,u_y)$\, into two summands, one with a higher order of
regularity and one that is small with respect to \,$\|\,\cdot\,\|_{\Pot}$\,. We will derive asymptotics for these two kinds of potentials separately, in Lemmas~\ref{L:fasymp:fourier-more} and \ref{L:fasymp:fourier-small}, 
and finally combine the two to obtain the proof of Theorem~\ref{T:fasymp:fourier}. 

In the proof, we will make  use of the Banach spaces \,$\ell^p(\Z)$\, and \,$\ell^p(\Z^2)$\,. We note that for 
\,$1 \leq p < q \leq \infty$\,, we have \,$\ell^p \subset \ell^q$\, and for any \,$(a_k)\in \ell^p$\,: 
\begin{equation}
\label{eq:fasymp:fourier-small:lplq}
\|a_k\|_{\ell^q} \leq \|a_k\|_{\ell^p} \; . 
\end{equation}

{\footnotesize
\begin{proof}[Proof of Equation~\eqref{eq:fasymp:fourier-small:lplq}.]
Let \,$(a_k)\in \ell^p$\, be given. Without loss of generality we may suppose \,$\|a_k\|_{\ell^p} = 1$\,, then we have in particular \,$|a_k|\leq 1$\, for all \,$k$\,. Therefore we have
\,$|a_k|^q \leq |a_k|^p$\, because of \,$p <q$\,, and hence 
$$ \|a_k\|_{\ell^q}^q = \sum_k |a_k|^q \leq \sum_k |a_k|^p = \|a_k\|_{\ell^p}^p = 1\;, $$
whence \,$\|a_k\|_{\ell^q} \leq 1 = \|a_k\|_{\ell^p}$\, follows.
\end{proof}

}

There are two important
operations for such sequences:

First, the component-wise multiplication \,$(a_k\cdot b_k)$\, of two sequences \,$a_k$\, and \,$b_k$\,. The H\"older inequality states that whenever \,$a_k\in \ell^p$\, and \,$b_k\in \ell^{p'}$\,, where
\,$p,p'\geq 1$\, and \,$\tfrac{1}{p}+\tfrac{1}{p'}=1$\,, we have \,$(a_k\cdot b_k)\in \ell^1$\, and \,$\|a_k\cdot b_k\|_{\ell^1} \leq \|a_k\|_{\ell^p} \cdot \|b_k\|_{\ell^{p'}}$\,. 

The second important operation is the convolution of two sequences \,$(a_k)$\, and \,$(b_k)$\,. The convolution \,$a*b$\, is defined by
\begin{equation}
\label{eq:fasymp:convolution}
(a * b)_k := \sum_{j\in \Z} a_{k-j}\,b_j = (b* a)_k \;,
\end{equation}
for all sequences \,$(a_k)$\, and \,$(b_k)$\, such that the above infinite sum converges for all \,$k$\,. Young's inequality states that for \,$p,q,r \geq 1$\, with \,$\tfrac{1}{r}+1 = \tfrac{1}{p}+\tfrac{1}{q}$\,,
and any sequences \,$(a_k)\in \ell^p$\, and \,$(b_k)\in \ell^q$\,, the convolution \,$a*b$\, is well-defined, we have \,$a*b \in \ell^r$\, and \,$\|a*b\|_{\ell^r} \leq \|a\|_{\ell^p} \cdot \|b\|_{\ell^q}$\,
(\cite{Bennett-Sharpley:1988}, Theorem~4.2.4, p.~199).
Moreover, the convolution with the sequences \,$(\tfrac{1}{k})_{k\in\Z}$\, and \,$(\tfrac{1}{|k|})_{k\in\Z}$\, (which are only weakly \,$\ell^1$\,), where we assign \,$\tfrac{1}{0}$\, an arbitrary value, will be important for us. 
In relation to this, we have for any \,$1<p\leq 2$\, (not for \,$p=1$\,, however!) and any \,$a_k \in \ell^p(\Z)$\, also \,$(\tfrac{1}{k}*a),(\tfrac{1}{|k|}*a)\in \ell^p(\Z)$\, and 
\begin{equation}
\label{eq:fasymp:fourier-small:weakyoung}
\left\|\tfrac{1}{k}*a\right\|_{\ell^p}, \left\|\tfrac{1}{|k|} * a\right\|_{\ell^p} \;\leq\; C_{Y,p}\cdot \|a\|_{\ell^p}
\end{equation}
with a constant \,$C_{Y,p}>0$\, depending only on \,$p$\,.

{\footnotesize
\begin{proof} 
Let \,$T$\, be either of the linear operators \,$(a_k) \mapsto (\tfrac{1}{k} * a_k)$\, or \,$(a_k) \mapsto (\tfrac{1}{|k|} * a_k)$\,. Then \,$T$\, maps \,$\ell^1(\Z)$\, into the Lorentz space
\,$\ell^{1,\infty}(\Z)$\,, and it maps \,$\ell^2(\Z)$\, into \,$\ell^2(\Z)$\,. The reason for the latter inclusion is that up to a factor, \,$\tfrac{1}{k}$\, is the Fourier transform of \,$g(x):=x\in L^2([0,1])$\,. 
For a given \,$(a_k)\in \ell^2(\Z)$\, there exists some \,$f\in L^2([0,1])$\, whose Fourier transform is \,$(a_k)$\,. Then the \,$L^2$-function \,$f\cdot g$\, has the Fourier transform
\,$(\tfrac{1}{k}*a_k) = T((a_k))$\,, and therefore we have \,$T((a_k))\in \ell^2(\Z)$\,. The claim now follows by application of the Marcinkiewicz interpolation theorem (see \cite{Bennett-Sharpley:1988}, Corollary~4.4.14, p.~226)
to the operator \,$T$\,. 
\end{proof}
}



In the sequel, we will consider the Fourier transform for integrable functions \,$f \in L^1([0,1])$\, on the circle resp.~functions \,$g \in L^1([0,1]^2)$\, on the torus. For \,$k\in \Z$\, resp.~for \,$k_1,k_2\in \Z$\,
we define the Fourier transform
\begin{equation}
\label{eq:fasymp:fourier}
\widehat{f}(k) := \int_0^1 f(t)\,e^{-2\pi i k t}\,\mathrm{d}t \qmq{resp.} \widehat{g}(k_1,k_2) := \int_0^1 g(t_1,t_2)\,e^{-2\pi i (k_1t_1+k_2t_2)}\,\mathrm{d}^2 t \; .
\end{equation}
It is important to note that for \,$f\in L^2([0,1])$\,, we have \,$\widehat{f} \in \ell^2(\Z)$\,, and the Fourier map \,$L^2([0,1]) \to \ell^2(\Z)$\, is an isomorphism of Banach (actually Hilbert) spaces, i.e.~it is bijective
and preserves the norm by Plancherel's identity (\cite{Grafakos:2008}, Proposition~3.1.16(1), p.~170): 
\,$\|\widehat{f}\|_{\ell^2(\Z)} = \|f\|_{L^2([0,1])}$\, for any \,$f\in L^2([0,1])$\,. Moreover, the Fourier series \,$\sum_{k\in \Z} \wh{f}(k)\cdot e^{2\pi i m x}$\,
converges to \,$f$\, with respect to the \,$L^2$-norm (\cite{Grafakos:2008}, Proposition~3.1.16(2), p.~170). Analogous statements hold for functions on the torus \,$[0,1]^2$\,.

We also note that the multiplication of functions corresponds under Fourier transformation to the convolution of the corresponding Fourier series: For any \,$f_1,f_2\in L^2([0,1])$\,, we have \,$\widehat{f_1\cdot f_2} = \widehat{f_1}
* \widehat{f_2}$\,. If \,$f\in W^{1,2}([0,1])$\, is periodic, we also have \,$\widehat{f'}(k) = 2\pi i k \, \widehat{f}(k)$\, for \,$k\in \Z$\,. The following Lemma concerns the inverse of the latter operation, namely
the Fourier coefficients of the anti-derivative of a function \,$f$\,. 

\label{not:fasymp:ellp-n}
For this, and in the sequel, we will have use besides the sequence spaces \,$\ell^p(N)$\, (with \,$N\subset \Z$\,, \,$1\leq p \leq \infty$\,)
also for the spaces \,$\ell^p_n(N)$\, (\,$n\in \Z$\,). We say that a sequence \,$(a_k)_{k\in \N}$\, is in \,$\ell^p_n(N)$\, if and only if the sequence \,$(k^n\,a_k)_{k\in \N}$\, is in \,$\ell^p(N)$\,; 
here (and also in associated formulas in the sequel) we put \,$0^n := 1$\, for all \,$n\in \Z$\, to simplify notation. Of course, we equip \,$\ell^p_n(N)$\, with the norm 
$$ \|a_k\|_{\ell^p_n(N)} := \|k^n\,a_k\|_{\ell^p(N)} \; . $$
In this way, \,$\ell^p_n(N)$\, becomes a Banach space (a Hilbert space for \,$p=2$\,). 

\begin{lem}
\label{L:fasymp:integrate}
Let \,$f\in L^1([0,1])$\, be such that \,$\widehat{f} \in \ell^1_{-1}(\Z)$\,. We put \,$F(x) := \int_0^x f(x)\,\mathrm{d}x$\,. Then we have 
\begin{equation}
\label{eq:fasymp:integrate:widehatF}
\widehat{F}(k) = \begin{cases} \frac{\widehat{f}(k)-\widehat{f}(0)}{2\pi i k} & \text{for \,$k\neq 0$\,} \\ \frac12 \widehat{f}(0) - \sum_{j\in \Z\setminus\{0\}} \frac{\widehat{f}(j)}{2\pi i j} & \text{for \,$k=0$\,} \end{cases} \; . 
\end{equation}
\end{lem}

\begin{proof}
Note that \,$F(0)=0$\, and \,$F(1)=\widehat{f}(0)$\, holds.
For \,$k\in \Z\setminus \{0\}$\, we thus have
\begin{align*}
\widehat{F}(k) & = \int_0^1 F(t)\,e^{-2\pi i k t}\,\mathrm{d}t \\
& = \left. -\frac{1}{2\pi i k}\,F(t)\,e^{-2\pi i k t} \right|_{t=0}^{t=1} - \int_0^1 f(t)\,\left(-\frac{1}{2\pi i k}\right)e^{-2\pi i k t}\,\mathrm{d}t \\
& = -\frac{1}{2\pi i k} \widehat{f}(0) + \frac{1}{2\pi i k} \int_0^1 f(t)\,e^{-2\pi i k t}\,\mathrm{d}t \\
& = \frac{\widehat{f}(k)-\widehat{f}(0)}{2\pi i k} \;.
\end{align*}
To calculate \,$\widehat{F}(0)$\,, we consider for \,$j\in \Z$\, the functions \,$f_j(x) := e^{2\pi i j x}$\, and 
$$ F_j(x) := \int_0^x f_j(t)\,\mathrm{d}t = \begin{cases} \frac{e^{2\pi i j x}-1}{2\pi i j} & \text{for \,$j\neq 0$\,} \\ x & \text{for \,$j=0$\,} \end{cases} \; . $$
Then we have for \,$j\neq 0$\,
$$ \widehat{F}_j(0) = \int_0^1 F_j(x)\,\mathrm{d}x = \left. \frac{\tfrac{1}{2\pi i j}e^{2\pi i j x}-x}{2\pi i j} \right|_{x=0}^{x=1} = -\frac{1}{2\pi i j} \;, $$
and for \,$j=0$\,
$$ \widehat{F}_0(0) = \int_0^1 F_0(x)\,\mathrm{d}x = \left. \frac12 x^2 \right|_{x=0}^{x=1} = \frac12 \; . $$
It follows that we have
$$ \widehat{F}(0) = \frac12 \widehat{f}(0) - \sum_{j\in \Z\setminus \{0\}} \frac{\widehat{f}(j)}{2\pi i j} \; . $$
\end{proof}

The following two lemmas each prove a sort of a Fourier asymptotic in a special case with respect to the potential \,$(u,u_y)$\,: Lemma~\ref{L:fasymp:fourier-small} deals with the case where
\,$\|(u,u_y)\|_{\Pot}$\, is small, and Lemma~\ref{L:fasymp:fourier-more} deals with the case where \,$(u,u_y)\in \Pot^1_{np}$\, holds, i.e.~\,$(u,u_y)$\, is once more differentiable. 
Because every \,$(u,u_y)\in \Pot_{np}$\, can be decomposed as a sum of one potential each of these two kinds, we will be able to combine these two lemmas into a full proof of Theorem~\ref{T:fasymp:fourier}.

Both lemmas operate in the setting of the proof of Theorem~\ref{T:asymp:basic}. In particular they concern the regauged and modified
monodromy \,$E$\, introduced in that proof and its series representation \,$E(x) = \sum_{n=0}^\infty E_n(x)$\,, see Equation~\eqref{eq:asymp:E-power-series}. 
In the following Lemma~\ref{L:fasymp:fourier-small}, note that the constant term \,$\unity$\, equals \,$(-1)^k\,E_{0,\lambda_{k,0}}(1)$\, and the Fourier term
\,$\tfrac12 \left( \begin{smallmatrix} -a_k & b_k \\ b_k & a_k \end{smallmatrix} \right)$\, equals \,$(-1)^k\,E_{1,\lambda_{k,0}}(1)$\,. Therefore the main point
of Lemma~\eqref{L:fasymp:fourier-small} is to estimate the sum \,$\sum_{n\geq 2} E_{n,\lambda_{k,0}}(1)$\, of the remaining terms in \,$E_{\lambda_{k,0}}(1)$\,.
It is relatively easy to see that each individual of these summands is in \,$\ell^p(k,\C^{2\times 2})$\,, but it is a remarkable fact that the behavior is not worsened
by taking the infinite sum over all \,$n\geq 2$\,. 

\begin{lem}
\label{L:fasymp:fourier-small}
Suppose \,$1<p\leq 2$\,. 
There exists a constant \,$R_0>0$\, such that for all \,$(u,u_y)\in \Pot_{np}$\, with \,$\|(u,u_y)\|_\Pot \leq R_0$\,, we have
$$ (-1)^k\,E_{\lambda_{k,0}}(1) = \unity + \frac12 \begin{pmatrix} -a_k & b_k \\ b_k & a_k \end{pmatrix} + \ell^p(k,\C^{2\times 2}) \; , $$
and the \,$\ell^p$-norm of the \,$\ell^p$-sequence is \,$\leq C\cdot \|(u,u_y)\|_\Pot^2$\, with some constant \,$C>0$\, (that depends on \,$p$\,). 
Here \,$a_k$\, and \,$b_k$\, are the Fourier coefficients of \,$u_z$\, as in Equation~\eqref{eq:fasymp:fourier:akbk}.
\end{lem}

\begin{proof}
We continue to work in the setting of the proof of Theorem~\ref{T:asymp:basic}. At first we do not restrict the norm of \,$(u,u_y)$\,.
We set \,$\lambda=\lambda_{k,0}$\, with \,$k\geq 1$\,, and omit the subscript \,${}_{\lambda_{k,0}}$\, from \,$E$\, in the sequel. 

The strategy we follow is to prepare the situation for a proof by induction on \,$n$\, that \,$E_n(1) \in \ell^p(k,\C^{2\times 2})$\, holds with a specific estimate for \,$\|E_n(1)\|_{\ell^p(k,\C^{2\times 2})}$\,. 
Via this estimate we will show that \,$\sum_{n\geq 2} E_n(1)$\, is bounded in \,$\ell^p(k,\C^{2\times 2})$\, by a geometric sum. This geometric sum converges (only) if \,$\|(u,u_y)\|_{\Pot}$\, is small,
yielding the result of the lemma.

For \,$n\geq 1$\, we define functions \,$h_n(s_0,s_1)$\, and \,$g_n(s)$\, by
$$ h_1(s_0,s_1) := -\frac12\,u_z(s_1) \;, $$
for \,$n\geq 2$\,
$$ h_n(s_0,s_1) := \left( -\frac{1}{2} \right)^n \int_{s_2=0}^{s_0} \!\!\!\!\!\! u_z(s_1+s_2) \int_{s_3=0}^{s_1} \!\!\!\!\!\! u_z(s_2+s_3) \cdots \int_{s_n=0}^{s_{n-2}} \!\!\!\!\!\!\!\! u_z(s_{n-1}+s_n)\,u_z(s_n)\,\mathrm{d}^{n-1}s \;, $$
and for \,$n\geq 1$\,
$$ g_n(s) := h_n(1-s,s) \; . $$
Then we have by Equations~\eqref{eq:asymp:Gn-def}, \eqref{eq:asymp:H1-def} and \eqref{eq:asymp:Hn-def}
$$ G_n(s) = g_n(s)\,L^n \qmq{and} H_n(s_0,s_1) = h_n(s_0,s_1)\,L^n \; , $$
where we put
$$ L := \begin{pmatrix} 1 & 0 \\ 0 & -1 \end{pmatrix} \qmq{and therefore have} 
L^n = \begin{cases} L & \text{for \,$n$\, odd} \\ \unity & \text{for \,$n$\, even} \end{cases} \; . $$

For \,$n=0$\, we have by Equation~\eqref{eq:asymp:E0}
$$ E_0(1) = \begin{pmatrix} \cos(k\pi\cdot 1) & -\sin(k\pi\cdot 1) \\ \sin(k\pi\cdot 1) & \cos(k\pi\cdot 1) \end{pmatrix} = (-1)^k\cdot \unity \;. $$
Moreover we have \,$\zeta(\lambda)=k\pi$\,, \,$w(\lambda)=1$\,, and therefore for \,$n\geq 1$\,  by Equations~\eqref{eq:asymp:EnGn} and \eqref{eq:asymp:E0}
\begin{align}
E_n(1) & = \int_0^1 E_0(1-2s)\,G_n(s)\,\mathrm{d}s = \int_0^1 \begin{pmatrix} \cos((1-2s)k\pi) & -\sin((1-2s)k\pi) \\ \sin((1-2s)k\pi) & \cos((1-2s)k\pi) \end{pmatrix} \,G_n(s)\,\mathrm{d}s \notag \\
& = (-1)^k\cdot \int_0^1 \begin{pmatrix} \cos(2k\pi s) & \sin(2k\pi s) \\ -\sin(2k\pi s) & \cos(2k\pi s) \end{pmatrix} \,G_n(s)\,\mathrm{d}s \notag \\
\label{eq:fasymp:fourier-small:En}
& = (-1)^k \cdot \int_0^1 \begin{pmatrix} \cos(2k\pi s) & \sin(2k\pi s) \\ -\sin(2k\pi s) & \cos(2k\pi s) \end{pmatrix} \,g_n(s)\,\mathrm{d}s \cdot L^n \;.
\end{align}
Specifically for \,$n=1$\, we have
$$ g_1(s) = h_1(1-s,s) = -\frac12 u_z(s) \;, $$
hence
\begin{align*}
E_1(1) & = -\frac{(-1)^k}{2} \int_0^1 \begin{pmatrix} \cos(2k\pi s) & \sin(2k\pi s) \\ -\sin(2k\pi s) & \cos(2k\pi s) \end{pmatrix} \,u_z(s)\,\mathrm{d}s\cdot L \\
& = -\frac{(-1)^k}{2} \begin{pmatrix} a_k & b_k \\ -b_k & a_k \end{pmatrix} \cdot L = \frac{(-1)^k}{2} \begin{pmatrix} -a_k & b_k \\ b_k & a_k \end{pmatrix} \; . 
\end{align*}
Thus it remains to show that for \,$(u,u_y)\in \Pot_{np}$\, with \,$\|(u,u_y)\|_\Pot \leq R_0$\, with some \,$R_0>0$\,,
$$ \sum_{n=2}^\infty E_n(1) \in \ell^p(k) $$
holds, where \,$\|\sum_{n=2}^\infty E_n(1)\|_{\ell^p(k)} \leq C\cdot \|(u,u_y)\|_\Pot^2$\, holds with a constant \,$C>0$\,;
and because of Equations~\eqref{eq:fasymp:fourier-small:En} and \eqref{eq:asymp:cossinexp}, for this it is in turn sufficient
to show that for the discrete Fourier coefficients \,$\widehat{g}_n$\, of \,$g_n$\,
we have
\begin{equation}
\label{eq:fasymp:fourier-small:gn-sum}
\sum_{n=2}^\infty \widehat{g}_n(k) \in \ell^p(k) \qmq{with} \left\|\sum_{n=2}^\infty \widehat{g}_n(k)\right\|_{\ell^p(k)} \leq C \cdot \|(u,u_y)\|_\Pot^2
\end{equation}
with some constant \,$C>0$\,. 

Because of \,$h_n \in L^2([0,1]^2)$\,, we can express the bivariate function \,$h_n$\, on the torus by its Fourier series:
$$ h_n(s_0,s_1) = \sum_{k_0,k_1\in \Z} \widehat{h}_n(k_0,k_1)\,e^{2\pi i (k_0s_0 + k_1s_1)} \; . $$
Thereby we obtain
\begin{align}
\widehat{g}_n(k) & = \int_0^1 g_n(s)\,e^{-2\pi i k s}\,\mathrm{d}s = \int_0^1 h_n(1-s,s)\,e^{-2\pi i k s}\,\mathrm{d}s \notag \\
& = \sum_{k_0,k_1\in \Z} \widehat{h}_n(k_0,k_1) \int_0^1 e^{2\pi i(k_0(1-s)+k_1 s)}\,e^{-2\pi i k s}\,\mathrm{d}s \notag\\
& = \sum_{k_0,k_1\in \Z} \widehat{h}_n(k_0,k_1) \int_0^1 e^{2\pi i (k_1-k_0-k)s}\,\mathrm{d}s \notag \\
\label{eq:fasymp:fourier-small:gnhn}
& = \sum_{k_0,k_1\in \Z} \widehat{h}_n(k_0,k_1) \delta_{k_1,k+k_0} = \sum_{\ell\in \Z} \widehat{h}_n(\ell-k,\ell) \; . 
\end{align}

For the sake of brevity, we put \,$c_k := -\tfrac12 \widehat{u_z}(k)$\, for \,$k\in \Z$\,; then we have
\,$(c_k) \in \ell^2(\Z)$\, and \,$\|c\|_{\ell^2} \leq \|(u,u_y)\|_\Pot$\,. 

We now make the following claim:
\begin{quote}
There exists a constant \,$C_1>0$\, such that for every \,$n\geq 2$\,, we have
\begin{gather}
\label{eq:fasymp:fourier-small:claim-l1}
\forall \,k\in \Z \; : \; \wh{h}_n(\ell-k,\ell) \in \ell^1(\ell) \\
\label{eq:fasymp:fourier-small:claim-lp}
\|\wh{h}_n(\ell-k,\ell)\|_{\ell^1(\ell)} \in \ell^p(k) \\
\label{eq:fasymp:fourier-small:claim-lpnorm}
\left\|\;\|\wh{h}_n(\ell-k,\ell)\|_{\ell^1(\ell)}\;\right\|_{\ell^p(k)} \leq (C_1 \cdot \|c\|_{\ell^2})^n \; . 
\end{gather}
\end{quote}

The claims \eqref{eq:fasymp:fourier-small:claim-l1}, \eqref{eq:fasymp:fourier-small:claim-lp} imply that \,$\wh{h}_n$\, lies in a Bochner space \,$\ell^1\otimes \ell^p$\,.

Before we prove the claim \eqref{eq:fasymp:fourier-small:claim-l1}--\eqref{eq:fasymp:fourier-small:claim-lpnorm}, we would like to note that the statement \eqref{eq:fasymp:fourier-small:gn-sum} which we need to show
follows from \eqref{eq:fasymp:fourier-small:claim-l1}--\eqref{eq:fasymp:fourier-small:claim-lpnorm}. Indeed, we then have for \,$n\geq 2$\,
$$ |g_n(k)| \overset{\eqref{eq:fasymp:fourier-small:gnhn}}{\leq} \|\wh{h}_n(\ell-k,\ell)\|_{\ell^1(\ell)} $$
and therefore by \eqref{eq:fasymp:fourier-small:claim-lpnorm}
$$ \|g_n(k)\|_{\ell^p(k)} \leq (C_1\cdot \|c\|_{\ell^2})^n \; . $$
If we now put \,$R_0 := \tfrac{1}{2C_1}$\, and \,$C := 2\,C_1^2$\,, and consider \,$(u,u_y)\in \Pot$\, with \,$\|c\|_{\ell^2} \leq \|(u,u_y)\|_\Pot \leq R_0$\,, then we thus see
$$ \sum_{n=2}^\infty \|g_n(k)\|_{\ell^p(k)} \leq \frac{(C_1\cdot \|c\|_{\ell^2})^2}{1-C_1\cdot \|c\|_{\ell^2}} \leq 2\,(C_1\cdot \|c\|_{\ell^2})^2 \leq C\cdot \|(u,u_y)\|_\Pot^2 \; . $$
Therefore \eqref{eq:fasymp:fourier-small:gn-sum} follows from the claim \eqref{eq:fasymp:fourier-small:claim-l1}--\eqref{eq:fasymp:fourier-small:claim-lpnorm}.

We will now prove the claim \eqref{eq:fasymp:fourier-small:claim-l1}--\eqref{eq:fasymp:fourier-small:claim-lpnorm} by induction on \,$n\geq 2$\,. First, let us look at the case \,$n=2$\,. Then we have
$$ h_2(s_0,s_1) = \int_{s_2=0}^{s_0} \vi_2(s_2,s_1)\,\mathrm{d}s_2 \qmq{with} \vi_2(s_0,s_1) = (-\tfrac12 u_z(s_0))\cdot (-\tfrac12 u_z(s_0+s_1)) \; . $$
For \,$(k_0,k_1)\in \Z^2$\, we thus have
\begin{equation}
\label{eq:fasymp:fourier-small:vi2}
\widehat{\vi}_2(k_0,k_1)
= (\delta_{j_1,0}\,c_{j_0} * \delta_{j_0,j_1}\,c_{j_0})_{(k_0,k_1)} = \sum_{j_0,j_1\in \Z} \delta_{j_1,0}\,c_{j_0}\,\delta_{(k_0-j_0),(k_1-j_1)}\,c_{k_0-j_0} 
= c_{k_0-k_1}\,c_{k_1} \; , 
\end{equation}
and therefore by Lemma~\ref{L:fasymp:integrate}
\begin{align}
\wh{h}_2(k_0,k_1) & \;\;=\;\; \begin{cases} \frac{1}{2\pi i k_0}\,\wh{\vi}_2(k_0,k_1) - \frac{1}{2\pi i k_0}\,\wh{\vi}_2(0,k_1) & \text{for \,$k_0 \neq 0$\,} \\ \frac12 \wh{\vi}_2(0,k_1) - \sum_{\ell\neq 0} \frac{\wh{\vi}_2(\ell,k_1)}{2\pi i \ell} & \text{for \,$k_0=0$\,} \end{cases} \notag \\
\label{eq:fasymp:fourier-small:whh2}
& \overset{\eqref{eq:fasymp:fourier-small:vi2}}{=} \begin{cases} \frac{1}{2\pi i k_0}\,c_{k_0-k_1}\,c_{k_1} - \frac{1}{2\pi i k_0}\,c_{-k_1}\,c_{k_1} & \text{for \,$k_0 \neq 0$\,} \\ \frac12 \,c_{-k_1}\,c_{k_1} - c_{k_1}\,\sum_{\ell\neq 0} \frac{c_{\ell-k_1}}{2\pi i \ell} & \text{for \,$k_0=0$\,} \end{cases} \; . 
\end{align}

Therefore we obtain for \,$k\in \Z$\,
\begin{align*}
& \|\wh{h}_2(\ell-k,\ell)\|_{\ell^1(\ell)} = \sum_{\ell\neq k} |\wh{h}_2(\ell-k,\ell)| + |\wh{h}_2(0,k)| \\
\overset{\eqref{eq:fasymp:fourier-small:whh2}}{\leq} & \frac{1}{2\pi} |c_{-k}| \sum_{\ell\neq k} \frac{|c_\ell|}{|\ell-k|} + \frac{1}{2\pi} \sum_{\ell\neq k} \frac{|c_{-\ell}|\cdot|c_\ell|}{|\ell-k|} + \frac12 |c_{-k}|\cdot |c_k| + \frac{1}{2\pi} \,|c_k| \sum_{\ell\neq 0} \frac{|c_{\ell-k}|}{|\ell|} \\
\;\;=\;\; & \frac{1}{2\pi} |c_{-k}| \cdot ( \frac{1}{|\ell|} * |c_\ell|)_k + \frac{1}{2\pi} \,(\frac{1}{|\ell|} * |c_{-\ell}\,c_\ell|)_k + \frac12 |c_{-k}|\cdot |c_k| + \frac{1}{2\pi} |c_k| \cdot (\frac{1}{|\ell|} * |c_{-\ell}|)_k \\
\;\;<\;\; & \infty \; . 
\end{align*}
Thus we have \,$\wh{h}_2(\ell-k,\ell)\in \ell^1(\ell)$\, and moreover, by use of the H\"older inequality and Young's inequality in the variant for weakly \,$\ell^1$-sequences of \eqref{eq:fasymp:fourier-small:weakyoung}:
\begin{align*}
& \left\|\;\|\wh{h}_2(\ell-k,\ell)\|_{\ell^1(\ell)}\;\right\|_{\ell^p(k)} \\
\leq & \frac{1}{2\pi} \,\|c\|_{\ell^2} \, C_{Y,2}\, \|c\|_{\ell^2} + \frac{1}{2\pi} \, C_{Y,p} \, C_{1,p} \, \|c\|_{\ell^2}\,\|c\|_{\ell^2} + \frac{1}{2}\, \|c\|_{\ell^2}\,\|c\|_{\ell^2} + \frac{1}{2\pi}\, \|c\|_{\ell^2} \, C_{Y,2} \, \|c\|_{\ell^2} \\
\leq & (C_1\cdot \|c\|_{\ell^2})^2
\end{align*}
with a suitable \,$C_1>0$\,. This shows that \eqref{eq:fasymp:fourier-small:claim-l1}--\eqref{eq:fasymp:fourier-small:claim-lpnorm} hold for \,$n=2$\,.

We now tackle the induction step. We suppose that \,$n\geq 3$\, is given such that \eqref{eq:fasymp:fourier-small:claim-l1}--\eqref{eq:fasymp:fourier-small:claim-lpnorm} hold for \,$n-1$\,. We then have
$$ h_n(s_0,s_1) = \int_{s_2=0}^{s_0} \vi_n(s_2,s_1)\,\mathrm{d}s_2 
\qmq{with}
\vi_n(s_0,s_1) := -\frac12\,{u_z}(s_0+s_1)\,h_{n-1}(s_1,s_0) \; . $$
Therefore we have
\begin{align*}
\widehat{\vi}_n(k_0,k_1) & = (\delta_{j_0,j_1}\,c_{j_0} * \widehat{h}_{n-1}(j_1,j_0))_{(k_0,k_1)} = \sum_{j_0,j_1\in\Z} \delta_{j_0,j_1}\,c_{j_0}\,\widehat{h}_{n-1}(k_1-j_1,k_0-j_0) \\
& = \sum_{j\in \Z} c_j\,\widehat{h}_{n-1}(k_1-j,k_0-j) \; .
\end{align*}
and hence by Lemma~\ref{L:fasymp:integrate}
{\footnotesize
\begin{align}
\wh{h}_n(k_0,k_1) & = \begin{cases} \frac{1}{2\pi i k_0}\,\wh{\vi}_n(k_0,k_1) - \frac{1}{2\pi i k_0}\,\wh{\vi}_n(0,k_1) & \text{for \,$k_0 \neq 0$\,} \\ \frac12 \wh{\vi}_n(0,k_1) - \sum_{\ell\neq 0} \frac{\wh{\vi}_n(\ell,k_1)}{2\pi i \ell} & \text{for \,$k_0=0$\,} \end{cases} \notag \\
\label{eq:fasymp:fourier-small:whhn}
& = \begin{cases} \frac{1}{2\pi i k_0}\,\sum_{j} c_j\,\wh{h}_{n-1}(k_1-j,k_0-j) - \frac{1}{2\pi i k_0}\,\sum_j c_j\,\wh{h}_{n-1}(k_1-j,-j) & \text{for \,$k_0 \neq 0$\,} \\ \frac12 \,\sum_j c_j \wh{h}_{n-1}(k_1-j,-j)-\tfrac{1}{2\pi i} \sum_{\ell\neq 0,j\in\Z} \frac{1}{\ell}\,c_j\,\wh{h}_{n-1}(k_1-j,\ell-j) & \text{for \,$k_0=0$\,} \end{cases} \; . 
\end{align}
}
We therefore obtain
\begin{align}
& \|\wh{h}_n(\ell-k,\ell)\|_{\ell^1(\ell)} = \sum_{\ell\neq k} |\wh{h}_n(\ell-k,\ell)| + |\wh{h}_n(0,k)| \notag \\
\leq & \frac{1}{2\pi} \underbrace{\sum_{\ell\neq k,j\in\Z} \frac{|c_j|\,|\wh{h}_{n-1}(\ell-j,\ell-k-j)|}{|\ell-k|}}_{(A)}
+ \frac{1}{2\pi} \underbrace{\sum_{\ell\neq k,j\in\Z} \frac{|c_j|\,|\wh{h}_{n-1}(\ell-j,-j)}{|\ell-k|}}_{(B)} \notag \\
& \qquad\qquad\qquad 
+ \frac{1}{2}\,\underbrace{\sum_{j\in\Z} |c_j|\,|\wh{h}_{n-1}(k-j,-j)|}_{(C)}
+ \frac{1}{2\pi} \underbrace{\sum_{\ell\neq 0,j\in \Z} \frac{|c_j|\,|\wh{h}_{n-1}(k-j,\ell-j)|}{|\ell|}}_{(D)}
\label{eq:fasymp:fourier-small:whhn-masteresti}
\end{align}
We estimate the four sums labeled \,$(A)$\,, \,$(B)$\,, \,$(C)$\, and \,$(D)$\, separately. First, we have by Cauchy-Schwarz's inequality and the variant \eqref{eq:fasymp:fourier-small:weakyoung} of Young's inequality 
\begin{align*}
(A) & = \sum_{\ell\neq k,j\in\Z} \frac{|c_j|\,|\wh{h}_{n-1}(\ell-j,\ell-k-j)|}{|\ell-k|} = \sum_{\ell\neq-j,j\in \Z} \frac{|c_j|\,|\wh{h}_{n-1}(\ell+k,\ell)|}{|\ell+j|} \\
& = \sum_{j\in\Z} |c_j|\,\left( \tfrac{1}{|\ell|} * |\wh{h}_{n-1}(-\ell+k,-\ell)|\right)_j \leq \|c\|_{\ell^2} \cdot \left\|\tfrac{1}{|\ell|} * |\wh{h}_{n-1}(\ell+k,\ell)|\right\|_{\ell^2(\ell)} \\
& \leq \|c\|_{\ell^2} \cdot C_{Y,2} \cdot \|\wh{h}_{n-1}(\ell+k,\ell)\|_{\ell^2(\ell)} \overset{\eqref{eq:fasymp:fourier-small:lplq}}{\leq} \|c\|_{\ell^2} \cdot C_{Y,2} \cdot \|\wh{h}_{n-1}(\ell+k,\ell)\|_{\ell^1(\ell)} \; . 
\end{align*}
By the induction hypothesis, it thus follows that \,$(A) < \infty$\,, \,$(A) \in \ell^p(k)$\, and
$$ \|(A)\|_{\ell^p(k)} \leq \|c\|_{\ell^2} \cdot C_{Y,2} \cdot \left\| \|\wh{h}_{n-1}(\ell+k,\ell)\|_{\ell^1(\ell)} \right\|_{\ell^p(k)} \leq (C_1\cdot \|c\|_{\ell^2})^n \;, $$
when \,$C_1>0$\, is chosen with \,$C_1 \geq C_{Y,2}$\,. 

Next, we have by Cauchy-Schwarz's inequality and \eqref{eq:fasymp:fourier-small:lplq}
\begin{align*}
(B) & = \sum_{\ell\neq k,j\in\Z} \frac{|c_j|\,|\wh{h}_{n-1}(\ell-j,-j)}{|\ell-k|} \leq \sum_{\ell\neq k} \frac{1}{|\ell-k|} \|c\|_{\ell^2} \, \|\wh{h}_{n-1}(\ell-j,-j)\|_{\ell^2(j)} \\
& \leq \sum_{\ell\neq k} \frac{1}{|\ell-k|} \|c\|_{\ell^2} \, \|\wh{h}_{n-1}(\ell-j,-j)\|_{\ell^1(j)} = \|c\|_{\ell^2} \cdot (\tfrac{1}{|\ell|} * \|\wh{h}_{n-1}(\ell-j,-j)\|_{\ell^1(j)})_k \; . 
\end{align*}
By the induction  hypothesis, it thus follows that \,$(B) < \infty$\,. Because of the induction hypothesis
\,$\|\wh{h}_{n-1}(\ell-j,-j)\|_{\ell^1(j)} \in \ell^p(\ell)$\,, we also have 
\,$\tfrac{1}{|\ell|} * \|\wh{h}_{n-1}(\ell-j,-j)\|_{\ell^1(j)} \in \ell^p$\, by the variant \eqref{eq:fasymp:fourier-small:weakyoung} of Young's inequality,
and therefore \,$(B) \in \ell^p(k)$\, and
$$ \|(B)\|_{\ell^p(k)} \leq \|c\|_{\ell^2} \cdot C_{Y,p} \cdot \left\| \|\wh{h}_{n-1}(\ell-j,-j)\|_{\ell^1(j)} \right\|_{\ell^p(\ell)} \leq (C_1\cdot \|c\|_{\ell^2})^n \;, $$
when \,$C_1>0$\, is chosen with \,$C_1 \geq C_{Y,p}$\,. 

Third, we estimate
\begin{align*}
(C) & = \sum_{j\in\Z} |c_j|\,|\wh{h}_{n-1}(k-j,-j)| \leq \|c\|_{\ell^2} \cdot \|\wh{h}_{n-1}(k-j,-j)\|_{\ell^2(j)} \\
& \overset{\eqref{eq:fasymp:fourier-small:lplq}}{\leq} \|c\|_{\ell^2} \cdot \|\wh{h}_{n-1}(k-j,-j)\|_{\ell^1(j)} \; . 
\end{align*}
By the induction hypothesis, it thus follows that \,$(C) < \infty$\,, \,$(C) \in \ell^p(k)$\, and
$$ \|(C)\|_{\ell^p(k)} \leq \|c\|_{\ell^2} \cdot \left\| \|\wh{h}_{n-1}(k-j,-j)\|_{\ell^1(j)} \right\|_{\ell^p(k)} \leq (C_1\cdot \|c\|_{\ell^2})^n \;, $$
when \,$C_1>0$\, is chosen with \,$C_1 \geq 1$\,. 

Finally, we obtain
\begin{align*}
(D) & = \sum_{\ell\neq 0,j\in \Z} \frac{|c_j|\,|\wh{h}_{n-1}(k-j,\ell-j)|}{|\ell|} = \sum_{\ell\neq 0,j\in \Z} \frac{|c_{\ell-j}|\,|\wh{h}_{n-1}(j+(k-\ell),j)|}{|\ell|} \\
& \leq \sum_{\ell\neq 0} \frac{1}{|\ell|}\,\|c\|_{\ell^2} \, \|\wh{h}_{n-1}(j+(k-\ell),j)\|_{\ell^2(j)} \overset{\eqref{eq:fasymp:fourier-small:lplq}}{\leq} \|c\|_{\ell^2} \cdot \sum_{\ell\neq 0} \frac{1}{|\ell|}\,\|\wh{h}_{n-1}(j+(k-\ell),j)\|_{\ell^1(j)} \\
& = \|c\|_{\ell^2} \cdot \left( \tfrac{1}{|\ell|} *_\ell \|\wh{h}_{n-1}(j+\ell,j)\|_{\ell^1(j)} \right)_k \;. 
\end{align*}
By the induction  hypothesis, it thus follows that \,$(D) < \infty$\,. Because of the induction hypothesis
\,$\|\wh{h}_{n-1}(j+\ell,j)\|_{\ell^1(j)} \in \ell^p(\ell)$\,, we also have 
\,$\tfrac{1}{|\ell|} *_\ell \|\wh{h}_{n-1}(j+\ell,j)\|_{\ell^1(j)} \in \ell^p$\, by the variant \eqref{eq:fasymp:fourier-small:weakyoung} of Young's inequality,
and therefore \,$(D) \in \ell^p(k)$\, and
$$ \|(D)\|_{\ell^p(k)} \leq \|c\|_{\ell^2} \cdot C_{Y,p} \cdot \left\| \|\wh{h}_{n-1}(j+\ell,j)\|_{\ell^1(j)} \right\|_{\ell^p(\ell)} \leq (C_1\cdot \|c\|_{\ell^2})^n \;, $$
when \,$C_1>0$\, is chosen with \,$C_1 \geq C_{Y,p}$\,. 

By applying these estimates for (A), (B), (C) and (D) to \eqref{eq:fasymp:fourier-small:whhn-masteresti}, we see
that the claim \eqref{eq:fasymp:fourier-small:claim-l1}--\eqref{eq:fasymp:fourier-small:claim-lpnorm} indeed holds for \,$n$\,, where we choose \,$C_1 := \max\{1,C_{Y,2},C_{Y,p}\}>0$\, (independently of \,$n$\,).
\end{proof}

\begin{lem}
\label{L:fasymp:fourier-more}
For every \,$(u,u_y)\in \Pot^1_{np}$\, there exist constants \,$\sigma_u \in \C$\, and \,$\rho_u>0$\, so that
for \,$k\geq 1$\,, the matrices \,$A_k \in \C^{2\times 2}$\, defined by the equation
$$ (-1)^k \, E_{\lambda_{k,0}}(1) = \unity + \frac12 \begin{pmatrix} -a_k & b_k \\ b_k & a_k \end{pmatrix} + \frac{\sigma_u}{k}\begin{pmatrix} 0 & -1 \\ 1 & 0 \end{pmatrix} + A_k $$
satisfy \,$|A_k| \leq \frac{\rho_u}{k^2}$\,. 
Here both \,$\sigma_u$\, and \,$\rho_u$\, are bounded on any closed ball in \,$\Pot^1_{np}$\,. 
\end{lem}

\begin{proof}
We work in the situation of the proof of Proposition~\ref{P:asymp:more}. When we fix \,$\lambda=\lambda_{k,0}$\, and \,$x=1$\,, we have
$$ \overline{E}_0(1) = (-1)^k\cdot \unity $$
and
\begin{align*}
\overline{E}_1(1) & = \overline{E}_0(1) \cdot \int_0^1 \overline{E}_0(-t) \cdot \overline{\beta}(t) \cdot \overline{E}_0(t) \,\mathrm{d}t \\
& = (-1)^k \,\lambda_{k,0}^{-1/2}\,\left( \int_0^1 \overline{\beta}_+(t)\,\mathrm{d}t + \int_0^1 \overline{E}(2t)\,\overline{\beta}_-(t)\,\mathrm{d}t \right) \\
& = (-1)^k \,\lambda_{k,0}^{-1/2}\,\biggr( \int_0^1 (-\tfrac12\,u_z^2 + \tfrac14\,\cosh(u)-1)\,\mathrm{d}t \, \begin{pmatrix} 0 & 1 \\ -1 & 0 \end{pmatrix} \\
& \quad\qquad\qquad\qquad + \int_0^1 \begin{pmatrix} \cos(2\pi k t) & -\sin(2\pi k t) \\ \sin(2\pi k t) & \cos(2\pi k t) \end{pmatrix} \,\begin{pmatrix} 0 & 1 \\ 1 & 0 \end{pmatrix} \,(-u_{zx}-\tfrac14 \sinh(u))\,\mathrm{d}t \biggr) \\
& = \frac{(-1)^k}{4\pi k}\,\left( \overline{\sigma}_u\, \begin{pmatrix} 0 & 1 \\ -1 & 0 \end{pmatrix} - \int_0^1 \begin{pmatrix} -\sin(2\pi k t) & \cos(2\pi k t) \\ \cos(2\pi k t) & \sin(2\pi k t) \end{pmatrix} \,u_{zx}(t)\,\mathrm{d}t \right) + O(k^{-3})
\end{align*}
with \,$\overline{\sigma}_u := \int_0^1 (-\tfrac12\,u_z^2 + \tfrac14\,\cosh(u)-1)\,\mathrm{d}t$\,; note that \,$\sinh(u)\in W^{2,2}([0,1])$\, holds, and therefore the Fourier coefficients of this function are \,$O(k^{-2})$\,; see
also Equation~\eqref{eq:vacuum:lambdak0-asymp}. 
\end{proof}

\begin{proof}[Proof of Theorem~\ref{T:fasymp:fourier}.]
%
We first consider the case \,$k>0$\,, and we regauge \,$\alpha$\, as in the proof of Theorem~\ref{T:asymp:basic}. 
By the analogous argument as in the proof of Theorem~\ref{T:asymp:basic} we see that the contribution of \,$\gamma$\,
to \,$\wt{F}$\, is of order \,$O(|\lambda_{k,0}|^{-1/2}) = O(\tfrac{1}{k}) \in \ell^p(k)$\,, and therefore can again be neglected.
We have \,$E_0(1)=(-1)^k\,\unity$\, and \,$E_1(1) = \tfrac{(-1)^k}{2}\left( \begin{smallmatrix} -a_k & b_k \\ b_k & a_k \end{smallmatrix} \right)$\,.
In view of the regauging function \,$g$\,, it therefore suffices to show that \,$\sum_{n=2}^\infty E_n(1) \in \ell^p(k,\C^{2\times 2})$\,. 

Because \,$\Pot^1_{np}$\, is dense in \,$\Pot_{np}$\,, there exist \,$(u^{[1]},u^{[1]}_y)\in \Pot^1_{np}$\, and \,$(u^{[2]},u^{[2]}_y)\in \Pot_{np}$\, with \,$\|(u^{[2]},u^{[2]}_y)\|_\Pot \leq R_0$\, (where \,$R_0$\, is the constant 
from Lemma~\ref{L:fasymp:fourier-small}), such that \,$(u,u_y) = (u^{[1]},u^{[1]}_y) + (u^{[2]},u^{[2]}_y)$\, holds. For \,$\nu\in\{1,2\}$\,, we denote the quantities associated to \,$(u^{[\nu]},u^{[\nu]}_y)$\, by the superscript \,${}^{[\nu]}$\,. 
Then by Lemma~\ref{L:fasymp:fourier-more} we have \,$\sum_{n=2}^\infty E_n^{[1]}(1) \in \ell^p(k,\C^{2\times 2})$\,. We also note that by Proposition~\ref{P:asymp:more}, there exists a constant \,$C_1>0$\, such that we have
for all \,$x\in [0,1]$\, 
$$ |E^{[1]}(x)-E_0(x)| \leq \frac{C_1}{\sqrt{|\lambda|}} \;; $$
note that \,$w(\lambda_{k,0})=1$\, holds. We also have \,$|E_0(x)|,|E_0(x)^{-1}|\leq 2$\,, and therefore \,$|E^{[1]}(x)| \leq C_1+2=:C_2$\,. We moreover obtain
\begin{align*}
|(E^{[1]})^{-1}| & \leq |E_0^{-1}| + |(E^{[1]})^{-1} - E_0^{-1}| \leq 2 + |E_0^{-1}|\,|(E^{[1]})^{-1}|\,|E_0-E^{[1]}| \\
& \leq 2 + 2\,C_1\,|\lambda|^{-1/2}\,|(E^{[1]})^{-1}| 
\end{align*}
and therefore
$$ |(E^{[1]})^{-1}| \leq 2\,(1-2\,C_1\,|\lambda|^{-1/2})^{-1} \;. $$
It follows that for \,$|\lambda|$\, sufficiently large, we have \,$|(E^{[1]})^{-1}| \leq C_3$\, for some \,$C_3>0$\,. Then we also have
$$ |(E^{[1]})^{-1}-E_0^{-1}| \leq |E_0^{-1}|\,|(E^{[1]})^{-1}|\,|E_0-E^{[1]}| \leq 2 \cdot C_3 \cdot C_1\,|\lambda|^{-1/2} = C_4\,|\lambda|^{-1/2} $$
with \,$C_4 := 2\,C_1\,C_3$\,. We put \,$A := \max\{2,C_1,C_2,C_3,C_4\}$\,. 

We also note that 
we have \,$\sum_{n=2}^\infty E^{[2]}_n(1) \in \ell^p(k,\C^{2\times 2})$\, and \,$\|\sum_{n=2}^\infty E^{[2]}_n(1)\|_{\ell^p} \leq C \cdot \|(u^{[2]},u_y^{[2]})\|_\Pot$\,  by Lemma~\ref{L:fasymp:fourier-small}.

To compare \,$E_n$\, with \,$E_n^{[1]}$\,, we now apply Lemma~\ref{L:asymp:compare}(1) with \,$\alpha = \wt{\alpha}_0 + \wt{\beta}^{[1]}$\, and \,$\beta = \wt{\beta}^{[2]}$\,.
Therefrom we obtain
\begin{align*}
E_n(1) & = E^{[1]}(1) \cdot \int_{t_1=0}^x \int_{t_2=0}^{t_1} \cdots \int_{t_n=0}^{t_{n-1}} \prod_{j=1}^n E^{[1]}(t_j)^{-1}\,\wt{\beta}^{[2]}(t_j)\,E^{[1]}(t_j)\,\mathrm{d}^nt \\
& = E_0(1) \cdot \int_{t_1=0}^x \int_{t_2=0}^{t_1} \cdots \int_{t_n=0}^{t_{n-1}} \prod_{j=1}^n E_0(t_j)^{-1}\,\wt{\beta}^{[2]}(t_j)\,E_0(t_j)\,\mathrm{d}^nt + D_n(1) \\
& = E^{[2]}(1) + D_n(1)
\end{align*}
with
{\tiny
\begin{align*}
D_n(x) & = (E^{[1]}(x)-E_0(x)) \cdot \int_{t_1=0}^x \cdots \int_{t_n=0}^{t_{n-1}} \prod_{j=1}^n E_0(t_j)^{-1}\,\wt{\beta}^{[2]}(t_j)\,E_0(t_j)\,\mathrm{d}^nt \\
& \qquad + E^{[1]}(x) \cdot \sum_{\ell=1}^n \int_{t_1=0}^x \cdots \int_{t_n=0}^{t_{n-1}} \left(\prod_{j=1}^{\ell-1} E^{[1]}(t_j)^{-1}\,\wt{\beta}^{[2]}(t_j)\,E^{[1]}(t_j)\right) \\
&       \hspace{5cm} \cdot \left( (E^{[1]}(t_\ell)^{-1}-E_0(t_\ell)^{-1})\,\wt{\beta}^{[2]}(t_\ell)\,E_0(t_\ell) + E^{[1]}(t_\ell)^{-1}\,\wt{\beta}^{[2]}(t_\ell)\,(E^{[1]}(t_\ell)-E_0(t_\ell))\right) \\
&       \hspace{5cm} \cdot \left(\prod_{j=\ell+1}^{n} E_0(t_j)^{-1}\,\wt{\beta}^{[2]}(t_j)\,E_0(t_j)\right) \mathrm{d}^nt \; ,
\end{align*}
}
whence it follows
$$ |D_n(x)| \leq (2n+1) \cdot \frac{x^n}{n!} \cdot \|\wt{\beta}^{[2]}\|^n\cdot A^{2n} \cdot A\,|\lambda_{k,0}|^{-1/2} \leq \frac{(\wt{A}\,\|(u^{[2]},u^{[2]}_y)\|_\Pot\,x)^{n}}{n!} \cdot \frac{1}{k} $$
with some \,$\wt{A}>0$\,. If \,$\|(u^{[2]},u^{[2]}_y)\|_\Pot$\, is sufficiently small, we therefore see that 
$$ \sum_{n=2}^\infty |D_n(x)| \leq C \cdot \|(u^{[2]},u^{[2]}_y)\|_\Pot^2 \cdot \frac{1}{k} \in \ell^p(k) $$
holds with some \,$C>0$\,. 

We now simply combine
$$ \sum_{n=2}^\infty |E_n(1)| \leq \sum_{n=2}^\infty |E^{[1]}_n(1)| + \sum_{n=2}^\infty |D_n(1)| \in \ell^p(k) \; , $$
completing the proof of Theorem~\ref{T:fasymp:fourier} for the case \,$k>0$\,. 

Finally, we use Proposition~\ref{P:mono:symmetry}(1) to derive the case \,$k<0$\, from the case \,$k>0$\,. 
We let \,$\wt{u}(z) := -u(\overline{z})$\, and \,$g(\lambda) := \left( \begin{smallmatrix} 1 & 0 \\ 0 & \lambda \end{smallmatrix} \right)$\,.
Then we have \,$\wt{u}_z(z) = -u_{\overline{z}}(\overline{z})$\,, and therefore \,$\wt{a}_k$\, resp.~\,$\wt{b}_k$\, (defined in Equations~\eqref{eq:fasymp:fourier:akbk}) 
are the cosine resp.~the sine Fourier coefficients of \,$\wt{u}_z$\,. 
Also, we have \,$e^{-\wt{u}(0)/2} = e^{u(0)/2}=\tau^{-1}$\,. 
Thus we obtain from Proposition~\ref{P:mono:symmetry}(1) and the case \,$k>0$\, of the present theorem, remembering \,$\lambda_{-k,0}=\lambda_{k,0}^{-1}$\,:
\begin{align*}
(-1)^k\,M_u(\lambda_{-k,0})
& = g(\lambda_{-k,0}^{-1})^{-1} \cdot (-1)^k\,M_{\wt{u}}(\lambda_{-k,0}^{-1}) \cdot g(\lambda_{-k,0}^{-1}) \\
& = g(\lambda_{k,0})^{-1} \cdot (-1)^k\,M_{\wt{u}}(\lambda_{k,0}) \cdot g(\lambda_{k,0}) \\
& = g(\lambda_{k,0})^{-1} \cdot \left( \unity + \frac{1}{2}\begin{pmatrix} -\wt{a}_k & -\lambda_{k,0}^{-1/2}\,\tau\,\wt{b}_k \\ \lambda_{k,0}^{1/2}\,\tau^{-1}\,\wt{b}_k & \wt{a}_k \end{pmatrix} \right. \\
& \hspace{3cm} \left. + \begin{pmatrix} \wt{r}_{11,k} & \lambda_{k,0}^{-1/2}\,\wt{r}_{12,k} \\ \lambda_{k,0}^{1/2}\,\wt{r}_{21,k} & \wt{r}_{22,k} \end{pmatrix} \right) \cdot g(\lambda_{k,0}) \\
& = \unity + \frac{1}{2}\begin{pmatrix} -\wt{a}_k & -\lambda_{k,0}^{1/2}\,\tau\,\wt{b}_k \\ \lambda_{k,0}^{-1/2}\,\tau^{-1}\,\wt{b}_k & \wt{a}_k \end{pmatrix}
+ \begin{pmatrix} \wt{r}_{11,k} & \lambda_{k,0}^{1/2}\,\wt{r}_{12,k} \\ \lambda_{k,0}^{-1/2}\,\wt{r}_{21,k} & \wt{r}_{22,k} \end{pmatrix} \\
& = \unity + \frac{1}{2}\begin{pmatrix} -\wt{a}_k & -\lambda_{-k,0}^{-1/2}\,\tau\,\wt{b}_k \\ \lambda_{-k,0}^{1/2}\,\tau^{-1}\,\wt{b}_k & \wt{a}_k \end{pmatrix}
+ \begin{pmatrix} \wt{r}_{11,k} & \lambda_{-k,0}^{-1/2}\,\wt{r}_{12,k} \\ \lambda_{-k,0}^{1/2}\,\wt{r}_{21,k} & \wt{r}_{22,k} \end{pmatrix} \; .
\end{align*}
By putting \,$r_{ij,-k} := \wt{r}_{ij,k}$\,, the proof of the case \,$k<0$\, is completed.
\end{proof}

\section{The consequences of the Fourier asymptotic for the spectral data}
\label{Se:asympdiv}

We can use the Fourier asymptotics of the previous section to improve our description of the asymptotic behavior of the (classical) spectral divisor \,$D=\{(\lambda_k,\mu_k)\}$\, of a potential \,$(u,u_y)$\,. This is similar to the analogous refinement in the treatment of the 1-dimensional Schr\"odinger equation
found in \cite{Poeschel-Trubowitz:1987}, Theorem~2.4, p.~35.

We would expect a similar refinement for the asymptotic behavior of the branch points \,$\vkap_{k,\nu}$\, of the spectral curve \,$\Sigma$\,.
However, we postpone the investigation of this refinement until Section~\ref{Se:asympfinal}, to avoid technical difficulties.

The following Proposition~\ref{P:asympdiv:asympdiv-neu} 
describes a relationship between the distance between the corresponding points of two spectral divisors, and the function values at \,$\lambda_{k,0}$\, of the functions \,$c$\, and \,$a$\,
in the corresponding monodromies. We will then use this relationship in Corollary~\ref{C:asympdiv:asympdiv-neu} to derive the (improved) asymptotic behavior of the spectral data from the result of Theorem~\ref{T:fasymp:fourier}.

Note that Proposition~\ref{P:asympdiv:asympdiv-neu} gives more information than what would be needed for Corollary~\ref{C:asympdiv:asympdiv-neu}: First, we compare two arbitrary spectral divisors
with each other (instead of one spectral divisor to the divisor of the vacuum), and second, we give in Proposition~\ref{P:asympdiv:asympdiv-neu}(1) two different variants of the estimate,
where only the second one (involving the sequence \,$\wt{\eps}_k$\,) would be needed for Corollary~\ref{C:asympdiv:asympdiv-neu}. We do this to facilitate a further application of Proposition~\ref{P:asympdiv:asympdiv-neu}
in the proof of Lemma~\ref{L:finite:eta} (where the zeros of the function \,$\Delta'(\lambda)$\, are studied). 

We use the notations of the latter part of Section~\ref{Se:excl}.

\begin{prop}
\label{P:asympdiv:asympdiv-neu}
\begin{enumerate}
\item
Let \,$c^{[1]},c^{[2]}: \C^* \to \C$\, be holomorphic functions which satisfy the basic asymptotics of Theorem~\ref{T:asymp:basic}, i.e.~there exist numbers
\,$\tau^{[1]},\tau^{[2]}\in\C^*$\, so that for every \,$\eps>0$\, there exists \,$R>0$\, so that we have
\begin{align*}
|c^{[\nu]}(\lambda)-\tau^{[\nu]}\,c_0(\lambda)| & \leq \eps \,|\lambda|^{1/2}\,w(\lambda) \text{ for \,$|\lambda|\geq R$\,} \\
\qmq{and} |c^{[\nu]}(\lambda)-(\tau^{[\nu]})^{-1}\,c_0(\lambda)| & \leq \eps \,|\lambda|^{1/2}\,w(\lambda) \text{ for \,$|\lambda|\leq \tfrac{1}{R}$\,,} 
\end{align*}
where \,$\nu\in\{1,2\}$\,. Then let \,$(\lambda_k^{[\nu]})_{k\in \Z}$\, be the sequence of zeros of \,$c^{[\nu]}$\, as in Proposition~\ref{P:excl:basic}(2).

In this setting we have for \,$k>0$\,
\begin{align*}
& \left| \bigr(\lambda_k^{[1]}-\lambda_k^{[2]}\bigr) \,-\, 8\,(-1)^{k}\,\bigr( (\tau^{[1]})^{-1}\,c^{[1]}(\lambda_{k,0})-(\tau^{[2]})^{-1}\,c^{[2]}(\lambda_{k,0})\bigr) \right| \\
\leq\;& \eps_k \cdot |\lambda_k^{[1]}-\lambda_k^{[2]}| + |\lambda_k^{[1]}-\lambda_{k,0}| \cdot \frac{C}{k}\cdot \max_{\lambda\in \overline{U_{k,\delta}}} \left| (\tau^{[2]})^{-1}\,c^{[2]}{}(\lambda) - (\tau^{[1]})^{-1}\,c^{[1]}{}(\lambda) \right| \\
\leq\;& \wt{\eps}_k \cdot \max_{\lambda\in \overline{U_{k,\delta}}} \left| (\tau^{[2]})^{-1}\,c^{[2]}{}(\lambda) - (\tau^{[1]})^{-1}\,c^{[1]}{}(\lambda) \right| 
\end{align*}
and
\begin{align*}
& \left| \bigr(\lambda_{-k}^{[1]}-\lambda_{-k}^{[2]}\bigr) \,-\, 8\,(-1)^{k+1}\,\lambda_{-k,0}\,\bigr( \tau^{[1]}\,c^{[1]}(\lambda_{-k,0})-\tau^{[2]}\,c^{[2]}(\lambda_{-k,0})\bigr) \right|  \\
\leq\;& \eps_{-k} \cdot |\lambda_{-k}^{[1]}-\lambda_{-k}^{[2]}| + |\lambda_{-k}^{[1]}-\lambda_{-k,0}| \cdot C\cdot k \cdot \max_{\lambda\in \overline{U_{-k,\delta}}} \left| \tau^{[2]}\,c^{[2]}(\lambda) - \tau^{[1]}\,c^{[1]}(\lambda) \right|  \\
\leq\;& \wt{\eps}_{-k} \cdot \frac{1}{k^2}\cdot \max_{\lambda\in \overline{U_{-k,\delta}}} \left| \tau^{[2]}\,c^{[2]}(\lambda) - \tau^{[1]}\,c^{[1]}(\lambda)\bigr) \right| \;, 
\end{align*}
where \,$(\eps_k)_{k\in\Z}$\,, \,$(\wt{\eps}_k)_{k\in\Z}$\, are sequences converging towards zero for \,$k\to\pm\infty$\,, and \,$C>0$\, is a constant.

\item
Suppose that additionally holomorphic functions \,$a^{[1]},a^{[2]}: \C^*\to\C$\, are given that satisfy the corresponding basic asymptotics of
Theorem~\ref{T:asymp:basic}, i.e.~there exist numbers
\,$\upsilon^{[1]},\upsilon^{[2]}\in\C^*$\, so that for every \,$\eps>0$\, there exists \,$R>0$\, so that we have
\begin{align*}
|a^{[\nu]}(\lambda)-\upsilon^{[\nu]}\,a_0(\lambda)| & \leq \eps \,w(\lambda) \text{ for \,$|\lambda|\geq R$\,} \\
\qmq{and} |a^{[\nu]}(\lambda)-(\upsilon^{[\nu]})^{-1}\,a_0(\lambda)| & \leq \eps \,w(\lambda) \text{ for \,$|\lambda|\leq \tfrac{1}{R}$\,.} 
\end{align*}
Then put \,$\mu_k^{[\nu]} := a^{[\nu]}(\lambda_k^{[\nu]})$\, for \,$k\in \Z$\,, see Proposition~\ref{P:excl:basic}(3).

In this setting there exist \,$C_1,C_2>0$\, so that we have for \,$k>0$\,
\begin{align*}
& \left| \bigr( (\upsilon^{[2]})^{-1}\,\mu_k^{[2]} - (\upsilon^{[1]})^{-1}\,\mu_k^{[1]}\bigr) \;-\; \bigr( (\upsilon^{[2]})^{-1}\,a^{[2]}(\lambda_{k,0}) - (\upsilon^{[1]})^{-1}\,a^{[1]}(\lambda_{k,0})\bigr) \right| \\
\leq\; &  C_1\cdot \frac{|\lambda_k^{[1]}-\lambda_{k,0}|}{k} \cdot \max_{\lambda\in\overline{U_{k,\delta}}} \left| (\upsilon^{[2]})^{-1}\,a^{[2]}(\lambda) - (\upsilon^{[1]})^{-1}\,a^{[1]}(\lambda)\right| \\
& \qquad\qquad + C_2 \cdot \frac{|\lambda_k^{[2]}-\lambda_k^{[1]}|}{k} \cdot \left( 1+\max_{\lambda\in\overline{U_{k,\delta}}} |(\upsilon^{[2]})^{-1}\,a^{[2]}{}(\lambda)-a_0(\lambda)| \right)  
\end{align*}
and
\begin{align*}
& \left| \bigr( \upsilon^{[2]}\,\mu_{-k}^{[2]} - \upsilon^{[1]}\,\mu_{-k}^{[1]}\bigr) \;-\; \bigr( \upsilon^{[2]}\,a^{[2]}(\lambda_{-k,0}) - \upsilon^{[1]}\,a^{[1]}(\lambda_{-k,0})\bigr) \right|  \\
\leq\; &  C_{1}\cdot |\lambda_{-k}^{[1]}-\lambda_{-k,0}|\cdot k^3 \cdot \max_{\lambda\in\overline{U_{-k,\delta}}} \left| \upsilon^{[2]}\,a^{[2]}(\lambda) - \upsilon^{[1]}\,a^{[1]}(\lambda)\right|  \\
& \qquad\qquad + C_{2} \cdot |\lambda_{-k}^{[2]}-\lambda_{-k}^{[1]}| \cdot k^3 \cdot \left( 1+\max_{\lambda\in\overline{U_{-k,\delta}}} |\upsilon^{[2]}\,a^{[2]}(\lambda)-a_0(\lambda)| \right) \; . 
\end{align*}
\end{enumerate}
\end{prop}

\begin{proof}
We first show in both (1) and (2) the estimates for \,$\lambda_k$\, and \,$\mu_k$\, with \,$k$\, positive. 

\emph{For (1).} 
The strategy of the proof is to express the quantity \,$(\tau^{[1]})^{-1}\,c^{[1]}(\lambda_{k,0})-(\tau^{[2]})^{-1}\,c^{[2]}(\lambda_{k,0})$\, in terms of integrals, and then use the given asymptotic behavior
of the \,$c^{[\nu]}$\, (also in the form of Corollary~\ref{C:asymp:Mprime}, which follows by the application of Cauchy's inequality) to relate these integrals to \,$\lambda_k^{[1]}-\lambda_k^{[2]}$\,.

By definition, we have \,$c^{[\nu]}(\lambda_k^{[\nu]}) =0$\, and therefore
\begin{align}
& (\tau^{[1]})^{-1}\,c^{[1]}(\lambda_{k,0})-(\tau^{[2]})^{-1}\,c^{[2]}(\lambda_{k,0}) \notag \\
= \;& -(\tau^{[1]})^{-1}\,\bigr(c^{[1]}(\lambda_k^{[1]})- c^{[1]}(\lambda_{k,0})\bigr)+(\tau^{[2]})^{-1}\,\bigr(c^{[2]}(\lambda_{k}^{[2]})-c^{[2]}(\lambda_{k,0})\bigr) \notag \\
= \; & (\tau^{[2]})^{-1}\,\int_{\lambda_{k,0}}^{\lambda_{k}^{[2]}} c^{[2]}{}'(\lambda)\,\mathrm{d}\lambda -(\tau^{[1]})^{-1}\,\int_{\lambda_{k,0}}^{\lambda_{k}^{[1]}} c^{[1]}{}'(\lambda)\,\mathrm{d}\lambda \notag \\
= \; & \int_{\lambda_k^{[1]}}^{\lambda_k^{[2]}} c_0'(\lambda)\,\mathrm{d}\lambda + \int_{\lambda_{k,0}}^{\lambda_k^{[2]}} \bigr( (\tau^{[2]})^{-1}\,c^{[2]}{}'(\lambda)-c_0'(\lambda)\bigr) \mathrm{d}\lambda \notag \\ 
\label{eq:asympdiv:asympdiv-neu:lambda-master}
& \qquad\qquad + \int_{\lambda_k^{[1]}}^{\lambda_k^{[2]}} \bigr( (\tau^{[2]})^{-1}\,c^{[2]}{}'(\lambda) - c_0'(\lambda) \bigr) \mathrm{d}\lambda \; . 
\end{align}
We now handle the three resulting integrals separately.

By Equation~\eqref{eq:asymp:c0'} we have
\begin{align*}
c_0'(\lambda) & = \frac{\lambda-1}{8\lambda}\,\cos(\zeta(\lambda)) + \frac{1}{2\,\sqrt{\lambda}}\,\sin(\zeta(\lambda)) \\
& = \frac{(-1)^k}{8} + \frac{1}{8}\,\bigr(\cos(\zeta(\lambda))-(-1)^k\bigr) - \frac{1}{8\lambda}\,\cos(\zeta(\lambda)) + \frac{1}{2\,\sqrt{\lambda}}\,\sin(\zeta(\lambda)) \; . 
\end{align*}
It follows that there exists a constant \,$C_1>0$\, so that for all \,$k\geq 1$\, and all \,$\lambda \in U_{k,\delta}$\, we have
$$ \left| c_0'(\lambda)-\frac{(-1)^k}{8} \right| \leq \frac{C_1}{k} \; , $$
and therefore
\begin{equation}
\label{eq:asympdiv:asympdiv-neu:lambda-A}
\left| \int_{\lambda_k^{[1]}}^{\lambda_k^{[2]}} c_0'(\lambda)\,\mathrm{d}\lambda - \frac{(-1)^k}{8}\,(\lambda_k^{[1]}-\lambda_k^{[2]}) \right| \leq \frac{C_1}{k} \cdot |\lambda_k^{[1]}-\lambda_k^{[2]}| \; . 
\end{equation}

Next, there exists a constant \,$C_2>0$\, so that we have
\begin{align}
& \left| \int_{\lambda_{k,0}}^{\lambda_k^{[1]}} \bigr( (\tau^{[2]})^{-1}\,c^{[2]}{}'(\lambda) - (\tau^{[1]})^{-1}\,c^{[1]}{}'(\lambda)\bigr) \mathrm{d}\lambda \right| \notag \\
\leq \; & |\lambda_k^{[1]}-\lambda_{k,0}| \cdot \max_{\lambda\in[\lambda_{k,0},\lambda_k^{[1]}]} \left| (\tau^{[2]})^{-1}\,c^{[2]}{}'(\lambda) - (\tau^{[1]})^{-1}\,c^{[1]}{}'(\lambda) \right| \notag \\
\label{eq:asympdiv:asympdiv-neu:lambda-B}
\leq \; & |\lambda_k^{[1]}-\lambda_{k,0}| \cdot \frac{C_2}{k} \cdot \max_{\lambda\in \overline{U_{k,\delta}}} \left| (\tau^{[2]})^{-1}\,c^{[2]}{}(\lambda) - (\tau^{[1]})^{-1}\,c^{[1]}{}(\lambda) \right|\;, 
\end{align}
where we denote for \,$\lambda_1,\lambda_2 \in \C^*$\, by \,$[\lambda_1,\lambda_2]$\, the straight line from \,$\lambda_1$\, to \,$\lambda_2$\, 
in the complex plane,
and where the second \,$\leq$-sign follows from an application of Cauchy's inequality similar to the one in the proof of Corollary~\ref{C:asymp:Mprime}.

Moreover by Corollary~\ref{C:asymp:Mprime}(1), the sequence
$$ \eps_k^{[1]} := \max_{\lambda\in[\lambda_{k}^{[1]},\lambda_k^{[2]}]} \left| (\tau^{[2]})^{-1}\,c^{[2]}{}'(\lambda) - c_0'(\lambda) \right| $$
converges to zero for \,$k\to\infty$\, (note that \,$w(\lambda)$\, is uniformly bounded on \,$[\lambda_{k}^{[1]},\lambda_k^{[2]}] \subset U_{k,\delta}$\,),
and therefore we have
\begin{align}
\left| \int_{\lambda_k^{[1]}}^{\lambda_k^{[2]}} \bigr( (\tau^{[2]})^{-1}\,c^{[2]}{}'(\lambda) - c_0'(\lambda) \bigr) \mathrm{d}\lambda \right|
& \leq |\lambda_k^{[2]}-\lambda_k^{[1]}| \cdot \max_{\lambda\in[\lambda_{k}^{[1]},\lambda_k^{[2]}]} \left| (\tau^{[2]})^{-1}\,c^{[2]}{}'(\lambda) - c_0'(\lambda) \right| \notag \\
\label{eq:asympdiv:asympdiv-neu:lambda-C}
& \leq \eps_k^{[1]} \cdot |\lambda_k^{[2]}-\lambda_k^{[1]}| \; . 
\end{align}

By applying the estimates \eqref{eq:asympdiv:asympdiv-neu:lambda-A}, \eqref{eq:asympdiv:asympdiv-neu:lambda-B} and \eqref{eq:asympdiv:asympdiv-neu:lambda-C}
to Equation~\eqref{eq:asympdiv:asympdiv-neu:lambda-master}, we obtain
\begin{align}
& \left| \bigr( (\tau^{[1]})^{-1}\,c^{[1]}(\lambda_{k,0})-(\tau^{[2]})^{-1}\,c^{[2]}(\lambda_{k,0}) \bigr) \;-\; \frac{(-1)^k}{8}(\lambda_k^{[1]}-\lambda_k^{[2]}) \right| \notag \\
\label{eq:asympdiv:asympdiv-neu:lambda-master2}
\leq\; & \eps_k^{[2]} \cdot |\lambda_k^{[1]}-\lambda_k^{[2]}| + |\lambda_k^{[1]}-\lambda_{k,0}| \cdot \frac{C_2}{k}\cdot \max_{\lambda\in \overline{U_{k,\delta}}} \left| (\tau^{[2]})^{-1}\,c^{[2]}{}(\lambda) - (\tau^{[1]})^{-1}\,c^{[1]}{}(\lambda) \right|
\end{align}
with the sequence \,$\eps_k^{[2]} := \tfrac{C_1}{k} + \eps_k^{[1]}$\,, which converges to zero for \,$k\to\infty$\,, and therefore the first claimed estimate
for \,$|\lambda_k^{[1]}-\lambda_k^{[2]}|$\,.

For the second claimed estimate we note that 
by Proposition~\ref{P:excl:basic}(2) there exists a sequence \,$(\eps_k^{[3]})$\, which converges to zero for \,$k\to\infty$\, so that
\,$|\lambda_k^{[1]}-\lambda_{k,0}| \leq k\cdot \eps_k^{[3]}$\, holds, 
and therefore it follows from \eqref{eq:asympdiv:asympdiv-neu:lambda-master2} that
\begin{align}
& \left| \bigr( (\tau^{[1]})^{-1}\,c^{[1]}(\lambda_{k,0})-(\tau^{[2]})^{-1}\,c^{[2]}(\lambda_{k,0}) \bigr) \;-\; \frac{(-1)^k}{8}(\lambda_k^{[1]}-\lambda_k^{[2]}) \right| \notag \\
\label{eq:asympdiv:asympdiv-neu:lambda-master3}
\leq\; & \eps_k^{[2]} \cdot |\lambda_k^{[1]}-\lambda_k^{[2]}| + C_2\cdot \eps_k^{[3]}\cdot \max_{\lambda\in \overline{ U_{k,\delta}}} \left| (\tau^{[2]})^{-1}\,c^{[2]}{}(\lambda) - (\tau^{[1]})^{-1}\,c^{[1]}{}(\lambda) \right|
\end{align}
holds.

It follows from \eqref{eq:asympdiv:asympdiv-neu:lambda-master3} that there exists \,$C_3>0$\, so that 
$$ |\lambda_k^{[1]}-\lambda_k^{[2]}| \leq C_3 \cdot \max_{\lambda\in \overline{U_{k,\delta}}} \left| (\tau^{[2]})^{-1}\,c^{[2]}{}(\lambda) - (\tau^{[1]})^{-1}\,c^{[1]}{}(\lambda) \right| $$
holds, and that we therefore have
\begin{align*}
& \left| \bigr( (\tau^{[1]})^{-1}\,c^{[1]}(\lambda_{k,0})-(\tau^{[2]})^{-1}\,c^{[2]}(\lambda_{k,0}) \bigr) \;-\; \frac{(-1)^k}{8}(\lambda_k^{[1]}-\lambda_k^{[2]}) \right| \\
\leq\; & \eps_k^{[4]}\cdot \max_{\lambda\in \overline{U_{k,\delta}}} \left| (\tau^{[2]})^{-1}\,c^{[2]}{}(\lambda) - (\tau^{[1]})^{-1}\,c^{[1]}{}(\lambda) \right|
\end{align*}
with the sequence \,$\eps_k^{[4]} := \eps_k^{[2]}\cdot C_3 + C_2 \cdot \eps_k^{[3]}$\,, which converges to zero for \,$k\to\infty$\,. Hence the
second claimed estimate for \,$|\lambda_k^{[1]}-\lambda_k^{[2]}|$\, follows.

\emph{For (2).}
We have \,$\mu_k^{[\nu]}=a^{[\nu]}(\lambda_k^{[\nu]})$\, by definition and therefore
\begin{align}
& \bigr( (\upsilon^{[2]})^{-1}\,\mu_k^{[2]} - (\upsilon^{[1]})^{-1}\,\mu_k^{[1]}\bigr) \;-\; \bigr( (\upsilon^{[2]})^{-1}\,a^{[2]}(\lambda_{k,0}) - (\upsilon^{[1]})^{-1}\,a^{[1]}(\lambda_{k,0})\bigr) \notag \\
=\; & (\upsilon^{[2]})^{-1}\cdot \bigr(a^{[2]}(\lambda_k^{[2]})- a^{[2]}(\lambda_{k,0})\bigr) - (\upsilon^{[1]})^{-1}\cdot \bigr(a^{[1]}(\lambda_k^{[1]})- a^{[1]}(\lambda_{k,0})\bigr) \notag \\
= \; & (\upsilon^{[2]})^{-1}\cdot \int_{\lambda_{k,0}}^{\lambda_k^{[2]}} a^{[2]}{}'(\lambda)\,\mathrm{d}\lambda - (\upsilon^{[1]})^{-1}\cdot \int_{\lambda_{k,0}}^{\lambda_k^{[1]}} a^{[1]}{}'(\lambda)\,\mathrm{d}\lambda \notag \\
= \; & \int_{\lambda_{k,0}}^{\lambda_k^{[1]}} \bigr( (\upsilon^{[2]})^{-1}\,a^{[2]}{}'(\lambda) - (\upsilon^{[1]})^{-1}\,a^{[1]}{}'(\lambda)\bigr)\,\mathrm{d}\lambda
\label{eq:asympdiv:asympdiv-neu:mu-master}
+ (\upsilon^{[2]})^{-1}\cdot \int_{\lambda_{k}^{[1]}}^{\lambda_k^{[2]}} a^{[2]}{}'(\lambda)\,\mathrm{d}\lambda \; . 
\end{align}
There exist constants \,$C_4,C_5,C_6>0$\, so that we have
\begin{align*}
& \left| \int_{\lambda_{k,0}}^{\lambda_k^{[1]}} \bigr( (\upsilon^{[2]})^{-1}\,a^{[2]}{}'(\lambda) - (\upsilon^{[1]})^{-1}\,a^{[1]}{}'(\lambda)\bigr)\,\mathrm{d}\lambda \right| \\
\leq\; & |\lambda_k^{[1]}-\lambda_{k,0}| \cdot \max_{\lambda\in[\lambda_{k,0},\lambda_k^{[1]}]} \left| (\upsilon^{[2]})^{-1}\,a^{[2]}{}'(\lambda) - (\upsilon^{[1]})^{-1}\,a^{[1]}{}'(\lambda) \right| \\
\leq\; & |\lambda_k^{[1]}-\lambda_{k,0}| \cdot \frac{C_4}{k}\cdot \max_{\lambda\in \overline{U_{k,\delta}}} \left| (\upsilon^{[2]})^{-1}\,a^{[2]}{}(\lambda) - (\upsilon^{[1]})^{-1}\,a^{[1]}{}(\lambda) \right| 
\end{align*}
(where the second inequality follows from Cauchy's inequality), and also
\begin{align*}
\left| (\upsilon^{[2]})^{-1}\cdot \int_{\lambda_{k}^{[1]}}^{\lambda_k^{[2]}} a^{[2]}{}'(\lambda)\,\mathrm{d}\lambda \right| 
& \leq |\lambda_k^{[2]}-\lambda_k^{[1]}| \cdot |\upsilon^{[2]}|^{-1}\cdot \max_{\lambda\in[\lambda_k^{[1]},\lambda_k^{[2]}]} |a^{[2]}{}'(\lambda)| \\
& \leq |\lambda_k^{[2]}-\lambda_k^{[1]}| \cdot \left( \frac{C_5}{k} + \max_{\lambda\in[\lambda_k^{[1]},\lambda_k^{[2]}]} |(\upsilon^{[2]})^{-1}\,a^{[2]}{}'(\lambda)-a_0'(\lambda)| \right) \\
& \leq \frac{|\lambda_k^{[2]}-\lambda_k^{[1]}|}{k} \cdot \left( C_5 + C_6\cdot \max_{\lambda\in\overline{U_{k,\delta}}} |(\upsilon^{[2]})^{-1}\,a^{[2]}{}(\lambda)-a_0(\lambda)| \right)  
\end{align*}
(the last estimate again by Cauchy's inequality). 
By taking the absolute value of Equation~\eqref{eq:asympdiv:asympdiv-neu:mu-master} and then applying the preceding estimates, we obtain the 
result claimed for \,$\mu_k^{[1]}-\mu_k^{[2]}$\, in the Proposition.

\medskip

The estimates for \,$\lambda_{-k}^{[\nu]}$\, and \,$\mu_{-k}^{[\nu]}$\, in (1) and (2) can now be derived by applying the previously shown results
to the holomorphic functions \,$\wt{c}^{[\nu]},\wt{a}^{[\nu]}: \C^* \to \C$\, for \,$\nu\in\{1,2\}$\, given by
\begin{equation*}
\wt{c}^{[\nu]}(\lambda) := \lambda\cdot c^{[\nu]}(\lambda^{-1}) \qmq{and} \wt{a}^{[\nu]}(\lambda) := a^{[\nu]}(\lambda^{-1}) \;. 
\end{equation*}
We denote the quantities associated to \,$\wt{c}^{[\nu]}$\, resp.~to \,$\wt{a}^{[\nu]}$\, by attaching a tilde to the associated symbol.
\,$\wt{c}^{[\nu]}$\, and \,$\wt{a}^{[\nu]}$\, satisfy the asymptotic hypotheses required in (1) resp.~(2) of the Proposition, 
with the constants 
\begin{equation*}
\wt{\tau}^{[\nu]}=(\tau^{[\nu]})^{-1} \qmq{and} \wt{\upsilon}^{[\nu]}=(\upsilon^{[\nu]})^{-1} \; . 
\end{equation*} 
Moreover, we have 
\begin{equation*}
\wt{\lambda}_{k}^{[\nu]} = (\lambda_{-k}^{[\nu]})^{-1} \qmq{and} \wt{\mu}_k^{[\nu]} = \mu_{-k}^{[\nu]} \; .
\end{equation*}
The estimates for \,$\lambda_{-k}^{[\nu]}$\, and \,$\mu_{-k}^{[\nu]}$\, now follow by applying the previous results to \,$\wt{c}^{[\nu]}$\, and \,$\wt{a}^{[\nu]}$\,. 
\end{proof}

\label{not:asympdiv:ellp-nm}
To simplify notations in the sequel, we consider besides the space \,$\ell^p_n$\, also \,$\ell^p_{n,m} := \ell^p_{n,m}(\Z)$\,. We define the corresponding norm for sequences \,$(a_k)_{k\in \Z}$\, by
$$ \|a_k\|_{\ell^p_{n,m}} := \|a_k\|_{\ell^p_n(k>0)} + |a_0| + \|a_{-k}\|_{\ell^p_m(k>0)} \; . $$
Of course, we put \,$\ell^p_{n,m} := \Mengegr{(a_k)_{k\in\Z}}{\|a_k\|_{\ell^p_{n,m}}<\infty}$\,; in this way, \,$\ell^p_{n,m}$\, becomes a Banach space.

\begin{cor}
\label{C:asympdiv:asympdiv-neu}
Let \,$(u,u_y) \in \Pot_{np}$\, (or \,$M(\lambda)$\, a monodromy matrix satisfying the asymptotic properties of Theorems~\ref{T:asymp:basic} and \ref{T:fasymp:fourier}), 
and \,$D$\, the classical spectral divisor of \,$(u,u_y)$\, (or of \,$M(\lambda)$\,). We enumerate \,$D=\{(\lambda_k,\mu_k)\}_{k\in \Z}$\, as in Proposition~\ref{P:excl:basic}(2),(3), then we have
$$ \lambda_k-\lambda_{k,0} \in \ell^2_{-1,3}(k) \qmq{and} \left\{ \begin{matrix} \mu_k-\upsilon\,\mu_{k,0} & \text{if \,$k\geq 0$\,} \\ \mu_k-\upsilon^{-1}\,\mu_{k,0} & \text{if \,$k<0$\,} \end{matrix} \right\} \in \ell^2_{0,0}(k) \; , $$
where \,$\upsilon := e^{(u(1)-u(0))/4} \in \C^*$\,. 
\end{cor}

\begin{proof}
We write the monodromy \,$M(\lambda)$\, of \,$(u,u_y)$\, as \,$M(\lambda)=\left( \begin{smallmatrix} a(\lambda) & b(\lambda) \\ c(\lambda) & d(\lambda) \end{smallmatrix} \right)$\,, and put
\,$\tau := e^{-(u(0)+u(1))/4}$\,. 
Then the hypotheses of 
Proposition~\ref{P:asympdiv:asympdiv-neu} are satisfied for \,$c^{[1]}=c$\,, \,$a^{[1]}=a$\, and \,$c^{[2]}=c_0$\,, \,$a^{[2]}=a_0$\, by Theorem~\ref{T:asymp:basic}.

By Proposition~\ref{P:asympdiv:asympdiv-neu}(1) 
there exists \,$C_1>0$\, and a sequence \,$(\eps_k)_{k\in \Z}$\, with \,$\eps_k \to 0$\, for \,$k\to\pm\infty$\, such that we have for \,$k>0$\,
\begin{align*}
& \left| \bigr(\lambda_k-\lambda_{k,0}\bigr) \,-\, 8\,(-1)^{k}\,\bigr( \tau^{-1}\,c(\lambda_{k,0})-c_0(\lambda_{k,0})\bigr) \right| \\
\leq\;& \eps_k^{[1]} \cdot |\lambda_k-\lambda_{k,0}| + |\lambda_k-\lambda_{k,0}| \cdot \frac{C_1}{k}\cdot \max_{\lambda\in \overline{U_{k,\delta}}} \left| c_0(\lambda) - \tau^{-1}\,c(\lambda) \right| 
\end{align*}
and
\begin{align*}
& \left| \bigr(\lambda_{-k}-\lambda_{-k,0}\bigr) \,-\, 8\,(-1)^{k+1}\,\lambda_{-k,0}\,\bigr( \tau\,c(\lambda_{-k,0})-c_0(\lambda_{-k,0})\bigr) \right|  \\
\leq\;& \eps_{-k}^{[1]} \cdot |\lambda_{-k}-\lambda_{-k,0}| + |\lambda_{-k}-\lambda_{-k,0}| \cdot C_1\cdot k \cdot \max_{\lambda\in \overline{U_{-k,\delta}}} \left| c_0(\lambda) - \tau\,c(\lambda) \right|  \;. 
\end{align*}
For \,$k\to \infty$\,, 
$$ \qmq{both} \frac{C_1}{k}\cdot \max_{\lambda\in \overline{U_{k,\delta}}} \left| c_0(\lambda) - \tau^{-1}\,c(\lambda) \right|  \qmq{and} C_1\cdot k \cdot \max_{\lambda\in \overline{U_{-k,\delta}}} \left| c_0(\lambda) - \tau\,c(\lambda) \right|  $$
converge to zero for \,$k\to\infty$\, by Theorem~\ref{T:asymp:basic}, and therefore it follows that there exists \,$C_2>0$\, so that 
\begin{equation}
\label{eq:asympdiv:asympdiv-neu:C-lambda-pre}
|\lambda_k-\lambda_{k,0}| \leq C_2 \cdot |\tau^{-1}\,c(\lambda_{k,0})-c_0(\lambda_{k,0})| \qmq{and} |\lambda_{-k}-\lambda_{-k,0}| \leq \frac{C_2}{k^2} \cdot |\tau\,c(\lambda_{-k,0})-c_0(\lambda_{-k,0})|
\end{equation}
holds. It follows from Theorem~\ref{T:fasymp:fourier} that
$$ |\tau^{-1}\,c(\lambda_{k,0})-c_0(\lambda_{k,0})| \in \ell^2_{-1}(k>0) \qmq{and} |\tau\,c(\lambda_{-k,0})-c_0(\lambda_{-k,0})| \in \ell^2_1(k>0) $$
holds, and therefore we obtain from \eqref{eq:asympdiv:asympdiv-neu:C-lambda-pre} that
$$ |\lambda_k-\lambda_{k,0}| \in \ell^2_{-1}(k>0) \qmq{and} |\lambda_{-k}-\lambda_{-k,0}| \leq \ell^2_3(k>0) $$
and therefore \,$\lambda_k-\lambda_{k,0} \in \ell^2_{-1,3}(k)$\, holds.

Moreover, by Proposition~\ref{P:asympdiv:asympdiv-neu}(2) there exist \,$C_2,C_3>0$\, so that we have for \,$k>0$\,
\begin{align*}
& \left| \bigr( \mu_{k,0} - \upsilon^{-1}\,\mu_k\bigr) \;-\; \bigr( a_0(\lambda_{k,0}) - \upsilon^{-1}\,a(\lambda_{k,0})\bigr) \right| \\
\leq\; &  C_2\cdot \frac{|\lambda_k-\lambda_{k,0}|}{k} \cdot \max_{\lambda\in\overline{U_{k,\delta}}} \left| a_0(\lambda) - \upsilon^{-1}\,a(\lambda)\right| 
+ C_3 \cdot \frac{|\lambda_k-\lambda_{k,0}|}{k} \cdot 1
\end{align*}
and
\begin{align*}
& \left| \bigr( \mu_{-k,0} - \upsilon\,\mu_{-k}\bigr) \;-\; \bigr( a_0(\lambda_{-k,0}) - \upsilon\,a(\lambda_{-k,0})\bigr) \right|  \\
\leq\; &  C_{2}\cdot |\lambda_{-k}-\lambda_{-k,0}|\cdot k^3 \cdot \max_{\lambda\in\overline{U_{-k,\delta}}} \left| a_0(\lambda) - \upsilon\,a(\lambda)\right|  
+ C_{3} \cdot |\lambda_{-k,0}-\lambda_{-k}| \cdot k^3 \cdot 1 \; . 
\end{align*}
Because we have previously shown \,$\lambda_k-\lambda_{k,0}\in \ell^2_{-1,3}(k)$\,, and for \,$k\to \infty$\,
$$ \qmq{both} \max_{\lambda\in\overline{U_{k,\delta}}} \left| a_0(\lambda) - \upsilon^{-1}\,a(\lambda)\right| \qmq{and} \max_{\lambda\in\overline{U_{-k,\delta}}} \left| a_0(\lambda) - \upsilon\,a(\lambda)\right|  $$
converge to zero by Theorem~\ref{T:asymp:basic}, and we have
$$ \bigr| a_0(\lambda_{k,0}) - \upsilon^{-1}\,a(\lambda_{k,0})\bigr| \in \ell^2_0(k>0) \qmq{and} \bigr| a_0(\lambda_{-k,0}) - \upsilon\,a(\lambda_{-k,0})\bigr| \in \ell^2_0(k>0) $$
by Theorem~\ref{T:fasymp:fourier}, it follows that
$$ \mu_k-\upsilon\,\mu_{k,0} \in \ell^2_0(k>0) \qmq{and} \mu_{-k}-\upsilon^{-1}\,\mu_{-k,0} \in \ell^2_0(k>0) $$
holds.
\end{proof}

\begin{Def}
\label{D:asympdiv:asympdiv}
We call a positive generalized divisor \,$\mathcal{D}$\, or the underlying classical
divisor \,$D$\, on a spectral curve \,$\Sigma \subset \C^*\times \C$\, (or a discrete multi-set \,$D$\, of points in \,$\C^* \times \C$\,) \emph{non-periodic asymptotic}, 
if the support of \,$D$\, is enumerated by a sequence \,$(\lambda_k,\mu_k)_{k\in \Z}$\,, and if there exists a number \,$\upsilon\in\C^*$\,
so that we have
$$ \lambda_k-\lambda_{k,0} \in \ell^2_{-1,3}(k) \qmq{and} \left\{ \begin{matrix} \mu_k-\upsilon\,\mu_{k,0} & \text{if \,$k\geq 0$\,} \\ \mu_k-\upsilon^{-1}\,\mu_{k,0} & \text{if \,$k<0$\,} \end{matrix} \right\} \in \ell^2_{0,0}(k) \; . $$
We call a non-periodic asymptotic divisor \,$\mathcal{D}$\, resp.~\,$D$\, \emph{asymptotic}, 
if the above holds with \,$\upsilon=1$\, i.e.~if we have
$$ \lambda_k-\lambda_{k,0} \in \ell^2_{-1,3}(k) \qmq{and} \mu_k-\mu_{k,0} \in \ell^2_{0,0}(k) \; . $$
We denote the \emph{space of asymptotic classical divisors} (on any spectral curve, regarded as point multi-sets in \,$\C^* \times \C^*$\,) by \,$\Div$\,. 
\end{Def}

\label{not:asympdiv:Div}
Corollary~\ref{C:asympdiv:asympdiv-neu} shows that if \,$(\Sigma,\mathcal{D})$\, are the spectral data
belonging to a non-periodic potential \,$(u,u_y)\in \Pot_{np}$\, (or to a monodromy \,$M(\lambda)$\, satisfying the asymptotic properties of
Theorems~\ref{T:asymp:basic} and \ref{T:fasymp:fourier})
then \,$\mathcal{D}$\, is an non-periodic asymptotic divisor on \,$\Sigma$\,. If \,$(u,u_y) \in \Pot$\, holds, i.e.~if the potential is periodic,
or if \,$\upsilon=1$\, holds in the setting of the present section, then \,$\mathcal{D}$\, is in fact an asymptotic divisor on \,$\Sigma$\,.

In view of Definition~\ref{D:asympdiv:asympdiv} it is tempting to identify the space \,$\Div$\, of asymptotic classical divisors with the Banach space 
\,$\ell^2_{-1,3} \oplus \ell^2_{0,0}$\,. But we need to be careful, because for the points of an asymptotic classical divisor \,$D$\, that lie in the ``compact part'' of \,$\Sigma$\,
(i.e.~for those finitely many points of \,$D$\, which do not need to lie in their excluded domains, see Proposition~\ref{P:excl:basic}(2),(3))
there is no canonical enumeration. Consequently,
the space of asymptotic divisors has a different structure from \,$\ell^2_{-1,3} \oplus \ell^2_{0,0}$\, near those asymptotic divisors \,$D$\,
which contain a point of higher multiplicity (which can occur only for the finitely many divisor points in the ``compact part'' of \,$\Sigma$\,). 

To describe the structure of the space of asymptotic divisors, we consider the group \,$P(\Z)$\, of finite permutations of \,$\Z$\,, i.e.~of
permutations \,$\sigma: \Z \to \Z$\, for which there exists \,$N\in \N$\, (dependent on \,$\sigma$\,) with \,$\sigma(k)=k$\, for all
\,$k\in \Z$\, with \,$|k|>N$\,. \,$P(\Z)$\, acts on \,$\ell^2_{-1,3} \oplus \ell^2_{0,0}$\, by permuting the elements of sequences, and the
quotient space \,$(\ell^2_{-1,3} \oplus \ell^2_{0,0})/P(\Z)$\, is isomorphic to the space \,$\Div$\, of asymptotic divisors. 

Letting \,$D^{[1]},D^{[2]}\in \Div$\, be two asymptotic divisors, 
represented in the usual form \,$D^{[\nu]}=\{(\lambda_k^{[\nu]},\mu_k^{[\nu]})\}_{k\in \Z}$\, for \,$\nu\in\{1,2\}$\,,
we define a distance between them by
$$ \|D^{[1]}-D^{[2]}\|_{\Div} := \inf_{\sigma_1,\sigma_2 \in P(\Z)} \left( \left\|\lambda_{\sigma_1(k)}^{[1]}-\lambda_{\sigma_2(k)}^{[2]}\right\|_{\ell^2_{-1,3}}^2
+ \left\|\mu_{\sigma_1(k)}^{[1]}-\mu_{\sigma_2(k)}^{[2]} \right\|_{\ell^2_{0,0}}^2 \right)^{1/2} \; . $$
Then the topology induced by this distance is the quotient topology of \,$(\ell^2_{-1,3} \oplus \ell^2_{0,0})/P(\Z)$\,; in the sequel we will regard \,$\Div$\, with this topology and this distance function. 

Near any divisor \,$D \in \Div$\, that does not contain any points of higher multiplicity, the space \,$\Div$\, is locally isomorphic to 
\,$\ell^2_{-1,3} \oplus \ell^2_{0,0}$\, (with the product Hilbert norm).
However, if \,$D$\, does contain points of higher multiplicity, this is no longer the case.
Coordinates of \,$\Div$\, near such a point are obtained in the following way: Suppose \,$D=\{(\lambda_k,\mu_k)\}_{k\in \Z}$\,, where
points \,$(\lambda_k,\mu_k)$\, occur more than once according to multiplicity. Choose \,$N\in \N$\, such that \,$(\lambda_k,\mu_k) \in \wh{U}_{k,\delta}$\,
holds for all \,$k$\, with \,$|k|>N$\, (in particular, none of these \,$(\lambda_k,\mu_k)$\, have multiplicity \,$>1$\,). Then the \,$2N+1$\, 
coefficients of a polynomial with zeros in all the \,$\lambda_k$\,, \,$|k|\geq N$\, (these coefficients are the elementary symmetric polynomials in \,$\lambda_k$\,, \,$|k|\leq N$\,), 
together with the \,$\lambda_k$\, with \,$|k|>N$\, provide coordinates for \,$\Div$\, near \,$D$\,. 

\part{The inverse problem for the monodromy}

\section{Asymptotic spaces of holomorphic functions}
\label{Se:As}

We now introduce ``asymptotic spaces'' of holomorphic functions on \,$\C^*$\, or on \,$\Sigma$\, which have prescribed descent to \,$0$\, for \,$\lambda\to\infty$\,
and/or for \,$\lambda\to 0$\,; for this purpose we will cover \,$\C^*$\, by a sequence \,$(S_k)_{k\in \Z}$\, of annuli, the descent of the functions described by these spaces
will be uniform on each of these annuli \,$S_k$\, (up to a factor \,$w(\lambda)^s$\,). 

The purpose of these spaces is to describe the asymptotic behavior of the monodromy of a potential \,$(u,u_y)\in \Pot$\,; more specifically
we will see in Sections~\ref{Se:interpolate} and \ref{Se:asympfinal} that the differences between the functions comprising the monodromy matrix and the corresponding
functions for the vacuum are members of such asymptotic spaces as we are about to introduce. This formulation of the asymptotic behaviour of the monodromy 
(which we will regard as the final asymptotics) will liberate us from the special role that the points \,$\lambda_{k,0}$\, play in the Fourier asymptotic of Theorem~\ref{T:fasymp:fourier}. 

It is a remarkable fact that we will be able to control the asymptotic behavior of the monodromy by way of these asymptotic spaces not only on the excluded domains, but on all of \,$\C^*$\,.
This is similar to the basic asymptotics of Theorem~\ref{T:asymp:basic}, but it is very interesting that the addition of the Fourier asymptotics of Theorem~\ref{T:fasymp:fourier}
which states only the \,$\ell^2$-summability of the asymptotic difference at the single sequence of points \,$(\lambda_{k,0})$\, already implies control of an \,$\ell^2$-summability-type
not only in excluded domains (i.e.~small neighborhoods of \,$\lambda_{k,0}$\,) but on entire annuli in \,$\C^*$\, via the asymptotic spaces defined below.

The underlying reason why this control is possible, and the fact that inspired the definition of the asymptotic spaces, is that the holomorphic functions comprising the monodromy
can be described by infinite sums resp.~products (as we will see in Section~\ref{Se:interpolate}), and these infinite sums and products can be estimated on the \,$S_k$\,
essentially by the statements of Proposition~\ref{P:inf:sumholo}(2),(3) resp.~Proposition~\ref{P:inf:prodholo}(2).

In this section, we suppose that a spectral curve \,$\Sigma$\, is given only where we define or investigate objects involving \,$\Sigma$\,. 

For \,$k\in \Z$\, we put
\begin{equation}
\label{eq:As:Sk}
S_k := \begin{cases}
\Mengegr{\lambda\in\C^*}{(k-\tfrac12)\pi \leq |\zeta(\lambda)| \leq (k+\tfrac12)\pi, |\lambda| > 1} & \text{for \,$k>0$\,} \\
\Mengegr{\lambda\in\C^*}{|\zeta(\lambda)|\leq \tfrac{\pi}{2}} & \text{for \,$k=0$\,} \\
\Mengegr{\lambda\in\C^*}{(-k-\tfrac12)\pi \leq |\zeta(\lambda)| \leq (-k+\tfrac12)\pi, |\lambda| < 1} & \text{for \,$k<0$\,} 
\end{cases}
\end{equation}
and
\begin{equation}
\label{eq:As:whSk}
\wh{S}_k := \Menge{(\lambda,\mu)\in \Sigma}{\lambda \in S_k} \; .
\end{equation}
We call \,$S_k$\, resp.~\,$\wh{S}_k$\, the \emph{annulus near \,$\lambda_{k,0}$\,} in \,$\C^*$\, resp.~in \,$\Sigma$\,.  
Note that for \,$0<\delta<\pi-\tfrac12$\,, \,$U_{k,\delta} \subset S_k$\, holds, that different \,$S_k$\, resp.~\,$\wh{S}_k$\, intersect at most at their boundary, 
and that we have \,$\bigcup_{k\in \Z} S_k = \C^*$\, resp.~\,$\bigcup_{k\in \Z} \wh{S}_k = \Sigma$\,. 

The reader is reminded of the function \,$w(\lambda)$\, introduced in Equation~\eqref{eq:asymp:w}. 

We now define the asymptotic spaces which we will use in Section~\ref{Se:asympfinal} to describe the asymptotic behavior of the holomorphic functions constituting the monodromy \,$M(\lambda)$\, of a potential \,$(u,u_y)\in \Pot$\,,
and of many other functions in the course of this work. 
We introduce three versions of the asymptotic spaces: One that controls the function both near \,$\lambda=\infty$\, and near \,$\lambda=0$\,, and one each where only one of the ``ends'' near \,$\lambda=\infty$\,
resp.~near \,$\lambda=0$\, is controlled.

\begin{Def}
\label{D:As:As}
Let \,$0<p\leq \infty$\,, \,$n,m\in \Z$\,, \,$s\geq 0$\, and a domain \,$G\subset \C^*$\, resp.~\,$\wh{G}\subset \Sigma$\, be given.%
\footnote{In most of our applications, we will have \,$G \in \{\C^*,V_\delta\}$\, resp.~\,$\wh{G}\in \{\Sigma,\wh{V}_\delta\}$\,, \,$p\in \{2,\infty\}$\, and \,$s\in \{0,1\}$\,.}

Then we say that a holomorphic function \,$f: G \to \C$\, resp.~\,$f: \wh{G} \to \C$\, 
has \emph{\,$\ell^p_{n,m}$-asymptotic of type \,$s$\,} if there exists a sequence \,$(a_k)_{k\in \Z} \in \ell^p_{n,m}(k)$\, of non-negative numbers such that 
$$ \forall \, k\in \Z \; \forall\,\lambda \in G \cap S_k \; : \; |f(\lambda)| \leq a_k \cdot w(\lambda)^s $$
resp.~
$$ \forall \, k\in \Z \; \forall\,\lambda \in \wh{G} \cap \wh{S}_k \; : \; |f(\lambda)| \leq a_k \cdot w(\lambda)^s $$
holds. We call any such sequence \,$(a_k)$\, a \emph{bounding sequence} for \,$f$\,. We denote the space of all \,$\ell^p_{n,m}$-asymptotic functions \,$f: G \to \C$\, resp.~\,$f: \wh{G}\to\C$\, of type \,$s$\, 
by \,$\As(G,\ell^p_{n,m},s)$\, resp.~\,$\As(\wh{G},\ell^p_{n,m},s)$\,. These spaces become Banach spaces via the norm
$$ \|f\|_{\As(G,\ell^p_{n,m},s)} := \inf \|a_k\|_{\ell^p_{n,m}} \qmq{resp.~} \|f\|_{\As(\wh{G},\ell^p_{n,m},s)} := \inf \|a_k\|_{\ell^p_{n,m}} \;, $$
where the infimum is taken over all bounding sequences \,$(a_k)$\, for \,$f$\,. 
\end{Def}

\begin{Def}
\label{D:As:As0inf}
Let \,$0<p\leq \infty$\,, \,$n,m\in \Z$\,, \,$s\geq 0$\, and a domain \,$G\subset \C^*$\, resp.~\,$\wh{G}\subset \Sigma$\, be given.
We put 
$$ G_\infty := \Mengegr{\lambda\in \C^*}{|\lambda|>1} \cap G \qmq{and} G_0 := \Mengegr{\lambda\in \C^*}{|\lambda|<1} \cap G $$
resp.~
$$ \wh{G}_\infty := \Mengegr{(\lambda,\mu)\in \Sigma}{|\lambda|>1} \cap \wh{G} \qmq{and} \wh{G}_0 := \Mengegr{(\lambda,\mu)\in \Sigma}{|\lambda|<1} \cap \wh{G} \; . $$

\begin{enumerate}
\item
We say that a holomorphic function \,$f: G \to \C$\, resp.~\,$f: \wh{G} \to \C$\, 
has \emph{\,$\ell^p_{n}$-asymptotic of type \,$s$\, near \,$\lambda=\infty$\,}, if \,$f|G_\infty \in \As(G_\infty,\ell^p_{n,0},s)$\, resp.~\,$f|\wh{G}_\infty \in \As(\wh{G}_\infty,\ell^p_{n,0},s)$\, holds.
We denote the space of these functions with \,$\As_\infty(G,\ell^p_n,s)$\, resp.~\,$\As_\infty(\wh{G},\ell^p_n,s)$\,, it becomes a Banach space via the norm
$$ \|f\|_{\As_\infty(G,\ell^p_n,s)} := \bigr\| \,f|G_\infty\,\bigr\|_{\As(G_\infty,\ell^p_{n,0},s)} \qmq{resp.} \|f\|_{\As_\infty(\wh{G},\ell^p_n,s)} := \bigr\| \,f|\wh{G}_\infty\,\bigr\|_{\As(\wh{G}_\infty,\ell^p_{n,0},s)} \; . $$
\item
We say that a holomorphic function \,$f: G \to \C$\, resp.~\,$f: \wh{G} \to \C$\, 
has \emph{\,$\ell^p_{m}$-asymptotic of type \,$s$\, near \,$\lambda=0$\,}, if \,$f|G_0 \in \As(G_0,\ell^p_{0,m},s)$\, resp.~\,$f|\wh{G}_0 \in \As(\wh{G}_0,\ell^p_{0,m},s)$\, holds.
We denote the space of these functions with \,$\As_0(G,\ell^p_m,s)$\, resp.~\,$\As_0(\wh{G},\ell^p_m,s)$\,, it becomes a Banach space via the norm
$$ \|f\|_{\As_0(G,\ell^p_m,s)} := \bigr\| \,f|G_0\,\bigr\|_{\As(G_0,\ell^p_{0,m},s)} \qmq{resp.} \|f\|_{\As_0(\wh{G},\ell^p_m,s)} := \bigr\| \,f|\wh{G}_0\,\bigr\|_{\As(\wh{G}_0,\ell^p_{0,m},s)} \; . $$
\end{enumerate}
\end{Def}


\begin{rem}
An entire function \,$f$\, with \,$f \in \As_\infty(\C,\ell^\infty_{0},s)$\, satisfies 
$$ |f(\lambda)| \leq C_1\cdot \exp(s\cdot |\IM(\zeta(\lambda))|) \leq C_2\cdot \exp(\tfrac{s}{4}\cdot \sqrt{|\lambda|}) \qmq{for \,$\lambda\in \C$\, with \,$|\lambda|$\, large,} $$
with constants \,$C_1,C_2>0$\, and therefore is an entire function of order \,$\rho=\tfrac12$\, and type \,$\sigma=\tfrac{s}{4}$\,. 
\end{rem}

The following proposition states several important facts on asymptotic functions. 
Proposition~\ref{P:interpolate:l2asymp}(1) shows that a holomorphic function on \,$\C^*$\, which is asymptotic
on some \,$V_\delta$\,, i.e.~on the area outside the excluded domains, is in fact asymptotic on all of \,$\C^*$\,; this will be especially useful in several instances.

\begin{prop}
\label{P:interpolate:l2asymp}
Let \,$0<p\leq \infty$\,, \,$n,m \in \Z$\, and \,$s \geq 0$\,. 
\begin{enumerate}
\item Let \,$f: \C^* \to \C$\, be a holomorphic function so that 
\,$f|V_\delta \in \As(V_\delta,\ell^p_{n,m},s)$\, holds for some \,$\delta>0$\,, and let \,$(a_k) \in \ell^p_{n,m}(k)$\, be a bounding sequence for \,$f$\,.

Then already \,$f \in \As(\C^*,\ell^p_{n,m},s)$\, holds and \,$((4e^\delta)^s \cdot a_k)$\, is a bounding sequence for \,$f$\,, in particular we have
$$ \|f\|_{\As(\C^*,\ell^p_{n,m},s)} \leq (4e^\delta)^s \cdot \bigr\|\,f|V_\delta\,\bigr\|_{\As(V_\delta,\ell^p_{n,m},s)} \; . $$
\item Let \,$G\subset \C^*$\, be a domain, \,$f\in \As(G,\ell^p_{n,m},s)$\, and \,$j \in \Z$\,. Then we have 
\,$\lambda^{j/2}\cdot f \in \As(G,\ell^p_{n-j,m+j},s)$\,, and for any bounding sequence \,$(a_k)$\, for \,$f$\,, 
$$ b_k := \begin{cases} (4\pi (k+1))^j\cdot a_k & \text{if \,$k>0$\,} \\ (4\pi)^j\cdot a_0 & \text{if \,$k=0$\,} \\ \left( \tfrac{1}{4\pi (|k|-1)} \right)^j\cdot a_k & \text{if \,$k<0$\,} \end{cases} $$
is a bounding sequence for \,$\lambda^{j/2}\cdot f$\,. In particular, we have
$$ \|\lambda^{j/2}\cdot f\|_{\As(G,\ell^p_{n-j,m+j},s)} \leq (8\pi)^{|j|} \cdot \|f\|_{\As(G,\ell^p_{n,m},s)} \; . $$
\item Let \,$f\in \As(\C^*,\ell^p_{n,m},s)$\, and \,$j\in \N$\,. Then for the \,$j$-th derivative \,$f^{(j)}$\, of \,$f$\,, we have \,$f^{(j)} \in \As(\C^*,\ell^p_{n+j,m-3j},s)$\,,
and if \,$(a_k)$\, is a bounding sequence for \,$f$\,, then \,$b_k := \tfrac{12^s\cdot j!}{r_k^j} \cdot \max\{a_{k-1},a_k,a_{k+1}\}$\, is a bounding sequence for \,$f^{(j)}$\,, where
$$ r_k := \begin{cases} k & \text{if \,$k>0$\,} \\ 1 & \text{if \,$k=0$\,} \\ \tfrac{1}{16\pi^2\,|k|^3} & \text{if \,$k<0$\,} \end{cases} \; . $$
In particular, we have
$$ \|f^{(j)}\|_{\As(\C^*,\ell^p_{n+j,m-3j},s)} \leq C \cdot \|f\|_{\As(\C^*,\ell^p_{n,m},s)} $$
with a constant \,$C>0$\,. 
\end{enumerate}
In (1)--(3), analogous statements hold for the spaces \,$\As_\infty(\C^*,\ell^p_{n},s)$\, and \,$\As_0(\C^*,\ell^p_{m},s)$\,.
\end{prop}

\begin{proof}
\emph{For (1).} 
For \,$k\in \Z$\, and \,$\lambda\in S_k \cap V_\delta$\, we have \,$|f(\lambda)| \leq a_k\,w(\lambda)^s$\, by definition.
We need to show that for \,$\lambda\in U_{k,\delta} \subset S_k$\,, 
$$ |f(\lambda)| \leq (4\,e^\delta)^s\,a_k \cdot w(\lambda)^s $$ 
holds. Indeed we have by Proposition~\ref{P:vac2:excldom-M0}(1) for any \,$\lambda'\in \overline{U_{k,\delta}}$\, 
$$ \frac12 \leq w(\lambda') \leq 2\,e^\delta $$
and therefore for \,$\lambda\in U_{k,\delta}$\, by the maximum principle for holomorphic functions
$$ |f(\lambda)| \leq \max_{\lambda' \in \partial U_{k,\delta}} |f(\lambda')| \leq a_k \cdot \left( \max_{\lambda' \in \partial U_{k,\delta}} w(\lambda')\right)^s \leq a_k \cdot (2\,e^{\delta})^s \leq (4\,e^{\delta})^s\,a_k \cdot w(\lambda)^s \; . $$

\emph{For (2).}
We let \,$(a_k)$\, be a bounding sequence for \,$f$\,, and define \,$(b_k)$\, as in the Proposition. Note that \,$(a_k) \in \ell^p_{n,m}(k)$\, implies \,$(b_k) \in \ell^p_{n-j,m+j}(k)$\,. 
For \,$k\in \Z$\, and \,$\lambda\in S_k$\,, we have
$$ \pi\cdot (k-\tfrac12) \leq |\zeta(\lambda)| \leq \pi \cdot (k+\tfrac12) $$
and therefore
$$ \begin{cases}
16\,\pi^2\,(k-1)^2 \leq |\lambda| \leq 16\,\pi^2\,(k+1)^2 & \text{if \,$k>0$\,} \\
4\,\pi^2 \leq |\lambda| \leq 16\,\pi^2 & \text{if \,$k=0$\,} \\
\tfrac{1}{16\,\pi^2\,(|k|+1)^2} \leq |\lambda| \leq \tfrac{1}{16\,\pi^2\,(|k|-1)^2} & \text{if \,$k<0$\,}
\end{cases} \;, $$
hence
$$ \begin{cases}
\left(4\,\pi\,(k-1)\right)^j \leq |\lambda|^{j/2} \leq \left(4\,\pi\,(k+1)\right)^j & \text{if \,$k>0$\,} \\
(2\pi)^j \leq |\lambda|^{j/2} \leq (4\pi)^j & \text{if \,$k=0$\,} \\
\left( \tfrac{1}{4\,\pi\,(|k|+1)} \right)^j \leq |\lambda|^{j/2} \leq \left( \tfrac{1}{4\,\pi\,(|k|-1)} \right)^j & \text{if \,$k<0$\,} 
\end{cases} \;. $$
For \,$\lambda \in S_k \cap G$\,, we therefore obtain
$$ |\lambda^{j/2}\cdot f(\lambda)| \leq |\lambda|^{j/2} \cdot a_k\cdot w(\lambda)^s \leq b_k \cdot w(\lambda)^s \; , $$
whence the claimed statements follow.

\emph{For (3).}
We let \,$(a_k)$\, be a bounding sequence for \,$f$\,, and define \,$(r_k)$\, and \,$(b_k)$\, as in the Proposition. Note that \,$(a_k) \in \ell^p_{n,m}(k)$\, implies \,$(b_k) \in \ell^p_{n+j,m-3j}(k)$\,. 
Let \,$k\in \Z$\, be given. \,$r_k$\, is chosen such that for any \,$\lambda \in S_k$\, and \,$\lambda'\in \C^*$\, with \,$|\lambda'-\lambda|=r_k$\,, we have
\begin{equation}
\label{eq:interpolate:l2asymp:3zeta}
|\zeta(\lambda')-\zeta(\lambda)| \leq 1 \; .
\end{equation}
It follows from \eqref{eq:interpolate:l2asymp:3zeta} that we have
\begin{equation}
\label{eq:interpolate:l2asymp:3Sk3}
\lambda' \in S_{k-1}\cup S_k \cup S_{k+1} \; .
\end{equation}
We also obtain from \eqref{eq:interpolate:l2asymp:3zeta} by Proposition~\ref{P:asymp:w}(1):
\begin{align*}
w(\lambda') & \leq 2\,e^{|\IM(\zeta(\lambda'))|} \leq 2\,e^{|\IM(\zeta(\lambda))|}\,e^{|\IM(\zeta(\lambda')-\zeta(\lambda))|} \leq 2\,e^{|\IM(\zeta(\lambda))|}\,e^{|\zeta(\lambda')-\zeta(\lambda)|} \\
& \leq 2\cdot 2\,w(\lambda)\cdot e \leq 12\,w(\lambda) \; . 
\end{align*}
We now obtain by Cauchy's inequality, applied to the holomorphic function \,$f$\, on the disk \,$\overline{B(\lambda,r_k)}$\,:
\begin{align*}
|f^{(j)}(\lambda)| 
& \leq \frac{j!}{r_k^j}\cdot \max_{|\lambda'-\lambda|=r_k} |f(\lambda')| \leq \frac{j!}{r_k^j}\cdot \max\{a_{k-1},a_k,a_{k+1}\}\cdot \left( \max_{|\lambda'-\lambda|=r_k} w(\lambda')\right)^s \\
& \leq \frac{j!}{r_k^j}\cdot \max\{a_{k-1},a_k,a_{k+1}\}\cdot (12\,w(\lambda))^s = b_k\cdot w(\lambda)^s \; ,
\end{align*}
whence the claimed statements follow.
\end{proof}

\section{Interpolating holomorphic functions}
\label{Se:interpolate}

One of the main results of this text is that the monodromy \,$M(\lambda)$\, can be reconstructed uniquely (up to a sign in the off-diagonal entries) from the spectral data \,$(\Sigma,\calD)$\,.
More specifically, the given spectral data \,$(\Sigma,\calD)$\, provide the following information on the holomorphic functions comprising the monodromy
$$ M(\lambda) = \begin{pmatrix} a(\lambda) & b(\lambda) \\ c(\lambda) & d(\lambda) \end{pmatrix} $$
to which they belong: The components \,$\lambda_*$\, of the points \,$(\lambda_*,\mu_*)$\, in the support of \,$\calD$\, give all the zeros of \,$c$\, (with multiplicity), and moreover the components \,$\mu_*$\,
provide information on function values of \,$a$\, resp.~\,$d$\, by the equations \,$a(\lambda_*)=\mu_*$\,, \,$d(\lambda_*)=\mu_*^{-1}$\,. (If \,$(\lambda_*,\mu_*)$\, is a point of degree \,$m\geq 2$\, in the divisor
\,$\calD$\,, then the generalized spectral divisor \,$\calD$\, also determines the derivatives \,$d'(\lambda_*),\dotsc,d^{(m-1)}(\lambda_*)$\, in a manner that will be explained in detail in Lemma~\ref{L:special:sectional-lemma}.)

So the problem at hand is to reconstruct a holomorphic function on \,$\C^*$\, from either the knowledge of its zeros, or from its function values at a sequence of points, in both cases together with the knowledge
of the asymptotic behavior of the function near \,$\lambda=\infty$\, and near \,$\lambda=0$\,. We will address these two problems in the present section.

Proposition~\ref{P:interpolate:lambda} is concerned with the reconstruction of a holomorphic function on \,$\C^*$\, ``with the asymptotic behavior of \,$c$\,'' from the knowledge of its zeros.
In a way, this proposition is an adaption of Hadamard's Factorization Theorem (see for example \cite{Conway:1978}, Theorem~XI.3.4, p.~289)
to our specific situation. The most significant difference between our situation and the classical Theorem 
is that whereas the classical Theorem concerns entire functions with zeros accumulating near \,$\lambda=\infty$\,, we are interested in holomorphic functions on \,$\C^*$\,,
whose zeros accumulate both near \,$\lambda=\infty$\, and near \,$\lambda=0$\,. Notice that similarly to Hadamard's Theorem, we obtain an explicit representation of \,$c$\, as an infinite product.

Thereafter we study in Corollary~\ref{C:interpolate:cdivlin} the behavior of the function \,$\tfrac{c(\lambda)}{\lambda-\lambda_k}$\, on \,$U_{k,\delta}$\,; here \,$\lambda_k$\, is a root of \,$c$\, and \,$U_{k,\delta}$\,
is the excluded domain associated to this root. 
We need to understand the behavior of such functions both for the proof of Proposition~\ref{P:interpolate:mu} and on several further
occasions in the course of this work.

Finally, Proposition~\ref{P:interpolate:mu} concerns the reconstruction of a holomorphic function on \,$\C^*$\, ``with the asymptotic behavior of \,$a$\, or \,$d$\,'' from the knowledge of its function values
at the zeros of \,$c$\,; if \,$c$\, has zeros of higher order, then also values of the derivatives of \,$a$\, resp.~\,$d$\, at these points need to be known. We obtain an explicit description of \,$a$\,
resp.~\,$d$\, as an infinite series.

In Section~\ref{Se:special}, we will use Propositions~\ref{P:interpolate:lambda} and \ref{P:interpolate:mu} to reconstruct the monodromy \,$M(\lambda)$\, from its spectral data \,$(\Sigma,\calD)$\,. 

We mention that in the proofs of the present section (and nowhere else) we use the results on infinite sums and products collected in Appendix~\ref{Ap:inf}. 

\begin{prop}[Interpolation by the zeros.]
\label{P:interpolate:lambda}
\begin{enumerate}
\item
Let a sequence \,$(\lambda_k)_{k\in \Z}$\, in \,$\C^*$\, be given, such that we have
$$ \lambda_k-\lambda_{k,0} \in \ell^2_{-1,3}(k) \; . $$
Then the infinite product%
\footnote{Because of the freedom of choice of the branch of the square root function, \,$\tau$\, is determined only up to sign. 
In view of the fact that the off-diagonal entries of the monodromy are determined by their zeros only up to a sign (see Proposition~\ref{P:excl:unique}(1)), this is to be expected.}
\begin{equation}
\label{eq:interpolate:lambda:tau}
\tau := \left( \prod_{k\in \Z} \frac{\lambda_{k,0}}{\lambda_k} \right)^{1/2}
\end{equation}
converges absolutely in \,$\C^*$\,, and the infinite product
\begin{equation}
\label{eq:interpolate:lambda:c}
c(\lambda) = \frac14\,\tau\,(\lambda-\lambda_0) \cdot \prod_{k=1}^\infty \frac{\lambda_{k}-\lambda}{16\,\pi^2\,k^2} \cdot \prod_{k=1}^{\infty} \frac{\lambda-\lambda_{-k}}{\lambda}
\end{equation}
converges locally uniformly to a holomorphic function \,$c=c(\lambda):\C^* \to \C$\,. \,$c$\, has zeros in all the \,$\lambda_k$\, (with the appropriate multiplicity, if some of the \,$\lambda_k$\, coincide) and no others.
Moreover, we have 
\begin{equation}
\label{eq:interpolate:lambda:As}
c-\tau\,c_0 \in \As_\infty(\C^*,\ell^2_{-1},1) \qmq{and} c-\tau^{-1}\,c_0 \in \As_0(\C^*,\ell^2_1,1) \; .
\end{equation}
\item
Let \,$R_0>0$\, be given. Then there exists a constant \,$C>0$\, (depending only on \,$R_0$\,), such that for any pair of sequences \,$(\lambda_k^{[1]})_{k\in \Z}, (\lambda_k^{[2]})_{k\in \Z} \in \ell^2_{-1,3}(k)$\,
with \,$\|\lambda_k^{[\nu]}-\lambda_{k,0}\|_{\ell^2_{-1,3}} \leq R_0$\,, the quantities \,$\tau^{[\nu]}$\, and \,$c^{[\nu]}$\, corresponding to \,$(\lambda_k^{[\nu]})$\, as in (1) satisfy:
\begin{enumerate}
\item \,$|\tau^{[1]}-\tau^{[2]}|, |(\tau^{[1]})^{-1}-(\tau^{[2]})^{-1}| \leq C \cdot \|\lambda_k^{[1]}-\lambda_k^{[2]}\|_{\ell^2_{-1,3}}$\,
\item \,$\tau^{[2]}\,c^{[1]} - \tau^{[1]}\,c^{[2]} \in \As_\infty(\C^*,\ell^2_{-1},1)$\, \\
\,$(\tau^{[2]})^{-1}\,c^{[1]}(\lambda)-(\tau^{[1]})^{-1}\,c^{[2]}(\lambda) \in \As_0(\C^*,\ell^2_1,1)$\,
\item
With the notation \,$0^{-1} := 1$\,, 
\begin{equation}
\label{eq:interpolate:lambda:akdef}
a_k := \begin{cases} k^{-1}\cdot |\lambda_k^{[1]}-\lambda_k^{[2]}| & \text{for \,$k\geq 0$\,} \\ |k|^3\cdot |\lambda_k^{[1]}-\lambda_k^{[2]}| & \text{for \,$k<0$\,} \end{cases}
\end{equation}
and 
\begin{equation}
\label{eq:interpolate:lambda:rkdef}
r_k := C \cdot \begin{cases} k\cdot \left(a_k * \frac{1}{|k|} \right) & \text{for \,$k\geq 0$\,} \\ |k|^{-1}\cdot \left(a_k * \frac{1}{|k|} \right) & \text{for \,$k<0$\,} \end{cases} \;, 
\end{equation}
we have \,$r_k \in \ell^2_{-1,1}(k)$\,, and \,$(r_k)_{k>0}$\, resp.~\,$(r_k)_{k<0}$\, is a bounding sequence for \,$\tau^{[2]}\,c^{[1]} - \tau^{[1]}\,c^{[2]}$\, for \,$\lambda\to\infty$\, resp.~for
\,$(\tau^{[2]})^{-1}\,c^{[1]}(\lambda)-(\tau^{[1]})^{-1}\,c^{[2]}(\lambda)$\, for \,$\lambda\to 0$\,.
\item In particular, we have
$$
\left.
\begin{matrix}
\left\| \tau^{[2]}\,c^{[1]} - \tau^{[1]}\,c^{[2]} \right\|_{\As_\infty(\C^*,\ell^2_{-1},1)} \\
\left\| (\tau^{[2]})^{-1}\,c^{[1]}-(\tau^{[1]})^{-1}\,c^{[2]} \right\|_{\As_0(\C^*,\ell^2_1,1)}
\end{matrix}  \right\}
\;\leq\; C \cdot \left\|\lambda_k^{[1]}-\lambda_k^{[2]} \right\|_{\ell^2_{-1,3}} \; . $$
\end{enumerate}

\end{enumerate}
\end{prop}

\begin{proof}
We begin by looking at the case \,$\lambda_k = \lambda_{k,0}$\,. The corresponding number \,$\tau$\, defined by Equation~\eqref{eq:interpolate:lambda:tau} is clearly \,$\tau=1$\,, 
and we show that the function \,$c$\, defined by Equation~\eqref{eq:interpolate:lambda:c} equals
\,$c_0(\lambda) = \sqrt{\lambda}\,\sin(\zeta(\lambda))$\,. Indeed, by virtue of the product expansion of the sine function
$$ \sin(z) = z \cdot \prod_{k=1}^\infty \left( 1- \frac{z^2}{k^2\,\pi^2} \right) \;,$$
the following formulas for the \,$\lambda_{k,0}$\,
$$ \lambda_{0,0} = -1 \;, \quad \lambda_{k,0} + \lambda_{-k,0} = 16\,\pi^2\,k^2-2 \qmq{and} \lambda_{k,0} \cdot \lambda_{-k,0} = 1 \;, $$
and the equation 
$$ \zeta(\lambda)^2 = \frac{\lambda^2+2\lambda+1}{16\,\lambda} \;, $$
we have 
\begin{align}
c_0(\lambda) & = \sqrt{\lambda}\,\sin(\zeta(\lambda)) = \sqrt{\lambda} \,\zeta(\lambda) \,\prod_{k=1}^\infty \left( 1 - \frac{\zeta(\lambda)^2}{k^2\,\pi^2} \right) \notag \\
& = \frac14\,(\lambda+1) \,\prod_{k=1}^\infty \frac{16\pi^2k^2\,\lambda - (\lambda^2 + 2\lambda+1)}{16\pi^2k^2\cdot \lambda} \notag \\
& = \frac14\,(\lambda-\lambda_{0,0})\,\prod_{k=1}^\infty \frac{-\lambda^2+(\lambda_{k,0}+\lambda_{-k,0})\,\lambda-\lambda_{k,0}\,\lambda_{-k,0}}{16\pi^2k^2 \cdot \lambda} \notag \\
& = \frac14\,(\lambda-\lambda_{0,0})\,\prod_{k=1}^\infty \frac{(\lambda_{k,0}-\lambda)\cdot(\lambda-\lambda_{-k,0})}{16\pi^2k^2 \cdot \lambda} \notag \\
\label{eq:interpolate:lambda:c0}
& = \frac14\,(\lambda-\lambda_{0,0})\,\prod_{k=1}^\infty \frac{\lambda_{k,0}-\lambda}{16\,\pi^2\,k^2}\,\prod_{k=1}^\infty \frac{\lambda-\lambda_{-k,0}}{\lambda} \; . 
\end{align}
Note that the final expression in \eqref{eq:interpolate:lambda:c0} is the infinite product of Equation~\eqref{eq:interpolate:lambda:c}, where \,$\lambda_k=\lambda_{k,0}$\,. 

To prove the proposition, we will use the results from Appendix~\ref{Ap:inf} to show the convergence and estimation of the infinite products involved; we will apply the results from Appendix~\ref{Ap:inf}
both for products over \,$k\geq 1$\, and for products over \,$k\geq 0$\,, see Remark~\ref{R:inf:startingindex}. The numbers \,$C_k>0$\, occurring in the following estimates are constants that depend only on \,$R_0$\,. 
Moreover we choose \,$0 < \delta_0 <\delta < \pi-\tfrac12$\,, then all but finitely many of the \,$\lambda_k$\, resp.~\,$\lambda_k^{[\nu]}$\, are in \,$U_{k,\delta}$\,. 

We note that because of \,$\lambda_k-\lambda_{k,0} \in \ell^2_{-1,3}(k)$\,, we have both \,$\lambda_k-\lambda_{k,0} \in \ell^2_{-1}(k\geq 1)$\, and \,$\lambda_{-k}^{-1}-\lambda_{k,0} \in \ell^2_{-1}(k\geq 0)$\,,
and 
\begin{equation}
\label{eq:interpolate:lambda:normsplitting}
\left\| \lambda_k-\lambda_{k,0} \right\|_{\ell^2_{-1}(k\geq 1)} \;,\; \left\|\lambda_{-k}^{-1}-\lambda_{k,0}\right\|_{\ell^2_{-1}(k\geq 0)} \leq C_1 \cdot R_0 \; . 
\end{equation}
Inspired by this ``decomposition'', we write
$$ c(\lambda) = \frac{1}{4}\, \tau\cdot \rho \cdot c_+(\lambda)\cdot c_-\left(\frac{1}{\lambda}\right) \cdot \lambda $$
with
$$ c_+(\lambda) := \prod_{k=1}^\infty \left( 1-\frac{\lambda}{\lambda_k} \right)\;,\quad c_-(\lambda) := \prod_{k=0}^\infty \left( 1-\frac{\lambda}{\lambda_{-k}^{-1}} \right)\;,\quad \rho := \prod_{k=1}^\infty \frac{\lambda_k}{16\,\pi^2\,k^2} \;, $$
and 
$$ \tau = \left( \sigma_+ \cdot \sigma_-^{-1} \right)^{1/2} $$
with
$$ \sigma_+ := \prod_{k=1}^\infty \frac{\lambda_{k,0}}{\lambda_k} \qmq{and} \sigma_- := \prod_{k=0}^\infty \frac{\lambda_{k,0}}{\lambda_{-k}^{-1}} \; . $$
Then the infinite products defining \,$\sigma_\pm$\, and \,$\rho$\, converge absolutely in \,$\C^*$\, by Proposition~\ref{P:inf:prod}, 
and the infinite products defining \,$c_\pm(\lambda)$\, converge locally uniformly to holomorphic functions on \,$\C$\,
by Proposition~\ref{P:inf:prodholo}(1). It also follows from Proposition~\ref{P:inf:prod} that
$$ |\sigma_{\pm}|\,,\; |\sigma_\pm^{-1}|\,,\; |\rho|\,,\; |\rho^{-1}| \leq C_2 $$
and therefore also 
\begin{equation}
\label{eq:interpolate:lambda:taubound}
|\tau|\,,\; |\tau^{-1}| \leq C_2
\end{equation}
holds.

It follows that the product defining \,$\tau$\, converges absolutely in \,$\C^*$\,, and the infinite products in the definition of \,$c$\, in Equation~\eqref{eq:interpolate:lambda:c}
converge absolutely and locally uniformly to a well-defined holomorphic function \,$c$\,. It is clear that \,$c$\, has zeros in all the \,$\lambda_k$\, (with appropriate multiplicity), and no others. 

Before we show the asymptotic behavior of \,$c$\, given in \eqref{eq:interpolate:lambda:As}, we first work in the setting of (2). For this purpose, we let a 
pair of sequences \,$(\lambda_k^{[1]})_{k\in \Z}, (\lambda_k^{[2]})_{k\in \Z} \in \ell^2_{-1,3}(k)$\,
with \,$\|\lambda_k^{[\nu]}-\lambda_{k,0}\|_{\ell^2_{-1,3}} \leq R_0$\, be given, and denote the quantities associated to \,$(\lambda_k^{[\nu]})$\, by the superscript \,${}^{[\nu]}$\, (for \,$\nu \in \{1,2\}$\,). 
In relation to the splitting \eqref{eq:interpolate:lambda:normsplitting}, we note that we have
\begin{align}
C_3 \cdot \left\| \lambda_k^{[1]} - \lambda_k^{[2]} \right\|_{\ell^2_{-1,3}} 
& \leq \left\| \lambda_k^{[1]} - \lambda_k^{[2]} \right\|_{\ell^2_{-1}(k\geq 1)} + \left\| (\lambda_k^{[1]})^{-1} - (\lambda_k^{[2]})^{-1} \right\|_{\ell^2_{-1}(k\geq 0)} \notag \\
\label{eq:interpolate:lambda:normsplitting2}
& \leq C_4 \cdot \left\| \lambda_k^{[1]} - \lambda_k^{[2]} \right\|_{\ell^2_{-1,3}} \; . 
\end{align}

We now show (2)(a).
By Proposition~\ref{P:inf:prod} we have
$$ \left| \frac{\sigma^{[1]}_\pm}{\sigma^{[2]}_\pm} - 1 \right| \leq C_5 \cdot \left\| \lambda_k^{[1]} - \lambda_k^{[2]} \right\|_{\ell^2_{-1,3}} $$
and also (by exchanging the roles of \,$(\lambda_k^{[1]})$\, and \,$(\lambda_k^{[2]})$\,) 
$$ \left| \frac{\sigma^{[2]}_\pm}{\sigma^{[1]}_\pm} - 1 \right| \leq C_6 \cdot \left\| \lambda_k^{[1]} - \lambda_k^{[2]} \right\|_{\ell^2_{-1,3}} \;,$$
and therefore 
\begin{align*}
\left| \frac{\sigma_+^{[1]}}{\sigma_-^{[1]}} - \frac{\sigma_+^{[2]}}{\sigma_-^{[2]}} \right|
& \leq |\sigma_+^{[2]}| \cdot |\sigma_-^{[2]}|^{-1} \cdot \left( \left| \frac{\sigma_-^{[2]}}{\sigma_-^{[1]}} \right| \cdot \left| \frac{\sigma_+^{[1]}}{\sigma_+^{[2]}} - 1 \right| + \left| \frac{\sigma_-^{[2]}}{\sigma_-^{[1]}} - 1 \right| \right) \\
& \leq C_2^2 \cdot \left( C_2^2 \cdot C_5 \cdot \left\| \lambda_k^{[1]} - \lambda_k^{[2]} \right\|_{\ell^2_{-1,3}} + C_6 \cdot \left\| \lambda_k^{[1]} - \lambda_k^{[2]} \right\|_{\ell^2_{-1,3}} \right) \\
& = C_7 \cdot \left\| \lambda_k^{[1]} - \lambda_k^{[2]} \right\|_{\ell^2_{-1,3}} \; . 
\end{align*}
Because we have \,$|\tau^{[\nu]}|^2 \geq C_2^{-2} > 0 $\, by Equation~\eqref{eq:interpolate:lambda:taubound}, and 
the square root function is Lipschitz continuous on the interval \,$[C_2^{-2},\infty)$\,, it follows that we have
$$ |\tau^{[1]} - \tau^{[2]}| = \left| \left( \frac{\sigma_+^{[1]}}{\sigma_-^{[1]}} \right)^{1/2} - \left( \frac{\sigma_+^{[2]}}{\sigma_-^{[2]}} \right)^{1/2} \right|
\leq C_8 \cdot \left\| \lambda_k^{[1]} - \lambda_k^{[2]} \right\|_{\ell^2_{-1,3}} \; . $$
By exchanging the roles of \,$(\lambda_k^{[1]})$\, and \,$(\lambda_k^{[2]})$\, we also obtain 
$$ |\left( \tau^{[1]}\right)^{-1} - \left( \tau^{[2]}\right)^{-1} | \leq C_9 \cdot \left\| \lambda_k^{[1]} - \lambda_k^{[2]} \right\|_{\ell^2_{-1,3}} \; , $$
completing the proof of (2)(a).

Continuing the proof of the remaining parts of (2) and of (1), we at first consider only the parts of the Proposition concerned with the asymptotic behaviour of the function \,$c$\, for \,$\lambda\to\infty$\,.
We begin by investigating the function \,$\frac{\tau^{[2]}\,c^{[1]}(\lambda)}{\tau^{[1]}\,c^{[2]}(\lambda)}$\,. 
We have
\begin{align}
\frac{\tau^{[2]}\,c^{[1]}(\lambda)}{\tau^{[1]}\,c^{[2]}(\lambda)}-1 & = \frac{\tau^{[2]} \cdot \tfrac{1}{4}\, \tau^{[1]}\, \rho^{[1]} \, c_+^{[1]}(\lambda)\, c_-^{[1]}\left(\lambda^{-1}\right) \, \lambda}{\tau^{[1]} \cdot \tfrac{1}{4}\, \tau^{[2]}\, \rho^{[2]} \, c_+^{[2]}(\lambda)\, c_-^{[2]}\left(\lambda^{-1}\right) \, \lambda} -1 
= \frac{\rho^{[1]} \, c_+^{[1]}(\lambda)\, c_-^{[1]}\left(\lambda^{-1}\right)}{\rho^{[2]} \, c_+^{[2]}(\lambda)\, c_-^{[2]}\left(\lambda^{-1}\right)} -1 \notag \\
& = \left( \frac{\rho^{[1]} \, c_+^{[1]}(\lambda)}{\rho^{[2]} \, c_+^{[2]}(\lambda)} -1 \right) \cdot \frac{c_-^{[1]}\left(\lambda^{-1}\right)}{c_-^{[2]}\left(\lambda^{-1}\right)} 
+ \left( \frac{c_-^{[1]}\left(\lambda^{-1}\right)}{c_-^{[2]}\left(\lambda^{-1}\right)}  -1 \right) \notag \\
\label{eq:interpolate:lambda:infinitysplit}
& = \left( g(\lambda)-1 \right) \cdot h(\lambda) + \left( h(\lambda)-1 \right)
\end{align}
with
$$ g(\lambda) := \frac{\rho^{[1]} \, c_+^{[1]}(\lambda)}{\rho^{[2]} \, c_+^{[2]}(\lambda)} 
= \prod_{k=1}^\infty \frac{\lambda_k^{[1]}-\lambda}{\lambda_k^{[2]}-\lambda} $$
and
$$ h(\lambda) := \frac{c_-^{[1]}\left(\lambda^{-1}\right)}{c_-^{[2]}\left(\lambda^{-1}\right)}
= \prod_{k=0}^\infty \frac{\lambda_{-k}^{[1]}-\lambda}{\lambda_{-k}^{[2]}-\lambda} \; . $$
By Proposition~\ref{P:inf:prodholo}(2) we have for \,$\lambda \in S_n \cap V_\delta$\,, \,$n\geq 1$\,,
$$ |g(\lambda)-1| \leq r_n^{[g]} \;, $$
where 
$$ r_n^{[g]} := C_{10} \cdot \left( \left\{ \begin{matrix} a_k & \text{for \,$k>0$\,} \\ 0 & \text{for \,$k\leq 0$\,} \end{matrix} \right\} * \frac{1}{|k|} \right)_n $$
with a constant \,$C_{10}>0$\,, and the sequence \,$(a_k)$\, is the one from Equation~\eqref{eq:interpolate:lambda:akdef}.

Likewise by Proposition~\ref{P:inf:prodholo}(2), applied to \,$(\lambda_{-k}^{[\nu]})^{-1}$\, in the place of \,$\lambda_k^{[\nu]}$\,,
we obtain for \,$\lambda \in S_n \cap V_\delta$\,, \,$n\geq 1$\,,
$$ |h(\lambda)-1| \leq r_n^{[h]} \;, $$
where 
$$ r_n^{[h]} := C_{11} \cdot \left( \left\{ \begin{matrix} \frac{|(\lambda_{-k}^{[1]})^{-1}-|(\lambda_{-k}^{[2]})^{-1}|}{|k|} & \text{for \,$k\geq 0$\,} \\ 0 & \text{for \,$k< 0$\,} \end{matrix} \right\} * \frac{1}{|k|} \right)_n 
\leq C_{12}\cdot \left( \left\{ \begin{matrix} a_k & \text{for \,$k\leq 0$\,} \\ 0 & \text{for \,$k> 0$\,} \end{matrix} \right\} * \frac{1}{|k|} \right)_{-n}$$
with constants \,$C_{11},C_{12}>0$\,, and again with the sequence \,$(a_k)$\, from Equation~\eqref{eq:interpolate:lambda:akdef}.

We have \,$r_n^{[g]},r_n^{[h]} \leq r_n$\,, where the sequence \,$(r_n)$\, is defined by Equation~\eqref{eq:interpolate:lambda:rkdef}
(and we choose \,$C$\, as the maximum of \,$C_{10}$\, and \,$C_{12}$\,), and thus we have
\begin{equation}
\label{eq:interpolate:lambda:gh-asymp}
|g(\lambda)-1|\,,\;|h(\lambda)-1|\; \leq \; r_n \; . 
\end{equation}
From this estimate it also follows that there exists \,$C_{13}>0$\, so that 
\begin{equation}
\label{eq:interpolate:lambda:gh-estim}
|g(\lambda)|\,,\;|h(\lambda)|\; \leq \; C_{13} \; . 
\end{equation}
By taking the absolute value in Equation~\eqref{eq:interpolate:lambda:infinitysplit}, and then applying 
the estimates \eqref{eq:interpolate:lambda:gh-asymp} and \eqref{eq:interpolate:lambda:gh-estim}, we obtain
\begin{equation}
\label{eq:interpolate:lambda:quot-asymp}
\left|\frac{\tau^{[2]}\,c^{[1]}(\lambda)}{\tau^{[1]}\,c^{[2]}(\lambda)} - 1 \right| \leq C_{14}\cdot r_n
\end{equation}
with a constant \,$C_{14}>0$\,. 

By applying this estimate in the setting of (1), i.e.~for \,$\lambda_k^{[1]}=\lambda_k$\, and \,$\lambda_k^{[2]}=\lambda_{k,0}$\,,
we obtain in particular
\begin{equation*}
\left|\frac{c(\lambda)}{\tau\,c_0(\lambda)} - 1 \right| \leq C_{14}\cdot r_n \; . 
\end{equation*}
and therefore
$$ |c(\lambda)-\tau\,c_0(\lambda)| \leq C_{15}\cdot r_n \cdot |c_0(\lambda)| \; . $$
Because \,$|c_0(\lambda)| = |\lambda|^{1/2}\cdot |\sin(\zeta(\lambda))|$\, is comparable on \,$V_\delta$\, to \,$|\lambda|^{1/2}
\cdot w(\lambda)$\, by Proposition~\ref{P:vac2:excldom-M0}(3), it follows that \,$(c-\tau\,c_0)|V_\delta \in \As_\infty(V_\delta,\ell^2_{-1},1)$\,
holds. By Proposition~\ref{P:interpolate:l2asymp}(1) we in fact have \,$c-\tau\,c_0 \in \As_\infty(\C^*,\ell^2_{-1},1)$\,. 
This shows the half of \eqref{eq:interpolate:lambda:As} concerned with \,$\lambda\to\infty$\,. Moreover, because of
\,$c_0 \in \As_\infty(\C^*,\ell^\infty_{-1},1)$\,, we conclude \,$c\in \As_\infty(\C^*,\ell^\infty_{-1},1)$\, and thus
\begin{equation}
\label{eq:interpolate:lambda:c-estim}
|c(\lambda)| \leq C_{15} \cdot n \cdot w(\lambda)
\end{equation}
with a constant \,$C_{15}>0$\,. 

We return to the setting of (2). Again for \,$\lambda \in S_n \cap V_\delta$\, with \,$n\geq 1$\,, we have
$$ |\tau^{[2]}\,c^{[1]}(\lambda) - \tau^{[1]}\,c^{[2]}(\lambda)| = |\tau^{[1]}| \cdot |c^{[2]}(\lambda)| \cdot 
\left|\frac{\tau^{[2]}\,c^{[1]}(\lambda)}{\tau^{[1]}\,c^{[2]}(\lambda)} - 1 \right| \; . $$
From the estimates \eqref{eq:interpolate:lambda:taubound}, \eqref{eq:interpolate:lambda:c-estim}
and \eqref{eq:interpolate:lambda:quot-asymp} we thus conclude
$$ |\tau^{[2]}\,c^{[1]}(\lambda) - \tau^{[1]}\,c^{[2]}(\lambda)| \leq C_2 \cdot C_{15} \cdot w(\lambda) \cdot C_{14}\cdot r_n \; . $$
This shows that \,$(\tau^{[2]}\,c^{[1]}-\tau^{[1]}\,c^{[2]})|V_\delta \in \As_\infty(V_\delta,\ell^2_{-1},1)$\, and
therefore by Proposition~\ref{P:interpolate:l2asymp}(1)
\,$\tau^{[2]}\,c^{[1]}-\tau^{[1]}\,c^{[2]} \in \As_\infty(\C^*,\ell^2_{-1},1)$\, holds.
Moreover, we see that \,$(r_k)$\, is a bounding sequence for \,$\tau^{[2]}\,c^{[1]}-\tau^{[1]}\,c^{[2]}$\, 
(where the constant \,$C>0$\, occurring in the definition \eqref{eq:interpolate:lambda:rkdef} of \,$r_k$\, is enlarged
appropriately). This shows the half of (2)(b) and (2)(c) concerned with \,$\lambda\to\infty$\,, and the half
of (2)(d) concerned with \,$\lambda\to\infty$\, follows by application of the version of Young's inequality
for weakly \,$\ell^1$-sequences (see Equation~\eqref{eq:fasymp:fourier-small:weakyoung}). 

It remains to show the asymptotic statements in (1) and (2)(b)--(d) also for \,$\lambda\to 0$\,. We do this
by reducing the situation for \,$\lambda\to 0$\, to the previously proven situation for \,$\lambda\to\infty$\,. For this
purpose, we first note
$$ c_0(\lambda^{-1}) = \lambda^{-1}\cdot c_0(\lambda) \; . $$
Then we put \,$\wt{\lambda}_k := \lambda_{-k}^{-1}$\,, and denote the quantities associated to \,$(\wt{\lambda}_k)$\, by
a tilde \,$\widetilde{\ }$\,. With \,$\lambda_k-\lambda_{k,0} \in \ell^2_{-1,3}(k)$\,, we also have 
\,$\wt{\lambda}_k-\lambda_{k,0} \in \ell^2_{-1,3}(k)$\,, and for the corresponding sequences \,$(a_k)$\, and \,$(r_k)$\,
defined by Equation~\eqref{eq:interpolate:lambda:akdef} resp.~\eqref{eq:interpolate:lambda:rkdef}, we have
that \,$\wt{a}_k$\, is comparable to \,$a_{-k}$\,, and therefore \,$\wt{r}_k$\, is comparable to \,$r_{-k}$\,. 
Explicit calculations yield
$$ \sigma_+ = -\wt{\lambda}_0\cdot \wt{\sigma}_-\;,\quad \sigma_- = -\wt{\lambda}_0 \cdot \wt{\sigma}_+\;,\quad 
\tau = \wt{\tau}^{-1}\;,\quad \rho = -\wt{\lambda}_0 \cdot \wt{\tau}^2 \cdot \wt{\rho} $$
and for \,$\lambda \in \C^*$\,
$$ c_+(\lambda) = \frac{1}{1-\wt{\lambda}_0\cdot \lambda}\cdot \wt{c}_-(\lambda)
\qmq{and} c_-(\lambda) = \frac{\wt{\lambda}_0-\lambda}{\wt{\lambda}_0} \cdot \wt{c}_+(\lambda) $$
and therefore
\begin{align*}
c(\lambda^{-1}) & = \frac14\,\tau\cdot \rho\cdot c_+(\lambda^{-1}) \cdot c_-(\lambda) \cdot \lambda^{-1} \\
& = \frac14\,\wt{\tau}^{-1} \cdot (-\wt{\lambda}_0\,\wt{\tau}^2\,\wt{\rho}) \cdot \left( \frac{1}{1-\wt{\lambda}_0\, \lambda^{-1}}\, \wt{c}_-(\lambda^{-1}) \right) \cdot \left( \frac{\wt{\lambda}_0-\lambda}{\wt{\lambda}_0} \, \wt{c}_+(\lambda) \right) \cdot \lambda^{-1} \\
& = \frac14\cdot \wt{\tau} \cdot \wt{\rho} \cdot \wt{c}_+(\lambda) \cdot \wt{c}_-(\lambda^{-1}) \cdot \lambda \cdot \lambda^{-1} \\
& = \wt{c}(\lambda) \cdot \lambda^{-1} \; . 
\end{align*}
Using these formulas, the asymptotic estimates for \,$c$\, resp.~for \,$c^{[\nu]}$\, for \,$\lambda\to 0$\, follow from the previously
shown asymptotic estimates for \,$\lambda\to\infty$\, applied to \,$\wt{c}$\, resp.~to \,$\wt{c}^{[\nu]}$\,. 
\end{proof}

\begin{cor}
\label{C:interpolate:c-bestasymp}
Suppose that \,$c:\C^*\to\C$\, is a holomorphic function which satisfies the following two asymptotic properties 
for some \,$\tau\in\C^*$\,:
\begin{itemize}
\item[(a)] For every \,$\eps>0$\, there exists \,$R>0$\, such that we have
$$ \text{for all \,$\lambda\in \C^*$\, with \,$|\lambda|\geq R$\,: } \qquad |c(\lambda)-\tau\,c_0(\lambda)| \leq \eps\,|\lambda|^{1/2}\,w(\lambda) $$
and
$$ \text{for all \,$\lambda\in \C^*$\, with \,$|\lambda|\leq \tfrac{1}{R}$\,: } \qquad |c(\lambda)-\tau^{-1}\,c_0(\lambda)| \leq \eps\,|\lambda|^{1/2}\,w(\lambda) \;. $$
\item[(b)]
\,$(c(\lambda_{k,0}))_{k\in \Z} \in \ell^2_{-1,1}(k)$\,.
\end{itemize}
Then we already have
\begin{equation}
\label{eq:interpolate:c-bestasymp:asymp}
c-\tau\,c_0 \in \As_\infty(\C^*,\ell^2_{-1},1) \qmq{and} c-\tau^{-1}\,c_0 \in \As_0(\C^*,\ell^2_1,1)
\end{equation}
and
$$ \tau = \pm \left( \prod_{k\in\Z} \frac{\lambda_{k,0}}{\lambda_k} \right)^{1/2} \; , $$
where \,$(\lambda_k)_{k\in \Z}$\, is the sequence of zeros of \,$c$\, as in Proposition~\ref{P:excl:basic}(2).

\emph{Addendum.} If \,$L$\, is a set of holomorphic functions on \,$\C^*$\, which satisfy (a) and (b) in such a way that \,$R=R(\eps)>0$\, in (a) can be chosen uniformly for all \,$c\in L$\,, and in (b)
there is a uniform bound \,$(z_k)_{k\in \Z} \in \ell^2_{-1,1}(k)$\, for the sequences \,$(c(\lambda_{k,0}))$\, for all \,$c\in L$\,, then there exists a uniform bounding sequence for the asymptotics
in \eqref{eq:interpolate:c-bestasymp:asymp} that applies for all \,$c\in L$\,. 
\end{cor}

\begin{proof}
By Proposition~\ref{P:excl:basic}(2) the zeros of \,$c$\, are enumerated by a sequence \,$(\lambda_k)_{k\in \Z}$\,
such that for every \,$\delta>0$\,, there exists \,$N\in \N$\, so that \,$\lambda_k \in U_{k,\delta}$\, holds for \,$k\in \Z$\,
with \,$|k|\geq N$\,. By Corollary~\ref{C:asympdiv:asympdiv-neu} we have \,$\lambda_k-\lambda_{k,0}\in \ell^2_{-1,3}(k)$\,. 
Therefore Proposition~\ref{P:interpolate:lambda}(1) shows that there exists a holomorphic function \,$\wt{c}:\C^*\to\C$\,
that has zeros in all the \,$\lambda_k$\, (with the appropriate multiplicity) and no others, and which satisfies
$$ \wt{c}-\wt{\tau}\,c_0 \in \As_\infty(\C^*,\ell^2_{-1},1) \qmq{and} \wt{c}-\wt{\tau}^{-1}\,c_0 \in \As_0(\C^*,\ell^2_1,1) $$
with
$$ \wt{\tau} = \left( \prod_{k\in\Z} \frac{\lambda_{k,0}}{\lambda_k} \right)^{1/2} \; . $$
By Proposition~\ref{P:excl:unique}(1) we have \,$(c,\tau) = \pm(\wt{c},\wt{\tau})$\,, whence the claimed statement
follows.

In the situation of the addendum, the hypotheses ensure that there is a uniform sequence \,$(w_k)_{k\in \Z} \in \ell^2_{-1,1}(k)$\, of non-negative real numbers so that we have \,$|\lambda_k-\lambda_{k,0}| \leq w_k$\, for every sequence
\,$(\lambda_k)_{k\in \Z}$\, that is the sequence of zeros of a \,$c\in L$\,. Let
$$ a_k := \begin{cases} k^{-1}\cdot w_k & \text{for \,$k\geq 0$\,} \\ |k|^3\cdot w_k & \text{for \,$k<0$\,} \end{cases}
\qmq{and}
r_k := C  \cdot \begin{cases} k \cdot \left(a_k * \tfrac{1}{|k|}\right) & \text{for \,$k\geq 0$\,} \\ |k|^{-1} \cdot \left(a_k * \tfrac{1}{|k|}\right) & \text{for \,$k< 0$\,} \end{cases} $$
with the constant \,$C>0$\, from Proposition~\ref{P:interpolate:lambda}(2)(c). 
By applying that proposition with \,$\lambda_k^{[1]}=\lambda_k$\,, \,$\lambda_k^{[2]}=\lambda_{k,0}$\,,
we see that with the sequence \,$(r_k)_{k>0}$\, resp.~\,$(r_k)_{k<0}$\, is a bounding sequence for \,$\wt{c}-\wt{\tau}\,c_0$\, resp.~for \,$\wt{c}-\wt{\tau}^{-1}\,c_0$\, (that is independent of \,$c\in L$\,). 
\end{proof}

\begin{cor}
\label{C:interpolate:cdivlin}
\begin{enumerate}
\item
In the setting of Proposition~\ref{P:interpolate:lambda}(1), \,$\tfrac{c(\lambda)}{\lambda-\lambda_k}$\, is a holomorphic function on \,$\C^*$\, for every \,$k\in \Z$\,. 
There exists a constant \,$C>0$\, (depending only on \,$R_0$\,) and a sequence \,$(r_k)\in \ell^{2}_{0,-2}(k)$\, such that for every \,$k>0$\, and every \,$\lambda\in U_{k,\delta}$\, we have
$$ \left| \frac{c(\lambda)}{\lambda-\lambda_k} - \tau\,\frac{(-1)^k}{8} \right| \leq C\,\frac{|\lambda-\lambda_k|}{k} + r_k $$
and for every \,$k<0$\, and every \,$\lambda\in U_{k,\delta}$\, we have
$$ \left| \frac{c(\lambda)}{\lambda-\lambda_k} - \left( -\tau^{-1}\,\frac{(-1)^k}{8}\,\lambda_{k,0}^{-1} \right) \right| \leq C\,|\lambda-\lambda_k|\,k^5 + r_k \; . $$
Here we can choose
$$ r_k := \frac{1}{8}\,\lambda_{k,0}^{-1} + \begin{cases} C \cdot \left(a_k * \frac{1}{|k|} \right) & \text{for \,$k>0$\,} \\  C\,k^2 \cdot \left(a_k * \frac{1}{|k|} \right)& \text{for \,$k<0$\,}  \end{cases} 
\qmq{with}
a_k := \begin{cases} k^{-1}\cdot |\lambda_k-\lambda_{k,0}| & \text{for \,$k>0$\,} \\ k^3\cdot |\lambda_k-\lambda_{k,0}| & \text{for \,$k<0$\,} \end{cases} \; . $$
\item
In the setting of Proposition~\ref{P:interpolate:lambda}(2)
there exists a constant \,$C>0$\, (depending only on \,$R_0$\,)  and 
a sequence \,$r_k \in \ell^{2}_{0,-2}(k)$\, with \,$\|r_k\|_{\ell^{2}_{0,-2}} \leq C \cdot \|\lambda_k^{[1]}-\lambda_{k}^{[2]}\|_{\ell^{2}_{-1,3}}$\,
such that we have for all \,$k\in \Z$\, and \,$\lambda\in S_k$\, if \,$k>0$\,
$$ \left| \frac{c^{[1]}(\lambda)}{\tau^{[1]}\cdot(\lambda-\lambda_{k}^{[1]})} - \frac{c^{[2]}(\lambda)}{\tau^{[2]}\cdot(\lambda-\lambda_{k}^{[2]})} \right| \leq r_k $$
and if \,$k<0$\,
$$ \left| \frac{c^{[1]}(\lambda)}{(\tau^{[1]})^{-1}\cdot(\lambda-\lambda_{k}^{[1]})} - \frac{c^{[2]}(\lambda)}{(\tau^{[2]})^{-1}\cdot(\lambda-\lambda_{k}^{[2]})} \right| \leq r_k \; . $$

More specifically, we can choose
$$ r_k := \begin{cases} C \cdot \left(a_k * \frac{1}{|k|} \right) & \text{for \,$k>0$\,} \\  C\,k^2 \cdot \left(a_k * \frac{1}{|k|} \right)& \text{for \,$k<0$\,}  \end{cases}
\qmq{with}
a_k := \begin{cases} k^{-1}\cdot |\lambda_k^{[1]}-\lambda_k^{[2]}| & \text{for \,$k>0$\,} \\ k^3\cdot |\lambda_k^{[1]}-\lambda_k^{[2]}| & \text{for \,$k<0$\,} \end{cases} \; . $$
\end{enumerate}
\end{cor}

\begin{proof}
We consider only the case \,$k>0$\,, and prove (2) before (1).

\emph{For (2).}
We note that because 
the holomorphic function \,$c^{[\nu]}$\, has a zero at \,$\lambda=\lambda_k^{[\nu]}$\,, \,$\tfrac{c^{[\nu]}(\lambda)}{\lambda-\lambda_k^{[\nu]}}$\, is a holomorphic function on \,$\C^*$\,,
and we have
\begin{align}
& \frac{c^{[1]}(\lambda)}{\tau^{[1]}\cdot(\lambda-\lambda_{k}^{[1]})} - \frac{c^{[2]}(\lambda)}{\tau^{[2]}\cdot(\lambda-\lambda_{k}^{[2]})} \notag \\
\label{eq:interpolate:cdivlin:c12}
= & \frac{1}{\tau^{[1]}\,\tau^{[2]}} \left( (\tau^{[2]}\,c^{[1]}(\lambda)-\tau^{[1]}\,c^{[2]}(\lambda))\cdot \frac{1}{\lambda-\lambda_{k}^{[1]}} + \tau^{[1]}\,c^{[2]}(\lambda)\cdot \left( \frac{1}{\lambda-\lambda_{k}^{[1]}} - \frac{1}{\lambda-\lambda_{k}^{[2]}} \right) \right) \; . 
\end{align}

We choose \,$\delta>0$\, so large that \,$\lambda_k^{[\nu]}\in U_{k,\delta}$\, holds for all \,$k$\,, and at first consider
\,$\lambda \in S_k \cap V_{2\delta}$\,. Thus there exists \,$C_1>0$\, with \,$|\lambda-\lambda_k^{[\nu]}| \geq C_1\cdot k$\,
for all \,$k>0$\,, \,$\nu\in \{1,2\}$\, and all \,$\lambda\in S_k \cap V_{2\delta}$\,. We also have
$$ \left| \frac{1}{\lambda-\lambda_{k}^{[1]}} - \frac{1}{\lambda-\lambda_{k}^{[2]}} \right| = \frac{|\lambda_{k}^{[1]}-\lambda_k^{[2]}|}{|\lambda-\lambda_k^{[1]}| \cdot |\lambda-\lambda_k^{[2]}|}
\leq \frac{|\lambda_{k}^{[1]}-\lambda_k^{[2]}|}{(C_1\,k)^2} \; . $$
Moreover, by Proposition~\ref{P:interpolate:lambda}(2) we have \,$(\tau^{[2]}\,c^{[1]}(\lambda)-\tau^{[1]}\,c^{[2]}(\lambda)) \in \As_\infty(\C^*,\ell^2_{-1},1)$\, and \,$c^{[\nu]}(\lambda) \in \As_\infty(\C^*,\ell^\infty_{-1},1)$\,.
By plugging these results into Equation~\eqref{eq:interpolate:cdivlin:c12}, we see that there is a sequence \,$r_k \in \ell^2_0(k>0)$\, with \,$\|r_k\|_{\ell^2_0} \leq C \cdot \|\lambda_k^{[1]}-\lambda_k^{[2]}\|_{\ell^2_{-1,3}}$\, for some \,$C>0$\,
such that for all \,$k>0$\, and all \,$\lambda \in S_k\cap V_{2\delta}$\, we have
$$ \left| \frac{c^{[1]}(\lambda)}{\tau^{[1]}\cdot(\lambda-\lambda_{k}^{[1]})} - \frac{c^{[2]}(\lambda)}{\tau^{[2]}\cdot(\lambda-\lambda_{k}^{[2]})} \right| \leq r_k \; . $$
Because of the maximum principle for holomorphic functions, this inequality then also holds for \,$\lambda \in U_{k,2\delta} \subset S_k$\,.

The ``more specific'' description of \,$r_k$\, follows from the corresponding description in Proposition~\ref{P:interpolate:lambda}(2).

\emph{For (1).}
We consider the holomorphic functions \,$g_k(\lambda) := \tfrac{c(\lambda)}{\lambda-\lambda_k}$\, and \,$g_{k,0}(\lambda) := \tfrac{c_0(\lambda)}{\lambda-\lambda_{k,0}}$\, on \,$U_{k,\delta}$\,. 
For \,$\lambda\in U_{k,\delta}$\,, we have
\begin{equation}
\label{eq:interpolate:cdivlin:cc0-split}
\frac{c(\lambda)}{\lambda-\lambda_k} - \tau\,\frac{(-1)^k}{8}
= g_k(\lambda) - g_k(\lambda_{k,0}) + g_k(\lambda_{k,0}) - \tau\,g_{k,0}(\lambda_{k,0}) + \tau\,\left(g_{k,0}(\lambda_{k,0}) - \frac{(-1)^k}{8} \right) \; . 
\end{equation}
We estimate the three differences on the right hand side separately.

First, for \,$\lambda\in U_{k,\delta}$\,, we have the Taylor approximation of \,$c$\, near its zero \,$\lambda=\lambda_k$\,
$$ |c(\lambda)-c'(\lambda_k)\cdot(\lambda-\lambda_k)| \leq \frac{1}{2}\,|c''(\xi)| \cdot |\lambda-\lambda_k|^2 $$
with some \,$\xi \in U_{k,\delta}$\, (dependent on \,$\lambda$\,), and therefore
$$ \left| g_k'(\lambda) \right| = \left| \frac{c'(\lambda)\,(\lambda-\lambda_{k})-c(\lambda)}{(\lambda-\lambda_k)^2} \right| \leq \frac12\,|c''(\xi)| \; . $$
Because of \,$c'' \in \As(\C^*,\ell^\infty_{1,-5},1)$\,, there exists \,$C>0$\, with \,$\bigr|g_k'|U_{k,\delta}\bigr| \leq \tfrac{C}{k}$\,, and thus 
$$ |g_k(\lambda)-g_k(\lambda_{k,0})| \leq \frac{C}{k}\,|\lambda-\lambda_{k,0}| \leq \frac{C}{k}\,|\lambda-\lambda_{k}| + \frac{C}{k}|\lambda_k-\lambda_{k,0}| \; . $$

Second, we have \,$g_k(\lambda_{k,0}) - \tau\,g_{k,0}(\lambda_{k,0}) \in \ell^2_0(k)$\, by (2) (applied with \,$\lambda_k^{[1]}=\lambda_k$\, and \,$\lambda_k^{[2]} = \lambda_{k,0}$\,). And third, we have
$$ g_{k,0}(\lambda_{k,0}) - \frac{(-1)^k}{8} = c_0'(\lambda_{k,0})-\frac{(-1)^k}{8} = -\frac{(-1)^k}{8}\,\lambda_{k,0}^{-1} \; . $$

By plugging the preceding results into Equation~\eqref{eq:interpolate:cdivlin:cc0-split}, we obtain the claimed statement.
\end{proof}

\begin{prop}[Interpolation by the values.]
\label{P:interpolate:mu}
Let sequences \,$(\lambda_k)_{k\in \Z}$\, and \,$(\mu_k)_{k\in \Z}$\, and a number \,$\upsilon \in \C^*$\, be given, such that we have
$$ \lambda_k - \lambda_{k,0} \in \ell^2_{-1,3}(k) \qmq{and} \left. \begin{cases} \mu_k-\upsilon\,\mu_{k,0} & \text{if \,$k\geq 0$\,} \\ \mu_k-\upsilon^{-1}\,\mu_{k,0} & \text{if \,$k<0$\,} \end{cases} \right\} \in \ell^2_{0,0}(k) \; . $$
We further require that whenever \,$\lambda_k = \lambda_{\wt{k}}$\, holds for some \,$k,\wt{k}\in \Z$\,, we have \,$\mu_k = \mu_{\wt{k}}$\,. 

Then let \,$c:\C^*\to \C$\, be the holomorphic function with zeros at the \,$\lambda_k$\, defined in Proposition~\ref{P:interpolate:lambda},
and put for \,$k\in \Z$\,
$$ d_k := \# \Menge{\wt{k}\in \Z}{\lambda_{\wt{k}} = \lambda_k} = \ord_{\lambda_k}(c) \; . $$

\begin{enumerate}


\item
There are only finitely many \,$k\in \Z$\, with \,$d_k\geq 2$\,,
the infinite sum%
\footnote{Note that if \,$d_k>1$\, holds, then the summand \,$\tfrac{\mu_k \cdot (d_k-1)!\cdot c(\lambda)}{c^{(d_k)}(\lambda_k) \cdot (\lambda-\lambda_k)^{d_k}}$\, will occur \,$d_k$\, times in the sum defining \,$a$\,. This is the reason why the factor \,$(d_k-1)!$\, (instead of \,$d_k!$\,) occurs in the numerator.}
\begin{equation}
\label{eq:f-val:interpolate:g}
a(\lambda) := \sum_{k\in \Z} \frac{\mu_k \cdot (d_k-1)!\cdot c(\lambda)}{c^{(d_k)}(\lambda_k) \cdot (\lambda-\lambda_k)^{d_k}}
\end{equation}
converges absolutely and locally uniformly on \,$\C^*$\,, 
we have
\begin{equation}
\label{eq:interpolate:mu:a-asymp}
a-\upsilon\,a_0 \in \As_\infty(\C^*,\ell^2_0,1) \qmq{and} a-\upsilon^{-1}\,a_0 \in \As_0(\C^*,\ell^2_0,1) \;,
\end{equation}
and
\begin{equation}
\label{eq:interpolate:mu:a-values}
a(\lambda_k)=\mu_k \qmq{for all \,$k\in \Z$\,.}
\end{equation}

\item 
Let \,$\Lambda := \Menge{\lambda_k}{k\in\Z,\;d_k \geq 2}$\, be the set of multiple zeros of \,$c$\,,
and for each \,$\lambda_*\in \Lambda$\, choose 
\,$t_{\lambda_*,1},\dotsc,t_{\lambda_*,\ord_{\lambda_*}(c)-1}\in \C$\,. With the function \,$a$\, defined by Equation~\eqref{eq:f-val:interpolate:g}, \,$\wt{a}:\C^*\to\C$\,,
\begin{equation}
\label{eq:f-val:interpolate:wta}
\wt{a}(\lambda) := a(\lambda) + \sum_{\lambda_*\in\Lambda} \sum_{j=1}^{\ord_{\lambda_*}(c)-1} t_{\lambda_*,j}\cdot \frac{c(\lambda)}{(\lambda-\lambda_*)^j}
\end{equation}
is another holomorphic function with
\begin{equation}
\label{eq:interpolate:mu:wta-asymp}
\wt{a}-\upsilon\,a_0 \in \As_\infty(\C^*,\ell^2_0,1) \qmq{and} \wt{a}-\upsilon^{-1}\,a_0 \in \As_0(\C^*,\ell^2_0,1) \;,
\end{equation}
and
\begin{equation}
\label{eq:interpolate:mu:wta-values}
\wt{a}(\lambda_k)=\mu_k \qmq{for all \,$k\in \Z$\,.}
\end{equation}
Moreover, any holomorphic function \,$\wt{a}$\, that satisfies \eqref{eq:interpolate:mu:wta-asymp} and 
\eqref{eq:interpolate:mu:wta-values} is obtained in this way.

\item
Now let two pairs of sequences \,$(\lambda_k^{[\nu]})_{k\in \Z}$\,, \,$(\mu_k^{[\nu]})_{k\in \Z}$\, and numbers \,$\upsilon^{[\nu]}\in\C^*$\, (where \,$\nu\in\{1,2\}$\,) that satisfy the hypotheses of (1) be given, and denote the corresponding holomorphic function of (1)
with  \,$a^{[\nu]}: \C^* \to \C$\,. Then for every \,$R_0>0$\, there exists a constant \,$C>0$\, (depending only on \,$R_0$\,),
such that if we have 
$$ \left\|\lambda_k^{[\nu]}-\lambda_{k,0}\right\|_{\ell^2_{-1,3}} \leq R_0 \qmq{and} \left\| \left. \begin{cases} \mu_k-\upsilon\,\mu_{k,0} & \text{if \,$k\geq 0$\,} \\ \mu_k-\upsilon^{-1}\,\mu_{k,0} & \text{if \,$k<0$\,} \end{cases} \right\} \right\|_{\ell^2_{0,0}} \leq R_0 \;, $$
then the following holds:
\begin{enumerate}
\item \,$\upsilon^{[2]}\,a^{[1]}-\upsilon^{[1]}\,a^{[2]} \in \As_\infty(\C^*,\ell^2_0,1)$\,, \\
\,$(\upsilon^{[2]})^{-1}\,a^{[1]}-(\upsilon^{[1]})^{-1}\,a^{[2]} \in \As_0(\C^*,\ell^2_0,1)$\,
\item
With the notation \,$\tfrac10:=1$\,,
\begin{align}
a_k & := \begin{cases} k^{-1}\cdot |\lambda_k^{[1]}-\lambda_k^{[2]}| & \text{for \,$k\geq 0$\,} \\ |k|^3\cdot |\lambda_k^{[1]}-\lambda_k^{[2]}| & \text{for \,$k<0$\,} \end{cases} \notag \\
\label{eq:interpolate:mu:akbkdef}
\qmq{and} b_k & := \begin{cases} |\upsilon^{[2]}\,\mu_k^{[1]}-\upsilon^{[1]}\,\mu_k^{[2]}| & \text{for \,$k\geq 0$\,} \\ |(\upsilon^{[2]})^{-1}\,\mu_k^{[1]}-(\upsilon^{[1]})^{-1}\,\mu_k^{[2]}| & \text{for \,$k<0$\,} \end{cases} 
\end{align}
and 
\begin{align}
r_k & := C \cdot \left( \left( a_k * \frac{1}{|k|} \right) * \frac{1}{|k|} + \left( b_k + \frac{|\tau^{[1]}-\tau^{[2]}|}{|k|} \right) * \frac{1}{|k|} \right. \notag \\
\label{eq:interpolate:mu:rkdef}
& \qquad\qquad\qquad \left. + \frac{|\upsilon^{[1]}-(\upsilon^{[1]})^{-1}|+|\upsilon^{[2]}-(\upsilon^{[2]})^{-1}|}{|k|} \right) \;,
\end{align}
we have \,$r_k \in \ell^2_{0,0}(k)$\,, 
and \,$(r_k)_{k>0}$\, resp.~\,$(r_k)_{k<0}$\, is a bounding sequence for 
\,$\upsilon^{[2]}\,a^{[1]}-\upsilon^{[1]}\,a^{[2]}$\, for \,$\lambda\to\infty$\, resp.~for
\,$(\upsilon^{[2]})^{-1}\,a^{[1]}-(\upsilon^{[1]})^{-1}\,a^{[2]}$\, for \,$\lambda\to 0$\,.
\item In particular, we have
\begin{gather*}
\left.
\begin{matrix}
\left\| \upsilon^{[2]}\,a^{[1]}-\upsilon^{[1]}\,a^{[2]} \right\|_{\As_\infty(\C^*,\ell^2_{-1},1)} \\
\left\| (\upsilon^{[2]})^{-1}\,a^{[1]}-(\upsilon^{[1]})^{-1}\,a^{[2]} \right\|_{\As_0(\C^*,\ell^2_1,1)}
\end{matrix}  \right\} \\
\qquad\qquad \;\leq\; C \cdot \biggr( \left\|\lambda_k^{[1]}-\lambda_k^{[2]} \right\|_{\ell^2_{-1,3}} + \left\| b_k \right\|_{\ell^2_{0,0}} + |\tau^{[1]}-\tau^{[2]}| \\
\qquad\qquad\qquad\qquad\qquad\qquad\qquad +\; |\upsilon^{[1]}-(\upsilon^{[1]})^{-1}| + |\upsilon^{[2]}-(\upsilon^{[2]})^{-1}| \biggr) \; . 
\end{gather*}
\end{enumerate}
\end{enumerate}
\end{prop}

\begin{rem}
\label{R:interpolate:mu}
If no two \,$\lambda_k$\, coincide in the setting of Proposition~\ref{P:interpolate:mu}, then we have \,$d_k=1$\, for all \,$k$\,, and therefore
\begin{equation}
\label{eq:f-val:interpolate:g1}
a(\lambda) = \sum_{k\in \Z} \frac{\mu_k \cdot c(\lambda)}{c'(\lambda_k) \cdot (\lambda-\lambda_k)} 
\end{equation}
is the \emph{only} holomorphic function \,$a: \C^*\to\C$\, that satisfies \eqref{eq:interpolate:mu:a-asymp}
and \eqref{eq:interpolate:mu:a-values}.

In the general setting, Proposition~\ref{P:interpolate:mu}(1),(2) shows that for every \,$\lambda_k$\, that occurs more than once, i.e.~\,$d_k\geq 2$\,,
we may prescribe the values of \,$a'(\lambda_k),\dotsc,a^{(d_k-1)}(\lambda_k)$\, arbitrarily; then there exists one and only one holomorphic function \,$a:\C^*\to\C$\,
that satisfies \eqref{eq:interpolate:mu:a-asymp} and \eqref{eq:interpolate:mu:a-values}, and has the prescribed values of its derivatives.
\end{rem}

\begin{proof}[Proof of Proposition~\ref{P:interpolate:mu}.]
\emph{For (1) and (3).}
Because of the requirement \,$\lambda_k-\lambda_{k,0} \in \ell^2_{-1,3}(k)$\,, every excluded domain \,$U_{k,\delta}$\, contains asymptotically
and totally exactly one of the \,$\lambda_k$\,. 
It follows that only finitely many \,$\lambda_k$\, can coincide with one another, 
and therefore we have \,$d_k=1$\, for all  \,$k$\, with only finitely many exceptions. 

Every single summand \,$\frac{\mu_k \cdot (d_k-1)!\cdot c(\lambda)}{c^{(d_k)}(\lambda_k) \cdot (\lambda-\lambda_k)^{d_k}}$\, with \,$d_k \geq 2$\,
is by itself in \,$\As(\C^*,\ell^\infty_{2d_k-1,1},1)$\, and therefore in \,$\As(\C^*,\ell^2_{0,0},1)$\,. For the proof of the convergence results and the
claims on the asymptotic behavior, we may therefore suppose without loss of generality that \,$d_k=1$\, holds for all \,$k$\,.
Then \,$a(\lambda)$\, is given by Equation~\eqref{eq:f-val:interpolate:g1}. 

Because of the hypotheses, the sequence \,$(\mu_k)$\, is bounded, we have \,$c'-\tau\,c_0' \in \As(\C^*,\ell^2_{0,-2},1)$\,
(where \,$\tau$\, is the number defined by Equation~\eqref{eq:interpolate:lambda:tau}) and therefore \,$\tfrac{1}{c'(\lambda_k)}
\in \ell^\infty_{0,2}(k)$\,, and \,$\tfrac{c(\lambda)}{\lambda-\lambda_k} \in \ell^\infty_{2,0}(k)$\, locally uniformly with
respect to \,$\lambda\in \C^*$\,. It follows that 
$$ \frac{\mu_k\cdot c(\lambda)}{c'(\lambda_k) \cdot (\lambda-\lambda_k)} \in \ell^\infty_{2,2}(k) \subset \ell^1_{0,0}(k) $$
holds (again uniformly in \,$\lambda$\,). This shows that the infinite sum defining the function \,$a$\, in
Equation~\eqref{eq:f-val:interpolate:g1} indeed converges absolutely and locally uniformly, and thus 
Equation~\eqref{eq:f-val:interpolate:g1} defines a holomorphic function \,$a$\,. For \,$k,\ell \in \Z$\, we have
$$ \left. \frac{(d_k-1)!\cdot c(\lambda)}{c^{(d_k)}(\lambda_k) \cdot (\lambda-\lambda_k)^{d_k}} \right|_{\lambda=\lambda_\ell} = 
\begin{cases} \frac{1}{d_k} & \text{for \,$k=\ell$\,} \\ 0 & \text{for \,$k\neq \ell$\,} \end{cases} \; . $$
Because the term \,$\frac{\mu_k \cdot (d_k-1)!\cdot c(\lambda)}{c^{(d_k)}(\lambda_k) \cdot (\lambda-\lambda_k)^{d_k}}$\,
occurs \,$d_k$\, times in the infinite sum in Equation \eqref{eq:f-val:interpolate:g}, this shows that \,$a(\lambda_k)=\mu_k$\, holds. 

Before we show the claims in \eqref{eq:interpolate:mu:wta-asymp}
on the asymptotic behavior of the function \,$a$\,, we work in the setting of (3). As usual, we denote the objects associated
with \,$(\lambda_k^{[\nu]})$\, and \,$(\mu_k^{[\nu]})$\, by the superscript \,${}^{[\nu]}$\, (for \,$\nu\in \{1,2\}$\,). We
will first show the asymptotic statements in (3) for the case \,$\lambda\to\infty$\,. For this purpose, we will again
use statements from Appendix~\ref{Ap:inf}. The constants \,$C_k>0$\, occurring in the sequel only depend on \,$R_0$\, and
upper and lower bounds for \,$|\upsilon^{[\nu]}|$\,. 

We write 
$$ a^{[\nu]}(\lambda) = a_+^{[\nu]}(\lambda) + a_-^{[\nu]}(\lambda) $$
with
$$ a_+^{[\nu]}(\lambda) = \sum_{k=1}^\infty \frac{\mu_k^{[\nu]} \cdot c^{[\nu]}(\lambda)}{(c^{[\nu]})'(\lambda_k^{[\nu]}) \cdot (\lambda-\lambda_k^{[\nu]})} 
\qmq{and} a_-^{[\nu]}(\lambda) = \sum_{k=0}^\infty \frac{\mu_{-k}^{[\nu]} \cdot c^{[\nu]}(\lambda)}{(c^{[\nu]})'(\lambda_{-k}^{[\nu]}) \cdot (\lambda-\lambda_{-k}^{[\nu]})} \; . $$
To estimate \,$\upsilon^{[2]}\,a^{[1]}- \upsilon^{[1]}\,a^{[2]}$\,,
we will look at the functions \,$\upsilon^{[2]}\,a^{[1]}_+- \upsilon^{[1]}\,a^{[2]}_+$\, and 
\,$\upsilon^{[2]}\,a^{[1]}_-- \upsilon^{[1]}\,a^{[2]}_-$\, separately. The handling
of the two functions is quite dissimilar, owing to the different asymptotic behavior of \,$\lambda_k^{[\nu]}$\,, \,$\mu_k^{[\nu]}$\, and
\,$(c^{[\nu]})'(\lambda_{-k}^{[\nu]})$\, for \,$k\to\infty$\, and for \,$k\to-\infty$\,.

Let us first look at \,$\upsilon^{[2]}\,a^{[1]}_+- \upsilon^{[1]}\,a^{[2]}_+$\,. For \,$\lambda\in \C^*$\, we have
\begin{align}
\allowdisplaybreaks
& \upsilon^{[2]}\,a^{[1]}_+(\lambda)- \upsilon^{[1]}\,a^{[2]}_+(\lambda) \notag \\
= \; & \sum_{k=1}^\infty \left( \frac{\upsilon^{[2]}\,\mu_k^{[1]} \cdot (\tau^{[1]})^{-1}\,c^{[1]}(\lambda)}{(\tau^{[1]})^{-1}\,(c^{[1]})'(\lambda_k^{[1]}) \cdot (\lambda-\lambda_k^{[1]})} 
- \frac{\upsilon^{[1]}\,\mu_k^{[2]} \cdot (\tau^{[2]})^{-1}\,c^{[2]}(\lambda)}{(\tau^{[2]})^{-1}\,(c^{[2]})'(\lambda_k^{[2]}) \cdot (\lambda-\lambda_k^{[2]})}  \right) \notag \\
\label{eq:interpolate:mu:a+1}
= \; & (\tau^{[1]})^{-1}\,c^{[1]}(\lambda) \cdot \sum_{k=1}^\infty \bigr(\upsilon^{[2]}\,\mu_k^{[1]}- \upsilon^{[1]}\,\mu_k^{[2]}\bigr) \cdot 
\frac{1}{(\tau^{[1]})^{-1}\,(c^{[1]})'(\lambda_k^{[1]}) \cdot (\lambda-\lambda_k^{[1]})} \\
\label{eq:interpolate:mu:a+2}
& \qquad\qquad + ((\tau^{[1]})^{-1}\,c^{[1]}(\lambda)- (\tau^{[2]})^{-1}\,c^{[2]}(\lambda)) \cdot \sum_{k=1}^\infty \frac{\upsilon^{[1]}\,\mu_k^{[2]}}{(\tau^{[1]})^{-1}\,(c^{[1]})'(\lambda_k^{[1]})\cdot (\lambda-\lambda_k^{[1]})} \\
\label{eq:interpolate:mu:a+3}
& \qquad\qquad + (\tau^{[2]})^{-1}\,c^{[2]}(\lambda) \cdot \sum_{k=1}^\infty \left( \frac{1}{(\tau^{[1]})^{-1}\,(c^{[1]})'(\lambda_k^{[1]})}-\frac{1}{(\tau^{[2]})^{-1}\,(c^{[2]})'(\lambda_k^{[2]})} \right) \cdot \frac{\upsilon^{[1]}\,\mu_k^{[2]}}{\lambda-\lambda_k^{[1]}} \\
\label{eq:interpolate:mu:a+4}
& \qquad\qquad + (\tau^{[2]})^{-1}\,c^{[2]}(\lambda) \cdot \sum_{k=1}^\infty \left( \frac{1}{\lambda-\lambda_k^{[1]}}-\frac{1}{\lambda-\lambda_k^{[2]}} \right) \cdot \frac{\upsilon^{[1]}\,\mu_k^{[2]}}{(\tau^{[2]})^{-1}\,(c^{[2]})'(\lambda_k^{[2]})} \; . 
\end{align}
We will estimate the four expressions \eqref{eq:interpolate:mu:a+1}--\eqref{eq:interpolate:mu:a+4} individually. For this
purpose, we fix \,$n\in \N$\, and consider \,$\lambda \in S_n \cap V_\delta$\,. 

For the expression \eqref{eq:interpolate:mu:a+1}, we note that \,$(c^{[1]})'(\lambda_k^{[1]})$\, is bounded away from zero, and that we
have \,$b_k \in \ell^2(k\geq 1)$\, (where \,$b_k$\, is defined by Equation~\eqref{eq:interpolate:mu:akbkdef}), whence we obtain
\begin{gather*}
\left| \sum_{k=1}^\infty \bigr(\upsilon^{[2]}\,\mu_k^{[1]}- \upsilon^{[1]}\,\mu_k^{[2]}\bigr) \cdot 
\frac{1}{(\tau^{[1]})^{-1}\,(c^{[1]})'(\lambda_k^{[1]}) \cdot (\lambda-\lambda_k^{[1]})} \right| 
\leq C_1\cdot \sum_{k=1}^\infty \frac{b_k}{|\lambda-\lambda_k^{[1]}|}
\end{gather*}
with a constant \,$C_1>0$\,. By Proposition~\ref{P:inf:sumholo}(3), the latter sum is \,$\leq \frac{C_2}{n} \cdot \left( b_k * \frac{1}{|k|} \right)_n\leq \frac{C_3}{n}\,r_n$\, with constants \,$C_2,C_3>0$\, and \,$(r_k)$\, defined in 
Equation~\eqref{eq:interpolate:mu:rkdef}.
Because we have \,$c^{[1]} \in \As_\infty(\C^*,\ell^\infty_{-1},1)$\,, and hence \,$|c^{[1]}(\lambda)| \leq C_4\cdot n \cdot w(\lambda)$\,,
it follows that the expression \eqref{eq:interpolate:mu:a+1} is \,$\leq C_5 \cdot r_n \cdot w(\lambda)$\,.

For the expression \eqref{eq:interpolate:mu:a+2}, we note that \,$\upsilon^{[1]}\,\mu_k^{[2]} = \upsilon^{[1]}\,\upsilon^{[2]}\,(-1)^k + \ell^2(k\geq 1)$\, and moreover \,$(\tau^{[1]})^{-1}\,(c^{[1]})'(\lambda_k^{[1]}) = \tfrac18\,(-1)^k+\ell^2(k\geq 1)$\, holds, whence
it follows that there exists a sequence \,$t_k \in \ell^2(k\geq 1)$\, with 
$$ \frac{\upsilon^{[1]}\,\mu_k^{[2]}}{(\tau^{[1]})^{-1}\,(c^{[1]})'(\lambda_k^{[1]})} = 8\,\upsilon^{[1]}\,\upsilon^{[2]} + t_k $$
and hence
$$ \sum_{k=1}^\infty \frac{\upsilon^{[1]}\,\mu_k^{[2]}}{(\tau^{[1]})^{-1}\,(c^{[1]})'(\lambda_k^{[1]}) \cdot (\lambda-\lambda_k^{[1]})} 
= 8\,\upsilon^{[1]}\,\upsilon^{[2]}\,\sum_{k=1}^\infty \frac{1}{\lambda-\lambda_k^{[1]}} + \sum_{k=1}^\infty \frac{t_k}{\lambda-\lambda_k^{[1]}} \; . $$
By Proposition~\ref{P:inf:sumholo}(1),(3) it follows that there exists a constant \,$C_6>0$\, with 
\begin{equation}
\label{eq:interpolate:mu:a+2-sum}
\left| \sum_{k=1}^\infty \frac{\upsilon^{[1]}\,\mu_k^{[2]}}{(\tau^{[1]})^{-1}\,(c^{[1]})'(\lambda_k^{[1]}) \cdot (\lambda-\lambda_k^{[1]})} \right| \leq C_6 \cdot \frac{1}{n} \; . 
\end{equation}
On the other hand, by Proposition~\ref{P:interpolate:lambda}(2) we have 
\,$(\tau^{[1]})^{-1}\,c^{[1]}(\lambda) - (\tau^{[2]})^{-1}\,c^{[2]}(\lambda) \in \As_\infty(\C^*,\ell^2_{-1},1)$\, with
\,$C_7 \cdot k\cdot (a_k * \tfrac{1}{k})$\, being a bounding sequence. 
We have
$$ C_7 \cdot k\cdot (a_k * \tfrac{1}{k}) 
\leq C_8 \cdot k\cdot (a_k * (\tfrac1{|k|} * \tfrac1{|k|})) = C_8 \cdot k \cdot ((a_k *\tfrac1{|k|})*\tfrac1{|k|}) \leq k\cdot r_k\; . $$
Together with \eqref{eq:interpolate:mu:a+2-sum},
it follows that the expression \eqref{eq:interpolate:mu:a+2} is \,$\leq C_9\cdot r_n \cdot w(\lambda)$\,. 

For the expression \eqref{eq:interpolate:mu:a+3}, we note that by Proposition~\ref{P:interpolate:lambda}(2), we have
$$ (\tau^{[1]})^{-1}\,c^{[1]}(\lambda)-(\tau^{[2]})^{-1}\,c^{[2]}(\lambda) = (\tau^{[1]})^{-1}\,(\tau^{[2]})^{-1}\,(\tau^{[2]}\,c^{[1]}(\lambda)-\tau^{[1]}\,c^{[2]}(\lambda)) \in \As_\infty(\C^*,\ell^2_{-1},1) $$
and \,$C_{10}\cdot k\cdot (a_k * \tfrac1{|k|})$\, is a bounding sequence.
By Proposition~\ref{P:interpolate:l2asymp}(3), therefore \,$(\tau^{[1]})^{-1}\,(c^{[1]})'-(\tau^{[2]})^{-1}\,(c^{[2]})'\in \As_\infty(\C^*,\ell^2_0,1)$\, holds,
and \,$C_{11}\cdot (a_k*\tfrac{1}{|k|}) \in \ell^2(k\geq 1)$\, is a bounding sequence with a constant \,$C_{11}>0$\,. Moreover, \,$\upsilon^{[1]}\,\mu_k^{[2]}$\, is bounded.
Thus we have
$$ \left| \sum_{k=1}^\infty \left( \frac{1}{(\tau^{[1]})^{-1}\,(c^{[1]})'(\lambda_k^{[1]})}-\frac{1}{(\tau^{[2]})^{-1}\,(c^{[2]})'(\lambda_k^{[2]})} \right) \cdot \frac{\upsilon^{[1]}\,\mu_k^{[2]}}{\lambda-\lambda_k^{[1]}} \right|
\leq C_{12}\cdot \sum_{k=1}^\infty \frac{(a_k*\tfrac{1}{|k|})}{|\lambda-\lambda_k^{[1]}|} $$
with a constant \,$C_{12}>0$\,. By Proposition~\ref{P:inf:sumholo}(3) it follows that the sum to the right in the inequality above
is \,$\leq \tfrac{C_{13}}{n}\cdot ((a_k*\tfrac{1}{|k|})*\tfrac{1}{|k|})_n \leq \tfrac{C_{13}}{n}\cdot r_n$\,. Together with the
fact that \,$c^{[2]} \in \As_\infty(\C^*,\ell^\infty_{-1},1)$\, holds, it follows that the 
expression \eqref{eq:interpolate:mu:a+3} is \,$\leq C_{14}\cdot r_n \cdot w(\lambda)$\,. 

Finally, for the expression \eqref{eq:interpolate:mu:a+4} we again note that 
\,$\frac{\upsilon^{[1]}\,\mu_k^{[2]}}{(\tau^{[2]})^{-1}\,(c^{[2]})'(\lambda_k^{[2]})}$\, is bounded, and thus we have
$$ \left| \sum_{k=1}^\infty \left( \frac{1}{\lambda-\lambda_k^{[1]}}-\frac{1}{\lambda-\lambda_k^{[2]}} \right) \cdot \frac{\upsilon^{[1]}\,\mu_k^{[2]}}{(\tau^{[2]})^{-1}\,(c^{[2]})'(\lambda_k^{[2]})}\right|
\leq C_{15} \cdot \sum_{k=1}^\infty \frac{|\lambda_k^{[1]}-\lambda_k^{[2]}|}{|\lambda-\lambda_k^{[1]}| \cdot |\lambda-\lambda_k^{[2]}|} \; . $$
By Proposition~\ref{P:inf:sumholo}(4), the sum on the right hand side of the preceding inequality is
\,$\leq \tfrac{C_{16}}{n} \cdot (a_k*\tfrac{1}{|k|})_n \leq \tfrac{C_{16}}{n}\cdot r_n$\,. Similarly as for the expression
expression \eqref{eq:interpolate:mu:a+3} above, it follows that the 
expression \eqref{eq:interpolate:mu:a+4} is \,$\leq C_{17}\cdot r_n \cdot w(\lambda)$\,. 

The preceding estimates show that for \,$\lambda \in S_n \cap V_\delta$\,, 
\begin{equation}
\label{eq:interpolate:mu:a+}
\left( \upsilon^{[2]}\,a^{[1]}_+(\lambda)- \upsilon^{[1]}\,a^{[2]}_+(\lambda)\right) \leq C_{18} \cdot r_n \cdot w(\lambda)
\end{equation}
holds with a constant \,$C_{18}>0$\,.

We now attend to \,$\upsilon^{[2]}\,a^{[1]}_-- \upsilon^{[1]}\,a^{[2]}_-$\,. We let \,$\lambda\in \C^*$\, be given and write
\begin{align}
& \upsilon^{[2]}\,a^{[1]}_-(\lambda)- \upsilon^{[1]}\,a^{[2]}_-(\lambda) \notag \\
\label{eq:interpolate:mu:a-step1}
=\; & \bigr(\upsilon^{[2]}-(\upsilon^{[2]})^{-1}\bigr)\,a^{[1]}_-(\lambda) + \bigr( (\upsilon^{[2]})^{-1}\,a^{[1]}_-(\lambda) - (\upsilon^{[1]})^{-1}\,a^{[2]}_-(\lambda) \bigr) 
+ \bigr((\upsilon^{[1]})^{-1}-\upsilon^{[1]}\bigr)\,a^{[2]}_-(\lambda) \; . 
\end{align}
For estimating the various parts of Equation~\eqref{eq:interpolate:mu:a-step1}, we keep the following inequalities in mind: 
We fix \,$n\in \N$\, and let \,$\lambda \in S_n$\,. Then we have for \,$k\geq 0$\,
$$ \frac{1}{|\lambda-\lambda_{-k}^{[\nu]}|} \leq \frac{C_{19}}{|n^2-k^{-2}|} = \frac{C_{19}\cdot k^2}{|n^2\,k^2-1|} = \frac{C_{19}\cdot k^2}{|nk+1| \cdot |nk-1|} \leq \frac{C_{19}\cdot k^2}{nk \cdot |n+k|} = \frac{C_{19}\cdot k}{n\cdot |n+k|} $$
with a constant \,$C_{19}>0$\,. We may also assume without loss of generality that \,$|\lambda| \geq 2\,|\lambda_{-k}^{[\nu]}|$\, holds for all \,$k\geq 0$\, and \,$\nu\in\{1,2\}$\,, then we have
\begin{equation}
\label{eq:interpolate:mu:lambda-lambdak}
\frac{1}{|\lambda-\lambda_{-k}^{[\nu]}|} \leq \frac{1}{2\,|\lambda|} \; . 
\end{equation}
We also have \,$c^{[\nu]}\in \As_\infty(\C^*,\ell^\infty_{-1},1)$\, and therefore
$$ |c^{[\nu]}(\lambda)| \leq C_{20}\cdot n\cdot w(\lambda) \; . $$
Moreover, from the asymptotic for \,$c^{[\nu]}$\, for \,$\lambda\to 0$\,, namely \,$c^{[\nu]}-(\tau^{[\nu]})^{-1}\,c_0\in \As_0(\C^*,\ell^2_1,1)$\,
(Proposition~\ref{P:interpolate:lambda}(2)(b)), it follows by Proposition~\ref{P:interpolate:l2asymp}(3) that \,$(c^{[\nu]})'-(\tau^{[\nu]})^{-1}\,c_0'\in \As_0(\C^*,\ell^2_{-2},1)$\, holds, and therefore there exists
a constant \,$C_{21}>0$\, such that 
$$ \frac{1}{|(c^{[\nu]})'(\lambda_{-k}^{[\nu]})|} \leq C_{21}\cdot \frac{1}{k^2} $$
holds for \,$k\geq 0$\,. We also note that \,$\mu_{-k}^{[\nu]}$\, is bounded.

In the first instance, we obtain using these estimates
$$ |a^{[\nu]}_-(\lambda)| = \left| \sum_{k=0}^\infty \frac{\mu_{-k}^{[\nu]} \cdot c^{[\nu]}(\lambda)}{(c^{[\nu]})'(\lambda_{-k}^{[\nu]}) \cdot (\lambda-\lambda_{-k}^{[\nu]})} \right|
\leq C_{22}\cdot \frac{|c(\lambda)|}{|\lambda|} \cdot \sum_{k=0}^\infty \frac{1}{k^2} \leq C_{23}\cdot \frac{|c(\lambda)|}{|\lambda|} \; , $$
and therefore 
\begin{equation}
\label{eq:interpolate:mu:a-estim}
a^{[\nu]}_- \in \As_\infty(\C^*,\ell^\infty_1,1) \; .
\end{equation}
This shows that two of the terms occurring in \eqref{eq:interpolate:mu:a-step1}, namely \,$\bigr(\upsilon^{[2]}-(\upsilon^{[2]})^{-1}\bigr)\,a^{[1]}_-$\, and \,$\bigr((\upsilon^{[1]})^{-1}-\upsilon^{[1]}\bigr)\,a^{[2]}_-$\,
are in \,$\As_\infty(\C^*,\ell^\infty_1,1)$\, with \,$C_{23}\cdot \bigr|\upsilon^{[2]}-(\upsilon^{[2]})^{-1}\bigr|\cdot \tfrac{1}{k} \leq C_{23}\cdot r_k$\, 
resp.~\,$C_{23}\cdot \bigr|\upsilon^{[1]}-(\upsilon^{[1]})^{-1}\bigr|\cdot \tfrac{1}{k} \leq C_{23}\cdot r_k$\, being a bounding sequence.

We now estimate the remaining term in \eqref{eq:interpolate:mu:a-step1}, \,$(\upsilon^{[2]})^{-1}\,a^{[1]}_-(\lambda) - (\upsilon^{[1]})^{-1}\,a^{[2]}_-$\,. To do so, we split this expression similarly as we did for
\,$\upsilon^{[2]}\,a^{[1]}_+- \upsilon^{[1]}\,a^{[2]}_+$\,:
{\allowdisplaybreaks
\begin{align}
& (\upsilon^{[2]})^{-1}\,a^{[1]}_-(\lambda) - (\upsilon^{[1]})^{-1}\,a^{[2]}_- \notag \\
=\; & \sum_{k=0}^\infty \left( \frac{(\upsilon^{[2]})^{-1}\,\mu_{-k}^{[1]} \cdot (\tau^{[1]})^{-1}\,c^{[1]}(\lambda)}{(\tau^{[1]})^{-1}\,(c^{[1]})'(\lambda_{-k}^{[1]}) \cdot (\lambda-\lambda_{-k}^{[1]})} 
- \frac{(\upsilon^{[1]})^{-1}\,\mu_{-k}^{[2]} \cdot (\tau^{[2]})^{-1}\,c^{[2]}(\lambda)}{(\tau^{[2]})^{-1}\,(c^{[2]})'(\lambda_{-k}^{[2]}) \cdot (\lambda-\lambda_{-k}^{[2]})}  \right) \notag \\
\label{eq:interpolate:mu:a-1}
= \; & (\tau^{[1]})^{-1}\,c^{[1]}(\lambda) \cdot \sum_{k=0}^\infty \bigr((\upsilon^{[2]})^{-1}\,\mu_{-k}^{[1]}- (\upsilon^{[1]})^{-1}\,\mu_{-k}^{[2]}\bigr) \cdot 
\frac{1}{(\tau^{[1]})^{-1}\,(c^{[1]})'(\lambda_{-k}^{[1]}) \cdot (\lambda-\lambda_{-k}^{[1]})} \\
\label{eq:interpolate:mu:a-2}
& \qquad\qquad + ((\tau^{[1]})^{-1}\,c^{[1]}(\lambda)- (\tau^{[2]})^{-1}\,c^{[2]}(\lambda)) \cdot \sum_{k=0}^\infty \frac{(\upsilon^{[1]})^{-1}\,\mu_{-k}^{[2]}}{(\tau^{[1]})^{-1}\,(c^{[1]})'(\lambda_{-k}^{[1]})\cdot (\lambda-\lambda_{-k}^{[1]})} \\
\label{eq:interpolate:mu:a-3}
& \qquad\qquad + (\tau^{[2]})^{-1}\,c^{[2]}(\lambda) \cdot \sum_{k=0}^\infty \left( \frac{1}{(\tau^{[1]})^{-1}\,(c^{[1]})'(\lambda_{-k}^{[1]})}-\frac{1}{(\tau^{[2]})^{-1}\,(c^{[2]})'(\lambda_{-k}^{[2]})} \right) \cdot \frac{(\upsilon^{[1]})^{-1}\,\mu_{-k}^{[2]}}{\lambda-\lambda_{-k}^{[1]}} \\
\label{eq:interpolate:mu:a-4}
& \qquad\qquad + (\tau^{[2]})^{-1}\,c^{[2]}(\lambda) \cdot \sum_{k=0}^\infty \left( \frac{1}{\lambda-\lambda_{-k}^{[1]}}-\frac{1}{\lambda-\lambda_{-k}^{[2]}} \right) \cdot \frac{(\upsilon^{[1]})^{-1}\,\mu_{-k}^{[2]}}{(\tau^{[2]})^{-1}\,(c^{[2]})'(\lambda_{-k}^{[2]})} \; . 
\end{align}
}
We again treat the expressions \eqref{eq:interpolate:mu:a-1}--\eqref{eq:interpolate:mu:a-4} separately:

The expression~\eqref{eq:interpolate:mu:a-1} satisfies
\begin{align*}
|\eqref{eq:interpolate:mu:a-1}| & \leq C_{24}\cdot n\cdot w(\lambda) \cdot \sum_{k=0}^\infty \frac{b_{-k}}{k^2} \cdot \frac{k}{n\cdot |n+k|} \leq C_{24}\cdot w(\lambda)\cdot \sum_{k=0}^\infty b_{-k} \cdot \frac{1}{|n+k|} \\
& \leq C_{24}\cdot w(\lambda) \cdot \left( b_k * \frac{1}{|k|} \right)_n \leq C_{24}\cdot r_n \cdot w(\lambda) \; . 
\end{align*}

The expression~\eqref{eq:interpolate:mu:a-2} satisfies
\begin{align*}
|\eqref{eq:interpolate:mu:a-2}| & \leq C_{26}\cdot \left| (\tau^{[1]})^{-1}\,c^{[1]}(\lambda)- (\tau^{[2]})^{-1}\,c^{[2]}(\lambda)\right| \cdot \frac{1}{|\lambda|}\cdot \sum_{k=0}^\infty \frac{1}{k^2} \\
& \leq C_{27}\cdot \left| (\tau^{[1]})^{-1}\,c^{[1]}(\lambda)- (\tau^{[2]})^{-1}\,c^{[2]}(\lambda)\right| \cdot \frac{1}{|\lambda|} \; . 
\end{align*}
We have \,$(\tau^{[1]})^{-1}\,c^{[1]}- (\tau^{[2]})^{-1}\,c^{[2]} \in \As_\infty(\C^*,\ell^2_{-1},1)$\, with bounding sequence \,$C_{28}\cdot k\cdot (a_k*\tfrac{1}{|k|}) \leq C_{29}\cdot k\cdot r_k$\, by Proposition~\ref{P:interpolate:lambda}(2)(b),(c)
and therefore it follows that 
$$ |\eqref{eq:interpolate:mu:a-2}| \leq C_{29} \cdot r_n \cdot w(\lambda) \; . $$

For the expression~\eqref{eq:interpolate:mu:a-3} we have
\begin{align}
|\eqref{eq:interpolate:mu:a-3}| 
& \leq C_{30}\cdot n\cdot w(\lambda) \cdot \sum_{k=0}^\infty \left| \frac{1}{(\tau^{[1]})^{-1}\,(c^{[1]})'(\lambda_{-k}^{[1]})}-\frac{1}{(\tau^{[2]})^{-1}\,(c^{[2]})'(\lambda_{-k}^{[2]})} \right| \cdot \frac{k}{n\cdot|n+k|} \notag \\
\label{eq:interpolate:mu:a-3pre}
& \leq C_{31}\cdot w(\lambda) \cdot \sum_{k=0}^\infty \frac{\left| (\tau^{[2]})^{-1}\,(c^{[2]})'(\lambda_{-k}^{[2]})-(\tau^{[1]})^{-1}\,(c^{[1]})'(\lambda_{-k}^{[1]})\right|}{k^4} \cdot \frac{k}{|n+k|} \; . 
\end{align}
We now note that we have
\begin{align}
& (\tau^{[2]})^{-1}\,(c^{[2]})'(\lambda_{-k}^{[2]})-(\tau^{[1]})^{-1}\,(c^{[1]})'(\lambda_{-k}^{[1]}) \notag \\
=\; & \bigr( (\tau^{[2]})^{-1}-(\tau^{[1]})^{-1}\bigr)\,(c^{[2]})'(\lambda_{-k}^{[2]}) + \bigr( (\tau^{[1]})^{-1}\,(c^{[2]})'(\lambda_{-k}^{[2]})-(\tau^{[2]})^{-1}\,(c^{[1]})'(\lambda_{-k}^{[1]}) \bigr) \notag\\
& \qquad\qquad + \bigr( (\tau^{[2]})^{-1}-(\tau^{[1]})^{-1}\bigr)\,(c^{[1]})'(\lambda_{-k}^{[1]}) \notag \\
\notag
=\; & \tau^{[1]}\,\tau^{[2]}\,(\tau^{[1]}-\tau^{[2]})\,((c^{[1]})'(\lambda_{-k}^{[1]})+(c^{[2]})'(\lambda_{-k}^{[2]})) \notag \\
\notag
& \qquad\qquad + \bigr( (\tau^{[1]})^{-1}\,(c^{[2]})'(\lambda_{-k}^{[2]})-(\tau^{[2]})^{-1}\,(c^{[1]})'(\lambda_{-k}^{[1]}) \bigr) \; ,
\end{align}
and because of \,$(c^{[\nu]})' \in \As_0(\C^*,\ell^\infty_{-2},1)$\,
$$ \left| \tau^{[1]}\,\tau^{[2]}\,(\tau^{[1]}-\tau^{[2]})\,((c^{[1]})'(\lambda_{-k}^{[1]})+(c^{[2]})'(\lambda_{-k}^{[2]})) \right| \leq C_{32}\cdot k^2\cdot \left| \tau^{[1]}-\tau^{[2]} \right| \; . $$
Moreover, by Proposition~\ref{P:interpolate:lambda}(2)(b),(c) we have \,$(\tau^{[1]})^{-1}\,c^{[2]}- (\tau^{[2]})^{-1}\,c^{[1]} \in \As_0(\C^*,\ell^2_{1},1)$\, with \,$C_{33}\cdot \tfrac{1}{k} (a_k*\tfrac{1}{|k|})_{-k}$\, a bounding sequence,
and from this it follows by Proposition~\ref{P:interpolate:l2asymp}(3) that we have \,$(\tau^{[1]})^{-1}\,(c^{[2]})'- (\tau^{[2]})^{-1}\,(c^{[1]})' \in \As_0(\C^*,\ell^2_{-2},1)$\, with \,$C_{34}\cdot k^2\,(a_k*\tfrac{1}{|k|})_{-k}$\, 
a bounding sequence. Thus we have
$$ \left| (\tau^{[1]})^{-1}\,(c^{[2]})'(\lambda_{-k}^{[2]})-(\tau^{[2]})^{-1}\,(c^{[1]})'(\lambda_{-k}^{[1]}) \right| \leq C_{34}\cdot k^2\cdot \left(a_k * \frac{1}{|k|}\right)_{-k} \; , $$
and therefore
\begin{gather*}
\left| \tau^{[1]}\,\tau^{[2]}\,(\tau^{[1]}-\tau^{[2]})\,((c^{[1]})'(\lambda_{-k}^{[1]})+(c^{[2]})'(\lambda_{-k}^{[2]})) \right| \\
\leq k^2 \cdot \left(C_{32}\,\left| \tau^{[1]}-\tau^{[2]} \right| + C_{34}\, \left(a_k * \frac{1}{|k|}\right)_{-k} \right) \; . 
\end{gather*}
By plugging this result into \eqref{eq:interpolate:mu:a-3pre}, we obtain
\begin{align*}
|\eqref{eq:interpolate:mu:a-3}| 
& \leq C_{35}\cdot w(\lambda) \cdot \sum_{k=0}^\infty \frac{\left| \tau^{[1]}-\tau^{[2]} \right| + \left(a_k * \tfrac{1}{|k|}\right)_{-k}}{k}\cdot \frac{1}{|n+k|} \\
& \leq C_{35}\cdot w(\lambda) \cdot \left( \left( \frac{\left| \tau^{[1]}-\tau^{[2]} \right|}{|k|} + \left(a_k * \frac{1}{|k|}\right) \right) * \frac{1}{|k|} \right) \\
& \leq C_{36}\cdot r_n \cdot w(\lambda) \; . 
\end{align*}

Finally, for the expression~\eqref{eq:interpolate:mu:a-4} we have
\begin{align*}
|\eqref{eq:interpolate:mu:a-4}| 
& \leq C_{37}\cdot n\cdot w(\lambda) \cdot \sum_{k=0}^\infty \left| \frac{1}{\lambda-\lambda_{-k}^{[1]}}- \frac{1}{\lambda-\lambda_{-k}^{[2]}} \right| \cdot \frac{1}{k^2} \\
& \leq C_{37}\cdot n\cdot w(\lambda) \cdot \sum_{k=0}^\infty \frac{ |\lambda_{-k}^{[1]} - \lambda_{-k}^{[2]}|}{|\lambda-\lambda_{-k}^{[1]}| \cdot |\lambda-\lambda_{-k}^{[2]}|} \cdot \frac{1}{k^2} \\
& \leq C_{38}\cdot n\cdot w(\lambda) \cdot \sum_{k=0}^\infty |\lambda_{-k}^{[1]} - \lambda_{-k}^{[2]}| \cdot \frac{k^2}{n^2 \cdot |n+k|^2} \cdot \frac{1}{k^2} \\
& \leq C_{38}\cdot \frac{1}{n}\cdot w(\lambda) \cdot \sum_{k=0}^\infty \frac{k^{-3}\cdot a_{-k}}{|n+k|^2} \\
& \leq C_{38}\cdot w(\lambda) \cdot \sum_{k=0}^\infty a_{-k}\cdot \frac{1}{|n+k|} \leq C_{38}\cdot w(\lambda) \cdot \left( a_k * \frac{1}{|k|} \right)_{n} \\
& \leq C_{39}\cdot r_n \cdot w(\lambda) \; . 
\end{align*}

From the preceding estimates, we obtain that for all \,$n\in \N$\, and \,$\lambda\in S_n$\, we have
\begin{equation}
\label{eq:interpolate:mu:a-}
\left( \upsilon^{[2]}\,a^{[1]}_-(\lambda)- \upsilon^{[1]}\,a^{[2]}_-(\lambda)\right) \leq C_{40} \cdot r_n \cdot w(\lambda)
\end{equation}
holds with a constant \,$C_{40}>0$\,.

By combining \eqref{eq:interpolate:mu:a+} and \eqref{eq:interpolate:mu:a-} we obtain that
$$ \left( \upsilon^{[2]}\,a^{[1]}(\lambda)- \upsilon^{[1]}\,a^{[2]}(\lambda)\right) \leq C_{41} \cdot r_n \cdot w(\lambda) $$
holds with a constant \,$C_{41}>0$\,, and therefore we have \,$\upsilon^{[2]}\,a^{[1]}- \upsilon^{[1]}\,a^{[2]} \in \As_\infty(\C^*,\ell^2,1)$\, 
and that \,$r_n$\, is a bounding sequence if the constant \,$C>0$\, in Equation~\eqref{eq:interpolate:mu:rkdef} is chosen appropriately.
This shows the parts of (3)(a) and (3)(b) concerning \,$\lambda\to\infty$\,. 
The part of (3)(c) concerning \,$\lambda\to\infty$\, follows by an application of the version \eqref{eq:fasymp:fourier-small:weakyoung} of Young's inequality. 

We next show the part of \eqref{eq:interpolate:mu:a-asymp} in (1) concerning \,$\lambda\to \infty$\,, i.e.~we show
that in the situation of (1), \,$a-\nu\,a_0 \in \As_\infty(\C^*,\ell^2_0,1)$\, holds. To do so, we apply the situation
of (3) with \,$\lambda_k^{[1]} = \lambda_k$\,, \,$\mu_k^{[1]}=\mu_k$\, and \,$\lambda_k^{[2]}=\lambda_{k,0}$\,, 
\,$\mu_k^{[2]}=\mu_{k,0}$\,. Then \,$\upsilon^{[1]}=\upsilon$\, and \,$\upsilon^{[2]}=1$\, holds. By (3)(a) we then
have \,$a^{[1]}-\upsilon\,a^{[2]} \in \As_\infty(\C^*,\ell^2_0,1)$\,. We therefore only need to show that
\,$a^{[2]}$\, (i.e.~the function \,$a$\, obtained from Equation~\eqref{eq:f-val:interpolate:g1} by setting
\,$\lambda_k=\lambda_{k,0}$\, and \,$\mu_k=\mu_{k,0}$\,) equals \,$a_0=\cos(\zeta(\lambda))$\,.

We base our proof of this fact on the partial fraction expansion of the cotangent function: For \,$z\in \C\setminus \Menge{k\pi}{k\in \Z}$\,,
\begin{equation}
\label{eq:f-val:interpolate:cot}
\cot(z) = \frac{1}{z} + \sum_{k=1}^\infty \frac{2z}{z^2 - (k\pi)^2}
\end{equation}
holds. We also remember
\begin{equation}
\label{eq:f-zero:f0'}
c_0'(\lambda_{k,0}) = \frac{\lambda_{k,0}-1}{8\,\lambda_{k,0}}\,\mu_{k,0} \; . 
\end{equation}
Therefore we have for \,$\lambda \in \C^* \setminus \Menge{\lambda_{k,0}}{k\in \Z}$\,
\begingroup
\allowdisplaybreaks
\begin{align}
a_0(\lambda) 
& \;\;\;\;=\;\;\;\; \cos(\zeta(\lambda)) = \sin(\zeta(\lambda)) \cdot \cot(\zeta(\lambda)) = \frac{c_0(\lambda)}{\sqrt{\lambda}}\cdot \cot(\zeta(\lambda)) \notag \\
& \,\;\overset{\eqref{eq:f-val:interpolate:cot}}{=}\; \frac{c_0(\lambda)}{\sqrt{\lambda}}\cdot \left( \frac{1}{\zeta(\lambda)} + \sum_{k=1}^\infty \frac{2\,\zeta(\lambda)}{\zeta(\lambda)^2 - (k\pi)^2} \right) \notag \\
& \;\;\;\;=\;\;\;\; c_0(\lambda) \cdot \left( \frac{4}{\lambda+1} + \sum_{k=1}^\infty \frac{\lambda+1}{2\,\lambda\,(\zeta(\lambda)^2 - \zeta(\lambda_{k,0})^2)} \right) \notag \\
& \;\;\;\;=\;\;\;\; c_0(\lambda) \cdot \left( \frac{4}{\lambda+1} + \sum_{k=1}^\infty \frac{8\,(\lambda+1)\,\lambda_{k,0}}{(\lambda+1)^2\,\lambda_{k,0} - (\lambda_{k,0}+1)^2\,\lambda} \right) \notag \\
& \;\;\;\;=\;\;\;\; c_0(\lambda) \cdot \left( \frac{4}{\lambda+1} + \sum_{k=1}^\infty \frac{8\,(\lambda+1)}{(\lambda-\lambda_{-k,0})\,(\lambda-\lambda_{k,0})} \right) \notag \\
& \;\;\;\;=\;\;\;\; c_0(\lambda) \cdot \left( \frac{4}{\lambda+1} + \sum_{k=1}^\infty \frac{8}{\lambda_{k,0}-\lambda_{-k,0}}\left( \frac{\lambda_{k,0}+1}{\lambda-\lambda_{k,0}} - \frac{\lambda_{-k,0}+1}{\lambda-\lambda_{-k,0}} \right) \right) \notag \\
& \;\;\;\;=\;\;\;\; c_0(\lambda) \cdot \left( \frac{4}{\lambda+1} + \sum_{k\in\Z\setminus\{0\}} \frac{8}{\lambda_{k,0}-\lambda_{-k,0}}\,\frac{\lambda_{k,0}+1}{\lambda-\lambda_{k,0}} \right) \notag \\
& \;\overset{\eqref{eq:f-zero:f0'}}{=}\; c_0(\lambda) \cdot \left( \frac{4}{\lambda+1}\,\frac{(\lambda_{0,0}-1)\,\mu_{0,0}}{8\,\lambda_{0,0}\,f_0'(\lambda_{k,0})} + \sum_{k\in\Z\setminus\{0\}} \frac{8}{\lambda_{k,0}-\lambda_{-k,0}}\,\frac{\lambda_{k,0}+1}{\lambda-\lambda_{k,0}}\,\frac{(\lambda_{k,0}-1)\,\mu_{k,0}}{8\,\lambda_{k,0}\,c_0'(\lambda_{k,0})} \right) \notag \\
& \overset{\lambda_{0,0}=-1}{=} c_0(\lambda) \cdot \sum_{k\in \Z} \frac{\mu_{k,0}}{c_0'(\lambda_{k,0}) \cdot (\lambda-\lambda_{k,0})} \; . \notag
\end{align}
\endgroup
Therefore the desired equality holds for all \,$\lambda\in\C^* \setminus \Menge{\lambda_{k,0}}{k\in\Z}$\,. Because its both sides are holomorphic in \,$\lambda$\,, 
that equality in fact holds for \,$\lambda\in \C^*$\,. 

Next we show the parts of (1) and (3) concerning \,$\lambda\to 0$\,. Similarly as we did in the proof of 
Proposition~\ref{P:interpolate:lambda}, we do so by reducing the situation for \,$\lambda\to 0$\, to the previously
proven situation for \,$\lambda\to\infty$\,. For this purpose, we put \,$\wt{\lambda}_k := \lambda_{-k}^{-1}$\, 
and \,$\wt{\mu}_k := \mu_{-k}$\,. Mutatis mutandis, these sequences again satisfy the hypotheses of part (3) of the proposition.
Again we denote the quantities associated to \,$(\wt{\lambda}_k)$\, and \,$\wt{\mu}_k$\, by a tilde \,$\wt{\ }$\,. 
We note that \,$\wt{a}_k$\,, \,$\wt{b}_k$\, and \,$\wt{r}_k$\, is comparable to \,$a_{-k}$\,, \,$b_{-k}$\, and \,$r_{-k}$\,,
respectively. We may again suppose without loss of generality that \,$d_k=1$\, holds for all \,$k\in \Z$\,. Then we obtain by explicit calculation
$$ \tau = \wt{\tau}^{-1} \qmq{and} \upsilon = \wt{\upsilon}^{-1} \;, $$
and
$$\quad c(\lambda^{-1}) = \wt{c}(\lambda) \cdot \lambda^{-1} \qmq{and} c'(\lambda_k) = -\wt{\lambda}_{-k}\cdot \wt{c}'(\wt{\lambda}_{-k}) $$
and therefore
\begin{align*}
a(\lambda^{-1}) & = \sum_{k\in \Z} \frac{\mu_k \cdot c(\lambda^{-1})}{c'(\lambda_k) \cdot (\lambda^{-1}-\lambda_k)} 
= \sum_{k\in \Z} \frac{\wt{\mu}_{-k} \cdot \wt{c}(\lambda)\cdot \lambda^{-1}}{(-\wt{\lambda}_{-k})\cdot c'(\wt{\lambda}_{-k}) \cdot (\lambda^{-1}-\wt{\lambda}_{-k}^{-1})} \\
& = \sum_{k\in \Z} \frac{\wt{\mu}_{-k} \cdot \wt{c}(\lambda)}{c'(\wt{\lambda}_{-k}) \cdot (\lambda-\wt{\lambda}_{-k})} 
= \sum_{k\in \Z} \frac{\wt{\mu}_{k} \cdot \wt{c}(\lambda)}{c'(\wt{\lambda}_{k}) \cdot (\lambda-\wt{\lambda}_{k})} 
= \wt{a}(\lambda) \; . 
\end{align*}
Using these formulas, one derives the asymptotic estimates for \,$\lambda\to 0$\, in (1) and (3) from the corresponding
estimates for \,$\lambda\to\infty$\, that have been shown before.

\emph{For (2).}
For any \,$\lambda_* \in \Lambda$\, we put \,$d(\lambda_*) := \ord_{\lambda_*}(c) \geq 2$\,. For 
\,$1\leq j \leq d(\lambda_*)-1$\,, we have
\,$\vi_{\lambda_*,j}(\lambda) := \tfrac{c(\lambda)}{(\lambda-\lambda_*)^j} \in \As(\C^*,\ell^\infty_{2j-1,1},1) \subset \As(\C^*,\ell^2_{0,0},1)$\,, and
therefore \eqref{eq:interpolate:mu:a-asymp} implies \eqref{eq:interpolate:mu:wta-asymp}. Moreover, we have
\,$\vi_{\lambda_*,j}(\lambda_k)=0$\, for all \,$k$\, (including the case
\,$\lambda_k=\lambda_*$\,), and therefore \,$\wt{a}(\lambda_k)=a(\lambda_k)=\mu_k$\,. 

Now suppose that \,$\wh{a}$\, is another holomorphic function that satisfies \eqref{eq:interpolate:mu:a-asymp}
and \eqref{eq:interpolate:mu:a-values}. Because \,$\vi_{\lambda_*,j}$\, has at every \,$\lambda_k$\, with \,$\lambda_k\neq \lambda_*$\,
a zero of order \,$d_k$\,, and because the $i$-th derivative of \,$\vi_{\lambda_*,j}$\, at \,$\lambda_*$\, satisfies 
$$ \vi_{\lambda_*,j}^{(i)}(\lambda_*) = 0 \qmq{for \,$0 \leq i \leq d(\lambda_*)-j-1$\,} $$
and 
$$ \vi_{\lambda_*,j}^{(d(\lambda_*)-j)}(\lambda_*) \neq 0 \;, $$
the numbers \,$t_{\lambda_*,j}$\, can be chosen such that for every \,$\lambda_* \in \Lambda$\, and every \,$1 \leq j \leq d(\lambda_*)-1$\,,
\,$\wt{a}^{(j)}(\lambda_*) = \wh{a}^{(j)}(\lambda_*)$\, holds, i.e.~\,$\wt{a}-\wh{a}$\, has at \,$\lambda_*$\, a zero of order
(at least) \,$d(\lambda_*)$\,. By Proposition~\ref{P:excl:unique}(2), \,$\wt{a}=\wh{a}$\, follows.
\end{proof}

\begin{cor}
\label{C:interpolate:a-bestasymp}
Suppose that \,$a:\C^*\to\C$\, is a holomorphic function which satisfies the following two asymptotic properties 
for some \,$\upsilon\in\C^*$\,:
\begin{itemize}
\item[(a)] For every \,$\eps>0$\, there exists \,$R>0$\, such that we have
$$ \text{for all \,$\lambda\in \C^*$\, with \,$|\lambda|\geq R$\,: } \qquad |a(\lambda)-\upsilon\,a_0(\lambda)| \leq \eps\,w(\lambda) $$
and
$$ \text{for all \,$\lambda\in \C^*$\, with \,$|\lambda|\leq \tfrac{1}{R}$\,: } \qquad |a(\lambda)-\upsilon^{-1}\,a_0(\lambda)| \leq \eps\,w(\lambda) \;. $$
\item[(b)]
\,$\left. \begin{cases} a(\lambda_{k,0})-\upsilon\,a_0(\lambda_{k,0}) & \text{if \,$k\geq 0$\,} \\ a(\lambda_{k,0})-\upsilon^{-1}\,a_0(\lambda_{k,0}) & \text{if \,$k< 0$\,} \end{cases} \right\}  \in \ell^2_{0,0}(k)$\,.
\end{itemize}
Then we already have
\begin{equation}
\label{eq:interpolate:a-bestasymp:asymp}
a-\upsilon\,a_0 \in \As_\infty(\C^*,\ell^2_{0},1) \qmq{and} a-\upsilon^{-1}\,a_0 \in \As_0(\C^*,\ell^2_0,1) \;. 
\end{equation}

\emph{Addendum.} If \,$L$\, is a set of holomorphic functions which satisfy (a) and (b) in such a way that the values of \,$\upsilon$\, corresponding to \,$a\in L$\, are bounded and bounded away from zero, that
\,$R=R(\eps)>0$\, in (a) can be chosen uniformly for all \,$a\in L$\,, and that in (b)
there is a uniform bound \,$(z_k)_{k\in \Z} \in \ell^2_{0,0}(k)$\, for the sequences in (b) for all \,$a\in L$\,, then there exists a uniform bounding sequence for the asymptotics
in \eqref{eq:interpolate:a-bestasymp:asymp} that applies for all \,$a\in L$\,. 
\end{cor}

\begin{proof}
Let \,$\lambda_k := \lambda_{k,0}$\, and \,$\mu_k := a(\lambda_{k,0})$\,. Then we have by hypothesis (b)
$$ \lambda_k - \lambda_{k,0} \in \ell^2_{-1,3}(k) \qmq{and} \left. \begin{cases} \mu_k-\upsilon\,\mu_{k,0} & \text{if \,$k\geq 0$\,} \\ \mu_k-\upsilon^{-1}\,\mu_{k,0} & \text{if \,$k<0$\,} \end{cases} \right\} \in \ell^2_{0,0}(k) \; . $$
By Proposition~\ref{P:interpolate:mu}(1) it follows that there exists a holomorphic function \,$\wt{a}:\C^*\to\C$\, that
satisfies
$$ \wt{a}-\upsilon\,a_0 \in \As_\infty(\C^*,\ell^2_0,1) \qmq{and} \wt{a}-\upsilon^{-1}\,a_0 \in \As_0(\C^*,\ell^2_0,1) $$
and \,$\wt{a}(\lambda_k)=\mu_k$\, for all \,$k\in \Z$\,. Because no two of the \,$\lambda_k$\, coincide (i.e.~the function \,$c_0$\,
with zeros at the \,$\lambda_{k,0}$\, has no zeros of higher order), it follows from Proposition~\ref{P:excl:unique}(2) 
that \,$\wt{a}=a$\, holds, and thus the claimed statement follows.

In the setting of the addendum, the hypothesis that the numbers \,$\upsilon$\, corresponding to \,$a\in L$\, are bounded and bounded away from zero ensures that there exists a constant \,$C_1>0$\, so that
\,$|\upsilon-\upsilon^{-1}|\leq C_1$\, holds for all these \,$\upsilon$\,. Let
$$ r_k := C \cdot \left( z_k * \frac{1}{|k|} + \frac{C_1}{|k|} \right) $$
with the constant \,$C>0$\, from Proposition~\ref{P:interpolate:mu}(3)(b). 
By applying that proposition with \,$\lambda_k^{[1]}=\lambda_k^{[2]}=\lambda_{k,0}$\,, \,$\mu_k^{[1]} = a(\lambda_{k,0})$\, and \,$\mu_k^{[2]} = \mu_{k,0}$\,, 
we see that with the sequence \,$(r_k)_{k>0}$\, resp.~\,$(r_k)_{k<0}$\, is a bounding sequence for \,$a-\upsilon\,a_0$\, resp.~for \,$a-\upsilon^{-1}\,a_0$\, (that is independent of \,$a\in L$\,). 
\end{proof}

\section{Final description of the asymptotic of the monodromy}
\label{Se:asympfinal}

Resulting from the interpolation theorems from Section~\ref{Se:interpolate} and by using the asymptotic spaces introduced in Section~\ref{Se:As}, 
we can now describe the asymptotics of the monodromy, its discriminant \,$\Delta^2-4$\, and also of the extended frame
in their final form. Moreover, we now prove the refined asymptotics for the branch points \,$\vkap_{k,\nu}$\, of the spectral curve,
which are analogous to Corollary~\ref{C:asympdiv:asympdiv-neu} for the spectral divisor, and which had been postponed in Section~\ref{Se:asympdiv}.

\begin{Def}
\label{D:asympfinal:monodromy-asymptotic}
We say that a \,$(2\times 2)$-matrix of holomorphic functions \,$\C^*\to \C$\,
$$ M(\lambda) = \left( \begin{matrix} a(\lambda) & b(\lambda) \\ c(\lambda) & d(\lambda) \end{matrix} \right) $$
with \,$\det(M(\lambda))=1$\,
has \emph{non-periodic monodromy asymptotic} or is an \emph{non-periodic (asymptotic) monodromy} with respect to numbers \,$\tau,\upsilon \in \C^*$\, if we have
\begin{itemize}
\item[(a)]
\begin{align*}
a(\lambda)-\upsilon\,a_0(\lambda) & \in \As_\infty(\C^*,\ell^2_{0},1) \\
b(\lambda)-\tau^{-1}\,b_0(\lambda) & \in \As_\infty(\C^*,\ell^2_{1},1) \\
c(\lambda)-\tau\,c_0(\lambda) & \in \As_\infty(\C^*,\ell^2_{-1},1) \\
d(\lambda)-\upsilon^{-1}\,d_0(\lambda) & \in \As_\infty(\C^*,\ell^2_{0},1) 
\end{align*}
\item[(b)]
\begin{align*}
a(\lambda)-\upsilon^{-1}\,a_0(\lambda) & \in \As_0(\C^*,\ell^2_{0},1) \\
b(\lambda)-\tau\,b_0(\lambda) & \in \As_0(\C^*,\ell^2_{-1},1) \\
c(\lambda)-\tau^{-1}\,c_0(\lambda) & \in \As_0(\C^*,\ell^2_{1},1) \\
d(\lambda)-\upsilon\,d_0(\lambda) & \in \As_0(\C^*,\ell^2_{0},1) \;.
\end{align*}
\end{itemize}
If this is the case with \,$\upsilon=1$\,, then we say that \,$M(\lambda)$\, has \emph{monodromy asymptotic} or is an \emph{(asymptotic) monodromy} with respect to 
\,$\tau\in \C^*$\,.

\label{not:asympfinal:Mon}
We denote the space of non-periodic monodromy asymptotic matrices with respect to \,$\tau,\upsilon$\, by \,$\Mon_{\tau,\upsilon}$\,, and put
\,$\Mon_{\tau} := \Mon_{\tau,\upsilon=1}$\,. \,$\Mon_{np} := \bigcup_{\tau,\upsilon\in\C^*} \Mon_{\tau,\upsilon}$\, is the space of all
non-periodic monodromy asymptotic matrices, and \,$\Mon := \bigcup_{\tau\in\C^*} \Mon_\tau$\, is the space of all
monodromy asymptotic matrices.
\end{Def}

Note that a monodromy \,$M(\lambda) \in \Mon_{np}$\, in particular has the properties shown for spectral monodromies in Theorems~\ref{T:asymp:basic} and \ref{T:fasymp:fourier}. 
Therefore Corollary~\ref{C:asympdiv:asympdiv-neu} concerning the asymptotic behavior of a spectral divisor is applicable to the spectral divisor \,$D$\, associated to a monodromy \,$M(\lambda)\in \Mon_{np}$\,.

We next show that the monodromy \,$M(\lambda)$\, associated to a (non-periodic) potential \,$(u,u_y) \in \Pot_{np}$\, indeed has monodromy asymptotic in the sense of Definition~\ref{D:asympfinal:monodromy-asymptotic}.

\newpage

\begin{thm}
\label{T:asympfinal:monodromy}
Let \,$(u,u_y)\in \Pot_{np}$\,,
\,$M(\lambda)$\, be the monodromy associated to \,$(u,u_y)$\,, \,$\Delta(\lambda)$\, its trace function,
and \,$\tau := e^{-(u(0)+u(1))/4}$\, and \,$\upsilon := e^{(u(1)-u(0))/4}$\,.

Then we have:
\begin{enumerate}
\item
\,$M(\lambda) \in \Mon_{\tau,\upsilon}$\,. \\
Moreover, if \,$P$\, is a relatively compact subset of \,$\Pot_{np}$\,, then there exist uniform bounding sequences for the entries of the monodromy that apply to the monodromies of the 
potentials \,$(u,u_y)\in P$\,. 
\item
If \,$(u,u_y)\in\Pot$\, holds (i.e.~\,$(u,u_y)$\, is periodic), then we in fact have \,$M(\lambda)\in \Mon_\tau$\,.
\item
We have the following \emph{trace formula}: If \,$D=\{(\lambda_k,\mu_k)\}$\, is the classical spectral divisor corresponding to \,$M(\lambda)$\,, then we have
$$ e^{(u(0)+u(1))/2} = \prod_{k\in \Z} \frac{\lambda_k}{\lambda_{k,0}} \; . $$
\end{enumerate}
\end{thm}

\begin{rem}
It is tempting to try to replace the condition ``\,$P$\, is a relatively compact subset of \,$\Pot_{np}$\,'' in 
Theorem~\ref{T:asympfinal:monodromy}(1) by ``\,$P$\, is a small closed ball in \,$\Pot_{np}$\,''. 
However with this variant of the condition, the theorem would not be true, as the condition that \,$P$\, 
needs to be relatively compact is ultimately derived from the requirement that the set \,$N\subset L^1([a,b])$\, 
is compact in the version of the Riemann-Lebesgue Lemma described in Lemma~\ref{L:asymp:riemannlebesgue}
(which had been used to prove the basic asymptotic for the monodromy, Theorem~\ref{T:asymp:basic}). It is clear
that one does not obtain a uniform estimate in the Riemann-Lebesgue Lemma if \,$N$\, is replaced by a small ball
in \,$L^1([a,b])$\,, and thus we also cannot relax the condition on \,$P$\, in Theorem~\ref{T:asympfinal:monodromy}(1)
in a similar way.
\end{rem}

\begin{proof}[Proof of Theorem~\ref{T:asympfinal:monodromy}.]
We write
$$ M(\lambda)= \left( \begin{matrix} a(\lambda) & b(\lambda) \\ c(\lambda) & d(\lambda) \end{matrix} \right)
\qmq{and} \Delta(\lambda) = a(\lambda)+d(\lambda) \; . $$

\emph{For (1).}
We note the facts we already know  from Theorems~\ref{T:asymp:basic} and \ref{T:fasymp:fourier} about the asymptotics of the monodromy of \,$(u,u_y)$\,.
From them, the statements on the functions \,$a$\, and \,$d$\, follow from Corollary~\ref{C:interpolate:a-bestasymp},
and the statements on \,$c$\, follow from Corollary~\ref{C:interpolate:c-bestasymp}. The statements on \,$b$\,
follow by applying Corollary~\ref{C:interpolate:c-bestasymp} with \,$-\lambda\cdot b(\lambda)$\, in the place of \,$c(\lambda)$\,.

If \,$P\subset \Pot_{np}$\, is a relatively compact set, then Theorems~\ref{T:asymp:basic} and \ref{T:fasymp:fourier} show
that the entries of the monodromies of \,$(u,u_y)\in P$\, satisfy the hypotheses of the Addendum in 
Corollary~\ref{C:interpolate:c-bestasymp} resp.~Corollary~\ref{C:interpolate:a-bestasymp}. These Addenda
imply the existence of uniform bounding sequences in the sense described in the statement.

\emph{For (2).}
This follows from the fact that the periodicity of \,$(u,u_y)$\, is equivalent to \,$\upsilon=1$\,. 

\emph{For (3).} 
The trace formula follows from calculating the parameter \,$\tau$\, of the asymptotic of the monodromy in two different ways: On one hand, we have
$$ \tau = e^{-(u(0)+u(1))/4} $$
by Theorem~\ref{T:asymp:basic}. On the other hand, we have 
$$ \tau = \pm \left( \prod_{k\in \Z} \frac{\lambda_{k,0}}{\lambda_k} \right)^{1/2} $$
by Proposition~\ref{P:interpolate:lambda}(1). By combining these two equations, we obtain
$$ e^{(u(0)+u(1))/2} = \tau^{-2} = \prod_{k\in \Z} \frac{\lambda_{k}}{\lambda_{k,0}} \; . $$
\end{proof}

\begin{cor}
\label{C:asympfinal:monodromy2}
Let \,$M(\lambda) \in \Mon_{\tau,\upsilon}$\, with \,$\tau,\upsilon \in \C^*$\, be given; this applies for example in the 
setting of Theorem~\ref{T:asympfinal:monodromy} (where \,$M(\lambda)$\, is the monodromy of a non-periodic potential \,$(u,u_y) \in \Pot_{np}$\,). 
We write \,$ M(\lambda)= \left( \begin{smallmatrix} a(\lambda) & b(\lambda) \\ c(\lambda) & d(\lambda) \end{smallmatrix} \right)$\,
and \,$ \Delta(\lambda) = \tr(M(\lambda))=a(\lambda)+d(\lambda)$\,. Then we have
\begin{enumerate}
\item \,$a,d,\Delta \in \As(\C^*,\ell^\infty_{0,0},1)$\,, \,$b \in \As(\C^*,\ell^\infty_{1,-1},1)$\, and \,$c \in \As(\C^*,\ell^\infty_{-1,1},1)$\,.
\item \,$\Delta(\lambda)-\tfrac{\upsilon+\upsilon^{-1}}{2}\,\Delta_0(\lambda) \in \As(\C^*,\ell^2_{0,0},1)$\,. 
\item \,$\Delta^2 - \left( \tfrac{\upsilon+\upsilon^{-1}}{2}\,\Delta_0\right)^2, bc-b_0c_0 \in \As(\C^*,\ell^2_{0,0},2)$\,.
\end{enumerate}
We now suppose that \,$\upsilon=1$\, holds, and let \,$\Sigma$\, be the spectral curve associated to \,$M(\lambda)$\,. Then we also have
\begin{enumerate}
\addtocounter{enumi}{3}
\item \,$\mu-\mu_0, \mu^{-1}-\mu_0^{-1} \in \As(\Sigma,\ell^2_{0,0},1)$\,, where we set \,$\mu_0 := e^{i\,\zeta(\lambda)}$\,. 
\item \,$\mu,\mu^{-1} \in \As(\Sigma,\ell^\infty_{0,0},1)$\,. 
\item \,$\left( \tfrac{\mu-d}{c} - \tfrac{i}{\tau\,\sqrt{\lambda}} \right)|\wh{V}_\delta \in \As_\infty(\wh{V}_\delta,\ell^2_{1},0)$\, 
and \,$\left( \tfrac{\mu-d}{c} - \tfrac{i}{\tau^{-1}\,\sqrt{\lambda}} \right)|\wh{V}_\delta \in \As_0(\wh{V}_\delta,\ell^2_{-1},0)$\,, 
where \,$\sqrt{\lambda}$\, is chosen as the holomorphic function on \,$\wh{V}_\delta$\, with the correct sign.
\end{enumerate}
\end{cor}

\begin{proof}
\emph{For (1).}
These statements follow from Definition~\ref{D:asympfinal:monodromy-asymptotic} 
together with the fact that 
$$ a_0,d_0,\Delta_0 \in \As(\C^*,\ell^\infty_{0,0},1)\;,\quad b_0 \in \As(\C^*,\ell^\infty_{1,-1},1) \qmq{and} c_0 \in \As(\C^*,\ell^\infty_{-1,1},1) $$
holds.

\emph{For (2).}
Because of \,$\Delta=a+d$\, and \,$\Delta_0=2a_0=2d_0$\, we have
\begin{align*}
\Delta(\lambda)-\tfrac{\upsilon+\upsilon^{-1}}{2}\,\Delta_0(\lambda) 
& = \bigr( a(\lambda)-\upsilon\,a_0(\lambda) \bigr) + \bigr( d(\lambda)-\upsilon^{-1}\,d_0(\lambda) \bigr) \\
& = \bigr( a(\lambda)-\upsilon^{-1}\,a_0(\lambda) \bigr) + \bigr( d(\lambda)-\upsilon\,d_0(\lambda) \bigr) \;. 
\end{align*}
The claimed statement follows from this equation
together with the asymptotics for \,$a$\, and \,$d$\, in Definition~\ref{D:asympfinal:monodromy-asymptotic}. 

\emph{For (3).}
This follows from the preceding results by the equations
\begin{align*}
\Delta^2 - \left( \tfrac{\upsilon+\upsilon^{-1}}{2}\,\Delta_0\right)^2
& = \left( \Delta + \tfrac{\upsilon+\upsilon^{-1}}{2}\,\Delta_0\right)\cdot \left( \Delta - \tfrac{\upsilon+\upsilon^{-1}}{2}\,\Delta_0\right) \\ 
\intertext{and}
b\,c-b_0\,c_0 & = (b-\tau^{-1}\,b_0)\cdot c + \tau^{-1}\,b_0 \cdot (c-\tau\,c_0)\\
&  = (b-\tau\,b_0)\cdot c + \tau\,b_0 \cdot (c-\tau^{-1}\,c_0) \;. 
\end{align*}

\emph{For (4).}
We have \,$\mu=\tfrac12(\Delta + \sqrt{\Delta^2-4})$\, and \,$\mu_0 = \tfrac12(\Delta_0 + \sqrt{\Delta_0^2-4})$\,. Because we have \,$\Delta-\Delta_0\in \As(\C^*,\ell^2_{0,0},1)$\, by (2), it remains to show 
that
\begin{equation}
\label{eq:asympfinal:monodromy2:sqrtDelta2-claim}
\left( \sqrt{\Delta^2-4}-\sqrt{\Delta_0^2-4} \right) \in \As(\Sigma,\ell^2_{0,0},1)
\end{equation}
holds. To show this, we use the formula
$$ \sqrt{\Delta^2-4}-\sqrt{\Delta_0^2-4} = \frac{\Delta^2-\Delta_0^2}{\sqrt{\Delta^2-4}+\sqrt{\Delta_0^2-4}} \; . $$
Because both \,$\sqrt{\Delta^2-4}$\, and \,$\sqrt{\Delta_0^2-4}$\, are comparable to \,$w(\lambda)$\, on \,$\wh{V}_\delta$\,,  it follows from (3) that 
$$ \left. \left( \sqrt{\Delta^2-4}-\sqrt{\Delta_0^2-4} \right) \right|\wh{V}_\delta \in \As(\wh{V}_\delta,\ell^2_{0,0},1) $$
holds, whence \eqref{eq:asympfinal:monodromy2:sqrtDelta2-claim} follows by Proposition~\ref{P:interpolate:l2asymp}(1);
this shows \,$\mu-\mu_0 \in \As(\Sigma,\ell^2_{0,0},1)$\,. 

Similarly we obtain \,$\mu^{-1}-\mu_0^{-1} \in \As(\Sigma,\ell^2_{0,0},1)$\, from the equations \,$\mu^{-1}=\tfrac12(\Delta - \sqrt{\Delta^2-4})$\, and \,$\mu_0^{-1} = \tfrac12(\Delta_0 - \sqrt{\Delta_0^2-4})$\,,

\emph{For (5).}
Because of \,$\mu_0,\mu_0^{-1} \in \As(\Sigma,\ell^\infty_{0,0},1)$\,, 
it follows from (4) that \,$\mu,\mu^{-1} \in \As(\Sigma,\ell^\infty_{0,0},1)$\, holds.

\emph{For (6).} We have \,$\tfrac{i}{\sqrt{\lambda}} = \tfrac{\mu_0-d_0}{c_0}$\, and therefore
\begin{equation}
\label{eq:asympfinal:monodromy2:mudc}
\frac{\mu-d}{c}-\frac{i}{\tau\,\sqrt{\lambda}} = \frac{1}{\tau\,c_0}\bigr( \,(\mu-\mu_0) + (d-d_0)\,\bigr) - \frac{\mu-d}{c\cdot \tau\,c_0}\,(c-\tau\,c_0) \; .
\end{equation}
Because both \,$c$\, and \,$c_0$\, are comparable to \,$\sqrt{|\lambda|}\cdot w(\lambda)$\, on \,$\wh{V}_\delta$\,, and we have \,$(\mu-\mu_0)|\wh{V}_\delta \in \As(\wh{V}_\delta,\ell^2_{0,0},1)$\,,
\,$d-d_0 \in \As(\C^*,\ell^2_{0,0},1)$\,, \,$\mu,d \in \As(\Sigma,\ell^\infty_{0,0},1)$\, and \,$c-\tau\,c_0 \in \As_\infty(\C^*,\ell^2_{-1,1},1)$\,, it follows that 
\,$\left( \tfrac{\mu-d}{c}-\tfrac{i}{\tau\,\sqrt{\lambda}} \right)|\wh{V}_\delta \in \As_\infty(\wh{V}_\delta,\ell^2_{1},0)$\, holds.

The other asymptotic assessment \,$\left( \tfrac{\mu-d}{c} - \tfrac{i}{\tau^{-1}\,\sqrt{\lambda}} \right)|\wh{V}_\delta \in \As_0(\wh{V}_\delta,\ell^2_{-1},0)$\, is similarly obtained by replacing
\,$\tau$\, with \,$\tau^{-1}$\, in Equation~\eqref{eq:asympfinal:monodromy2:mudc}.
\end{proof}

In the sequel, we will use the asymptotic behavior of the monodromy described in Definition~\ref{D:asympfinal:monodromy-asymptotic} resp.~Theorem~\ref{T:asympfinal:monodromy} 
as well as the consequences in Corollary~\ref{C:asympfinal:monodromy2} without referencing these statements explicitly every time.

\smallskip

The following proposition concerns the zeros of the discriminant \,$\Delta^2-4$\, of an asymptotic monodromy \,$M(\lambda)\in \Mon$\,. It shows that the zeros of \,$\Delta^2-4$\, (corresponding
to the branch points and singularities of the associated spectral curve \,$\Sigma$\,) have the same kind of asymptotic behavior as was shown for the divisor points in
Corollary~\ref{C:asympdiv:asympdiv-neu}. Note that this statement is true only for periodic asymptotic monodromies (i.e.~for the case \,$\upsilon=1$\,). 

\begin{prop}
\label{P:asympfinal:vkap}
Let \,$M(\lambda) \in \Mon$\, be given (or let \,$M(\lambda)$\, be the monodromy of a potential \,$(u,u_y) \in \Pot$\,), 
\,$\Delta(\lambda)$\, be the trace function of \,$M(\lambda)$\,.
We consider the two sequences
\,$(\vkap_{k,1})$\, and \,$(\vkap_{k,2})$\, of zeros of \,$\Delta^2-4$\, (i.e.~of branch points and singularities of
the spectral curve \,$\Sigma$\, associated to \,$\Delta$\,) defined in Proposition~\ref{P:excl:basic}(1). 

Then we have \,$\vkap_{k,\nu}-\lambda_{k,0} \in \ell^2_{-1,3}(k)$\, for \,$\nu\in \{1,2\}$\,. 
\end{prop}

\begin{proof}
By Corollary~\ref{C:asympfinal:monodromy2}(2) we have \,$\Delta-\Delta_0 \in \As(\C^*,\ell^2_{0,0},1)$\,. We let \,$(\Sigma,\calD)$\, be the spectral
data of \,$M(\lambda)$\,. By Proposition~\ref{P:excl:basic}(2),(3) and Corollary~\ref{C:asympdiv:asympdiv-neu}, the points in the support of \,$\calD$\, are enumerated by a sequence \,$(\lambda_k,\mu_k)_{k\in \Z}$\,
on \,$\Sigma$\,, such that \,$\lambda_k-\lambda_{k,0}\in \ell^2_{-1,3}(k)$\, and \,$\mu_k-\mu_{k,0}\in \ell^2_{0,0}(k)$\, holds. 

Using the fact that \,$\Delta_0(\lambda_{k,0})=2(-1)^k=\Delta(\vkap_{k,\nu})$\, holds,
we write 
\begin{align}
& \Delta_0(\vkap_{k,\nu})-\Delta_0(\lambda_{k,0}) \notag \\
=\; & \Delta_0(\vkap_{k,\nu})-\Delta(\vkap_{k,\nu}) \notag \\
\label{eq:monore:direct:Delta-Delta0}
=\; & \bigr(\Delta_0(\lambda_k)-\Delta_0(\lambda_{k,0})\bigr) + \bigr( \Delta_0(\lambda_{k,0}) - \Delta(\lambda_k) \bigr) + \int_{\lambda_k}^{\vkap_{k,\nu}} (\Delta_0'-\Delta')(\lambda)\,\mathrm{d}\lambda \; . 
\end{align}
To evaluate the terms involved in this equation, we first note  the following Taylor expansion for \,$\Delta_0(\lambda) = 2\cos(\zeta(\lambda))$\, for \,$\lambda\in U_{k,\delta}$\,:
\begin{equation}
\label{eq:monore:direct:taylor}
\Delta_0(\lambda) = \Delta_0(\lambda_{k,0}) + (-1)^k\,(\zeta(\lambda)-\zeta(\lambda_{k,0}))^2 + O\bigr( (\zeta(\lambda)-\zeta(\lambda_{k,0}))^4 \bigr) \; .
\end{equation}
As a consequence, we see that  we have for the left hand side of \eqref{eq:monore:direct:Delta-Delta0}
\begin{equation}
\label{eq:monore:direct:vkap-lhs}
\Delta_0(\vkap_{k,\nu})-\Delta_0(\lambda_{k,0}) = (-1)^k\cdot (\zeta(\vkap_{k,\nu})-\zeta(\lambda_{k,0}))^2 \cdot (1+O((\zeta(\vkap_{k,\nu})-\zeta(\lambda_{k,0}))^2)) \; . 
\end{equation}

We now investigate the three summands on the right hand side of \eqref{eq:monore:direct:Delta-Delta0} individually.

First, we note that because of \,$\lambda_k-\lambda_{k,0} \in \ell^2_{-1,3}(k)$\, we have \,$\zeta(\lambda_k)-\zeta(\lambda_{k,0})\in \ell^2_{0,0}(k)$\, by Proposition~\ref{P:vac2:excldom-new}(1), 
whence it follows by the Taylor expansion
Equation~\eqref{eq:monore:direct:taylor} that
\begin{equation}
\label{eq:monore:direct:vkap-rhs1}
\Delta_0(\lambda_k)-\Delta_0(\lambda_{k,0}) \in \ell^1_{0,0}(k)
\end{equation}
holds.

Second, we have
\begin{align*}
\Delta(\lambda_k) - \Delta_0(\lambda_{k,0})
& = \mu_k + \mu_k^{-1} - 2\,(-1)^k = (\mu_k^{1/2}-(-1)^k\,\mu_k^{-1/2})^2 \notag \\
& = \left( \frac{\mu_k-\mu_k^{-1}}{\mu_k^{1/2}+(-1)^k\,\mu_k^{-1/2}} \right)^2 = \left( \frac{(\mu_k-\mu_{k,0})-(\mu_k^{-1}-\mu_{k,0}^{-1})}{\mu_k^{1/2}+(-1)^k\,\mu_k^{-1/2}} \right)^2 \; . 
\end{align*}
We have \,$\mu_k-\mu_{k,0} \in \ell^2_{0,0}(k)$\,, and therefore also \,$\mu_k^{-1}-\mu_{k,0}^{-1} \in \ell^2_{0,0}(k)$\, (because \,$z\mapsto z^{-1}$\, is locally Lipschitz continuous near \,$z=\pm 1$\,). 
Because the denominator \,$\mu_k^{1/2}+(-1)^k\,\mu_k^{-1/2}$\, in the expression above is bounded away from zero, it follows that we have
\begin{equation}
\label{eq:monore:direct:vkap-rhs2}
\Delta(\lambda_k)-\Delta_0(\lambda_{k,0}) \in \ell^1_{0,0}(k) \; . 
\end{equation}

Third, because of \,$\Delta-\Delta_0 \in \As(\C^*,\ell^2_{0,0},1)$\, (Corollary~\ref{C:asympfinal:monodromy2}(2)),
we have \,$\Delta'-\Delta_0'\in \As(\C^*,\ell^2_{1,-3},1)$\, by Proposition~\ref{P:interpolate:l2asymp}(3), and thus we obtain
\begin{align}
\left| \int_{\lambda_k}^{\vkap_{k,\nu}} (\Delta_0'-\Delta')(\lambda)\,\mathrm{d}\lambda \right|
& \leq |\vkap_{k,\nu}-\lambda_k| \cdot \ell^2_{1,-3}(k) = |\zeta(\vkap_{k,\nu})-\zeta(\lambda_k)| \cdot \ell^2_{0,0}(k) \notag \\
& \leq \bigr( |\zeta(\vkap_{k,\nu})-\zeta(\lambda_{k,0})| + \underbrace{|\zeta(\lambda_k)-\zeta(\lambda_{k,0})|}_{\in \ell^2_{0,0}(k)} \bigr) \cdot \ell^2_{0,0}(k) \notag \\
\label{eq:monore:direct:vkap-rhs3}
& = |\zeta(\vkap_{k,\nu})-\zeta(\lambda_{k,0})| \cdot \ell^2_{0,0}(k) + \ell^1_{0,0}(k) \; . 
\end{align}

By plugging the Equations~\eqref{eq:monore:direct:vkap-lhs}, \eqref{eq:monore:direct:vkap-rhs1}, \eqref{eq:monore:direct:vkap-rhs2} and \eqref{eq:monore:direct:vkap-rhs3} into Equation~\eqref{eq:monore:direct:Delta-Delta0}, we obtain
that there exist sequences \,$a_k \in \ell^1_{0,0}(k)$\, and \,$b_k \in \ell^2_{0,0}(k)$\, so that with \,$r_k := |\zeta(\vkap_{k,\nu})-\zeta(\lambda_{k,0})|$\, we have
$$ r_k^2 = a_k + b_k\cdot r_k $$
and therefore
$$ r_k = \frac{b_k}{2} \pm \sqrt{\frac{b_k^2}{4}-a_k} \in \ell^2_{0,0}(k) \; . $$
It follows that \,$\vkap_{k,\nu}-\lambda_{k,0} \in \ell^2_{-1,3}(k)$\, holds.
\end{proof}

In the next proposition, the discriminant function \,$\Delta^2-4$\, of a monodromy \,$M(\lambda)\in \Mon$\, is studied. Proposition~\ref{P:asympfinal:Delta2}(1),(2) gives a representation of \,$\Delta^2-4$\, 
analogous to the one of the monodromy function \,$c$\, in Proposition~\ref{P:interpolate:lambda}; more specifically, the infinite product in Proposition~\ref{P:asympfinal:Delta2}(1) is 
analogous to Equation~\eqref{eq:interpolate:lambda:c}, and the formula in Proposition~\ref{P:asympfinal:Delta2}(2) is analogous to the trace formula in Equation~\eqref{eq:interpolate:lambda:tau}. 
Proposition~\ref{P:asympfinal:Delta2}(3) gives a more detailed asymptotic description of \,$\Delta^2-4$\, on the excluded domains than is provided by Corollary~\ref{C:asympfinal:monodromy2}(3);
we will use this asymptotic description in Section~\ref{Se:jacobiprep} to prove certain estimates in preparation for the construction of the Jacobi variety for the spectral curve \,$\Sigma$\,.  

\begin{prop}
\label{P:asympfinal:Delta2}
Let \,$M(\lambda) \in \Mon$\, and \,$\Delta := \tr(M(\lambda))$\,; we enumerate the zeros \,$\vkap_{k,\nu}$\, of \,$\Delta^2-4$\, as in Proposition~\ref{P:excl:basic}(1).

\begin{enumerate}
\item We have
$$ \Delta(\lambda)^2-4 = - \frac{1}{4\lambda} \cdot (\lambda-\vkap_{0,1})(\lambda-\vkap_{0,2}) \cdot \prod_{k=1}^\infty \frac{(\lambda-\vkap_{k,1})\,(\lambda-\vkap_{k,2})}{(16\,\pi^2\,k^2)^2} \cdot \prod_{k=1}^\infty \frac{(\lambda-\vkap_{-k,1})\,(\lambda-\vkap_{-k,2})}{\lambda^2} \; . $$
\item We have
$$ \prod_{k\in \Z} \frac{\vkap_{k,1}\cdot \vkap_{k,2}}{\lambda_{k,0}^2} = 1 \; . $$
\item There exists a \,$C>0$\, and a sequence \,$a_k \in \ell^{2}_{0,0}(k)$\, such that for all \,$k>0$\, and all \,$\lambda\in U_{k,\delta}$\, we have
$$ \left| \lambda_{k,0}\cdot \frac{\Delta(\lambda)^2-4}{(\lambda-\vkap_{k,1})\cdot(\lambda-\vkap_{k,2})} - \left( -\frac{1}{16} \right) \right| \leq \frac{C}{k}\cdot \left( |\lambda-\vkap_{k,1}| + |\lambda-\vkap_{k,2}| \right) + a_k \;, $$
whereas for \,$k<0$\, and \,$\lambda\in U_{k,\delta}$\, we have
$$ \left| \lambda_{k,0}^3\cdot \frac{\Delta(\lambda)^2-4}{(\lambda-\vkap_{k,1})\cdot(\lambda-\vkap_{k,2})} - \left( -\frac{1}{16} \right) \right| \leq C\,k^3\,\cdot \left( |\lambda-\vkap_{k,1}| + |\lambda-\vkap_{k,2}| \right) + a_k \; . $$
\end{enumerate}
\end{prop}

\begin{proof}
\emph{For (1) and (2).}
Let us put \,$f_0(\lambda) := \Delta_0(\lambda)^2-4=-4\,\sin(\zeta(\lambda))^2$\, and \,$f(\lambda) := \Delta(\lambda)^2-4$\,. 
Then we have \,$f-f_0 \in \As(\C^*,\ell^2_{0,0},2)$\, by Corollary~\ref{C:asympfinal:monodromy2}(3).

Moreover we put
\begin{align*}
\wt{f}(\lambda) & := - \frac{1}{4\lambda} \cdot (\lambda-\vkap_{0,1})(\lambda-\vkap_{0,2}) \cdot \prod_{k=1}^\infty \frac{(\lambda-\vkap_{k,1})\,(\lambda-\vkap_{k,2})}{(16\,\pi^2\,k^2)^2} \cdot \prod_{k=1}^\infty \frac{(\lambda-\vkap_{-k,1})\,(\lambda-\vkap_{-k,2})}{\lambda^2} \\
& = -\frac{4}{\lambda} \cdot \frac{1}{\tau_1\cdot \tau_2}\cdot c_1(\lambda)\cdot c_2(\lambda) 
\end{align*}
with
$$ c_\nu(\lambda) := \frac{1}{4}\,\tau_\nu\,(\lambda-\vkap_{0,\nu})\,\prod_{k=1}^\infty \frac{\vkap_{k,\nu}-\lambda}{16\,\pi^2\,k^2}\,\prod_{k=1}^\infty \frac{\lambda-\vkap_{-k,\nu}}{\lambda} $$
and 
$$ \tau_\nu := \left( \prod_{k\in \Z} \frac{\lambda_{k,0}}{\vkap_{k,\nu}} \right)^{1/2} $$
for \,$\nu\in\{1,2\}$\,. 

By Proposition~\ref{P:interpolate:lambda}(1) we have 
$$ c_\nu-\tau_\nu\,c_0 \in \As_\infty(\C^*,\ell^2_{-1},1) \qmq{and} c_\nu-\tau_\nu^{-1}\,c_0 \in \As_0(\C^*,\ell^2_1,1) $$
and therefore
$$ \wt{f}-f_0 \in \As_\infty(\C^*,\ell^2_0,2) \qmq{and} \wt{f}-\tau_1^{-2}\,\tau_2^{-2}\,f_0 \in \As_0(\C^*,\ell^2_0,2) \; . $$

The holomorphic functions \,$f$\, and \,$\wt{f}$\, have the same zeros with multiplicities, and by the preceding asymptotic estimates they satisfy the hypotheses of Proposition~\ref{P:excl:unique}(3)
with \,$\tau_\pm =1=\wt{\tau}_+$\, and \,$\wt{\tau}_- = \tau_1^{-2}\,\tau_2^{-2}$\,. It follows from that proposition that there exists a constant \,$A\in\C^*$\, with \,$f=A\cdot \wt{f}$\,, \,$1=A\cdot 1$\,
and \,$1=A\cdot \tau_1^{-1}\,\tau_2^{-2}$\,. Thus we have \,$\wt{f}=f$\, and 
$$ 1 = A = \tau_1^2\,\tau_2^2 = \prod_{k\in \Z} \frac{\lambda_{k,0}^2}{\vkap_{k,1}\,\vkap_{k,2}} \;. $$

\emph{For (3).}
We first consider the case \,$k>0$\,. 
By Corollary~\ref{C:interpolate:cdivlin}(1), there exists for \,$\nu\in\{1,2\}$\, a constant \,$C_\nu>0$\, and a sequence \,$a_k^{[\nu]} \in \ell^{2}_{0,-2}(k)$\, such that we have for all \,$k>0$\, and all \,$\lambda\in U_{k,\delta}$\, 
$$ \left| \frac{c_\nu(\lambda)}{\tau_\nu^{-1}\cdot (\lambda-\vkap_{k,\nu})} - \frac{(-1)^k}{8} \right| \leq C_\nu \frac{|\lambda-\vkap_{k,\nu}|}{k} + a_k^{[\nu]} \; , $$
and therefore 
\begin{align}
& \left| \frac{c_1(\lambda)\cdot c_2(\lambda)}{\tau_1\,\tau_2\,(\lambda-\vkap_{k,1})\,(\lambda-\vkap_{k,2})} - \frac{1}{64} \right| \notag \\ 
\leq & \left| \frac{c_1(\lambda)}{\tau_1\,(\lambda-\vkap_{k,1})} \right| \cdot \left| \frac{c_2(\lambda)}{\tau_2\,(\lambda-\vkap_{k,2})} - \frac{(-1)^k}{8} \right| + \left| \frac{c_1(\lambda)}{\tau_1\,(\lambda-\vkap_{k,1})} - \frac{(-1)^k}{8} \right| \cdot \left| \frac{(-1)^k}{8} \right| \notag \\
\label{eq:transdgl:Delta2:asymp-c1c2}
\leq & C\cdot \left( \frac{|\lambda-\vkap_{k,1}|}{k} + \frac{|\lambda-\vkap_{k,2}|}{k} \right) + a_k
\end{align}
with \,$C>0$\, and \,$a_k \in \ell^{2}_0(k>0)$\,. 
Furthermore, we have
\begin{equation}
\label{eq:transdgl:Delta2:asymp-lambda}
\left|\frac{\lambda_{k,0}}{\lambda} - 1 \right| = \frac{|\lambda_{k,0}-\lambda|}{|\lambda|} = O(k^{-1}) \; .
\end{equation}

By (1), we have
$$ \Delta(\lambda)^2-4 = -\frac{4}{\lambda} \cdot \frac{c_1(\lambda)\cdot c_2(\lambda)}{\tau_1\,\tau_2} $$
and therefore by the estimates Equations~\eqref{eq:transdgl:Delta2:asymp-c1c2} and \eqref{eq:transdgl:Delta2:asymp-lambda} we obtain

\begin{align*}
& \left| \lambda_{k,0}\cdot \frac{\Delta(\lambda)^2-4}{(\lambda-\vkap_{k,1})\cdot(\lambda-\vkap_{k,2})} - \left( -\frac{1}{16} \right) \right| 
= \left| 4\,\frac{\lambda_{k,0}}{\lambda}\,\frac{c_1(\lambda)\cdot c_2(\lambda)}{\tau_1\,\tau_2\,(\lambda-\vkap_{k,1})\,(\lambda-\vkap_{k,2})} -\frac{1}{16} \right| \\
\leq & \underbrace{\left| \frac{\lambda_{k,0}}{\lambda}-1 \right|}_{=O(k^{-1})} \cdot \underbrace{\left| 4\,\frac{c_1(\lambda)\cdot c_2(\lambda)}{\tau_1\,\tau_2\,(\lambda-\vkap_{k,1})\,(\lambda-\vkap_{k,2})} \right|}_{=O(1)} 
+ 4\,\underbrace{\left| \frac{c_1(\lambda)\cdot c_2(\lambda)}{\tau_1\,\tau_2\,(\lambda-\vkap_{k,1})\,(\lambda-\vkap_{k,2})} -\frac{1}{64} \right|}_{\leq C\cdot \left( \frac{|\lambda-\vkap_{k,1}|}{k} + \frac{|\lambda-\vkap_{k,2}|}{k} \right) + a_k} \\
\leq & \frac{C'}{k}\cdot \left( |\lambda-\vkap_{k,1}| + |\lambda-\vkap_{k,2}| \right) + a_k'
\end{align*}
with a new \,$C'>0$\, and a new sequence \,$a_k' \in \ell^{2}_0(k>0)$\,.
This proves (3) for the case \,$k>0$\,. 

In the case \,$k<0$\, we proceed similarly. The estimate analogous to Equation~\eqref{eq:transdgl:Delta2:asymp-c1c2} is
\begin{align*}
\left| \frac{c_1(\lambda)\cdot c_2(\lambda)}{\tau_1^{-1}\,\tau_2^{-1}\,(\lambda-\vkap_{k,1})\,(\lambda-\vkap_{k,2})} - \frac{1}{64}\,\lambda_{k,0}^{-2} \right| & \leq C\cdot \left( |\lambda-\vkap_{k,1}| + |\lambda-\vkap_{k,2}| \right)\cdot k^7 + a_k 
\end{align*}
with \,$C>0$\, and \,$a_k \in \ell^{2}_{-4}(k)$\,, and the estimate Equation~\eqref{eq:transdgl:Delta2:asymp-lambda} also holds in the case \,$k<0$\,. 

We thus obtain (note that we have \,$\tau_1\,\tau_2 = \tau_1^{-1}\,\tau_2^{-1}$\, by (2))
\begin{align*}
& \left| \lambda_{k,0}^3\cdot \frac{\Delta(\lambda)^2-4}{(\lambda-\vkap_{k,1})\cdot(\lambda-\vkap_{k,2})} - \left( -\frac{1}{16} \right) \right| 
= \left| 4\,\frac{\lambda_{k,0}^3}{\lambda}\,\frac{c_1(\lambda)\cdot c_2(\lambda)}{\tau_1\,\tau_2\,(\lambda-\vkap_{k,1})\,(\lambda-\vkap_{k,2})} -\frac{1}{16} \right| \\
\leq & \underbrace{\left| \frac{\lambda_{k,0}}{\lambda}-1 \right|}_{=O(k^{-1})} \cdot \underbrace{\left| 4\,\lambda_{k,0}^2\,\frac{c_1(\lambda)\cdot c_2(\lambda)}{\tau_1\,\tau_2\,(\lambda-\vkap_{k,1})\,(\lambda-\vkap_{k,2})} \right|}_{=O(1)} \\
& \qquad \qquad + \underbrace{4\,\lambda_{k,0}^2}_{=O(k^{-4})}\,\underbrace{\left| \frac{c_1(\lambda)\cdot c_2(\lambda)}{\tau_1^{-1}\,\tau_2^{-1}\,(\lambda-\vkap_{k,1})\,(\lambda-\vkap_{k,2})} -\frac{1}{64}\,\lambda_{k,0}^{-2} \right|}_{\leq C\cdot \left( |\lambda-\vkap_{k,1}| + |\lambda-\vkap_{k,2}| \right)\cdot k^7 + a_k } \\
\leq & C'\,k^3\,\cdot \left( |\lambda-\vkap_{k,1}| + |\lambda-\vkap_{k,2}| \right) + a_k'
\end{align*}
with new constants \,$C'>0$\, and \,$a_k' \in \ell^{2}_0(k<0)$\,. This proves also the case \,$k<0$\,. 
\end{proof}

Finally we consider the extended frame \,$F(x,\lambda):=F_\lambda(x)$\, associated to a potential \,$(u,u_y)\in \Pot_{np}$\,, see Section~\ref{Se:minimal}.  For the construction of the Darboux coordinates on the space of 
(periodic) potentials in Section~\ref{Se:darboux}, we need an asymptotic estimate for \,$F(x,\lambda)$\, that is uniform in \,$x\in [0,1]$\,,
and which is analogous to the asymptotics we have for the monodromy \,$M(\lambda)=F(1,\lambda)$\,. The following proposition provides the necessary result:

\begin{prop}
\label{P:asympfinal:frame}
Let \,$(u,u_y)\in \Pot_{np}$\, and
\,$F(x,\lambda) := \left( \begin{smallmatrix} a(x,\lambda) & b(x,\lambda) \\ c(x,\lambda) & d(x,\lambda) \end{smallmatrix} \right)$\, 
be the extended frame associated to \,$(u,u_y)$\,. We also consider the extended frame of the vacuum
$$ F_0(x,\lambda) := \begin{pmatrix} a_0(x,\lambda) & b_0(x,\lambda) \\ c_0(x,\lambda) & d_0(x,\lambda) \end{pmatrix} := \begin{pmatrix}
\cos(x\,\zeta(\lambda)) & -\lambda^{-1/2}\,\sin(x\,\zeta(\lambda)) \\ \lambda^{1/2}\,\sin(x\,\zeta(\lambda)) & \cos(x\,\zeta(\lambda)) \end{pmatrix} $$
and put \,$\tau(x) := e^{-(u(0)+u(x))/4}$\,, \,$\upsilon(x) := e^{(u(x)-u(0))/4}$\,. 

Then for every \,$x\in [0,1]$\, we have 
\begin{itemize}
\item[(1)]
\begin{align*}
a(x,\lambda)-\upsilon(x)\,a_0(x,\lambda) & \in \As_\infty(\C^*,\ell^2_{0},1) \\
b(x,\lambda)-\tau(x)^{-1}\,b_0(x,\lambda) & \in \As_\infty(\C^*,\ell^2_{1},1) \\
c(x,\lambda)-\tau(x)\,c_0(x,\lambda) & \in \As_\infty(\C^*,\ell^2_{-1},1) \\
d(x,\lambda)-\upsilon(x)^{-1}\,d_0(x,\lambda) & \in \As_\infty(\C^*,\ell^2_{0},1) 
\end{align*}
\item[(2)]
\begin{align*}
a(x,\lambda)-\upsilon(x)^{-1}\,a_0(x,\lambda) & \in \As_0(\C^*,\ell^2_{0},1) \\
b(x,\lambda)-\tau(x)\,b_0(x,\lambda) & \in \As_0(\C^*,\ell^2_{-1},1) \\
c(x,\lambda)-\tau(x)^{-1}\,c_0(x,\lambda) & \in \As_0(\C^*,\ell^2_{1},1) \\
d(x,\lambda)-\upsilon(x)\,d_0(x,\lambda) & \in \As_0(\C^*,\ell^2_{0},1) \;,
\end{align*}
\end{itemize}
and these asymptotic estimates are uniform in the sense that there
exist bounding sequences that are suitable for all \,$x\in [0,1]$\,. 
\end{prop}


\begin{proof}
For given \,$x_0 \in [0,1]$\,, we consider the functions \,$\wt{u},\wt{u}_y$\, given for \,$t\in [0,1]$\, by
$$ \wt{u}(t) = \begin{cases} u(t) & \text{for \,$0\leq t \leq x_0$\,} \\ u(x_0) & \text{for \,$x_0 < t\leq 1$\,} \end{cases}
\qmq{and}
\wt{u}_y(t) = \begin{cases} u_y(t) & \text{for \,$0\leq t \leq x_0$\,} \\ 0 & \text{for \,$x_0 < t\leq 1$\,} \end{cases} \; . $$
Then we have \,$(\wt{u},\wt{u}_y)\in \Pot_{np}$\,.

Let us denote by \,$\wt{F}(t,\lambda)$\, the extended frame of the potential \,$(\wt{u},\wt{u}_y)$\,, and let us denote
by \,$\wh{F}(t,\lambda)$\, the extended frame of the constant potential \,$(\wh{u}=u(x_0),\wh{u}_y=0)\in\Pot$\,. 
(For the geometric meaning of this constant potential, refer to Remark~\ref{R:excl:constpot}.)
Then we have for \,$t\in[0,1]$\,
$$ \wt{F}(t,\lambda)  = \begin{cases}
F(t,\lambda) & \text{for \,$0\leq t \leq x_0$\,} \\ \wh{F}(t-x_0,\lambda)\cdot F(x_0,\lambda) & \text{for \,$x_0 < t \leq 1$\,} 
\end{cases} \;, $$
in particular
$$ F(x_0,\lambda) = \wh{F}(1-x_0,\lambda)^{-1} \cdot \wt{F}(1,\lambda) \; . $$

By Theorem~\ref{T:asympfinal:monodromy}, we have \,$\wt{F}(1,\lambda) = \wt{M}(\lambda) \in \Mon_{\tau=e^{-(u(0)+u(x_0))/4},\upsilon=e^{(u(x_0)-u(0))/4}}$\,, and because 
the set \,$P$\, of potentials \,$(\wt{u},\wt{u}_y)$\,, where \,$x_0$\, runs through \,$[0,1]$\,, is relatively compact in \,$\Pot$\,,
there exists bounding sequences corresponding to this asymptotic that are independent of \,$x_0\in[0,1]$\,.

Concerning the extended frame \,$\wh{F}(t,\lambda)$\, of the constant potential, we note that the corresponding \,$1$-form 
$$ \wh{\alpha}_\lambda = \frac{1}{4}\,\begin{pmatrix} 0 & -e^{u(x_0)/2}-\lambda^{-1}\,e^{-u(x_0)/2} \\ e^{u(x_0)/2}+\lambda\,e^{-u(x_0)/2} & 0 \end{pmatrix}\mathrm{d}x $$
is independent of \,$t$\,, and therefore can be calculated explicitly as
\begin{align*}
\wh{F}(t,\lambda) & = \exp(t\cdot \wh{\alpha}_\lambda) 
= \begin{pmatrix} \cos(t\cdot \xi(\lambda)) & -\beta(\lambda)^{-1}\cdot \sin(t\cdot \xi(\lambda)) \\ \beta(\lambda) \cdot \sin(t\cdot \xi(\lambda)) & \cos(t\cdot \xi(\lambda)) \end{pmatrix}  
\end{align*}
with
\begin{align*}
\xi(\lambda) & := \frac14\,\sqrt{(\lambda\,e^{-u(x_0)/2}+e^{u(x_0)/2})\cdot(\lambda^{-1}\,e^{-u(x_0)/2}+e^{u(x_0)/2})}  \\
\qmq{and} \beta(\lambda) & :=  \sqrt{\frac{e^{u(x_0)/2}+\lambda\,e^{-u(x_0)/2}}{e^{u(x_0)/2}+\lambda^{-1}\,e^{-u(x_0)/2}}}
\end{align*}
(compare the proof of Proposition~\ref{P:excl:basic}(2)). An explicit calculation shows that \,$\wt{F}(t,\lambda)$\, is asymptotically close to
$$ \begin{pmatrix} a_0(t,\lambda) & e^{u(x_0)/2}\,b_0(t,\lambda) \\ e^{-u(x_0)/2}\,c_0(t,\lambda) & d_0(t,\lambda) \end{pmatrix}
\qmq{resp.~to} \begin{pmatrix} a_0(t,\lambda) & e^{-u(x_0)/2}\,b_0(t,\lambda) \\ e^{u(x_0)/2}\,c_0(t,\lambda) & d_0(t,\lambda) \end{pmatrix} $$
for \,$\lambda\to\infty$\, resp.~for \,$\lambda\to 0$\,. 

It follows that for \,$\lambda\to\infty$\,, \,$F(x_0,\lambda)=\wh{F}(1-x_0,\lambda)^{-1}\cdot \wt{F}(1,\lambda)$\, is for \,$\lambda\to\infty$\, asymptotically close to
(where we abbreviate \,$\zeta := \zeta(\lambda)$\,)
{\tiny
\begin{align*}
& \begin{pmatrix} a_0(1-x_0,\lambda) & e^{u(x_0)/2}\,b_0(1-x_0,\lambda) \\ e^{-u(x_0)/2}\,c_0(1-x_0,\lambda) & d_0(1-x_0,\lambda) \end{pmatrix}^{-1} 
\cdot \begin{pmatrix} e^{(u(x_0)-u(0))/4}\,a_0(1,\lambda) & e^{(u(0)+u(x_0))/4}\,b_0(1,\lambda) \\ e^{-(u(0)+u(x_0))/4}\,c_0(1,\lambda) & e^{-(u(x_0)-u(0))/4}\,d_0(1,\lambda) \end{pmatrix} \\
=\; & \begin{pmatrix} d_0(1-x_0,\lambda) & -e^{u(x_0)/2}\,b_0(1-x_0,\lambda) \\ -e^{-u(x_0)/2}\,c_0(1-x_0,\lambda) & a_0(1-x_0,\lambda) \end{pmatrix}
\cdot \begin{pmatrix} e^{(u(x_0)-u(0))/4}\,a_0(1,\lambda) & e^{(u(0)+u(x_0))/4}\,b_0(1,\lambda) \\ e^{-(u(0)+u(x_0))/4}\,c_0(1,\lambda) & e^{-(u(x_0)-u(0))/4}\,d_0(1,\lambda) \end{pmatrix} \\
=\; & \left( \begin{matrix} e^{(u(x_0)-u(0))/4}\,\cos((1-x_0)\zeta)\,\cos(\zeta) + e^{(u(x_0)-u(0))/4}\,\sin((1-x_0)\,\zeta)\,\sin(\zeta) \\
-e^{-(u(0)+u(x_0))/4}\,\lambda^{1/2}\,\sin((1-x_0)\,\zeta)\,\cos(\zeta) + e^{-(u(0)+u(x_0))/4}\,\cos((1-x_0)\,\zeta)\,\lambda^{1/2}\,\sin(\zeta) \end{matrix} \right. \\
& \qquad\qquad\qquad \left. \begin{matrix} -e^{(u(0)+u(x_0))/4}\,\cos((1-x_0)\,\zeta)\,\lambda^{-1/2}\,\sin(\zeta) + e^{(u(x_0)+u(0))/4}\,\lambda^{-1/2}\,\sin((1-x_0)\,\zeta)\,\cos(\zeta) \\
e^{(u(0)-u(x_0))/4}\,\lambda^{1/2}\,\sin((1-x_0)\,\zeta)\,\lambda^{-1/2}\,\sin(\zeta) + e^{-(u(x_0)-u(0))/4}\,\cos((1-x_0)\,\zeta)\,\cos(\zeta) \end{matrix} \right) \\
=\; & \begin{pmatrix} e^{(u(x_0)-u(0))/4}\,\cos(x_0\,\zeta) & -e^{(u(0)+u(x_0))/4}\,\lambda^{-1/2}\,\sin(x_0\,\zeta) \\
e^{-(u(0)+u(x_0))/4}\,\lambda^{1/2}\,\sin(x_0\,\zeta) & e^{-(u(x_0)-u(0))/4}\,\cos(x_0\,\zeta) \end{pmatrix} \\
=\; & \begin{pmatrix} \upsilon(x_0)\,a_0(x_0,\lambda) & \tau(x_0)^{-1}\,b_0(x_0,\lambda) \\ \tau(x_0)\,c_0(x_0,\lambda) & \upsilon(x_0)^{-1}\,d_0(x_0,\lambda) \end{pmatrix} \; . 
\end{align*}

}

An analogous calculation shows that for \,$x\to 0$\,, \,$F(x_0,\lambda)$\, is asymptotically close to
$$ \begin{pmatrix} \upsilon(x_0)^{-1}\,a_0(x_0,\lambda) & \tau(x_0)\,b_0(x_0,\lambda) \\ \tau(x_0)^{-1}\,c_0(x_0,\lambda) & \upsilon(x_0)\,d_0(x_0,\lambda) \end{pmatrix} \; . $$
\end{proof}

\section{Non-special divisors and the inverse problem for \,$\Mon\to\Div$\,}
\label{Se:special}

In this section we characterize the pairs \,$(\Sigma,\mathcal{D})$\, which uniquely define a monodromy \,$M(\lambda) \in \Mon$\, such that \,$(\Sigma,\mathcal{D})$\, is the spectral data of \,$M(\lambda)$\,.
Remember that in Section~\ref{Se:spectrum}, we associated to any \,$M(\lambda)\in \Mon$\,
(henceforth, we will no longer consider non-periodically asymptotic objects, i.e.~from here on we will always suppose \,$\upsilon=1$\,)
the trace function \,$\Delta(\lambda) := \tr(M(\lambda))$\, and thereby the spectral curve 
\begin{equation}
\label{eq:special:spectralcurve}
\Sigma := \Menge{(\lambda,\mu)\in\C^*\times \C}{\mu^2-\Delta(\lambda)\cdot \mu+1=0} \;,
\end{equation}
and moreover the spectral divisor, i.e.~the generalized divisor \,$\calD$\, on \,$\Sigma$\, generated by \,$1$\, and \,$\tfrac{\mu-d}{c}$\,. \,$\calD$\, is positive and asymptotic (see Definition~\ref{D:asympdiv:asympdiv}). 

We now suppose that any pair \,$(\Sigma,\calD)$\, is given, and ask if there exists a monodromy \,$M(\lambda)\in\Mon$\, such that \,$(\Sigma,\calD)$\, is the spectral data of \,$M(\lambda)$\,.
Here we suppose that \,$\Sigma$\, is a spectral curve ``of the type we consider'', i.e.~that there exists a holomorphic function \,$\Delta:\C^*\to\C$\, with \,$\Delta-\Delta_0\in \As(\C^*,\ell^2_{0,0},1)$\, 
so that the complex curve \,$\Sigma$\, is given by Equation~\eqref{eq:special:spectralcurve}. Moreover, we suppose that \,$\calD$\, is a generalized divisor on \,$\Sigma$\,,
i.e.~a subsheaf of the sheaf of meromorphic functions on \,$\Sigma$\, that is finitely generated over the sheaf \,$\calO$\, of holomorphic functions on \,$\Sigma$\,; moreover we require
\,$\calD$\, to be positive (i.e.~\,$\calO \subset \calD$\,) and asymptotic in the sense of Definition~\ref{D:asympdiv:asympdiv}. 

However it turns out that  two further conditions need
to be imposed on \,$\calD$\, so that \,$(\Sigma,\calD)$\, is the spectral data of a monodromy \,$M(\lambda)\in\Mon$\,. The first of these conditions, which we call
``\,$\calD$\, is \emph{compatible}'' below, is a condition concerning the local structure of spectral divisors at singularities of \,$\Sigma$\,, namely precisely the property that has been shown
for spectral divisors in Proposition~\ref{P:spectrum:dims}(2). Thus every spectral divisor \,$\calD$\, also is compatible.

The second condition that we need to impose is one also concerning the asymptotic behavior of \,$\calD$\,, namely that \,$\calD$\, is
non-special (which is also required in the case of finite type divisors). We have not yet shown that all spectral divisors are non-special (because the definition of specialty involves the asymptotic
behavior near \,$\lambda=\infty$\, resp.~\,$\lambda=0$\,, we could not have done so in Section~\ref{Se:spectrum}), but this follows from Theorem~\ref{T:special:special} below, which states 
that a compatible, asymptotic, positive generalized divisor \,$\calD$\, is a spectral divisor corresponding to some monodromy \,$M(\lambda)$\, if and only if it is non-special.

As in Section~\ref{Se:spectrum}, we let \,$\mathcal{O}$\, be the sheaf of holomorphic functions on \,$\Sigma$\,, \,$\wh{\Sigma}$\, be the normalization of \,$\Sigma$\,, and \,$\wh{\calO}$\, be the 
direct image onto \,$\Sigma$\, of the sheaf of holomorphic functions on \,$\wh{\Sigma}$\,. 
For a positive generalized divisor \,$\calD$\, we call the classical divisor
$$ D: \Sigma\to\N_0,\;(\lambda,\mu)\mapsto \dim(\mathcal{D}_{(\lambda,\mu)}/\mathcal{O}_{(\lambda,\mu)}) $$
the underlying classical divisor of \,$\Sigma$\,, and in Definition~\ref{D:asympdiv:asympdiv} we described what it means for \,$\calD$\, resp.~\,$D$\, to be asymptotic. 
In the sequel, we will also use other notations from Section~\ref{Se:spectrum}.

\begin{Def}
\label{D:special:special}
Let \,$\calD$\, be an asymptotic, positive generalized divisor on \,$\Sigma$\,.
\begin{enumerate}
\item
We call \,$\calD$\, \emph{compatible}, if for every singular point \,$(\lambda_*,\mu_*)\in \Sigma$\, that is in the support of \,$\calD$\,, say with degree \,$m := \dim(\calD_{(\lambda_*,\mu_*)}/\calO_{(\lambda_*,\mu_*)}) \geq 1$\,,
there exists \,$g \in \calD_{(\lambda_*,\mu_*)}$\,
so that \,$\eta := g-\tfrac{\mu-\mu^{-1}}{(\lambda-\lambda_*)^m}$\, is a meromorphic function germ in \,$\lambda$\, alone, with \,$\polord^{\Sigma}(\eta) \leq \polord^{\Sigma}(g)$\,.
\item
We call \,$\mathcal{D}$\, \emph{special}, if \,$\mathcal{D}$\, contains a section \,$f$\,
that is a non-constant meromorphic function on \,$\Sigma$\, which is bounded on \,$\wh{V}_\delta$\, for some \,$\delta>0$\,.
Otherwise we say that \,$\mathcal{D}$\, is \emph{non-special}.
\end{enumerate}
\end{Def}

We prepare the proof of the theorem on the reconstruction of the monodromy from the spectral data (Theorem~\ref{T:special:special}) with the following lemma:

\begin{lem}
\label{L:special:sectional-lemma}
Let \,$\calD$\, be a compatible, positive generalized divisor on \,$\Sigma$\,, and \,$(\lambda_*,\mu_*)\in \Sigma$\, is in the support of \,$\calD$\,, i.e.~\,$m:=\dim(\calD_{(\lambda_*,\mu_*)}/\calO_{(\lambda_*,\mu_*)}) \geq 1$\, holds.
We further suppose that \,$(\lambda-\lambda_*)^{-1}$\, is \emph{not} a section of \,$\calD$\,.

Moreover, we suppose that holomorphic function germs \,$c,d$\, in \,$\lambda$\, at \,$\lambda=\lambda_*$\, are given with \,$\ord^{\C}_{\lambda_*}(c)=m$\,. 
In the sequel we give a characterization of the condition that \,$\tfrac{\mu-d}{c}\in \mathcal{D}_{(\lambda_*,\mu)}$\, holds for ``all'' (one or two) values of \,$\mu\in \C^*$\, with \,$(\lambda_*,\mu)\in \Sigma$\,; for this we distinguish
between regular and singular points \,$(\lambda_*,\mu_*)$\, of \,$\Sigma$\,. 

\begin{enumerate}
\item 
Suppose that \,$(\lambda_*,\mu_*)$\, is a regular point of \,$\Sigma$\,. If \,$(\lambda_*,\mu_*)$\, is not a branch point (i.e.~\,$\Delta(\lambda_*)^2-4 \neq 0$\, holds), then the holomorphic function \,$\mu^{-1}$\, on \,$\Sigma$\,
can near \,$(\lambda_*,\mu_*)$\, be represented as a holomorphic function \,$\vartheta(\lambda)$\, in \,$\lambda$\,, and we have
\,$\tfrac{\mu-d}{c} \in \calD_{(\lambda_*,\mu_*)}$\, and \,$\tfrac{\mu-d}{c} \in \calD_{(\lambda_*,\mu_*^{-1})}$\, if and only if we have
\begin{equation}
\label{eq:special:sectional-lemma:regular-nobranch}
d^{(k)}(\lambda_*) = \vartheta^{(k)}(\lambda_*) \qmq{for \,$k\in\{0,\dotsc,m-1\}$\,.}
\end{equation}
If \,$(\lambda_*,\mu_*)$\, is a branch point (i.e.~\,$\Delta^2-4$\, has a simple zero at \,$\lambda_*$\,) then we have \,$m=1$\,, and \,$\tfrac{\mu-d}{c}\in \calD_{(\lambda_*,\mu_*)}$\, holds if and only if we have
\begin{equation}
\label{eq:special:sectional-lemma:regular-branch}
d(\lambda_*) = \mu_*^{-1} \; . 
\end{equation}

\item
Suppose that \,$(\lambda_*,\mu_*)$\, is a singular point of \,$\Sigma$\, (i.e.~\,$\Delta^2-4$\, has a zero of order \,$\wh{n} \geq 2$\, at \,$\lambda_*$\,); we let
\,$g \in \calD_{(\lambda_*,\mu_*)}$\, be so that \,$\eta := g-\tfrac{\mu-\mu^{-1}}{(\lambda-\lambda_*)^m}$\, is a meromorphic function germ in \,$\lambda$\, alone, with \,$\polord^{\Sigma}(\eta) \leq \polord^{\Sigma}(g)$\,.
Then \,$m\leq \wh{n}-j_0$\, holds, \,$\wt{\eta} := \frac12\bigr( \Delta-(\lambda-\lambda_*)^m \cdot \eta\bigr)$\, 
is a holomorphic function germ in \,$\lambda$\,, and \,$\tfrac{\mu-d}{c} \in \calD_{(\lambda_*,\mu_*)}$\, holds
if and only if we have
\begin{equation}
\label{eq:special:sectional-lemma:singular}
d^{(k)}(\lambda_*) = \wt{\eta}^{(k)}(\lambda_*) \qmq{for \,$k\in\{0,\dotsc,m-1\}$\,.}
\end{equation}
\end{enumerate}
\emph{Addendum.} If the condition of Equation~\eqref{eq:special:sectional-lemma:regular-nobranch}, \eqref{eq:special:sectional-lemma:regular-branch} or \eqref{eq:special:sectional-lemma:singular}
holds in the setting of (1) resp.~(2), then the function germ \,$b := \tfrac{(\Delta-d)\cdot d - 1}{c}$\, is holomorphic at \,$\lambda_*$\,.
\end{lem}

\begin{proof}
\emph{For (1).}
If \,$(\lambda_*,\mu_*)$\, is a regular point of \,$\Sigma$\,, then \,$\calD_{(\lambda_*,\mu_*)}$\, is completely characterized by the underlying classical divisor, therefore a meromorphic function germ \,$f$\,
is in \,$\calD_{(\lambda_*,\mu_*)}$\, if and only if it has at \,$(\lambda_*,\mu_*)$\, at most a pole of order \,$m$\,. 
Let us now first suppose that \,$(\lambda_*,\mu_*)$\, is not a branch point of \,$\Sigma$\,. We then have \,$\mu_*^{-1} \neq \mu_*$\,, and because of the hypothesis
that \,$(\lambda-\lambda_*)^{-1}$\, is not a section of \,$\calD$\,, we see that \,$(\lambda_*,\mu_*^{-1})$\, is not in the support of \,$\calD$\,. Therefore \,$f\in \calD_{(\lambda_*,\mu_*^{-1})}$\, holds 
if and only if \,$f$\, is holomorphic at \,$(\lambda_*,\mu_*)$\,. Moreover, \,$\calO_{(\lambda_*,\mu_*)}$\, is generated by \,$\lambda$\, as
a ring, and therefore the holomorphic function \,$\mu^{-1}$\, can near \,$(\lambda_*,\mu_*)$\, be represented as a holomorphic function \,$\vartheta(\lambda)$\, in \,$\lambda$\,. 
In any event, \,$\tfrac{\mu-d}{c}$\, has at \,$(\lambda_*,\mu_*)$\, a pole of order at most \,$m$\,, and therefore \,$\tfrac{\mu-d(\lambda)}{c(\lambda)} \in \calD_{(\lambda_*,\mu_*)}$\, holds. In this situation,
Equation~\eqref{eq:special:sectional-lemma:singular} is equivalent to the fact that the holomorphic function germ \,$\mu-d(\lambda)$\, has at \,$(\lambda_*,\mu_*^{-1})$\, a zero of order at least \,$m$\,,
and this is in turn equivalent to the fact that \,$\tfrac{\mu-d(\lambda)}{c(\lambda)}$\, is holomorphic at \,$(\lambda_*,\mu_*^{-1})$\,, which is equivalent to \,$\tfrac{\mu-d(\lambda)}{c(\lambda)} \in 
\calD_{(\lambda_*,\mu_*)}$\,. Moreover, we have
\begin{equation}
\label{eq:special:sectional-lemma:regular-Deltadd-1}
(\Delta-d)\cdot d - 1 = (\mu+\mu^{-1}-d)\cdot d - \mu\cdot \mu^{-1} = -(\mu-d)\cdot (\mu^{-1}-d) \;, 
\end{equation}
therefore the function \,$(\Delta-d)\cdot d-1$\, has a zero of order at least \,$m$\, at \,$\lambda_*$\,. Because \,$c$\, also has a zero of order \,$m$\, at \,$\lambda_*$\,, it follows that \,$b$\, is holomorphic
at \,$\lambda_*$\,.

Now suppose that \,$(\lambda_*,\mu_*)$\, is a regular branch point of \,$\Sigma$\,. Then \,$\mu_*=\mu_*^{-1}$\, holds, and because \,$\mathrm{d}\lambda$\, has a simple zero at \,$\lambda_*$\,, the function
\,$(\lambda-\lambda_*)^{-1}$\, has a pole of order \,$2$\, when regarded as a meromorphic function on \,$\Sigma$\,. Therefore \,$m\geq 2$\, would imply that \,$(\lambda-\lambda_*)^{-1}$\, is a section
of \,$\calD$\,, in contradiction to the hypothesis. Hence \,$m=1$\, holds. Therefore \,$c$\, has a simple zero as a function in \,$\lambda$\,, but a zero of order \,$2$\, as a holomorphic function germ
on \,$\Sigma$\,. In this situation, the equation \,$d(\lambda_*)=\mu_*^{-1}=\mu_*$\, is equivalent to the fact that \,$\mu-d(\lambda)$\, has (at least) a simple zero at \,$(\lambda_*,\mu_*)$\, 
(because \,$\mu$\, is a coordinate on \,$\Sigma$\, near \,$(\lambda_*,\mu_*)$\,); this is equivalent to \,$\tfrac{\mu-d}{c}$\, having at most a pole of order \,$1$\, at \,$(\lambda_*,\mu_*)$\, and therefore to
\,$\tfrac{\mu-d}{c}\in\calD_{(\lambda_*,\mu_*)}$\,.
Moreover, because \,$\mu^{-1}-d$\, also has a zero at \,$(\lambda_*,\mu_*)$\,, Equation~\eqref{eq:special:sectional-lemma:regular-Deltadd-1} shows that the function \,$(\Delta-d)\cdot d-1$\, has as function on \,$\Sigma$\,
a zero of order at least \,$2$\,, and therefore as function in \,$\lambda$\, a zero of order at least \,$1$\,. Because \,$c$\, also has a zero of order \,$m=1$\, as function in \,$\lambda$\,, it follows again
that \,$b$\, is holomorphic.

\emph{For (2).} 
Because \,$\calD$\, is compatible, there exists \,$g \in \calD_{(\lambda_*,\mu_*)}$\,
so that \,$\eta := g-\tfrac{\mu-\mu^{-1}}{(\lambda-\lambda_*)^m}$\, is a meromorphic function germ in \,$\lambda$\, alone, with \,$\polord^{\Sigma}(\eta) \leq \polord^{\Sigma}(g)$\,.
Here the pole order of \,$\eta$\, (as a function of \,$\lambda$\,) is at most \,$m$\,, and therefore \,$(\lambda-\lambda_*)^m\cdot \eta$\, and hence also \,$\wt{\eta}$\, is holomorphic. 
Moreover, because \,$c$\, has at \,$\lambda_*$\, a zero of order \,$m$\,, there exists an invertible \,$\gamma\in \calO_{\lambda_*}(\lambda)$\, with \,$c=\gamma\cdot (\lambda-\lambda_*)^m$\,. We now calculate
\begin{align*}
\frac{\mu-d}{c} - \frac{1}{2\gamma}\,g
& = \frac{\mu-d}{c} - \frac{1}{2\gamma}\left( \eta + \frac{\mu-\mu^{-1}}{(\lambda-\lambda_*)^m} \right) \\
& = \frac{\mu-d}{c} - \frac{\mu-\mu^{-1}}{2\gamma\,(\lambda-\lambda_*)^m} - \frac{(\lambda-\lambda_*)^m\cdot \eta}{2\gamma\,(\lambda-\lambda_*)^m} \\
& = \frac{2\mu-2d-\mu+\mu^{-1}-(\lambda-\lambda_*)^m\cdot \eta}{2\,c} 
= \frac{\mu+\mu^{-1}-2d-(\lambda-\lambda_*)^m\cdot \eta}{2\,c} \\
& = \frac{\Delta-2d-(\lambda-\lambda_*)^m\cdot \eta}{2\,c} 
= \frac{\wt{\eta}-d}{c} \; . 
\end{align*}
Because we have \,$\tfrac{1}{2\gamma}g \in \calD_{(\lambda_*,\mu_*)}$\,, this calculation shows that \,$\tfrac{\mu-d}{c} \in \calD_{(\lambda_*,\mu_*)}$\, holds if and only if 
\,$\tfrac{\wt{\eta}-d}{c} \in \calD_{(\lambda_*,\mu_*)}$\, holds. Because the latter germ is meromorphic in \,$\lambda$\, alone, and we have \,$(\lambda-\lambda_*)^{-1}\not\in\calD_{(\lambda_*,\mu_*)}$\, by
hypothesis, \,$\tfrac{\wt{\eta}-d}{c} \in \calD_{(\lambda_*,\mu_*)}$\, holds if and only if \,$\tfrac{\wt{\eta}-d}{c}$\, is holomorphic, and this is true if and only if \,$\wt{\eta}-d$\, has
a zero of order at least \,$m$\, (with respect to \,$\lambda$\,). The latter condition is equivalent to \eqref{eq:special:sectional-lemma:singular}. 


For the proof of the Addendum, we let \,$\wh{n}$\, be the order of the singularity \,$(\lambda_*,\mu_*)$\, of \,$\Sigma$\,, i.e.~\,$\wh{n}$\, is the order of the zero of \,$\Delta^2-4$\, at \,$\lambda_*$\,.
In the sequel, we will distinguish between the cases that \,$\wh{n}=2n+1$\, is odd, and that \,$\wh{n}=2n$\, is even. Moreover let \,$s\geq 0$\, be the maximal pole order that occurs in \,$\calD_{(\lambda_*,\mu_*)}$\,. 
As a first stage for the proof, we will show that 
\begin{equation}
\label{eq:special:sectional-lemma:s-dim}
s \leq \wh{n} - 2j_0 \; 
\end{equation}
holds, where \,$j_0\in\{0,\dotsc,n\}$\, is the number from Proposition~\ref{P:spectrum:locally-free}. Let \,$f\in \calD_{(\lambda_*,\mu_*)}$\, be of pole order \,$s$\,. By the fact that \,$\calD$\, is a \,$\calO$-module
and by Proposition~\ref{P:spectrum:locally-free} it then follows that
\begin{equation}
\label{eq:special:sectional-lemma:s-dim-pre}
(\lambda-\lambda_*)^{n-j_0}\cdot f \;,\; \frac{\mu-\mu^{-1}}{(\lambda-\lambda_*)^{j_0}}\cdot f \;\in\; \calD_{(\lambda_*,\mu_*)} \; .
\end{equation}
If \,$\wh{n}=2n+1$\, is odd, then \eqref{eq:special:sectional-lemma:s-dim-pre} means that
$$ \sqrt{\lambda-\lambda_*}^{2n-2j_0}\cdot f\;,\; \sqrt{\lambda-\lambda_*}^{2n-2j_0+1}\cdot f \;\in\; \calD_{(\lambda_*,\mu_*)} $$
and therefore 
$$ \sqrt{\lambda-\lambda_*}^{-s+2n-2j_0} \cdot \wh{\calO}_{(\lambda_*,\mu_*)} \;\subset\; \calD_{(\lambda_*,\mu_*)} $$
holds, where \,$\wh{\calO}$\, denotes the direct image of the sheaf of holomorphic functions on the normalization of \,$\Sigma$\,. 
If \,$-s+2n-2j_0\leq -2$\, were true, this would imply \,$(\lambda-\lambda_*)^{-1}\in\calD_{(\lambda_*,\mu_*)}$\,, in contradiction to the hypothesis. Thus we have \,$-s+2n-2j_0 \geq -1$\, and therefore
\eqref{eq:special:sectional-lemma:s-dim} holds in this case. 

If \,$\wh{n}=2n$\, is even, then we represent \,$f=(f_1,f_2)$\, as a pair of germs of meromorphic functions in \,$\lambda$\, at the two points above \,$(\lambda_*,\mu_*)$\, in the normalization of \,$\Sigma$\,,
and denote the pole order of \,$f_\nu$\, by \,$s_\nu$\, (\,$\nu\in \{1,2\}$\,), ordered such that \,$s_1\geq s_2$\, holds. Then we have \,$s_1 +s_2=s$\,, and 
\eqref{eq:special:sectional-lemma:s-dim-pre} means that
$$ (\lambda-\lambda_*)^{-s_1+n-j_0} \cdot \wh{\calO}_{(\lambda_*,\mu_*)} \;\subset\; \calD_{(\lambda_*,\mu_*)} $$
holds. If \,$-s_1+n-j_0 \leq -1$\, were true, this would imply \,$(\lambda-\lambda_*)^{-1}\in\calD_{(\lambda_*,\mu_*)}$\,, in contradiction to the hypothesis. Therefore we have \,$-s_1+n-j_0 \geq 0$\,, whence
\eqref{eq:special:sectional-lemma:s-dim} follows also in this case.

It follows from \eqref{eq:special:sectional-lemma:s-dim} by Proposition~\ref{P:spectrum:msj0} that 
\begin{equation}
\label{eq:special:sectional-lemma:m-dim}
m = s+j_0 \leq \wh{n} - j_0 
\end{equation}
holds.

To show that \,$b$\, is holomorphic in \,$\lambda$\,, we note that by \eqref{eq:special:sectional-lemma:singular}, \,$d-\wt{\eta}$\, has at \,$\lambda_*$\, a zero of order \,$m$\, (with respect to \,$\lambda$\,),
so we can represent \,$d$\, in the form \,$d=\wt{\eta}+(\lambda-\lambda_*)^m\cdot \zeta$\, with a holomorphic function germ \,$\zeta$\, in \,$\lambda$\,. We now calculate:
\begin{align}
b & = \frac{(\Delta-d)\cdot d - 1}{c} = \frac{(\Delta-\wt{\eta}-(\lambda-\lambda_*)^m\cdot\zeta)\cdot (\wt{\eta}+(\lambda-\lambda_*)^m\cdot\zeta) - 1}{\gamma\cdot (\lambda-\lambda_*)^m} \notag \\
\label{eq:special:sectional-lemma:calc1}
& = \frac{1}{\gamma}\cdot \left( \frac{(\Delta-\wt{\eta})\cdot \wt{\eta} - 1}{(\lambda-\lambda_*)^m} + \Delta\,\zeta-2\,\wt{\eta}\,\zeta-(\lambda-\lambda_*)^m\,\zeta^2 \right) \; . 
\end{align}
Because \,$\Delta\,\zeta-2\,\wt{\eta}\,\zeta-(\lambda-\lambda_*)^m\,\zeta^2$\, is obviously holomorphic, it remains to show that \,$\tfrac{(\Delta-\wt{\eta})\cdot \wt{\eta} - 1}{(\lambda-\lambda_*)^m}$\,
is holomorphic. For this purpose, we continue to calculate:
\begin{align}
\frac{(\Delta-\wt{\eta})\cdot \wt{\eta} - 1}{(\lambda-\lambda_*)^m}
& = \frac{\bigr( \Delta-\tfrac12(\Delta-(\lambda-\lambda_*)^m\eta) \bigr) \cdot \tfrac12\bigr( \Delta-(\lambda-\lambda_*)^m\eta \bigr) -1}{(\lambda-\lambda_*)^m} \notag \\
& = \frac14\,\frac{\bigr( \Delta+(\lambda-\lambda_*)^m \eta\bigr) \cdot \bigr( \Delta-(\lambda-\lambda_*)^m\eta \bigr) -4}{(\lambda-\lambda_*)^m} \notag \\
\label{eq:special:sectional-lemma:calc2}
& = \frac14\,\frac{\Delta^2-\bigr((\lambda-\lambda_*)^m\,\eta\bigr)^2-4}{(\lambda-\lambda_*)^m} = \frac14\left( \frac{\Delta^2-4}{(\lambda-\lambda_*)^m} - (\lambda-\lambda_*)^m\cdot \eta^2 \right) \; . 
\end{align}
\,$\Delta^2-4$\, has (with respect to \,$\lambda$\,) a zero of order \,$\wh{n}$\,, whereas \,$(\lambda-\lambda_*)^m$\, has by \eqref{eq:special:sectional-lemma:m-dim} a zero of order at most \,$\wh{n}-j_0 \leq \wh{n}$\,. 
Therefore \,$\tfrac{\Delta^2-4}{(\lambda-\lambda_*)^m}$\, is holomorphic in \,$\lambda$\,. Moreover, we have
$$ \polord^{\C}(\eta^2) = \polord^\Sigma(\eta) \leq \polord^\Sigma(g) \leq s = m-j_0 \leq m $$
(where the last equals sign again follows from Proposition~\ref{P:spectrum:msj0}), and hence \,$(\lambda-\lambda_*)^m\cdot \eta^2$\, is also holomorphic. Therefore it follows from Equation~\eqref{eq:special:sectional-lemma:calc2}
that \,$\tfrac{(\Delta-\wt{\eta})\cdot \wt{\eta} - 1}{(\lambda-\lambda_*)^m}$\, is holomorphic, and thus we see from Equation~\eqref{eq:special:sectional-lemma:calc1} that \,$b$\, is holomorphic.
\end{proof}

Now we are ready to prove that any non-special, compatible, asymptotic, positive generalized divisor on \,$\Sigma$\, is the spectral divisor of an asymptotic monodromy \,$M(\lambda) \in \Mon$\,.
This is the content of the implication (a) $\Rightarrow$ (c) of the following theorem. Also note the equivalence (a) $\Leftrightarrow$ (b) of the theorem, which shows that for a
compatible, asymptotic, positive divisor \,$\calD$\, the property that \,$\calD$\, is non-special (which is by its definition related to the behavior of \,$\calD$\, near \,$\lambda=0$\, and \,$\lambda=\infty$\,)
can be characterized in terms of the local behavior of \,$\calD$\, at the points of its support. 

\newpage

\begin{thm}
\label{T:special:special}
Let \,$\mathcal{D}$\, be an compatible, asymptotic, positive generalized divisor on \,$\Sigma$\,. Then the following statements are equivalent:
\begin{itemize}
\item[(a)] \,$\mathcal{D}$\, is non-special.
\item[(b)] There does not exist any \,$\lambda_* \in \C^*$\, so that the meromorphic function \,$(\lambda-\lambda_*)^{-1}$\, on \,$\Sigma$\, is a section of \,$\mathcal{D}$\,. 
\item[(c)] There exists a monodromy \,$M(\lambda)=\left( \begin{smallmatrix} a(\lambda) & b(\lambda)\\c(\lambda) & d(\lambda) \end{smallmatrix} \right) \in \Mon$\, with \,$a+d=\Delta$\,, such that \,$\mathcal{D}$\, is the spectral divisor of \,$M(\lambda)$\,, i.e.~\,$\mathcal{D}$\, is generated by \,$1$\, and \,$\tfrac{\mu-d}{c}$\, over \,$\mathcal{O}$\,. Moreover, \,$M(\lambda)$\, is determined uniquely
up to a joint change of sign of the functions \,$b$\, and \,$c$\,. 
\item[(d)] There is one and only one global meromorphic function \,$f$\, on \,$\Sigma$\, that is a section of \,$\mathcal{D}$\, and 
such that 
\,$\left( f-\tfrac{i}{\tau\,\sqrt{\lambda}} \right)|\wh{V}_\delta \in \As_\infty(\wh{V}_\delta,\ell^2_{1},0)$\,
and \,$\left( f-\tfrac{i}{\tau^{-1}\sqrt{\lambda}} \right)|\wh{V}_\delta \in \As_0(\wh{V}_\delta,\ell^2_{-1},0)$\,
holds for some \,$\delta>0$\,.
(Here \,$\sqrt{\lambda}$\, is a fixed holomorphic branch of the square root function on \,$\wh{V}_\delta$\,, and \,$\tau\in\C^*$\, is defined by Equation~\eqref{eq:interpolate:lambda:tau} via
the underlying classical divisor \,$D=\{(\lambda_k,\mu_k)\}$\, of \,$\calD$\,.)
\end{itemize}
\emph{Addendum.} If \,$\mathcal{D}$\, is non-special, then the following statements hold concerning the monodromy \,$M(\lambda)$\, 
in (c), where we let \,$D=\{(\lambda_k,\mu_k)\}_{k\in \Z}$\, be the underlying classical divisor of \,$\mathcal{D}$\,:
\begin{itemize}
\item[(1)] \,$c$\, is characterized up to sign by the fact that it has zeros in all the \,$\lambda_k$\, (with the correct multiplicity) and no others; this function is given by Equation~\eqref{eq:interpolate:lambda:c} 
in Proposition~\ref{P:interpolate:lambda}.
\item[(2)] The functions \,$a$\, resp.~\,$d$\, satisfy \,$a(\lambda_k)=\mu_k$\, resp.~\,$d(\lambda_k)=\mu_k^{-1}$\, for all \,$k\in \Z$\,; these functions are given by Equations~\eqref{eq:f-val:interpolate:g}
and \eqref{eq:f-val:interpolate:wta} in Proposition~\ref{P:interpolate:mu} (where the \,$t_{\lambda_*,j}$\, are uniquely determined from the spectral data \,$(\Sigma,\mathcal{D})$\,).
\item[(3)] The function \,$b$\, is characterized by \,$\det(M)=ad-bc=1$\,.
\item[(4)] We have \,$M(\lambda) \in \Mon_\tau$\, with \,$\tau = \pm \left( \prod_{k\in\Z} \frac{\lambda_{k,0}}{\lambda_k} \right)^{1/2}$\,.
\item[(5)] \,$M(\lambda)$\, is uniquely determined by the spectral data \,$(\Sigma,\mathcal{D})$\, up to a joint change of sign of \,$b$\,, \,$c$\, and \,$\tau$\,. 
\end{itemize}
Moreover, the uniquely determined function \,$f$\, from (d) then equals \,$\pm \tfrac{\mu-d}{c}$\,. If we choose a sign for the holomorphic function \,$\sqrt{\lambda}$\, on \,$\wh{V}_\delta$\,, the condition
that \,$f=\tfrac{\mu-d}{c}$\, holds in (d) determines the sign of \,$c$\, in (1) and the sign of \,$\tau$\, in (4), and therefore determines \,$M(\lambda)$\, in (5).
\end{thm}

\begin{proof}
\emph{For (a) $\Rightarrow$ (b).}
Assume to the contrary that there exists \,$\lambda_* \in \C^*$\, so that \,$f := (\lambda-\lambda_*)^{-1}$\, is a section of \,$\mathcal{D}$\,.
Then we can choose \,$\delta>0$\, so large that \,$\lambda_* \not\in V_\delta$\, holds, and with this choice of \,$\delta$\,,
the non-constant meromorphic function \,$f$\, is bounded on \,$\wh{V}_\delta$\,. Therefore \,$\mathcal{D}$\, would be special,
in contradiction to the hypothesis (a).

\smallskip




\emph{For (b) $\Rightarrow$ (c).}
Because \,$\mathcal{D}$\, is asymptotic, the underlying classical divisor \,$D$\, is of the form 
\,$D=\{(\lambda_k,\mu_k)\}_{k\in \Z}$\, with 
\,$\lambda_k-\lambda_{k,0} \in \ell^2_{-1,3}(k)$\, and \,$\mu_k-\mu_{k,0} \in \ell^2_{0,0}(k)$\,. 

By Proposition~\ref{P:interpolate:lambda}(1),
$$ \tau := \pm \left( \prod_{k\in\Z} \frac{\lambda_{k,0}}{\lambda_k} \right)^{1/2} $$
then converges in \,$\C^*$\,, and after fixing the sign of \,$\tau$\,,
there exists a holomorphic function \,$c(\lambda)$\, on \,$\C^*$\, with zeros in all the \,$\lambda_k$\, (with multiplicity) and no others, and which satisfies the asymptotic property
$$ c-\tau\,c_0 \in \As_\infty(\C^*,\ell^2_{-1},1) \qmq{and} c-\tau^{-1}\,c_0 \in \As_0(\C^*,\ell^2_1,1) \; ;$$
\,$c$\, is unique by Proposition~\ref{P:excl:unique}(1), and \,$c(\lambda)$\, is given explicitly by Equation~\eqref{eq:interpolate:lambda:c}.


The next step of the construction is to obtain a holomorphic function \,$d$\, with \,$d-d_0 \in \As(\C^*,\ell^2_{0,0},1)$\, such that \,$\tfrac{\mu-d}{c}$\, is a global section of \,$\mathcal{D}$\,. 
For this purpose, we apply Lemma~\ref{L:special:sectional-lemma}, which is applicable because of the hypothesis (b).
This lemma shows that the condition that \,$\tfrac{\mu-d}{c}$\, be a section of \,$\mathcal{D}$\, is equivalent to the prescription of the values of
\,$d^{(\nu)}(\lambda_k)$\, for every \,$k\in \Z$\, and \,$\nu\in \{0,\dotsc,m-1\}$\,, where \,$m$\, is the degree of the point \,$(\lambda_k,\mu_k)$\, in \,$\calD$\,, and \,$m$\, also equals the multiplicity
of the zero \,$\lambda_k$\, of \,$c$\,. Therefore Proposition~\ref{P:interpolate:mu}(1), (2) shows (see Remark~\ref{R:interpolate:mu}) 
that there exists one and only one holomorphic function \,$d$\, with \,$d-d_0\in\As(\C^*,\ell^2_{0,0},1)$\,
so that the conditions from Lemma~\ref{L:special:sectional-lemma} are satisfied for all points in the support of \,$\calD$\,, and therefore \,$\tfrac{\mu-d}{c}$\, is a global section of \,$\calD$\,.



We let \,$a := \Delta-d$\,, this is also a holomorphic function in \,$\lambda$\, and we have \,$a-a_0 \in \As(\C^*,\ell^2_{0,0},1)$\,. Moreover, \,$b := \tfrac{ad-1}{c}$\, is holomorphic by the Addendum of 
Lemma~\ref{L:special:sectional-lemma}.
From the equation
$$ b - \tau^{-1}\,b_0 = \frac{ad-1}{c}-\frac{a_0\,d_0-1}{\tau\,c_0} = \frac{1}{c} \cdot \bigr( (a-a_0)\,d+a_0\,(d-d_0) \bigr) - \frac{b_0}{\tau\,c}\cdot (c-\tau\,c_0) $$
one obtains \,$(b-\tau^{-1}\,b_0)|V_\delta \in \As_\infty(\C^*,\ell^2_1,1)$\, and then by Proposition~\ref{P:interpolate:l2asymp}(1) \,$b-\tau^{-1}\,b_0 \in \As_\infty(\C^*,\ell^2_1,1)$\,.
Similarly one also shows \,$b-\tau\,b_0 \in \As_0(\C^*,\ell^2_{-1},1)$\,. 

We now consider
$$ M(\lambda) := \begin{pmatrix} a(\lambda) & b(\lambda) \\ c(\lambda) & d(\lambda) \end{pmatrix} \; . $$
By the preceding construction, we have \,$\det(M(\lambda))=1$\,, \,$M(\lambda) \in \Mon_\tau$\, and \,$a+d=\Delta$\,. Moreover, both \,$1$\, and \,$\tfrac{\mu-d}{c}$\, are global sections of \,$\calD$\,, so the
spectral divisor \,$\calD_M$\, of \,$M(\lambda)$\, (which is generated by \,$1$\, and \,$\tfrac{\mu-d}{c}$\,) is contained in \,$\calD$\,. Both \,$\calD$\, and \,$\calD_M$\, are asymptotic divisors, i.e.~they have
asymptotically and totally the same degree. Therefore \,$\calD_M \subset \calD$\, in fact implies \,$\calD_M = \calD$\,.

\smallskip

\emph{For (c) $\Rightarrow$ (d).}
We let \,$M(\lambda)=\left( \begin{smallmatrix} a(\lambda) & b(\lambda)\\c(\lambda) & d(\lambda) \end{smallmatrix} \right)$\,
be the monodromy from (c). Then \,$f:=\tfrac{\mu-d}{c}$\, is a section of \,$\mathcal{D}$\,, and 
\,$\left( f-\tfrac{i}{\tau\,\sqrt{\lambda}} \right)|\wh{V}_\delta \in \As_\infty(\wh{V}_\delta,\ell^2_{1},0)$\,
and \,$\left( f-\tfrac{i}{\tau^{-1}\sqrt{\lambda}} \right)|\wh{V}_\delta \in \As_0(\wh{V}_\delta,\ell^2_{-1},0)$\,
hold by Corollary~\ref{C:asympfinal:monodromy2}(6). 

For the proof of the uniqueness of \,$f$\,, we suppose that another section \,$\wt{f}$\, of \,$\mathcal{D}$\, with 
\,$\left( \wt{f}-\tfrac{i}{\tau\,\sqrt{\lambda}} \right)|\wh{V}_\delta \in \As_\infty(\wh{V}_\delta,\ell^2_{1},0)$\,
and \,$\left( \wt{f}-\tfrac{i}{\tau^{-1}\sqrt{\lambda}} \right)|\wh{V}_\delta \in \As_0(\wh{V}_\delta,\ell^2_{-1},0)$\,
is given. We first show
that the odd parts of \,$f$\, and \,$\wt{f}$\, are equal. Because \,$f$\, and \,$\wt{f}$\, are sections in \,$\calD$\,, 
\,$g := c\cdot (f-\wt{f})$\, is a holomorphic function on \,$\Sigma$\,, and we have \,$(f-\wt{f})|\wh{V}_\delta \in \As(\wh{V}_\delta,\ell^2_{1,-1},0)$\,,
and therefore \,$g|\wh{V}_\delta \in \As(\wh{V}_\delta,\ell^2_{0,0},1)$\,. We now represent \,$g=g_+ + (\mu-\mu^{-1})\cdot g_-$\, with 
holomorphic functions \,$g_+,g_-$\, in \,$\lambda\in \C^*$\,. Along with \,$g$\,, we also have \,$((\mu-\mu^{-1})\cdot g_-)|\wh{V}_\delta
\in \As(\wh{V}_\delta,\ell^2_{0,0},1)$\,. Because \,$\mu-\mu^{-1}$\, is comparable to \,$w(\lambda)$\, on \,$\wh{V}_\delta$\,,
this implies \,$g_-|V_\delta \in \As(V_\delta,\ell^2_{0,0},0)$\, and therefore by Proposition~\ref{P:interpolate:l2asymp}(1) 
\,$g_- \in \As(\C^*,\ell^2_{0,0},0)$\,. This shows that the holomorphic function \,$g_-$\, on \,$\C^*$\, can be extended to a holomorphic
function on \,$\PP^1$\, by setting \,$g_-(0)=g_-(\infty)=0$\,. It follows that \,$g_-=0$\, holds, and thus the odd parts
of \,$f$\, and \,$\wt{f}$\, are equal.

Because we have
$$ f = \frac{\mu-d}{c} = \frac{(\Delta-2d)+(\mu-\mu^{-1})}{2c} \;, $$
it follows that there exists a meromorphic function \,$\wt{d}$\, in \,$\lambda$\, so that
$$ \wt{f} = \frac{(\Delta-2\wt{d})+(\mu-\mu^{-1})}{2c} = \frac{\mu-\wt{d}}{c} $$
holds. Because \,$g=c\cdot (\wt{f}-f)$\, is holomorphic on \,$\Sigma$\,, \,$\wt{d}$\, is in fact holomorphic in \,$\lambda$\,. Moreover, \,$(\wt{f}-f)|\wh{V}_\delta \in \As(\wh{V}_\delta,\ell^2_{1,-1},0)$\, implies
$$ \left. -\frac{d-\wt{d}}{c} \right|_{V_\delta} = \left.\left( \frac{\Delta-2\wt{d}}{2c}- \frac{\Delta-2d}{2c} \right) \right|_{V_\delta} \;\in\; \As(V_\delta,\ell^2_{1,-1},0) $$
and therefore \,$(d-\wt{d})|V_\delta  \in \As(V_\delta,\ell^2_{0,0},1)$\,, hence \,$d-\wt{d} \in \As(\C^*,\ell^2_{0,0},1)$\, by Proposition~\ref{P:interpolate:l2asymp}(1). Lemma~\ref{L:special:sectional-lemma}
shows that the condition \,$\tfrac{\mu-d}{c} \in H^0(\Sigma,\calD)$\, prescribes for any \,$(\lambda_*,\mu_*)$\, in the support of \,$\calD$\,, say of degree \,$m$\,, the values of
\,$d(\lambda_*)$\,, \,$d'(\lambda_*)$\,, \,$\dotsc$\,, \,$d^{(m-1)}(\lambda_*)$\,, and therefore we have \,$\ord_{\lambda_*}(d-\wt{d}) \geq \ord_{\lambda_*}(c)$\, for every \,$\lambda_*\in \C$\, with \,$c(\lambda_*)=0$\,.
Proposition~\ref{P:excl:unique}(2) shows that \,$d=\wt{d}$\, and therefore \,$f=\wt{f}$\, holds.

\smallskip

\emph{For (d) $\Rightarrow$ (a).}
Suppose that (d) holds, i.e.~there exists one and only one section \,$f$\, of \,$\mathcal{D}$\, such that
\,$\left( f-\tfrac{i}{\tau\,\sqrt{\lambda}} \right)|\wh{V}_\delta \in \As_\infty(\wh{V}_\delta,\ell^2_{1},0)$\,
and \,$\left( f-\tfrac{i}{\tau^{-1}\sqrt{\lambda}} \right)|\wh{V}_\delta \in \As_0(\wh{V}_\delta,\ell^2_{-1},0)$\,
holds for some \,$\delta>0$\,.

Assume to the contrary that \,$\mathcal{D}$\, is special. Thus there exists a section
\,$g$\, of \,$\mathcal{D}$\, which is non-constant and bounded on \,$\wh{V}_\delta$\, for some \,$\delta>0$\,. Then we have in particular
\,$\tfrac{1}{\sqrt{\lambda}}\cdot g|\wh{V}_\delta \in \As(\wh{V}_\delta,\ell^2_{1,-1},0)$\,, and therefore for every
\,$s\in \C^*$\,, \,$f_s := f + s \cdot \tfrac{1}{\sqrt{\lambda}}\,g$\, is a section of \,$\mathcal{D}$\, different from \,$f$\,,
for which 
\,$\left( f_s-\tfrac{i}{\tau\,\sqrt{\lambda}} \right)|\wh{V}_\delta \in \As_\infty(\wh{V}_\delta,\ell^2_{1},0)$\,
and \,$\left( f_s-\tfrac{i}{\tau^{-1}\sqrt{\lambda}} \right)|\wh{V}_\delta \in \As_0(\wh{V}_\delta,\ell^2_{-1},0)$\, holds.
This contradicts the uniqueness statement of (d).
\end{proof}

With this theorem, we have clarified the equivalence between monodromies \,$M(\lambda) \in \Mon$\, and spectral data \,$(\Sigma,\calD)$\,. 

From this point onward, we would like to avoid technical complications by restricting the class of divisors in such a way that the generalized divisor \,$\calD$\, is uniquely determined by the
underlying classical divisor \,$D$\,. We will achieve this by requiring that the degree of \,$\calD$\, above any \,$\lambda_* \in \C^*$\, is at most \,$1$\,; we will call such divisors \emph{tame}.
The set of tame divisors is open and dense in \,$\Div$\,; also note that the divisor \,$\calD_0$\, of the vacuum is tame (see Section~\ref{Se:vacuum}). 
Because tame generalized divisors are uniquely characterized by their underlying classical divisors, we will work mostly with the classical
divisors in the case of tame divisors.

The concept is made precise by the following definition and proposition:

\begin{Def}
\label{D:special:tame}
\begin{enumerate}
\item We say that an asymptotic classical divisor \,$D=\{(\lambda_k,\mu_k)\} \in \Div$\, (viewed as a multi-set of points in \,$\C^*\times\C^*$\,) is \emph{tame} if the \,$\lambda_k$\, are pairwise unequal.
We denote the subset of tame divisors in \,$\Div$\, by \,$\Div_{tame}$\,.
\item We say that a generalized divisor \,$\calD$\, on the spectral curve \,$\Sigma$\,  is \emph{tame}, if it is compatible, asymptotic and positive, and if the total degree of \,$\calD$\, above any
\,$\lambda_* \in \C^*$\, is at most \,$1$\,. 
\item We say that a potential \,$(u,u_y)\in \Pot$\, is \emph{tame}, if its associated spectral divisor is tame. We denote the subset of tame potentials in \,$\Pot$\, by \,$\Pot_{tame}$\,. 
\end{enumerate}
\end{Def}

At the end of Section~\ref{Se:asympdiv} we identified the space \,$\Div$\, of asymptotic divisors with the quotient space \,$(\ell^2_{-1,3}\oplus \ell^2_{0,0})/P(\Z)$\,, and viewed in this way,
\,$\Div_{tame}$\, is an open and dense subset of \,$\Div$\,. 
The singular points of \,$\Div$\, are those \,$D\in\Div$\, which contain a point of multiplicity \,$\geq 2$\,, and  therefore the points of \,$\Div_{tame}$\, are regular points of \,$\Div$\,. 
Thus \,$\Div_{tame}$\, has the structure of an (infinite-dimensional) smooth manifold.

The following Proposition~\ref{P:special:phoenixfromashes}(1) shows that every tame generalized divisor \,$\calD$\, is non-special, and therefore is the spectral divisor of a monodromy \,$M(\lambda)$\,;
moreover it shows that if a point \,$(\lambda_k,\mu_k)$\, in the support of \,$\calD$\, is a singular point of \,$\Sigma$\, for \,$|k|$\, large, then \,$\calD$\, looks at \,$(\lambda_k,\mu_k)$\, like
the spectral divisor of the vacuum at any of its points, see Section~\ref{Se:vacuum}. Moreover Proposition~\ref{P:special:phoenixfromashes}(2) shows that the classical tame divisors \,$D$\,
(regarded as multi-sets of points in \,$\C^* \times \C^*$\,) are in one-to-one correspondence with the spectral data \,$(\Sigma,\calD)$\,, where \,$\calD$\, is a tame generalized divisor on the
spectral curve \,$\Sigma$\,. The latter fact justifies our approach of investigating tame classical (rather than generalized) divisors in the sequel. 

\begin{prop}
\label{P:special:phoenixfromashes}
\begin{enumerate}
\item
Let \,$\calD$\, be a tame generalized divisor on \,$\Sigma$\,. Then the following holds:
\begin{enumerate}
\item The underlying classical divisor \,$D$\, of \,$\calD$\, is tame.
\item \,$\calD$\, is non-special.
\item We write \,$D=\{(\lambda_k,\mu_k)\}$\, as in Proposition~\ref{P:excl:basic}(3). Then there exists \,$N\in \N$\, such that for every \,$k\in \Z$\, with \,$|k|> N$\, where \,$(\lambda_k,\mu_k)$\, 
is a singular point of \,$\Sigma$\,, this singular point is an ordinary double point and we have \,$\calD_{(\lambda_k,\mu_k)} = \wh{\calO}_{(\lambda_k,\mu_k)}$\,. Here \,$\wh{\calO}$\, denotes the direct image in \,$\Sigma$\, 
of the sheaf of holomorphic functions on the normalization of \,$\Sigma$\,. 
\end{enumerate}
\item
Let \,$D=\{(\lambda_k,\mu_k)\} \in \Div_{tame}$\,.
Then there exists one and only one holomorphic function \,$\Delta:\C^* \to \C$\, with \,$\Delta-\Delta_0 \in \As(\C^*,\ell^2_{0,0},1)$\, so that the hyperelliptic complex curve
\begin{equation}
\label{eq:special:phoenixfromashes:phoenix}
\Sigma := \Mengegr{(\lambda,\mu)\in \C^*\times \C^*}{\mu^2 - \Delta(\lambda)\cdot \mu+1=0} 
\end{equation}
contains \,$D$\,, and there exists one and only one tame generalized divisor \,$\calD$\, on \,$\Sigma$\, so that \,$D$\, is the support of \,$\calD$\,.
\end{enumerate}
\end{prop}

\begin{proof}
\emph{For (1)(a).} This is obvious.

\emph{For (1)(b).} For any \,$\lambda_* \in \C^*$\,, the meromorphic function \,$(\lambda-\lambda_*)^{-1}$\, on \,$\Sigma$\, has total pole order \,$2$\, above \,$\lambda_*$\, (regardless of whether \,$\Sigma$\,
is regular or singular above \,$\lambda_*$\,), and therefore cannot be a section of the tame divisor \,$\calD$\,. Thus \,$\calD$\, is non-special by Theorem~\ref{T:special:special}(b)$\Rightarrow$(a).

\emph{For (1)(c).} We enumerate the zeros of the discriminant function \,$\Delta^2-4$\, by two sequences \,$(\vkap_{k,1})$\, and \,$(\vkap_{k,2})$\, as in Proposition~\ref{P:excl:basic}(1),
then there exists \,$N\in \N$\, so that \,$\vkap_{k,\nu} \in U_{k,\delta}$\, and \,$(\lambda_k,\mu_k) \in \wh{U}_{k,\delta}$\, 
holds for all \,$k\in \Z$\, with \,$|k|>N$\, and \,$\nu \in \{1,2\}$\,. Because the excluded domains \,$U_{k,\delta}$\, are disjoint, it follows
that \,$\Delta^2-4$\, can have only zeros of order \,$\leq 2$\, on \,$U_{k,\delta}$\, with \,$|k|>N$\,, and therefore any singularities of \,$\Sigma$\, in \,$\wh{U}_{k,\delta}$\, are then of order \,$2$\,,
i.e.~ordinary double points.

Now let \,$k\in \Z$\, be given with \,$|k|>N$\, so that the divisor point \,$(\lambda_k,\mu_k)$\, is a singular point of \,$\Sigma$\,. By the preceding argument, this singularity of \,$\Sigma$\, is then 
an ordinary double point, and hence the \,$\delta$-invariant of \,$\Sigma$\, at \,$(\lambda_k,\mu_k)$\, is \,$1$\,. Moreover, \,$\calD$\, is non-special by (1)(b), and therefore is the spectral
divisor of a monodromy \,$M(\lambda) \in \Mon$\, by Theorem~\ref{T:special:special}(a)$\Rightarrow$(c). Here the degree of \,$(\lambda_k,\mu_k)$\, in \,$\calD$\, is \,$1$\,, because \,$\calD$\, is tame,
and therefore it follows from Proposition~\ref{P:spectrum:dims}(3) that \,$j_0=1$\, and \,$s=0$\, holds, where \,$j_0$\, is the number from Proposition~\ref{P:spectrum:locally-free} and 
\,$s$\, is the maximal pole order occurring in \,$\calD_{(\lambda_k,\mu_k)}$\,. Because \,$j_0$\, thus equals the \,$\delta$-invariant of \,$(\lambda_k,\mu_k)$\,, we have \,$\calR_{(\lambda_k,\mu_k)} = \wh{\calO}_{(\lambda_k,\mu_k)}$\,
for the ring \,$\calR_{(\lambda_k,\mu_k)}$\, over which \,$\calD_{(\lambda_k,\mu_k)}$\, is locally free (see Proposition~\ref{P:spectrum:locally-free}), and because of \,$s=0$\, we have
\,$\calD_{(\lambda_k,\mu_k)} = \calR_{(\lambda_k,\mu_k)}$\,. Thus \,$\calD_{(\lambda_k,\mu_k)} = \wh{\calO}_{(\lambda_k,\mu_k)}$\, holds.

\emph{For (2).}
We have
$$ \frac12\,\left(\mu_k+\mu_k^{-1}\right)-\mu_{k,0} = \frac12\,\left( 1 - \frac{1}{\mu_k\,\mu_{k,0}} \right)\cdot (\mu_k-\mu_{k,0}) \;. $$
\,$\mu_k-\mu_{k,0} \in \ell^2_{0,0}(k)$\, and \,$\mu_{k,0}=(-1)^k$\, implies \,$\lim_{k\to\pm\infty} \tfrac12\bigr( 1-\tfrac{1}{\mu_k\,\mu_{k,0}} \bigr) =0$\, and therefore by the above calculation
\,$\frac12\,\left(\mu_k+\mu_k^{-1}\right)-\mu_{k,0} \in \ell^2_{0,0}(k)$\,. 
Because \,$D$\, is tame, the \,$\lambda_k$\, are pairwise different, therefore 
Proposition~\ref{P:interpolate:mu}(1),(2) shows that there exists one and only one holomorphic function \,$\tfrac12\,\Delta:\C^* \to \C$\, with \,$\tfrac12\,\Delta-\tfrac12\,\Delta_0 \in \As(\C^*,\ell^2_{0,0},1)$\, and
$$ \frac12\,\Delta(\lambda_k) = \frac12\,\left( \mu_k+\mu_k^{-1} \right) \qmq{for all \,$k\in \Z$\,.} $$
If we now define the complex curve \,$\Sigma$\, by Equation~\eqref{eq:special:phoenixfromashes:phoenix}, then we have \,$D \subset \Sigma$\,, and there is no other way to choose \,$\Delta$\, so that
this inclusion holds. We let \,$\calO$\, be the sheaf of holomorphic functions on \,$\Sigma$\,, and \,$\wh{\calO}$\, be the direct image of the sheaf of holomorphic functions on the normalization \,$\wh{\Sigma}$\, 
of \,$\Sigma$\,. 

By Proposition~\ref{P:interpolate:lambda}(1) there exists (up to sign) one and only one holomorphic function \,$c: \C^* \to \C$\, with \,$c-\tau\,c_0 \in \As_\infty(\C^*,\ell^2_{-1},1)$\, and
\,$c-\tau^{-1}\,c_0 \in \As_0(\C^*,\ell^2_1,1)$\, for some \,$\tau\in\C^*$\,, which has zeros at all the \,$\lambda_k$\,, \,$k\in \Z$\, and no others. 
And again because the \,$\lambda_k$\, are pairwise different, 
there exists one and only one holomorphic function \,$d: \C^*\to \C$\, with \,$d-d_0 \in \As(\C^*,\ell^2_{0,0},1)$\, such that \,$d(\lambda_k)=\mu_k^{-1}$\, holds for all \,$k\in \Z$\,, 
by Proposition~\ref{P:interpolate:mu}(1),(2).

We let \,$\calD$\, be the generalized divisor on \,$\Sigma$\, generated by \,$1$\, and \,$\tfrac{\mu-d}{c}$\,
over \,$\calO$\,. Clearly \,$\calD$\, is positive, and the support of \,$\calD$\, is \,$D$\,, hence \,$\calD$\, is asymptotic. 
Because \,$D$\, is tame, the total degree of \,$\calD$\, above any \,$\lambda_* \in \C^*$\, is at most \,$1$\,. 
Moreover \,$\calD$\, is compatible: If \,$(\lambda_*,\mu_*)\in D$\,
is a singular point of \,$\Sigma$\,, let \,$g := 2\,\tfrac{\mu-d}{\lambda-\lambda_*}$\,, then \,$g \in H^0(\Sigma,\calD)$\, holds and \,$\eta := g-\tfrac{\mu-\mu^{-1}}{\lambda-\lambda_*} = \tfrac{\Delta-2d}{\lambda-\lambda_*}$\,
is a meromorphic function in \,$\lambda$\, with a pole of order at most \,$1$\,. Therefore it follows that \,$\calD$\, is tame.
%
\end{proof}

By virtue of Proposition~\ref{P:special:phoenixfromashes}, the tame classical divisors are in one-to-one correspondence with spectral data \,$(\Sigma,\calD)$\,, where \,$\Sigma$\, is a spectral curve and \,$\calD$\, is a tame,
generalized divisor on \,$\Sigma$\,. 
In the sequel, we will therefore speak of the spectral curve \,$\Sigma$\, and trace function \,$\Delta=\mu+\mu^{-1}$\,, and of the tame generalized divisor \,$\calD$\, on \,$\Sigma$\, 
associated to a tame classical divisor \,$D\in\Div_{tame}$\, (regarded as a multi-set of points in \,$\C^* \times \C^*$\,).
In this manner, any tame classical divisor \,$D$\, gives rise to a monodromy \,$M(\lambda)$\, with spectral divisor \,$D$\, by Theorem~\ref{T:special:special}.

\part{The inverse problem for periodic potentials (Cauchy data)}

\section{Divisors of finite type}
\label{Se:finite}

An asymptotic divisor \,$\calD$\, on a spectral curve \,$\Sigma$\, is said to be \emph{of finite type}, if the following conditions hold:
\begin{enumerate}
\item The spectral curve \,$\Sigma$\, has finite geometric genus (i.e.~only finitely many of the double points of the 
spectral curve of the vacuum have ``opened up'' into a pair of branch points with positive distance).
\item All but finitely many of the points \,$(\lambda_*,\mu_*)$\, in the support of \,$\calD$\, lie in double points of \,$\Sigma$\,, and we have \,$\calD_{(\lambda_*,\mu_*)} = \wh{\calO}_{(\lambda_*,\mu_*)}$\,,
where \,$\wh{\calO}$\, denotes the direct image in \,$\Sigma$\, of the sheaf of holomorphic functions on the normalization of \,$\Sigma$\,. 
\end{enumerate}
Spectral data \,$(\Sigma,\calD)$\, of finite type thus look like the spectral data of the vacuum near all but finitely many of the divisor points, see Section~\ref{Se:vacuum}. 
In particular the equation \,$\calD_{(\lambda_*,\mu_*)} = \wh{\calO}_{(\lambda_*,\mu_*)}$\, implies that the ring \,$\calR_{(\lambda_*,\mu_*)}$\, from Proposition~\ref{P:spectrum:locally-free} over which
\,$\calD_{(\lambda_*,\mu_*)}$\, is locally free equals \,$\wh{\calO}_{(\lambda_*,\mu_*)}$\,; this shows that a generalized divisor \,$\calD$\, of finite type is locally free in the normalization \,$\wh{\Sigma}$\,
of \,$\Sigma$\, with the exception of finitely many points.

Among the asymptotic divisors, those that are of finite type play a special role, because the spectral data associated to solutions \,$u$\, of the sinh-Gordon equation that are doubly periodic,
meaning that there is a number \,$\tau\in\C$\, with \,$\IM(\tau)>0$\, such that \,$u(z+1)=u(z+\tau)=u(z)$\, holds for all \,$z\in \C$\,, are of finite type. (Conversely though, not every solution \,$u$\,
associated to spectral data of finite type is doubly periodic, because an additional periodicity condition has to be fulfilled.) With an additional closing condition, doubly periodic solutions of the
sinh-Gordon equation correspond to constant mean curvature tori in \,$\R^3$\,, or to minimal tori in \,$S^3$\,.
The finite type solutions \,$u$\, of the sinh-Gordon equation, and consequently the minimal tori, have been famously classified by \textsc{Pinkall}/\textsc{Sterling} 
(\cite{Pinkall/Sterling:1989}), and independently by \textsc{Hitchin} (\cite{Hitchin:1990}).

For many other integrable systems, where one can similarly define the concept of a finite type potential, the finite type potentials are known to be dense in the space of all potentials.
For example, this is true for the integrable system associated to the KdV equation (where the ``finite type potentials'' are called finite-gap solutions), 
see \cite{Kappeler/Poeschel:2003}, Section~11. Therefore we expect the finite type potentials to be dense
in the space of all potentials also in our setting for the sinh-Gordon equation, and likewise that the set of divisors of finite type is dense in the space of all asymptotic potentials. 

It is the objective of the present section to prove the latter statement, i.e.~that the divisors of finite type are indeed dense in \,$\Div$\,. 
Later, in Corollary~\ref{C:diffeo:dense}(1), we will see that also the set of finite type potentials in \,$\Pot_{tame}$\, is dense in \,$\Pot_{tame}$\,.

We will see that it suffices to consider divisors of finite type that are tame, because the tame divisors themselves comprise an open and dense set in \,$\Div$\,. 
Proposition~\ref{P:special:phoenixfromashes}(1)(c) shows that for a tame generalized divisor
\,$\calD$\, and a point \,$(\lambda_k,\mu_k)$\, in the support of \,$\calD$\, with \,$|k|$\, large that is a singular point of the associated spectral curve \,$\Sigma$\,, we already have
that the singularity \,$(\lambda_k,\mu_k)$\, is an ordinary double point of \,$\Sigma$\, and that \,$\calD_{(\lambda_k,\mu_k)} = \wh{\calO}_{(\lambda_k,\mu_k)}$\, holds. Moreover by Proposition~\ref{P:special:phoenixfromashes}
the generalized tame divisors are in one-to-one correspondence with classical tame divisors. Therefore the definition of classical finite type divisors given below is for tame
divisors equivalent to the the definition of finite type above.

\begin{Def}
\label{D:finite:finite}
\begin{enumerate}
\item Let \,$D=\Menge{(\lambda_k,\mu_k)}{k\in \Z}$\, be an asymptotic divisor on a spectral curve \,$\Sigma$\, with trace function \,$\Delta = \mu+\mu^{-1}$\,, and let \,$\vkap_{k,\nu}$\, be the zeros of
\,$\Delta^2-4$\, (see Proposition~\ref{P:excl:basic}(1)). We say that \,$D$\, is \emph{of finite type}, if \,$\vkap_{k,1}=\vkap_{k,2}=\lambda_k$\, holds for all \,$k\in \Z$\, with at most finitely 
many exceptions.
\item We say that a potential \,$(u,u_y)\in\Pot$\, is \emph{of finite type}, if its associated spectral divisor \,$D$\, is of finite type. 
\end{enumerate}
\end{Def}

As promised we will show in the present section that the set of finite type divisors is dense in \,$\Div$\,. In fact, the following theorem states that the set of tame divisors of finite type
is dense in \,$\Div_{tame}$\,. 
Because \,$\Div_{tame}$\, is open and dense in \,$\Div$\,, it follows from this theorem that the finite type divisors are dense in \,$\Div$\,. 

\begin{thm}
\label{T:finite:finite}
Let \,$D=\{(\lambda_k,\mu_k)\}\in\Div_{tame}$\, be given. Then for every \,$\eps>0$\, there exists a divisor of finite type \,$D^*=\{(\lambda_k^*,\mu_k^*)\} \in \Div_{tame}$\, with
$$ \|D^*-D\|_{\Div} \leq \eps \; . $$
Moreover, for given \,$N\in \N$\,, \,$D^*$\, can be chosen such that
\begin{equation}
\label{eq:finite:finite:fixed-cond}
\forall \, k\in \Z, \; |k|\leq N \; : \; \lambda_k^*=\lambda_k\;,\;\;\mu_k^* = \mu_k
\end{equation}
holds. 
\end{thm}

We prepare the proof of Theorem~~\ref{T:finite:finite} with two lemmas. 
For a given holomorphic function \,$\Delta$\, 
with \,$\Delta-\Delta_0 \in \As(\C^*,\ell^2_{0,0},1)$\, (intended to be the trace function of a spectral curve \,$\Sigma$\,), the first lemma (Lemma~\ref{L:finite:eta}) studies the asymptotic behavior 
of the sequence \,$(\eta_n)$\, of zeros of \,$\Delta'$\,. These zeros are of interest in relation to potentials of finite type because if \,$\Sigma$\, has an ordinary double point 
(or, in fact, any other singularity) above some \,$\lambda_*\in \C^*$\,, then \,$\Delta'(\lambda_*)=0$\, holds.
The second lemma (Lemma~\ref{L:finite:Delta}) provides a Lipschitz estimate for the quantity \,$\Delta(\eta_n)-\Delta(\lambda_n)$\,
where \,$(\lambda_k)$\, is a fixed sequence, but \,$\Delta$\, and correspondingly \,$(\eta_n)$\, varies. Both lemmas will be used in the proof of Theorem~\ref{T:finite:finite} to set up a
fixed point equation to describe the property of a divisor being of finite type, and to show that the Banach fixed point theorem applies to this equation.

\newpage

\begin{lem}
\label{L:finite:eta}
\begin{enumerate}
\item
Let \,$\Delta:\C^* \to \C$\, be a holomorphic function with \,$\Delta-\Delta_0 \in \As(\C^*,\ell^2_{0,0},1)$\,. Then \,$\Delta'$\, 
has asymptotically and totally a zero in each excluded domain \,$U_{k,\delta}$\,, and besides these exactly one additional zero.

If \,$\eta_*\in \C^*$\, denotes a zero of \,$\Delta'$\, for which \,$|\eta_*-1|$\, becomes minimal, and 
if we denote by \,$(\eta_k)_{k\in \Z}$\, the sequence of the remaining zeros of \,$\Delta'$\,, where \,$\eta_k \in U_{k,\delta}$\, holds for \,$|k|$\, large,
then \,$\eta_k - \lambda_{k,0} \in \ell^2_{-1,3}(k)$\, holds. 

In the sequel, we will denote the sequence of all zeros of \,$\Delta'$\, by \,$(\eta_k)_{k\in \Z\cup \{*\}}$\,. 

\item
Let \,$R_0>0$\,, let  \,$\Delta^{[1]},\Delta^{[2]}:\C^* \to \C$\, be two holomorphic functions with
\,$\Delta^{[\nu]}-\Delta_0 \in \As(\C^*,\ell^2_{0,0},1)$\, and \,$\|\Delta^{[\nu]}-\Delta_0\|_{\As(\C^*,\ell^2_{0,0},1)} \leq R_0$\, for \,$\nu\in\{1,2\}$\,.
We denote the zeros of \,$\Delta^{[\nu]}{}'$\, as in (1) by \,$(\eta_k^{[\nu]})_{k\in \Z\cup \{*\}}$\,. 
Then there exists a constant \,$C>0$\,, dependent only on \,$R_0$\,, such that 
$$ |\eta_k^{[1]}-\eta_k^{[2]}| \leq C \cdot \left\{ \begin{matrix} k & \text{if \,$k>0$\,} \\ |k|^{-3} & \text{if \,$k<0$\,} \end{matrix}  \right\}  \cdot  \max\{b_{k-1},b_k,b_{k+1}\} $$
holds, where \,$(b_k)_{k\in \Z}\in \ell^2_{0,0}(k)$\, is a bounding sequence for \,$\Delta^{[1]}-\Delta^{[2]}$\, in \,$\As(\C^*,\ell^2_{0,0},1)$\,. 
\end{enumerate}
\end{lem}

\begin{proof}
%
\emph{For (1).}
Because of \,$\Delta-\Delta_0 \in \As(\C^*,\ell^2_{0,0},1)$\,, we have by Proposition~\ref{P:interpolate:l2asymp}(3)
\begin{equation}
\label{eq:finite:eta:Delta-asymp}
\Delta'-\Delta_0' \in \As(\C^*,\ell^2_{1,-3},1) \qmq{and} \Delta''-\Delta_0'' \in \As(\C^*,\ell^2_{2,-6},1)
\end{equation}
Moreover, we have
\begin{equation}
\label{eq:finite:eta:Delta0-estim}
\Delta_0' \in \As(\C^*,\ell^\infty_{1,-3},1) \qmq{and} \Delta_0'' \in \As(\C^*,\ell^\infty_{2,-6},1)
\end{equation}
and therefore also 
\begin{equation}
\label{eq:finite:eta:Delta-estim}
\Delta' \in \As(\C^*,\ell^\infty_{1,-3},1) \qmq{and} \Delta'' \in \As(\C^*,\ell^\infty_{2,-6},1) \; . 
\end{equation}

We begin by showing that \,$\Delta'$\, has at least one zero \,$\eta_*$\,. 
It follows from \eqref{eq:finite:eta:Delta-asymp}, \eqref{eq:finite:eta:Delta0-estim} and \eqref{eq:finite:eta:Delta-estim} that we have
$$ \left. \left( \frac{\Delta''}{\Delta'}-\frac{\Delta_0''}{\Delta_0'}\right)\right|_{V_\delta} = \left. \frac{(\Delta''-\Delta_0'')\cdot \Delta_0' + \Delta_0'' \cdot (\Delta_0'-\Delta')}{\Delta' \cdot \Delta_0'} \right|_{V_\delta}
\;\in\; \As(\C^*,\ell^2_{1,-3},1) \;; $$
because the sequence of circumferences of \,$U_{k,\delta}$\, is in \,$\ell^\infty_{-1,3}(k)$\,, this implies
$$ \int_{\partial U_{k,\delta}} \left( \frac{\Delta''}{\Delta'}-\frac{\Delta_0''}{\Delta_0'} \right)\,\mathrm{d}\lambda \;\in\; \ell^2_{0,0}(k) \; . $$
Because \,$\int_{\partial U_{k,\delta}} \left( \frac{\Delta''}{\Delta'}-\frac{\Delta_0''}{\Delta_0'} \right)\,\mathrm{d}\lambda$\, is the difference of the ``zero counting integrals'' for \,$\Delta'$\, and for \,$\Delta_0'$\,,
this integral is an integer multiple of \,$2\pi i$\,, and therefore equals zero for large \,$|k|$\,. Thus we see that for \,$|k|$\, large, \,$\Delta'$\, and \,$\Delta_0'$\, have the same number of zeros on \,$U_{k,\delta}$\,,
i.e.~one. In particular, \,$\Delta'$\, has at least one zero, and we let \,$\eta_* \in \C^*$\, be the zero of \,$\Delta'$\, for which \,$|\eta_*-1|$\, becomes minimal. (If there are several such zeros, we make an
arbitrary choice.)

We now put
\begin{equation}
\label{eq:finite:eta:cdef}
\tau := \sqrt{\eta_*} \in \C^* \qmq{and} c(\lambda) := -8\tau\,\frac{\lambda^2}{\lambda-\eta_*}\,\Delta'(\lambda) \;.
\end{equation}
Because \,$\eta_*$\, is a zero of \,$\Delta'$\,, \,$c(\lambda)$\, is a holomorphic function on \,$\C^*$\,, and we claim that we have
\begin{equation}
\label{eq:finite:eta:c-asymp}
c-\tau\,c_0 \in \As_\infty(\C^*,\ell^2_{-1},1) \qmq{and} c-\tau^{-1}\,c_0 \in \As_0(\C^*,\ell^2_1,1) \; . 
\end{equation}
Indeed, because of \,$\Delta_0(\lambda) = \cos(\zeta(\lambda))$\, and \,$c_0(\lambda)=\sqrt{\lambda}\,\sin(\zeta(\lambda))$\,, we have
$$ \Delta_0'(\lambda) = -\frac{1}{8}\,\left( \lambda^{-1/2}-\lambda^{-3/2} \right)\,\sin\left(\zeta(\lambda)\right) = -\frac{\lambda-1}{8\lambda^2}\,c_0(\lambda) $$
and therefore
$$ c_0(\lambda) = -\frac{8\lambda^2}{\lambda-1}\,\Delta_0'(\lambda) \; . $$
We thus obtain via \eqref{eq:finite:eta:Delta-asymp} and \eqref{eq:finite:eta:Delta0-estim}
$$ c-\tau\,c_0 = \frac{-8\,\tau\,\lambda^2}{\lambda-\eta_*}\cdot(\Delta'-\Delta_0') + \frac{-8\,\tau\,(\eta_*-1)\,\lambda^2}{(\lambda-\eta_*)\cdot(\lambda-1)}\cdot \Delta_0' \;\in\; \As_\infty(\C^*,\ell^2_{-1},1) $$
and
$$ c-\tau^{-1}\,c_0 = \frac{-8\,\tau^{-1}\,\eta_*\,\lambda^2}{\lambda-\eta_*}\cdot (\Delta'-\Delta_0') + \frac{-8\,\tau^{-1}\,(\eta_*-1)\,\lambda^3}{(\lambda-\eta_*)\cdot(\lambda-1)}\cdot \Delta_0'
\;\in\; \As_0(\C^*,\ell^2_1,1) \; . $$

By Proposition~\ref{P:excl:basic}(2) it follows from \eqref{eq:finite:eta:c-asymp} that the holomorphic function \,$c(\lambda)$\, has asymptotically and totally exactly one zero in every excluded domain \,$U_{k,\delta}$\,,
and that the sequence \,$(\eta_k)_{k\in \Z}$\, of these zeros satisfies \,$\eta_k - \lambda_{k,0} \in \ell^2_{-1,3}(k)$\,. Because of \,$\Delta'(\lambda)=\tfrac{\lambda-\eta_*}{-8\,\tau\,\lambda^2}\cdot c(\lambda)$\,,
it follows that \,$\eta_*$\, and the \,$(\eta_k)_{k\in \Z}$\, are all the zeros of \,$\Delta'$\,. 


\emph{For (2).} In the setting of (2), we denote the objects defined above in relation to \,$\Delta^{[\nu]}$\, by the superscript \,${}^{[\nu]}$\, (\,$\nu\in\{1,2\}$\,). If \,$(b_k)$\, is a bounding sequence
for \,$\Delta^{[1]}-\Delta^{[2]} \in \As(\C^*,\ell^2_{0,0},1)$\,, then by Proposition~\ref{P:interpolate:l2asymp}(3) there exists a constant \,$C_1>0$\, so that 
$$ \begin{cases} \tfrac{C_1}{k}\cdot \max\{b_{k-1},b_k,b_{k+1}\} & \text{for \,$k>0$\,} \\ C_1 \cdot k^3 \cdot \max\{b_{k-1},b_k,b_{k+1}\} & \text{for \,$k<0$\,} \end{cases} $$
is a bounding sequence for \,$\Delta^{[1]}{}'-\Delta^{[2]}{}'\in \As(\C^*,\ell^2_{1,-3},1)$\,, and therefore with another constant \,$C_2>0$\, and the sequence
$$ \wt{b}_k := \begin{cases} C_2\cdot k\cdot \max\{b_{k-1},b_k,b_{k+1}\} & \text{for \,$k>0$\,} \\ \tfrac{C_2}{k} \cdot \max\{b_{k-1},b_k,b_{k+1}\} & \text{for \,$k<0$\,} \end{cases} \;,$$
\,$(\wt{b}_k)_{k>0}$\, is a bounding sequence for \,$\tau^{[2]}\,c^{[1]}-\tau^{[1]}\,c^{[2]}\in \As_\infty(\C^*,\ell^2_{-1},1)$\, and 
\,$(\wt{b}_k)_{k<0}$\, is a bounding sequence for \,$\tau^{[1]}\,c^{[1]}-\tau^{[2]}\,c^{[2]}\in \As_0(\C^*,\ell^2_{1},1)$\,. By application of Proposition~\ref{P:asympdiv:asympdiv-neu}(1), the claimed statement follows.
%
\end{proof}

\newpage

\begin{lem}
\label{L:finite:Delta}
Let \,$R_0>0$\,, and let \,$(\lambda_k)_{k\in \Z}$\, with \,$\lambda_k-\lambda_{k,0} \in \ell^2_{-1,3}(k)$\, and \,$\|\lambda_k-\lambda_{k,0}\|_{\ell^2_{-1,3}} \leq R_0$\, be given.
Further suppose that for \,$\nu\in\{1,2\}$\,,
\,$\Delta^{[\nu]}: \C^* \to \C$\, is a holomorphic function with \,$\Delta^{[\nu]}-\Delta_0 \in \As(\C^*,\ell^2_{0,0},1)$\, and \,$\|\Delta^{[\nu]}-\Delta_0\|_{\As(\C^*,\ell^2_{0,0},1)} \leq R_0$\,. 
We put \,$w_k^{[\nu]} := \Delta^{[\nu]}(\lambda_k)$\, and let 
\,$(\eta_k^{[\nu]})_{k\in \Z\cup \{*\}}$\, be the sequence of zeros of \,$(\Delta^{[\nu]})'$\, as in Lemma~\ref{L:finite:eta}(1). Then there exists a constant \,$C>0$\,, depending only on \,$R_0$\,, 
so that we have for \,$n\in \Z$\,
\begin{gather*}
\left| \left( \Delta^{[1]}(\eta_n^{[1]}) - \Delta^{[1]}(\lambda_n)\right) \;-\; \left( \Delta^{[2]}(\eta_n^{[2]}) - \Delta^{[2]}(\lambda_n)\right) \right|   \\
\leq C \cdot \left( |\zeta(\eta_n^{[1]})-\zeta(\lambda_n)|\cdot \left( |w_k^{[1]}-w_k^{[2]}| * \frac{1}{|k|} \right)_n + |\zeta(\eta_n^{[1]})-\zeta(\eta_n^{[2]})|^2 \right) \; .
\end{gather*}
\end{lem}

\begin{proof}
In the sequel, we always have \,$\nu\in\{1,2\}$\,, and all \,$C_k$\, are constants \,$>0$\, that depend only on \,$R_0$\,. We have
\begin{align}
& \left( \Delta^{[1]}(\eta_n^{[1]}) - \Delta^{[1]}(\lambda_n)\right) \;-\; \left( \Delta^{[2]}(\eta_n^{[2]}) - \Delta^{[2]}(\lambda_n)\right) \notag \\
= & \int_{\lambda_n}^{\eta_n^{[1]}}(\Delta^{[1]})'(\lambda)\,\mathrm{d}\lambda - \int_{\lambda_n}^{\eta_n^{[2]}}(\Delta^{[2]})'(\lambda)\,\mathrm{d}\lambda \notag \\
\label{eq:finite:Delta:master}
= & \int_{\lambda_n}^{\eta_n^{[1]}}(\Delta^{[1]}-\Delta^{[2]})'(\lambda)\,\mathrm{d}\lambda + \int_{\eta_n^{[2]}}^{\eta_n^{[1]}}(\Delta^{[2]})'(\lambda)\,\mathrm{d}\lambda \; . 
\end{align}
We will estimate the two summands in the last expression individually.

For the first term in \eqref{eq:finite:Delta:master}, we have
\begin{align}
\left| \int_{\lambda_n}^{\eta_n^{[1]}}(\Delta^{[1]}-\Delta^{[2]})'(\lambda)\,\mathrm{d}\lambda \right| 
& \leq |\eta_n^{[1]}-\lambda_n| \cdot \max_{\lambda\in U_{n,\delta}} \left| \left( \Delta^{[1]}-\Delta^{[2]} \right)'(\lambda) \right| \notag \\
& \overset{(*)}{\leq} |\eta_n^{[1]}-\lambda_n|\cdot \ell^\infty_{1,-3}(n)\cdot \max_{\lambda\in U_{n,\delta}} \left| \left( \Delta^{[1]}-\Delta^{[2]} \right)(\lambda) \right| \notag \\
\label{eq:finite:Delta:Apre}
& \overset{(\dagger)}{\leq} C_1 \cdot |\zeta(\eta_n^{[1]})-\zeta(\lambda_n)| \cdot \max_{\lambda\in U_{n,\delta}} \left| \left( \Delta^{[1]}-\Delta^{[2]} \right)(\lambda) \right| \; , 
\end{align}
where the estimate marked $(*)$ follows from Cauchy's inequality, 
and the estimate marked $(\dagger)$ follows from Proposition~\ref{P:vac2:excldom-new}(1).

To estimate \,$\left( \Delta^{[1]}-\Delta^{[2]} \right)(\lambda)$\,, we apply Proposition~\ref{P:interpolate:mu}(3): In the setting described there, we choose \,$\lambda_k^{[1]}=\lambda_k^{[2]}=\lambda_k$\,
and \,$\mu_k^{[\nu]}=\tfrac12 w_k^{[\nu]}$\,. Then we have \,$\tau^{[1]}=\tau^{[2]}$\,, \,$\upsilon^{[1]}=\upsilon^{[2]}=1$\,, \,$a_k=0$\, and \,$b_k = \tfrac12\,|w_k^{[1]}-w_k^{[2]}|$\,. It follows 
from Proposition~\ref{P:interpolate:mu}(3)(b) that there exists \,$C_2>0$\, so that
$$ r_n := C_2\cdot \left( |w_k^{[1]}-w_k^{[2]}| * \frac{1}{|k|} \right)_n \in \ell^2_{0,0}(n) $$
is a bounding sequence for \,$\Delta^{[1]}-\Delta^{[2]}$\, in \,$\As(\C^*,\ell^2_{0,0},1)$\,; here we set \,$\tfrac10:=1$\, as usual. Because \,$w(\lambda)$\, is bounded on \,$U_\delta$\,, there exists \,$C_3>0$\,
so that 
$$ \max_{\lambda\in U_{n,\delta}} \left| \left( \Delta^{[1]}-\Delta^{[2]} \right)(\lambda) \right| \leq C_3 \cdot r_n \; . $$
By plugging this estimate into \eqref{eq:finite:Delta:Apre} we obtain
\begin{equation}
\label{eq:finite:Delta:A}
\left| \int_{\lambda_n}^{\eta_n^{[1]}}(\Delta^{[1]}-\Delta^{[2]})'(\lambda)\,\mathrm{d}\lambda \right| 
\leq C_4 \cdot |\zeta(\eta_n^{[1]})-\zeta(\lambda_n)|\cdot \left( |w_k^{[1]}-w_k^{[2]}| * \frac{1}{|k|} \right)_n 
\end{equation}
with \,$C_4 := C_3 \cdot C_2 \cdot C_1$\,.

We now turn our attention to the second summand in \eqref{eq:finite:Delta:master}. 
We have
\begin{equation}
\label{eq:finite:Delta:Bpre}
\int_{\eta_n^{[2]}}^{\eta_n^{[1]}} (\Delta^{[2]})'(\lambda)\,\mathrm{d}\lambda = \int_{\eta_n^{[2]}}^{\eta_n^{[1]}} g_n(\lambda)\cdot (\lambda-\eta_n^{[2]})\,\mathrm{d}\lambda \;,
\end{equation}
where the function 
$$ g_n(\lambda) := \frac{(\Delta^{[2]})'(\lambda)}{\lambda-\eta_n^{[2]}} $$
is holomorphic because \,$(\Delta^{[2]})'$\, has a zero at \,$\eta_n^{[2]}$\,. We now define (compare the proof of Lemma~\ref{L:finite:eta}(1))
$$ \tau := \sqrt{\eta^{[2]}_*} \in \C^* \qmq{and} c(\lambda) := -8\tau\,\frac{\lambda^2}{\lambda-\eta_*^{[2]}}\,(\Delta^{[2]})'(\lambda) \;, $$
then we have \,$c-\tau\,c_0 \in \As_\infty(\C^*,\ell^2_{-1},1)$\,, \,$c-\tau^{-1}\,c_0 \in \As_0(\C^*,\ell^2_1,1)$\, and 
$$ g_n(\lambda) = -\frac{1}{8\tau}\,\frac{\lambda-\eta^{[2]}_*}{\lambda^2}\,\frac{c(\lambda)}{\lambda-\eta_n^{[2]}} \; . $$
For \,$\lambda \in U_{n,\delta}$\, we have by Corollary~\ref{C:interpolate:cdivlin}(1): \,$\left| \tfrac{c(\lambda)}{\lambda-\eta_n^{[2]}} \right| \in \ell^\infty_{0,-2}(n)$\,,
and we also have \,$\left| \tfrac{\lambda-\eta^{[2]}_*}{\lambda^2} \right| \in \ell^\infty_{2,-4}(n)$\,. Thus we obtain
$$ |g_n(\lambda)| \in \ell^\infty_{2,-6}(n) \qmq{for \,$\lambda\in U_{n,\delta}$\,} \;, $$
and therefore by Equation~\eqref{eq:finite:Delta:Bpre}
\begin{align}
\left| \int_{\eta_n^{[2]}}^{\eta_n^{[1]}} (\Delta^{[2]})'(\lambda)\,\mathrm{d}\lambda \right|
& \leq |\eta_n^{[2]}-\eta_n^{[1]}| \cdot \ell^\infty_{2,-6}(n) \cdot \max_{\lambda\in [\eta_n^{[1]},\eta_n^{[2]}]} |\lambda-\eta_n^{[2]}| \notag \\
\label{eq:finite:Delta:B}
& \leq \ell^\infty_{2,-6}(n) \cdot |\eta_n^{[2]}-\eta_n^{[1]}|^2 \overset{(*)}{\leq} C_5 \cdot |\zeta(\eta_n^{[2]})-\zeta(\eta_n^{[1]})|^2  \; ,
\end{align}
where $(*)$ again follows from Proposition~\ref{P:vac2:excldom-new}(1).

By taking the absolute value in Equation~\eqref{eq:finite:Delta:master} and then applying the estimates \eqref{eq:finite:Delta:A} and \eqref{eq:finite:Delta:B}, we obtain the claimed statement.
\end{proof}

We are now ready to prove Theorem~\ref{T:finite:finite}:

\begin{proof}[Proof of Theorem~\ref{T:finite:finite}.]
The following construction depends on a number \,$N\in \N$\,. We will see that if \,$N$\, is chosen large enough, then this construction will yield a divisor of finite type \,$D^*$\, so that
\eqref{eq:finite:finite:fixed-cond} holds. 

The idea of the proof is as follows: We need to construct the trace function \,$\Delta(\lambda)$\, (with the asymptotic \,$\Delta-\Delta_0 \in \As(\C^*,\ell^2_{0,0},1)$\,)
for a spectral curve \,$\Sigma$\, such that \,$(\lambda_k,\mu_k)\in \Sigma$\, holds for \,$|k|\leq N$\, and so that \,$\Sigma$\, has a double point in each excluded domain \,$\wh{U}_{k,\delta}$\, with \,$|k|>N$\,.
The latter condition means: We have \,$\Delta(\eta_k)=2(-1)^k$\, for \,$|k|>N$\,, where \,$(\eta_k)_{k\in \Z \cup \{*\}}$\, is the sequence of zeros of \,$\Delta'$\, as in Lemma~\ref{L:finite:eta}(1).
Because the \,$\lambda_k$\, are pairwise unequal (\,$D$\, being tame), Proposition~\ref{P:interpolate:mu} shows that we can uniquely determine
a trace function \,$\Delta$\, with \,$\Delta-\Delta_0 \in \As(\C^*,\ell^2_{0,0},1)$\, by prescribing the values of \,$\Delta$\, at the points \,$\lambda_k$\, as
$$ \Delta(\lambda_k) = 2(-1)^k + z_k $$
with a sequence \,$(z_k) \in \ell^2_{0,0}(k)$\,. \,$\Delta'$\, has a zero in each excluded domain \,$U_{k,\delta}$\, with \,$|k|$\, large, and therefore \,$\Delta$\, is approximately constant on 
\,$U_{k,\delta}$\,. If some \,$\Delta$\, defined via any sequence \,$(z_k)\in \ell^2_{0,0}(k)$\, is given, we can therefore expect to decrease \,$|\Delta(\eta_k)-2(-1)^k|$\, by passing
from \,$(z_k)$\, to \,$(\wt{z}_k)$\, defined by
$$ \wt{z}_k = z_k - \big(\Delta(\eta_k)-2(-1)^k\big) = \Delta(\lambda_k)-\Delta(\eta_k) \; . $$
In particular the desired equality \,$\Delta(\eta_k)=2(-1)^k$\, is equivalent to the fixed point equation \,$\wt{z}_k = z_k$\,. The following proof shows that (for \,$N$\, chosen suitably large), the iteration map
\,$(z_k) \mapsto (\wt{z}_k)$\, defines a contraction with respect to the \,$\ell^2$-norm, and thus has a unique fixed point by Banach's Fixed Point Theorem. 

To carry out the described idea, we consider the
Banach space \,$\mathfrak{B}_N$\, of sequences \,$(z_k)$\,, where the index \,$k$\, runs through all the integers with \,$|k|>N$\,, equipped with the \,$\ell^2$-norm. To each member
\,$(z_k)\in \mathfrak{B}_N$\, we associate the holomorphic function \,$\Delta: \C^*\to \C$\, with \,$\Delta-\Delta_0 \in \As(\C^*,\ell^2_{0,0},1)$\, and
$$ \Delta(\lambda_k) = \begin{cases}
\mu_k + \mu_k^{-1} & \text{for \,$|k|\leq N$\,} \\
2(-1)^k + z_k & \text{for \,$|k|>N$\,} \end{cases} $$
for all \,$k\in \Z$\,; existence of \,$\Delta$\, follows from Proposition~\ref{P:interpolate:mu}(1) (applied to \,$\tfrac12 \Delta$\, with \,$\upsilon=1$\,), and \,$\Delta$\, is uniquely determined by these
equations because \,$D$\, is tame, and hence the \,$\lambda_k$\, are pairwise unequal (Proposition~\ref{P:interpolate:mu}(2)).
In this situation, we also have the
sequence \,$(\eta_k)_{k\in \Z \cup \{*\}}$\, of the zeros of \,$\Delta'$\, as in Lemma~\ref{L:finite:eta}(1); the sequence satisfies \,$\eta_k-\lambda_{k,0}\in \ell^2_{-1,3}(k)$\,. 

We now define for given \,$(z_k)\in \mathfrak{B}_N$\, and the associated objects \,$\Delta$\,, \,$(\eta_k)$\, a new sequence \,$(\wt{z}_k)_{|k|>N}$\, by
$$ \wt{z}_k := z_k - \bigr( \Delta(\eta_k) - 2\,(-1)^k \bigr) = \Delta(\lambda_k)-\Delta(\eta_k) \; . $$
We then have for \,$|k|>N$\,
$$ |\wt{z}_k| = |\Delta(\lambda_k)-\Delta(\eta_k)| \leq \int_{\eta_k}^{\lambda_k} \underbrace{|\Delta'(\lambda)|}_{\in \As(\C^*,\ell^\infty_{1,-3},1)} \,\mathrm{d}\lambda \leq \ell^\infty_{1,-3}(k) \cdot \underbrace{|\lambda_k-\eta_k|}_{\in \ell^2_{-1,3}(k)} \;\in\; \ell^2_{0,0}(k) \;, $$
and therefore we have \,$(\wt{z}_k) \in \mathfrak{B}_N$\,. 
Hence, the map
$$ \Phi: \mathfrak{B}_N \to \mathfrak{B}_N,\; (z_k) \mapsto (\wt{z}_k) $$
is well-defined. 

Before we show that for sufficiently large \,$N$\,, \,$\Phi$\, is a contraction on some small closed ball in \,$\mathfrak{B}_N$\,, we look at what happens if we apply \,$\Phi$\, to the sequence \,$z_k=0 \in \mathfrak{B}_N$\,.
It should be noted that the function \,$\Delta^{[0]}$\, associated to this sequence depends on the choice of \,$N$\,, as do the associated sequences
\,$(\eta_k^{[0]})$\, and \,$(\wt{z}_n) := \Phi(0) \in \mathfrak{B}_N$\,. First we note that \,$\Delta^{[0]}$\, is characterized by
$$ \frac12\,\Delta^{[0]}(\lambda_k) = \mu_k^{[0]} := \begin{cases}
\tfrac12(\mu_k + \mu_k^{-1}) & \text{for \,$|k|\leq N$\,} \\
\tfrac12(\mu_{k,0}+\mu_{k,0}^{-1})=(-1)^k & \text{for \,$|k|>N$\,} \end{cases} \;, $$
and therefore, if we put \,$R_0 := \max\{\|\lambda_k-\lambda_{k,0}\|_{\ell^2_{-1,3}},\,\|\tfrac12(\mu_k+\mu_k^{-1})-\mu_{k,0}\|_{\ell^2_{0,0}}\}$\, 
(the value of this constant depends on the divisor \,$D$\,, but not on \,$N$\,), we have
$$ \|\lambda_k-\lambda_{k,0}\|_{\ell^2_{-1,3}} \leq R_0 \qmq{and}  \|\mu_k^{[0]}-\mu_{k,0}\|_{\ell^2_{0,0}} \leq R_0 \;, $$
therefore it follows from Proposition~\ref{P:interpolate:mu}(3)(c) that there exists a constant \,$C_1>0$\, (again depending on \,$D$\,, but not on \,$N$\,) so that
$$ \|\Delta^{[0]}-\Delta_0\|_{\As(\C^*,\ell^2_{0,0},1)} \leq C_1 \cdot R_0 $$
holds, moreover by Proposition~\ref{P:interpolate:mu}(3)(b) (note that we have \,$\upsilon^{[\nu]}=1$\, and \,$\tau^{[2]}=1$\, in the application of that proposition),
the \,$\ell^2_{0,0}$-sequence
$$ C_1 \cdot \left( \left(a_k * \frac{1}{|k|}\right)*\frac{1}{|k|} + \left( |\mu_k^{[0]}-\mu_{k,0}| + \frac{|\tau-1|}{|k|} \right) * \frac{1}{|k|} \right) $$
with
$$ a_k := \begin{cases} k^{-1}\cdot |\lambda_k-\lambda_{k,0}| & \text{for \,$k\geq 0$\,} \\ |k|^3\cdot |\lambda_k-\lambda_{k,0}| & \text{for \,$k<0$\,} \end{cases}  $$
is a bounding sequence for \,$\Delta^{[0]}-\Delta_0$\,. Because we have \,$|\mu_k^{[0]}-\mu_{k,0}| \leq |\tfrac12(\mu_k+\mu_k^{-1})-\mu_{k,0}|$\, for all \,$k$\,, in fact 
$$ b_k := C_1 \cdot \left( \left(a_k * \frac{1}{|k|}\right)*\frac{1}{|k|} + \left( \left|\tfrac12(\mu_k+\mu_k^{-1})-\mu_{k,0}\right| + \frac{|\tau-1|}{|k|} \right) * \frac{1}{|k|} \right) $$
is another bounding sequence in \,$\ell^2_{0,0}(k)$\, for \,$\Delta^{[0]}-\Delta_0$\,; this sequence does not depend on \,$N$\,.
Next we estimate \,$\|\eta_k^{[0]}-\lambda_{k,0}\|_{\ell^2_{-1,3}(|k|>N)}$\, by applying Lemma~\ref{L:finite:eta}(2) in the setting \,$\Delta^{[1]}=\Delta^{[0]}$\,, \,$\Delta^{[2]} = \Delta_0$\,: 
Because \,$\|\Delta^{[0]}-\Delta_0\|_{\As(\C^*,\ell^2_{0,0},1)}$\, is bounded independently of \,$N$\,, there exists a constant \,$C_2>0$\,
which is independent of \,$N$\,, so that we have for \,$k\in \Z$\,
\begin{align*}
|\eta_k^{[0]}-\lambda_{k,0}| 
& \leq \tfrac13\,C_2 \cdot \left\{ \begin{matrix} k & \text{if \,$k>0$\,} \\ |k|^{-3} & \text{if \,$k<0$\,} \end{matrix} \right\} \cdot \max\{b_{k-1},b_k,b_{k+1}\} \; . 
\end{align*}
Therefore we have
\begin{equation}
\label{eq:finite:finite:etak0}
\|\eta_k^{[0]}-\lambda_{k,0}\|_{\ell^2_{-1,3}(|k|>N)} \leq C_2 \cdot \|b_k\|_{\ell^2_{0,0}(|k|>N)} \; .
\end{equation}
Finally, we have
\begin{align*}
\left|\wt{z}_k^{[0]}\right| & = \left| \Delta^{[0]}(\lambda_k)-\Delta^{[0]}(\eta_k) \right| 
\leq \left( |\eta_k^{[0]}-\lambda_{k,0}| + |\lambda_k-\lambda_{k,0}| \right) \cdot \max_{\lambda\in U_{k,\delta}} |\Delta^{[0]}{}'(\lambda)| \; .
\end{align*}
Because we have \,$\Delta^{[0]}-\Delta_0\in \As(\C^*,\ell^2_{0,0},1)$\,, where \,$(b_k)$\, is a bounding sequence independent of \,$N$\,, we also have 
\,$\Delta^{[0]}{}'-\Delta_0' \in \As(\C^*,\ell^2_{1,-3},1)$\, with a bounding sequence independent of \,$N$\, by Proposition~\ref{P:interpolate:l2asymp}(3). Because we also have
\,$\Delta_0' \in \As(\C^*,\ell^\infty_{1,-3},1)$\,, it follows that we have \,$\Delta^{[0]}{}' \in \As(\C^*,\ell^\infty_{1,-3},1)$\, with a bounding sequence \,$(c_k)\in \ell^\infty_{1,-3}(k)$\,
that is independent of \,$N$\,. With this sequence, we have
\begin{equation*}
\left|\wt{z}_k^{[0]}\right|  \leq \left( |\eta_k^{[0]}-\lambda_{k,0}| + |\lambda_k-\lambda_{k,0}| \right) \cdot c_k \; . 
\end{equation*}
Hence there exist constants \,$C_3,C_4>0$\,, independent of \,$N$\,, so that 
\begin{align}
\left\|\wt{z}_k^{[0]}\right\|_{\mathfrak{B}_N}  & \;\;\,\leq\;\;\, C_3 \cdot \left( \|\eta_k^{[0]}-\lambda_{k,0}\|_{\ell^2_{-1,3}} + \|\lambda_k-\lambda_{k,0}\|_{\ell^2_{-1,3}(|k|>N)} \right) \notag\\
& \overset{\eqref{eq:finite:finite:etak0}}{\leq} C_3 \cdot \left( C_2\cdot \|b_k\|_{\ell^2_{0,0}(|k|>N)} + \|\lambda_k-\lambda_{k,0}\|_{\ell^2_{-1,3}(|k|>N)} \right) \notag\\
\label{eq:finite:finite:wtz0}
& \;\;\,\leq\;\; C_4 \cdot \left( \|b_k\|_{\ell^2_{0,0}(|k|>N)} + \|\lambda_k-\lambda_{k,0}\|_{\ell^2_{-1,3}(|k|>N)} \right) \; . 
\end{align}

We now fix besides \,$N\in \N$\, also \,$\delta>0$\,, and consider the closed ball
$$ \mathfrak{B}_{N,\delta} := \Mengegr{(z_k)\in \mathfrak{B}_N}{\|z_k\|_{\mathfrak{B}_N} \leq \delta} $$
in \,$\mathfrak{B}_N$\,. We will show that if we choose \,$\delta$\, small enough and \,$N$\, large enough, then \,$\Phi$\, is a contracting self-mapping on \,$\mathfrak{B}_{N,\delta}$\,. 

For this purpose, we let two sequences \,$(z_k^{[1]}),(z_k^{[2]})\in \mathfrak{B}_N$\, be given, and we denote the objects associated to \,$(z_k^{[\nu]})$\, by 
\,$\Delta^{[\nu]}$\,, \,$(\eta_k^{[\nu]})$\, and \,$(\wt{z}_k^{[\nu]})$\,. Moreover, we denote the objects associated to the zero sequence \,$(z_k=0)$\, by 
\,$\Delta^{[0]}$\,, \,$(\eta_k^{[0]})$\, and \,$(\wt{z}_k^{[0]})$\, as before.
By Lemma~\ref{L:finite:Delta} (for its application, note that for the \,$w_k^{[\nu]}$\, from the Lemma we have \,$w_k^{[1]}-w_k^{[2]} = \Delta^{[1]}(\lambda_k)-\Delta^{[2]}(\lambda_k) = z_k^{[1]}-z_k^{[2]}$\,)
we have for \,$|n|> N$\,
\begin{align*}
|\wt{z}_n^{[1]}-\wt{z}_n^{[2]}|
& = \left| \left( \Delta^{[1]}(\eta_n^{[1]}) - \Delta^{[1]}(\lambda_n)\right) \;-\; \left( \Delta^{[2]}(\eta_n^{[2]}) - \Delta^{[2]}(\lambda_n)\right) \right|  \\
& \leq C_5 \cdot \left( |\zeta(\eta_n^{[1]})-\zeta(\lambda_n)|\cdot \left( \frac{1}{|k|} * |z_k^{[1]}-z_k^{[2]}| \right)_n + |\zeta(\eta_n^{[1]})-\zeta(\eta_n^{[2]})|^2 \right) \; ,
\end{align*}
where \,$C_5$\, and all \,$C_k$\, (\,$k>5$\,) occurring in the sequel are positive constants, which apply uniformly for all \,$(z_k)\in \mathfrak{B}_{N,\delta}$\,, all sufficiently large \,$N\in \N$\,,
and all \,$\delta>0$\, which are smaller than some arbitrarily fixed upper bound.
From this estimate, we obtain by Cauchy-Schwarz's inequality, the variant \eqref{eq:fasymp:fourier-small:weakyoung} of Young's inequality for weakly \,$\ell^1$-sequences and Proposition~\ref{P:vac2:excldom-new}(1)

{\footnotesize
\begin{align}
\|\wt{z}_n^{[1]}-\wt{z}_n^{[2]}\|_{\mathfrak{B}_N}
& \leq C_6 \cdot \|\wt{z}_n^{[1]}-\wt{z}_n^{[2]}\|_{\ell^1_{0,0}(|n|>N)} \notag \\
& \leq C_7 \cdot \left( \|\zeta(\eta_n^{[1]})-\zeta(\lambda_n)\|_{\ell^2_{0,0}(|n|>N)} \cdot \|z_n^{[1]}-z_n^{[2]}\|_{\ell^2_{0,0}(|n|>N)} + \|\zeta(\eta_n^{[1]})-\zeta(\eta_n^{[2]})\|_{\ell^2_{0,0}(|n|>N)}^2 \right) \notag \\
\label{eq:finite:finite:wtzn12-norm-pre}
& \leq C_8 \cdot \left( \|\eta_n^{[1]}-\lambda_n\|_{\ell^2_{-1,3}(|n|>N)} \cdot \|z_n^{[1]}-z_n^{[2]}\|_{\ell^2_{0,0}(|n|>N)} + \|\eta_n^{[1]}-\eta_n^{[2]}\|_{\ell^2_{-1,3}(|n|>N)}^2 \right) \; .
\end{align}
}

We now note that we have by Lemma~\ref{L:finite:eta}(2) and Proposition~\ref{P:interpolate:mu}(2)(c) (in the application of the latter proposition, we have \,$\lambda_k^{[1]}=\lambda_k^{[2]}$\,, \,$\tau^{[1]}=\tau^{[2]}$\,
and \,$\upsilon^{[\nu]}=1$\,)
\begin{equation}
\label{eq:finite:finite:eta12-norm}
\|\eta_n^{[1]}-\eta_n^{[2]}\|_{\ell^2_{-1,3}(|n|>N)} \leq C_9 \cdot \|\Delta^{[1]}-\Delta^{[2]}\|_{\As(\C^*,\ell^2_{0,0},1)} \leq C_{10} \cdot \|z_n^{[1]}-z_n^{[2]}\|_{\mathfrak{B}_N} \;,
\end{equation}
and for \,$\nu \in \{1,2\}$\, we also have
{\footnotesize
\begin{align}
\|\eta_n^{[\nu]}-\lambda_n\|_{\ell^2_{-1,3}(|n|>N)} & \;\;\leq\;\; \|\eta_n^{[\nu]}-\eta_n^{[0]}\|_{\ell^2_{-1,3}(|n|>N)} + \|\eta_n^{[0]}-\lambda_{n,0}\|_{\ell^2_{-1,3}(|n|>N)} + \|\lambda_n-\lambda_{n,0}\|_{\ell^2_{-1,3}(|n|>N)} \notag \\
\label{eq:finite:finite:etalambda-norm}
& \overset{\eqref{eq:finite:finite:eta12-norm}}{\overset{\eqref{eq:finite:finite:etak0}}{\leq}} C_{10}\cdot \|z_n^{[\nu]}\|_{\mathfrak{B}_N} + C_2 \cdot \|b_n\|_{\ell^2_{0,0}(|n|>N)} + \|\lambda_n-\lambda_{n,0}\|_{\ell^2_{-1,3}(|n|>N)} \; . 
\end{align}
}
By applying these estimates to \eqref{eq:finite:finite:wtzn12-norm-pre}, we obtain
{\footnotesize
\begin{align}
\|\wt{z}_n^{[1]}-\wt{z}_n^{[2]}\|_{\mathfrak{B}_N}
& \leq \biggr( C_{10}\cdot \|z_n^{[1]}\|_{\mathfrak{B}_N} + C_2 \cdot \|b_k\|_{\ell^2_{0,0}(|n|>N)} \notag \\
& \qquad\qquad + \|\lambda_n-\lambda_{n,0}\|_{\ell^2_{-1,3}(|n|>N)} + C_{10}^2 \cdot \|z_n^{[1]}-z_n^{[2]}\|_{\mathfrak{B}_N} \biggr) \cdot \|z_n^{[1]}-z_n^{[2]}\|_{\mathfrak{B}_N} \notag\\
\label{eq:finite:finite:wtzn12-norm}
& \leq C_{11} \cdot \left( \|z_n^{[1]}\|_{\mathfrak{B}_N} + \|z_n^{[2]}\|_{\mathfrak{B}_N} + \|b_k\|_{\ell^2_{0,0}(|n|>N)} + \|\lambda_n-\lambda_{n,0}\|_{\ell^2_{-1,3}(|n|>N)} \right) \cdot \|z_n^{[1]}-z_n^{[2]}\|_{\mathfrak{B}_N} \;.
\end{align}
}
Inequality~\eqref{eq:finite:finite:wtzn12-norm} shows that \,$\Phi$\, is Lipschitz continuous on \,$\mathfrak{B}_{N,\delta}$\,. 

If we now choose 
$$ \delta := \min \left\{ \frac{1}{8\,C_{11}}, \frac{\eps}{C_{10}+C_2+2} \right\} $$
and then \,$N\in \N$\, so large that the following inequalities hold:
\begin{align*}
\|b_n\|_{\ell^2_{0,0}(|n|>N)} & \leq \min\left\{ \delta, \frac{\delta}{4\,C_4} \right\} \\
\|\lambda_n-\lambda_{n,0}\|_{\ell^2_{-1,3}(|n|>N)} & \leq \min\left\{ \delta, \frac{\delta}{4\,C_4} \right\} \\
\|\mu_n-\mu_{n,0}\|_{\ell^2_{0,0}(|n|>N)} & \leq \delta \;, 
\end{align*}
then it follows from Equation~\eqref{eq:finite:finite:wtzn12-norm} that for \,$(z_k^{[\nu]}) \in \mathfrak{B}_{N,\delta}$\, we have
\begin{equation}
\label{eq:finite:finite:contraction}
\|\wt{z}_n^{[1]}-\wt{z}_n^{[2]}\|_{\mathfrak{B}_N} \leq C_{11} \cdot (\delta+\delta+\delta+\delta) \cdot \|z_n^{[1]}-z_n^{[2]}\|_{\mathfrak{B}_N} \leq \frac12 \|z_n^{[1]}-z_n^{[2]}\|_{\mathfrak{B}_N} \; .
\end{equation}
Moreover, by setting \,$z_k^{[1]}=z_k \in \mathfrak{B}_{N,\delta}$\, and \,$z_k^{[2]}=0$\, in this inequality, we also see
$$ \|\wt{z}_n-\wt{z}_n^{[0]}\|_{\mathfrak{B}_N} \leq \frac12 \|z_n\|_{\mathfrak{B}_N} $$
and therefore
\begin{align}
\label{eq:finite:finite:selfmap}
\|\wt{z}_n\|_{\mathfrak{B}_N} & \;\;\,\leq\;\; \frac12 \,\|z_n\|_{\mathfrak{B}_N} + \|\wt{z}_n^{[0]}\|_{\mathfrak{B}_N} \notag \\
& \overset{\eqref{eq:finite:finite:wtz0}}{\leq} \frac12 \,\|z_n\|_{\mathfrak{B}_N} + C_4\cdot \left( \|b_n\|_{\ell^2_{0,0}(|n|>N)} + \|\lambda_n-\lambda_{n,0}\|_{\ell^2_{-1,3}(|n|>N)} \right) \notag \\
& \;\;\,\leq\;\; \frac12\,\delta + C_4\cdot \left(\frac{\delta}{4\,C_4}+\frac{\delta}{4\,C_4}\right) = \delta \; . 
\end{align}
Inequality~\eqref{eq:finite:finite:selfmap} shows that \,$\Phi$\, maps the complete metric space \,$\mathfrak{B}_{N,\delta}$\ into itself, and inequality~\eqref{eq:finite:finite:contraction} shows that 
\,$\Phi$\, is a contraction (with Lipschitz constant \,$\tfrac12$\,) on this space.

The Banach Fixed Point Theorem therefore shows
that \,$\Phi$\, has exactly one fixed point \,$(z_n^*)$\, in \,$\mathfrak{B}_{N,\delta}$\,. We let \,$\Delta^*$\, and \,$(\eta_k^*)$\, be the objects associated to \,$(z_n^*)$\,. Then we define the
divisor \,$D^*=\{(\lambda_k^*,\mu_k^*)\}$\, by
$$ \lambda_k^* := \begin{cases} \lambda_k & \text{for \,$|k|\leq N$\,} \\ \eta_k^* & \text{for \,$|k|>N$\,} \end{cases}
\qmq{and}
\mu_k^* := \begin{cases} \mu_k & \text{for \,$|k|\leq N$\,} \\ (-1)^k & \text{for \,$|k|>N$\,} \end{cases} \; . $$
By construction, \,$D^*$\, is an asymptotic divisor on the spectral curve \,$\Sigma^*$\, corresponding to \,$\Delta^*$\,. Because \,$(z_k^*)$\, is a fixed point of \,$\Phi$\,, we have for \,$|k|>N$\,
$$ z_k^* = z_k^* - (\Delta^*(\eta_k^*)-2(-1)^k) $$
and therefore \,$\Delta^*(\eta_k^*)=2(-1)^k$\,. Because we have \,$\Delta^*{}'(\eta_k^*)=0$\, by the definition of \,$\eta_k^*$\,, it follows that the spectral curve 
\,$\Sigma^*$\, associated to \,$\Delta^*$\, has a double point at \,$(\eta_k^*,2(-1)^k)=(\lambda_k^*,\mu_k^*)$\,.
Hence \,$D^*$\, is of finite type.

Moreover we have
\begin{align*}
\|D^*-D\|_{\Div}
& \;\;\,\leq\;\; \|\lambda_n^*-\lambda_n\|_{\ell^2_{-1,3}(n)} + \|\mu_n^*-\mu_n\|_{\ell^2_{0,0}(n)}  \\
& \;\;\,=\;\; \|\eta_n^*-\lambda_n\|_{\ell^2_{-1,3}(|n|>N)} + \|(-1)^n-\mu_n\|_{\ell^2_{0,0}(|n|>N)} \\
& \overset{\eqref{eq:finite:finite:etalambda-norm}}{\leq} C_{10}\cdot \|z_n^{*}\|_{\mathfrak{B}_N} + C_2 \cdot \|b_n\|_{\ell^2_{0,0}(|n|>N)} \\
& \qquad\qquad + \|\lambda_n-\lambda_{n,0}\|_{\ell^2_{-1,3}(|n|>N)} + \|(-1)^n-\mu_n\|_{\ell^2_{0,0}(|n|>N)} \\
& \;\;\,\leq\;\; C_{10}\cdot\delta + C_2\cdot \delta + \delta + \delta = (C_{10}+C_2+2)\cdot \delta \leq \eps \; .
\end{align*}

Finally, if we choose \,$N$\, large enough that \,$\lambda_k\in U_{k,\delta}$\, holds for \,$|k|>N$\,, 
and suppose that \,$\eps$\, is small enough that the condition \,$\|D^*-D\|_{\Div}<\eps$\, implies that then also \,$\lambda_k^*\in U_{k,\delta}$\,
holds for \,$|k|>N$\,, then the \,$\lambda_k^*$\, are necessarily pairwise unequal. Therefore \,$D^*$\, is then tame.
\end{proof}

\section{Darboux coordinates for the space of potentials}
\label{Se:darboux}

In the present section we will view the space of potentials \,$\Pot$\, as an (infinite-dimensional) symplectic manifold, and construct coordinates for this manifold which are adapted to the symplectic structure.
By analogy with the finite-dimensional situation, we will call these coordinates \emph{Darboux coordinates}. 

\label{not:darboux:deltaf}
More specifically, we consider the Hilbert space \,$\Pot$\, (a hyperplane in \,$\Pot_{np}=W^{1,2}([0,1])\times L^2([0,1])$\,). For \,$(u,u_y)\in\Pot$\,, the tangent space \,$T_{(u,u_y)}\Pot$\, is canonically 
isomorphic to \,$\Pot$\,, and we will typically denote elements of \,$T_{(u,u_y)}\Pot$\, by \,$(\delta u,\delta u_y)$\, and \,$(\wt{\delta} u,\wt{\delta} u_y)$\,. 
For any function \,$f$\, defined on \,$\Pot$\,, we let \,$\delta f := \tfrac{\partial f}{\partial (u,u_y)} \cdot (\delta u,\delta u_y)$\, be the variation of \,$f$\, in the direction 
of \,$(\delta u,\delta u_y) \in T_{(u,u_y)}\Pot$\,.

In this setting we define a bilinear form \,$\Omega$\, on each tangent space \,$T_{(u,u_y)}\Pot$\, by
\begin{equation}
\label{eq:darboux:Omega-defined}
\Omega: T_{(u,u_y)}\Pot \times T_{(u,u_y)}\Pot \to \C,\; \bigr( \, (\delta u,\delta u_y)\,,\, (\wt{\delta} u,\wt{\delta} u_y) \, \bigr) \mapsto \int_0^1 \bigr( \delta u \cdot \wt{\delta} u_y - \wt{\delta}u \cdot \delta u_y \bigr) \,\mathrm{d}x \; .
\end{equation}
\,$\Omega$\, defines a non-degenerate symplectic form on \,$\Pot$\,. 

\newpage

We recall a theorem for the \emph{finite-dimensional} situation due to \textsc{Darboux}:

\begin{thm}[Darboux \cite{Darboux:1882}]
Any symplectic manifold \,$(M,\Omega)$\, of dimension \,$2n<\infty$\, is \emph{locally} symplectomorphic to an open subset of \,$(\R^{2n},\Omega_0)$\,, where \,$\Omega_0$\, is the canonical symplectic form on \,$\R^{2n}$\,. 
\end{thm}

In the sequel, we will obtain coordinates on the symplectic space \,$\Pot$\, near any tame potential \,$(u,u_y)$\,, which are analogous to the coordinates promised by Darboux's
theorem for the finite-dimensional case. These coordinates are an excellent and important instrument for understanding the structure of \,$\Pot$\,; we will base our proof that the map
\,$\Pot_{tame}\to\Div_{tame}$\, is a diffeomorphism (in Section~\ref{Se:diffeo}) on the use of these coordinates. Concerning the analogous coordinates for the space of potentials for the
1-dimensional Schr\"odinger equation (see \cite{Poeschel-Trubowitz:1987},
Theorem~2.8, p.~44), \textsc{P\"oschel/Trubowitz} write: ``It is only a slight exaggeration to say that [these coordinates are] the basis of almost everything else we are going to do.''. In our situation,
we will not use the Darboux coordinates quite as intensely, but they will still prove to be very useful. The application of the Darboux coordinates for the 1-dimensional Schr\"odinger equation
to the spectral theory for the KdV equation (for which the 1-dimensional Schr\"odinger operator is the Lax operator) is described in \cite{Kappeler/Poeschel:2003}, Section~8.

In our situation concerning the sinh-Gordon equation, \textsc{M.~Knopf} has constructed in \cite{Knopf:2015} a symplectic basis (with respect to \,$\Omega$\,) for a certain subspace \,$U$\, of 
\,$T_{(u,u_y)}\Pot$\,, and thereby effectively Darboux coordinates for this subspace. Knopf did not show that this ``subspace'' is actually all of \,$T_{(u,u_y)}\Pot$\,, but we will combine
his representation with
the asymptotic estimates of the present work (especially the asymptotic description of the extended frame given in Proposition~\ref{P:asympfinal:frame}, giving rise to the asymptotic descriptions
of \,$\delta M(\lambda)$\, in Lemma~\ref{L:diffeo:delM-asymp} and of \,$\delta \lambda_k$\,, \,$\delta \mu_k$\, in Lemma~\ref{L:diffeo:delD-asymp} below) 
to prove that in fact \,$U=T_{(u,u_y)}\Pot$\, holds. Therefore Knopf's paper \cite{Knopf:2015} actually provides Darboux coordinates for \,$\Pot$\, (near tame potentials). 

Because we will base our construction of Darboux coordinates on \cite{Knopf:2015}, we now report on the results of that paper, giving his results in the notations of the present work.
When adapting the results from \cite{Knopf:2015}, one should note that the notations of that paper differ from ours in several 
important points: (1) The norming of the potential \,$u$\, differ, as is evidenced by Knopf's sinh-Gordon equation \,$\Delta u + 2\,\sinh(2u)=0$\, (\cite{Knopf:2015}, Equation~(1.1)) in comparison to 
our sinh-Gordon equation \,$\Delta u + \sinh(u)=0$\, (Equation~\eqref{eq:minimal:sinh-gordon}). (2) Knopf's flat connection form \,$\alpha_\lambda$\, differs from ours by a factor \,$2$\,, and by an additional 
conjugation with a constant matrix (compare Equation~\eqref{eq:mono:alphaxy} to \cite{Knopf:2015}, Equations~(2.1)).
(3) Knopf normalizes the eigenvector field of the monodromy as \,$(1,\tfrac{\mu-a}{b})$\,, whereas we use
\,$(\tfrac{\mu-d}{c},1)$\,. In the definition of the (classical) spectral divisor, the pair of holomorphic functions \,$(c,a)$\, is therefore replaced by \,$(b,d)$\,. Moreover, quantities have been renamed in the following report to avoid clashes with the 
conventions of the present work.

We fix \,$(u,u_y)\in \Pot_{tame}$\, and let \,$D=\{(\lambda_k,\mu_k)\}$\, be the spectral divisor of \,$(u,u_y)$\,. By hypothesis, \,$D$\, is tame, i.e.~\,$D$\, does not contain any double points.
Therefore \,$\lambda_k$\, and \,$\mu_k$\, can be interpreted as well-defined smooth functions on \,$\Pot$\, near \,$(u,u_y)$\,, and \,$c'(\lambda_k)\neq 0$\, holds for all \,$k\in \Z$\,. 

In the sequel we consider the extended frame \,$F(x,\lambda)$\, corresponding to the potential \,$(u,u_y)$\,, we write it in the form
\begin{equation}
\label{eq:darboux:extended-frame}
F(x,\lambda) = \left( \begin{matrix} a(x,\lambda) & b(x,\lambda) \\ c(x,\lambda) & d(x,\lambda) \end{matrix} \right)
\end{equation}
like we already did in Proposition~\ref{P:asympfinal:frame}, and also consider the extended frame of the vacuum
$$ F_0(x,\lambda) := \begin{pmatrix} a_0(x,\lambda) & b_0(x,\lambda) \\ c_0(x,\lambda) & d_0(x,\lambda) \end{pmatrix} := \begin{pmatrix}
\cos(x\,\zeta(\lambda)) & -\lambda^{-1/2}\,\sin(x\,\zeta(\lambda)) \\ \lambda^{1/2}\,\sin(x\,\zeta(\lambda)) & \cos(x\,\zeta(\lambda)) \end{pmatrix} \; . $$


We now define for \,$k\in \Z$\,
\begin{align}
\label{eq:darboux:vkwk}
v_k := (v_{k,1},v_{k,2}) \qmq{and} w_k := (w_{k,1},w_{k,2})
\end{align}
with
\begin{align*}
v_{k,1}(x) & := a(x,\lambda_k)\, c(x,\lambda_k) \\
v_{k,2}(x) & := \tfrac{i}{4}\left( (e^{u(x)/2}-\lambda_k\,e^{-u(x)/2})\,a(x,\lambda_k)^2 + (e^{u(x)/2}-\lambda_k^{-1}\,e^{-u(x)/2})\,c(x,\lambda_k)^2  \right) \\
w_{k,1}(x) & := a(x,\lambda_k)\, d(x,\lambda_k)+b(x,\lambda_k)\, c(x,\lambda_k) \\
w_{k,2}(x) & := \tfrac{i}{2}\big( (e^{u(x)/2}-\lambda_k\,e^{-u(x)/2})\,a(x,\lambda_k)\, b(x,\lambda_k) \\
& \qquad\qquad + (e^{u(x)/2}-\lambda_k^{-1}\,e^{-u(x)/2})\,c(x,\lambda_k)\, d(x,\lambda_k) \big) 
\end{align*}
and we also put
\begin{equation}
\label{eq:darboux:vartheta}
\vartheta_k := \int_0^1 \left( \lambda_k\,a(x,\lambda_k)^2 + \lambda_k^{-1}\,c(x,\lambda_k)^2 \right)\cdot e^{u(x)/2}\,\mathrm{d}x \;. 
\end{equation}
In \cite{Knopf:2015}, it is shown in the proof of Theorem~5.5 that \,$c'(\lambda_k) = -i\,\tfrac{\mu_k}{2\,\lambda_k}\cdot \vartheta_k$\, holds; because we have \,$c'(\lambda_k)\neq 0$\, it follows from that
equation that
$$ \vartheta_k \neq 0 $$
holds.

We can now state the main result of \cite{Knopf:2015}:
\begin{thm}[Knopf \cite{Knopf:2015}]
\label{T:darboux:knopf}
Let \,$(u,u_y)\in\Pot_{tame}$\,. 
\begin{enumerate}
\item
For \,$k\in \Z$\, we have \,$v_k,w_k \in T_{(u,u_y)}\Pot$\,, and up to the factor \,$\vartheta_k$\,, \,$(v_k,w_k)$\, is a 
system of symplectic vectors with respect to \,$\Omega$\,, i.e.~we have for all \,$k,\ell \in \Z$\,
\begin{align*}
\Omega(v_k,v_\ell) & = 0 \\
\Omega(v_k,w_\ell) & = \vartheta_k \cdot \delta_{k\ell} \quad \text{(\,$\delta_{k\ell}$\,: Kronecker delta)}\\
\Omega(w_k,w_\ell) & = 0 \; . 
\end{align*}
\item
For all \,$(\delta u,\delta u_y)\,,\,(\wt{\delta} u,\wt{\delta} u_y) \in T_{(u,u_y)}\Pot$\, that are finite linear combinations of the \,$v_k$\, and \,$w_k$\,, we have
\begin{equation}
\label{eq:darboux:Omega}
\Omega\bigr( \, (\delta u,\delta u_y)\,,\,(\wt{\delta} u,\wt{\delta} u_y) \,\bigr) = \frac{i}{2}\,\sum_{k\in \Z} \left( \frac{\delta \lambda_k}{\lambda_k} \cdot \frac{\wt{\delta} \mu_k}{\mu_k} - \frac{\wt{\delta} \lambda_k}{\lambda_k} \cdot \frac{\delta \mu_k}{\mu_k} \right) \; .
\end{equation}
\end{enumerate}
\end{thm}

\newpage

In the sequel, we liberate Knopf's preceding result from the stated restrictions. More specifically, we prove the following two statements:
\begin{enumerate}
\item The symplectic system \,$(v_k,w_k)$\, defined by Equations~\eqref{eq:darboux:vkwk} 
spans \,$T_{(u,u_y)}\Pot$\, (in the Hilbert space sense, i.e.~\,$\overline{\mathrm{span}\Menge{v_k,w_k}{k\in \Z}}=T_{(u,u_y)}\Pot$\,).
\item The infinite sum on the right-hand side of Equation~\eqref{eq:darboux:Omega} converges absolutely for arbitrary \,$(\delta u,\delta u_y)\,,\,(\wt{\delta} u,\wt{\delta} u_y) \in T_{(u,u_y)}\Pot$\,.
\end{enumerate}

In this way we will generalize Knopf's Theorem~\ref{T:darboux:knopf} to the result of the following theorem. This result shows that the \,$(v_k,w_k)$\, essentially define Darboux coordinates on \,$\Pot$\, (near any tame potential),
and the functions \,$\lambda_k^{-1}\,\delta\lambda_k$\, and \,$\mu_k^{-1}\,\delta\mu_k$\, define Darboux coordinates on \,$\Div$\, (near any tame divisor, where \,$\Div$\, is equipped with the symplectic structure induced by the map
\,$\Pot\to\Div$\, associating to each potential its associated spectral divisor).

\begin{thm}
\label{T:darboux:darboux}
Let \,$(u,u_y)\in\Pot_{tame}$\,. 
\begin{enumerate}
\item
The symplectic system \,$(v_k,w_k)$\, defined by Equations~\eqref{eq:darboux:vkwk} is a symplectic basis of \,$T_{(u,u_y)} \Pot$\, (as before, up to the factor \,$\vartheta_k$\,). 
\item
For all \,$(\delta u,\delta u_y)\,,\,(\wt{\delta} u,\wt{\delta} u_y) \in T_{(u,u_y)}\Pot$\,,
\begin{equation}
\label{eq:darboux:darboux:Omega}
\Omega\bigr( \, (\delta u,\delta u_y)\,,\,(\wt{\delta} u,\wt{\delta} u_y) \,\bigr) = \frac{i}{2}\,\sum_{k\in \Z} \left( \frac{\delta \lambda_k}{\lambda_k} \cdot \frac{\wt{\delta} \mu_k}{\mu_k} - \frac{\wt{\delta} \lambda_k}{\lambda_k} \cdot \frac{\delta \mu_k}{\mu_k} \right) 
\end{equation}
holds, and the sum on the right hand side of this equation converges.
\end{enumerate}
\end{thm}

\begin{rem}
In Theorem~\ref{T:darboux:darboux} we are still restricted to tame potentials resp.~divisors, i.e.~spectral divisors without multiple points. To lift this restriction, we would need to deal with the complication
that if \,$(\lambda_k,\mu_k)=(\lambda_{k'},\mu_{k'})$\, holds for \,$k\neq k'$\,, then also \,$v_k=v_{k'}$\, and \,$w_k=w_{k'}$\, holds, and therefore we do not get ``enough'' independent coordinates by
the \,$(v_k,w_k)$\, in this case. To obtain the ``missing'' coordinates, we would need to consider appropriate coordinates on \,$\Div$\, in such points; they are given by the elementary symmetric polynomials
in the \,$\lambda_k$\, resp.~\,$\mu_k$\,, see the discussion at the end of Section~\ref{Se:asympdiv}.
\end{rem}


The remainder of the section is dedicated to the proof of Theorem~\ref{T:darboux:darboux}. 

\begin{lem}
\label{L:darboux:basis}
Let \,$(u,u_y)\in \Pot_{tame}$\, be given. 
The symplectic system \,$(v_k,w_k)$\, defined by Equations~\eqref{eq:darboux:vkwk} is a basis of \,$T_{(u,u_y)} \Pot$\,. 
\end{lem}

\begin{proof}
By Theorem~\ref{T:darboux:knopf}(1) (=\cite{Knopf:2015}, Theorem~5.4), 
\,$(v_k,w_k)$\, is a  symplectic system (up to the factor \,$\vartheta_k$\,)  with respect to the non-degenerate symplectic form \,$\Omega$\,, and therefore these vectors are linear independent.
It remains to show that they span \,$T_{(u,u_y)} \Pot\cong \Mengegr{(u,u_y)\in W^{1,2}([0,1]) \times L^2([0,1])}{u(0)=u(1)}$\, in the Hilbert space sense, 
i.e.~that the set \,$\mathrm{span}\Menge{v_k,w_k}{k\in\Z}$\, is dense in \,$T_{(u,u_y)}\Pot$\,. 
Because \,$\Mengegr{u\in W^{1,2}([0,1])}{u(0)=u(1)}$\, is dense in \,$L^2([0,1])$\,, the preceding claim is equivalent to the statement that the set \,$\mathrm{span}\Menge{v_k,w_k}{k\in\Z}$\, is dense in
\,$L^2([0,1]) \times L^2([0,1])$\,. We will prove the latter statement
by applying the following fact from functional analysis, which we cite from \cite{Poeschel-Trubowitz:1987}, Theorem~D.3, p.~163f.~(where its proof can
also be found):

\begin{quote}
Let \,$(h_i)_{i\in I}$\, be an orthonormal basis of a Hilbert space \,$H$\,, and suppose that \,$(\wt{h}_i)_{i\in I}$\, is another sequence of vectors in \,$H$\, that is linear independent. If in addition
$$ \sum_{i\in I} \|h_i-\wt{h}_i\|^2 < \infty \;, $$
then \,$(\wt{h}_i)_{i\in I}$\, is also a basis of \,$H$\,. 
\end{quote}

To apply this theorem to the present situation (with \,$H:=L^2([0,1]) \times L^2([0,1])$\,), 
we need to investigate the asymptotic behavior of the \,$v_k$\, and the \,$w_k$\,. In relation to this, we denote for \,$1\leq p\leq \infty$\, and \,$n\in \Z$\, 
by \,$\ell^p_n(L^2([0,1]))$\, the space of sequences \,$(f_k)_{k\geq 1}$\, in \,$L^2([0,1])$\, with 
\,$(\|f_k\|_{L^2([0,1])}) \in \ell^p_n(k\geq 1)$\,.


As an exemplar case we look at \,$v_{k,1}(x)=a(x,\lambda_k)\cdot c(x,\lambda_k)$\, for \,$k>0$\,. By Proposition~\ref{P:asympfinal:frame} we have
\begin{gather*}
a-\upsilon\,a_0 \in \As_\infty(\C^*,\ell^2_0,1)\;,\quad c-\tau\,c_0 \in \As_\infty(\C^*,\ell^2_{-1},1) \;, \\
c \in \As_\infty(\C^*,\ell^\infty_{-1},1) \qmq{and} \upsilon\,a_0 \in \As_\infty(\C^*,\ell^\infty_0,1) \;,
\end{gather*}
where all these asymptotic assessments are uniform in \,$x\in [0,1]$\, and where we put \,$\tau(x) := e^{-(u(0)+u(x))/4}$\,, \,$\upsilon(x) := e^{(u(x)-u(0))/4}$\,. As a consequence we have
\begin{align*}
a(x,\lambda_k)-\upsilon(x)\,a_0(x,\lambda_k) & \in \ell^2_0(L^2([0,1])) \\
c(x,\lambda_k)-\tau(x)\,c_0(x,\lambda_k) & \in \ell^2_{-1}(L^2([0,1])) \\
c(x,\lambda_k) & \in \ell^\infty_{-1}(L^2([0,1])) \\
\upsilon(x)\,a_0(x,\lambda_k) & \in \ell^\infty_0(L^2([0,1])) 
\end{align*}
and therefore
\begin{align*}
& v_{k,1}(x) - \frac{1}{2}\,\sqrt{\lambda_k}\,e^{-u(0)/2}\,\sin(2\zeta(\lambda_k)x) \\
=\; & a(x,\lambda_k)\cdot c(x,\lambda_k)-\upsilon(x)\,a_0(x,\lambda_k)\,\tau(x)\,c_0(x,\lambda_k) \in \ell^2_{-1}(L^2([0,1])) \; . 
\end{align*}
Moreover, it is easy to check that
$$ \sqrt{\lambda_k}\,\sin(2\zeta(\lambda_k)x) - 4\pi k\,\sin(2k\pi\,x) \in \ell^2_{-1}(L^2([0,1])) $$
holds, and thus we obtain
$$ v_{k,1}(x) - 2\pi k \,e^{-u(0)/2}\,\sin(2k\pi\,x) \in \ell^2_{-1}(L^2([0,1])) \; . $$
By carrying out analogous calculations for \,$v_{k,2}$\, and \,$w_{k,\nu}$\,, one obtains for \,$k>0$\,
\begin{align*}
v_{k,1}(x) & - k\cdot 2\pi \,e^{-u(0)/2}\cdot\sin(2k\pi\,x) & & \in \ell^2_{-1}(L^2([0,1])) \\
v_{k,2}(x) & - k^2 \cdot (-4)\pi^2\,i\,e^{-u(0)/2}\cdot\cos(2k\pi\, x) & &\in \ell^2_{-2}(L^2([0,1])) \\
w_{k,1}(x) & - \cos(2k\pi\,x) & & \in \ell^2_0(L^2([0,1])) \\
w_{k,2}(x) & -  k\cdot 2\pi i \cdot \sin(2k\pi\,x) & & \in \ell^2_{-1}(L^2([0,1])) 
\end{align*}
and also 
\begin{align*}
v_{-k,1}(x) & - k^{-1}\cdot \tfrac{1}{8\pi}\, e^{u(0)/2}\cdot \sin(2k\pi\,x) & & \in \ell^2_{1}(L^2([0,1])) \\
v_{-k,2}(x) & - \tfrac{i}{4}\,e^{u(0)/2}\cdot \cos(2k\pi\, x) & &\in \ell^2_{0}(L^2([0,1])) \\
w_{-k,1}(x) & - \cos(2k\pi\,x) & & \in \ell^2_0(L^2([0,1])) \\
w_{-k,2}(x) & - k\cdot (-2)\pi i \cdot \sin(2k\pi\,x) & & \in \ell^2_{-1}(L^2([0,1])) 
\end{align*}
We now define elements of \,$L^2([0,1]) \times L^2([0,1])$\, for \,$k\geq 1$\, by:
%
%
\begin{align}
\wt{h}_{s1,k} & := \sqrt{2}\cdot \left( \frac{1}{4\pi}\,e^{u(0)/2}\cdot k^{-1}\cdot v_k + 4\pi\,e^{-u(0)/2}\cdot k \cdot v_{-k} \right) \notag \\
\wt{h}_{s2,k} & := \frac{\sqrt{2}}{4\pi i} \cdot k^{-1} \cdot \left( w_k-w_{-k} \right) \notag \\
\wt{h}_{c1,k} & := \frac{\sqrt{2}}{2} \cdot \left(w_k + w_{-k} \right) \notag \\
\label{eq:darboux:basis:wth}
\wt{h}_{c2,k} & := \sqrt{2} \cdot \left( \frac{1}{8\pi^2}\,i\,e^{u(0)/2}\cdot k^{-2}\cdot v_k - 2i\,e^{-u(0)/2}\cdot v_{-k} \right) 
%
\end{align}
Moreover, we put
$$ \wt{h}_{c1,0} := v_0 \qmq{and} \wt{h}_{c2,0} := w_0 \; . $$
Because the vectors in \,$\Menge{v_k,w_k}{k\in \Z}$\, are linear independent, the above definitions show that the vectors in \,$\Menge{\wt{h}_{c1,k},\wt{h}_{c2,k}}{k\geq 0} \cup \Menge{\wt{h}_{s1,k},\wt{h}_{s2,k}}{k\geq 1}$\,
are also linear independent. Moreover, from the preceding asymptotic assessments of the \,$v_{\pm k,\nu}$\, and the \,$w_{\pm k,\nu}$\, one can calculate that 
\begin{align}
\wt{h}_{c1,k}-h_{c1,k}\,,\, \wt{h}_{c2,k}-h_{c2,k} & \;\in\; \ell^2_0(L^2([0,1])) \qmq{for \,$k\geq 0$\,} \notag \\
\label{eq:darboux:basis:asymp}
\qmq{and} \wt{h}_{s1,k}-h_{s1,k}\,,\, \wt{h}_{s2,k}-h_{s2,k} & \;\in\; \ell^2_0(L^2([0,1])) \qmq{for \,$k\geq 1$\,} 
\end{align}
holds, where we put for \,$k\geq 0$\,
$$ h_{c1,k}(x) := \bigr(\sqrt{2}\,\cos(2k\pi\, x)\,,\,0\bigr) \qmq{and} h_{c2,k}(x) := \bigr(0\,,\,\sqrt{2}\,\cos(2k\pi\, x)\bigr) $$
and for \,$k\geq 1$\,
$$ h_{s1,k}(x) := \bigr(\sqrt{2}\,\sin(2k\pi\, x)\,,\,0\bigr) \qmq{and} h_{s2,k}(x) := \bigr(0\,,\,\sqrt{2}\,\sin(2k\pi\, x)\bigr) \; . $$
\,$\Menge{h_{c1,k},h_{c2,k}}{k\geq 0} \cup \Menge{h_{s1,k},h_{s2,k}}{k\geq 1}$\, is an orthonormal basis of the Hilbert space \,$H := L^2([0,1]) \times L^2([0,1])$\,, the vectors in
\,$\Menge{\wt{h}_{c1,k},\wt{h}_{c2,k}}{k\geq 0} \cup \Menge{\wt{h}_{s1,k},\wt{h}_{s2,k}}{k\geq 1}$\, are linear independent, and 
by \eqref{eq:darboux:basis:asymp} we have 
$$ \sum_{k=0}^\infty \|\wt{h}_{c1,k}-h_{c1,k}\|_H^2 + \sum_{k=0}^\infty \|\wt{h}_{c2,k}-h_{c2,k}\|_H^2 + \sum_{k=1}^\infty \|\wt{h}_{s1,k}-h_{s1,k}\|_H^2 + \sum_{k=1}^\infty \|\wt{h}_{s2,k}-h_{s2,k}\|_H^2
\;<\; \infty \; . $$
It follows by the theorem cited at the beginning of the proof that 
\,$\Menge{\wt{h}_{c1,k},\wt{h}_{c2,k}}{k\geq 0} \cup \Menge{\wt{h}_{s1,k},\wt{h}_{s2,k}}{k\geq 1}$\, is a basis of \,$H$\,. By solving the Equations~\eqref{eq:darboux:basis:wth} defining the
\,$\wt{h}_{\dotsc}$\, for \,$v_{\pm k}$\, and \,$w_{\pm k}$\, it follows that \,$\Menge{v_k,w_k}{k\in \Z}$\, also is a basis of \,$H$\,. This completes the proof.
\end{proof}

After the question of Theorem~\ref{T:darboux:darboux}(1) has been settled with the preceding Lemma, we next need to show the convergence of the infinite sum occurring in Theorem~\ref{T:darboux:darboux}(2).
To do so, we will have to describe the asymptotic behavior of the summand \,$\tfrac{\delta \lambda_k}{\lambda_k} \cdot \tfrac{\wt{\delta} \mu_k}{\mu_k} - \tfrac{\wt{\delta} \lambda_k}{\lambda_k} \cdot \tfrac{\delta \mu_k}{\mu_k}$\,
for \,$k\to\pm\infty$\,, and thus we will need to find out about the asymptotic behaviour of the functions \,$\delta \lambda_k$\, and \,$\delta \mu_k$\,. 

As preparation for the latter task, the following Lemma describes the asymptotic behaviour of the variation \,$\delta M(\lambda)$\, of the monodromy \,$M(\lambda)$\,. 

\begin{lem}
\label{L:diffeo:delM-asymp}
Let \,$(u,u_y)\in \Pot_{tame}$\, be given. We denote the monodromy of \,$(u,u_y)$\, by \,$M(\lambda) = \left( \begin{smallmatrix} a(\lambda) & b(\lambda) \\ c(\lambda) & d(\lambda) \end{smallmatrix} \right)$\, 
as usual.
\begin{enumerate}
\item \emph{(The case \,$\lambda\to\infty$\,.)} There exist functions depending on \,$\lambda\in \C^*$\, and \,$t\in [0,1]$\, with
\begin{align*}
f_{1,a} & \in \As_\infty(\C^*,\ell^2_0,2)\,, & f_{2,a} & \in \As_\infty(\C^*,\ell^\infty_1,2)\,, \\
f_{1,b} & \in \As_\infty(\C^*,\ell^2_1,2)\,, & f_{2,b} & \in \As_\infty(\C^*,\ell^\infty_2,2)\,, \\
f_{1,c} & \in \As_\infty(\C^*,\ell^2_{-1},2)\,, & f_{2,c} & \in \As_\infty(\C^*,\ell^\infty_0,2)\,, \\
f_{1,d} & \in \As_\infty(\C^*,\ell^2_0,2)\,, & f_{2,d} & \in \As_\infty(\C^*,\ell^\infty_1,2)\,,
\end{align*}
where all these memberships in asymptotic spaces are uniform in \,$t\in [0,1]$\,, so that we have
{\allowdisplaybreaks
\begin{align*}
\delta a(\lambda) & = -\tfrac12\,a(\lambda)\,\int_0^1 \delta u_z(t)\,\cos(2\zeta(\lambda)t)\,\mathrm{d}t + \tfrac12\,\tau\,\lambda^{1/2}\,b(\lambda)\,\int_0^1 \delta u_z(t)\,\sin(2\zeta(\lambda)t)\,\mathrm{d}t \\
& \qquad\qquad + \int_0^1 \delta u_z(t)\,f_{1,a}(t)\,\mathrm{d}t + \int_0^1 \delta u(t)\,f_{2,a}(t)\,\mathrm{d}t \\
\delta b(\lambda) & = \tfrac12\,b(\lambda)\,\left( \int_0^1 \delta u_z(t)\,\cos(2\zeta(\lambda)t)\,\mathrm{d}t - \delta u(0) \right) + \tfrac12\,\tau^{-1}\,\lambda^{-1/2}\,a(\lambda)\,\int_0^1 \delta u_z(t)\,\sin(2\zeta(\lambda)t)\,\mathrm{d}t \\
& \qquad\qquad + \int_0^1 \delta u_z(t)\,f_{1,b}(t)\,\mathrm{d}t + \int_0^1 \delta u(t)\,f_{2,b}(t)\,\mathrm{d}t \\
\delta c(\lambda) & = -\tfrac12\,c(\lambda)\left( \int_0^1 \delta u_z(t)\,\cos(2\zeta(\lambda)t)\,\mathrm{d}t-\delta u(0) \right) + \tfrac12\,\tau\,\lambda^{1/2}\,d(\lambda)\,\int_0^1 \delta u_z(t)\,\sin(2\zeta(\lambda)t)\,\mathrm{d}t \\
& \qquad\qquad + \int_0^1 \delta u_z(t)\,f_{1,c}(t)\,\mathrm{d}t + \int_0^1 \delta u(t)\,f_{2,c}(t)\,\mathrm{d}t \\
\delta d(\lambda) & = \tfrac12\,d(\lambda)\,\int_0^1 \delta u_z(t)\,\cos(2\zeta(\lambda)t)\,\mathrm{d}t + \tfrac12\,\tau^{-1}\,\lambda^{-1/2}\,c(\lambda)\,\int_0^1 \delta u_z(t)\,\sin(2\zeta(\lambda)t)\,\mathrm{d}t \\
& \qquad\qquad + \int_0^1 \delta u_z(t)\,f_{1,d}(t)\,\mathrm{d}t + \int_0^1 \delta u(t)\,f_{2,d}(t)\,\mathrm{d}t \; . 
\end{align*}
}
\item 
\emph{(The case \,$\lambda\to0$\,.)}
There exist functions depending on \,$\lambda\in \C^*$\, and \,$t\in [0,1]$\, with
\begin{align*}
g_{1,a} & \in \As_0(\C^*,\ell^2_0,2)\,, & g_{2,a} & \in \As_0(\C^*,\ell^\infty_1,2)\,, \\
g_{1,b} & \in \As_0(\C^*,\ell^2_{-1},2)\,, & g_{2,b} & \in \As_0(\C^*,\ell^\infty_0,2)\,, \\
g_{1,c} & \in \As_0(\C^*,\ell^2_{1},2)\, & g_{2,c} & \in \As_0(\C^*,\ell^\infty_2,2)\, \\
g_{1,d} & \in \As_0(\C^*,\ell^2_0,2)\,, & g_{2,d} & \in \As_0(\C^*,\ell^\infty_1,2)\,, 
\end{align*}
where all these memberships in asymptotic spaces are uniform in \,$t\in [0,1]$\,, so that we have
{\allowdisplaybreaks
\begin{align*}
\delta a(\lambda) & = \tfrac12\,a(\lambda)\,\int_0^1 \delta u_{\overline{z}}(t)\,\cos(2\zeta(\lambda)t)\,\mathrm{d}t - \tfrac12\,\tau^{-1}\,\lambda^{1/2}\,b(\lambda)\,\int_0^1 \delta u_{\overline{z}}(t)\,\sin(2\zeta(\lambda)t)\,\mathrm{d}t \\
& \qquad\qquad + \int_0^1 \delta u_{\overline{z}}(t)\,g_{1,a}(t)\,\mathrm{d}t + \int_0^1 \delta u(t)\,g_{2,a}(t)\,\mathrm{d}t \\
\delta b(\lambda) & = \tfrac12\,b(\lambda)\,\left( -\int_0^1 \delta u_{\overline{z}}(t)\,\cos(2\zeta(\lambda)t)\,\mathrm{d}t + \delta u(0) \right) - \tfrac12\,\tau\,\lambda^{-1/2}\,a(\lambda)\,\int_0^1 \delta u_{\overline{z}}(t)\,\sin(2\zeta(\lambda)t)\,\mathrm{d}t \\
& \qquad\qquad + \int_0^1 \delta u_{\overline{z}}(t)\,g_{1,b}(t)\,\mathrm{d}t + \int_0^1 \delta u(t)\,g_{2,b}\,\mathrm{d}t \\
\delta c(\lambda) & = -\tfrac12\,c(\lambda)\left( -\int_0^1 \delta u_{\overline{z}}(t)\,\cos(2\zeta(\lambda)t)\,\mathrm{d}t+\delta u(0) \right) - \tfrac12\,\tau^{-1}\,\lambda^{1/2}\,d(\lambda)\,\int_0^1 \delta u_{\overline{z}}(t)\,\sin(2\zeta(\lambda)t)\,\mathrm{d}t \\
& \qquad\qquad + \int_0^1 \delta u_{\overline{z}}(t)\,g_{1,c}(t)\,\mathrm{d}t + \int_0^1 \delta u(t)\,g_{2,c}(t)\,\mathrm{d}t \\
\delta d(\lambda) & = -\tfrac12\,d(\lambda)\,\int_0^1 \delta u_{\overline{z}}(t)\,\cos(2\zeta(\lambda)t)\,\mathrm{d}t - \tfrac12\,\tau\,\lambda^{-1/2}\,c(\lambda)\,\int_0^1 \delta u_{\overline{z}}(t)\,\sin(2\zeta(\lambda)t)\,\mathrm{d}t \\
& \qquad\qquad + \int_0^1 \delta u_{\overline{z}}(t)\,g_{1,d}(t)\,\mathrm{d}t + \int_0^1 \delta u(t)\,g_{2,d}(t)\,\mathrm{d}t \; . 
\end{align*}
}
\end{enumerate}
\end{lem}

\begin{rem}
If we apply Lemma~\ref{L:diffeo:delM-asymp} to the vacuum potential, i.e.~\,$(u,u_y)=(0,0)$\,, then all the functions \,$f_{\dotsc}$\, and \,$g_{\dotsc}$\, vanish (see the proof below). It follows that up to factors
which do not depend on \,$(\delta u,\delta u_y)$\,, the variations \,$\delta a_0(\lambda_{k,0})\cdot (\delta u,\delta u_y),\dotsc,\delta d_0(\lambda_{k,0})\cdot (\delta u,\delta u_y)$\, are the \,$|k|$-th Fourier coefficients of
\,$\delta u_z$\, (for \,$k>0$\,) or of \,$\delta u_{\overline{z}}$\, (for \,$k<0$\,). This is a linearised version of the observation that for the monodromy \,$M(\lambda)$\, corresponding to a potential \,$(u,u_y)\in\Pot$\,,
a first order approximation of \,$M(\lambda_{k,0})$\, is given by the \,$|k|$-th Fourier coefficients of \,$u_z$\, resp.~of \,$u_{\overline{z}}$\,, see Theorem~\ref{T:fasymp:fourier}.
\end{rem}

\begin{proof}[Proof of Lemma~\ref{L:diffeo:delM-asymp}.]
\emph{For (1).}
Let \,$(u,u_y) \in \Pot_{tame}$\, and \,$(\delta u,\delta u_y)\in T_{(u,u_y)}\Pot$\, be given.


We once again need to use the regauging of the extended frame which we described in the proof of Theorem~\ref{T:asymp:basic}. In the sequel, we will use the notations taken from that proof,
especially with regard to the regauging map \,$g$\, defined in Equation~\eqref{eq:asymp:expand:g}, the regauged 1-form \,$\wt{\alpha}=\wt{\alpha}_0+\beta+\gamma$\, from Equation~\eqref{eq:asymp:basic:wtalpha},
the regauged extended frame resp.~monodromy \,$\wt{F}(x)$\, resp.~\,$\wt{M}$\, defined by Equation~\eqref{eq:asymp:regauge}, and the corresponding vacuum solution \,$\wt{F}_0$\, given by Equation~\eqref{eq:asymp:E0}.


We begin by calculating \,$\delta \wt{F}(x)$\,. From the fact that \,$\wt{F}$\, satisfies the homogeneous linear initial value problem
$$ \wt{F}'(x) = \wt{\alpha}\,\wt{F}(x)\qmq{with} \wt{F}(0)=\unity\;, $$
it follows that \,$\delta\wt{F}$\, solves the inhomogeneous linear initial value problem
$$ (\delta\wt{F})'(x) = \wt{\alpha}\cdot \delta\wt{F} + \delta\wt{\alpha}\cdot \wt{F} \qmq{with} \delta\wt{F}(0)=0 \;, $$
and therefore \,$\delta\wt{F}$\, is given by
\begin{equation}
\label{eq:diffeo:delM-asymp:deltawtF-pre}
\delta\wt{F}(x) = \wt{F}(x)\,\int_{0}^x \wt{F}(t)^{-1}\cdot \delta\wt{\alpha}(t)\cdot \wt{F}(t)\,\mathrm{d}t \; .
\end{equation}
Now we have
$$ \wt{\alpha} = \wt{\alpha}_0 -\frac12\,u_z\,L + \gamma \qmq{with} \gamma = \frac{1}{4}\,\lambda^{-1/2}\,\begin{pmatrix} 0 & -e^u+1 \\ e^{-u}-1 & 0 \end{pmatrix} \qmq{and} 
L := \left( \begin{matrix} 1 & 0 \\ 0 & -1 \end{matrix} \right) $$
and therefore
$$ \delta \wt{\alpha} = -\frac12\,\delta u_z\,L + \delta\gamma \qmq{with} \delta\gamma = \frac14\,\lambda^{-1/2}\,\delta u\, \begin{pmatrix} 0 & -e^u \\ -e^{-u} & 0 \end{pmatrix} \; . $$
It therefore follows from Equation~\eqref{eq:diffeo:delM-asymp:deltawtF-pre} that we have
\begin{equation}
\label{eq:diffeo:delM-asymp:deltawtF}
\delta\wt{F}(x) = \wt{F}(x)\,\int_{0}^x \left( \delta u_z(t)\cdot \wt{X}(t) + \delta u_z(t)\cdot \wt{R}_{1}(t) + \delta u(t) \cdot \wt{R}_{2}(t) \right)\,\mathrm{d}t
\end{equation}
with
\begin{align}
\label{eq:diffeo:delM-asymp:wtX}
\wt{X}(t) & := -\frac12\,\wt{F}_0(t)^{-1}\,L\,\wt{F}_0(t) \overset{\eqref{eq:asymp:E0}}{=} \frac12 \begin{pmatrix} -\cos(2\,\zeta(\lambda)\,t) & \sin(2\,\zeta(\lambda)\,t) \\ \sin(2\,\zeta(\lambda)\,t) & \cos(2\,\zeta(\lambda)\,t) \end{pmatrix} \\
\label{eq:diffeo:delM-asymp:wtR1}
\wt{R}_{1}(t) & := -\frac12\,\bigr( (\wt{F}(t)^{-1}-\wt{F}_0(t)^{-1})\,L\,\wt{F}_0(t) + \wt{F}(t)^{-1}\,L\,(\wt{F}(t)-\wt{F}_0(t)) \bigr) \\
\label{eq:diffeo:delM-asymp:wtR2}
\wt{R}_{2}(t) & := \frac14\,\lambda^{-1/2}\,\wt{F}(t)^{-1}\,\begin{pmatrix} 0 & -e^{u(t)} \\ -e^{-u(t)} & 0 \end{pmatrix} \,\wt{F}(t) \; . 
\end{align}

We now undo the regauging by \,$g$\,, to obtain \,$\delta F$\, from \,$\delta \wt{F}$\,: 
By Equation~\eqref{eq:asymp:regauge} we have
$$ F(x) = g(x)\cdot \wt{F}(x) \cdot g(0)^{-1} $$
and therefore 
\begin{align*}
\delta F(x) & = g(x)\cdot \delta \wt{F}(x) \cdot g(0)^{-1} + \delta g(x)\cdot \wt{F}(x) \cdot g(0)^{-1} + g(x)\cdot \wt{F}(x) \cdot \delta(g(0)^{-1}) \\
& = g(x)\cdot \delta \wt{F}(x) \cdot g(0)^{-1} + \delta g(x)\cdot \wt{F}(x) \cdot g(0)^{-1} - g(x)\cdot \wt{F}(x) \cdot g(0)^{-1} \cdot \delta g(0) \cdot g(0)^{-1} \\
& = g(x)\cdot \delta \wt{F}(x) \cdot g(0)^{-1} + \delta g(x)\cdot g(x)^{-1}\cdot F(x) - F(x) \cdot \delta g(0) \cdot g(0)^{-1} \; . 
\end{align*}
When we set \,$x=1$\, in the above equation, we obtain (remember that the potential \,$(u,u_y)\in\Pot_{tame}$\, is periodic, and therefore \,$g(1)=g(0)$\, holds)
\begin{align}
\delta F(1) & = g(0)\cdot \delta \wt{F}(1) \cdot g(0)^{-1} + \delta g(0)\cdot g(0)^{-1}\cdot F(1) - F(1) \cdot \delta g(0) \cdot g(0)^{-1} \notag \\
& = g(0)\cdot \delta \wt{F}(1) \cdot g(0)^{-1} + [g(0)^{-1}\cdot \delta g(0), F(1)] \notag\\
& \overset{(*)}{=} g(0)\cdot \delta \wt{F}(1) \cdot g(0)^{-1} + \frac{1}{4}\,\delta u(0)\, [L, F(1)] \notag \\
& = g(0)\cdot \delta \wt{F}(1) \cdot g(0)^{-1} + \frac{1}{2}\,\delta u(0)\, \begin{pmatrix} 0 & -b(1,\lambda) \\ c(1,\lambda) & 0 \end{pmatrix} \notag \\
\label{eq:diffeo:delM-asymp:delF1}
& =  F(1)\,\int_{0}^1 \left( \delta u_z(t)\cdot X(t) + \delta u_z(t)\cdot R_{1}(t) + \delta u(t) \cdot R_{2}(t) \right)\,\mathrm{d}t \notag \\
& \qquad\qquad + \frac{1}{2}\,\delta u(0)\, \begin{pmatrix} 0 & -b(1,\lambda) \\ c(1,\lambda) & 0 \end{pmatrix} 
\end{align}
with (compare Equations~\eqref{eq:diffeo:delM-asymp:wtX}--\eqref{eq:diffeo:delM-asymp:wtR2})
\begin{align*}
X(t) & := g(0)\,\wt{X}(t)\,g(0)^{-1} = \frac12 \begin{pmatrix} -\cos(2\,\zeta(\lambda)\,t) & \tau^{-1}\,\lambda^{-1/2}\,\sin(2\,\zeta(\lambda)\,t) \\ \tau\,\lambda^{1/2}\,\sin(2\,\zeta(\lambda)\,t) & \cos(2\,\zeta(\lambda)\,t) \end{pmatrix} \\
R_1(t) & := g(0)\,\wt{R}_{1}(t)\,g(0)^{-1} = -\frac12\,\bigr( (F(t)^{-1}-\wh{F}_0(t)^{-1})\,L\,\wh{F}_0(t) + F(t)^{-1}\,L\,(F(t)-\wh{F}_0(t)) \bigr) \\
R_2(t) & := g(0)\,\wt{R}_{2}(t)\,g(0)^{-1} = \frac14\,\lambda^{-1/2}\,F(t)^{-1}\,\begin{pmatrix} 0 & -\lambda^{-1/2}\,\tau^{-1}\,e^{u(t)} \\ -\lambda^{1/2}\,\tau\,e^{-u(t)} & 0 \end{pmatrix} \,F(t) \\
\wh{F}_0(t) & := g(0)\,\wt{F}_0(t)\,g(0)^{-1} = \begin{pmatrix} a_0 & \tau^{-1}\,b_0 \\ \tau\,c_0 & d_0 \end{pmatrix} \; ;
\end{align*}
for the equals sign marked $(*)$ we note that we obtain from Equation~\eqref{eq:asymp:expand:g} 
\begin{equation*}
\delta g(x) = \frac14\,\delta u(x)\,\begin{pmatrix} e^{u/4} & 0 \\ 0 & -\lambda^{1/2}\,e^{-u/4} \end{pmatrix} \qmq{and therefore} g(x)^{-1} \cdot \delta g(x) = \frac14\,\delta u(x)\,L \; .
\end{equation*}

We now write the extended frame \,$F$\, in the form
$$ F = \begin{pmatrix} a & b \\ c & d \end{pmatrix} $$
(note that the functions \,$a,\dotsc,d$\, here depend not only on \,$\lambda\in \C^*$\,, but on \,$x\in [0,1]$\, as well). Because of \,$\det(F)=\det(\wh{F}_0)=1$\, we then have
$$ F^{-1} = \begin{pmatrix} d & -b \\ -c & a \end{pmatrix} \qmq{and} \wh{F}_0^{-1} = \begin{pmatrix} d_0 & -\tau^{-1}\,b_0 \\ -\tau\,c_0 & a_0 \end{pmatrix} \;. $$
By carrying out the matrix multiplications in Equation~\eqref{eq:diffeo:delM-asymp:delF1}, and using these formulae, we see that 
the component functions of \,$\delta M = \delta F(1)$\, are indeed of the form claimed in part (1) of the Theorem, where 
\begin{align*}
f_{a,1} & = -\tfrac12\bigr( (d-d_0)a_0 + (b-\tau^{-1}\,b_0)c_0 \bigr) & f_{a,2} & = \tfrac14\,\bigr(\tau\,a\,b\,e^{-u(t)} - \tau^{-1}\,c\,d\,\lambda^{-1}\,e^{u(t)}\bigr) \\
f_{b,1} & = -\tfrac12\bigr( (d-d_0)b_0 + (b-\tau^{-1}\,b_0)d_0 \bigr) & f_{b,2} & = \tfrac14\,\bigr( \tau\,b^2\,e^{-u(t)} -\tau^{-1}\,d^2\,\lambda^{-1}\,e^{u(t)} \bigr) \\
f_{c,1} & = -\tfrac12\bigr( -(c-\tau\,c_0)a_0 + (a-a_0)c_0 \bigr) & f_{c,2} & = \tfrac14\,\bigr(\tau\,a^2\,e^{-u(t)} + \tau^{-1}\,c^2\,\lambda^{-1}\,e^{u(t)} \bigr) \\
f_{d,1} & = -\tfrac12\bigr( -(c-\tau\,c_0)b_0+(a-a_0)d_0 \bigr) & f_{d,2} & = \tfrac14\,\bigr( -\tau\,a\,b\,e^{-u(t)} + \tau^{-1}\,c\,d\,\lambda^{-1}\,e^{u(t)} \bigr) \; . 
\end{align*}
It remains to show that these functions have the asymptotic behavior stated in part (1) of the Theorem, and this follows from the asymptotic behavior of the extended frame \,$F$\, described in 
Proposition~\ref{P:asympfinal:frame}, together with the fact that the functions \,$e^{\pm u(t)}$\, are (continuous and therefore) bounded for \,$t\in [0,1]$\,.

\emph{For (2).}
We consider besides the given potential \,$(u,u_y)$\, also the potential
\,$(\wt{u},\wt{u}_y) := (-u,u_y)$\,, and denote the quantities associated to the latter potential by a tilde. We will reduce (2) for \,$(u,u_y)$\, to the statement of (1) applied to \,$(\wt{u},\wt{u}_y)$\,. 

We have \,$\wt{\tau}=\tau^{-1}$\,,
$$ (\delta \wt{u},\delta \wt{u}_y) = (-\delta u,\delta u_y) $$
and therefore
$$ \delta \wt{u} = -\delta u \qmq{and} \delta \wt{u}_z = -\delta u_{\overline{z}} \; . $$
By Proposition~\ref{P:mono:symmetry}(1) we have
$$ M(\lambda^{-1}) = g^{-1} \cdot \wt{M}(\lambda) \cdot g \qmq{with} g := \begin{pmatrix} 1 & 0 \\ 0 & \lambda \end{pmatrix} \;, $$
i.e.
$$ \begin{pmatrix} a(\lambda^{-1}) & b(\lambda^{-1}) \\ c(\lambda^{-1}) & d(\lambda^{-1}) \end{pmatrix}
= \begin{pmatrix} \wt{a}(\lambda) & \lambda\cdot \wt{b}(\lambda) \\ \lambda^{-1}\cdot \wt{c}(\lambda) & \wt{d}(\lambda) \end{pmatrix} \; . $$
By differentiation of this equation with respect to \,$(u,u_y)$\, it follows:
$$ \begin{pmatrix} \delta a(\lambda^{-1}) & \delta b(\lambda^{-1}) \\ \delta c(\lambda^{-1}) & \delta d(\lambda^{-1}) \end{pmatrix}
= \begin{pmatrix} \delta \wt{a}(\lambda) & \lambda\cdot \delta \wt{b}(\lambda) \\ \lambda^{-1}\cdot \delta \wt{c}(\lambda) & \delta \wt{d}(\lambda) \end{pmatrix} \; . $$
Thus we obtain via the result of (1) 
\begin{align*}
\delta a(\lambda^{-1}) 
& = \delta\wt{a}(\lambda) \\
& = -\tfrac12\,\wt{a}(\lambda)\,\int_0^1 \delta \wt{u}_z(t)\,\cos(2\zeta(\lambda)t)\,\mathrm{d}t + \tfrac12\,\wt{\tau}\,\lambda^{1/2}\,\wt{b}(\lambda)\,\int_0^1 \delta \wt{u}_z(t)\,\sin(2\zeta(\lambda)t)\,\mathrm{d}t \\
& \qquad\qquad + \int_0^1 \delta \wt{u}_z(t)\,\wt{f}_{1,a}(t)\,\mathrm{d}t + \int_0^1 \delta \wt{u}(t)\,\wt{f}_{2,a}(t)\,\mathrm{d}t \\
& = \tfrac12\,a(\lambda^{-1})\,\int_0^1 \delta u_{\overline{z}}(t)\,\cos(2\zeta(\lambda^{-1})t)\,\mathrm{d}t \\
& \qquad\qquad - \tfrac12\,\tau^{-1}\,\lambda^{1/2}\,\lambda^{-1}\,b(\lambda^{-1})\,\int_0^1 \delta u_{\overline{z}}(t)\,\sin(2\zeta(\lambda^{-1})t)\,\mathrm{d}t \\
& \qquad\qquad + \int_0^1 \delta u_{\overline{z}}(t)\,g_{1,a}(t)\bigr|_{\lambda^{-1}}\,\mathrm{d}t + \int_0^1 \delta u(t)\,g_{2,a}(t)\bigr|_{\lambda^{-1}}\,\mathrm{d}t 
\end{align*}
with the functions \,$g_{1,a} := -\wt{f}_{1,a}\bigr|_{\lambda^{-1}}$\, and \,$g_{2,a} := -\wt{f}_{2,a}\bigr|_{\lambda^{-1}}$\,.
The uniform asymptotic properties \,$\wt{f}_{1,a} \in \As_\infty(\C^*,\ell^2_0,2)$\, and \,$\wt{f}_{2,a} \in \As_\infty(\C^*,\ell^\infty_1,2)$\,
imply \,$g_{1,a} \in \As_0(\C^*,\ell^2_0,2)$\, and \,$g_{2,a} \in \As_0(\C^*,\ell^\infty_1,2)$\,. By substituting \,$\lambda$\, for \,$\lambda^{-1}$\, in the above equation we obtain the result for \,$\delta a$\, in (2). 
The asymptotic equations for \,$\delta b,\delta c,\delta d$\, are shown in the same way.
%
\end{proof}

We next apply Lemma~\ref{L:diffeo:delM-asymp} to describe the asymptotic behavior of \,$\delta \lambda_k$\, and \,$\delta \mu_k$\,:

\begin{lem}
\label{L:diffeo:delD-asymp}
Let \,$(u,u_y)\in \Pot_{tame}$\, be given, and let \,$D = \{(\lambda_k,\mu_k)\} \in \Div_{tame}$\, be the spectral divisor of \,$(u,u_y)$\,. 
Because \,$D$\, does not contain any double points, the functions \,$\lambda_k$\, and \,$\mu_k$\,
are well-defined and smooth on neighborhoods of \,$(u,u_y)$\, in \,$\Pot$\,. In this setting, there exist sequences
\,$(a_k) \in \ell^2_{-1,3}(k)$\, and \,$(b_k) \in \ell^2_{0,0}(k)$\, of non-negative real numbers, 
so that we have for \,$(\delta u,\delta u_y)\in T_{(u,u_y)}\Pot$\,:
\begin{enumerate}
\item For \,$k>0$\,:
\begin{align*}
\left| \delta \lambda_k - (-4)\,(-1)^k\,\lambda_k^{1/2}\,\mu_k^{-1}\,\int_0^1 \delta u_z(t)\,\sin(2k\pi t)\,\mathrm{d}t \right|
& \leq a_k \cdot \|(\delta u,\delta u_y)\|_{\Pot} \\
\left| \delta \mu_k - \left(-\frac12\right) \,\mu_k\,\int_0^1 \delta u_z(t)\,\cos(2k\pi t)\,\mathrm{d}t \right|
& \leq b_k \cdot \|(\delta u,\delta u_y)\|_{\Pot} 
\end{align*}
\item For \,$k<0$\,:
\begin{align*}
\left| \delta \lambda_k - 4\,(-1)^k\,\lambda_k^{3/2}\,\mu_k^{-1}\,\int_0^1 \delta u_{\overline{z}}(t)\,\sin(2k\pi t)\,\mathrm{d}t \right|
& \leq a_k \cdot \|(\delta u,\delta u_y)\|_{\Pot} \\
\left| \delta \mu_k - \frac12 \,\mu_k\,\int_0^1 \delta u_{\overline{z}}(t)\,\cos(2k\pi t)\,\mathrm{d}t \right|
& \leq b_k \cdot \|(\delta u,\delta u_y)\|_{\Pot} 
\end{align*}
\end{enumerate}
Moreover, \,$\delta \lambda_k$\, and \,$\delta \mu_k$\, depend continuously on \,$(\delta u,\delta u_y)$\,. 
\end{lem}

\begin{proof}
\emph{For (1).} We have \,$k>0$\, here, and use the notations of Lemma~\ref{L:diffeo:delM-asymp}.
Because \,$D$\, is tame, we have \,$c'(\lambda_k)\neq 0$\,, and therefore
the implicit function theorem shows that the function \,$\lambda_k$\, is smooth, and that we have
\begin{equation}
\label{eq:diffeo:delD-asymp:deltalambda-pre}
\delta \lambda_k = -\frac{1}{c'(\lambda_k)}\,(\delta c)(\lambda_k) \; .
\end{equation}
We have \,$c'(\lambda_k) = \left. \tfrac{c(\lambda)}{\lambda-\lambda_k} \right|_{\lambda=\lambda_k}$\,, and therefore
Corollary~\ref{C:interpolate:cdivlin}(1) shows that 
$$ c'(\lambda_k) - \tau\,\frac{(-1)^k}{8} \in \ell^2_0(k) $$
holds; because the function \,$z\mapsto z^{-1}$\, is locally Lipschitz continuous near \,$z=\tau\,\tfrac{(-1)^k}{8}$\,, it follows that
we also have
\begin{equation}
\label{eq:diffeo:delD-asymp:c'inverse}
\frac{1}{c'(\lambda_k)} = \frac{8\,(-1)^k}{\tau} + r_k^{[1]} 
\end{equation}
with a sequence \,$(r_k^{[1]}) \in \ell^2_0(k)$\,. 

Moreover, we estimate \,$(\delta c)(\lambda_k)$\, via Lemma~\ref{L:diffeo:delM-asymp}(1): Because we have \,$c(\lambda_k)=0$\,
and \,$d(\lambda_k)=\mu_k^{-1}$\,, we obtain 
\begin{equation}
\label{eq:diffeo:delD-asymp:deltac}
(\delta c)(\lambda_k) = \frac12\,\tau\,\lambda_k^{1/2}\,\mu_k^{-1}\,\int_{0}^1 \delta u_z(t)\,\sin(2\zeta(\lambda_k)t)\,\mathrm{d}t
 + \int_0^1 \delta u_z(t)\,f_{1,c}(t)\,\mathrm{d}t + \int_0^1 \delta u(t)\,f_{2,c}(t)\,\mathrm{d}t 
\end{equation}
with functions \,$f_{1,c} \in \As_\infty(\C^*,\ell^2_{-1},2)$\, and \,$f_{2,c}\in \As_\infty(\C^*,\ell^\infty_0,2)$\, uniformly
in \,$t\in[0,1]$\,. We let \,$(r_k^{[f_{c,1}]})\in \ell^2_{-1}(k)$\, resp.~\,$(r_k^{[f_{c,2}]})\in \ell^\infty_0(k)$\, be bounding
sequences for \,$f_{c,1}(t)$\, resp.~\,$f_{c,2}(t)$\, for all \,$t\in [0,1]$\,. Then we have
\begin{equation}
\label{eq:diffeo:delD-asymp:deltac-estim}
|(\delta c)(\lambda_k)| \leq \left( |\tau|\,|\lambda_k|^{1/2}\,|\mu_k|^{-1}+r_k^{[f_{1,c}]}\right) \cdot \|\delta u_z\|_{L^2([0,1])} 
+ r_k^{[f_{2,c}]}\cdot \|\delta u\|_{L^2([0,1])} 
\end{equation}
by Cauchy-Schwarz's inequality, and the fact that \,$|\sin(2\zeta(\lambda_k)t)|\leq 2$\, holds 
for all \,$k\in \Z$\, and \,$t\in [0,1]$\,. 

We also note that we have \,$\lambda_k-\lambda_{k,0}\in \ell^2_{-1}(k>0)$\, and therefore \,$\zeta(\lambda_k)-\pi k = \zeta(\lambda_k)-\zeta(\lambda_{k,0}) \in \ell^2_0(k>0)$\,; because \,$\sin(z)$\, is Lipschitz continuous
on any horizontal strip in the complex plane, it follows that there exists a sequence \,$(r_k^{[2]})\in \ell^2(k)$\, such that
\begin{equation}
\label{eq:diffeo:delD-asymp:sinsin}
|\sin(2\,\zeta(\lambda_k)\,t) - \sin(2\pi k t)| \leq r_k^{[2]}
\end{equation}
holds for all \,$k>0$\, and \,$t\in [0,1]$\,. 


We now apply the preceding estimates to Equation~\eqref{eq:diffeo:delD-asymp:deltalambda-pre} to obtain the asymptotic assessment for \,$\delta \lambda_k$\, claimed in the Lemma. Indeed, we have
\begin{align*}
\delta \lambda_k 
& \overset{\eqref{eq:diffeo:delD-asymp:deltalambda-pre}}{=} -\frac{1}{c'(\lambda_k)}\,(\delta c)(\lambda_k) \\
& \;\;\,=\;\;\, -\frac{8\,(-1)^k}{\tau}\,(\delta c)(\lambda_k) - \left( \frac{1}{c'(\lambda_k)} - \frac{8\,(-1)^k}{\tau} \right)\,(\delta c)(\lambda_k) \\
& \overset{\eqref{eq:diffeo:delD-asymp:deltac}}{=} 
-4\,(-1)^k\,\lambda_k^{1/2}\,\mu_k^{-1}\,\int_{0}^1 \delta u_z(t)\,\sin(2\zeta(\lambda_k)t)\,\mathrm{d}t \\
& \qquad\qquad - \frac{8\,(-1)^k}{\tau}\left( \int_0^1 \delta u_z(t)\,f_{1,c}(t)\,\mathrm{d}t + \int_0^1 \delta u(t)\,f_{2,c}(t)\,\mathrm{d}t \right) \\
& \qquad\qquad - \left( \frac{1}{c'(\lambda_k)} - \frac{8\,(-1)^k}{\tau} \right)\,(\delta c)(\lambda_k) 
\end{align*}
and therefore
\begin{align*}
& \delta \lambda_k - (-4)\,(-1)^k\,\lambda_k^{1/2}\,\mu_k^{-1}\,\int_0^1 \delta u_z(t)\,\sin(2k\pi t)\,\mathrm{d}t \\
= & -4\,(-1)^k\,\lambda_k^{1/2}\,\mu_k^{-1}\,\int_{0}^1 \delta u_z(t)\,\bigr(\sin(2\zeta(\lambda_k)t)-\sin(2k\pi t)\bigr) \,\mathrm{d}t \\
& \qquad\qquad - \frac{8\,(-1)^k}{\tau}\left( \int_0^1 \delta u_z(t)\,f_{1,c}(t)\,\mathrm{d}t + \int_0^1 \delta u(t)\,f_{2,c}(t)\,\mathrm{d}t \right) \\
& \qquad\qquad - \left( \frac{1}{c'(\lambda_k)} - \frac{8\,(-1)^k}{\tau} \right)\,(\delta c)(\lambda_k)  \; . 
\end{align*}
By use of the triangle inequality, Cauchy-Schwarz's inequality and the estimates \eqref{eq:diffeo:delD-asymp:c'inverse}, 
\eqref{eq:diffeo:delD-asymp:deltac-estim} and \eqref{eq:diffeo:delD-asymp:sinsin}, we obtain
\begin{align*}
& \left| \delta \lambda_k - (-4)\,(-1)^k\,\lambda_k^{1/2}\,\mu_k^{-1}\,\int_0^1 \delta u_z(t)\,\sin(2k\pi t)\,\mathrm{d}t \right| \\
\leq\; & a_k^{[1]} \cdot \|\delta u_z\|_{L^2([0,1])} + a_k^{[2]} \cdot \|\delta u\|_{L^2([0,1])} \\
\leq\; & (a_k^{[1]} + a_k^{[2]}) \cdot \|(\delta u,\delta u_y)\|_{\Pot}
\end{align*}
with
\begin{align*}
a_k^{[1]} & := 4\,|\lambda_k|^{1/2}\,|\mu_k|^{-1}\,r_k^{[2]} + 8\,|\tau|^{-1}\,r_k^{[f_{1,c}]} + r_k^{[1]}\cdot \left( |\tau|\,|\lambda_k|^{1/2}\,|\mu_k|^{-1}+r_k^{[f_{1,c}]}\right) \;, \\
a_k^{[2]} & := 8\,|\tau|^{-1}\,r_k^{[f_{2,c}]} + r_k^{[1]} \cdot r_k^{[f_{2,c}]} \;.
\end{align*}
Because we have \,$a_k^{[1]},a_k^{[2]} \in \ell^2_{-1}(k)$\,, the asymptotic estimate claimed for \,$\delta \lambda_k$\, in part (1) of
the lemma follows.

For the estimate for \,$\delta \mu_k$\,, we first note that \,$\mu_k = a(\lambda_k)$\, is a smooth function near \,$(u,u_y)$\,. However, the computation of \,$\delta \mu_k$\, becomes slightly easier if
we base our calculation on \,$\mu_k^{-1} = d(\lambda_k)$\, instead. The reason for this is that the expressions for \,$\delta d(\lambda_k)$\, in Lemma~\ref{L:diffeo:delM-asymp} are somewhat simpler than
those for \,$\delta a(\lambda_k)$\,, owing to the fact that \,$c(\lambda_k)=0$\, holds, whereas we do not know anything about \,$b(\lambda_k)$\,. 

Thus we calculate:
\begin{equation}
\label{eq:diffeo:delD-asymp:deltamu-pre}
\delta \mu_k = \delta\bigr( (\mu_k^{-1})^{-1} \bigr) = -\mu_k^{2}\cdot \delta(\mu_k^{-1}) = -\mu_k^{2}\cdot \delta(d(\lambda_k)) = -\mu_k^{2} \cdot \left( \delta d(\lambda_k) + d'(\lambda_k)\cdot \delta \lambda_k \right) \; .
\end{equation}
For \,$\delta d(\lambda_k)$\,, we have by Lemma~\ref{L:diffeo:delM-asymp}(1) (again noting that \,$c(\lambda_k)=0$\, and 
\,$d(\lambda_k)=\mu_k^{-1}$\, holds):
\begin{equation}
\label{eq:diffeo:delD-asymp:deltad}
\delta d(\lambda_k) = \tfrac12\,\mu_k^{-1}\,\int_0^1 \delta u_z(t)\,\cos(2\zeta(\lambda_k)t)\,\mathrm{d}t 
+ \int_0^1 \delta u_z(t)\,f_{1,d}(t)\,\mathrm{d}t + \int_0^1 \delta u(t)\,f_{2,d}(t)\,\mathrm{d}t
\end{equation}
with functions \,$f_{1,d} \in \As_\infty(\C^*,\ell^2_0,2)$\, and \,$f_{2,d}\in \As_\infty(\C^*,\ell^\infty_1,2)$\, uniformly
in \,$t\in[0,1]$\,. We let \,$(r_k^{[f_{d,1}]})\in \ell^2_{0}(k)$\, resp.~\,$(r_k^{[f_{d,2}]})\in \ell^\infty_1(k)$\, be uniform bounding
sequences for \,$f_{d,1}(t)$\, resp.~\,$f_{d,2}(t)$\, for all \,$t\in [0,1]$\,.

Next, we estimate \,$d'(\lambda_k)$\,. We have \,$d_0'(\lambda_{k,0})=0$\, and therefore
$$ d'(\lambda_k) = d'(\lambda_k)-d_0'(\lambda_{k,0}) = d'(\lambda_k)-d_0'(\lambda_k) + \int_{\lambda_{k,0}}^{\lambda_k} d_0''(\lambda)\,\mathrm{d}\lambda \; . $$
We have \,$d'-d_0' \in \As(\C^*,\ell^2_{1,-3},1)$\, and therefore \,$d'(\lambda_k)-d_0'(\lambda_k)\in \ell^2_{1,-3}(k)$\,. Moreover,
we have \,$d_0'' \in \As(\C^*,\ell^\infty_{2,-6},1)$\, and \,$\lambda_k-\lambda_{k,0} \in \ell^2_{-1,3}(k)$\, and therefore 
\,$\int_{\lambda_{k,0}}^{\lambda_k} d_0''(\lambda)\,\mathrm{d}\lambda \in \ell^2_1(k)$\,. Thus we obtain
\begin{equation}
\label{eq:diffeo:delD-asymp:d'-estim}
|d'(\lambda_k)| =: r_k^{[3]} \in \ell^2_1(k) \; . 
\end{equation}
Moreover, it follows from the asymptotic estimate for \,$\delta \lambda_k$\, shown above that we have
\begin{equation}
\label{eq:diffeo:delD-asymp:deltalambda-estim}
|\delta \lambda_k| \leq \wt{a}_k \cdot \|(\delta u,\delta u_y)\|_{\Pot}
\end{equation}
with the sequence
$$ \wt{a}_k := 4\,|\lambda_k|^{1/2}\,|\mu_k|^{-1} + a_k \in \ell^\infty_{-1}(k) \; . $$

Similarly as Equation~\eqref{eq:diffeo:delD-asymp:sinsin}, we note moreover that there exists
a sequence \,$(r_k^{[4]})\in \ell^2(k)$\, such that
\begin{equation}
\label{eq:diffeo:delD-asymp:coscos}
|\cos(2\,\zeta(\lambda_k)\,t) - \cos(2\pi k t)| \leq r_k^{[4]}
\end{equation}
holds for all \,$k>0$\, and \,$t\in [0,1]$\,. 

We now obtain
\begin{align*}
\delta \mu_k
& \overset{\eqref{eq:diffeo:delD-asymp:deltamu-pre}}{=} 
-\mu_k^{2} \cdot \left( \delta d(\lambda_k) + d'(\lambda_k)\cdot \delta \lambda_k \right) \\
& \overset{\eqref{eq:diffeo:delD-asymp:deltad}}{=} 
-\tfrac12\,\mu_k\,\int_0^1 \delta u_z(t)\,\cos(2\zeta(\lambda_k)t)\,\mathrm{d}t 
-\mu_k^{2} \int_0^1 \delta u_z(t)\,f_{1,d}(t)\,\mathrm{d}t \\
& \qquad\qquad -\mu_k^{2} \int_0^1 \delta u(t)\,f_{2,d}(t)\,\mathrm{d}t -\mu_k^{2} \cdot d'(\lambda_k)\cdot \delta \lambda_k 
\end{align*}
and therefore
\begin{align*}
& \delta \mu_k - \left(-\tfrac12\right) \,\mu_k\,\int_0^1 \delta u_z(t)\,\cos(2k\pi t)\,\mathrm{d}t \\
= & -\tfrac12\,\mu_k^{2}\,\int_0^1 \delta u_z(t)\,\bigr( \cos(2\zeta(\lambda_k)t)-\cos(2k\pi t)\bigr) \,\mathrm{d}t 
-\mu_k^{2} \int_0^1 \delta u_z(t)\,f_{1,d}(t)\,\mathrm{d}t \\
& \qquad\qquad -\mu_k^{2} \int_0^1 \delta u(t)\,f_{2,d}(t)\,\mathrm{d}t -\mu_k^{2} \cdot d'(\lambda_k)\cdot \delta \lambda_k \; . 
\end{align*}
By use of the triangle inequality, Cauchy-Schwarz's inequality and the estimates \eqref{eq:diffeo:delD-asymp:d'-estim}, 
\eqref{eq:diffeo:delD-asymp:deltalambda-estim} and \eqref{eq:diffeo:delD-asymp:coscos}, we obtain
\begin{align*}
& \left| \delta \mu_k - \left(-\tfrac12\right) \,\mu_k\,\int_0^1 \delta u_z(t)\,\cos(2k\pi t)\,\mathrm{d}t \right| \\
\leq \; & b_k^{[1]} \cdot \|\delta u_z\|_{L^2([0,1])} + b_k^{[2]} \cdot \|\delta u\|_{L^2([0,1])} + r_k^{[3]} \cdot \wt{a}_k\cdot \|(\delta u,\delta u_y)\|_{\Pot} \\
\leq \; & (b_k^{[1]} + b_k^{[2]} + r_k^{[3]} \cdot \wt{a}_k) \cdot \|(\delta u,\delta u_y)\|_{\Pot}
\end{align*}
with
$$ b_k^{[1]} := -\frac12\,|\mu_k|^{2}\,r_k^{[4]} + |\mu_k|^{2}\,r_k^{[f_{1,d}]} \qmq{and}
b_k^{[2]} := |\mu_k|^2\,r_k^{[f_{2,d}]} \; . $$
Because we have \,$b_k^{[1]},b_k^{[2]}, r_k^{[3]}\cdot \wt{a}_k \in \ell^2_{0}(k)$\,, and moreover \,$\mu_k-\mu_k^{-3} \in \ell^2_0(k)$\, holds, the asymptotic estimate claimed for \,$\delta \mu_k$\, 
in part (1) of the lemma follows.


\emph{For (2).} This is shown similarly as (1), using the asymptotic estimates for \,$\delta M$\, for \,$\lambda\to 0$\,
from Lemma~\ref{L:diffeo:delM-asymp}(2).

Equations~\eqref{eq:diffeo:delD-asymp:deltalambda-pre} and \eqref{eq:diffeo:delD-asymp:deltamu-pre} also show that \,$\delta \lambda_k$\, and \,$\delta \mu_k$\, depend continuously on \,$(\delta u,\delta u_y)$\,. 
\end{proof}

\begin{lem}
\label{L:darboux:convergence}
Let \,$(u,u_y) \in \Pot_{tame}$\, and \,$(\delta u,\delta u_y), (\wt{\delta} u,\wt{\delta} u_y) \in T_{(u,u_y)}\Pot$\, be given. 
Then the sum
$$ \frac12\,\sum_{k\in \Z} \left( \frac{\delta \lambda_k}{\lambda_k} \cdot \frac{\wt{\delta} \mu_k}{\mu_k} - \frac{\wt{\delta} \lambda_k}{\lambda_k} \cdot \frac{\delta \mu_k}{\mu_k} \right)  $$
converges absolutely, and it depends continuously on \,$(\delta u,\delta u_y)$\, and \,$(\wt{\delta} u,\wt{\delta} u_y)$\,. 
\end{lem}

\begin{proof}
Because of \,$(\delta u,\delta u_y) \in T_{(u,u_y)}\Pot \subset W^{1,2}([0,1]) \times L^2([0,1])$\,, we have \,$\delta u_z, \delta u_{\overline{z}} \in L^2([0,1])$\,, and therefore the Fourier coefficients
$$ \int_0^1 \delta u_z(t)\,\cos(2k\pi t)\,\mathrm{d}t \qmq{and} \int_0^1 \delta u_z(t)\,\sin(2k\pi t)\,\mathrm{d}t $$
and 
$$ \int_0^1 \delta u_{\overline{z}}(t)\,\cos(2k\pi t)\,\mathrm{d}t \qmq{and} \int_0^1 \delta u_{\overline{z}}(t)\,\sin(2k\pi t)\,\mathrm{d}t $$
are in \,$\ell^2(k)$\,. Lemma~\ref{L:diffeo:delD-asymp} therefore shows that 
$$ \delta \lambda_k \in \ell^2_{-1,3}(k) \qmq{and} \delta \mu_k \in \ell^2_{0,0}(k) $$
holds. Moreover, because of \,$\lambda_k = \lambda_{k,0} + \ell^2_{-1,3}(k)$\, we have \,$\frac{1}{\lambda_k} \in \ell^\infty_{2,-2}(k)$\,, and because of \,$\mu_k = \mu_{k,0} + \ell^2_{0,0}(k)$\,, we have
\,$\tfrac{1}{\mu_k} \in \ell^\infty_{0,0}(k)$\,. Thus we obtain
$$ \frac{\delta \lambda_k}{\lambda_k} \in \ell^2_{1,1}(k) \qmq{and} \frac{\delta \mu_k}{\mu_k} \in \ell^2_{0,0}(k) \; . $$
By applying this result both for \,$(\delta u,\delta u_y)$\, and \,$(\wt{\delta} u,\wt{\delta} u_y)$\, we see that
$$ \frac{\delta \lambda_k}{\lambda_k} \cdot \frac{\wt{\delta} \mu_k}{\mu_k} - \frac{\wt{\delta} \lambda_k}{\lambda_k} \cdot \frac{\delta \mu_k}{\mu_k} \in \ell^1_{1,1}(k) \subset \ell^1(k) $$
holds, which means that the sum under investigation converges absolutely in \,$\C$\,.

Moreover, the summands depend continuously on \,$(\delta u,\delta u_y)$\, and \,$(\wt{\delta} u,\wt{\delta} u_y)$\, by Lem\-ma~\ref{L:diffeo:delD-asymp}. Because the above estimates show
that the convergence of the sum is locally uniform, it follows that also the sum depends continuously on \,$(\delta u,\delta u_y)$\, and \,$(\wt{\delta} u,\wt{\delta} u_y)$\,. 
\end{proof}

We now have all the pieces in place to complete the proof of Theorem~\ref{T:darboux:darboux}.

\begin{proof}[Proof of Theorem~\ref{T:darboux:darboux}.]
\emph{For (1).} This is Lemma~\ref{L:darboux:basis}.

\emph{For (2).}
By Theorem~\ref{T:darboux:knopf}(2) (=\cite{Knopf:2015}, Theorem~5.5), Equation~\eqref{eq:darboux:darboux:Omega} holds for all \,$(\delta u,\delta u_y), (\wt{\delta} u,\wt{\delta} u_y) \in T_{(u,u_y)}\Pot$\, that are
\emph{finite} linear combinations of the \,$v_k$\, and the \,$w_k$\,. Because the vectors \,$(\delta u,\delta u_y)$\, with this property are dense in \,$T_{(u,u_y)}\Pot$\, by Lemma~\ref{L:darboux:basis},
and both sides of Equation~\eqref{eq:darboux:darboux:Omega} are well-defined and continuous for all \,$(\delta u,\delta u_y), (\wt{\delta} u,\wt{\delta} u_y) \in T_{(u,u_y)}\Pot$\,
(the right-hand side by Lemma~\ref{L:darboux:convergence}), it follows that Equation~\eqref{eq:darboux:darboux:Omega} holds for all \,$(\delta u,\delta u_y), (\wt{\delta} u,\wt{\delta} u_y) \in T_{(u,u_y)}\Pot$\,.
\end{proof}

\section{The inverse problem for \,$\Pot\to\Div$\,}
\label{Se:diffeo}

In Section~\ref{Se:special} we showed that the monodromy \,$M(\lambda)$\, of a tame potential \,$(u,u_y)\in \Pot$\, is uniquely determined by its (tame) spectral divisor \,$D$\,, 
at least up to a change of sign of the off-diagonal entries. 
We are now ready to prove the important result that also the potential \,$(u,u_y)$\, itself is uniquely determined by \,$D$\,, at least if \,$D$\, is tame. In fact, we will show that for the map
\,$\Phi: \Pot \to \Div$\, maps each potential \,$(u,u_y)\in\Pot$\, onto its spectral divisor, \,$\Phi|\Pot_{tame}: \Pot_{tame}\to\Div_{tame}$\, is a diffeomorphism onto an open and dense subset of
\,$\Div_{tame}$\,. We will obtain this result relatively easily by using the Darboux coordinates for \,$\Pot_{tame}$\, constructed in the preceding section.

\newpage

We will prove
\begin{thm}
\label{T:diffeo:diffeo}
\,$\Phi|\Pot_{tame}:\Pot_{tame} \to \Div_{tame}$\, is a diffeomorphism onto an open and dense subset of \,$\Div_{tame}$\,. 
\end{thm}

\begin{rem}
\,$\Phi$\, is not immersive at potentials \,$(u,u_y) \in \Pot\setminus \Pot_{tame}$\,. This is to be expected if the classical spectral divisor \,$D\in \Div$\, corresponding to \,$(u,u_y)$\, contains singular points
of the associated spectral curve \,$\Sigma$\, with multiplicity \,$\geq 2$\,, because then not even the monodromy \,$M(\lambda)$\, is uniquely determined by \,$D$\,, as we need the additional information
contained in the generalized spectral divisor \,$\calD$\, to reconstruct \,$M(\lambda)$\, in this case (see Section~\ref{Se:special}). 
However, even if all the points occurring in \,$D$\, with multiplicity \,$\geq 2$\, are regular points of \,$\Sigma$\, (then the monodromy \,$M(\lambda)$\, is uniquely determined by \,$D$\,), there is an
entire family of integral curves of \,$x$-translation in \,$\Div$\, that intersect in such a point, and therefore \,$\Phi$\, cannot be immersive.
To make \,$\Phi$\, an immersion (and consequently a local diffeomorphism) near such points, we
would therefore need to replace the range \,$\Div$\, of \,$\Phi$\, by a suitable blow-up at its singularities (see~e.g.~\cite{Hartshorne:1977}, p.~163ff.). We do not carry out such a construction here.

The reason why the image of \,$\Phi|\Pot_{tame}$\, is not all of \,$\Div_{tame}$\, is that even though any tame divisor \,$D\in\Div_{tame}$\, is non-special, it is possible for \,$D$\, to become special under 
\,$x$-translation. 
If this occurs, the potential corresponding to \,$D$\, has a singularity for the corresponding value of \,$x$\,, and thus \,$D$\, cannot correspond to a potential \,$(u,u_y)\in \Pot$\, in our sense.
The investigation of sinh-Gordon potentials with singularities, corresponding to divisors \,$D\in\Div$\, which become special under \,$x$-translation for some value of \,$x$\,, would be extremely interesting
for the study of compact constant mean curvatures; we discuss this perspective in the concluding Section~\ref{Se:perspectives}.
\end{rem}

For the proof of Theorem~\ref{T:diffeo:diffeo}, we will use a well-known fact from the theory of solutions of the sinh-Gordon equation of finite type: For any divisor \,$D\in\Div$\, of finite type
(see Definition~\ref{D:finite:finite}),
such that the \,$x$-translations \,$D(x)$\, are non-special for every \,$x\in [0,1]$\,, there exists one and only one potential \,$(u,u_y)\in \Pot$\, with \,$\Phi((u,u_y))=D$\,. 

By the \,$x$-translation of a divisor \,$D\in \Div$\, we here mean the following: Suppose that a potential \,$(u,u_y)\in \Pot$\, and \,$x\in \R$\, are given, and denote the translation by \,$x$\, on \,$\R$\,
by \,$L_x(x') := x'+x$\,. Then \,$(u\circ L_x,u_y \circ L_x)$\, is another (periodic) potential in \,$\Pot$\,, and if \,$(\Sigma,D)$\, are the spectral data of \,$(u,u_y)$\,, then the spectral curve
of \,$(u\circ L_x,u_y\circ L_x)$\, will also be \,$\Sigma$\, for all \,$x\in \R$\,, whereas we denote the spectral divisor 
of \,$(u\circ L_x,u_y \circ L_x)$\, by \,$D(x)\in\Div$\,; this is the \,$x$-translation of the divisor \,$D$\,. 

At least if \,$D\in\Div$\, is tame, we can characterize the \,$x$-translation of \,$D$\, by differential equations, without recourse to a potential \,$(u,u_y)$\, of which \,$D$\, is the spectral divisor:
We write \,$D=\{(\lambda_k,\mu_k)\}$\,, and regard \,$\lambda_k$\, and \,$\mu_k$\, as functions in \,$x\in \R$\,. Because \,$D(x)$\, is a divisor on the spectral curve \,$\Sigma$\, associated to \,$D$\, for all \,$x\in \R$\,,
it suffices to describe the motion of \,$\lambda_k(x)$\, under \,$x$\,. Let us write the monodromy corresponding to \,$D(x)$\, by Theorem~\ref{T:special:special} as 
\,$M(\lambda,x)=\left( \begin{smallmatrix} a(\lambda,x) & b(\lambda,x) \\ c(\lambda,x) & d(\lambda,x) \end{smallmatrix} \right)$\,,%
\footnote{Note that the functions \,$a,b,c,d$\, comprising this monodromy are (for \,$x\neq 1$\,) different from the functions comprising the extended frame, which we also denoted by \,$a,b,c,d$\, 
e.g.~in Proposition~\ref{P:asympfinal:frame} and in Section~\ref{Se:darboux}.}
then the characteristic equation \,$c(\lambda_k(x),x)=0$\, together 
with the differential equation \eqref{eq:mono:monodromy-dgl} for \,$\tfrac{\partial\ }{\partial x}M(\lambda,x)$\, yields a differential equation for \,$\lambda_k$\,:
\begin{equation}
\label{eq:diffeo:xtrans-ode}
\frac{\partial \lambda_k}{\partial x} = -\frac{1}{4\,c'(\lambda_k)}\cdot (\lambda\,\tau+\tau^{-1})\cdot (\mu_k-\mu_k^{-1}) \qmq{with} \tau := \left( \prod_{k\in \Z} \frac{\lambda_{k,0}}{\lambda_k} \right)^{1/2}
\end{equation}
(see the proof of Proposition~\ref{P:jacobitrans:dgl-V1} for the detailed calculation). For given \,$x\in \R$\,, if there is a solution \,$\lambda_k$\, of this differential equation on \,$[0,x]$\, (resp.~on \,$[x,0]$\,)
for all \,$k\in \Z$\,, and \,$D(x) := \{(\lambda_k(x),\mu_k(x))\} \in \Div$\, holds, we call \,$D(x)$\, the \,$x$-translation of \,$D$\,; 
here \,$\mu_k(x)$\, is determined by the condition that \,$D(x)$\, moves smoothly on the spectral curve \,$\Sigma$\,
of \,$D$\,. Note that it is not clear a priori that \,$D(x)$\, exists for any \,$x\neq 0$\, (when one does not know that \,$D$\, stems from a potential \,$(u,u_y)\in \Pot$\,). As long as \,$D(x)$\, is defined,
it depends smoothly on \,$D$\, and on \,$x$\,. 

We will study the \,$x$-translation (and the \,$y$-translation as well) for general tame divisors \,$D$\, in Section~\ref{Se:jacobitrans} by using the Jacobi variety associated to the spectral curve corresponding to \,$D$\,.
At the moment however, we are interested in the \,$x$-translation only of finite type divisors \,$D$\,. If \,$D=\{(\lambda_k,\mu_k)\}\in\Div_{tame}$\,
is given and \,$(\lambda_k,\mu_k)$\, is a double point of the spectral curve \,$\Sigma$\, corresponding to \,$D$\, for some \,$k\in \Z$\,, then the corresponding \,$\lambda_k$\, will not move under \,$x$-translation,
i.e.~we will have \,$(\lambda_k(x),\mu_k(x))=(\lambda_k(0),\mu_k(0))$\, for all \,$x\in \R$\,. If \,$D$\, is of finite type, the system of differential equations \eqref{eq:diffeo:xtrans-ode} contains in fact only
finitely many nontrivial equations; it follows that solutions \,$\lambda_k(x)$\, exist and \,$D(x)=(\lambda_k(x),\mu_k(x))$\, remains in \,$\Div$\, as long as \,$D(x)$\, does not run into a non-tame (special) divisor,
and this is the case at least for small \,$|x|$\,. Note that the \,$x$-translation of a finite type divisor \,$D$\, is again of finite type.

The following fact is a well-known result from the theory of finite type solutions of the sinh-Gordon equation. 

\begin{lem}[Bobenko]
\label{L:diffeo:finite}
For any \,$D\in \Div$\, that is of finite type, and so that the \,$x$-translation \,$D(x)$\, exists for all \,$x\in [0,1]$\, and is non-special, 
there exists one and only one potential \,$(u,u_y) \in \Pot$\, with \,$\Phi((u,u_y))=D$\,.
\end{lem}

\begin{proof}
The potential \,$(u,u_y)$\, corresponding to \,$D$\, has been constructed explicitly in terms of theta functions by \textsc{Bobenko}, see \cite{Bobenko:1991a}, Theorem~4.1,
also see the construction by \textsc{Babich} in \cite{Babich:1991} and \cite{Babich:1992}. 
In the case where \,$D$\, satisfies a condition of reality that is equivalent to the condition that the corresponding potential \,$(u,u_y)$\, is real-valued
(this reality condition already implies that \,$D(x)$\, is non-special for all \,$x$\,), a nice explicit construction of the potential \,$(u,u_y)$\, in terms of
vector-valued Baker-Akhiezer functions on the spectral curve \,$\Sigma$\, corresponding to \,$D$\, was described by \textsc{Knopf} in 
\cite{Knopf:2013}, Proposition~4.34.
\end{proof}

The instrument for the application of Lemma~\ref{L:diffeo:finite} to our situation, where the divisor \,$D$\, is in general of infinite type, is the following lemma:

\begin{lem}
\label{L:diffeo:xnsp-dense}
The set of divisors \,$D\in \Div_{tame}$\, that are of finite type, and such that \,$D(x)$\, exists for all \,$x\in [0,1]$\, and is tame (in particular, non-special),
is dense in \,$\Div_{tame}$\,. 
\end{lem}

\begin{proof}
For \,$N\in \N$\, we let \,$\Div_N$\, be the set of those divisors \,$D=\{(\lambda_k,\mu_k)\} \in \Div_{tame}$\, for which \,$(\lambda_k,\mu_k)$\, is a double point of the spectral curve \,$\Sigma$\, associated to \,$D$\,
for all \,$k\in \Z$\, with \,$|k|> N$\,; clearly the divisors in \,$\Div_N$\, are of finite type. Moreover, we let \,$\Div_{N,xtame}$\, be the set of divisors \,$D\in \Div_N$\, so that \,$D(x)$\, exists and is tame
for all \,$x\in [0,1]$\,. We will show that \,$\Div_{N,xtame}$\, is open and dense in \,$\Div_N$\,.

As discussed above, the only obstacle for the existence of \,$D(x)$\, for a finite type divisor \,$D$\, is that the \,$x$-translation runs into a special divisor. It follows that for \,$D\in \Div_N \setminus \Div_{N,xtame}$\,
there exists \,$x\in (0,1]$\, such that \,$D(x)$\, is non-tame, i.e.~such that there exist \,$k,\ell\in \Z$\, with \,$|k|,|\ell|\leq N$\, and \,$k\neq \ell$\, so that \,$\lambda_k(x)=\lambda_\ell(x)$\, holds; 
this value of \,$x$\, is then the maximal value for which \,$D(x)$\, is defined. This shows that
$$ \Div_N \setminus \Div_{N,xtame} = \bigcup_{\substack{|k|,|\ell|\leq N\\k\neq \ell}} \bigcup_{x\in [0,1]} \Mengegr{D=\{(\lambda_k,\mu_k)\}\in \Div_N}{\lambda_k(x)=\lambda_\ell(x)} $$
holds. The set \,$\Mengegr{D=\{(\lambda_k,\mu_k)\}\in \Div_N}{\lambda_k(x)=\lambda_\ell(x)}$\, is a complex hypersurface in \,$\Div_N$\, for every fixed value of \,$k,\ell,x$\, because \,$D(x)$\, depends smoothly on \,$D$\,, 
and therefore  \,$\bigcup_{x\in [0,1]} \Mengegr{D=\{(\lambda_k,\mu_k)\}\in \Div_N}{\lambda_k(x)=\lambda_\ell(x)}$\, is for every fixed value of \,$k,\ell$\, contained in a real hypersurface in \,$\Div_N$\, because
\,$D(x)$\, also depends smoothly on \,$x$\,. Thus \,$\Div_N \setminus \Div_{N,xtame}$\, is contained in a union of finitely many real hypersurfaces in \,$\Div_N$\,, and therefore \,$\Div_{N,xtame}$\, is open and dense
in \,$\Div_N$\,.

Therefore the set \,$\bigcup_{N\in \N} \Div_{N,xtame}$\, of divisors \,$D\in\Div_{tame}$\, that are of finite type and such that \,$D(x)$\, exists and is tame for all \,$x\in [0,1]$\, is open and dense in the set
\,$\bigcup_{N\in \N} \Div_N$\, of all divisors of finite type in \,$\Div_{tame}$\,. Because the set of finite type divisors in \,$\Div_{tame}$\, is dense in \,$\Div_{tame}$\, by Theorem~\ref{T:finite:finite},
we see that the set of divisors \,$D\in \Div_{tame}$\, that are of finite type, and such that \,$D(x)$\, exists for all \,$x\in [0,1]$\, and is tame,
is dense in \,$\Div_{tame}$\,. 
\end{proof}

To further prepare the proof of Theorem~\ref{T:diffeo:diffeo}, we show that the map \,$\Phi: \Pot\to\Div$\, is smooth,
and calculate the action of its derivative on the symplectic basis \,$(v_k,w_k)$\, of \,$T_{(u,u_y)}\Pot$\, from 
Section~\ref{Se:darboux}:

\begin{lem}
\label{L:diffeo:Phi-smooth}
\begin{enumerate}
\item
The map \,$\Phi:\Pot\to\Div,\; (u,u_y) \mapsto (\lambda_k,\mu_k)_{k\in \Z}$\, is smooth at \,$(u,u_y)\in \Pot_{tame}$\,
in the ``weak'' sense that every component function \,$\lambda_k:\Pot\to\C$\, or \,$\mu_k:\Pot\to\C$\, (with \,$k\in \Z$\,) 
depends smoothly on \,$(u,u_y)$\,.%
\footnote{We will see in the proof of Theorem~\ref{T:diffeo:diffeo} that \,$\Phi$\, is in fact smooth in a ``stronger''
sense, namely as a map into the Banach space \,$\Div$\,.}
\item
Let \,$(u,u_y) \in \Pot_{tame}$\, be given, let \,$F(x,\lambda)$\, be the extended frame corresponding to \,$(u,u_y)$\,
written as in Equation~\eqref{eq:darboux:extended-frame},
let \,$(v_k,w_k)$\, be the symplectic basis of \,$T_{(u,u_y)}\Pot$\, defined in \eqref{eq:darboux:vkwk}
(see also Theorem~\ref{T:darboux:darboux}(1)), let \,$\vartheta_k$\, be the numbers defined in Equation~\eqref{eq:darboux:vartheta},
and put
$$ \wt{\vartheta}_k := \int_0^1 \left( \lambda_k\,a(x,\lambda_k)\,b(x,\lambda_k) + \lambda_k^{-1}\,c(x,\lambda_k)\,d(x,\lambda_k) \right)\cdot e^{u(x)/2}\,\mathrm{d}x \; . $$
For \,$k\in \Z$\, we moreover let \,$e_k=(\delta_{jk})_{j\in \Z}$\, be the \,$\Z$-sequence whose \,$k$-th member is \,$1$\, and
all other members are \,$0$\,; then we define elements of \,$T_{\Phi((u,u_y))}\Div \cong \ell^2_{-1,3} \oplus \ell^2_{0,0}$\,
by \,$e_{\lambda,k} := (e_k,0)$\, and \,$e_{\mu,k} := (0,e_k)$\,. Clearly \,$(e_{\lambda,k},e_{\mu,k})_{k\in \Z}$\, is a basis of
\,$T_{\Phi((u,u_y))}\Div$\,. 

In this setting we have for \,$k\in \Z$\, (where we interpret \,$\Phi'((u,u_y))$\, according to the ``weak'' differentiability
of (1)):
\begin{align*}
\Phi'((u,u_y))v_k & = -i\,\vartheta_k\,\mu_k\,e_{\mu,k} \;, \\
\Phi'((u,u_y))w_k & = -2\,\lambda_k\,e_{\lambda,k} + i\,\mu_k\,\wt{\vartheta}_k\,e_{\mu,k} \; . 
\end{align*}
\end{enumerate}
\end{lem}

\begin{proof}
\emph{For (1).}
We begin by showing that \,$\Phi$\, is smooth in the ``weak'' sense described in the lemma:
The \,$1$-form \,$\alpha_\lambda$\, defined by Equation~\eqref{eq:mono:alphaxy} depends smoothly on \,$(u,u_y)\in \Pot$\,. Therefore the extended frame \,$F_\lambda(x)$\,, characterized by
the initial value problem
$$ \mathrm{d}F_\lambda = \alpha_\lambda\,F_\lambda \qmq{with} F_\lambda(0)=\unity $$
also depends smoothly on \,$(u,u_y)$\,, and hence, the monodromy \,$M(\lambda) = F_\lambda(1)=\left( \begin{smallmatrix} a(\lambda) & b(\lambda) \\ c(\lambda) & d(\lambda) \end{smallmatrix} \right)$\, depends smoothly on \,$(u,u_y)$\,;
here the entries of the monodromy and the entries of the extended frame are related by \,$a(\lambda)=a(1,\lambda)$\, and likewise for the functions \,$b$\,, \,$c$\, and \,$d$\,. 

In the spectral divisor \,$D=\{(\lambda_k,\mu_k)\}$\, corresponding to the monodromy \,$M(\lambda)$\,, the \,$\lambda_k$\, are the zeros of the function \,$c(\lambda)$\,, which depends smoothly on \,$(u,u_y)$\,. Therefore the \,$\lambda$-coordinates of the divisor points depend
smoothly on \,$(u,u_y)$\,; this is true even in the event of zeros of \,$c$\, of higher order by virtue of our construction of \,$\Div$\, as the quotient \,$(\ell^2_{-1,3} \oplus \ell^2_{0,0})/P(\Z)$\,, 
see the end of Section~\ref{Se:asympdiv}. Finally, the \,$\mu$-components of the divisor points are obtained as \,$\mu_k = a(\lambda_k)$\,, where the function \,$a$\, depends smoothly on \,$(u,u_y)$\,.
Thus the \,$\mu$-components also depend smoothly on \,$(u,u_y)$\,. Hence, the map \,$\Phi: \Pot\to\Div$\, is smooth.

\emph{For (2).}
We let \,$(u,u_y)\in \Pot_{tame}$\, be given, and let \,$D:=\Phi((u,u_y)) \in \Div_{tame}$\, be the spectral divisor of \,$(u,u_y)$\,. We equip \,$T_{(u,u_y)}\Pot$\, with the non-degenerate symplectic form \,$\Omega$\,
from Equation~\eqref{eq:darboux:Omega-defined}. Moreover, because \,$D$\, is tame, \,$D$\, is a regular point of \,$\Div$\, and therefore we can consider the tangent space \,$T_D\Div \cong \ell^2_{-1,3} \oplus \ell^2_{0,0}$\,, 
which we equip with the non-degenerate symplectic form
$$ \wt{\Omega}: T_D\Div \times T_D\Div \to \C,\; \bigr(\, (\delta \lambda_k,\delta \mu_k)\,,(\wt{\delta} \lambda_k,\wt{\delta} \mu_k)\,\bigr) \mapsto 
\frac{i}{2}\,\sum_{k\in \Z} \left( \frac{\delta \lambda_k}{\lambda_k} \cdot \frac{\wt{\delta} \mu_k}{\mu_k} - \frac{\wt{\delta} \lambda_k}{\lambda_k} \cdot \frac{\delta \mu_k}{\mu_k} \right)  \;; $$
the sum converges by Lemma~\ref{L:darboux:convergence}. By Theorem \ref{T:darboux:darboux}(2) the derivative \,$\Phi'((u,u_y)): T_{(u,u_y)}\Pot \to T_D\Div$\, is a symplectomorphism from
\,$(T_{(u,u_y)}\Pot,\Omega)$\, to \,$(T_D\Div,\wt{\Omega})$\,. Moreover we have \,$\Phi'((u,u_y))(\delta u,\delta u_y)=(\delta \lambda_k,\delta \mu_k)$\,; from this equation and the definitions of 
\,$\wt{\Omega}$\,, \,$e_{\lambda,k}$\, and \,$e_{\mu,k}$\, it follows that
\begin{align*}
\delta \lambda_k & = -\frac{2}{i}\,\lambda_k\,\mu_k\,\wt{\Omega}\left(e_{\mu,k}, \Phi'((u,u_y))(\delta u,\delta u_y)\right) \\
\qmq{and} \delta \mu_k & = \frac{2}{i}\,\lambda_k\,\mu_k\,\wt{\Omega}\left(e_{\lambda,k}, \Phi'((u,u_y))(\delta u,\delta u_y)\right) 
\end{align*}
holds. 


Now let \,$(\delta u,\delta u_y) \in T_{(u,u_y)}\Pot$\, be given. 
\textsc{Knopf} has shown in \cite{Knopf:2015} (in the proof of Theorem~5.5,
essentially by evaluating an ungauged analogue to our Equation~\eqref{eq:diffeo:delM-asymp:deltawtF-pre})
that the corresponding variations \,$\delta a$\, and \,$\delta c$\, of entries of the extended frame satisfy
$$ \delta a(\lambda_k) = -i\,\mu_k\,\Omega\left(w_k,(\delta u,\delta u_y)\right) \qmq{and}
\delta c(\lambda_k) = -i\,\mu_k\,\Omega\left(v_k,(\delta u,\delta u_y)\right) \;.$$
By an analogous calculation he also showed that the derivatives \,$a'(\lambda)$\, and \,$c'(\lambda)$\, of \,$a$\, resp.~\,$c$\, 
with respect to \,$\lambda$\, satisfy
\begin{equation*}
a'(\lambda_k) = -i \frac{\mu_k}{2\,\lambda_k}\,\wt{\vartheta}_k \qmq{and}
c'(\lambda_k) = -i \frac{\mu_k}{2\,\lambda_k}\,{\vartheta}_k \;.
\end{equation*}

By the implicit function theorem for \,$\lambda_k$\, and differentiation of the equality \,$\mu_k=a(\lambda_k)$\, we have
$$ \delta \lambda_k = -\frac{\delta c(\lambda_k)}{c'(\lambda_k)}
\qmq{and} \delta \mu_k = \delta a(\lambda_k) + a'(\lambda_k) \cdot \delta \lambda_k \; . $$

By combining the preceding equations, and the fact that \,$\Phi'((u,u_y))$\, is a symplectomorphism with respect to 
\,$\Omega$\, and \,$\wt{\Omega}$\,, we obtain
\begin{align*}
& -\frac{2}{i}\,\lambda_k\,\mu_k\,\wt{\Omega}\left(e_{\mu,k}, \Phi'((u,u_y))(\delta u,\delta u_y)\right) \\
=\; & \delta \lambda_k = -\frac{\delta c(\lambda_k)}{c'(\lambda_k)} = -\frac{-i\,\mu_k\,\Omega\left(v_k,(\delta u,\delta u_y)\right)}{-i \frac{\mu_k}{2\,\lambda_k}\,{\vartheta}_k} = -2\,\frac{\lambda_k}{\vartheta_k}\,\Omega\left(v_k,(\delta u,\delta u_y)\right) \\
=\; & -2\,\frac{\lambda_k}{\vartheta_k}\,\wt{\Omega}\left(\Phi'((u,u_y))v_k,\Phi'((u,u_y))(\delta u,\delta u_y)\right) 
\end{align*}
and therefore, because \,$(\delta u,\delta u_y)\in T_{(u,u_y)}\Pot$\, was arbitrary and \,$\Phi'((u,u_y))$\, is a symplectomorphism,
$$ -\frac{2}{i}\,\lambda_k\,\mu_k\,e_{\mu,k} = -2\,\frac{\lambda_k}{\vartheta_k}\,\Phi'((u,u_y))v_k \;, $$
whence 
$$ \Phi'((u,u_y))v_k = -i\,\vartheta_k\,\mu_k\,e_{\mu,k} $$
follows.

Similarly we have
\begin{align*}
& \frac{2}{i}\,\lambda_k\,\mu_k\,\wt{\Omega}\left(e_{\lambda,k}, \Phi'((u,u_y))(\delta u,\delta u_y)\right) \\
=\; & \delta \mu_k = \delta a(\lambda_k) + a'(\lambda_k) \cdot \delta \lambda_k 
= -i\,\mu_k\,\Omega\left(w_k,(\delta u,\delta u_y)\right) + (-i) \frac{\mu_k}{2\,\lambda_k}\,\wt{\vartheta}_k \cdot \delta\lambda_k \\
=\; & -i\,\mu_k\,\wt{\Omega}\left(\Phi'((u,u_y))w_k,\Phi'((u,u_y))(\delta u,\delta u_y)\right) + \mu_k^2\,\wt{\vartheta}_k \,\wt{\Omega}\left(e_{\mu,k}, \Phi'((u,u_y))(\delta u,\delta u_y)\right) \\
\end{align*}
and therefore
$$ \frac{2}{i}\,\lambda_k\,\mu_k\,e_{\lambda,k} = -i\,\mu_k\,\Phi'((u,u_y))w_k + \mu_k^2\,\wt{\vartheta}_k \,e_{\mu,k} \;, $$
whence
$$ \Phi'((u,u_y))w_k = i\,\mu_k\,\wt{\vartheta}_k\,e_{\mu,k} - 2\,\lambda_k\,e_{\lambda,k} $$
follows.
\end{proof}

\newpage

We are now ready for the proof of Theorem~\ref{T:diffeo:diffeo}:

\begin{proof}[Proof of Theorem~\ref{T:diffeo:diffeo}.]
We let \,$(u,u_y)\in \Pot_{tame}$\, be given, and let \,$D:=\Phi((u,u_y)) \in \Div_{tame}$\, be the spectral divisor of \,$(u,u_y)$\,. 
We use the notations from Lemma~\ref{L:diffeo:Phi-smooth} and its proof, in particular we equip \,$T_{(u,u_y)}\Pot$\, and
\,$T_D\Div$\, with the symplectic forms \,$\Omega$\, resp.~\,$\wt{\Omega}$\,.

To prove that \,$\Phi$\, is a local diffeomorphism near \,$(u,u_y)$\,, we apply the inverse function theorem.
We therefore need to show that \,$\Phi:\Pot\to\Div$\, is smooth at \,$(u,u_y)$\, as a map between Banach spaces,
and that the derivative \,$\Phi'((u,u_y)): T_{(u,u_y)}\Pot \to T_D\Div$\, of \,$\Phi$\, 
is an invertible linear operator of Banach spaces.

Because \,$\Phi$\, is smooth in the ``weak'' sense of Lemma~\ref{L:diffeo:Phi-smooth}(1), we have the linear
map \,$\Phi'((u,u_y)): T_{(u,u_y)}\Pot \to T_D\Div$\,. This linear map is a symplectomorphism with respect to
\,$\Omega$\, and \,$\wt{\Omega}$\, by Theorem~\ref{T:darboux:darboux}(2), and is therefore bijective. 
It remains to show that this linear map is an invertible operator of Banach spaces, i.e.~that it is continuous
and that its inverse is also continuous. 

For this we need to show that there exist constants \,$C_1,C_2>0$\, so that 
\begin{align*}
 C_1\cdot \left\|v_k\right\|_{\Pot} & \leq \left\| \Phi'((u,u_y))v_k \right\|_{\Div} \leq C_2 \cdot \left\|v_k\right\|_{\Pot} \\
\qmq{and} C_1\cdot \left\|w_k\right\|_{\Pot} & \leq \left\| \Phi'((u,u_y))w_k \right\|_{\Div} \leq C_2 \cdot \left\|w_k\right\|_{\Pot} 
\end{align*}
holds for all \,$k\in \Z$\,, where \,$(v_k,w_k)$\, is the basis of \,$T_{(u,u_y)}\Pot$\, introduced in \eqref{eq:darboux:vkwk},
see also Theorem~\ref{T:darboux:darboux}(1). It follows from the asymptotic assessment for \,$v_k$\, and \,$w_k$\,
in the proof of Lemma~\ref{L:darboux:basis} that there exist constants \,$C_3,\dotsc,C_6>0$\, so that 
\begin{align*}
\|v_k\|_{\Pot} & = \left.\begin{cases} C_3\, k^2 & \text{for \,$k\geq 0$\,} \\ C_4 & \text{for \,$k<0$\,} \end{cases}\right\} + \ell^2_{-2,0}(k) \\
\qmq{and}
\|w_k\|_{\Pot} & = \left.\begin{cases} C_5\, k & \text{for \,$k\geq 0$\,} \\ C_6\,|k| & \text{for \,$k<0$\,} \end{cases}\right\} + \ell^2_{-1,-1}(k)
\end{align*}
holds. To show that \,$\Phi'((u,u_y))$\, is an invertible operator of Banach spaces it therefore suffices to show that
there exist constants \,$C_7,C_8>0$\, so that 
\begin{align}
\label{eq:diffeo:diffeo:Phi'inv-claim-vk}
\left.\begin{cases} C_7\,k^2 & \text{for \,$k\geq 0$\,} \\ C_7 & \text{for \,$k<0$\,} \end{cases} \right\} 
& \leq \left\| \Phi'((u,u_y))v_k \right\|_{\Div}
\leq \left.\begin{cases} C_8\,k^2 & \text{for \,$k\geq 0$\,} \\ C_8 & \text{for \,$k<0$\,} \end{cases} \right\} \\
\label{eq:diffeo:diffeo:Phi'inv-claim-wk}
\qmq{and}
\left.\begin{cases} C_7\,k & \text{for \,$k\geq 0$\,} \\ C_7\,|k| & \text{for \,$k<0$\,} \end{cases} \right\} 
& \leq \left\| \Phi'((u,u_y))w_k \right\|_{\Div}
\leq \left.\begin{cases} C_8\,k & \text{for \,$k\geq 0$\,} \\ C_8\,|k| & \text{for \,$k<0$\,} \end{cases} \right\} 
\end{align}
holds for all \,$k\in \Z$\, with \,$|k|$\, large. 

For \eqref{eq:diffeo:diffeo:Phi'inv-claim-vk}, we note that \,$\Phi'((u,u_y))v_k=-i\,\vartheta_k\,\mu_k\,e_{\mu,k}$\,
holds by Lemma~\ref{L:diffeo:Phi-smooth}(2). 
The asymptotics for the extended frame \,$F(x,\lambda)$\, in Proposition~\ref{P:asympfinal:frame} show that
\,$\vartheta_k = \vartheta_{k,0}+\ell^2_{-2,0}(k)$\, holds, where
$$ \vartheta_{k,0} := \int_{0}^1 \left( \lambda_{k,0}\,\cos(k\pi\,x)^2 + \sin(k\pi\,x)^2 \right) \,\mathrm{d}x = \frac{1}{2}(\lambda_{k,0}+1) $$
denotes the quantity corresponding to \,$\vartheta_k$\, for the vacuum. Moreover \,$\mu_k = (-1)^k + \ell^2_{0,0}(k)$\,
and \,$\|e_{\mu,k}\|_{\Div}=1$\, holds. Thus we obtain
$$ \left\| \Phi'((u,u_y))v_k\right\|_{\Div} = |\vartheta_k|\cdot|\mu_k|\cdot \|e_{\mu,k}\|_{\Div} 
= \frac{1}{2}\,(\lambda_{k,0}+1) + \ell^2_{-2,0}(k) \;, $$
which implies \eqref{eq:diffeo:diffeo:Phi'inv-claim-vk}. 

For \eqref{eq:diffeo:diffeo:Phi'inv-claim-wk}, we similarly 
note that \,$\Phi'((u,u_y))w_k=-2\,\lambda_k\,e_{\lambda,k} + i\,\mu_k\,\wt{\vartheta}_k\,e_{\mu,k}$\, 
holds by Lemma~\ref{L:diffeo:Phi-smooth}(2). In addition to the information on \,$\mu_k$\, and \,$e_{\mu,k}$\, already given,
we have \,$\lambda_k = \lambda_{k,0} + \ell^2_{-1,3}(k)$\,, 
$$ \|e_{\lambda,k}\|_{\Div} = \begin{cases} k^{-1} & \text{for \,$k\geq 0$\,} \\ |k|^3 & \text{for \,$k<0$\,} \end{cases} \;, $$
and the asymptotics for the extended frame \,$F(x,\lambda)$\, in Proposition~\ref{P:asympfinal:frame} show that
\,$\wt{\vartheta}_k \in \ell^2_{-1,-1}(k)$\, holds, because the quantity \,$\wt{\vartheta}_{k,0}$\, corresponding to \,$\wt{\vartheta}_k$\,
for the vacuum vanishes:
$$ \wt{\vartheta}_{k,0} = \int_0^1 \left( -\lambda_{k,0}^{1/2}\,\cos(k\pi\,x)\,\sin(k\pi\,x) + \lambda_{k,0}^{-1/2}\,\sin(k\pi\,x)\,\cos(k\pi\,x) \right)\,\mathrm{d}x=0 \; . $$
Thus we obtain 
\begin{align*}
\left\| \Phi'((u,u_y))w_k- (-2)\,\lambda_k\,e_{\lambda,k}\right\|_{\Div} 
& = \left\| i\,\mu_k\,\wt{\vartheta}_k\,e_{\mu,k} \right\|_{\Div} \\
& = (1+\ell^2_{0,0}(k)) \cdot \ell^2_{-1,-1}(k)\cdot 1 \;\in\; \ell^2_{-1,-1}(k)
\end{align*}
with
\begin{align*}
\left\| (-2)\,\lambda_k\,e_{\lambda,k} \right\|_{\Div}
& = 2\cdot (\lambda_{k,0}+\ell^2_{-1,3}(k)) \cdot \left.\begin{cases} k^{-1} & \text{for \,$k\geq 0$\,} \\ |k|^3 & \text{for \,$k<0$\,} \end{cases} \right\} \\
& = \left. \begin{cases} C_{9}\,k & \text{for \,$k\geq 0$\,} \\ C_{10}\,|k| & \text{for \,$k<0$\,} \end{cases} \right\} + \ell^2_{0,0}(k) \;, 
\end{align*}
which implies \eqref{eq:diffeo:diffeo:Phi'inv-claim-wk}. 

Therefore \,$\Phi'((u,u_y))$\, is an invertible operator
of Banach spaces. Hence it follows from the inverse function theorem that \,$\Phi$\, is a local diffeomorphism near \,$(u,u_y)$\,. 


Because \,$\Phi|\Pot_{tame}$\, thus is a local diffeomorphism, its image is an open subset of \,$\Div_{tame}$\,. Because all finite type divisors \,$D\in \Div_{tame}$\, for which \,$D(x)$\, exists and is tame 
for all \,$x\in [0,1]$\, are contained in this image by Lemma~\ref{L:diffeo:finite},
and the set of these divisors is dense in \,$\Div_{tame}$\, by Lemma~\ref{L:diffeo:xnsp-dense}, the image of \,$\Phi|\Pot_{tame}$\, is also dense in \,$\Div_{tame}$\,. 

It remains to show that \,$\Phi|\Pot_{tame}$\, is injective. For this purpose let \,$(u^{[1]},u_y^{[1]})$\,, \,$(u^{[2]},u_y^{[2]}) \in \Pot_{tame}$\, be given with \,$\Phi((u^{[1]},u_y^{[1]}))=\Phi((u^{[2]},u_y^{[2]})) =: D \in \Div_{tame}$\,. 
Because \,$\Phi$\, is a diffeomorphism near \,$(u^{[\nu]},u_y^{[\nu]})$\, (for \,$\nu\in\{1,2\}$\,), there exist neighborhoods \,$U^{[\nu]}$\, of \,$(u^{[\nu]},u_y^{[\nu]})$\, in \,$\Pot_{tame}$\, and \,$V^{[\nu]}$\, of \,$D$\,
in \,$\Div_{tame}$\, so that \,$\Phi|U^{[\nu]}: U^{[\nu]}\to V^{[\nu]}$\, is a diffeomorphism.

By Lemma~\ref{L:diffeo:xnsp-dense} there exists a sequence \,$(D_n)_{n\geq 1}$\, of divisors of finite type, such that \,$D_n(x)$\, exists and is tame for all \,$x\in[0,1]$\,, which converges to \,$D$\,. 
Without loss of generality, we may suppose 
\,$D_n \in V^{[1]}\cap V^{[2]}$\, for all \,$n\geq 1$\,, and then \,$(u^{[\nu,n]},u^{[\nu,n]}_y) := (\Phi|U^{[\nu]})^{-1}(D_n)$\, is a sequence in \,$U^{[\nu]}$\, which converges to 
\,$(\Phi|U^{[\nu]})^{-1}(D_n)=(u^{[\nu]},u^{[\nu]}_y)$\,. On the other hand, \,$D_n$\, has only one pre-image under \,$\Phi$\, by the uniqueness statement in Lemma~\ref{L:diffeo:finite}. 
Therefore \,$(u^{[1,n]},u^{[1,n]}_y) = (u^{[2,n]},u^{[2,n]}_y)$\, holds for all \,$n\geq 1$\,, whence \,$(u^{[1]},u^{[1]}_y) = (u^{[2]},u^{[2]}_y)$\, follows by taking the limit \,$n\to\infty$\,. 
\end{proof}

\begin{cor}
\label{C:diffeo:dense}
\begin{enumerate}
\item The set of potentials of finite type in \,$\Pot_{tame}$\, is dense in \,$\Pot_{tame}$\,.
\item The set of divisors \,$D\in\Div_{tame}$\, such that \,$D(x)$\, exists and is tame for all \,$x\in [0,1]$\, is open and dense in \,$\Div_{tame}$\,. 
\end{enumerate}
\end{cor}

\begin{proof}
\emph{For (1).}
\label{not:diffeo:Div-fin}
The set \,$\Div_{fin}$\, of divisors of finite type in \,$\Div_{tame}$\, is dense in \,$\Div_{tame}$\, by Theorem~\ref{T:finite:finite}; because \,$\Phi[\Pot_{tame}]$\, is open in \,$\Div_{tame}$\, by Theorem~\ref{T:diffeo:diffeo},
it follows that \,$\Div_{fin} \cap \Phi[\Pot_{tame}]$\, is dense in \,$\Phi[\Pot_{tame}]$\,. Because \,$\Phi: \Pot_{tame}\to\Phi[\Pot_{tame}]$\, is a diffeomorphism again by Theorem~\ref{T:diffeo:diffeo}, it follows
that \,$\Phi^{-1}[\Div_{fin} \cap \Phi[\Pot_{tame}]]$\, is dense in \,$\Pot_{tame}$\,; this is the set of finite type potentials in \,$\Pot_{tame}$\,.

\emph{For (2).}
The set of divisors \,$D\in\Div_{tame}$\, such that \,$D(x)$\, exists and is tame for all \,$x\in [0,1]$\, is open in \,$\Div_{tame}$\, because \,$D(x)$\, depends continuously on \,$D$\, and on \,$x$\,.
This set is dense in \,$\Div_{tame}$\, because already the set of finite type divisors such that \,$D(x)$\, exists and is tame for all \,$x\in [0,1]$\, is dense in \,$\Div_{tame}$\, by 
Lemma~\ref{L:diffeo:xnsp-dense}.
\end{proof}

\part{The Jacobi variety of the spectral curve}

\section{Estimate of certain integrals}
\label{Se:jacobiprep}


Our next big objective is the construction of Jacobi coordinates for the spectral curve \,$\Sigma$\, corresponding to some potential \,$(u,u_y)\in \Pot$\,. Similarly as for the construction of the Jacobi
variety for compact Riemann surfaces, our construction of the Jacobi variety for \,$\Sigma$\, in Section~\ref{Se:jacobi} will involve path integrals on the spectral curve \,$\Sigma$\, which are of the form
$$ \int_{P_0}^{P_1} \frac{\Phi(\lambda)}{\mu-\mu^{-1}}\,\mathrm{d}\lambda \;, $$
where \,$\Phi$\, is a holomorphic function in \,$\lambda$\,. 

To make the construction tractable, we will need to do two things: In the present section, we will estimate \,$\int_{P_0}^{P_1} \tfrac{1}{\mu-\mu^{-1}}\,\mathrm{d}\lambda$\, 
and \,$\int_{P_0}^{P_1} \tfrac{1}{|\mu-\mu^{-1}|}\,\mathrm{d}\lambda$\, (Proposition~\ref{P:jacobiprep:Deltaint-neu}). In the former integral, the integrand is holomorphic, and therefore the value of the integral depends only on the 
homology type of the path of integration. However the second integrand is not holomorphic, and thus the value of the second integral depends on the specific choice of the path of integration.
We will choose a special class of \emph{admissible} paths of integration, which permit to connect each pair of points contained in an excluded domain \,$\wh{U}_{k,\delta} \subset \Sigma$\,, 
but which are simple enough (they are composed of lifts of straight lines and circular arcs in the \,$\lambda$-plane) so that the integral in question can be estimated relatively explicitly.

The other preparation needed for the construction of the Jacobi coordinates on \,$\Sigma$\, is the study of the asymptotic behavior of holomorphic \,$1$-forms on \,$\Sigma$\, which are square-integrable.
This is the topic of Section~\ref{Se:holo1forms}.

\medskip

We suppose that \,$\Sigma$\, is the spectral curve associated to some potential \,$(u,u_y)\in \Pot$\, (or defined by a holomorphic function \,$\Delta:\C^*\to\C$\, with \,$\Delta-\Delta_0 \in \As(\C^*,\ell^2_{0,0},1)$\,).
As before, we have the holomorphic functions \,$\lambda,\mu: \Sigma\to\C^*$\,, and we denote the zeros of \,$\Delta^2-4$\, (corresponding to the branch points resp.~singularities of \,$\Sigma$\,)
by \,$\vkap_{k,\nu}$\, as in Proposition~\ref{P:excl:basic}(1);
we have \,$\vkap_{k,\nu}-\lambda_{k,0} \in \ell^2_{-1,3}(k)$\, by Proposition~\ref{P:asympfinal:vkap}.

Because \,$\vkap_{k,\nu}$\, is a zero of \,$\Delta^2-4$\, for \,$k\in \Z$\, and \,$\nu\in \{1,2\}$\,, 
there exists one and only one point in \,$\Sigma$\, above \,$\vkap_{k,\nu}$\, (which is a fixed point of the hyperelliptic involution of \,$\Sigma$\,). In the sequel 
we will take the liberty of denoting this point of \,$\Sigma$\, also by \,$\vkap_{k,\nu}$\,. 

To prepare for the evaluation of the mentioned integrals, we begin 
by showing that the function \,$\tfrac{1}{\mu-\mu^{-1}}$\, behaves on \,$\wh{U}_{k,\delta}$\, 
essentially like \,$\tfrac{1}{\sqrt{(\lambda-\vkap_{k,1})\cdot(\lambda-\vkap_{k,2})}}$\,, provided that \,$|k|$\, is so large
that \,$\vkap_{k,1}$\, and \,$\vkap_{k,2}$\, are the only zeros of \,$\Delta^2-4$\, (branch points or singularities of \,$\Sigma$\,) 
that occur in \,$\wh{U}_{k,\delta}$\,. 
This statement is made precise by the following lemma:

\begin{lem}
\label{L:jacobiprep:Psi-intro}
Let \,$k\in \Z$\,. There exists \,$N\in \N$\, so that for all \,$k\in \Z$\, with \,$|k|> N$\, 
there exists a holomorphic function \,$\Psi_k$\, on \,$\wh{U}_{k,\delta}$\, with
\begin{equation}
\label{eq:jacobiprep:Psi:psidef}
\Psi_k(\lambda,\mu)^2 = (\lambda-\vkap_{k,1})\cdot (\lambda-\vkap_{k,2}) \;.
\end{equation}
We have \,$\Psi_k \circ \sigma = -\Psi_k$\,, where \,$\sigma$\, is the hyperelliptic involution of \,$\Sigma$\,. 
\,$\Psi_k$\, has the following asymptotic property, which also fixes the sign of \,$\Psi_k$\,:

There exists \,$C>0$\, and a sequence \,$a_k \in \ell^2_{-1,3}(|k|> N)$\, such that for all \,$k\in \Z$\, with \,$|k|> N$\, 
and all \,$(\lambda,\mu) \in \wh{U}_{k,\delta}\setminus\{\vkap_{k,1},\vkap_{k,2}\}$\, we have
if \,$k>0$\, 
$$ \left| \frac{1}{\mu-\mu^{-1}} - 4i(-1)^k\,\lambda_{k,0}^{1/2}\,\frac{1}{\Psi_k(\lambda,\mu)} \right|
\leq \frac{C \cdot \left( |\lambda-\vkap_{k,1}| + |\lambda-\vkap_{k,2}| \right) + a_k}{|\Psi_k(\lambda,\mu)|} \;, $$
and if \,$k<0$\, 
$$ \left| \frac{1}{\mu-\mu^{-1}} - 4i(-1)^k\,\lambda_{k,0}^{3/2}\,\frac{1}{\Psi_k(\lambda,\mu)} \right|
\leq \frac{C \cdot \left( |\lambda-\vkap_{k,1}| + |\lambda-\vkap_{k,2}| \right) + a_k}{|\Psi_k(\lambda,\mu)|} \; . $$
\end{lem}

\begin{proof}
By Proposition~\ref{P:excl:basic}(1) there exists \,$N\in \N$\, so that for every \,$k$\, with \,$|k|> N$\,, \,$\vkap_{k,1}$\,
and \,$\vkap_{k,2}$\, are all the zeros of \,$\Delta^2-4$\, in \,$U_{k,\delta}$\,, and therefore
\,$\vkap_{k,1}$\, and \,$\vkap_{k,2}$\, are all the branch points of \,$\wh{U}_{k,\delta}$\,. For \,$|k|>N$\, it is therefore clear 
that there exists a holomorphic function \,$\Psi_k$\, on \,$\wh{U}_{k,\delta}$\, so that
Equation~\eqref{eq:jacobiprep:Psi:psidef} holds; we then also have \,$\Psi_k \circ \sigma = -\Psi_k$\,.
Equation~\eqref{eq:jacobiprep:Psi:psidef} determines \,$\Psi_k$\, only up to sign on each of the (one or two) connected components of \,$\wh{U}_{k,\delta} \setminus \{\vkap_{k,1},\vkap_{k,2}\}$\,.
The sign is fixed by the stated asymptotic property for \,$\Psi_k$\,, which we now prove: 

We let \,$\rho=1$\, if \,$k>0$\,, \,$\rho=3$\, if \,$k<0$\,. By Proposition~\ref{P:asympfinal:Delta2}(3) there exists \,$C_1>0$\, 
such that we have for \,$\lambda\in U_{k,\delta}$\,
\begin{gather}
\left| \lambda_{k,0}^\rho \cdot \frac{\Delta(\lambda)^2-4}{(\lambda-\vkap_{k,1})\cdot(\lambda-\vkap_{k,2})} - \left( -\frac{1}{16} \right) \right| \notag \\
\label{eq:transdgl:Delta2:asymp:Delta}
\leq C_1 \cdot \left\{ \begin{matrix} k^{-1} & \text{if \,$k>0$\,} \\ k^3 & \text{if \,$k<0$\,} \end{matrix} \right\}  \cdot \left( |\lambda-\vkap_{k,1}| + |\lambda-\vkap_{k,2}| \right) + \ell^{2}_{0,0}(k) \; .
\end{gather}
Because of the equation \,$\mu-\mu^{-1}=\sqrt{\Delta(\lambda)^2-4}$\, and the fact that \,$z\mapsto z^{-1/2}$\, is Lipschitz continuous near \,$z=-\tfrac{1}{16}$\,, there exists \,$C_2>0$\, and \,$\eps \in \{-1,1\}$\, with 
$$ \left| \lambda_{k,0}^{-\rho/2} \cdot \frac{\Psi_k(\lambda,\mu)}{\mu-\mu^{-1}} - 4i\,\eps \right|
\leq C_2 \cdot \left\{ \begin{matrix} k^{-1} & \text{if \,$k>0$\,} \\ k^3 & \text{if \,$k<0$\,} \end{matrix} \right\}  \cdot \left( |\lambda-\vkap_{k,1}| + |\lambda-\vkap_{k,2}| \right) + \ell^{2}_{0,0}(k) \; , $$
whence
\begin{equation*}
\left| \frac{1}{\sqrt{\mu-\mu^{-1}}} - 4i\,\eps\,\lambda_{k,0}^{\rho/2}\,\frac{1}{\Psi_k(\lambda,\mu)} \right|
\leq \frac{C_3 \cdot \left( |\lambda-\vkap_{k,1}| + |\lambda-\vkap_{k,2}| \right) + \ell^{2}_{-1,3}(k)}{|\Psi_k(\lambda,\mu)|}
\end{equation*}
follows. By choosing the proper sign for \,$\Psi_k$\,, we can arrange \,$\eps=1$\,, and therefore obtain the claimed asymptotic property.
\end{proof}

For \,$k\in \Z$\, we let \,$\vkap_{k,*}\in \C^*$\, be the midpoint between \,$\vkap_{k,1}$\, and \,$\vkap_{k,2}$\,, i.e.
\begin{equation}
\label{eq:jacobiprep:vkap-k-star}
\vkap_{k,*} := \frac12(\vkap_{k,1}+\vkap_{k,2}) \; .
\end{equation}
Because of \,$\vkap_{k,\nu}-\lambda_{k,0} \in \ell^2_{-1,3}(k)$\, we have \,$\vkap_{k,*}-\lambda_{k,0}\in \ell^2_{-1,3}(k)$\,, 
and for \,$|k|$\, large, we have \,$\vkap_{k,*} \in U_{k,\delta}$\,.

If \,$k\in \Z$\, is such that \,$\vkap_{k,1}$\, and \,$\vkap_{k,2}$\, are the only zeros of \,$\Delta^2-4$\, in \,$\wh{U}_{k,\delta}$\,
(i.e.~\,$|k|>N$\, in the notations of Lemma~\ref{L:jacobiprep:Psi-intro}), then we put 
\begin{equation}
\label{eq:jacobiprep:Ukd'}
\wh{U}_{k,\delta}' := \begin{cases} \wh{U}_{k,\delta} & \text{if \,$\vkap_{k,1}\neq \vkap_{k,2}$\,} \\ \wh{U}_{k,\delta} \setminus \{\vkap_{k,*}\} & \text{if \,$\vkap_{k,1}=\vkap_{k,2}$\,} \end{cases}
\qmq{and}
{U}_{k,\delta}' := \begin{cases} {U}_{k,\delta} & \text{if \,$\vkap_{k,1}\neq \vkap_{k,2}$\,} \\ {U}_{k,\delta} \setminus \{\vkap_{k,*}\} & \text{if \,$\vkap_{k,1}= \vkap_{k,2}$\,} \end{cases} \; .
\end{equation}
Note that \,$U_{k,\delta}'$\, is connected in any case; whereas \,$\wh{U}_{k,\delta}'$\, is connected only if
\,$\vkap_{k,1}\neq \vkap_{k,2}$\,, for \,$\vkap_{k,1}=\vkap_{k,2}$\, it has two connected components.

In the sequel, we will estimate the integrals of \,$\tfrac{1}{\Psi_k}$\, and of \,$\tfrac{1}{|\Psi_k|}$\, along paths in the excluded domain; because the second integrand \,$\tfrac{1}{|\Psi_k|}$\, is not holomorphic,
we need to specify in detail which paths we permit for estimating the corresponding integral, as was explained above. The admissible paths introduced in the following definition are those paths we permit in this context:

\begin{Def}
\label{D:jacobiprep:admissible}
Let \,$N\in \N$\, be as in Lemma~\ref{L:jacobiprep:Psi-intro} and \,$k\in \Z$\, with \,$|k|>N$\,.
We call a path \,$\gamma$\, in \,$\wh{U}_{k,\delta}'$\, \emph{admissible}, if the following properties hold, where we denote by \,$(\lambda_k^o,\mu_k^o)$\, and \,$(\lambda_k,\mu_k)$\, the endpoints of \,$\gamma$\,:
\begin{enumerate}
\item \,$\gamma$\, is composed of the lifts to \,$\Sigma$\, of at most three each of the following two types of path segments in \,$U_{k,\delta}'$\,: 
\begin{enumerate}
\item
A segment of a line through \,$\vkap_{k,*}$\, and contained in \,$U_{k,\delta}'$\,, traversed in either direction,
\item
A segment of a circle with center \,$\vkap_{k,*}$\, and contained in \,$U_{k,\delta}'$\,, where the circle may be parameterized in either direction, and the angular length is at most \,$4\pi$\, (i.e.~at most two full revolutions of the circle).
\end{enumerate}
\item For all \,$(\lambda,\mu)$\, on the trace of \,$\gamma$\, and for \,$\nu \in \{1,2\}$\, we have
\begin{equation}
\label{eq:jacobiprep:admissible:estimate}
|\lambda-\vkap_{k,\nu}| \leq \max \bigr\{ \,2\,|\vkap_{k,1}-\vkap_{k,2}|\,,\, |\lambda_k^o-\vkap_{k,*}|\,,\, |\lambda_k-\vkap_{k,*}|\,\bigr\} \; .
\end{equation}
\item The trace of \,$\gamma$\, meets the points \,$\vkap_{k,1}$\, and \,$\vkap_{k,2}$\, at most in its endpoints.
\end{enumerate}
\end{Def}

\begin{prop}
\label{P:jacobiprep:admissible}
Let \,$N\in \N$\, be as in Lemma~\ref{L:jacobiprep:Psi-intro} and \,$k\in \Z$\, with \,$|k|>N$\,.
For every \,$(\lambda_k^o,\mu_k^o),(\lambda_k,\mu_k)\in \wh{U}_{k,\delta}'$\, (from the same connected component of \,$\wh{U}_{k,\delta}'$\, if \,$\vkap_{k,1}=\vkap_{k,2}$\,) there exists an admissible path \,$\gamma$\, 
in \,$\wh{U}_{k,\delta}'$\, which connects \,$(\lambda_k^o,\mu_k^o)$\, to \,$(\lambda_k,\mu_k)$\,.
\end{prop}

\begin{proof}
We distinguish cases depending on whether \,$|\lambda_k-\vkap_{k,*}|$\, or \,$|\lambda_k^o-\vkap_{k,*}|$\, is larger or smaller than \,$|\vkap_{k,1}-\vkap_{k,2}|$\,. 
In the sequel, ``line'' or ``ray'' always means a line or ray contained in \,$U_{k,\delta}$\, starting in \,$\vkap_{k,*}$\, and passing through another point, 
and ``circle (segment)'' always means a circle (segment) contained in \,$U_{k,\delta}$\, with center \,$\vkap_{k,*}$\,. 

If \,$|\lambda_k^o-\vkap_{k,*}|, |\lambda_k-\vkap_{k,*}| \leq |\vkap_{k,1}-\vkap_{k,2}|$\, holds, we construct \,$\gamma$\, as composition of lifts in \,$\wh{U}_{k,\delta}'$\, of the following segments in \,$U_{k,\delta}'$\,:
\begin{itemize}
\item If the ray through \,$\lambda_k^o$\, passes through either \,$\vkap_{k,1}$\, or \,$\vkap_{k,2}$\,, a short circle segment so that the ray through its endpoint passes through neither
\,$\vkap_{k,1}$\, nor \,$\vkap_{k,2}$\,; otherwise nothing.
\item The segment of the line through the previous endpoint from that endpoint to a point \,$\lambda\in U_{k,\delta}$\, with \,$|\lambda-\vkap_{k,*}| = |\vkap_{k,1}-\vkap_{k,2}|$\,.
\item A circle segment from the previous endpoint, which contains a full circle if and only if \,$(\lambda_k^o,\mu_k^o)$\, and \,$(\lambda_k,\mu_k)$\, are on different ``leaves'' of \,$\wh{U}_{k,\delta}'$\,,
and whose endpoint is chosen in the following way: If the ray through \,$\lambda_k$\, does not pass through either \,$\vkap_{k,1}$\, or \,$\vkap_{k,2}$\,, then the endpoint is such that the ray through the
endpoint passes through \,$\lambda_k$\,; otherwise the arc is extended beyond that point by a short length so that the ray through its endpoint does not pass through either \,$\vkap_{k,1}$\, or \,$\vkap_{k,2}$\,.
\item A line segment connecting the previous endpoint to a point \,$\lambda'\in U_{k,\delta}$\, with \,$|\lambda'-\vkap_{k,*}| = |\lambda_k-\vkap_{k,*}|$\,.
\item If the lift of the previous endpoint is not \,$(\lambda_k,\mu_k)$\,, then a short circle segment connecting the previous endpoint to \,$(\lambda_k,\mu_k)$\,. 
\end{itemize}

If either \,$|\lambda_k^o-\vkap_{k,*}| > |\vkap_{k,1}-\vkap_{k,2}|$\, or \,$|\lambda_k-\vkap_{k,*}| > |\vkap_{k,1}-\vkap_{k,2}|$\, holds, we may suppose without loss of generality that 
\,$|\lambda_k^o-\vkap_{k,*}| \leq |\lambda_k-\vkap_{k,*}|$\, and therefore 
\,$|\lambda_k-\vkap_{k,*}| >  |\vkap_{k,1}-\vkap_{k,2}|$\, holds. We make no assumptions on the relation between \,$|\lambda_k^o-\vkap_{k,*}|$\, and \,$|\vkap_{k,1}-\vkap_{k,2}|$\,. In this case, we 
construct \,$\gamma$\, as composition of lifts in \,$\wh{U}_{k,\delta}'$\, of the following segments in \,$U_{k,\delta}'$\,:
\begin{itemize}
\item If the segment of the ray through \,$\lambda_k^o$\, from \,$\lambda_k^o$\, to a point \,$\lambda\in U_{k,\delta}$\, with \,$|\lambda-\vkap_{k,*}|=|\lambda_k-\vkap_{k,*}|$\, would pass through either \,$\vkap_{k,1}$\, or \,$\vkap_{k,2}$\,,
a short circle segment  so that the ray through its endpoint passes through neither \,$\vkap_{k,1}$\, nor \,$\vkap_{k,2}$\,; otherwise nothing.
\item A line segment from the previous endpoint to a point \,$\lambda$\, with \,$|\lambda-\vkap_{k,*}|=|\lambda_k-\vkap_{k,*}|$\,.
\item A circle segment from the previous endpoint, which contains a full circle if and only if \,$(\lambda_k^o,\mu_k^o)$\, and \,$(\lambda_k,\mu_k)$\, are on different ``leaves'' of \,$\wh{U}_{k,\delta}'$\,,
and which ends in \,$(\lambda_k,\mu_k)$\,. 
\end{itemize}
\end{proof}

We now estimate the path integral along admissible paths for the integrands \,$\tfrac{1}{\Psi_k}$\,, \,$\tfrac{1}{|\Psi_k|}$\, and \,$\tfrac{|\lambda-\xi_k|}{|\Psi_k|}$\,, and also the area integral
of \,$\left(\tfrac{|\lambda-\xi_k|}{|\Psi_k|}\right)^2$\, over an excluded domain. Note that while the holomorphic function \,$\Psi_k$\, is well-defined only on \,$\wh{U}_{k,\delta}$\, for \,$k$\,
sufficiently large (see Lemma~\ref{L:jacobiprep:Psi-intro}), the continuous function \,$|\Psi_k|$\, is well-defined on all of \,$\Sigma$\, and for all \,$k\in \Z$\,, namely by
$$ |\Psi_k(\lambda)| := \sqrt{|\lambda-\vkap_{k,1}|\cdot|\lambda-\vkap_{k,2}|} \;, $$
where \,$\sqrt{\;\;}$\, denotes the real square root function, moreover \,$|\Psi_k|$\, depends only on \,$\lambda\in\C^*$\,, and not on \,$\mu$\,. 

\begin{lem}
\label{L:jacobiprep:Psi-neu}
Let \,$N\in \N$\, be as in Lemma~\ref{L:jacobiprep:Psi-intro}, \,$k\in \Z$\,, and 
\,$(\lambda_k^o,\mu_k^o), (\lambda_k,\mu_k) \in \wh{U}_{k,\delta}'$\,. 
If \,$\vkap_{k,1}\neq \vkap_{k,2}$\, holds, we let \,$\xi_k \in U_{k,\delta}$\, be given, whereas for \,$\vkap_{k,1}=\vkap_{k,2}$\,
we put \,$\xi_k := \vkap_{k,*} \in U_{k,\delta}$\, and require that 
\,$(\lambda_k^o,\mu_k^o)$\, and \,$(\lambda_k,\mu_k)$\, are in the same connected component of \,$\wh{U}_{k,\delta}'$\,.
\begin{enumerate}


\item For \,$|k|>N$\, we have
\begin{align*}
\int_{(\lambda_k^o,\mu_k^o)} ^{(\lambda_k,\mu_k)} \frac{1}{\Psi_k}\,\mathrm{d}\lambda 
& = \ln\left( \frac{\lambda_k-\vkap_{k,*}+\Psi_k(\lambda_k,\mu_k)}{\lambda_k^o-\vkap_{k,*}+\Psi_k(\lambda_k^o,\mu_k^o)} \right) \; ,
\end{align*}
where we integrate along any path that is entirely contained in \,$\widehat{U}_{k,\delta}'$\,. 
Here \,$\ln(z)$\, is the branch of the complex logarithm function with \,$\ln(1)=2\pi i m$\,, 
where \,$m\in \Z$\, is the winding number of the path of integration around the pair of branch points \,$\vkap_{k,1}$\,, \,$\vkap_{k,2}$\,.
\item
There exists a constant \,$C>0$\, (depending neither on \,$k$\, nor on \,$\vkap_{k,\nu}$\,) so that for \,$|k|>N$\, we have
\begin{equation}
\label{eq:jacobiprep:Psi-neu:int-estimate-nouveau}
\int_{(\lambda_k^o,\mu_k^o)}^{(\lambda_k,\mu_k)} \frac{1}{|\Psi_k|}\,|\mathrm{d}\lambda|
\leq C \cdot \left| \int_{(\lambda_k^o,\mu_k^o)}^{(\lambda_k,\mu_k)} \frac{1}{\Psi_k}\,\mathrm{d}\lambda \right| \;,
\end{equation}
where we integrate along an admissible path (Definition~\ref{D:jacobiprep:admissible}).

\item
There exist constants \,$C_1,C_2>0$\, (depending neither on \,$k$\, nor on \,$\vkap_{k,\nu}$\, or \,$\xi_k$\,) so that we have for any \,$k\in \Z$\,
$$\int_{\vkap_{k,1}}^{(\lambda_k,\mu_k)} \frac{|\lambda-\xi_k|}{|\Psi_k|}\,|\mathrm{d}\lambda|
\leq \left( C_1\,\left| \frac{\xi_k-\vkap_{k,*}}{\vkap_{k,1}-\vkap_{k,2}} \right| + C_2 \right) \cdot |\lambda_k-\vkap_{k,1}| \;. $$
Here we integrate along the lift of the straight line from \,$\vkap_{k,1}$\, to \,$\lambda_k$\,, and the expression
\,$\left| \frac{\xi_k-\vkap_{k,*}}{\vkap_{k,1}-\vkap_{k,2}} \right|$\, is to be read as \,$0$\, in the case \,$\vkap_{k,1}=\vkap_{k,2}$\,. 

\item
There exists \,$C>0$\, (depending neither on \,$k$\, nor on \,$\vkap_{k,\nu}$\, or \,$\xi_k$\,) so that 
we have for any \,$k\in \Z$\,
\begin{equation*}
\int_{U_{k,\delta}} \frac{|\lambda-\xi_k|^2}{|\lambda-\vkap_{k,1}|\cdot |\lambda-\vkap_{k,2}|}\,\mathrm{d}^2\lambda \leq 
C\cdot \left( |\xi_k-\vkap_{k,*}|^2 + \frac{|\xi_k-\vkap_{k,*}|^4}{|\vkap_{k,1}-\vkap_{k,2}|^2} \right) + \ell^\infty_{-2,6}(k) \; ,
\end{equation*}
where \,$d^2\lambda$\, denotes the 2-dimensional Lebesgue measure on \,$\C$\,. 
In the case \,$\vkap_{k,1}=\vkap_{k,2}$\,, the expression \,$\frac{|\xi_k-\vkap_{k,*}|^4}{|\vkap_{k,1}-\vkap_{k,2}|^2}$\, is to be read as \,$0$\,. 
\end{enumerate}
\end{lem}

\begin{proof}[Proof of Lemma~\ref{L:jacobiprep:Psi-neu}.]
%
\emph{For (1).}
It is easily checked that \,$\ln(\lambda-\vkap_{k,*}+\Psi_k(\lambda,\mu))$\, is an anti-derivative of \,$\tfrac{1}{\Psi_k(\lambda,\mu)}$\,, and the claimed formula follows.


\emph{For (2).}
We first show that the inequality to be shown is invariant under change of the \,$\vkap_{k,\nu}$\,. For this purpose, we at first suppose that \,$\vkap_{k,1}\neq\vkap_{k,2}$\, holds,
let \,$\wh{U}$\, be the hyperelliptic surface above \,$\C$\, with branch points above \,$\vkap_1 := 1$\, and \,$\vkap_2 := -1$\,,
and let \,$\Psi(\lambda,\mu) := \sqrt{\lambda^2-1}$\, be the corresponding holomorphic function on \,$\wh{U}$\, with zeros in the branch points. We then have
\begin{align*}
\Psi_k(\lambda,\mu) & = \sqrt{(\lambda-\vkap_{k,1})\cdot (\lambda-\vkap_{k,2})} = \sqrt{(\lambda-\vkap_{k,*})^2-\bigr(\tfrac12(\vkap_{k,1}-\vkap_{k,2})\bigr)^2} \\
& = \frac12(\vkap_{k,1}-\vkap_{k,2}) \cdot \Psi\left( \frac{\lambda-\vkap_{k,*}}{\tfrac12(\vkap_{k,1}-\vkap_{k,2})} , \mu\right)
\end{align*}
and therefore, if we choose points \,$(\lambda^o,\mu^o), (\lambda,\mu) \in \wh{U}$\, with \,$\lambda^o = \tfrac{\lambda_k^o-\vkap_{k,*}}{\tfrac12(\vkap_{k,1}-\vkap_{k,2})}$\, 
resp.~\,$\lambda = \tfrac{\lambda_k-\vkap_{k,*}}{\tfrac12(\vkap_{k,1}-\vkap_{k,2})}$\,, we have
$$ \int_{(\lambda_k^o,\mu_k^o)}^{(\lambda_k,\mu_k)} \frac{1}{\Psi_k}\,\mathrm{d}\lambda' = \int_{(\lambda^o,\mu^o)}^{(\lambda,\mu)} \frac{1}{\Psi}\,\mathrm{d}\lambda' $$
and similarly
$$ \int_{(\lambda_k^o,\mu_k^o)}^{(\lambda_k,\mu_k)} \frac{1}{|\Psi_k|}\,|\mathrm{d}\lambda'| = \int_{(\lambda^o,\mu^o)}^{(\lambda,\mu)} \frac{1}{|\Psi|}\,|\mathrm{d}\lambda'| \;. $$
Therefore we see that the inequality~\eqref{eq:jacobiprep:Psi-neu:int-estimate-nouveau} is implied by the claim that 
there exists \,$C>0$\, so that we have for all \,$(\lambda^o,\mu^o),(\lambda,\mu)\in\wh{U}$\,
\begin{equation}
\label{eq:jacobiprep:Psi-neu:int-claim}
\int_{(\lambda^o,\mu^o)}^{(\lambda,\mu)} \frac{1}{|\Psi|}\,|\mathrm{d}\lambda'|
\leq C \cdot \left| \int_{(\lambda^o,\mu^o)}^{(\lambda,\mu)} \frac{1}{\Psi}\,\mathrm{d}\lambda' \right| \;. 
\end{equation}
For reasons of continuity this is true even in the case \,$\vkap_{k,1}=\vkap_{k,2}$\,.

For the proof of \eqref{eq:jacobiprep:Psi-neu:int-claim}, we note that \,$\vkap_{1}-\vkap_{2}=2$\, and \,$\vkap_* := \tfrac12(\vkap_{1}+\vkap_2) = 0$\, holds. 
We first consider the case where \,$|\lambda^o|,|\lambda|\leq 2$\, holds. Then the admissible path of integration from \,$(\lambda^o,\mu^o)$\, to \,$(\lambda,\mu)$\, runs entirely in
the pre-image of \,$\Mengegr{\lambda'\in \C}{|\lambda'|\leq 2}$\,. We note that 
$$ \int_{(\lambda^o,\mu^o)}^{(\lambda,\mu)} \frac{1}{\Psi}\,\mathrm{d}\lambda' = \ln\left( \frac{\lambda-\vkap_{*}+\Psi(\lambda,\mu)}{\lambda^o-\vkap_{*}+\Psi(\lambda^o,\mu^o)} \right) $$
equals \,$0$\, if and only if we choose the branch of the complex logarithm with \,$\ln(1)=0$\, (corresponding to integration paths that do not wind around the branch points \,$\vkap_\nu$\,), and moreover
$$ \lambda-\vkap_{*}+\Psi(\lambda,\mu) = \lambda^o-\vkap_{*}+\Psi(\lambda^o,\mu^o) $$
and therefore \,$(\lambda,\mu)=(\lambda^o,\mu^o)$\, holds. 

On the other hand, the function \,$\int_{(\lambda^o,\mu^o)}^{(\lambda,\mu)} \frac{1}{|\Psi|}\,|\mathrm{d}\lambda'|$\, also equals zero for \,$(\lambda,\mu)=(\lambda^o,\mu^o)$\,, and we will show below that
\begin{equation}
\label{eq:jacobiprep:Psi-neu:int-limko}
\lim_{(\lambda,\mu)\to(\lambda^o,\mu^o)} \frac{\int_{(\lambda^o,\mu^o)}^{(\lambda,\mu)} \frac{1}{|\Psi|}\,|\mathrm{d}\lambda'|}{\left| \int_{(\lambda^o,\mu^o)}^{(\lambda,\mu)} \frac{1}{\Psi}\,\mathrm{d}\lambda' \right|} \leq \sqrt{2}
\end{equation}
holds. Moreover, there exists a constant \,$C_1>0$\, so that \,$\int_{(\lambda^o,\mu^o)}^{(\lambda,\mu)} \frac{1}{|\Psi|}\,|\mathrm{d}\lambda'| \leq C_1$\, holds.
It follows from these facts by an argument of compactness that \eqref{eq:jacobiprep:Psi-neu:int-claim} holds for \,$|\lambda^o|,|\lambda|\leq 2$\,.

For the proof of \eqref{eq:jacobiprep:Psi-neu:int-limko}: If \,$(\lambda^o,\mu^o) \neq \vkap_{\nu}$\, holds for \,$\nu\in\{1,2\}$\,, then we have \,$\Psi(\lambda^o,\mu^o)\neq 0$\, and therefore by l'Hospital's rule and the Fundamental Theorem of Calculus
$$ \lim_{(\lambda,\mu)\to(\lambda^o,\mu^o)} \frac{\int_{(\lambda^o,\mu^o)}^{(\lambda,\mu)} \frac{1}{|\Psi|}\,|\mathrm{d}\lambda'|}{\left| \int_{(\lambda^o,\mu^o)}^{(\lambda,\mu)} \frac{1}{\Psi}\,\mathrm{d}\lambda' \right|}
= \frac{\frac{1}{|\Psi(\lambda^o,\mu^o)|}}{\left| \frac{1}{\Psi(\lambda^o,\mu^o)} \right|} = 1 \leq \sqrt{2} \;. $$
We now look at the case \,$(\lambda^o,\mu^o)=\vkap_\nu$\,, say with \,$\nu=1$\,. We may suppose that \,$(\lambda,\mu)$\, is close to \,$(\lambda^o,\mu^o)$\,, more specifically that 
\begin{equation}
\label{eq:jacobiprep:Psi-neu:lambda-vkapnu}
|\lambda-\vkap_1|\leq 1 \qmq{and therefore} |\lambda-\vkap_2|\geq 1
\end{equation}
holds. Moreover, we may suppose that the path of integration in \eqref{eq:jacobiprep:Psi-neu:int-limko} is a straight line from \,$(\lambda^o,\mu^o)=\vkap_1$\,
to \,$(\lambda,\mu)$\, (this is an admissible path). 
From \eqref{eq:jacobiprep:Psi-neu:lambda-vkapnu} it follows that we have
$$ \frac{1}{|\Psi(\lambda',\mu')|} \leq \frac{1}{\sqrt{|\lambda'-\vkap_1|}} $$
for \,$(\lambda',\mu')$\, on the path of integration, and therefore
\begin{align*}
\int_{\vkap_1}^{(\lambda,\mu)} \frac{1}{|\Psi|}\,|\mathrm{d}\lambda'| & \leq  \int_{\vkap_1}^{(\lambda,\mu)} \frac{1}{\sqrt{|\lambda'-\vkap_1|}}\,|\mathrm{d}\lambda'|
& \overset{(*)}{=} \left| \int_{\vkap_1}^{(\lambda,\mu)} \frac{1}{\sqrt{\lambda'-\vkap_1}}\,\mathrm{d}\lambda' \right| = 2\,\sqrt{|\lambda-\vkap_1|} \;, 
\end{align*}
where the equals sign marked $(*)$ follows from our choice of the path of integration, which implies that \,$\frac{1}{\sqrt{\lambda'-\vkap_1}}$\, has constant argument. On the other hand, we have
$$ \int_{\vkap_1}^{(\lambda,\mu)} \frac{1}{\Psi}\,\mathrm{d}\lambda' = \ln\left( \frac{\lambda-\vkap_{*}+\Psi(\lambda,\mu)}{\vkap_1-\vkap_{*}+\Psi(\vkap_1)} \right) 
= \ln\left( \lambda + \sqrt{\lambda^2-1} \right) = \arcosh(\lambda) \; . $$
Thus we obtain, again by the rule of l'Hospital: 
$$ \lim_{(\lambda,\mu)\to\vkap_1} \frac{\int_{\vkap_1}^{(\lambda,\mu)} \frac{1}{|\Psi|}\,|\mathrm{d}\lambda'|}{\left| \int_{\vkap_1}^{(\lambda,\mu)} \frac{1}{\Psi}\,\mathrm{d}\lambda' \right|}
\leq \lim_{(\lambda,\mu)\to\vkap_1} \frac{2\sqrt{|\lambda-\vkap_1|}}{|\arcosh(\lambda)|} = \lim_{(\lambda,\mu)\to\vkap_1} \frac{\tfrac{1}{\sqrt{|\lambda-\vkap_1|}}}{\left| \tfrac{1}{\sqrt{\lambda^2-1}} \right|} = \sqrt{2} \; .  $$
This completes the proof of \eqref{eq:jacobiprep:Psi-neu:int-limko}. 

We now prove \eqref{eq:jacobiprep:Psi-neu:int-claim} in the case where at least one of the inequalities \,$|\lambda|>2$\, and \,$|\lambda^o|>2$\, holds. 
The parts of the admissible path of integration that run within \,$\Mengegr{\lambda'\in\C}{|\lambda'|\leq 2}$\, are handled by the above 
argument, so we only consider the parts that run within \,$\Mengegr{\lambda'\in\C}{|\lambda'|\geq 2}$\, in the sequel. For such \,$\lambda'$\,, we have \,$|\lambda'-\vkap_{\nu}| \geq \tfrac12\,|\lambda'|$\, and therefore
\begin{equation}
\label{eq:jacobiprep:Psi-neu:ist-extestim}
\frac{1}{|\Psi(\lambda',\mu')|} \leq \frac{2}{|\lambda'|} \;.
\end{equation}
We also note that we have
\begin{align*}
 \int_{(\lambda^o,\mu^o)}^{(\lambda,\mu)} \frac{1}{\Psi}\,\mathrm{d}\lambda'
& = \ln\left( \frac{\lambda-\vkap_{*}+\Psi(\lambda,\mu)}{\lambda^o-\vkap_{*}+\Psi(\lambda^o,\mu^o)} \right)  \\
& = \ln \left| \frac{\lambda-\vkap_{*}+\Psi(\lambda,\mu)}{\lambda^o-\vkap_{*}+\Psi(\lambda^o,\mu^o)} \right| + i\cdot \arg\left( \frac{\lambda-\vkap_{*}+\Psi(\lambda,\mu)}{\lambda^o-\vkap_{*}+\Psi(\lambda^o,\mu^o)} \right) 
\end{align*}
and therefore
$$ \left| \int_{(\lambda^o,\mu^o)}^{(\lambda,\mu)} \frac{1}{\Psi}\,\mathrm{d}\lambda' \right| \geq \max\left\{ \ln \left| \frac{\lambda-\vkap_{*}+\Psi(\lambda,\mu)}{\lambda^o-\vkap_{*}+\Psi(\lambda^o,\mu^o)} \right|\;,\;
\left| \arg\left( \frac{\lambda-\vkap_{*}+\Psi(\lambda,\mu)}{\lambda^o-\vkap_{*}+\Psi(\lambda^o,\mu^o)} \right) \right| \right\} \;. $$
In the sequel, we will handle the segments of the part of the admissible path of integration contained in \,$\{|\lambda|\geq 2\}$\, individually, treating lines and circles separately. We will estimate 
\,$\int \tfrac{1}{|\Psi|}\,|\mathrm{d}\lambda'|$\, by a multiple of \,$\ln \left| \frac{\lambda-\vkap_{*}+\Psi(\lambda,\mu)}{\lambda^o-\vkap_{*}+\Psi(\lambda^o,\mu^o)} \right|$\, for lines,
and by a multiple of \,$\left| \arg\left( \frac{\lambda-\vkap_{*}+\Psi(\lambda,\mu)}{\lambda^o-\vkap_{*}+\Psi(\lambda^o,\mu^o)} \right) \right|$\, for circles. For this purpose we denote by \,$\gamma$\, the relevant
path, and by \,$(\lambda_1,\mu_1)$\,, \,$(\lambda_2,\mu_2)$\, its endpoints.

For lines, we have by \eqref{eq:jacobiprep:Psi-neu:ist-extestim} with a constant \,$C_2>0$\,
$$ \int_\gamma \frac{1}{|\Psi|}\,|\mathrm{d}\lambda'| \leq \int_\gamma \frac{2}{|\lambda'|}\,|\mathrm{d}\lambda'| = 2\,\ln\left( \frac{|\lambda_2|}{|\lambda_1|} \right) 
\leq C_2\,\ln \left| \frac{\lambda-\vkap_{*}+\Psi(\lambda,\mu)}{\lambda^o-\vkap_{*}+\Psi(\lambda^o,\mu^o)} \right| \;. $$

For circles, \,$|\lambda|$\, is constant along the path of integration, and therefore we have again by \eqref{eq:jacobiprep:Psi-neu:ist-extestim} with a constant \,$C_3>0$\,
$$ \int_\gamma \frac{1}{|\Psi|}\,|\mathrm{d}\lambda'| \leq \int_\gamma \frac{2}{|\lambda'|}\,|\mathrm{d}\lambda'| = 2\,\arg\left( \frac{\lambda_2}{\lambda_1} \right) 
\leq C_3\,\left| \arg\left(\frac{\lambda-\vkap_{*}+\Psi(\lambda,\mu)}{\lambda^o-\vkap_{*}+\Psi(\lambda^o,\mu^o)} \right) \right|\;. $$

This completes the proof of \eqref{eq:jacobiprep:Psi-neu:int-estimate-nouveau}.

\emph{For (3).}
In the case \,$\vkap_{k,1}=\vkap_{k,2}$\,, the integrand equals \,$1$\, because of our hypothesis \,$\xi_k=\vkap_{k,*}$\, for this case, 
and therefore we then have
$$ \int_{\vkap_{k,1}}^{(\lambda_k,\mu_k)} \frac{|\lambda-\xi_k|}{|\Psi_k|}\,|\mathrm{d}\lambda|  \leq |\lambda_k-\vkap_{k,1}| \; . $$

Thus we only need to consider the case \,$\vkap_{k,1}\neq \vkap_{k,2}$\, in the sequel. 
For this case we will show as an intermediate step
that there are constants \,$C_3,C_4>0$\, so that 
\begin{align}
& \int_{\vkap_{k,1}}^{(\lambda_k,\mu_k)} \frac{|\lambda-\xi_k|}{|\Psi_k|}\,|\mathrm{d}\lambda| \notag \\
\label{eq:jacobiprep:Psi-neu:3:intermediate}
\leq\; & C_3 \cdot |\lambda_k-\vkap_{k,1}| + C_4 \cdot (|\xi_k-\vkap_{k,*}|+|\vkap_{k,1}-\vkap_{k,2}|) \cdot \left| \arcosh\left( \frac{\lambda_k-\vkap_{k,*}}{\tfrac12\,(\vkap_{k,1}-\vkap_{k,2})} \right) \right| 
\end{align}
holds.

Indeed, it follows from an argument of compactness that \,$C_3$\,, \,$C_4$\, can be chosen such that \eqref{eq:jacobiprep:Psi-neu:3:intermediate} holds for \,$|k|\leq N$\,. Thus we only consider \,$|k|>N$\, in
the sequel. Then there exists a constant \,$C_5>0$\, so that we have
\begin{align}
\int_{\vkap_{k,1}}^{(\lambda_k,\mu_k)} \frac{1}{|\Psi_k|}\,|\mathrm{d}\lambda|
& \overset{(2)}{\leq} C_5\cdot \left| \int_{\vkap_{k,1}}^{(\lambda_k,\mu_k)} \frac{1}{\Psi_k}\,\mathrm{d}\lambda \right| \notag \\
& \overset{(1)}{=} C_5\cdot \left| \ln\left( \frac{\lambda_k-\vkap_{k,*}+\Psi_k(\lambda_k,\mu_k)}{\vkap_{k,1}-\vkap_{k,*}+\Psi_k(\vkap_{k,1})} \right) \right| \notag \\
\label{eq:jacobiprep:Psi-neu:3:Psi-est}
& \overset{(*)}{=} C_5\cdot \left| \arcosh\left( \frac{\lambda_k-\vkap_{k,*}}{\tfrac12\,(\vkap_{k,1}-\vkap_{k,2})} \right) \right| \; , 
\end{align}
where (2) and (1) refer to the respective parts of the present lemma, and the equality marked $(*)$ follows from the 
equation \,$\arcosh(z)=\ln(z+\sqrt{z^2-1})$\,, and therefore
\begin{align*}
& \int_{\vkap_{k,1}}^{(\lambda_k,\mu_k)} \frac{|\lambda-\xi_k|}{|\Psi_k|}\,|\mathrm{d}\lambda|\\
\leq\; &  \int_{\vkap_{k,1}}^{(\lambda_k,\mu_k)} \frac{|\lambda-\vkap_{k,*}|}{|\Psi_k|}\,|\mathrm{d}\lambda| + |\vkap_{k,*}-\xi_k|\cdot \int_{\vkap_{k,1}}^{(\lambda_k,\mu_k)} \frac{1}{|\Psi_k|}\,|\mathrm{d}\lambda| \\
\leq\; & \int_{\vkap_{k,1}}^{(\lambda_k,\mu_k)} \frac{|\lambda-\vkap_{k,*}|}{|\Psi_k|}\,|\mathrm{d}\lambda| + C_5\cdot |\vkap_{k,*}-\xi_k|\cdot \left| \arcosh\left( \frac{\lambda_k-\vkap_{k,*}}{\tfrac12\,(\vkap_{k,1}-\vkap_{k,2})} \right) \right| \; .
\end{align*}
This calculation shows that it suffices to consider the case \,$\xi_k = \vkap_{k,*}$\,
in the proof of \eqref{eq:jacobiprep:Psi-neu:3:intermediate}. In this case, \eqref{eq:jacobiprep:Psi-neu:3:intermediate} 
takes the form of the following inequality
\begin{equation}
\label{eq:jacobiprep:Psi-neu:3:intermediate2}
\int_{\vkap_{k,1}}^{(\lambda_k,\mu_k)} \frac{|\lambda-\vkap_{k,*}|}{|\Psi_k|}\,|\mathrm{d}\lambda|
\leq C_3 \cdot |\lambda_k-\vkap_{k,1}| + C_4 \cdot |\vkap_{k,1}-\vkap_{k,2}|\cdot \left| \arcosh\left( \frac{\lambda_k-\vkap_{k,*}}{\tfrac12\,(\vkap_{k,1}-\vkap_{k,2})} \right) \right| \; ,
\end{equation}
which we will now prove.

For the proof of \eqref{eq:jacobiprep:Psi-neu:3:intermediate2} we first look at the case \,$|\lambda_k-\vkap_{k,*}|\leq 
|\vkap_{k,1}-\vkap_{k,2}|$\,. Then we have for any \,$\lambda \in [\vkap_{k,1},\lambda_k]$\,
\begin{align*}
|\lambda-\vkap_{k,*}| & \leq |\lambda-\vkap_{k,1}| + |\vkap_{k,1}-\vkap_{k,*}| \leq |\lambda_k-\vkap_{k,1}| + |\vkap_{k,1}-\vkap_{k,*}|  \\
& \leq |\lambda_k-\vkap_{k,*}| + 2\,|\vkap_{k,1}-\vkap_{k,*}| \leq 2\,|\vkap_{k,1}-\vkap_{k,2}|
\end{align*}
and therefore by \eqref{eq:jacobiprep:Psi-neu:3:Psi-est}
\begin{align}
\int_{\vkap_{k,1}}^{(\lambda_k,\mu_k)} \frac{|\lambda-\vkap_{k,*}|}{|\Psi_k|}\,|\mathrm{d}\lambda|
& \leq 2\,|\vkap_{k,1}-\vkap_{k,2}| \cdot \int_{\vkap_{k,1}}^{(\lambda_k,\mu_k)} \frac{1}{|\Psi_k|}\,|\mathrm{d}\lambda| \notag \\
\label{eq:jacobiprep:Psi-neu:3:case1}
& \leq 2\,C_5\,|\vkap_{k,1}-\vkap_{k,2}| \cdot \left| \arcosh\left( \frac{\lambda_k-\vkap_{k,*}}{\tfrac12\,(\vkap_{k,1}-\vkap_{k,2})} \right) \right| \; ,
\end{align}
whence \eqref{eq:jacobiprep:Psi-neu:3:intermediate2} follows for this case.

We now consider the case \,$|\lambda_k-\vkap_{k,*}| > |\vkap_{k,1}-\vkap_{k,2}|$\,. Because the endpoints of the line segment
\,$[\vkap_{k,1},\lambda_k]$\, then satisfy
$$ |\vkap_{k,1}-\vkap_{k,*}| < |\vkap_{k,1}-\vkap_{k,2}| \qmq{and} |\lambda_k-\vkap_{k,*}| > |\vkap_{k,1}-\vkap_{k,2}| \;, $$
there is a unique intersection point between this line segment and the circle
\,$\Mengegr{\lambda\in\C}{|\lambda-\vkap_{k,*}|= |\vkap_{k,1}-\vkap_{k,2}|}$\,,
which we denote by \,$\lambda_* \in \C^*$\,. 
We further choose \,$\mu_*\in \C$\, so that \,$(\lambda_*,\mu_*)\in \Sigma$\, is in the appropriate leaf of \,$\Sigma$\,.
Then we have
\begin{equation}
\label{eq:jacobiprep:Psi-neu:3:intlin-estimate-split}
\int_{\vkap_{k,1}}^{(\lambda_k,\mu_k)} \frac{|\lambda-\vkap_{k,*}|}{|\Psi_k|}\,|\mathrm{d}\lambda| = 
\int_{\vkap_{k,1}}^{(\lambda_*,\mu_*)} \frac{|\lambda-\vkap_{k,*}|}{|\Psi_k|}\,|\mathrm{d}\lambda| + \int_{(\lambda_*,\mu_*)}^{(\lambda_k,\mu_k)} \frac{|\lambda-\vkap_{k,*}|}{|\Psi_k|}\,|\mathrm{d}\lambda| \;,
\end{equation}
where we integrate along straight lines in all three integrals.

We estimate the two integrals on the right hand side of Equation~\eqref{eq:jacobiprep:Psi-neu:3:intlin-estimate-split} separately. 
For the first integral, the estimate \eqref{eq:jacobiprep:Psi-neu:3:case1} is applicable because of \,$|\lambda_*-\vkap_{k,*}|
= |\vkap_{k,1}-\vkap_{k,2}|$\,, and therefore we have
\begin{align}
\int_{\vkap_{k,1}}^{(\lambda_*,\mu_*)} \frac{|\lambda-\vkap_{k,*}|}{|\Psi_k|}\,|\mathrm{d}\lambda|
& \leq 2\,C_5\,|\vkap_{k,1}-\vkap_{k,2}| \cdot \left| \arcosh\left( \frac{\lambda_*-\vkap_{k,*}}{\tfrac12\,(\vkap_{k,1}-\vkap_{k,2})} \right) \right| \notag \\
& \leq 2\,C_5\,|\vkap_{k,1}-\vkap_{k,2}| \cdot \left| \arcosh\left( \frac{\lambda_k-\vkap_{k,*}}{\tfrac12\,(\vkap_{k,1}-\vkap_{k,2})} \right) \right| \; . \notag
\end{align}
Concerning the second integral, we note that for \,$\lambda \in [\lambda_*,\lambda_k]$\, and \,$\nu\in\{1,2\}$\, we have 
$$ |\lambda-\vkap_{k,*}| \geq |\vkap_{k,1}-\vkap_{k,2}| = 2\,|\vkap_{k,*}-\vkap_{k,\nu}| $$
and therefore 
$$ |\lambda-\vkap_{k,\nu}| \geq |\lambda-\vkap_{k,*}| - |\vkap_{k,*}-\vkap_{k,\nu}| \geq \frac12\,|\lambda-\vkap_{k,*}| \;,$$
whence
$$ \frac{|\lambda-\vkap_{k,*}|}{|\Psi_k(\lambda,\mu)|} = \frac{|\lambda-\vkap_{k,*}|}{\sqrt{|\lambda-\vkap_{k,1}|\cdot |\lambda-\vkap_{k,2}|}} \leq 2 $$
follows. Thus we obtain
$$ \int_{(\lambda_*,\mu_*)}^{(\lambda_k,\mu_k)} \frac{|\lambda-\vkap_{k,*}|}{|\Psi_k|}\,|\mathrm{d}\lambda |
\leq \int_{(\lambda_*,\mu_*)}^{(\lambda_k,\mu_k)} 2\,|\mathrm{d}\lambda| = 2\,|\lambda_k-\lambda_*| \leq 2\,|\lambda_k-\vkap_{k,1}| \;. $$
By plugging these estimates into Equation~\eqref{eq:jacobiprep:Psi-neu:3:intlin-estimate-split}, 
we obtain the estimate \eqref{eq:jacobiprep:Psi-neu:3:intermediate2} also for the case 
of \,$|\lambda_k-\vkap_{k,*}| > |\vkap_{k,1}-\vkap_{k,2}|$\,.
This concludes the proof of \eqref{eq:jacobiprep:Psi-neu:3:intermediate2} and therefore of \eqref{eq:jacobiprep:Psi-neu:3:intermediate}.

To derive the estimate claimed in part (3) of the lemma from \eqref{eq:jacobiprep:Psi-neu:3:intermediate}, we use the fact that
there exists a constant \,$C_{\arcosh}>0$\, so that for every \,$z\in \C$\, we have 
\begin{equation}
\label{eq:jacobiprep:Psi-neu:arcosh-pre}
|\arcosh(z)| \leq C_{\arcosh} \cdot |z-1| \; .
\end{equation}
For the proof of \eqref{eq:jacobiprep:Psi-neu:arcosh-pre} we note that \,$\arcosh(z)$\, is differentiable on a suitably doubly slitted
plane with \,$\tfrac{\mathrm{d}\ }{\mathrm{d}z} \arcosh(z) = \tfrac{1}{\sqrt{z^2-1}}$\, and \,$\arcosh(1)=0$\,, and therefore we have
$$ |\arcosh(z)| = |\arcosh(z)-\arcosh(1)| = \left| \int_1^z \frac{1}{\sqrt{z^2-1}}\,\mathrm{d}z \right| \leq \int_1^z \frac{1}{\sqrt{|z^2-1|}}\,|\mathrm{d}z|  \;, $$
where the integration takes place along a path that is chosen in the slitted plane with a length commensurate with \,$|z-1|$\,. 
Because the integrand has only finitely many poles of order at most \,$\tfrac12$\,, the above integral is finite, 
and because the integrand is bounded for \,$|z|\to \infty$\,, there in fact exists a constant \,$C_{\arcosh}>0$\, with
$$ |\arcosh(z)| \leq \int_1^z \frac{1}{\sqrt{|z^2-1|}}\,|\mathrm{d}z| \leq C_{\arcosh}\cdot |z-1|\; . $$

By applying \eqref{eq:jacobiprep:Psi-neu:arcosh-pre} to the estimate \eqref{eq:jacobiprep:Psi-neu:3:intermediate}, we see
\begin{align*}
& \int_{\vkap_{k,1}}^{(\lambda_k,\mu_k)} \frac{|\lambda-\xi_k|}{|\Psi_k|}\,|\mathrm{d}\lambda| \\
\leq \; & C_3 \cdot |\lambda_k-\vkap_{k,1}| + C_4 \cdot (|\xi_k-\vkap_{k,*}|+|\vkap_{k,1}-\vkap_{k,2}|) \cdot C_{\arcosh} \cdot \left| \frac{\lambda_k-\vkap_{k,*}}{\tfrac12\,(\vkap_{k,1}-\vkap_{k,2})} -1 \right| \\
= \; & C_3 \cdot |\lambda_k-\vkap_{k,1}| + C_4 \cdot (|\xi_k-\vkap_{k,*}|+|\vkap_{k,1}-\vkap_{k,2}|) \cdot C_{\arcosh} \cdot \left| \frac{\lambda_k-\vkap_{k,1}}{\tfrac12\,(\vkap_{k,1}-\vkap_{k,2})} \right| \\
= \; & \left( \left( C_3+ 2\,C_4\,C_{\arcosh} \right) + 2\,C_4\,C_{\arcosh}\,\left| \frac{\xi_k-\vkap_{k,*}}{\vkap_{k,1}-\vkap_{k,2}} \right| \right) \cdot |\lambda_k-\vkap_{k,1}| \; . 
\end{align*}
Thus the estimate claimed in (3) holds (with \,$C_1 := 2\,C_4\,C_{\arcosh}$\, and \,$C_2 := C_3+ 2\,C_4\,C_{\arcosh}$\,).

\emph{For (4).}
If \,$\vkap_{k,1}=\vkap_{k,2}$\, holds, the integrand equals \,$1$\, (because we had required \,$\xi_k=\vkap_{k,*}$\, in this case)\,, and thus the integral
equals \,$\vol(U_{k,\delta}) \in \ell^\infty_{-2,6}(k)$\,. Thus we now suppose \,$\vkap_{k,1}\neq \vkap_{k,2}$\,, consider the case \,$k>0$\,, and choose
\,$\delta'>0$\, (independently of \,$k$\,) so that \,$U_{k,\delta} \subset B(\lambda_{k,0},k\delta')$\, holds. We then compute \,$\int_{B(\lambda_{k,0},k\delta')} \tfrac{|\lambda-\xi_k|^2}{|\lambda-\vkap_{k,1}|\cdot|\lambda-\vkap_{k,2}|}
\,\mathrm{d}^2\lambda$\,. We abbreviate \,$A:= \tfrac{|\xi_k-\vkap_{k,*}|}{|\vkap_{k,1}-\vkap_{k,2}|}$\, and 
split the domain of integration \,$B(\lambda_{k,0},k\delta')$\, into four parts:
\begin{align*}
M_1 & := B(\vkap_{k,1},\tfrac12|\vkap_{k,1}-\vkap_{k,2}|) \\
M_2 & := B(\vkap_{k,2},\tfrac12|\vkap_{k,1}-\vkap_{k,2}|) \\
M_3 & := B(\vkap_{k,*},\max\{|\xi_k-\vkap_{k,*}|,|\vkap_{k,1}-\vkap_{k,2}|\}) \setminus (M_1 \cup M_2) \\
M_4 & := B(\lambda_{k,0},k\delta') \setminus (M_1 \cup M_2 \cup M_3) \; . 
\end{align*}

For \,$\nu \in \{1,2\}$\, and \,$\lambda \in M_\nu$\, we have
\begin{align*}
|\lambda-\xi_k| & \leq |\lambda-\vkap_{k,\nu}| + |\vkap_{k,\nu}-\vkap_{k,*}| + |\vkap_{k,*}-\xi_k| \\
& \leq \frac12\,|\vkap_{k,1}-\vkap_{k,2}| + \frac12\,|\vkap_{k,1}-\vkap_{k,2}| + A\,|\vkap_{k,1}-\vkap_{k,2}| \\
& = (A+1)\cdot |\vkap_{k,1}-\vkap_{k,2}| 
\end{align*}
and
$$ |\lambda-\vkap_{k,3-\nu}| \geq |\vkap_{k,3-\nu}-\vkap_{k,\nu}| - |\lambda-\vkap_{k,\nu}| \geq \frac12\,|\vkap_{k,1}-\vkap_{k,2}| \; , $$
and therefore
$$ \frac{|\lambda-\xi_k|^2}{|\lambda-\vkap_{k,1}|\cdot|\lambda-\vkap_{k,2}|} \leq 2\,(A+1)^2\cdot |\vkap_{k,1}-\vkap_{k,2}| \cdot \frac{1}{|\lambda-\vkap_{k,\nu}|} \; . $$
Thus we obtain
\begin{align}
\int_{M_\nu} \frac{|\lambda-\xi_k|^2}{|\lambda-\vkap_{k,1}|\cdot|\lambda-\vkap_{k,2}|}\,\mathrm{d}^2\lambda
& \leq 2\,(A+1)^2\cdot |\vkap_{k,1}-\vkap_{k,2}| \cdot \int_{M_\nu} \frac{1}{|\lambda-\vkap_{k,\nu}|}\,\mathrm{d}^2\lambda \notag \\
& \leq 2\,(A+1)^2\cdot |\vkap_{k,1}-\vkap_{k,2}| \cdot \int_{r=0}^{\tfrac12|\vkap_{k,1}-\vkap_{k,2}|} \int_{\vi=0}^{2\pi} \frac{1}{r}\,r\,\mathrm{d}r\,\mathrm{d}\vi \notag \\
& = 2\pi\,(A+1)^2\,|\vkap_{k,1}-\vkap_{k,2}|^2 \notag\\ 
\label{eq:jacobiprep:Psi:6:xiM12}
& \leq  8\pi\,(|\xi_k-\vkap_{k,*}|^2+|\vkap_{k,1}-\vkap_{k,2}|^2) \; . 
\end{align}

For \,$\lambda \in M_3$\, we have 
$$ |\lambda-\xi_k| \leq |\lambda-\vkap_{k,*}| + |\vkap_{k,*}-\xi_k| \leq 2\,\max\{|\xi_k-\vkap_{k,*}|,|\vkap_{k,1}-\vkap_{k,2}|\} $$
and for \,$\nu \in \{1,2\}$\, 
$$ |\lambda-\vkap_{k,\nu}| \geq \frac12\,|\vkap_{k,1}-\vkap_{k,2}| \;, $$
and therefore
$$ \frac{|\lambda-\xi_k|^2}{|\lambda-\vkap_{k,1}|\cdot|\lambda-\vkap_{k,2}|} \leq 16\,\frac{(\max\{|\xi_k-\vkap_{k,*}|,|\vkap_{k,1}-\vkap_{k,2}|\})^2}{|\vkap_{k,1}-\vkap_{k,2}|^2} \; . $$
It follows that we have
\begin{align}
& \int_{M_3} \frac{|\lambda-\xi_k|^2}{|\lambda-\vkap_{k,1}|\cdot|\lambda-\vkap_{k,2}|}\,\mathrm{d}^2\lambda
\leq 16\,\frac{(\max\{|\xi_k-\vkap_{k,*}|,|\vkap_{k,1}-\vkap_{k,2}|\})^2}{|\vkap_{k,1}-\vkap_{k,2}|^2} \cdot \vol(M_3) \notag \\
\leq & 16\,\frac{(\max\{|\xi_k-\vkap_{k,*}|,|\vkap_{k,1}-\vkap_{k,2}|\})^2}{|\vkap_{k,1}-\vkap_{k,2}|^2} \cdot \vol(B(\vkap_{k,*},\max\{|\xi_k-\vkap_{k,*}|,|\vkap_{k,1}-\vkap_{k,2}|\})) \notag \\
= & 16\,\frac{(\max\{|\xi_k-\vkap_{k,*}|,|\vkap_{k,1}-\vkap_{k,2}|\})^2}{|\vkap_{k,1}-\vkap_{k,2}|^2} \cdot \pi \,(\max\{|\xi_k-\vkap_{k,*}|,|\vkap_{k,1}-\vkap_{k,2}|\})^2 \notag \\
\label{eq:jacobiprep:Psi:6:xiM3}
\leq & 16\,\pi\,\left( |\vkap_{k,1}-\vkap_{k,2}|^2 + \frac{|\xi_k-\vkap_{k,*}|^4}{|\vkap_{k,1}-\vkap_{k,2}|^2} \right) \; . 
\end{align}

Finally, for \,$\lambda\in M_4$\, we have
$$ |\lambda-\xi_k| \leq |\lambda-\vkap_{k,*}| + |\vkap_{k,*}-\xi_k| \leq |\lambda-\vkap_{k,*}| + |\lambda-\vkap_{k,*}| = 2\,|\lambda-\vkap_{k,*}| $$
and for \,$\nu \in \{1,2\}$\,
\begin{align*}
|\lambda-\vkap_{k,\nu}| & \geq |\lambda-\vkap_{k,*}| - |\vkap_{k,*}-\vkap_{k,\nu}| = |\lambda-\vkap_{k,*}| - \frac12\,|\vkap_{k,1}-\vkap_{k,2}| \geq \frac12\,|\lambda-\vkap_{k,*}|
\end{align*}
and hence
$$ \frac{|\lambda-\xi_k|^2}{|\lambda-\vkap_{k,1}|\cdot|\lambda-\vkap_{k,2}|} \leq 16 \;. $$
Thus we obtain
\begin{equation}
\label{eq:jacobiprep:Psi:6:xiM4}
\int_{M_4} \frac{|\lambda-\xi_k|^2}{|\lambda-\vkap_{k,1}|\cdot|\lambda-\vkap_{k,2}|}\,\mathrm{d}^2\lambda 
\leq 16\cdot \vol(M_4) \leq 16 \cdot \vol(B(\lambda_{k,0},k\delta')) = 16 \cdot \pi\,(k\delta')^2 \leq 64\,\pi\,k^2 
\end{equation}
It follows from Equations~\eqref{eq:jacobiprep:Psi:6:xiM12}, \eqref{eq:jacobiprep:Psi:6:xiM3} and \eqref{eq:jacobiprep:Psi:6:xiM4}, and the fact that \,$\vkap_{k,1}-\vkap_{k,2}\in \ell^2_{-1,3}(k)$\, holds,
that there exists \,$C>0$\, so that 
$$ \int_{U_{k,\delta}} \frac{|\lambda-\xi_k|^2}{|\lambda-\vkap_{k,1}|\cdot|\lambda-\vkap_{k,2}|}\,\mathrm{d}^2\lambda \leq C \cdot \max\left\{ 1, \frac{|\xi_k-\vkap_{k,*}|}{|\vkap_{k,1}-\vkap_{k,2}|} \right\}^4 \cdot (|\vkap_{k,1}-\vkap_{k,2}|^2 + k^2) $$
holds.
\end{proof}

\begin{prop}
\label{P:jacobiprep:Deltaint-neu}
Let \,$N\in \N$\, be as in Lemma~\ref{L:jacobiprep:Psi-intro}.
\begin{enumerate}
\item
There exist \,$C>0$\, and a sequence \,$b_k \in \ell^2_{-1,3}(k)$\, (depending on \,$\Sigma$\,) such that for any \,$k\in \Z$\, with \,$|k|> N$\, and for any \,$(\lambda_k^o,\mu_k^o),(\lambda_k,\mu_k) \in \wh{U}_{k,\delta}'$\, 
(in the case \,$\vkap_{k,1}=\vkap_{k,2}$\, we require that \,$(\lambda_k,\mu_k)$\, and \,$(\lambda_k^o,\mu_k^o)$\, lie in the same connected component of \,$\wh{U}_{k,\delta}'$\,)
we have
\begin{align*}
& \left| \int_{(\lambda_k^o,\mu_k^o)}^{(\lambda_k,\mu_k)} \frac{1}{\mu-\mu^{-1}}\,\mathrm{d}\lambda - 4i(-1)^k\,\lambda_{k,0}^{\rho/2}\,\ln \left( \frac{\lambda_k-\vkap_{k,*}+\Psi_k(\lambda_k,\mu_k)}{\lambda_k^o-\vkap_{k,*}+\Psi_k(\lambda_k^o,\mu_k^o)} \right) \right|\\
\leq\; & C \cdot \left( |\lambda_k-\lambda_k^o|+b_k \right) \cdot \left| \ln \left( \frac{\lambda_k-\vkap_{k,*}+\Psi_k(\lambda_k,\mu_k)}{\lambda_k^o-\vkap_{k,*}+\Psi_k(\lambda_k^o,\mu_k^o)} \right) \right| \; . 
\end{align*}
Here we integrate along any path that is contained in \,$\wh{U}_{k,\delta}'$\,, set \,$\rho=1$\, for \,$k>0$\, and \,$\rho=3$\, for \,$k<0$\,, 
and \,$\ln(z)$\, is the branch of the complex logarithm function with \,$\ln(1)=2\pi i m$\,, 
where \,$m\in \Z$\, is the winding number of the path of integration around the pair of branch points \,$\vkap_{k,1}$\,, \,$\vkap_{k,2}$\,.


\item
There exists  \,$C>0$\, so that for any \,$k\in \Z$\, with \,$|k|>N$\, and for any \,$(\lambda_k^o,\mu_k^o),(\lambda_k,\mu_k) \in \wh{U}_{k,\delta}'$\, as in (1)
we have
\begin{align}
\label{eq:jacobiprep:Deltaint-neu:absint}
\int_{(\lambda_k^o,\mu_k^o)}^{(\lambda_k,\mu_k)} \frac{1}{|\mu-\mu^{-1}|}\,|\mathrm{d}\lambda|
& \leq C\,\lambda_{k,0}^{\rho/2} \cdot \left| \ln \left( \frac{\lambda_k-\vkap_{k,*}+\Psi_k(\lambda_k,\mu_k)}{\lambda_k^o-\vkap_{k,*}+\Psi_k(\lambda_k^o,\mu_k^o)} \right) \right| \;, 
\end{align}
where we integrate along an admissible path (Definition~\ref{D:jacobiprep:admissible}), 
and set \,$\rho=1$\, for \,$k>0$\, and \,$\rho=3$\, for \,$k<0$\,.


Moreover, for all \,$k$\, with \,$\vkap_{k,1}\neq \vkap_{k,2}$\, we have
\begin{align}
\label{eq:transdgl:Delta2int}
\int_{\vkap_{k,1}}^{\vkap_{k,2}} \frac{1}{\mu-\mu^{-1}}\,\mathrm{d}\lambda & = -4(-1)^k\pi\,\lambda_{k,0}^{\rho/2} + \ell^{2}_{-1,3}(k) \; , \\
\label{eq:transdgl:Delta2intabs}
\int_{\vkap_{k,1}}^{\vkap_{k,2}} \frac{1}{|\mu-\mu^{-1}|}\,|\mathrm{d}\lambda| & = 4\pi\,\lambda_{k,0}^{\rho/2} + \ell^{2}_{-1,3}(k) \; , 
\end{align}
where the path of integration is a straight line from \,$\vkap_{k,1}$\, to \,$\vkap_{k,2}$\, in the \,$\lambda$-plane.

\item
There exist \,$C_1,C_2>0$\, so that for any \,$k\in \Z$\,, for any \,$(\lambda_k^o,\mu_k^o),(\lambda_k,\mu_k) \in \wh{U}_{k,\delta}$\, as in (1), and for arbitrary \,$\xi_k\in U_{k,\delta}$\, if \,$\vkap_{k,1}\neq \vkap_{k,2}$\,, whereas we set \,$\xi_k := \vkap_{k,*}$\, if \,$\vkap_{k,1}=\vkap_{k,2}$\,, we have
$$ \int_{\vkap_{k,1}}^{(\lambda_k,\mu_k)} \frac{|\lambda-\xi_k|}{|\mu-\mu^{-1}|}\,|\mathrm{d}\lambda|
\leq \left( C_1\,\left| \frac{\xi_k-\vkap_{k,*}}{\vkap_{k,1}-\vkap_{k,2}} \right| + C_2 \right) \cdot \lambda_{k,0}^{\rho/2}\cdot |\lambda_k-\vkap_{k,1}| \;, $$
where the path of integration is a straight line from \,$\vkap_{k,1}$\, to \,$\lambda_k$\, in the \,$\lambda$-plane,
we set \,$\rho=1$\, for \,$k>0$\, and \,$\rho=3$\, for \,$k<0$\,,
and the expression \,$\left| \frac{\xi_k-\vkap_{k,*}}{\vkap_{k,1}-\vkap_{k,2}} \right|$\, is to be read as \,$0$\, 
in the case \,$\vkap_{k,1}=\vkap_{k,2}$\,.  

\item
There exists \,$C>0$\, so that for any \,$k\in \Z$\, and \,$\xi_k$\, as in (3) we have
$$ \int_{U_{k,\delta}} \frac{|\lambda-\xi_k|^2}{|\Delta(\lambda)^2-4|}\,\mathrm{d}^2\lambda 
\leq C \,\lambda_{k,0}^\rho\, \cdot \left( |\xi_k-\vkap_{k,*}|^2+ \frac{|\xi_k-\vkap_{k,*}|^4}{|\vkap_{k,1}-\vkap_{k,2}|^2} + \ell^\infty_{-2,6}(k) \right) \; , $$
where \,$\mathrm{d}^2\lambda$\, again denotes the 2-dimensional Lebesgue measure on \,$\C$\,. 
Here we again put \,$\rho=1$\, for \,$k>0$\, and \,$\rho=3$\, for \,$k<0$\,; and in the case \,$\vkap_{k,1}=\vkap_{k,2}$\,, the expression \,$\frac{|\xi_k-\vkap_{k,*}|^4}{|\vkap_{k,1}-\vkap_{k,2}|^2}$\, is to be read as equal to \,$0$\,. 
\end{enumerate}
\end{prop}

\begin{proof}
\emph{For (1).}
We may suppose without loss of generality that the path of integration is composed of an admissible path connecting \,$(\lambda_k^o,\mu_k^o)$\, to \,$(\lambda_k,\mu_k)$\, and \,$|m|$\, circles around
\,$\vkap_{k,1}$\,, \,$\vkap_{k,2}$\, (with the orientation given by the sign of \,$m$\,). Note that each of the circles is also an admissible path.

If the point \,$(\lambda,\mu)$\, is on the trace of the path of integration, we then have because of Equation~\eqref{eq:jacobiprep:admissible:estimate} in Definition~\ref{D:jacobiprep:admissible}
\begin{align*}
|\lambda-\vkap_{k,\nu}| & \leq \max \bigr\{ \,2\,|\vkap_{k,1}-\vkap_{k,2}|\,,\, |\lambda_k^o-\vkap_{k,*}|\,,\, |\lambda_k-\vkap_{k,*}|\,\bigr\} \\
& \leq |\lambda_k-\lambda_k^o| + \wt{b}_k \qmq{with} \wt{b}_k := 2\,|\vkap_{k,1}-\vkap_{k,2}| + |\lambda_k^o-\vkap_{k,*}| \in \ell^2_{-1,3}(k) \; . 
\end{align*}
Therefore it follows from Lemma~\ref{L:jacobiprep:Psi-intro} that there exist constants \,$C_1,C_2>0$\, and \,$a_k \in \ell^2_{-1,3}(k)$\, so that we have for such points
\begin{align}
\left| \frac{1}{\mu-\mu^{-1}} - 4i(-1)^k\,\lambda_{k,0}^{\rho/2}\,\frac{1}{\Psi_k(\lambda,\mu)} \right|
\label{eq:transdgl:Delta2:asymp-Delta-12-better-pre}
& \leq \frac{C_1 \cdot \left( |\lambda-\vkap_{k,1}| + |\lambda-\vkap_{k,2}| \right) + a_k}{|\Psi_k(\lambda,\mu)|} \\
\label{eq:transdgl:Delta2:asymp-Delta-12-better}
& \leq C_2 \cdot \frac{|\lambda_k-\lambda_k^o| + b_k}{|\Psi_k(\lambda,\mu)|} \qmq{with} b_k := a_k + \wt{b}_k \; . 
\end{align}
We obtain by integration
\begin{gather*}
\left| \int_{(\lambda_k^o,\mu_k^o)}^{(\lambda_k,\mu_k)} \frac{1}{\mu-\mu^{-1}}\,\mathrm{d}\lambda - 4i(-1)^k\,\lambda_{k,0}^{\rho/2}\,\int_{(\lambda_k^o,\mu_k^o)}^{(\lambda_k,\mu_k)} \frac{1}{\Psi_k} \,\mathrm{d}\lambda \right| \hspace{3cm} \\
\hspace{3cm} \leq C_2 \cdot \left( |\lambda_k-\lambda_k^o| + b_k \right) \cdot \int_{(\lambda_k^o,\mu_k^o)}^{(\lambda_k,\mu_k)} \frac{1}{|\Psi_k|} \,\mathrm{d}\lambda \; . 
\end{gather*}
(1) now follows from Lemma~\ref{L:jacobiprep:Psi-neu}(1),(2).

\emph{For (2).}
It follows from Equation~\eqref{eq:transdgl:Delta2:asymp-Delta-12-better-pre} that there exists \,$C_3>0$\, so that we have for \,$(\lambda,\mu) \in \wh{U}_{k,\delta}$\,
$$ 
\frac{1}{|\mu-\mu^{-1}|} \leq C_3\,\lambda_{k,0}^{\rho/2} \,\frac{1}{|\Psi_k(\lambda,\mu)|} $$
and therefore
$$ 
\int_{(\lambda_k^o,\mu_k^o)}^{(\lambda_k,\mu_k)}\frac{1}{|\mu-\mu^{-1}|}\,|\mathrm{d}\lambda| \leq C_3\,\lambda_{k,0}^{\rho/2} \,\int_{(\lambda_k^o,\mu_k^o)}^{(\lambda_k,\mu_k)}\frac{1}{|\Psi_k(\lambda,\mu)|}\,|\mathrm{d}\lambda|\;.  $$
By Lemma~\ref{L:jacobiprep:Psi-neu}(2),(1) there exists \,$C_4>0$\, with
$$ \int_{(\lambda_k^o,\mu_k^o)}^{(\lambda_k,\mu_k)}\frac{1}{|\Psi_k(\lambda,\mu)|}\,|\mathrm{d}\lambda| \leq C_4 \cdot \left| \int_{(\lambda_k^o,\mu_k^o)}^{(\lambda_k,\mu_k)}\frac{1}{\Psi_k(\lambda,\mu)}\,\mathrm{d}\lambda \right| 
= C_4\cdot \left| \ln \left( \frac{\lambda_k-\vkap_{k,*}+\Psi_k(\lambda_k,\mu_k)}{\lambda_k^o-\vkap_{k,*}+\Psi_k(\lambda_k^o,\mu_k^o)} \right) \right| \;, $$
and Equation~\eqref{eq:jacobiprep:Deltaint-neu:absint} follows from these estimates.

Equations~\eqref{eq:transdgl:Delta2int} and \eqref{eq:transdgl:Delta2intabs}
follow from part (1) resp.~from Equation~\eqref{eq:jacobiprep:Deltaint-neu:absint} by plugging in \,$(\lambda_k^o,\mu_k^o)=\vkap_{k,1}$\, and
\,$(\lambda_k,\mu_k)=\vkap_{k,2}$\, and noting that \,$\vkap_{k,2}-\vkap_{k,1} \in \ell^2_{-1,3}(k)$\, and \,$\ln\left(\tfrac{\vkap_{k,2}-\vkap_{k,*}-\Psi_k(\vkap_{k,2})}{\vkap_{k,1}-\vkap_{k,*}-\Psi_k(\vkap_{k,1})}\right)=\ln(-1)=i\pi$\,
holds.

\emph{For (3).}
It follows from \eqref{eq:transdgl:Delta2:asymp-Delta-12-better} (applied with \,$\lambda_k^o=\vkap_{k,1}$\,) that we have
$$ \frac{|\lambda-\xi_k|}{|\mu-\mu^{-1}|} \leq C_5 \cdot \left( |\lambda_k-\vkap_{k,1}| + b_k \right) \cdot \frac{|\lambda-\xi_k|}{|\Psi_k(\lambda)|} \leq C_6 \cdot \lambda_{k,0}^{\rho/2} \cdot \frac{|\lambda-\xi_k|}{|\Psi_k(\lambda)|} $$
with constants \,$C_5,C_6>0$\, and a sequence \,$(b_k) \in \ell^2_{-1,3}(k)$\,. By integration we thus obtain 
$$ \int_{\vkap_{k,1}}^{(\lambda_k,\mu_k)} \frac{|\lambda-\xi_k|}{|\mu-\mu^{-1}|}\,|\mathrm{d}\lambda|
\leq C_6 \cdot \lambda_{k,0}^{\rho/2} \cdot \int_{\vkap_{k,1}}^{(\lambda_k,\mu_k)} \frac{|\lambda-\xi_k|}{|\Psi_k|} \,|\mathrm{d}\lambda| \; . $$
The claimed estimate now follows from the application of Lemma~\ref{L:jacobiprep:Psi-neu}(3).

\emph{For (4).}
From Equation~\eqref{eq:transdgl:Delta2:asymp:Delta} we obtain for \,$\lambda\in U_{k,\delta}$\,, because \,$z\mapsto z^{-1}$\, is Lipschitz continuous near \,$z=-\tfrac{1}{16}$\,,
\begin{gather*}
\left| \lambda_{k,0}^{-\rho} \cdot \frac{(\lambda-\vkap_{k,1})\cdot(\lambda-\vkap_{k,2})}{\Delta(\lambda)^2-4} - \left( -16 \right) \right|\\
\leq C_7 \cdot \left\{ \begin{matrix} k^{-1} & \text{if \,$k>0$\,} \\ k^3 & \text{if \,$k<0$\,} \end{matrix} \right\}  \cdot \left( |\lambda-\vkap_{k,1}| + |\lambda-\vkap_{k,2}| \right) + \ell^{2}_{0,0}(k) 
\end{gather*}
and therefore
\begin{gather*}
\left| \frac{1}{\Delta(\lambda)^2-4} - \left( - \frac{16\,\lambda_{k,0}^\rho}{(\lambda-\vkap_{k,1})\cdot(\lambda-\vkap_{k,2})} \right) \right| \\
\leq \frac{\lambda_{k,0}^\rho}{|\lambda-\vkap_{k,1}|\cdot|\lambda-\vkap_{k,2}|}\cdot \left( C_7\cdot \left\{ \begin{matrix} k^{-1} & \text{if \,$k>0$\,} \\ k^3 & \text{if \,$k<0$\,} \end{matrix} \right\}  \cdot \left( |\lambda-\vkap_{k,1}| + |\lambda-\vkap_{k,2}| \right) + \ell^{2}_{0,0}(k)\right) \;,
\end{gather*}
whence
$$ \left| \frac{1}{\Delta(\lambda)^2-4} \right| \leq C_8\,\lambda_{k,0}^\rho\,\frac{1}{|\lambda-\vkap_{k,1}|\cdot |\lambda-\vkap_{k,2}|} $$
follows. By integration and application of Lemma~\ref{L:jacobiprep:Psi-neu}(4), we obtain
\begin{align*}
\int_{U_{k,\delta}} \frac{|\lambda-\xi_k|^2}{|\Delta(\lambda)^2-4|}\,\mathrm{d}^2\lambda 
& \leq C_8\,\lambda_{k,0}^\rho\,\int_{U_{k,\delta}} \frac{|\lambda-\xi_k|^2}{|\lambda-\vkap_{k,1}|\cdot |\lambda-\vkap_{k,2}|}\,\mathrm{d}^2\lambda \\
& \leq C_{9}\,\lambda_{k,0}^\rho\, \left( |\xi_k-\vkap_{k,*}|^2 + \frac{|\xi_k-\vkap_{k,*}|^4}{|\vkap_{k,1}-\vkap_{k,2}|^2} + \ell^\infty_{-2,6}(k) \right) \; . 
\end{align*}
\end{proof}

\section{Asymptotic behavior of 1-forms on the spectral curve}
\label{Se:holo1forms}

\label{not:holo1forms:homology-basis}
We continue our preparations for the construction of the Jacobi coordinates on the spectral curve \,$\Sigma$\,. 
One important step in this construction is to obtain a basis \,$(\omega_n)_{n\in \Z}$\, of the space of square-integrable, holomorphic%
\footnote{The requirements of square-integrability and holomorphy need to be modified near singular points of \,$\Sigma$\,, see the precise
statements below.}
\,$1$-forms on \,$\Sigma$\, which is dual to a given homology basis \,$(A_n,B_n)_{n\in \Z}$\, on \,$\Sigma$\,
in the sense that \,$\int_{A_k} \omega_\ell = \delta_{k,\ell}$\, holds. Because we will need to assess the asymptotic
behavior of the \,$\omega_n$\,, a general statement of existence (like \cite{Feldman/Knoerrer/Trubowitz:2003}, Theorem~3.8, p.~28)
is not enough. Rather we will explicitly construct \,$\omega_n$\, in the form
\,$\omega_n = \tfrac{\Phi_n(\lambda)}{\mu-\mu^{-1}}\,\mathrm{d}\lambda$\,, where \,$\Phi_n$\, is a holomorphic function 
on \,$\C^*$\, which we will construct as a suitable linear combination of infinite products.
This presentation of \,$\Phi_n$\, will give us the asymptotic estimates we need. 

In the present section, we therefore study holomorphic \,$1$-forms on \,$\Sigma$\, and in particular the asymptotic behavior
of holomorphic \,$1$-forms \,$\omega= \tfrac{\Phi(\lambda)}{\mu-\mu^{-1}}\,\mathrm{d}\lambda$\,, where \,$\Phi(\lambda)$\, 
is an infinite product of the kind we are interested in for the construction of the \,$\omega_n$\,. 

More specifically, 
the following Proposition~\ref{P:jacobi:omega} shows that any square-integrable holomorphic 1-form \,$\omega$\, 
is anti-invariant with respect to the hyperelliptic involution and therefore of the form
\,$\omega=\tfrac{\Phi(\lambda)}{\mu-\mu^{-1}}\,\mathrm{d}\lambda$\, with some holomorphic function \,$\Phi: \C^*\to\C$\,. Conversely, the same Proposition gives a necessary condition and a (different) sufficient condition for 
the asymptotic behavior of the function \,$\Phi$\, so that the holomorphic \,$1$-form \,$\tfrac{\Phi(\lambda)}{\mu-\mu^{-1}}\,\mathrm{d}\lambda$\, is square-integrable.

Propositions~\ref{P:jacobi:asymp-Phi}--\ref{P:jacobi:asymp-omega} 
describe a specific construction method for square-integrable holomorphic 1-forms \,$\omega=\tfrac{\Phi(\lambda)}{\mu-\mu^{-1}}\,\mathrm{d}\lambda$\,,
where the function \,$\Phi$\, is given as an infinite product, so that \,$\Phi$\, has one zero in every excluded domain with a single exception;
we in particular describe the asymptotic behavior of such 1-forms. Note that the asymptotic behavior of such 1-forms 
is ``better'' than what can be shown for square-integrable 1-forms in the general case. (See Remark~\ref{R:jacobi:asymp-comparison}.) 

As in the previous section, we fix a holomorphic function \,$\Delta: \C^*\to \C$\, with \,$\Delta-\Delta_0 \in \As(\C^*,\ell^2_{0,0},1)$\,,
and thereby the associated spectral curve \,$\Sigma$\,. We continue to use the notations of Section~\ref{Se:jacobiprep}. 
In particular we let \,$\vkap_{k,\nu}$\, be the zeros of \,$\Delta^2-4$\, as before, interpret \,$\vkap_{k,\nu}$\, also
as points on \,$\Sigma$\,, and put \,$\vkap_{k,*} := \tfrac12(\vkap_{k,1}+\vkap_{k,2})$\,. We also continue to use the possibly punctured 
excluded domains \,$U_{k,\delta}'$\, and \,$\wh{U}_{k,\delta}'$\, defined by Equation~\eqref{eq:jacobiprep:Ukd'}. 

To avoid difficulties, we will suppose from now on: 
\begin{equation}
\label{eq:holo1forms:Delta-hypothesis} 
\begin{boxedminipage}{8.5cm}
\,$\Delta^2-4$\, does not have any zeros of order \,$\geq 3$\,.
\end{boxedminipage}
\end{equation}
Because \,$\Delta^2-4$\, has asymptotically and totally only two zeros in every excluded domain, this condition excludes
only the case where more than two zeros of \,$\Delta^2-4$\, combine in a single point in the ``compact part'' of \,$\Sigma$\,. It is a consequence
of this condition that \,$\Sigma$\, does not have any singularities other than ordinary double points. 

We then enumerate the zeros \,$\vkap_{k,\nu}$\, of \,$\Delta^2-4$\, in such a way that if \,$\vkap\in \C^*$\, is a zero of \,$\Delta^2-4$\,
of order \,$2$\,, there exists \,$k\in \Z$\, with \,$\vkap=\vkap_{k,1}=\vkap_{k,2}$\, (even for \,$|k|$\, small), and we
define
\begin{equation}
\label{eq:holo1forms:S-def} 
S := \Mengegr{k\in \Z}{\vkap_{k,1}=\vkap_{k,2}} \; .
\end{equation}
With this definition, \,$\Menge{\vkap_{k,\nu}}{k\in S}$\, is the set of singular points of \,$\Sigma$\,, and therefore 
\begin{equation}
\label{eq:holo1forms:Sigma'-def} 
\Sigma' := \Sigma \;\setminus\; \Menge{\vkap_{k,\nu}}{k\in S}
\end{equation}
is the regular set of \,$\Sigma$\,, a Riemann surface that is open and dense in \,$\Sigma$\,. 

\label{not:holo1forms:Omega-L2}
We denote the space of holomorphic 1-forms on \,$\Sigma$\, by \,$\Omega(\Sigma)$\, (in the singular points of \,$\Sigma$\,, the holomorphy of \,$\omega\in\Omega(\Sigma)$\, 
is considered in the normalization \,$\wh{\Sigma}$\, of \,$\Sigma$\,), 
and the space of square-integrable 1-forms on \,$\Sigma$\, by \,$L^2(\Sigma,T^*\Sigma)$\,.

The following proposition shows that any square-integrable, holomorphic \,$1$-form \,$\omega$\, on \,$\Sigma$\, is
anti-symmetric with respect to the hyperelliptic involution of \,$\Sigma$\,, and provides a sufficient criterion 
for the square-integrability of a holomorphic \,$1$-form. 

\begin{prop}
\label{P:jacobi:omega}
\begin{enumerate}
\item
Let \,$\omega \in \Omega(\Sigma) \cap L^2(\Sigma,T^*\Sigma)$\,. Then there exists a holomorphic function \,$\Phi: \C^* \to \C$\, so that
\begin{equation}
\label{eq:jacobi:omega:omega}
\omega = \frac{\Phi(\lambda)}{\mu-\mu^{-1}}\,\mathrm{d}\lambda
\end{equation}
holds. 
We have \,$\Phi \in \As(\C^*,\ell^2_{1,-3},1)$\,. 
\item
Let \,$\Phi \in \As(\C^*,\ell^2_{3/2,-5/2},1)$\, be given and put \,$\omega := \tfrac{\Phi}{\mu-\mu^{-1}}\,\mathrm{d}\lambda$\,. Then \,$\omega$\, is a meromorphic \,$1$-form on \,$\Sigma$\, that is holomorphic 
on \,$\Sigma \setminus \Menge{\vkap_{k,*}}{k\in S}$\,, and \,$\omega|\widehat{V}_\delta \in L^2(\wh{V}_\delta,T^*\wh{V}_\delta)$\, holds.
\item
In the setting of (2), suppose that the following condition additionally holds: 
There exists a finite  set \,$T\subset \Z$\, so that  \,$\Phi$\, has a zero \,$\xi_k$\, in \,$U_{k,\delta}$\, for every \,$k\in \Z\setminus T$\,, 
and there exists some \,$C_\xi>0$\, so that \,$|\xi_k-\vkap_{k,*}| \leq C_\xi \cdot |\vkap_{k,1}-\vkap_{k,2}|$\, holds for all \,$k\in \Z \setminus T$\,. 


Then \,$\omega$\, is holomorphic on \,$\Sigma \setminus \Menge{\vkap_{k,*}}{k\in T\cap S}$\,, and square-integrable
on \,$\Sigma \setminus \bigcup_{k\in T\cap S} \wh{U}_{k,\delta}$\,.
%
\end{enumerate}
\end{prop}

\begin{proof}
\emph{For (1).}
Let \,$\sigma: (\lambda,\mu) \mapsto (\lambda,\mu^{-1})$\, be the hyperelliptic involution of \,$\Sigma$\,. 
Then we have \,$\omega = \omega_+ + \omega_-$\, with \,$\omega_\pm := \tfrac12(\omega \pm \sigma^* \omega)$\,. 
\,$\omega_+$\, is \,$\sigma$-invariant, and therefore can be regarded as a square-integrable, holomorphic 1-form \,$f_+(\lambda)\,\mathrm{d}\lambda$\, on \,$\C^*$\,.
Because the holomorphic function \,$f_+$\, is square-integrable on \,$\C^*$\,, it is identically zero, and therefore
\,$\omega=\omega_-$\, is \,$\sigma$-anti-invariant. 

Hence \,$(\mu-\mu^{-1})\cdot \omega$\, is a \,$\sigma$-invariant, holomorphic \,$1$-form, and can thus be regarded
as a holomorphic \,$1$-form on \,$\C^*$\,. Therefore there exists a holomorphic function \,$\Phi: \C^* \to \C$\, with
\,$(\mu-\mu^{-1})\cdot \omega = \Phi(\lambda)\,\mathrm{d}\lambda$\,, and with this \,$\Phi$\, Equation~\eqref{eq:jacobi:omega:omega} holds.


It remains to show \,$\Phi \in \As(\C^*,\ell^2_{1,-3},1)$\,. For this we may suppose without loss of generality that \,$\delta>0$\, is chosen
so small that even the excluded domains \,$U_{k,3\delta}$\, do not overlap. 
Because of Proposition~\ref{P:interpolate:l2asymp}(1) it then suffices to show 
that \,$\Phi|V_{3\delta} \in \As(V_{3\delta},\ell^2_{1,-3},1)$\, holds. For each \,$k\in \Z$\,, let \,$\lambda_k \in S_k \cap V_{3\delta}$\, 
be such that \,$\tfrac{|\Phi(\lambda)|}{w(\lambda)}$\, attains its maximum on the compact set \,$S_k \cap V_{3\delta}$\, at \,$\lambda_k$\,
(where \,$S_k$\, is the annulus defined by Equation~\eqref{eq:As:Sk}). 
Let \,$M_k$\, be the topological annulus
$$ M_k := \begin{cases} \Menge{\lambda\in\C^*}{\delta k \leq |\lambda-\lambda_k|\leq 2\delta k} & \text{for \,$k>0$\,} \\ 
\Menge{\lambda\in\C^*}{\delta k^{-3} \leq |\lambda-\lambda_k|\leq 2\delta k^{-3}} & \text{for \,$k<0$\,} \
\end{cases} $$
for \,$|k|$\, sufficiently large. Then only the \,$M_k$\, of consecutive indices can intersect, and \,$M_k$\, is contained in \,$(S_{k-1}\cup S_k\cup S_{k+1})\cap V_\delta$\,. 

We now first look at \,$k>0$\,. 
Because of the mean value property for the holomorphic function \,$f(\lambda) := \bigr(\tfrac{\Phi}{\mu-\mu^{-1}}\bigr)^2 = \tfrac{\Phi^2}{\Delta^2-4}$\,, 
we have for any \,$r \in [\delta k,2\delta k]$\, 
$$ f(\lambda_k) = \frac{1}{2\pi} \int_{0}^{2\pi} f(\lambda_k + re^{it})\,\mathrm{d}t $$
and therefore
\begin{align*}
f(\lambda_k) & = \frac{1}{\delta k} \int_{r=\delta k}^{2\delta k} \frac{1}{2\pi} \int_{t=0}^{2\pi} f(\lambda_k + re^{it})\,\mathrm{d}t\,\mathrm{d}r \\
& = \frac{1}{2\pi \delta k} \int_{r=\delta k}^{2\delta k} \int_{t=0}^{2\pi} \frac{1}{r}\,f(\lambda_k + re^{it})\,\mathrm{d}t\,r\, \mathrm{d}r \\
& = \frac{1}{2\pi \delta k} \int_{M_k} \frac{1}{r}\,f(\lambda)\,\mathrm{d}^2\lambda \;, 
\end{align*}
where \,$\mathrm{d}^2\lambda$\, here and in the sequel denotes the Lebesgue measure on \,$\C$\,, and whence
$$ |f(\lambda_k)| \leq \frac{1}{2\pi \delta k} \cdot \frac{1}{\delta k} \cdot \int_{M_k} |f(\lambda)|\,\mathrm{d}^2\lambda $$
follows. Because of \,$\omega \in L^2(\Sigma,T^*\Sigma)$\,, we have in particular \,$f|V_\delta \in L^1(V_\delta)$\,. The \,$M_k$\, are contained
in \,$V_\delta$\,, and each point of \,$V_\delta$\, is covered by at most two of the \,$M_k$\,. Therefore we obtain:
$$ \left| f(\lambda_k) \right| \in \frac{1}{k^2} \cdot \ell^1(k) = \ell^1_2(k) $$
and hence
$$ \frac{\Phi(\lambda_k)}{\sqrt{\Delta(\lambda_k)^2-4}} \in \ell^2_1(k) $$ 
holds. Because of \,$\lambda_k \in V_{3\delta}$\,, \,$\sqrt{\Delta(\lambda_k)^2-4}$\, is comparable to \,$w(\lambda_k)$\,, and thus we obtain
$$ \frac{\Phi(\lambda_k)}{w(\lambda_k)} \in \ell^2_1(k) \; . $$
It follows from the definition of \,$\lambda_k$\, that \,$\Phi|V_{3\delta} \in \As_\infty(V_{3\delta},\ell^2_1,1)$\, holds. A similar argument yields \,$\Phi|V_{3\delta} \in \As_0(V_{3\delta},\ell^2_{-3},1)$\,, and thus we have
\,$\Phi|V_{3\delta} \in \As(V_{3\delta},\ell^2_{1,-3},1)$\,. By Proposition~\ref{P:interpolate:l2asymp}(1), \,$\Phi \in \As(\C^*,\ell^2_{1,-3},1)$\, follows.

\emph{For (2).}
It is clear that \,$\omega$\, is a meromorphic \,$1$-form on \,$\Sigma$\,, and that it is holomorphic except at the zeros of \,$\mu-\mu^{-1}$\,, i.e.~at the \,$\vkap_{k,\nu}$\,. If \,$\vkap_{k,\nu}$\, is a regular branch
point of \,$\Sigma$\,, i.e.~a simple zero of \,$\Delta^2-4$\,, then both \,$\mu-\mu^{-1}$\, and \,$\mathrm{d}\lambda$\, have a simple zero at this point, and therefore \,$\omega$\, is holomorphic also at these points.
It follows that \,$\omega$\, is holomorphic on \,$\Sigma \setminus \Menge{\vkap_{k,*}}{k\in S}$\,.

Because of \,$\Phi \in \As(\C^*,\ell^2_{3/2,-5/2},1)$\, and because \,$\mu-\mu^{-1}=\sqrt{\Delta(\lambda)^2-4}$\, is on \,$\widehat{V}_\delta$\, comparable to \,$w(\lambda)$\,, we have
\,$\tfrac{\Phi}{\mu-\mu^{-1}}|\wh{V}_\delta \in \As(\wh{V}_{\delta},\ell^2_{3/2,-5/2},0)$\, and therefore \,$\tfrac{\Phi^2}{\Delta^2-4}|V_\delta \in \As(V_\delta,\ell^1_{3,-5},0)$\,. Because of \,$\vol(S_k\cap V_\delta) \in \ell^\infty_{-3,5}$\,,
it follows that \,$\tfrac{\Phi^2}{\Delta^2-4} \in L^1(V_\delta)$\, and therefore \,$\omega|\wh{V}_\delta \in L^2(\wh{V}_\delta,T^*\wh{V}_\delta)$\, holds.

\emph{For (3).}
The additional condition ensures that \,$\xi_k=\vkap_{k,*}$\, holds for every \,$k\in (\Z\setminus T)\cap S$\,; this fact implies together with
\eqref{eq:holo1forms:Delta-hypothesis} that \,$\omega$\, is holomorphic on \,$\Sigma \setminus \Menge{\vkap_{k,*}}{k\in T\cap S}$\,.
Moreover, these hypotheses ensure that \,$\frac{\Phi^2}{\Delta^2-4}$\, has at most simple poles on every \,$U_{k,\delta}$\, with \,$k \in \Z \setminus (T \cap S)$\,, and therefore \,$\omega$\, is
square-integrable on each individual excluded domain \,$\wh{U}_{k,\delta}$\, with \,$k \in \Z \setminus (T\cap S)$\,. Because \,$\omega$\, is also square-integrable on \,$\wh{V}_\delta$\, by (2), it only remains
to show that the sum of the integrals of \,$\omega^2$\, on \,$\wh{U}_{k,\delta}$\,, \,$k\in \Z \setminus (T\cap S)$\, is finite.

For \,$\lambda\in U_{k,\delta}$\, we have
$$ |\Phi(\lambda)| \leq \max_{U_{k,\delta}} |\Phi'| \cdot |\lambda-\xi_k| \; . $$
We now note that \,$\Phi \in \As(\C^*,\ell^2_{3/2,-5/2},1)$\, implies \,$\Phi' \in \As(\C^*,\ell^2_{5/2,-11/2},1)$\, by Proposition~\ref{P:interpolate:l2asymp}(3), 
and therefore there exists a sequence \,$(a_k) \in \ell^2_{5/2,-11/2}(k)$\, so that for every \,$\lambda \in U_{k,\delta}$\, we have
$$ |\Phi(\lambda)| \leq a_k \cdot |\lambda-\xi_k| \; . $$
We now have by Proposition~\ref{P:jacobiprep:Deltaint-neu}(4)
\begin{align*}
\left| \int_{\wh{U}_{k,\delta}} \omega^2 \right|
& \leq 2\int_{U_{k,\delta}} \frac{|\Phi(\lambda)|^2}{|\Delta(\lambda)^2-4|}\,\mathrm{d}^2\lambda
\leq 2\,a_k^2 \cdot \int_{U_{k,\delta}} \frac{|\lambda-\xi_k|^2}{|\Delta^2-4|}\,\mathrm{d}^2\lambda \\
& \leq C\,a_k^2\,\lambda_{k,0}^\rho\cdot \left( |\xi_k-\vkap_{k,*}|^2 + \frac{|\xi_k-\vkap_{k,*}|^4}{|\vkap_{k,1}-\vkap_{k,2}|^2} + \ell^\infty_{-2,6}(k) \right) \;.
\end{align*}
We have \,$|\xi_k-\vkap_{k,*}|^2 \in \ell^\infty_{-2,6}(k)$\, and \,$\frac{|\xi_k-\vkap_{k,*}|^2}{|\vkap_{k,1}-\vkap_{k,2}|^2}\leq C_\xi^2$\,. It follows that we have
$$ \left| \int_{\wh{U}_{k,\delta}} \omega^2 \right| \leq C\,a_k^2\,\ell^\infty_{-2,6}(k) \cdot \ell^\infty_{-2,6}(k) \in \ell^1_{1,1}(k) \subset \ell^1(k) \; . $$
%
\end{proof}

In the following proposition we investigate holomorphic functions \,$\Phi:\C^*\to\C$\, which are infinite products and have one zero in every excluded domain \,$U_{k,\delta}$\, with a single exception \,$U_{n,\delta}$\,, \,$n\in \Z$\,. 
To show asymptotic estimates for such holomorphic functions, we will apply Proposition~\ref{P:interpolate:lambda}. (By adding one more linear factor corresponding a zero in the ``missing'' excluded domain \,$U_{n,\delta}$\,,
we obtain a function of the type of the holomorphic functions \,$c$\, studied in Proposition~\ref{P:interpolate:lambda}.) One more degree of freedom is provided by multiplying \,$\Phi$\, with 
\,$\lambda^\rho$\, for some \,$\rho\in \Z$\,. We will see in Proposition~\ref{P:jacobi:asymp-f} that there is exactly one choice of \,$\rho$\, so that the associated holomorphic
1-form \,$\tfrac{\Phi(\lambda)}{\mu-\mu^{-1}}\,\mathrm{d}\lambda$\, becomes square-integrable on \,$\Sigma$\, (or on \,$\Sigma \setminus \wh{U}_{n,\delta}$\, if \,$n\in S$\,). In Section~\ref{Se:jacobi}
we will use square-integrable holomorphic 1-forms of this type to construct the canonical basis of \,$\Omega(\Sigma)$\,. 

\newpage

\begin{prop}
\label{P:jacobi:asymp-Phi}
We fix \,$R_0>0$\,. Let 
\,$\rho\in \Z$\,, \,$n\in \Z$\,, and \,$(\xi_k)_{k\in \Z\setminus \{n\}}$\, with \,$\xi_k-\lambda_{k,0} \in \ell^{2}_{-1,3}(\Z \setminus \{n\})$\,, \,$\|\xi_k-\lambda_{k,0}\|_{\ell^{2}_{-1,3}(\Z \setminus \{n\})} \leq R_0$\,
and \,$\xi_k=\vkap_{k,*}$\, for all \,$k\in S\setminus \{n\}$\, be given.

The function
\begin{equation}
\label{eq:jacobi:asymp-Phi:Phi-def}
\Phi_{n,\xi,\rho}(\lambda) := 
\begin{cases}
\lambda^\rho \cdot (\lambda-\xi_0) \cdot \prod_{k\in \N\setminus\{n\}} \frac{\xi_{k}-\lambda}{16\,\pi^2\,k^2} \cdot \prod_{k\in \N} \frac{\lambda-\xi_{-k}}{\lambda} & \text{if \,$n> 0$\,} \\
\lambda^{\rho-1} \cdot (\lambda-\xi_0) \cdot \prod_{k\in \N} \frac{\xi_{k}-\lambda}{16\,\pi^2\,k^2} \cdot \prod_{k\in \N\setminus\{-n\}} \frac{\lambda-\xi_{-k}}{\lambda} & \text{if \,$n\leq 0$\,} 
\end{cases}
\end{equation}
is a holomorphic function on \,$\C^*$\, with the following asymptotic properties:
\begin{enumerate}
\item \emph{(Asymptotic behavior of \,$\Phi_{n,\xi,\rho}$\,.)} We have
\begin{align*}
\Phi_{n,\xi,\rho} - \Phi_{n,0,\rho} & \in \As_\infty(\C^*,\ell^{2}_{-2\rho+1},1) \\
\qmq{and} \Phi_{n,\xi,\rho} - \tau_\xi^{-2}\,\Phi_{n,0,\rho} & \in \As_0(\C^*,\ell^{2}_{2\rho+1},1) 
\end{align*}
with
\begin{equation}
\label{eq:jacobi:asymp-Phi:Phi0-def}
\Phi_{n,0,\rho} := \begin{cases}
4 \cdot \lambda^\rho \cdot \frac{16\pi^2 n^2}{\lambda_{n,0}-\lambda} \cdot c_0(\lambda) & \text{if \,$n>0$\,} \\
4 \cdot \lambda^{\rho} \cdot \frac{1}{\lambda-\lambda_{n,0}} \cdot c_0(\lambda) & \text{if \,$n\leq 0$\,} 
\end{cases} 
\end{equation}
and
$$ \tau_\xi := \left( \prod_{k\in \Z\setminus \{n\}} \frac{\lambda_{k,0}}{\xi_k} \right)^{1/2} \;. $$
In particular \,$\Phi_{n,\xi,\rho} \in \As(\C^*,\ell^\infty_{-2\rho+1,2\rho+1},1)$\, holds.

\item \emph{(Asymptotic comparison of two functions of the form \,$\tfrac{\Phi_{n,\xi,\rho}(\lambda)}{\lambda-\xi_k}$\, on \,$U_{k,\delta}$\,.)}
Let two sequences \,$(\xi_k^{[1]}),(\xi_k^{[2]})$\, of the kind of \,$(\xi_k)$\, above be given. 
We denote the quantities defined above associated to \,$(\xi_k^{[\nu]})$\, by the superscript \,${}^{[\nu]}$\, (for \,$\nu\in \{1,2\}$\,). 

Then there exists a constant \,$C>0$\, (dependent only on \,$R_0$\,) so that we have for every \,$k\in \Z \setminus \{n\}$\, and every \,$\lambda \in U_{k,\delta}$\, if \,$k>0$\, 
$$ \left| \frac{\Phi_{n,\xi,\rho}^{[1]}(\lambda)}{\lambda-\xi_k^{[1]}} - \frac{\Phi_{n,\xi,\rho}^{[2]}(\lambda)}{\lambda-\xi_k^{[2]}} \right| \leq C\,r_k $$
and if \,$k<0$\,
$$ \left| \frac{(\tau^{[1]})^2 \cdot\Phi_{n,\xi,\rho}^{[1]}(\lambda)}{\lambda-\xi_k^{[1]}} - \frac{(\tau^{[2]})^2\cdot \Phi_{n,\xi,\rho}^{[2]}(\lambda)}{\lambda-\xi_k^{[2]}} \right| \leq C\,r_k \; , $$
where the sequence \,$(r_k)\in \ell^{2}_{2-2\rho,2\rho-2}(k)$\, is defined for \,$n>0$\, by
$$ r_k := \begin{cases}
\frac{n^2\,k^{2\rho}}{|n^2-k^2|}\,(a_j * \frac{1}{|j|})_k & \text{if \,$k>0$\,} \\
k^{2-2\rho}\,(a_j * \frac{1}{|j|})_k & \text{if \,$k<0$\,} 
\end{cases} \;, $$
and for \,$n<0$\, by
$$ r_k := \begin{cases}
k^{2\rho-2}\,(a_j * \frac{1}{|j|})_k & \text{if \,$k>0$\,} \\
\frac{n^2\,k^{4-2\rho}}{|n^2-k^2|}\,(a_j * \frac{1}{|j|})_k & \text{if \,$k<0$\,} 
\end{cases} \;, $$
with \,$(a_k) \in \ell^2_{0,0}(k)$\, defined by
\begin{equation}
\label{eq:jacobi:asymp-Phi:ak-def}
a_k := \begin{cases} k^{-1}\,|\xi_k^{[1]}-\xi_k^{[2]}| & \text{if \,$k>0$\,} \\ k^{3}\,|\xi_k^{[1]}-\xi_k^{[2]}| & \text{if \,$k<0$\,} \end{cases} \;.
\end{equation}
We have \,$\|r_k\|_{\ell^{2}_{2-2\rho,2\rho-2}} \leq \|\xi_k^{[1]}-\xi_k^{[2]}\|_{\ell^{2}_{-1,3}}$\,.

\item \emph{(Asymptotic estimate of \,$\tfrac{\Phi_{n,\xi,\rho}(\lambda)}{\lambda-\xi_k}$\, on \,$U_{k,\delta}$\,.)}
Let \,$(\xi_k)$\, be as before and let \,$(r_k)$\, be as in (2) applied to \,$\xi_k^{[1]}=\xi_k$\, and \,$\xi_k^{[2]}=\lambda_{k,0}$\,. 

Then there exist \,$C_1,C_2>0$\, (dependent only on \,$R_0$\,)
such that for every \,$k\in \Z\setminus \{n\}$\, and \,$\lambda \in U_{k,\delta}$\, we have 
\begin{align*}
\text{if \,$n,k>0$\,:} \left| \frac{\Phi_{n,\xi,\rho}(\lambda)}{\lambda-\xi_k} - \frac{8(-1)^k\,\pi^2\,n^2\cdot \xi_k^\rho}{\lambda_{n,0}-\xi_k} \right| & \leq C_1\cdot \frac{n^2 \cdot k^{2\rho-1}}{|n^2-k^2|}\cdot |\lambda-\xi_k| + C_2\cdot r_k \\
\text{if \,$n>0,k<0$\,:} \left| \frac{\Phi_{n,\xi,\rho}(\lambda)}{\lambda-\xi_k} - \left( - \frac{8(-1)^k\,\pi^2\,n^2 \cdot \xi_k^{\rho-1}}{\tau_\xi^2\cdot (\lambda_{n,0}-\xi_k)} \right) \right| & \leq C_1 \cdot k^{5-2\rho}\cdot |\lambda-\xi_k| + C_2\cdot r_k \\
\text{if \,$n<0,k>0$\,:} \left| \frac{\Phi_{n,\xi,\rho}(\lambda)}{\lambda-\xi_k} - \left( - \frac{(-1)^k\cdot \xi_k^\rho}{2\,(\lambda_{n,0}-\xi_k)} \right) \right| & \leq C_1\cdot k^{2\rho-3}\cdot |\lambda-\xi_k|+C_2\cdot r_k  \\
\text{if \,$n,k<0$\,:} \left| \frac{\Phi_{n,\xi,\rho}(\lambda)}{\lambda-\xi_k} - \frac{(-1)^k\cdot \xi_k^{\rho-1}}{2\,\tau_\xi^2\cdot (\lambda_{n,0}-\xi_k)} \right| & \leq C_1\cdot \frac{n^2 \cdot k^{7-2\rho}}{|n^2-k^2|} \cdot |\lambda-\xi_k| + C_2\cdot r_k \; .  
\end{align*}
\end{enumerate}
\end{prop}

\begin{proof}
For the proof, we abbreviate \,$\Phi := \Phi_{n,\xi,\rho}$\, and \,$\Phi_0 := \Phi_{n,0,\rho}$\,.
Moreover, we put \,$\xi_n:=\lambda_{n,0} \in U_{n,\delta}$\,, and thereby obtain a sequence \,$\xi_k \in \ell^{2}_{-1,3}(\Z)$\,.

\emph{For (1).}
By Proposition~\ref{P:interpolate:lambda}
$$ c_\xi(\lambda) := \frac14\,\tau_\xi\,(\lambda-\xi_0) \cdot \prod_{k=1}^\infty \frac{\xi_{k}-\lambda}{16\,\pi^2\,k^2} \cdot \prod_{k=1}^{\infty} \frac{\lambda-\xi_{-k}}{\lambda} $$
defines a holomorphic function on \,$\C^*$\,, and we have
\begin{equation}
\label{eq:jacobi:asymp:Phi-c}
\Phi(\lambda) = 
\begin{cases}
\displaystyle
\frac{4}{\tau_\xi} \cdot \lambda^\rho \cdot \frac{16\pi^2 n^2}{\xi_n-\lambda}\cdot c_\xi(\lambda)
& \text{if \,$n>0$\,} \\
\displaystyle
\frac{4}{\tau_\xi} \cdot \lambda^{\rho} \cdot \frac{1}{\lambda-\xi_n}\cdot c_\xi(\lambda)
& \text{if \,$n<0$\,} 
\end{cases}
\; . 
\end{equation}
Because \,$c_\xi$\, has a zero at \,$\lambda=\xi_n$\,, it follows from this description that
\,$\Phi$\, is a holomorphic function on \,$\C^*$\,. Moreover, we have \,$c_0(\lambda)=\lambda^{1/2}\,\sin(\zeta(\lambda)) \in \As(\C^*,\ell^\infty_{-1,1},1)$\, and by Proposition~\ref{P:interpolate:lambda}
$$ c_\xi - \tau_\xi\,c_0 \in \As_\infty(\C^*,\ell^{2}_{-1},1) \qmq{and} c_\xi - \tau_\xi^{-1}\,c_0 \in \As_0(\C^*,\ell^{2}_1,1) \;, $$
and hence \,$c_\xi \in \As(\C^*,\ell^\infty_{-1,1},1)$\,. Because of Equation~\eqref{eq:jacobi:asymp:Phi-c}, it follows that we have
$$ \Phi - \Phi_0 \in \As_\infty(\C^*,\ell^{2}_{-2\rho+1},1) \qmq{and} \Phi - \tau_\xi^{-2}\,\Phi_0 \in \As_0(\C^*,\ell^{2}_{2\rho+1},1) \; . $$
Because of \,$\Phi_0 \in \As(\C^*,\ell^\infty_{-2\rho+1,2\rho+1},1)$\, it follows
in particular that \,$\Phi \in \As(\C^*,\ell^\infty_{-2\rho+1,2\rho+1},1)$\, holds.

\emph{For (2).}
By Equation~\eqref{eq:jacobi:asymp:Phi-c} we have 
\begin{equation}
\label{eq:jacobi:asymp:Phi12-quot-c}
\frac{\Phi^{[1]}(\lambda)}{\lambda-\xi_k^{[1]}}-\frac{\Phi^{[2]}(\lambda)}{\lambda-\xi_k^{[2]}} = 
\begin{cases}
\displaystyle
4\, \lambda^\rho \cdot \frac{16\pi^2 n^2}{\xi_n-\lambda}\cdot \left( \frac{c_\xi^{[1]}(\lambda)}{\tau_\xi^{[1]}\cdot (\lambda-\xi_k^{[1]})} - \frac{c_\xi^{[2]}(\lambda)}{\tau_\xi^{[2]}\cdot (\lambda-\xi_k^{[2]})} \right)
& \text{if \,$n>0$\,} \\
\displaystyle
4\, \lambda^{\rho} \cdot \frac{1}{\lambda-\xi_n}\cdot \left( \frac{c_\xi^{[1]}(\lambda)}{\tau_\xi^{[1]}\cdot (\lambda-\xi_k^{[1]})} - \frac{c_\xi^{[2]}(\lambda)}{\tau_\xi^{[2]}\cdot (\lambda-\xi_k^{[2]})} \right)
& \text{if \,$n<0$\,} 
\end{cases}
\; . 
\end{equation}
In the case \,$n>0$\, we note that we have for \,$k\in \Z$\, and \,$\lambda\in U_{k,\delta}$\,
\begin{equation}
\label{eq:jacobi:asymp:lambda-estimate}
|\lambda|^\rho \leq \begin{cases} C_3\cdot k^{2\rho} & \text{if \,$k>0$\,} \\ C_3 \cdot k^{-2\rho} & \text{if \,$k<0$\,} \end{cases}
\qmq{and}
\frac{1}{|\xi_n-\lambda|} \leq \begin{cases} C_4 \cdot \dfrac{1}{|n^2-k^2|} & \text{if \,$k>0$\,} \\ C_4 \cdot \dfrac{1}{n^2} & \text{if \,$k<0$\,} \end{cases} \; .
\end{equation}
Moreover, by Corollary~\ref{C:interpolate:cdivlin}(2) there exists \,$C_5>0$\, so that with \,$(a_k)\in \ell^2_{0,0}(k)$\, defined by Equation~\eqref{eq:jacobi:asymp-Phi:ak-def} we have for \,$k>0$\,
\begin{equation}
\label{eq:jacobi:asymp:c-quot-12+}
\left| \frac{c^{[1]}(\lambda)}{\tau^{[1]}\cdot(\lambda-\xi_{k}^{[1]})} - \frac{c^{[2]}(\lambda)}{\tau^{[2]}\cdot(\lambda-\xi_{k}^{[2]})} \right| \leq C_5\cdot \left(a_k * \frac{1}{|k|}\right)
\end{equation}
and for \,$k<0$\,
\begin{equation}
\label{eq:jacobi:asymp:c-quot-12-}
\left| \frac{c^{[1]}(\lambda)}{(\tau^{[1]})^{-1}\cdot(\lambda-\xi_{k}^{[1]})} - \frac{c^{[2]}(\lambda)}{(\tau^{[2]})^{-1}\cdot(\lambda-\xi_{k}^{[2]})} \right| \leq C_5\,k^2\cdot \left(a_k * \frac{1}{|k|}\right) \; .
\end{equation}

In the case of \,$k>0$\, we obtain for \,$\lambda\in U_{k,\delta}$\, by plugging the estimates \eqref{eq:jacobi:asymp:lambda-estimate} and \eqref{eq:jacobi:asymp:c-quot-12+} into Equation~\eqref{eq:jacobi:asymp:Phi12-quot-c}:
$$ \left| \frac{\Phi^{[1]}(\lambda)}{\lambda-\xi_k^{[1]}}-\frac{\Phi^{[2]}(\lambda)}{\lambda-\xi_k^{[2]}} \right| \leq 4 \cdot C_3\,k^{2\rho} \cdot 16\pi^2 n^2 \cdot \frac{C_4}{|n^2-k^2|} \cdot C_5\, \left(a_k * \frac{1}{|k|}\right) \leq C\cdot r_k$$
with \,$C := 16\,\pi^2\,C_3\,C_4\,C_5$\,, whereas for \,$k<0$\, we similarly obtain
$$ \left| (\tau^{[1]})^2 \cdot \frac{\Phi^{[1]}(\lambda)}{\lambda-\xi_k^{[1]}} - \frac{(\tau^{[2]})^2\cdot \Phi^{[2]}(\lambda)}{\lambda-\xi_k^{[2]}} \right| 
\leq 4\cdot C_3\,k^{-2\rho} \cdot 16\pi^2n^2 \cdot C_4\,\frac{1}{n^2} \cdot C_5\,k^2\, \left(a_k * \frac{1}{|k|}\right) \leq C\,r_k \;. $$

The case \,$n<0$\, is handled analogously.

\emph{For (3).}
Suppose that \,$k\in \Z \setminus \{n\}$\, is given. Here we only consider the case \,$n,k>0$\,; the other cases are shown in an analogous way. 
For given \,$\lambda \in U_{k,\delta}$\, we have
\begin{align}
& \left| \frac{\Phi(\lambda)}{\lambda-\xi_k} - \frac{8(-1)^k\,\pi^2\,n^2\cdot \xi_k^\rho}{\lambda_{n,0}-\xi_k} \right| \notag \\
\label{eq:jacobi:asymp:3}
\leq\; & \left| \frac{\Phi(\lambda)}{\lambda-\xi_k} - \Phi'(\xi_k)\right| 
+ \left| \Phi'(\xi_k) - \Phi_0'(\lambda_{k,0}) \right|
+ \left| \Phi_0'(\lambda_{k,0}) - \frac{8(-1)^k\,\pi^2\,n^2\cdot \xi_k^\rho}{\lambda_{n,0}-\xi_k} \right| \; .
\end{align}
Because of Equation~\eqref{eq:jacobi:asymp:Phi-c} we have
$$ \left( \frac{\Phi(\lambda)}{\lambda-\xi_k} \right)' = \frac{4\cdot 16\pi^2 n^2}{\tau_\xi} \cdot \left( \left( \frac{\lambda^\rho}{\xi_n-\lambda} \right)' \cdot \frac{c_\xi(\lambda)}{\lambda-\xi_k}
+ \frac{\lambda^\rho}{\xi_n-\lambda} \cdot \left( \frac{c_\xi(\lambda)}{\lambda-\xi_k} \right)' \right) \; . $$
Because \,$\tfrac{c_\xi(\lambda)}{\lambda-\xi_k}=O(1)$\, and \,$\left(\tfrac{c_\xi(\lambda)}{\lambda-\xi_k}\right)'=O(k^{-1})$\, uniformly on all the excluded domains, it follows that there exists \,$C_6>0$\, with
$$ \left| \left( \frac{\Phi(\lambda)}{\lambda-\xi_k} \right)' \right| \leq C_6\cdot \frac{n^2\cdot k^{2\rho-1}}{|n^2-k^2|} \; . $$
Because \,$\frac{\Phi(\lambda)}{\lambda-\xi_k}$\, is extended holomorphically at \,$\lambda=\xi_k$\, by the value \,$\Phi'(\xi_k)$\,, it follows that 
\begin{equation}
\label{eq:jacobi:asymp:3-1}
\left| \frac{\Phi(\lambda)}{\lambda-\xi_k} - \Phi'(\xi_k)\right| \leq C_6\cdot \frac{n^2\cdot k^{2\rho-1}}{|n^2-k^2|} \cdot |\lambda-\xi_k| 
\end{equation}
holds. 

Moreover, again because the holomorphic function \,$\tfrac{\Phi(\lambda)}{\lambda-\xi_k}$\,
resp.~\,$\tfrac{\Phi_0(\lambda)}{\lambda-\lambda_{k,0}}$\, is extended holomorphically at \,$\lambda=\xi_k$\, resp.~at \,$\lambda=\lambda_{k,0}$\, by the value \,$\Phi'(\xi_k)$\, resp.~\,$\Phi_0'(\lambda_{k,0})$\,,
it follows from (2) that 
\begin{equation}
\label{eq:jacobi:asymp:3-2}
\left| \Phi'(\xi_k) - \Phi_0'(\lambda_{k,0}) \right| \leq r_k
\end{equation}
holds.

Finally, we have by Equation~\eqref{eq:jacobi:asymp-Phi:Phi0-def}
$$ \Phi_{n,0,\rho}'(\lambda_{k,0}) = 4\,\lambda_{k,0}^\rho\cdot \frac{16\pi^2 n^2}{\lambda_{n,0}-\lambda_{k,0}} \cdot \underbrace{c_0'(\lambda_{k,0})}_{=(-1)^k/8}
= \frac{8(-1)^k\,\pi^2\,n^2\,\lambda_{k,0}^\rho}{\lambda_{n,0}-\lambda_{k,0}} $$
and therefore with constants \,$C_7,C_8>0$\,
\begin{equation}
\label{eq:jacobi:asymp:3-3}
\left| \Phi_{n,0,\rho}'(\lambda_{k,0}) - \frac{8(-1)^k\,\pi^2\,n^2\cdot \xi_k^\rho}{\lambda_{n,0}-\xi_k} \right| \leq C_7\cdot \frac{n^2\,k^{2\rho-2}}{|n^2-k^2|}\cdot |\xi_k-\lambda_{k,0}| 
= C_7\cdot \frac{n^2\,k^{2\rho-1}}{|n^2-k^2|}\cdot a_k \leq C_8 \cdot r_k \; . 
\end{equation}

By applying the estimates \eqref{eq:jacobi:asymp:3-1}, \eqref{eq:jacobi:asymp:3-2} and \eqref{eq:jacobi:asymp:3-3} to \eqref{eq:jacobi:asymp:3}, we obtain the claimed statement.
\end{proof}

As explained at the beginning of the present section, we are interested in holomorphic 1-forms \,$\omega = \tfrac{\Phi_{n,\xi,\rho}(\lambda)}{\mu-\mu^{-1}}\,\mathrm{d}\lambda$\,, where \,$\Phi_{n,\xi,\rho}$\, is as in
Proposition~\ref{P:jacobi:asymp-Phi}, and such that \,$\omega$\, is square-integrable on \,$\Sigma$\, (at least away from the singular points of \,$\Sigma$\,). The following proposition addresses the question
when such an \,$\omega$\, is in fact square-integrable. It turns out that a necessary condition for this to be the case is that 
the exponent \,$\rho$\, occurring in the definition of \,$\Phi_{n,\xi,\rho}$\, in Equation~\eqref{eq:jacobi:asymp-Phi:Phi-def} equals \,$\rho=-1$\,. For \,$\omega$\, to be actually square-integrable,
it is further necessary that the zeros \,$\xi_k$\, of \,$\Phi_{n,\xi,\rho}$\, are ``not too far removed'' from the center \,$\vkap_{k,*}$\, of the branch points, in comparison to the size \,$|\vkap_{k,1}-\vkap_{k,2}|$\,
of the ``handle'' on \,$\Sigma$\, defined by \,$\vkap_{k,1}$\, and \,$\vkap_{k,2}$\,. More specifically, the required condition is that there exist a constant \,$C_\xi>0$\, so that 
\,$|\xi_k-\vkap_{k,*}|\leq C_\xi \cdot |\vkap_{k,1}-\vkap_{k,2}|$\, holds for all \,$k \in \Z \setminus \{n\}$\,, compare Proposition~\ref{P:jacobi:omega}(3).

\begin{prop}
\label{P:jacobi:asymp-f}
In the setting of Proposition~\ref{P:jacobi:asymp-Phi},
the function \,$f_{n,\xi,\rho}$\, defined by 
$$ f_{n,\xi,\rho} := \frac{\Phi_{n,\xi,\rho}}{\mu-\mu^{-1}} $$
is a meromorphic function on \,$\Sigma$\, with the following asymptotic properties:
\begin{enumerate}
\item With 
$$ f_{n,0,\rho} := \begin{cases}
-2i\,\lambda^{\rho+1/2} \cdot \frac{16\pi^2 n^2}{\lambda_{n,0}-\lambda}  & \text{if \,$n>0$\,} \\
-2i\,\lambda^{\rho+1/2} \cdot \frac{1}{\lambda-\lambda_{n,0}} & \text{if \,$n<0$\,} 
\end{cases} 
$$
we have
\begin{align*}
(f_{n,\xi,\rho}-f_{n,0,\rho})|\widehat{V}_\delta & \in \As_\infty(\wh{V}_\delta,\ell^{2}_{-2\rho+1},0) \\
\qmq{and} (f_{n,\xi,\rho}-\tau_\xi^{-2}\,f_{n,0,\rho})|\widehat{V}_\delta & \in \As_0(\wh{V}_\delta,\ell^{2}_{2\rho+1},0) \;; 
\end{align*}
in particular \,$f_{n,\xi,\rho}|\widehat{V}_\delta \in \As(\wh{V}_\delta,\ell^\infty_{-2\rho+1,2\rho+1},0)$\,.
\item 
\,$f_{n,\xi,\rho}|\widehat{V}_\delta \in L^2(\widehat{V}_\delta)$\, if and only if \,$\rho=-1$\,. If \,$\rho=-1$\, holds and there exists \,$C_\xi>0$\, so that 
\,$|\xi_k-\vkap_{k,*}|\leq C_\xi \cdot |\vkap_{k,1}-\vkap_{k,2}|$\, holds for all \,$k \in \Z \setminus \{n\}$\,,
then we have \,$f_{n,\xi,\rho} \in L^2(\Sigma)$\, if \,$n\not\in S$\,, and \,$f_{n,\xi,\rho} \in L^2(\Sigma\setminus \wh{U}_{n,\delta})$\, if \,$n\in S$\,. 
%
\end{enumerate}
\end{prop}

\begin{proof}
We continue to use the notations from the proof of Proposition~\ref{P:jacobi:asymp-Phi}, and also abbreviate \,$f := f_{n,\xi,\rho}$\, and \,$f_0 := f_{n,0,\rho}$\,. 

\emph{For (1).}
We have
$$ \mu_0-\mu_0^{-1}=\sqrt{\Delta_0^2-4} = 2i\cdot \sin(\zeta(\lambda)) = 2i\,\lambda^{-1/2}\,c_0(\lambda) $$
and therefore
$$ f_0 = \frac{\Phi_0}{2i\,\lambda^{-1}\,c_0} = \frac{\Phi_0}{\mu_0-\mu_0^{-1}} \;, $$
whence
\begin{align*}
f-f_0 & = \frac{\Phi}{\mu-\mu^{-1}} - \frac{\Phi_0}{\mu_0-\mu_0^{-1}} \\
& = \Phi \cdot \left( \frac{1}{\mu-\mu^{-1}} - \frac{1}{\mu_0-\mu_0^{-1}} \right) + \left(\Phi-\Phi_0\right) \cdot \frac{1}{\mu_0-\mu_0^{-1}} \; . 
\end{align*}
follows. Because of
$$ \left. \left( (\mu-\mu^{-1}) - (\mu_0-\mu_0^{-1}) \right)\right|\widehat{V}_\delta \in \As(\wh{V}_\delta,\ell^{2}_{0,0},1) $$
and Proposition~\ref{P:jacobi:asymp-Phi}(1), we obtain \,$(f-f_0)|\widehat{V}_\delta \in \As_\infty(\wh{V}_\delta,\ell^{2}_{-2\rho+1},0)$\,. 
A similar calculation yields \,$(f-\tau_\xi^{-2}\,f_0)|\widehat{V}_\delta \in \As_0(\wh{V}_\delta,\ell^{2}_{2\rho+1},0)$\,. 


\emph{For (2).}
By Proposition~\ref{P:jacobi:omega}(1), a necessary condition for \,$f$\, to be square-integrable is that \,$\Phi \in \As(\C^*,\ell^2_{1,-3},1)$\, holds. We have \,$\Phi \in \As(\C^*,\ell^\infty_{-2\rho+1,2\rho+1},1)$\, 
by Proposition~\ref{P:jacobi:asymp-Phi}(1) (and no better asymptotic can hold). Because \,$\ell^\infty_{-2\rho+1,2\rho+1} \subset \ell^2_{1,-3}$\, holds
if and only if the inequalities \,$-2\rho+1 > 1 + \tfrac12$\,, i.e.~\,$\rho<-\tfrac14$\, and 
\,$2\rho+1 > -3 + \tfrac12$\,, i.e.~\,$\rho > -\tfrac74$\, holds, we see that \,$f$\, cannot be square-integrable for \,$\rho\neq -1$\,. 

Let us now suppose \,$\rho=-1$\,. Then we have \,$\Phi \in \As(\C^*,\ell^\infty_{3,-1},1)$\, by Proposition~\ref{P:jacobi:asymp-Phi}(1), and therefore in particular \,$\Phi \in \As(\C^*,\ell^2_{3/2,-5/2},1)$\,. 
Therefore \,$f|\wh{V}_\delta \in L^2(\wh{V}_\delta)$\, then holds. 

If moreover there exists \,$C_\xi>0$\, so that \,$|\xi_k-\vkap_{k,*}|\leq C_\xi \cdot |\vkap_{k,1}-\vkap_{k,2}|$\, holds for all \,$k \in \Z \setminus \{n\}$\,,
then we have \,$f\in L^2(\Sigma)$\, if \,$n\not\in S$\,, and \,$f \in L^2(\Sigma\setminus \wh{U}_{n,\delta})$\, if \,$n\in S$\, by Proposition~\ref{P:jacobi:omega}(3)
(applied with \,$T=\{n\}$\,).
\end{proof}

The preceding proposition showed that \,$\omega=\tfrac{\Phi_{n,\xi,\rho}}{\mu-\mu^{-1}}\,\mathrm{d}\lambda$\, is square-integrable only if \,$\rho=-1$\, holds and 
there exists \,$C_\xi>0$\, so that \,$|\xi_k-\vkap_{k,*}|\leq C_\xi \cdot |\vkap_{k,1}-\vkap_{k,2}|$\, holds for all \,$k \in \Z \setminus \{n\}$\,.
For this reason we will subsequently consider \,$\Phi_{n,\xi,\rho}$\, only where these requirements
are satisfied. In the following, final proposition on asymptotic estimates we study the asymptotic behavior of path integrals of holomorphic \,$1$-forms \,$\omega$\, of this type.

For this purpose we introduce the following notation:
Within the space \,$\Div$\, of all (classical) asymptotic divisors (regarded as point multi-sets in \,$\C^*\times \C^*$\,), 
we consider the subspace of those divisors whose support is contained in the spectral curve \,$\Sigma$\, resp.~in the
regular surface \,$\Sigma'$\, obtained from \,$\Sigma$\, by puncturing at the singularities 
(compare Equation~\eqref{eq:holo1forms:Sigma'-def}):
\begin{equation}
\label{eq:jacobi:DivSigma}
\Div(\Sigma) := \Mengegr{D \in \Div}{\mathrm{supp}(D) \subset \Sigma} \qmq{resp.} 
\Div(\Sigma') := \Mengegr{D \in \Div}{\mathrm{supp}(D) \subset \Sigma'} \; . 
\end{equation}
Then \,$\Div(\Sigma)$\, is an (infinite-dimensional) complex subvariety of \,$\Div$\,, and \,$\Div(\Sigma')$\, is an 
open and dense subset of \,$\Div(\Sigma)$\,. We also consider the corresponding open and dense subsets of tame
divisors (see Definition~\ref{D:special:tame}(1)):
\begin{equation}
\label{eq:jacobi:DivSigma-tame}
\Div_{tame}(\Sigma) := \Div(\Sigma) \cap \Div_{tame} \qmq{and} \Div_{tame}(\Sigma') := \Div(\Sigma') \cap \Div_{tame} \; .
\end{equation}

\begin{prop}
\label{P:jacobi:asymp-omega}
Let \,$C_\xi>0$\,. In the present proposition we use the notations of Propositions~\ref{P:jacobi:asymp-Phi} and \ref{P:jacobi:asymp-f}, and consider 1-forms 
\,$\omega = f_{n,\xi,\rho=-1}\,\mathrm{d}\lambda$\,, where \,$n\in \Z$\,, \,$\rho=-1$\, and \,$|\xi_k-\vkap_{k,*}|\leq C_\xi \cdot |\vkap_{k,1}-\vkap_{k,2}|$\, holds for all \,$k \in \Z \setminus \{n\}$\,.
By Proposition~\ref{P:jacobi:asymp-f}(2) any such \,$\omega$\, is square-integrable on \,$\Sigma$\, resp.~on \,$\Sigma \setminus \wh{U}_{n,\delta}$\, for \,$n\not\in S$\, resp.~for \,$n\in S$\,.
If \,$n\not\in S$\,, \,$\omega \in \Omega(\Sigma)$\, holds, whereas for \,$n\in S$\,, we have \,$\omega \in \Omega(\Sigma\setminus \{\vkap_{n,*}\})$\, and \,$\vkap_{n,*}$\, is a regular point of \,$\omega$\,
in the sense of \textsc{Serre} (\cite{Serre:1988}, IV.9, p.~68ff.).


For any path integral occurring in the sequel, we suppose that
the path of integration runs entirely in \,$\Sigma'$\,, 
and that for \,$|k|$\, large, the path runs entirely in \,$\wh{U}_{k,\delta}$\,. 

We let \,$N\in \N$\, be as in Lemma~\ref{L:jacobiprep:Psi-intro}.

\enlargethispage{2em}
\begin{enumerate}
\item
There exists a sequence \,$(b_n)\in \ell^2_{-1,-1}(n)$\, (depending only on \,$\Sigma$\,, \,$C_\xi$\, and \,$\delta$\,) so that for every \,$\omega=f_{n,\xi,\rho=-1}\,\mathrm{d}\lambda$\, of the type considered
in the present proposition with \,$|n|>N$\, and every \,$(\lambda_n^o,\mu_n^o),(\lambda_n,\mu_n) \in \wh{U}_{n,\delta}'$\, we have
\begin{align*}
\text{if \,$n>0$\,:} 
& \left| \int_{(\lambda_n^o,\mu_n^o)}^{(\lambda_n,\mu_n)} \omega - (-2i)\,\lambda_{n,0}^{1/2}\cdot \ln\left( \frac{\lambda_n-\vkap_{n,*}+\Psi_n(\lambda_n,\mu_n)}{\lambda_n^o-\vkap_{n,*}+\Psi_n(\lambda_n^o,\mu_n^o)} \right)\right| \\
& \qquad\qquad\quad \leq b_n \cdot \left| \ln\left( \frac{\lambda_n-\vkap_{n,*}+\Psi_n(\lambda_n,\mu_n)}{\lambda_n^o-\vkap_{n,*}+\Psi_n(\lambda_n^o,\mu_n^o)} \right)\right| \\
\text{if \,$n<0$\,:} 
& \left| \int_{(\lambda_n^o,\mu_n^o)}^{(\lambda_n,\mu_n)} \omega - \frac{2i}{\tau_\xi^2}\cdot \lambda_{n,0}^{-1/2}\cdot \ln\left( \frac{\lambda_n-\vkap_{n,*}+\Psi_n(\lambda_n,\mu_n)}{\lambda_n^o-\vkap_{n,*}+\Psi_n(\lambda_n^o,\mu_n^o)} \right) \right| \\
& \qquad\qquad\quad \leq b_n \cdot \left| \ln\left( \frac{\lambda_n-\vkap_{n,*}+\Psi_n(\lambda_n,\mu_n)}{\lambda_n^o-\vkap_{n,*}+\Psi_n(\lambda_n^o,\mu_n^o)} \right) \right| \; . 
\end{align*}
Here we integrate along any path that is contained in \,$\wh{U}_{n,\delta}$\,, and \,$\ln(z)$\, is the branch of the complex logarithm function with \,$\ln(1)=2\pi i m$\,, 
where \,$m\in \Z$\, is the winding number of the path of integration around the pair of branch points \,$\vkap_{k,1}$\,, \,$\vkap_{k,2}$\,.

\item 
Let \,$R_0>0$\, be given. 
There exists a constant \,$C>0$\,, depending only on \,$\Sigma$\,, \,$C_\xi$\, and \,$R_0$\,, so that
we have for every \,$\omega$\, of the type considered here, every divisor \,$D=\{(\lambda_k,\mu_k)\} \in \Div(\Sigma')$\, with \,$\|\lambda_k-\vkap_{k,1}\|_{\ell^2_{-1,3}} \leq R_0$\, and every \,$k \in \Z \setminus \{n\}$\,
\begin{align*}
\text{if \,$n,k>0$\,:} 
\left| \int_{\vkap_{k,1}}^{(\lambda_k,\mu_k)} \omega \right| & \leq C \cdot \frac{n^2}{|n^2-k^2|}\cdot a_k \\
\text{if \,$n>0$\,, \,$k<0$\,:} 
\left| \int_{\vkap_{k,1}}^{(\lambda_k,\mu_k)} \omega \right| & \leq C \cdot \frac{1}{k^2}\cdot a_k \\
\text{if \,$n<0$\,, \,$k>0$\,:} 
\left| \int_{\vkap_{k,1}}^{(\lambda_k,\mu_k)} \omega \right| & \leq C \cdot \frac{1}{k^2}\cdot a_k \\
\text{if \,$n,k<0$\,:} 
\left| \int_{\vkap_{k,1}}^{(\lambda_k,\mu_k)} \omega \right| & \leq C \cdot \frac{n^2}{|n^2-k^2|}\cdot a_k \; , 
\end{align*}
where the sequence \,$a_k \in \ell^2_{0,0}(k)$\, is given by 
\begin{equation}
\label{eq:jacobi:asymp-omega:ak-def}
a_k := \begin{cases} \tfrac{1}{k} |\lambda_k-\vkap_{k,1}| & \text{for \,$k>0$\,} \\ k^3 \,|\lambda_k-\vkap_{k,1}| & \text{for \,$k<0$\,} \end{cases} \;,
\end{equation}
and where the paths of integration of the above integrals run entirely in \,$\wh{U}_{k,\delta}$\, and do not wind around the pair of branch points \,$\vkap_{k,1}, \vkap_{k,2}$\,. 
\item
For every \,$\omega$\, of the type considered here and every divisor \,$D=\{(\lambda_k,\mu_k)\} \in \Div(\Sigma')$\,,
the sum \,$\sum_{k\in \Z\setminus \{n\}} \int_{\vkap_{k,1}}^{(\lambda_k,\mu_k)} \omega$\, converges absolutely, where the paths of
integration of the integrals are as in (2).
For every \,$R_0>0$\,, and with the associated sequence \,$(a_k) \in  \ell^2_{0,0}(k)$\, as in (2) there exist constants \,$C,C_1>0$\, so that we have if \,$\|\lambda_k-\vkap_{k,1}\|_{\ell^2_{-1,3}} \leq R_0$\, 
$$ \sum_{k\in \Z\setminus \{n\}} \left| \int_{\vkap_{k,1}}^{(\lambda_k,\mu_k)} \omega \right| \leq C\,n\,\left( a_k * \frac{1}{|k|} \right)_n + C_1 \in \ell^2_{-1,-1}(n)\; . $$
\item
Let \,$B$\, be any cycle or path of integration on \,$\Sigma$\, that avoids the branch points and the singularities of \,$\Sigma$\,. Then there exists a constant \,$C_B>0$\, (dependent on \,$\Sigma$\,, \,$B$\, and \,$C_\xi$\,) 
so that we have for every \,$\omega$\, of the type considered here
$$ \left| \int_{B} \omega \right| \leq C_B \;. $$ 
\end{enumerate}
\end{prop}

\begin{proof}
\emph{For (1).}
We again abbreviate \,$\Phi := \Phi_{\xi,n,\rho=-1}$\,, \,$f := f_{\xi,n,\rho=-1}= \tfrac{\Phi}{\mu-\mu^{-1}}$\, and
\,$\omega = f\,\mathrm{d}\lambda$\,, and use the quantities defined in Propositions~\ref{P:jacobi:asymp-Phi} and \ref{P:jacobi:asymp-f}.

By Corollary~\ref{C:interpolate:cdivlin}(1) we have for \,$\lambda\in U_{n,\delta}$\,:
\begin{align}
\text{if \,$n>0$\,:} \left| \frac{\Phi(\lambda)}{\lambda-\lambda_{n,0}} - \frac{(-1)^n\,\tau_\xi}{8} \right| & \leq C_1\,\frac{|\lambda-\lambda_{n,0}|}{n}+ r_n \notag \\
\label{eq:jacobi:asymp:intn-Phiasymp}
\text{if \,$n<0$\,:} \left| \frac{\Phi(\lambda)}{\lambda-\lambda_{n,0}} - \frac{(-1)^n\,\tau_\xi^{-1}\,\lambda_{n,0}^{-1}}{8} \right| & \leq C_1\,|\lambda-\lambda_{n,0}|\,n^5+ r_n 
\end{align}
with a sequence \,$(r_n)\in \ell^2_{0,-2}(n)$\,. One possible choice of \,$(r_n)$\, for given \,$(\xi_k)$\, is
$$ r_n = \frac{1}{8}\,\lambda_{n,0}^{-1} + \begin{cases} C_2 \cdot \left(a_k * \frac{1}{|k|} \right)_n & \text{for \,$n>0$\,} \\  C_2\,n^2 \cdot \left(a_k * \frac{1}{|k|} \right)_n & \text{for \,$n<0$\,}  \end{cases} $$
with
$$ a_k := \begin{cases} k^{-1}\cdot |\xi_k-\lambda_{k,0}| & \text{for \,$k>0$\,} \\ k^3\cdot |\xi_k-\lambda_{k,0}| & \text{for \,$k<0$\,} \end{cases} \; . $$
Because of
$$ |\xi_k-\lambda_{k,0}| \leq |\xi_k-\vkap_{k,*}| + |\vkap_{k,*}-\lambda_{k,0}| \leq C_\xi \cdot |\vkap_{k,1}-\vkap_{k,2}| + |\vkap_{k,*}-\lambda_{k,0}| \;, $$
it follows that the sequence \,$(a_k)$\, and therefore also the sequence \,$(r_n)$\, can be chosen such as to depend only on \,$C_\xi$\, and the \,$\vkap_{k,\nu}$\,. 

Moreover, there exists constants \,$C_3,C_4>0$\, (dependent only on \,$\delta$\,) so that we have for \,$\lambda \in U_{n,\delta}$\,
\begin{align}
\text{if \,$n>0$\,:} \left| \frac{16\,\pi^2\,n^2}{\lambda} - 1 \right| & = \left| \frac{\lambda-16\pi^2 n^2-\lambda}{\lambda} \right| \leq C_3 \cdot \frac{|\lambda-\lambda_{n,0}|}{\lambda_{n,0}} \leq C_4 \cdot \frac{|\lambda-\lambda_{n,0}|}{n} \notag\\
\label{eq:jacobi:asymp:intn-lambda}
\text{if \,$n<0$\,:} \left| \frac{1}{\lambda} - \frac{1}{\lambda_{n,0}} \right| & = \left| \frac{\lambda-\lambda_{n,0}}{\lambda\cdot\lambda_{n,0}} \right|
\leq C_3 \cdot \frac{|\lambda-\lambda_{n,0}|}{\lambda_{n,0}^2} \leq C_4 \cdot |\lambda-\lambda_{n,0}|\cdot n^5 \; . 
\end{align}

By applying the asymptotic estimates \eqref{eq:jacobi:asymp:intn-Phiasymp} and \eqref{eq:jacobi:asymp:intn-lambda} to Equation~\eqref{eq:jacobi:asymp:Phi-c},
we obtain with another constant \,$C_5>0$\,
\begin{align*}
\text{if \,$n>0$\,:} \left| \Phi(\lambda) - \left( -\frac{(-1)^n}{2} \right) \right|
& \leq C_5 \left( \frac{|\lambda-\lambda_{n,0}|}{n} + r_n \right) \\
\text{if \,$n<0$\,:} \left| \Phi(\lambda) - \frac{(-1)^n}{2\,\tau_\xi^2}\cdot \lambda_{n,0}^{-2}  \right|
& \leq C_5 \cdot \left( |\lambda-\lambda_{n,0}|\cdot n^7 + r_n \cdot n^2 \right) \; ,
\end{align*}
and therefore by multiplication with \,$\tfrac{1}{|\mu-\mu^{-1}|}$\,
{\footnotesize
\begin{align*}
\text{if \,$n>0$\,:} \left| \frac{\Phi(\lambda)}{\mu-\mu^{-1}} - \left( -\frac{(-1)^n}{2} \right)\,\frac{1}{\mu-\mu^{-1}} \right|
& \leq C_5 \left( \frac{|\lambda-\lambda_{n,0}|}{n} + r_n \right) \, \frac{1}{|\mu-\mu^{-1}|} \\
\text{if \,$n<0$\,:} \left| \frac{\Phi(\lambda)}{\mu-\mu^{-1}} - \frac{(-1)^n}{2\,\tau_\xi^2}\cdot \lambda_{n,0}^{-2}\cdot \frac{1}{\mu-\mu^{-1}}  \right|
& \leq C_5 \cdot \left( |\lambda-\lambda_{n,0}|\cdot n^7 + r_n \cdot n^2 \right)\,\frac{1}{|\mu-\mu^{-1}|} \; .
\end{align*}
}
By integration along a given path from \,$(\lambda_n^o,\mu_n^o)$\, to \,$(\lambda_n,\mu_n)$\, we obtain
{\footnotesize
\begin{align*}
\text{if \,$n>0$\,:} \left| \int_{(\lambda_n^o,\mu_n^o)}^{(\lambda_n,\mu_n)} \frac{\Phi(\lambda)}{\mu-\mu^{-1}}\,\mathrm{d}\lambda - \left( -\frac{(-1)^n}{2} \right)\,\int_{(\lambda_n^o,\mu_n^o)}^{(\lambda_n,\mu_n)} \frac{1}{\mu-\mu^{-1}}\,\mathrm{d}\lambda \right|
& \leq s_n \cdot \int_{(\lambda_n^o,\mu_n^o)}^{(\lambda_n,\mu_n)} \frac{1}{|\mu-\mu^{-1}|}\,\mathrm{d}\lambda \\
\text{if \,$n<0$\,:} \left| \int_{(\lambda_n^o,\mu_n^o)}^{(\lambda_n,\mu_n)} \frac{\Phi(\lambda)}{\mu-\mu^{-1}}\,\mathrm{d}\lambda - \frac{(-1)^n}{2\,\tau_\xi^2}\cdot \lambda_{n,0}^{-2}\cdot \int_{(\lambda_n^o,\mu_n^o)}^{(\lambda_n,\mu_n)} \frac{1}{\mu-\mu^{-1}}\,\mathrm{d}\lambda  \right|
& \leq s_n \cdot \int_{(\lambda_n^o,\mu_n^o)}^{(\lambda_n,\mu_n)} \frac{1}{|\mu-\mu^{-1}|}\,\mathrm{d}\lambda \; ,
\end{align*}
}
where we define the sequence \,$(s_n) \in \ell^2_{0,-4}(n)$\, by 
$$ s_n := \begin{cases}
C_5 \cdot \left( \frac{1}{n} \max_{\lambda\in U_{n,\delta}} |\lambda-\lambda_{n,0}| + r_n \right) & \text{if \,$n>0$\,} \\
C_5 \cdot \left( n^7 \cdot \max_{\lambda\in U_{n,\delta}} |\lambda-\lambda_{n,0}| + r_n \cdot n^2 \right) & \text{if \,$n<0$\,} 
\end{cases} \; . $$
Let \,$m\in \Z$\, be the winding number of the given path of integration around the pair of branch points \,$\vkap_{n,1}$\,, \,$\vkap_{n,2}$\,. Then we may suppose without loss of generality that the path of integration
is composed of an admissible path (Definition~\ref{D:jacobiprep:admissible}) and \,$|m|$\, circles (with the orientation given by the sign of \,$m$\,) around \,$\vkap_{n,1}$\,, \,$\vkap_{n,2}$\,. These circles are 
also admissible paths, and therefore it follows from Proposition~\ref{P:jacobiprep:Deltaint-neu}(1),(2)  that there exists another sequence \,$(t_k) \in \ell^2_{-1,3}(k)$\,  that depends only on \,$\Sigma$\, 
so that we have for \,$n>0$\, 
\begin{align*}
& \left| \int_{(\lambda_n^o,\mu_n^o)}^{(\lambda_n,\mu_n)} \frac{\Phi(\lambda)}{\mu-\mu^{-1}}\,\mathrm{d}\lambda - \left( -2i \right)\,\lambda_{n,0}^{1/2}\, \ln\left( \frac{\lambda_n-\vkap_{n,*}+\Psi_n(\lambda_n,\mu_n)}{\lambda_n^o-\vkap_{n,*}+\Psi_n(\lambda_n^o,\mu_n^o)} \right) \right| \\
\leq\;& C_6 \cdot (\lambda_{n,0}^{1/2}\,s_n + |\lambda_n-\vkap_{n,1}| + t_n) \cdot \left|\ln\left( \frac{\lambda_n-\vkap_{n,*}+\Psi_n(\lambda_n,\mu_n)}{\lambda_n^o-\vkap_{n,*}+\Psi_n(\lambda_n^o,\mu_n^o)} \right) \right|  \\
\leq\;& b_n \cdot \left|\ln\left( \frac{\lambda_n-\vkap_{n,*}+\Psi_n(\lambda_n,\mu_n)}{\lambda_n^o-\vkap_{n,*}+\Psi_n(\lambda_n^o,\mu_n^o)} \right) \right|  
\end{align*}
and if \,$n<0$\,
\begin{align*}
& \left| \int_{(\lambda_n^o,\mu_n^o)}^{(\lambda_n,\mu_n)} \frac{\Phi(\lambda)}{\mu-\mu^{-1}}\,\mathrm{d}\lambda - \frac{2i}{\tau_\xi^2}\cdot \lambda_{n,0}^{-1/2}\cdot \ln\left( \frac{\lambda_n-\vkap_{n,*}+\Psi_n(\lambda_n,\mu_n)}{\lambda_n^o-\vkap_{n,*}+\Psi_n(\lambda_n^o,\mu_n^o)} \right)  \right| \\
\leq\; & C_6 \cdot (s_n\,\lambda_{n,0}^{3/2} + (|\lambda_n-\vkap_{n,1}| + t_n)\,\lambda_{n,0}^{-2}) \cdot \left| \ln\left( \frac{\lambda_n-\vkap_{n,*}+\Psi_n(\lambda_n,\mu_n)}{\lambda_n^o-\vkap_{n,*}+\Psi_n(\lambda_n^o,\mu_n^o)} \right) \right| \\
\leq\; & b_n \cdot \left| \ln\left( \frac{\lambda_n-\vkap_{n,*}+\Psi_n(\lambda_n,\mu_n)}{\lambda_n^o-\vkap_{n,*}+\Psi_n(\lambda_n^o,\mu_n^o)} \right) \right| \; ,
\end{align*}
where \,$(b_n)\in \ell^2_{-1,-1}(n)$\, is defined by
$$ b_n := \begin{cases} C_6 \cdot (\lambda_{n,0}^{1/2}\,s_n + |\lambda_n-\vkap_{n,1}| + t_n) & \text{for \,$n>0$\,} \\ C_6 \cdot (s_n\,\lambda_{n,0}^{3/2} + (|\lambda_n-\vkap_{n,1}| + t_n)\,\lambda_{n,0}^{-2}) & \text{for \,$n<0$\,} \end{cases} \; . $$
Thus the claimed statement follows.

\emph{For (2).}
To simplify notation, we consider only the case \,$n,k>0$\, in the sequel. The other three cases with respect to the signs of \,$n$\, and \,$k$\, are handled analogously. 
Because of the hypothesis on the path of integration, we may suppose without loss of generality that the path of integration from \,$\vkap_{k,1}$\, to \,$(\lambda_k,\mu_k)$\, is
the lift of the straight line \,$[\vkap_{k,1},\lambda_k]$\, in \,$\C^*$\,; this is in particular an admissible path in the sense of Definition~\ref{D:jacobiprep:admissible}. 

By Proposition~\ref{P:jacobi:asymp-Phi}(3) there exist constants \,$C_1,C_2>0$\, (depending only on \,$C_\xi$\, and \,$R_0$\,) so that we have
for \,$\lambda \in [\vkap_{k,1},\lambda_k]$\, 
$$ \left| \frac{\Phi(\lambda)}{\lambda-\xi_k} - \frac{8(-1)^k\,\pi^2\,n^2}{\xi_k \cdot (\lambda_{n,0}-\xi_k)} \right| \leq 
\frac{n^2}{k^2 \cdot |n^2-k^2|}\cdot \left( C_1\cdot \frac{|\lambda-\xi_k|}{k} + C_2 \cdot \left( b_j * \frac{1}{|j|} \right)_k \right) $$
with
$$ b_j := \begin{cases} j^{-1}\,|\xi_j-\lambda_{j,0}| & \text{if \,$j>0$\,} \\ j^3\,|\xi_j-\lambda_{j,0}| & \text{if \,$j<0$\,} \end{cases} \; . $$
We have
\begin{align*}
|\lambda-\xi_k| & \leq |\lambda-\vkap_{k,1}| + |\vkap_{k,1}-\vkap_{k,*}| + |\vkap_{k,*}-\xi_k| \\
& \leq R_0\cdot k + \tfrac12 |\vkap_{k,1}-\vkap_{k,2}| + C_\xi\,|\vkap_{k,1}-\vkap_{k,2}| \\
& \leq R_0\cdot k + (C_\xi + \tfrac12)\cdot |\vkap_{k,1}-\vkap_{k,2}| 
\end{align*}
and
$$ |\xi_j-\lambda_{k,0}| \leq |\xi_j-\vkap_{k,*}|+|\vkap_{k,*}-\lambda_{k,0}| \leq C_\xi\cdot |\vkap_{k,1}-\vkap_{k,2}| + |\vkap_{k,*}-\lambda_{k,0}| \; . $$
Therefore it follows that there exists a constant \,$C_3>0$\, (dependent only on \,$R_0$\,, \,$C_\xi$\, and the \,$\vkap_{k,\nu}$\,) so that we have
$$ \left| \frac{\Phi(\lambda)}{\lambda-\xi_k}\right| \leq C_3\cdot \frac{n^2}{k^2 \cdot |n^2-k^2|} \; . $$

By multiplying this inequality with \,$\tfrac{\lambda-\xi_k}{\mu-\mu^{-1}}$\,, we obtain:
$$ \left| \frac{\Phi(\lambda)}{\mu-\mu^{-1}}\right| \leq C_3\cdot \frac{n^2}{k^2 \cdot |n^2-k^2|} \cdot \frac{|\lambda-\xi_k|}{|\mu-\mu^{-1}|} \; . $$
By integration along the straight line from \,$\vkap_{k,1}$\, to \,$(\lambda_k,\mu_k)$\, this yields
$$ \left| \int_{\vkap_{k,1}}^{(\lambda_k,\mu_k)} \frac{\Phi(\lambda)}{\mu-\mu^{-1}}\,\mathrm{d}\lambda \right| \leq C_3\cdot \frac{n^2}{k^2 \cdot |n^2-k^2|} \cdot \int_{\vkap_{k,1}}^{(\lambda_k,\mu_k)} \frac{|\lambda-\xi_k|}{|\mu-\mu^{-1}|} \,|\mathrm{d}\lambda| \; . $$
Because \,$\left| \tfrac{\xi_k-\vkap_{k,*}}{\vkap_{k,1}-\vkap_{k,2}} \right|\leq C_\xi$\, holds, 
Proposition~\ref{P:jacobiprep:Deltaint-neu}(3) implies that there exists a constant \,$C_4>0$\, (depending only on \,$C_\xi$\,) so that we have
$$ \int_{\vkap_{k,1}}^{(\lambda_k,\mu_k)} \frac{|\lambda-\xi_k|}{|\mu-\mu^{-1}|}\,|\mathrm{d}\lambda| \leq C_4 \cdot k \cdot |\lambda_k-\vkap_{k,1}| \;, $$
and from this estimate we obtain with \,$C_5 := C_3 \cdot C_4$\,
$$ \left| \int_{\vkap_{k,1}}^{(\lambda_k,\mu_k)} \frac{\Phi(\lambda)}{\mu-\mu^{-1}}\,\mathrm{d}\lambda \right| \leq C_5\cdot \frac{n^2}{k \cdot |n^2-k^2|} \cdot |\lambda_k-\vkap_{k,1}| \; . $$
In consideration of the definition of \,$(a_k)$\,, this gives the claimed estimate.

\emph{For (3).}
Let us consider the case \,$n>0$\,. Then we have by (2) with the quantities \,$C$\, and \,$(a_k)$\, defined there:
\begin{align*}
\sum_{k\in \Z\setminus \{n\}} \left| \int_{\vkap_{k,1}}^{(\lambda_k,\mu_k)} \omega \right|
& = \sum_{\substack{k\geq 0 \\ k\neq n}} \left| \int_{\vkap_{k,1}}^{(\lambda_k,\mu_k)} \omega \right| + \sum_{k<0} \left| \int_{\vkap_{k,1}}^{(\lambda_k,\mu_k)} \omega \right| \\
& \leq \sum_{\substack{k\geq 0 \\ k\neq n}} C\,\frac{n^2}{|n^2-k^2|}\, a_k + \sum_{k<0} \frac{C}{k^2}\, a_k \\
& = C\,n^2 \sum_{\substack{k\geq 0 \\ k\neq n}} \underbrace{\frac{1}{n+k}}_{\leq 1/n} \, a_k\, \frac{1}{|n-k|} + C\,\sum_{k<0} \frac{1}{k^2}\, a_k \\
& \leq C\,n\,\left( a_k * \frac{1}{|k|} \right)_n + C\,\left\| \frac{1}{k^2} \right\|_{\ell^2} \cdot\left\| a_k \right\|_{\ell^2} \\
& \leq C\,n\,\left( a_k * \frac{1}{|k|} \right)_n + C_6
\end{align*}
with a constant \,$C_6>0$\,. The case \,$n<0$\, is treated analogously.

 
\emph{For (4).}
We let \,$c_\xi$\, be as in the proof of Proposition~\ref{P:jacobi:asymp-Phi}(1). 
Then there exists a constant \,$C_7>0$\, (dependent only on \,$B$\, and \,$C_\xi$\,) so that 
$$ |c_\xi(\lambda)| \leq C_7 \cdot |c_0(\lambda)| $$
holds for all \,$(\lambda,\mu)$\, on the trace of \,$B$\,. 
By Equation~\eqref{eq:jacobi:asymp:Phi-c} it follows that there exists a constant \,$C_8>0$\, so that 
$$ |\Phi(\lambda)| \leq C_8\cdot |c_\xi(\lambda)| \leq C_8 \cdot C_7 \cdot |c_0(\lambda)| $$
holds for all \,$(\lambda,\mu)$\, on the trace of \,$B$\,. We thus have
$$ \left| \int_{B} \omega \right| \leq \int_{B} \frac{|\Phi(\lambda)|}{|\mu-\mu^{-1}|}\,|\mathrm{d}\lambda| \leq C_8\cdot C_7 \cdot C_9 \cdot C_{10} $$
with
$$ C_9 := \max_{\lambda\in \tr(B)} |c_0(\lambda)| \qmq{and} C_{10} := \int_{B} \frac{1}{|\mu-\mu^{-1}|}\,\mathrm{d}\lambda \;. $$
Both \,$C_9$\, and \,$C_{10}$\, are finite (positive) numbers that depend only on \,$B$\,, because the trace of \,$B$\, is compact and avoids the zeros of \,$\mu-\mu^{-1}$\,. 
This yields the claimed estimate with \,$C_B := C_7 \cdot C_8 \cdot C_9 \cdot C_{10}$\,. 
\end{proof}

\begin{rem}
\label{R:jacobi:asymp-comparison}
In the setting of the preceding propositions, suppose that \,$\tfrac{\Phi_{n,\xi,\rho}}{\mu-\mu^{-1}}\,\mathrm{d}\lambda$\,
is square-integrable, and therefore \,$\rho=-1$\, holds by Proposition~\ref{P:jacobi:asymp-f}(2). By Proposition~\ref{P:jacobi:asymp-Phi}(1) we then obtain
\,$\Phi_{n,\xi,\rho} \in \As(\C^*,\ell^\infty_{3,-1},1)$\,. Note that this is an improved statement for 1-forms constructed in the way of Propositions~\ref{P:jacobi:asymp-Phi}--\ref{P:jacobi:asymp-omega}
over the general statement  \,$\Phi \in \As(\C^*,\ell^2_{1,-3},1)$\, in Proposition~\ref{P:jacobi:omega}(1), because of \,$\ell^\infty_{3,-1}\subsetneq \ell^2_{1,-3}$\,.
\end{rem}

\section{Construction of the Jacobi variety for the spectral curve}
\label{Se:jacobi}

In the present section, we will construct Jacobi coordinates for the spectral curve \,$\Sigma$\,, which can have singularities, is non-compact, and has always infinite arithmetic genus and generally infinite geometric genus.
For the construction of the Jacobi variety, the estimates proven in the preceding two sections will play an important role.

\label{not:jacobi:homology-basis}
As an introduction,
we review the construction of the Jacobi variety for \emph{compact} Riemann surfaces 
(see for example \cite{Farkas/Kra:1992}, Section~III.6): 
Let \,$X$\, be a compact Riemann surface, say of genus \,$g \geq 1$\,, and let \,$(A_k,B_k)_{k=1,\dotsc,g}$\,
be a canonical homology basis of \,$X$\,, i.e.~\,$(A_k,B_k)$\, is a basis of the homology group \,$H_1(X,\Z)$\, 
with the intersection properties \,$A_k\times B_\ell = \delta_{k\ell}$\, (Kronecker delta), \,$A_k\times A_\ell = 0 = B_k \times B_\ell$\,
for \,$k,\ell\in \{1,\dotsc,g\}$\,. 
Then there exists a canonical basis \,$(\omega_k)_{k=1,\dotsc,g}$\, of the vector space \,$\Omega(X)$\, of holomorphic 1-forms on \,$X$\, 
that is dual to \,$(A_k)$\, in the sense that \,$\int_{A_k} \omega_\ell = \delta_{k\ell}$\, holds. To any given positive divisor \,$D=\{P_1,\dotsc,P_g\}$\, of degree \,$g$\, on \,$X$\,, we then associate
the quantity
$$ \wt{\vi}(D) := \left( \sum_{k=1}^g \int_{P_0}^{P_k} \omega_\ell \right)_{\ell=1,\dotsc,g} \;\in\; \C^g\;, $$
where \,$P_0\in X$\, is the ``origin point'', which we hold fixed. 
Because these integrals depend on the homology class of the paths of integration from \,$P_0$\, to \,$P_k$\, we choose, the quantity \,$\wt{\vi}(D)$\, is only defined modulo the \emph{period lattice}
$$ \Gamma := \left\langle \left( \int_{A_k} \omega_\ell \right)_{\ell=1,\dotsc,g} \;,\; \left( \int_{B_k} \omega_\ell \right)_{\ell=1,\dotsc,g} \right\rangle_{\Z} \subset \C^g \; . $$
Thus we obtain the \emph{Jacobi variety} \,$\Jac(X) := \C^g / \Gamma$\, of \,$X$\, and (by projecting the values of \,$\wt{\vi}$\, onto \,$\Jac(X)$\,) the \emph{Abel map} 
$$ \vi: \mathrm{Div}_g(X) \to \Jac(X) \; , $$
where \,$\mathrm{Div}_g(X)$\, denotes the space of positive divisors of degree \,$g$\, on \,$X$\,. 

To obtain a Jacobi variety and an Abel map for the spectral curve \,$\Sigma$\, by this strategy, we need to generalize the construction described above in two directions: 
First we need to deal with the fact that \,$\Sigma$\, is not compact and its homology group is generally infinite dimensional, and second we need to handle the fact that \,$\Sigma$\, can have singularities.

As consequence of the fact that the homology group of \,$\Sigma$\, is infinite dimensional, the space of square-integrable, holomorphic%
\footnote{Again the concepts of square-integrability and holomorphy need to be modified near singularities of \,$\Sigma$\,.}
1-forms on \,$\Sigma$\, is also infinite dimensional. Therefore the space \,$\C^g$\, occurring in the treatment of the compact case as the universal cover (or the tangent space) of the Jacobi variety needs to be
replaced by a suitable Banach space adapted to \,$\Sigma$\,, the period lattice \,$\Gamma$\, will also be infinite dimensional, the positive divisors of genus \,$g$\, will likewise be replaced
by positive divisors of infinitely many points, and the sum of integrals defining the Abel map will be a sum of infinitely many terms. To make sure that the latter sum converges absolutely, 
we need to impose an asymptotic condition on the divisors we consider, and that condition of course is the asymptotic condition for the space \,$\Div$\,, see Definition~\ref{D:asympdiv:asympdiv}. 
For such divisors, we will use the estimates of the preceding two sections to show that the Abel map, which is now defined for \,$D=\{(\lambda_k,\mu_k)\} \in \Div$\, via
$$ \wt{\vi}(D) := \left( \sum_{k\in \Z} \int_{(\lambda_k^o,\mu_k^o)}^{(\lambda_k,\mu_k)} \omega_\ell \right)_{\ell\in \Z} $$
with fixed points \,$(\lambda_k^o,\mu_k^o)\in \wh{U}_{k,\delta}$\,, is indeed convergent and maps into the right Banach space. 

Regarding singularities of \,$\Sigma$\, there is the problem that if \,$(\lambda_*,\mu_*)$\, is a singular point of \,$\Sigma$\, and \,$A_k$\, is a member of the homology basis of \,$\Sigma$\,
that encircles this point, then the corresponding member \,$\omega_k$\, of the canonical basis of the 1-forms on \,$\Sigma$\, will not be holomorphic at \,$(\lambda_*,\mu_*)$\,, but only regular
in the sense of \textsc{Serre} (\cite{Serre:1988}, IV.9, p.~68ff.). Moreover, \,$\omega_k$\, will not be square-integrable near \,$(\lambda_*,\mu_*)$\,. To handle this phenomenon, we follow
the approach of \textsc{Rosenlicht} in constructing generalized Jacobi varieties (see \cite{Rosenlicht:1954}, or Chapter~V of \cite{Serre:1988}) in admitting only divisors whose support is
contained in the regular set \,$\Sigma'$\, of \,$\Sigma$\,; in this sense we will in fact construct the Jacobi variety of \,$\Sigma'$\,. 
(Rosenlicht considered a more general setting in that he constructed a family of generalized Jacobi maps dependent on 
a modulus \,$\mathfrak{m}$\, on the singular set of \,$\Sigma$\,, but for our purposes the case \,$\mathfrak{m}=0$\, will suffice.) Considering only divisors whose support avoids the singular 
points of \,$\Sigma$\, is actually sufficient for our purposes: The main reason why we are interested in Jacobi coordinates for the spectral curve \,$\Sigma$\, is to describe the 
motion of spectral divisors under translation of the potential (this motion turns out to be linear in the Jacobi coordinates), and because divisor points that are singular points of \,$\Sigma$\,
remain stationary under translation, it suffices to consider the regular points of the divisor in this context.

\medskip

We continue to use the notations of the preceding two sections. In particular \,$\Sigma$\, is a spectral curve defined via Equation~\eqref{eq:spectral:Sigma2}
by a holomorphic function \,$\Delta:\C^*\to\C$\, with \,$\Delta-\Delta_0 \in \As(\C^*,\ell^2_{0,0},1)$\,,
we denote the zeros of \,$\Delta^2-4$\, by \,$\vkap_{k,\nu}$\, as before, interpret \,$\vkap_{k,\nu}$\, also as points on \,$\Sigma$\, and put \,$\vkap_{k,*}:=\tfrac12(\vkap_{k,1}+\vkap_{k,2})$\,. We continue to
exclude singularities of higher order on \,$\Sigma$\, as we did in Section~\ref{Se:holo1forms} by \eqref{eq:holo1forms:Delta-hypothesis}, i.e.~we require that
\begin{center}
\,$\Delta^2-4$\, does not have any zeros of order \,$\geq 3$\,.
\end{center}
As in Section~17 (see before Equation~\eqref{eq:holo1forms:S-def}) we then suppose that the zeros \,$\vkap_{k,\nu}$\, 
of \,$\Delta^2-4$\, (corresponding to the branch points and singularities of \,$\Sigma$\,) are numbered in such a way
that if \,$\vkap$\, is a zero of order \,$2$\, of \,$\Delta^2-4$\,, then we have \,$\vkap_{k,1}=\vkap_{k,2}=\vkap$\, for 
some \,$k\in \Z$\,. 

Moreover we again consider the possibly punctured
excluded domains \,$U_{k,\delta}'$\, and \,$\wh{U}_{k,\delta}'$\, defined by Equation~\eqref{eq:jacobiprep:Ukd'}, 
the set \,$S \subset \Z$\, of indices for which \,$\vkap_{k,1}=\vkap_{k,2}$\, is a double point of \,$\Sigma$\,
(Equation~\eqref{eq:holo1forms:S-def}) and the regular set \,$\Sigma' = \Sigma\setminus \Menge{\vkap_{k,\nu}}{k\in S}$\, of \,$\Sigma$\,. 
We also consider the spaces \,$\Div(\Sigma)$\, and \,$\Div(\Sigma')$\, 
of (classical) asymptotic divisors with support in \,$\Sigma$\, resp.~in \,$\Sigma'$\,. 



We now fix a basis of the homology of \,$\Sigma'$\,. For every \,$k\in \Z$\, there is a non-trivial cycle \,$A_k$\, in \,$\Sigma'$\, encircling the pair of points \,$(\vkap_{k,1},\vkap_{k,2})$\,. Note that (unlike in \,$\Sigma$\,)
this is true regardless of whether \,$k\not\in S$\, (then \,$\vkap_{k,1},\vkap_{k,2}$\, is a pair of branch points) or \,$k\in S$\, (then \,$\vkap_{k,1}=\vkap_{k,2}$\, is a puncture of \,$\Sigma'$\,) holds.
For \,$k\neq 0$\,, there exists only for \,$k\not\in S$\, another cycle \,$B_k$\, that encircles the pair of points \,$(\vkap_{k,1},\vkap_{k=0,1})$\,. A final cycle \,$B_0$\, comes from the observation
that \,$\sqrt{\lambda}$\, is a global parameter on \,$\wh{V}_\delta$\,. Because
the Riemann surface associated to \,$\sqrt{\lambda}$\, has branch points in \,$\lambda=0$\, and \,$\lambda=\infty$\,, we see that there is another non-trivial cycle \,$B_0$\, in \,$\wh{V}_\delta \subset \Sigma'$\, 
encircling these two branch points. Of these cycles we require that their intersection numbers satisfy 
\,$A_k\times A_\ell = 0 = B_k\times B_\ell$\, and \,$A_k\times B_\ell = \delta_{k\ell}$\, (Kronecker delta) for all \,$k,\ell$\,, then we call \,$(A_k,B_k)$\, a canonical basis of the homology of \,$\Sigma'$\,.
We can choose the cycles such that for \,$|k|$\, large, \,$A_k$\, runs entirely in \,$\wh{U}_{k,\delta}'$\,, and \,$B_k$\, runs in \,$\wh{U}_{k,\delta}' \cup \wh{U}_{k=0,\delta}' \cup \wh{V}_\delta$\,.

We note that with the canonical basis \,$(A_k,B_k)$\, defined in this way, we have 
for any 1-form \,$\omega \in \Omega(\Sigma)$\, and any \,$k\in \Z\setminus S$\, 
\begin{equation}
\label{eq:jacobi:intAk}
\int_{A_k} \omega = 2\cdot \int_{\vkap_{k,1}}^{\vkap_{k,2}} \omega \;,
\end{equation}
where on the right hand side the path of integration is the lift of the straight line \,$[\vkap_{k,1},\vkap_{k,2}]$\, to \,$\Sigma$\,. 

In the case \,$k\in S$\,, we have an analogue to the preceding equation: Consider a meromorphic \,$1$-form \,$\omega$\, on \,$\Sigma$\, that is regular (in the sense of \textsc{Serre}) at \,$\vkap_{k,*}$\,;
the latter means that regarded as a meromorphic \,$1$-form on the normalization \,$\wh{\Sigma}$\, of \,$\Sigma$\,, the pole order of \,$\omega$\, at each of the two points above \,$\vkap_{k,*}$\, in \,$\wh{\Sigma}$\,
is at most \,$1$\, (the order of the zero of \,$\mu-\mu^{-1}$\, there), and that the sum of the residues of \,$\omega$\, in these points is zero, see \cite{Serre:1988}, IV.9, p.~68ff.. 
Note that for \,$n\in S$\, the 1-form \,$\omega$\, constructed in Proposition~\ref{P:jacobi:asymp-omega} is regular at \,$\vkap_{n,*}$\, and holomorphic on \,$\Sigma \setminus \{\vkap_{n,*}\}$\,. 
If any meromorphic 1-form \,$\omega$\, is regular at \,$\vkap_{k,*}$\,, then we have
\begin{equation}
\label{eq:jacobi:intAk-singular}
\int_{A_k} \omega = 0 \;\Longleftrightarrow\; \text{\,$\omega$\, is holomorphic at \,$\vkap_{k,*}$\,} \;,
\end{equation}
where for the question of holomorphy at \,$\vkap_{k,*}$\,, \,$\omega$\, is again regarded as a \,$1$-form on the normalization \,$\wh{\Sigma}$\,. 

Our first task on the way to Jacobi coordinates for \,$\Sigma'$\,
is to obtain a canonical basis \,$(\omega_n)_{n\in \Z}$\, of holomorphic 1-forms \,$\omega_n \in \Omega(\Sigma')$\,
that is dual to the basis of the homology \,$(A_k,B_k)$\, in the sense that for all \,$k,n\in \Z$\,
$$ \int_{A_k} \omega_n = \delta_{kn} $$
holds. It will turn out that \,$\omega_n$\, is for \,$n\in \Z\setminus S$\, holomorphic on \,$\Sigma$\,, and for \,$n\in S$\,, \,$\omega_n$\, is holomorphic on \,$\Sigma \setminus \{\vkap_{n,*}\}$\, and
regular in \,$\vkap_{n,*}$\,. 
For \,$n\not\in S$\, we will have \,$\omega_n \in L^2(\Sigma,T^*\Sigma)$\,, whereas for \,$n\in S$\,, we will
have \,$\omega_n \in L^2(\Sigma \setminus \wh{U}_{n,\delta},T^*(\Sigma\setminus \wh{U}_{n,\delta}))$\, for every \,$\delta>0$\,. 

In the case \,$S=\varnothing$\,, 
the existence of such a basis follows from general results on open
Riemann surfaces (of infinite genus) which are parabolic in the sense of \textsc{Ahlfors} and \textsc{Nevanlinna}, as they are described in Chapter~1 of \cite{Feldman/Knoerrer/Trubowitz:2003}.
Indeed \,$\Sigma$\, is then parabolic (see the proof of Proposition~\ref{P:jacobi:abstractomega-neu} below), so the existence of the basis \,$(\omega_n)$\, follows from \cite{Feldman/Knoerrer/Trubowitz:2003}, Theorem~3.8, p.~28.
However, this argument does not apply to the case \,$S\neq\varnothing$\,, where \,$\Sigma$\, is no longer smooth. 
Moreover it turns out that even for \,$S=\varnothing$\, the information on the asymptotic
behavior of \,$\omega_k$\, we obtain from this general result (by Proposition~\ref{P:jacobi:omega}(1)) is not sufficient to prove that the infinite sum defining
the Jacobi coordinates converges. Therefore we need to carry out an explicit construction of the \,$\omega_n$\,. 


As preparation we need one consequence of the general results on parabolic Riemann surfaces:

\begin{prop}
\label{P:jacobi:abstractomega-neu}
Suppose that \,$\omega_1,\omega_2 \in \Omega(\Sigma) \cap L^2(\Sigma,T^*\Sigma)$\, are given so that
$$ \int_{A_k} \omega_1 = \int_{A_k} \omega_2 $$
holds for all \,$k\in \Z$\,. Then we have \,$\omega_1=\omega_2$\,. 
\end{prop}


\begin{rem}
This proposition generalizes to 1-forms that are holomorphic on \,$\Sigma$\, with the exception of finitely many poles in singular points of \,$\Sigma$\,, that are regular in these poles, and that are 
square-integrable on \,$\Sigma\setminus \bigcup_{k\in I} \wh{U}_{k,\delta}$\, (where \,$I \subset S$\, is the finite set of indices \,$k\in S$\, for which \,$\wh{U}_{k,\delta}$\, contains a pole of the 1-form). The reason is that 
if \,$\omega_1, \omega_2$\, are 1-forms of this kind with \,$\int_{A_k}(\omega_1-\omega_2)=0$\, for all \,$k\in \Z$\,, then this hypothesis implies that \,$\omega_1-\omega_2$\, is holomorphic on all of \,$\Sigma$\,
by \eqref{eq:jacobi:intAk-singular}; because \,$\bigcup_{k\in I} \wh{U}_{k,\delta}$\, is relatively compact in \,$\Sigma$\,, this also means that \,$\omega_1-\omega_2$\, is square-integrable on \,$\Sigma$\,. 
Proposition~\ref{P:jacobi:abstractomega-neu} then shows that \,$\omega_1-\omega_2=0$\, holds.

Note that this generalization is applicable to the members \,$\omega_n$\, of the canonical basis constructed in Theorem~\ref{T:jacobi:canonical} for \,$n\in S$\,,
because they have a regular pole in \,$\vkap_{n,*}$\, and are otherwise holomorphic on \,$\Sigma$\,, and they are square-integrable on \,$\Sigma \setminus \wh{U}_{n,\delta}$\,. 
\end{rem}

\begin{proof}[Proof of Proposition~\ref{P:jacobi:abstractomega-neu}.]
We first note that the normalization \,$\wh{\Sigma}$\, of \,$\Sigma$\, is parabolic in the sense of \textsc{Ahlfors} and \textsc{Nevanlinna}; this means that
\,$\wh{\Sigma}$\, has a harmonic exhaustion function \,$h$\,, i.e.~\,$h$\, is a continuous, proper, non-negative function on \,$\wh{\Sigma}$\,
that is harmonic on the complement of a compact subset (see \cite{Feldman/Knoerrer/Trubowitz:2003}, Definition~3.1, p.~25).
Indeed, \,$h := \bigr| \log(|\lambda|) \bigr|$\, is such a function on \,$\wh{\Sigma}$\,, so \,$\wh{\Sigma}$\, is parabolic. 

It follows by \cite{Feldman/Knoerrer/Trubowitz:2003}, Proposition~3.6, p.~27 that \,$\wh{\Sigma}$\, then also has an exhaustion function with finite charge, 
whence it follows by \cite{Feldman/Knoerrer/Trubowitz:2003}, Proposition~2.10, p.~22 that for any square-integrable, holomorphic \,$1$-form \,$\wh{\omega}$\, on \,$\wh{\Sigma}$\,,
the condition \,$\int_{A_k} \wh{\omega} = 0$\, for all \,$k\in \Z$\, implies \,$\wh{\omega}=0$\,. Because \,$\wh{\Sigma}\to \Sigma$\, is a one-sheeted
covering, we have the analogous statement for \,$\Sigma$\,: If \,$\omega$\, is a square-integrable, holomorphic \,$1$-form on \,$\Sigma$\, (where for the question of holomorphy
we again consider \,$\omega$\, on \,$\wh{\Sigma}$\, in the singular points of \,$\Sigma$\,) and \,$\int_{A_k} \omega=0$\, holds for all \,$k\in \Z$\,,
then we have \,$\omega=0$\,.

By applying this statement to \,$\omega := \omega_1-\omega_2$\,, the claimed result follows.
%
%
\end{proof}

We next proceed to the explicit construction of the canonical basis \,$(\omega_n)$\, of the holomorphic 1-forms. As explained above, we need this explicit representation to be able to apply the results
of Section~\ref{Se:jacobiprep} to obtain a finer description of the asymptotic behavior of the \,$\omega_n$\,; it is also needed for the treatment of the case \,$S\neq\varnothing$\,. 

The strategy for the explicit construction is to use the estimates of Proposition~\ref{P:jacobiprep:Deltaint-neu} and Propositions~\ref{P:jacobi:asymp-Phi}--\ref{P:jacobi:asymp-omega}
to construct holomorphic \,$1$-forms that have zeros in all excluded domains \,$\wh{U}_{k,\delta}$\, except for \,$k=n$\,,
and to use an argument based on the Banach Fixed Point Theorem to adjust the position of these zeros such that 
\,$\int_{A_k} \omega_n=0$\, holds for \,$k\neq n$\, with \,$|k|$\, large. Multiplication with an appropriate constant factor
ensures \,$\int_{A_n} \omega_n=1$\,, whereas an additional (finite) linear combination of 1-forms of the described kind is needed to achieve
\,$\int_{A_k} \omega_n=0$\, also for small values of \,$|k|$\,. 

\begin{thm}
\label{T:jacobi:canonical}
\begin{enumerate}
\item 
There exists \,$R_0>0$\, so that for every \,$C_\xi \geq R_0$\, there exists \,$N\in \N$\, with the following property:

For every \,$n\in \Z$\, and every given finite sequence \,$(\xi_{n,k})_{|k|\leq N, k\neq n}$\, 
with \,$|\xi_{n,k}-\vkap_{k,*}|\leq C_\xi \cdot |\vkap_{k,1}-\vkap_{k,2}|$\, there exists an extension of this sequence to an infinite sequence \,$(\xi_{n,k})_{k\in \Z\setminus \{n\}}$\, 
still with \,$|\xi_{n,k}-\vkap_{k,*}|\leq C_\xi \cdot |\vkap_{k,1}-\vkap_{k,2}|$\,, so that:
\begin{enumerate}
\item With \,$\wt{\Phi}_n := \Phi_{n,\xi_{n,k},\rho=-1}$\, (see Proposition~\ref{P:jacobi:asymp-Phi}) and \,$\wt{\omega}_n  := \tfrac{\wt{\Phi}_{n}(\lambda)}{\mu-\mu^{-1}}\, \mathrm{d}\lambda$\,, we have:
\begin{enumerate}
\item If \,$n\not\in S$\,, then \,$\wt{\omega}_n \in \Omega(\Sigma)\cap L^2(\Sigma,T^*\Sigma)$\, holds.
\item If \,$n\in S$\,, then \,$\wt{\omega}_n \in \Omega(\Sigma\setminus \{\vkap_{n,*}\}) \cap L^2(\Sigma\setminus \wh{U}_{n,\delta},T^*(\Sigma\setminus \wh{U}_{n,\delta}))$\, 
for every \,$\delta>0$\,, and \,$\wt{\omega}_n$\, is regular in \,$\vkap_{n,*}$\,. 
\end{enumerate}
\item We have
$$ \int_{A_k} \wt{\omega}_n = 0 \qmq{for every \,$k\in \Z\setminus\{n\}$\, with (\,$|k|>N$\, or \,$k\in S$\,)} \;. $$
\end{enumerate}

\item
Let \,$R_0>0$\, be as in (1), \,$C_\xi \geq R_0$\,, and \,$N\in \N$\, corresponding to \,$C_\xi$\, as in (1). 
We fix pairwise unequal 
\,$(\xi_k)_{|k|\leq N}$\, so that \,$|\xi_k-\vkap_{k,*}| \leq C_\xi\cdot |\vkap_{k,1}-\vkap_{k,2}|$\, holds for all \,$|k|\leq N$\,. We apply (1) for every \,$n\in \Z$\,, yielding a sequence
\,$(\xi_{n,k})_{k\in \Z\setminus\{n\}}$\, with \,$\xi_{n,k}=\xi_k$\, for \,$|k|\leq N$\,, \,$k\neq n$\, so that the associated function \,$\wt{\Phi}_n$\, and the 1-form \,$\wt{\omega}_n$\, have the properties (a) and (b) from (1).

Then for every \,$n\in \Z$\, there exist numbers \,$s_{n,\ell}\in \C$\, for \,$|\ell|\leq N$\, (where \,$s_{n,\ell}=0$\, for every \,$\ell \in S$\, with \,$|\ell|\leq N$\, and \,$\ell\neq n$\,) 
and for \,$|n|>N$\, also \,$s_{n,n}\in \C$\,, such that the 1-form
$$ \omega_n := \left. \begin{cases}
\displaystyle \sum_{\substack{|\ell|\leq N}}  s_{n,\ell}\,\wt{\omega}_\ell & \text{if \,$|n|\leq N$\,} \\
\displaystyle s_{n,n}\,\wt{\omega}_n + \sum_{\substack{|\ell|\leq N}}  s_{n,\ell}\,\wt{\omega}_\ell & \text{if \,$|n|>N$\,}
\end{cases} \right\} =  \frac{\Phi_n(\lambda)}{\mu-\mu^{-1}}\,\mathrm{d}\lambda $$
with
$$ \Phi_n := \begin{cases}
\displaystyle \sum_{\substack{|\ell|\leq N}}  s_{n,\ell}\cdot \wt{\Phi}_{\ell} & \text{if \,$|n|\leq N$\,} \\
\displaystyle s_{n,n}\cdot \wt{\Phi}_{n} + \sum_{\substack{|\ell|\leq N}}  s_{n,\ell}\cdot \wt{\Phi}_{\ell} & \text{if \,$|n|>N$\,} \end{cases}  $$
has the property (1)(a)(i)--(ii) and satisfies
$$ \int_{A_k} \omega_n = \delta_{kn} \qmq{for all \,$k,n\in \Z$\,.} $$
Here we have \,$\Phi_n \in \As(\C^*,\ell^\infty_{3,-1},1)$\,. 

The constants \,$s_{n,n}$\, have the following asymptotic behavior for \,$|n|\to\infty$\,:
\begin{equation}
\label{eq:jacobi:canonical:snn}
s_{n,n} = 
\begin{cases}
\frac{1}{4\pi} \cdot \lambda_{n,0}^{-1/2} + \ell^2_{1}(n) & \text{for \,$n>N$\,} \\
\frac{1}{4\pi} \cdot \tau_{\xi_{n,k}}^{2}\cdot \lambda_{n,0}^{1/2} + \ell^2_{1}(n) & \text{for \,$n<-N$\,}
\end{cases} \;,
\end{equation}
and the constants \,$s_{n,\ell}$\, with \,$|\ell|\leq N$\, are also of the order of \,$|n|^{-1}$\, for \,$n\to \pm\infty$\,. 
We have \,$s_{n,\ell}=0$\, for \,$\ell \in S$\,, \,$\ell\neq n$\,. 
\end{enumerate}
\end{thm}

\begin{rem}
The 1-forms \,$\omega_n$\, with \,$n\in \Z$\, from Theorem~\ref{T:jacobi:canonical}(2) are for \,$S=\varnothing$\, the canonical basis of \,$\Omega(\Sigma)\cap L^2(\Sigma,T^*\Sigma)$\, 
associated to the basis of the homology \,$(A_n,B_n)$\,
of \,$\Sigma$\,. If \,$n\in S$\, holds, then \,$\omega_n$\, is not holomorphic, but only regular in \,$\vkap_{n,*}$\,, but \,$(\omega_n)_{n\in \Z}$\, is still the appropriate replacement 
for the basis of \,$\Omega(\Sigma)\cap L^2(\Sigma,T^*\Sigma)$\, in the case where the spectral curve \,$\Sigma$\, has singularities.

In either case, Theorem~\ref{T:jacobi:canonical} shows that any \,$\omega_n$\, is a finite linear combination of 1-forms of the
form \,$\tfrac{\Phi_k(\lambda)}{\mu-\mu^{-1}}\,\mathrm{d}\lambda$\,, where the holomorphic functions \,$\Phi_k:\C^*\to\C$\, are infinite products of the type considered in Proposition~\ref{P:jacobi:asymp-Phi},
and therefore have asymptotically and totally exactly one zero in every excluded domain \,$\wh{U}_{j,\delta}$\,, with the exception of \,$\wh{U}_{k,\delta}$\,. 

A natural question is if \,$\omega_n$\, itself can be represented in the form \,$\tfrac{\Phi_n(\lambda)}{\mu-\mu^{-1}}\,\mathrm{d}\lambda$\, with the function \,$\Phi_n$\, being an infinite product 
as in Proposition~\ref{P:jacobi:asymp-Phi}, i.e.~without a linear combination.
If this were the case, then \,$\omega_n$\, would have asymptotically and totally exactly one zero in every excluded domain, except for \,$\wh{U}_{n,\delta}$\,. Unfortunately however, 
the answer to this question is negative in general, as the following argument shows:

It follows from Proposition~\ref{P:jacobi:asymp-Phi}(1)
that the holomorphic function \,$\Phi_n$\, from Theorem~\ref{T:jacobi:canonical}(2) has the following asymptotic behavior: There exist constants \,$z_{n,\infty},z_{n,0} \in \C$\, such that
\,$\Phi_n - z_{n,\infty}\,\Phi_0 \in \As_\infty(\C^*,\ell^2_3,1)$\, and \,$\Phi_n - z_{n,0} \Phi_0 \in \As_0(\C^*,\ell^2_{-1},1)$\, holds with the comparison function
\,$\Phi_0(\lambda) := \tfrac{c_0(\lambda)}{\lambda\cdot(\lambda-\lambda_{0,0})}$\,. If the constants \,$z_{n,\infty}$\, and \,$z_{n,0}$\, are both non-zero, then it can indeed be shown
that \,$\Phi_n(\lambda)$\, itself can be represented in the form \,$s_n \cdot \Phi_{n,\xi_{n,k},-1}$\, with a constant \,$s_n\in \C^*$\, and a sequence \,$(\xi_{n,k})_{k\in \Z\setminus\{n\}}$\, with 
\,$|\xi_{n,k}-\vkap_{k,*}| \leq C_\xi \cdot |\vkap_{k,1}-\vkap_{k,2}|$\,. 

However the case that either of the constants \,$z_{n,\infty}$\, or \,$z_{n,0}$\, vanishes can occur, depending on the specific shape of the spectral curve \,$\Sigma$\,. 
(For \,$n\to \infty$\,, \,$z_{n,\infty}$\, grows of the order \,$n^2$\,, and therefore is non-zero for sufficiently large \,$n$\,, but \,$z_{n,0}$\, remains bounded, and there is no general reason preventing it from becoming zero.
For \,$n\to -\infty$\,, the same is true with the roles of \,$z_{n,\infty}$\, and \,$z_{n,0}$\, reversed.)
This corresponds to the case that one or more zeros \,$\xi_{n,k}$\, from the compact part of \,$\Sigma$\, (i.e.~with \,$|k|$\, small) have moved out to \,$\infty$\, resp.~to \,$0$\,. It
is for this reason that a product representation of the form \,$s_n \cdot \Phi_{n,\xi_{n,k},-1}$\, is not possible for \,$\Phi_n$\, in general.
\end{rem}

\begin{proof}[Proof of Theorem~\ref{T:jacobi:canonical}.]
\emph{For (1).}
For an arbitrary sequence \,$(\xi_{n,k})$\, with \,$\xi_{n,k}-\lambda_{k,0} \in \ell^{2}_{-1,3}(k)$\,, \,$\Phi_{n,\xi_{n,k}} := \Phi_{n,\xi_{n,k},\rho=-1}$\, is a holomorphic function by Proposition~\ref{P:jacobi:asymp-Phi}, and
\,$\omega_{n,\xi_{n,k}} := \tfrac{\Phi_{n,\xi_{n,k}}(\lambda)}{\mu-\mu^{-1}}\,\mathrm{d}\lambda$\, is a holomorphic \,$1$-form on \,$\Sigma'$\, by Proposition~\ref{P:jacobi:asymp-omega}.
If moreover \,$|\xi_{n,k}-\vkap_{k,*}|\leq C_\xi\cdot |\vkap_{k,1}-\vkap_{k,2}|$\, holds for some \,$C_\xi>0$\,, then Proposition~\ref{P:jacobi:asymp-f}(2) yields that 
\,$\omega_{n,\xi_{n,k}} \in L^2(\Sigma,T^*\Sigma)$\, holds if \,$n\not\in S$\,, whereas \,$\omega_{n,\xi_{n,k}} \in L^2(\Sigma\setminus \wh{U}_{n,\delta},T^*(\Sigma\setminus \wh{U}_{n,\delta}))$\, holds 
for all \,$\delta>0$\, if \,$n\in S$\,. We also note that for \,$k\in S$\, with \,$k\neq n$\,, the condition \,$|\xi_{n,k}-\vkap_{k,*}|\leq C_\xi\cdot |\vkap_{k,1}-\vkap_{k,2}|$\, implies
\,$\xi_{n,k}=\vkap_{k,*}$\,; therefore both \,$\Phi_{n,\xi_{n,k}}(\lambda)\,\mathrm{d}\lambda$\, and \,$\mu-\mu^{-1}$\, have a zero of order \,$2$\, at \,$\vkap_{k,*}$\,, whence it follows that  
\,$\omega_{n,\xi_{n,k}}$\, is holomorphic in \,$\vkap_{k,*}$\,. Thus \,$\int_{A_k} \omega_{n,\xi_{n,k}}=0$\, holds automatically for \,$k\in S\setminus \{n\}$\, by \eqref{eq:jacobi:intAk-singular}. 


In the sequel we fix \,$R_0>0$\, and let \,$C_\xi\geq R_0$\, be given. We also fix \,$N\in \N$\,. It will 
turn out in the course of the proof how to choose \,$R_0$\, and \,$N$\, (the latter in dependence on \,$C_\xi$\,) so that the claimed statement holds.
We will always choose \,$N$\, at least as large as the \,$N$\, in Lemma~\ref{L:jacobiprep:Psi-intro}, so that Proposition~\ref{P:jacobiprep:Deltaint-neu}(1),(2) is applicable to all \,$k\in \Z$\, with \,$|k|>N$\,. 

Moreover, we let \,$n\in \Z$\, be given.  We then consider the Banach space
$$ \mathfrak{B} := \left\{(a_k)_{k\in \Z \setminus \{n\}} \;\left|\; \exists \,C\geq 0\;\forall\, k\in \Z\setminus\{n\} \; : \; |a_k| \leq C \cdot |\vkap_{k,1}-\vkap_{k,2}| \right. \right\} $$
with the norm
$$ \|a_k\|_{\mathfrak{B}} := \sup_{k\in \Z\setminus (\{n\} \cup S)} \left| \frac{a_k}{\vkap_{k,1}-\vkap_{k,2}} \right| \qmq{for \,$(a_k) \in \mathfrak{B}$\,.} $$
Note that for any \,$(a_k)\in \mathfrak{B}$\,, we have \,$a_k=0$\, for all \,$k\in S$\,. 
We also suppose that a finite sequence \,$(\xi_{n,k}^{[0]})_{|k|\leq N,k\neq n}$\, with \,$|\xi_{n,k}^{[0]}-\vkap_{k,*}| \leq C_\xi\cdot |\vkap_{k,1}-\vkap_{k,2}|$\, is given, and consider
the translated subspace disc in \,$\mathfrak{B}$\,
$$ \mathfrak{B}_{C_\xi} := \Mengegr{(\xi_{n,k})_{k\in \Z \setminus \{n\}}}{(\xi_{n,k}-\vkap_{k,*})\in \mathfrak{B}, \|\xi_{n,k}-\vkap_{k,*}\|_{\mathfrak{B}}\leq C_\xi, \; \text{\,$\xi_{n,k}=\xi_{n,k}^{[0]}$\, for \,$|k|\leq N$\,}} \;. $$

We now associate to each sequence \,$(\xi_{n,k}) \in \mathfrak{B}_{C_\xi}$\, 
a new sequence \,$(\wt{\xi}_{n,k})_{k\in \Z\setminus\{n\}}$\, in the following way: 
For \,$k\in \Z \setminus \{n\}$\,, we define the holomorphic function \,$\Phi_{n,\xi_{n,k};k}: \C^* \to \C$\, by
$$ \Phi_{n,\xi_{n,k};k}(\lambda) := \frac{\Phi_{n,\xi_{n,k}}(\lambda)}{\lambda-\xi_{n,k}} 
\; . $$
We then define \,$\wt{\xi}_{n,k}$\, in the following way: For \,$|k|\leq N$\, we put \,$\wt{\xi}_{n,k} := \xi_{n,k}=\xi_{n,k}^{[0]}$\,, for \,$|k|>N$\, with \,$k\in S$\, we put \,$\wt{\xi}_{n,k} := \vkap_{k,*}$\,,
and for \,$|k|>N$\, with \,$k\not\in S$\, we put 
\begin{equation}
\label{eq:jacobi:canonical:wtxidef}
\wt{\xi}_{n,k} := 
\frac{\displaystyle \int_{A_k} \lambda \cdot \frac{\Phi_{n,\xi_{n,k};k}(\lambda)}{\mu-\mu^{-1}} \,\mathrm{d}\lambda}{\displaystyle \int_{A_k} \frac{\Phi_{n,\xi_{n,k};k}(\lambda)}{\mu-\mu^{-1}} \,\mathrm{d}\lambda} 
= \vkap_{k,*} + \frac{\displaystyle \int_{A_k} (\lambda-\vkap_{k,*}) \cdot \frac{\Phi_{n,\xi_{n,k};k}(\lambda)}{\mu-\mu^{-1}} \,\mathrm{d}\lambda}{\displaystyle \int_{A_k} \frac{\Phi_{n,\xi_{n,k};k}(\lambda)}{\mu-\mu^{-1}} \,\mathrm{d}\lambda} \;.
\end{equation}

We will show below that (for suitable conditions on \,$R_0$\, and \,$N$\,) 
the map \,$(\xi_{n,k}) \mapsto (\wt{\xi}_{n,k})$\, we just defined has exactly one fixed point \,$(\xi_{n,k}^*)$\,. The sequence \,$(\xi_{n,k}^*)$\, has the following property:
For any \,$k\in \Z$\, with \,$|k|>N$\, and \,$k\not\in S$\, we have
\begin{align}
\int_{A_k} \frac{\Phi_{n,\xi_{n,k}^*}(\lambda)}{\mu-\mu^{-1}}\,\mathrm{d}\lambda
& = \int_{A_k} (\lambda-\xi^*_{n,k}) \cdot \frac{\Phi_{n,\xi_{n,k}^*;k}(\lambda)}{\mu-\mu^{-1}}\,\mathrm{d}\lambda \notag \\
& = \int_{A_k} \lambda \cdot \frac{\Phi_{n,\xi_{n,k}^*;k}(\lambda)}{\mu-\mu^{-1}}\,\mathrm{d}\lambda - \underbrace{\xi^*_{n,k}}_{=\wt{\xi}^*_{n,k}} \cdot \int_{A_k} \frac{\Phi_{n,\xi_{n,k}^*;k}(\lambda)}{\mu-\mu^{-1}}\,\mathrm{d}\lambda \notag \\
& = \int_{A_k} \lambda \cdot \frac{\Phi_{n,\xi_{n,k}^*;k}(\lambda)}{\mu-\mu^{-1}}\,\mathrm{d}\lambda - \frac{\displaystyle \int_{A_k} \lambda \cdot \frac{\Phi_{n,\xi_{n,k}^*;k}(\lambda)}{\mu-\mu^{-1}} \,\mathrm{d}\lambda}{\displaystyle \int_{A_k} \frac{\Phi_{n,\xi_{n,k}^*;k}(\lambda)}{\mu-\mu^{-1}} \,\mathrm{d}\lambda} \cdot \int_{A_k} \frac{\Phi_{n,\xi_{n,k}^*;k}(\lambda)}{\mu-\mu^{-1}}\,\mathrm{d}\lambda \notag \\
\label{eq:jacobi:canonical:fixiegoodie}
& = 0\;.
\end{align}
Moreover, for any \,$k\in \Z$\, with \,$|k|>N$\, and \,$k\in S$\,, the 1-form \,$\frac{\Phi_{n,\xi_{n,k}^*}(\lambda)}{\mu-\mu^{-1}}\,\mathrm{d}\lambda$\, extends holomorphically in \,$\vkap_{k,*}=\xi_{n,k}^*$\,,
and therefore we have
$$ \int_{A_k} \frac{\Phi_{n,\xi_{n,k}^*}(\lambda)}{\mu-\mu^{-1}}\,\mathrm{d}\lambda = 0 $$
also in this case. This shows the form \,$\wt{\omega}_n := \tfrac{\Phi_{n,\xi_{n,k}^*}(\lambda)}{\mu-\mu^{-1}}\,\mathrm{d}\lambda$\, then has the properties (a) and (b) from the theorem (1).

We will now show that if \,$N$\, is chosen large enough (depending on \,$C_\xi$\,), then \,$F: (\xi_{n,k}) \mapsto (\wt{\xi}_{n,k})$\, maps \,$\mathfrak{B}_{C_\xi}$\, into itself and is a contraction.
Therefore, by the Banach Fixed Point Theorem, there then exists one and only one fixed point of this map; this fixed point yields the solution to the problem of the part (1) of the theorem, as was just explained.

We begin by showing that the map \,$F$\, maps \,$\mathfrak{B}_{C_\xi}$\, into \,$\mathfrak{B}_{C_\xi}$\,. 
We let \,$(\xi_{n,k})\in \mathfrak{B}_{C_\xi}$\, be given and put \,$(\wt{\xi}_{n,k}) := F((\xi_{n,k}))$\,. 
Then we have to show that 
\begin{equation}
\label{eq:jacobi:canonical:F-bounded}
|\wt{\xi}_{n,k}-\vkap_{k,*}| \leq C_\xi\cdot |\vkap_{k,1}-\vkap_{k,2}|
\end{equation}
holds for all \,$k\in \Z$\, with \,$|k|> N$\, and \,$k\not\in S$\,.

In the sequel, we will consider the case where \,$n$\, and \,$k$\, are positive; the other cases are handled analogously. 
By Proposition~\ref{P:jacobiprep:Deltaint-neu}(2) and Equation~\eqref{eq:jacobi:intAk} we have
\begin{align*}
\int_{A_k} \frac{1}{\mu-\mu^{-1}}\,\mathrm{d}\lambda & = -8\pi\,(-1)^k\,\lambda_{k,0}^{1/2} + \ell^{2}_{-1}(k) \\
\int_{A_k} \frac{1}{|\mu-\mu^{-1}|}\,|\mathrm{d}\lambda| & = 8\pi\,\lambda_{k,0}^{1/2} + \ell^{2}_{-1}(k) \; ,
\end{align*}
therefore we may suppose that \,$N$\, is chosen large enough so that for \,$k>N$\, we have
\begin{equation}
\label{eq:jacobi:canonical:int1Ak-betragrein}
\frac12 \, \int_{A_k} \frac{1}{|\mu-\mu^{-1}|}\,|\mathrm{d}\lambda| \leq \left| \int_{A_k} \frac{1}{\mu-\mu^{-1}}\,\mathrm{d}\lambda \right| \leq \int_{A_k} \frac{1}{|\mu-\mu^{-1}|}\,|\mathrm{d}\lambda| \; . 
\end{equation}

Moreover, Proposition~\ref{P:jacobi:asymp-Phi}(3) shows that there exist constants \,$C_1,C_2>0$\, (depending only on \,$C_\xi$\,)
so that we have for \,$\lambda\in U_{k,\delta}$\,
\begin{equation}
\label{eq:jacobi:canonical:Phink++-estim-pre-pre}
\left| \Phi_{n,\xi_{n,k};k}(\lambda) - \frac{8(-1)^k\,\pi^2\,n^2}{\xi_{n,k} \cdot (\lambda_{n,0}-\xi_{n,k})} \right| \leq C_1 \cdot \frac{n^2}{k^2 \cdot |n^2-k^2|} \cdot \frac{|\lambda-\xi_{n,k}|}{k}+ C_2\cdot r_k
\end{equation}
with the sequence \,$(r_k) \in \ell^2_{0}(k)$\, given by
$$ r_k = \left(a_j * \frac{1}{|j|}\right)_k \qmq{with} a_k := \begin{cases} k^{-1}\cdot |\xi_{n,k}-\lambda_{k,0}| & \text{for \,$k>0$\,} \\ k^{3}\cdot |\xi_{n,k}-\lambda_{k,0}| & \text{for \,$k<0$\,}
\end{cases} \; . $$

We may suppose without loss of generality that the projection of the cycle \,$A_k$\, onto the \,$\lambda$-plane is the straight line \,$[\vkap_{k,1},\vkap_{k,2}] \subset U_{k,\delta}$\,, traversed twice
in opposite directions. Then we have for any \,$\lambda$\, in the projection of \,$A_k$\,, i.e.~for \,$\lambda \in [\vkap_{k,1},\vkap_{k,2}]$\,
$$ |\lambda-\xi_{n,k}| \leq |\lambda-\vkap_{k,*}| + |\vkap_{k,*}-\xi_{n,k}| \leq \left(\tfrac12+C_\xi\right)\cdot |\vkap_{k,1}-\vkap_{k,2}| \;, $$
and we also have
$$ |\xi_{n,k}-\lambda_{k,0}| \leq |\xi_{n,k}-\vkap_{k,*}| + |\vkap_{k,*}-\lambda_{k,0}| \leq C_\xi \cdot |\vkap_{k,1}-\vkap_{k,2}| + |\vkap_{k,*}-\lambda_{k,0}| \; . $$
We now define sequences \,$(a_k'), (r_k')\in \ell^2_0(k)$\, by 
$$ a_k' := \begin{cases} k^{-1}\cdot (C_\xi\cdot |\vkap_{k,1}-\vkap_{k,2}| + |\vkap_{k,*}-\lambda_{k,0}|) & \text{for \,$k>0$\,} \\ k^{3}\cdot (C_\xi\cdot |\vkap_{k,1}-\vkap_{k,2}| + |\vkap_{k,*}-\lambda_{k,0}|) & \text{for \,$k<0$\,} \end{cases} $$
and
$$ r_k' := C_1\,\left(C_\xi+\tfrac12\right)\,\frac{|\vkap_{k,1}-\vkap_{k,2}|}{k} + C_2\,\left(a_j' * \frac{1}{|j|}\right)_k \;. $$
It then follows from \eqref{eq:jacobi:canonical:Phink++-estim-pre-pre} that we have for \,$\lambda \in [\vkap_{k,1},\vkap_{k,2}]$\,
\begin{equation}
\label{eq:jacobi:canonical:Phink++-estim-pre}
\left| \Phi_{n,\xi_{n,k};k}(\lambda) - \frac{8(-1)^k\,\pi^2\,n^2}{\xi_{n,k} \cdot (\lambda_{n,0}-\xi_{n,k})} \right|
\leq \frac{n^2}{k^2 \cdot |n^2-k^2|} \cdot r_k' \; .
\end{equation}
Note that the sequence \,$(r_k') \in \ell^2_{0}(k)$\, depends only on \,$C_\xi$\, (and the spectral curve \,$\Sigma$\,), 
but neither on \,$n$\, nor on the specific sequence \,$(\xi_{n,k})\in \mathfrak{B}_{C_\xi}$\, under consideration. Therefore we may suppose that \,$N$\,
is chosen large enough so that \,$|r_k'|\leq \tfrac12$\, holds for \,$k>N$\,. 

From Equation~\eqref{eq:jacobi:canonical:Phink++-estim-pre} we obtain by multiplication with \,$\tfrac{1}{\mu-\mu^{-1}}$\,,
integration, and application of Equation~\eqref{eq:jacobi:canonical:int1Ak-betragrein} 
for \,$k>N$\,
\begin{gather}
\frac{2\,\pi^2\,n^2}{|\xi_{n,k}| \cdot |\lambda_{n,0}-\xi_{n,k}|} \cdot \int_{A_k} \frac{1}{|\mu-\mu^{-1}|}\,|\mathrm{d}\lambda|
\leq \left| \int_{A_k} \frac{\Phi_{n,\xi_{n,k};k}(\lambda)}{\mu-\mu^{-1}}\,\mathrm{d}\lambda \right| 
\hspace{3cm} \notag \\
\label{eq:jacobi:canonical:int-Phi}
\hspace{3cm}
\leq \int_{A_k} \frac{|\Phi_{n,\xi_{n,k};k}(\lambda)|}{|\mu-\mu^{-1}|}\,|\mathrm{d}\lambda|
\leq \frac{16\pi^2\,n^2}{|\xi_{n,k}| \cdot |\lambda_{n,0}-\xi_{n,k}|} \cdot \int_{A_k} \frac{1}{|\mu-\mu^{-1}|}\,|\mathrm{d}\lambda| \; . 
\end{gather}


For any \,$\lambda$\, in the projection of \,$A_k$\, to the \,$\lambda$-plane, i.e.~for \,$\lambda \in [\vkap_{k,1},\vkap_{k,2}]$\, we have \,$|\lambda-\vkap_{k,*}|\leq \tfrac12\,|\vkap_{k,1}-\vkap_{k,2}|$\,,
and therefore it follows from \eqref{eq:jacobi:canonical:int-Phi} that we have
\begin{equation}
\label{eq:jacobi:canonical:int-Phi-lambda}
\left| \int_{A_k} (\lambda-\vkap_{k,*})\cdot \frac{\Phi_{n,\xi_{n,k};k}(\lambda)}{\mu-\mu^{-1}}\,\mathrm{d}\lambda \right| \leq \frac12 \, |\vkap_{k,1}-\vkap_{k,2}| \cdot \frac{16\pi^2\,n^2}{|\xi_{n,k}| \cdot |\lambda_{n,0}-\xi_{n,k}|} \cdot \int_{A_k} \frac{1}{|\mu-\mu^{-1}|}\,|\mathrm{d}\lambda| \; . 
\end{equation}


From the estimates Equation~\eqref{eq:jacobi:canonical:int-Phi-lambda} and \eqref{eq:jacobi:canonical:int-Phi} we now obtain
$$ |\wt{\xi}_{n,k}-\vkap_{k,*}| \overset{\eqref{eq:jacobi:canonical:wtxidef}}{=} \left|\frac{\displaystyle \int_{A_k} (\lambda-\vkap_{k,*}) \cdot \frac{\Phi_{n,\xi_{n,k};k}(\lambda)}{\mu-\mu^{-1}} \,\mathrm{d}\lambda}{\displaystyle \int_{A_k} \frac{\Phi_{n,\xi_{n,k};k}(\lambda)}{\mu-\mu^{-1}} \,\mathrm{d}\lambda} \right|
\leq 4\cdot |\vkap_{k,1}-\vkap_{k,2}| \;. $$
By a similar argument for the case where \,$k$\, is negative, we also obtain for \,$k<-N$\,
$$ |\wt{\xi}_{n,k}-\vkap_{k,*}| \leq 4\cdot |\vkap_{k,1}-\vkap_{k,2}| \;, $$
and therefore Equation~\eqref{eq:jacobi:canonical:F-bounded} holds if \,$C_\xi \geq 4$\,. This shows that (at least) for \,$C_\xi \geq 4$\,, there exists a sufficiently large \,$N\in\N$\, 
such that \,$(\xi_{n,k})\in \mathfrak{B}_{C_\xi}$\, implies
\,$(\wt{\xi}_{n,k}) \in \mathfrak{B}_{C_\xi}$\,. The case \,$n<0$\, is handled in an analogous manner.
Therefore the map \,$F: (\xi_{n,k}) \mapsto (\wt{\xi}_{n,k})$\, indeed maps \,$\mathfrak{B}_{C_\xi}$\, into \,$\mathfrak{B}_{C_\xi}$\,. 


We now show that \,$F: \mathfrak{B}_{C_\xi} \to \mathfrak{B}_{C_\xi}$\, is a contraction, if \,$R_0$\, and \,$N$\, are suitably chosen. 
We let \,$(\xi_{n,k}^{[1]}), (\xi_{n,k}^{[2]}) \in \mathfrak{B}_{C_\xi}$\, be given, 
and put \,$\wt{\xi}_{n,k}^{[\nu]} := F(\xi_{n,k}^{[\nu]})$\, for \,$\nu \in \{1,2\}$\,. 
We again consider the case of positive \,$n$\, and \,$k$\,. 
Then we have for \,$k>N$\,, \,$k\not\in S$\, by Equation~\eqref{eq:jacobi:canonical:wtxidef}
\begin{align}
\wt{\xi}_{n,k}^{[1]} - \wt{\xi}_{n,k}^{[2]}
& = \frac{\displaystyle \int_{A_k} (\lambda-\vkap_{k,*}) \cdot \frac{\Phi_{n,\xi_{n,k}^{[1]};k}(\lambda)}{\mu-\mu^{-1}} \,\mathrm{d}\lambda}{\displaystyle \int_{A_k} \frac{\Phi_{n,\xi_{n,k}^{[2]};k}(\lambda)}{\mu-\mu^{-1}} \,\mathrm{d}\lambda} 
- \frac{\displaystyle \int_{A_k} (\lambda-\vkap_{k,*}) \cdot \frac{\Phi_{n,\xi_{n,k}^{[2]};k}(\lambda)}{\mu-\mu^{-1}} \,\mathrm{d}\lambda}{\displaystyle \int_{A_k} \frac{\Phi_{n,\xi_{n,k}^{[2]};k}(\lambda)}{\mu-\mu^{-1}} \,\mathrm{d}\lambda} \notag \\
& = \frac{\displaystyle \int_{A_k} (\lambda-\vkap_{k,*}) \cdot \frac{\Phi_{n,\xi_{n,k}^{[1]};k}(\lambda)-\Phi_{n,\xi_{n,k}^{[2]};k}(\lambda)}{\mu-\mu^{-1}} \,\mathrm{d}\lambda}{\displaystyle \int_{A_k} \frac{\Phi_{n,\xi_{n,k}^{[1]};k}(\lambda)}{\mu-\mu^{-1}} \,\mathrm{d}\lambda} \notag \\
\label{eq:jacobi:canonical:wtxik12-pre}
& \qquad\qquad + \frac{\displaystyle \int_{A_k} (\lambda-\vkap_{k,*}) \cdot \frac{\Phi_{n,\xi_{n,k}^{[2]};k}\lambda)}{\mu-\mu^{-1}} \,\mathrm{d}\lambda}{\displaystyle \int_{A_k} \frac{\Phi_{n,\xi_{n,k}^{[2]};k}(\lambda)}{\mu-\mu^{-1}} \,\mathrm{d}\lambda} 
\cdot \frac{\displaystyle \int_{A_k} \frac{\Phi_{n,\xi_{n,k}^{[2]};k}(\lambda)-\Phi_{n,\xi_{n,k}^{[1]};k}(\lambda)}{\mu-\mu^{-1}} \,\mathrm{d}\lambda}{\displaystyle \int_{A_k} \frac{\Phi_{n,\xi_{n,k}^{[1]};k}(\lambda)}{\mu-\mu^{-1}} \,\mathrm{d}\lambda} \; . 
\end{align}

We now note that by Proposition~\ref{P:jacobi:asymp-Phi}(2) we have for any \,$\lambda \in U_{k,\delta}$\,
$$ |\Phi_{n,\xi_{n,k}^{[1]};k}(\lambda)-\Phi_{n,\xi_{n,k}^{[2]};k}(\lambda)| \leq C_3\cdot r_k \;, $$
where \,$C_3>0$\, is a constant depending only on \,$C_\xi$\, and \,$(r_k)\in \ell^2_4(k)$\, is the sequence given by 
$$ r_k := \frac{n^2}{k^2\cdot |n^2-k^2|}\,\left(a_j * \frac{1}{|j|}\right)_k 
\qmq{with}
a_k := \begin{cases} k^{-1}\,|\xi_{n,k}^{[1]}-\xi_{n,k}^{[2]}| & \text{if \,$k>0$\,} \\ k^{3}\,|\xi_{n,k}^{[1]}-\xi_{n,k}^{[2]}| & \text{if \,$k<0$\,} \end{cases} \;. $$
By the definition of \,$\|\,\cdot\,\|_{\mathfrak{B}}$\,, we have for every \,$k$\,
$$ |\xi_{n,k}^{[1]}-\xi_{n,k}^{[2]}| \leq \|\xi_{n,k}^{[1]}-\xi_{n,k}^{[2]}\|_{\mathfrak{B}} \cdot |\vkap_{k,1}-\vkap_{k,2}| \;, $$
and thus we see that we have
$$ a_k \leq \|\xi_{n,k}^{[1]}-\xi_{n,k}^{[2]}\|_{\mathfrak{B}} \cdot a_k' \qmq{with} a_k' := \begin{cases} k^{-1}\,|\vkap_{k,1}-\vkap_{k,2}| & \text{if \,$k>0$\,} \\ k^{3}\,|\vkap_{k,1}-\vkap_{k,2}| & \text{if \,$k<0$\,} \end{cases} $$
and therefore
$$ r_k \leq \|\xi_{n,k}^{[1]}-\xi_{n,k}^{[2]}\|_{\mathfrak{B}} \cdot r_k' \qmq{with} r_k' := \frac{n^2}{k^2\cdot |n^2-k^2|}\,\left(a_j' * \frac{1}{|j|}\right)_k \; . $$
Note that again, the sequences \,$(a_k')$\, and therefore also \,$(r_k')$\, depend only on \,$C_\xi$\,, but not on \,$n$\, or on the specific sequences \,$(\xi_{n,k}^{[\nu]})\in \mathfrak{B}_{C_\xi}$\, under consideration. We now have
$$ |\Phi_{n,\xi_{n,k}^{[1]};k}(\lambda)-\Phi_{n,\xi_{n,k}^{[2]};k}(\lambda)| \leq C_3\,r_k \leq C_3\,\|\xi_{n,k}^{[1]}-\xi_{n,k}^{[2]}\|_{\mathfrak{B}} \cdot r_k' $$
and therefore by integration
\begin{equation}
\label{eq:jacobi:canonical:Phink12}
\left| \int_{A_k} \frac{\Phi_{n,\xi_{n,k}^{[2]};k}(\lambda)-\Phi_{n,\xi_{n,k}^{[1]};k}(\lambda)}{\mu-\mu^{-1}} \,\mathrm{d}\lambda \right| \leq C_3\,\|\xi_{n,k}^{[1]}-\xi_{n,k}^{[2]}\|_{\mathfrak{B}} \cdot r_k' \cdot \int_{A_k} \frac{1}{|\mu-\mu^{-1}|}\,|\mathrm{d}\lambda| \;,
\end{equation}
whence also 
\begin{gather}
\left| \int_{A_k} (\lambda-\vkap_{k,*})\cdot \frac{\Phi_{n,\xi_{n,k}^{[1]};k}(\lambda)-\Phi_{n,\xi_{n,k}^{[2]};k}(\lambda)}{\mu-\mu^{-1}} \,\mathrm{d}\lambda \right| \notag \\
\label{eq:jacobi:canonical:Phink12-lambda}
\leq \frac12\,C_3\,|\vkap_{k,1}-\vkap_{k,2}|\, \cdot \|\xi_{n,k}^{[1]}-\xi_{n,k}^{[2]}\|_{\mathfrak{B}} \cdot r_k' \cdot \int_{A_k} \frac{1}{|\mu-\mu^{-1}|}\,|\mathrm{d}\lambda| \; 
\end{gather}
follows.

By applying the estimates \eqref{eq:jacobi:canonical:int-Phi}, \eqref{eq:jacobi:canonical:int-Phi-lambda}, \eqref{eq:jacobi:canonical:Phink12} and \eqref{eq:jacobi:canonical:Phink12-lambda} 
to Equation~\eqref{eq:jacobi:canonical:wtxik12-pre}, and cancelling, we now obtain that there exists \,$C_4>0$\, (dependent only on \,$C_\xi$\,) with
\begin{equation}
\label{eq:jacobi:canonical:wtxik12}
|\wt{\xi}_{n,k}^{[1]}-\wt{\xi}_{n,k}^{[2]}| \leq C_4\cdot \left( a_j' * \frac{1}{|j|}\right)_k \cdot |\vkap_{k,1}-\vkap_{k,2}|\cdot \|\xi_{n,k}^{[1]}-\xi_{n,k}^{[2]}\|_{\mathfrak{B}} \;. 
\end{equation}
The analogous calculation applies for \,$k<-N$\,, yielding literally the same estimate as \eqref{eq:jacobi:canonical:wtxik12} for this case, and therefore we obtain
$$ \|\wt{\xi}_{n,k}^{[1]}-\wt{\xi}_{n,k}^{[2]}\|_{\mathfrak{B}} \leq L \cdot \|\xi_{n,k}^{[1]}-\xi_{n,k}^{[2]}\|_{\mathfrak{B}} \qmq{with} L := C_4 \cdot \max_{|k|>N} \left| \left( a_j' * \frac{1}{|j|}\right)_k \right| \; ; $$
note that \,$a_j' * \frac{1}{|j|} \in \ell^2_{0,0}(j)\subset \ell^\infty(j)$\, holds and therefore \,$L$\, is finite. 
This shows that the map \,$F$\, is Lipschitz continuous with the Lipschitz constant \,$L$\,. Because the sequence \,$a_j' * \frac{1}{|j|}$\, is an \,$\ell^2$-sequence which depends on 
\,$C_\xi$\, and nothing else, and also \,$C_4$\, depends only on \,$C_\xi$\,, 
we can choose \,$N$\, so large that \,$\max_{|k|>N} \left| \left( a_j' * \frac{1}{|j|}\right)_k \right| \leq \tfrac{1}{2\,C_4}$\, holds. Then we have \,$L \leq \tfrac12$\,, hence
with this choice of \,$N$\,, \,$F$\, is a contraction, completing the proof of (1).

\emph{For (2).}
By (1), the \,$1$-forms \,$\wt{\omega}_n$\, have the property that
\begin{equation}
\label{eq:jacobi:canonical:wtomega-prop-repeat}
\int_{A_k} \wt{\omega}_n = 0 \qmq{holds for all \,$k,n\in \Z$\, with (\,$|k|> N$\, or \,$k\in S$\,) and \,$k\neq n$\,.} 
\end{equation}

We consider the linear space
$$ V := \left\{ \omega \in \Omega(\Sigma) \cap L^2(\Sigma,T^*\Sigma) \left| \int_{A_k} \omega = 0 \text{ for \,$|k|>N$\, or \,$k\in S$\,} \right. \right\} \;. $$
Let \,$I := \Mengegr{k\in \Z \setminus S}{|k|\leq N}$\, and \,$m := \# I \leq 2N+1$\,. Then it follows from Proposition~\ref{P:jacobi:abstractomega-neu} that 
the linear map \,$f: V \to \C^m,\; \omega \mapsto \left( \int_{A_\ell} \omega \right)_{\ell \in I}$\, is injective and that therefore \,$\dim(V) \leq m$\, holds.
On the other hand, the holomorphic 1-forms \,$\wt{\omega}_\ell$\, with \,$\ell \in I$\, are in \,$\Omega(\Sigma)\cap L^2(\Sigma,T^*\Sigma)$\, and therefore in \,$V$\, by (1), and they are linear independent 
(because \,$\wt{\omega}_\ell$\, does not vanish at \,$\xi_\ell$\,, but all the \,$\wt{\omega}_{\ell'}$\, with \,$\ell'\neq \ell$\, do). Therefore we have \,$\dim(V)=m$\, and \,$(\wt{\omega}_\ell)_{\ell\in I}$\,
is a basis of \,$V$\,. The linear map \,$f: V \to \C^m$\, is therefore bijective, 
which means that for every finite sequence \,$(z_\ell)_{\ell\in I}$\, there exists one and only one \,$\omega\in V$\, with \,$\int_{A_\ell} \omega=z_\ell$\, for all \,$\ell \in I$\,.

Now let \,$n\in \Z$\, be given. If \,$|n|\leq N$\, and \,$n \not\in S$\, holds, then the preceding statement shows immediately that there exists \,$\omega_n \in V$\, such that
$$ \int_{A_k} \omega_n = \delta_{kn} \qmq{holds for all \,$k\in \Z$\,;} $$
because \,$(\wt{\omega}_\ell)_{\ell \in I}$\, is a basis of \,$V$\,, it follows that there exist numbers \,$s_{n,\ell}\in \C$\, for \,$\ell \in I$\, so that \,$\omega_n = \sum_{\ell\in I} s_{n,\ell}\,\wt{\omega}_\ell$\,
holds.

On the other hand, if either \,$n\in S$\, or \,$|n|>N$\, holds, we put \,$s_{n,n} := \left( \int_{A_n} \wt{\omega}_n \right)^{-1}$\,. There then exists \,$\wh{\omega}\in V$\, so that
$$ \int_{A_\ell} \wh{\omega} = -s_{n,n}\cdot \int_{A_\ell} \wt{\omega}_n \qmq{holds for all \,$\ell\in I$\,;} $$
because \,$(\wt{\omega}_\ell)_{\ell \in I}$\, is a basis of \,$V$\,, it follows that there exist numbers \,$s_{n,\ell}\in \C$\, for \,$\ell \in I$\, so that \,$\wh{\omega} = \sum_{\ell\in I} s_{n,\ell}\,\wt{\omega}_\ell$\,
holds; we then have
$$ \int_{A_k} \omega_n = \delta_{kn} \qmq{for all \,$k\in \Z$\,} $$
with
$$ \omega_n := s_{n,n}\cdot \wt{\omega}_n + \wh{\omega} = s_{n,n}\cdot \wt{\omega}_n + \sum_{\ell\in I} s_{n,\ell}\,\wt{\omega}_\ell \; . $$
In either case this shows the claimed representation of \,$\omega_n$\,. 
It follows from Proposition~\ref{P:jacobi:asymp-Phi}(1) that \,$\Phi_n \in \As(\C^*,\ell^\infty_{3,-1},1)$\, holds.

It remains to show the statements on the asymptotic behavior of \,$s_{n,n}$\, and \,$s_{n,\ell}$\,, and for this it suffices to consider the case \,$|n|>N$\,. We have
\begin{equation}
\label{eq:jacobi:canonical:sestim-Phi}
\Phi_{n,\xi_{n,k}}(\lambda) = 
\begin{cases}
-\frac{4}{\tau_{\xi_{n,k}}}\cdot 16\pi^2n^2 \cdot \lambda^{-1}\cdot \frac{c_{\xi_{n,k}}(\lambda)}{\lambda-\lambda_{n,0}} & \text{if \,$n>N$\,} \\
\frac{4}{\tau_{\xi_{n,k}}}\cdot \lambda^{-1}\cdot \frac{c_{\xi_{n,k}}(\lambda)}{\lambda-\lambda_{n,0}} & \text{if \,$n<-N$\,} 
\end{cases} 
\end{equation}
where \,$c_{\xi_{n,k}}$\, and \,$\tau_{\xi_{n,k}}$\, are as in the proof of Proposition~\ref{P:jacobi:asymp-Phi}(1). 
By Corollary~\ref{C:interpolate:cdivlin}(1), there exists a constant \,$C_5>0$\, and a sequence \,$(r_k)\in \ell^2_{0,-2}(k)$\, so that we have for \,$n,k\in \Z$\, with \,$|n|>N$\, and \,$\lambda \in U_{k,\delta}$\, 
$$ \begin{cases}
\left| \frac{c_{\xi_{n,k}}(\lambda)}{\tau_{\xi_{n,k}}\cdot (\lambda-\xi_{n,k})} - \frac{(-1)^k}{8} \right| \leq C_5\,\frac{|\lambda-\xi_{n,k}|}{k} + r_k & \text{if \,$k>0$\,} \\
\left| \frac{c_{\xi_{n,k}}(\lambda)}{\lambda-\xi_{n,k}} - \left( -\tau_{\xi_{n,k}}^{-1}\, \frac{(-1)^k}{8}\,\lambda_{n,0}^{-1} \right) \right| \leq C_5\,|\lambda-\xi_{n,k}|\,k^5 + r_k & \text{if \,$k<0$\,} 
\end{cases} \;, $$
here we put \,$\xi_{n,n} := \lambda_{n,0}$\, and the sequence \,$(r_k)$\, can be chosen independently of \,$n$\, because \,$\left| \tfrac{\xi_{n,k}-\vkap_{k,*}}{\vkap_{k,1}-\vkap_{k,2}} \right|$\, is bounded 
independently of \,$n$\,. By setting \,$k=n$\, in the preceding estimate we see that for \,$\lambda \in [\vkap_{k,1},\vkap_{k,2}]$\, we have
$$ \begin{cases}
\left| \frac{c_{\xi_{n,k}}(\lambda)}{\tau_{\xi_{n,k}}\cdot (\lambda-\lambda_{n,0})} - \frac{(-1)^k}{8} \right| \leq r_n' & \text{if \,$n>N$\,} \\
\left| \frac{c_{\xi_{n,k}}(\lambda)}{\lambda-\lambda_{n,0}} - \left( -\tau_{\xi_{n,k}}^{-1}\, \frac{(-1)^k}{8}\,\lambda_{n,0}^{-1} \right) \right| \leq r_n' & \text{if \,$n<-N$\,}
\end{cases} $$
with the sequence \,$(r_k') \in \ell^2_{0,-2}(k)$\, defined by
$$ r_k' := \begin{cases} C_5\,\left( C_\xi+\tfrac12 \right)\,\tfrac{|\vkap_{k,1}-\vkap_{k,2}|}{k} + r_k & \text{for \,$k>0$\,} \\ C_5\,\left( C_\xi+\tfrac12 \right)\,|\vkap_{k,1}-\vkap_{k,2}|\,k^5 + r_k & \text{for \,$k<0$\,} \end{cases} \;. $$
By applying this estimate, and also the fact that we have \,$\frac{16\pi^2 n^2}{\lambda} - 1 \in \ell^2_1(n)$\, (for the case \,$n>0$\,) resp.~\,$\frac{1}{\lambda} - \frac{1}{\lambda_{n,0}} \in \ell^2_{-1}(n)$\,
(for the case \,$n<0$\,) for \,$\lambda \in [\vkap_{k,1},\vkap_{k,2}]$\,, to Equation~\eqref{eq:jacobi:canonical:sestim-Phi}, we obtain
\begin{equation}
\label{eq:jacobi:canonical:sestim-Phiest}
\begin{cases}
\left| \Phi_{n,\xi_{n,k}}(\lambda)-\left( -\frac{(-1)^k}{2} \right) \right| \in \ell^2(n) & \text{for \,$n>N$\,} \\
\left| \Phi_{n,\xi_{n,k}}(\lambda)-\left( -\frac{(-1)^k}{2}\,\tau_{\xi_{n,k}}^{-2}\,\lambda_{n,0}^{-2} \right) \right| \in \ell^2_{-4}(n) & \text{for \,$n<-N$\,} \; . 
\end{cases} 
\end{equation}
By Proposition~\ref{P:jacobiprep:Deltaint-neu}(2) we have
$$ \int_{A_n} \frac{1}{\mu-\mu^{-1}}\,\mathrm{d}\lambda \overset{\eqref{eq:jacobi:intAk}}{=} 2\cdot \int_{\vkap_{n,1}}^{\vkap_{n,2}} \frac{1}{\mu-\mu^{-1}}\,\mathrm{d}\lambda = \begin{cases}
-8\pi(-1)^n\,\lambda_{n,0}^{1/2} + \ell^2_{-1}(n) & \text{for \,$n>N$\,} \\
-8\pi(-1)^n\,\lambda_{n,0}^{3/2} + \ell^2_{3}(n) & \text{for \,$n<-N$\,} 
\end{cases} $$
and therefore we derive from Equation~\eqref{eq:jacobi:canonical:sestim-Phiest}
$$ \int_{A_n} \frac{\Phi_{n,\xi_{n,k}}(\lambda)}{\mu-\mu^{-1}}\,\mathrm{d}\lambda = 
\begin{cases}
4\pi\cdot \lambda_{n,0}^{1/2} + \ell^2_{-1}(n) & \text{for \,$n>N$\,} \\
4\pi\cdot \tau_{\xi_{n,k}}^{-2}\cdot \lambda_{n,0}^{-1/2} + \ell^2_{-1}(n) & \text{for \,$n<-N$\,} 
\end{cases} $$
Thus we obtain
$$ s_{n,n} = \left( \int_{A_n} \frac{\Phi_{n,\xi_{n,k}}(\lambda)}{\mu-\mu^{-1}}\,\mathrm{d}\lambda \right)^{-1} =
\begin{cases}
\frac{1}{4\pi} \cdot \lambda_{n,0}^{-1/2} + \ell^2_{1}(n) & \text{for \,$n>N$\,} \\
\frac{1}{4\pi} \cdot \tau_{\xi_{n,k}}^{2}\cdot \lambda_{n,0}^{1/2} + \ell^2_{1}(n) & \text{for \,$n<-N$\,}.
\end{cases}  $$
We now calculate the \,$s_{n,\ell}$\,. For this purpose, let \,$P$\, be the ``restricted period matrix'' of the \,$\wt{\omega}_\ell$\,, i.e.~let \,$P=(p_{k,\ell})_{k,\ell\in I}$\, be the \,$(m\times m)$-matrix
with
$$ p_{k,\ell} := \int_{A_k} \wt{\omega}_\ell \; . $$
Note that the matrix \,$P$\, is invertible because of Proposition~\ref{P:jacobi:abstractomega-neu} and the fact that \,$(\wt{\omega}_\ell)_{\ell\in I}$\, is a basis of \,$V$\,. We have
$$ \omega_n = s_{n,n}\cdot \wt{\omega}_n + \sum_{\substack{\ell\in I}} s_{n,\ell}\cdot \wt{\omega}_\ell $$
and therefore for \,$k\in I$\,
$$ 0 = \int_{A_k} \omega_n = s_{n,n} \cdot \int_{A_k} \wt{\omega}_n + \sum_{\substack{\ell\in I}} p_{k,\ell} \cdot s_{n,\ell} \; . $$
Hence we have as an equality of \,$\C^m$-vectors
$$ (s_{n,\ell})_{\ell\in I} = -s_{n,n} \, P^{-1} \cdot \left( \int_{A_k} \wt{\omega}_n \right)_{k\in I} \; . $$
Because the \,$\Phi_{n,\xi_{n,k}}$\, are uniformly bounded on the ``compact part of \,$\Sigma$\,'' \,$\bigcup_{|k|\leq N} \wh{S}_k$\,, and \,$s_{n,n}$\, is of the order of \,$n^{-1}$\,, it follows from this equation
that also \,$s_{n,\ell}$\, is with respect to \,$n$\, of the order \,$n^{-1}$\,. 
\end{proof}

After having constructed the canonical basis \,$(\omega_n)$\, of the 1-forms in the preceding theorem, the next step for the construction of the Jacobi variety and Abel map for \,$\Sigma$\, is to
define the 
the Banach space \,$\wt{\Jac}(\Sigma)$\, which will serve as the covering space of the Jacobi variety (and also plays the role of a tangent space for \,$\Jac(\Sigma)$\,), 
and the map \,$\wt{\vi}$\, which will turn out to be the lift of the Abel map
\,$\vi:\Div(\Sigma') \to \Jac(\Sigma)$\,; we will also study the asymptotic behavior of these objects.

In the construction of the Abel map we only permit divisors whose support is contained in the regular set \,$\Sigma'$\, of
\,$\Sigma$\,. One way to handle divisor points in singularities would be to desingularize \,$\Sigma$\, at these points
(as \textsc{Hitchin} does when constructing his spectral curve in \cite{Hitchin:1990}). By this process, one would remove
the \,$\omega_n$\, and the coordinates of the Jacobi variety corresponding to singular divisor points, and obtain
a Jacobi variety associated to a Banach space for the remaining coordinates.

\begin{prop}
\label{P:jacobi:jaccoord}
We fix a divisor \,$D^o \in \Div(\Sigma')$\, (which will be used as the origin point for the definition of the Abel map) 
and denote by \,$\mathfrak{C}_{D^o}$\, the set of sequences \,$(\gamma_k)_{k\in \Z}$\, where each \,$\gamma_k$\, is a curve in \,$\Sigma'$\, running from
a point \,$(\lambda_k^o,\mu_k^o)\in \Sigma'$\, to another point \,$(\lambda_k,\mu_k) \in \Sigma'$\,, such that \,$(\lambda_k^o,\mu_k^o)_{k\in \Z}$\, equals the support of \,$D_0$\, and 
\,$D:=\Menge{(\lambda_k,\mu_k)}{k\in \Z}\in \Div(\Sigma')$\, holds; moreover for large \,$|k|$\, the curve \,$\gamma_k$\, runs entirely in \,$\wh{U}_{k,\delta}$\,,
and there is a number \,$m_\gamma\in \N$\, (depending on \,$\gamma$\, but not on \,$k$\,) so that the winding number of any \,$\gamma_k$\, around any branch point or puncture of \,$\Sigma'$\,
is at most \,$m_\gamma$\,. In this situation we call \,$D$\, the divisor induced by the sequence \,$(\gamma_k)_{k\in \Z}$\,. 

We also let \,$(\omega_n)_{n\in \Z}$\, be the basis of 1-forms from Theorem~\ref{T:jacobi:canonical}(2).
\begin{enumerate}
\item 
For \,$n\in \Z$\, and \,$(\gamma_k)_{k\in \Z} \in \mathfrak{C}_{D^o}$\, the infinite sum
\begin{equation}
\label{eq:jacobi:jaccoord:infinite-sum}
\sum_{k\in \Z} \int_{\gamma_k} \omega_n
\end{equation}
converges absolutely in \,$\C$\,, and we define the map
\begin{equation}
\label{eq:jacobi:jaccoord:wtvi-def}
\wt{\vi}_n : \mathfrak{C}_{D^o} \to \C,\;(\gamma_k)_{k\in \Z} \mapsto \sum_{k\in \Z} \int_{\gamma_k} \omega_n \;. 
\end{equation}

For every \,$(\gamma_k) \in \mathfrak{C}_{D^o}$\, there exist sequences \,$(b_n),(c_n) \in \ell^2(n)$\, so that for \,$n\in \Z$\, with \,$|n|>N$\, (where \,$N\in \N$\, is as in Theorem~\ref{T:jacobi:canonical}(2))
we have
\begin{equation}
\label{eq:jacobi:jaccoord:wtvi-asymp}
\wt{\vi}_n((\gamma_k)) = \frac{\sign(n)}{2\pi i}\cdot \ln\left( \frac{\lambda_n-\vkap_{n,*}+\Psi_n(\lambda_n,\mu_n)}{\lambda_n^o-\vkap_{n,*}+\Psi_n(\lambda_n^o,\mu_n^o)} \right) \cdot (1+b_n) + c_n \; ,
\end{equation}
where \,$\Psi_n$\, is as in Lemma~\ref{L:jacobiprep:Psi-intro}, and 
\,$\ln(z)$\, is the branch of the complex logarithm function with \,$\ln(1)=2\pi i m_n$\,
with \,$m_n\in \Z$\, being the winding number of \,$\gamma_n$\, around the pair of branch points \,$\vkap_{n,1}$\,, \,$\vkap_{n,2}$\,.


\item We consider the (if \,$S\neq \varnothing$\,, non-Hausdorff) topological vector space
$$ \wt{\Jac}(\Sigma) := \Mengegr{(a_n)_{n\in \Z}}{a_n \cdot (\vkap_{n,1}-\vkap_{n,2}) \in \ell^2_{-1,3}(n)}\,, $$
whose topology is induced by the semi-norm
$$ \|a_n\|_{\wt{\Jac}(\Sigma)} := \left\|a_n \cdot (\vkap_{n,1}-\vkap_{n,2})\right\|_{\ell^2_{-1,3}} \qmq{for \,$(a_n) \in \wt{\Jac}(\Sigma)$\,.} $$
If \,$S=\varnothing$\, holds, then \,$\|\,\cdot\,\|_{\wt{\Jac}(\Sigma)}$\, is a norm and \,$\wt{\Jac}(\Sigma)$\, is a Banach space. 
We have \,$\ell^\infty(\Z) \subset \wt{\Jac}(\Sigma)$\,,  and for any \,$(\gamma_k) \in \mathfrak{C}_{D^o}$\,
$$ \left( \wt{\vi}_n((\gamma_k)) \right)_{n\in\Z} \in \wt{\Jac}(\Sigma) \; . $$
Thereby we obtain the map
$$ \wt{\vi}: \mathfrak{C}_{D^o} \to \wt{\Jac}(\Sigma),\; (\gamma_k) \mapsto (\wt{\vi}_n(\gamma_k))_{n\in \Z} \; . $$
\item For any \,$k\in \Z\setminus S$\, we have
$$ \left( \int_{B_k} \omega_n \right)_{n\in \Z} \in \wt{\Jac}(\Sigma) \; . $$
\end{enumerate}
\end{prop}

\begin{rem}
The semi-norm and the topology of \,$\wt{\Jac}(\Sigma)$\, do not contain any information on the coordinates \,$a_k$\, 
in \,$\wt{\Jac}(\Sigma)$\, with \,$k\in S$\,. This means that the present construction of the Jacobi variety
is not very well suited for studying spectral curves with infinitely many singularities, when the corresponding
divisor points are outside these singularities (or at least, when infinitely many of them are). 
Such divisors occur, for example, as spectral divisors for minimal
surfaces in \,$S^3$\, consisting of infinitely many bubbletons. 
(On the spectral curve \,$\Sigma_0$\, of the vacuum, a bubbleton is induced by moving one divisor point that was
in a singularity of \,$\Sigma_0$\, away from this singularity.)
\end{rem}

\begin{proof}[Proof of Proposition~\ref{P:jacobi:jaccoord}.]
\emph{For (1).} 
By Theorem~\ref{T:jacobi:canonical}, we have
$$ \omega_n = \frac{\Phi_n(\lambda)}{\mu-\mu^{-1}}\,\mathrm{d}\lambda \;, $$
where the holomorphic function \,$\Phi_n$\, is a finite linear combination of functions \,$\wt{\Phi}_\ell$\,, 
\,$\ell\in \{-N,\dotsc,N\}\cup\{n\}$\, of the type studied in Proposition~\ref{P:jacobi:asymp-Phi}. 
In other words, \,$\omega_n$\, itself is a linear combination of 1-forms \,$\wt{\omega}_\ell$\, (with \,$\ell\in \{-N,\dotsc,N\} \cup \{n\}$\,) of the type studied in Proposition~\ref{P:jacobi:asymp-omega}. 
The corresponding sequences
\,$(\xi_{\ell,k})_{k\in\Z\setminus\{\ell\}}$\, of zeros of \,$\wt{\Phi}_\ell$\, satisfy the estimate \,$|\xi_{\ell,k}-\vkap_{k,*}|
\leq C_\xi\cdot |\vkap_{k,1}-\vkap_{k,2}|$\, with a constant \,$C_\xi>0$\, that is independent of \,$\ell$\, and \,$k$\,,
and therefore the estimates in Proposition~\ref{P:jacobi:asymp-Phi}--\ref{P:jacobi:asymp-omega} 
apply uniformly to all \,$\wt{\Phi}_\ell$\, resp.~\,$\wt{\omega}_\ell$\,. 

For the proof of the convergence of the infinite sum \eqref{eq:jacobi:jaccoord:infinite-sum} for some given \,$n\in \Z$\,,
we note that we have
$$ \int_{(\lambda_k^o,\mu_k^o)}^{(\lambda_k,\mu_k)} \omega_n = \int_{\vkap_{k,1}}^{(\lambda_k,\mu_k)} \omega_n - \int_{\vkap_{k,1}}^{(\lambda_k^o,\mu_k^o)} \omega_n \;, $$
where the paths of integration are chosen suitably on the right hand side. 
Therefore we may suppose without loss of generality that \,$(\lambda_k^o,\mu_k^o)=\vkap_{k,1}$\, holds for all \,$k\neq n$\,. 
Let \,$(\gamma_k)\in \mathfrak{C}_{D^o}$\,
be given. Because we have \,$\int_{A_k} \wt{\omega}_\ell=0$\, for \,$|k|$\, large, \,$k\neq n$\,, 
the question of the convergence of the sum \,$\sum_{k\in \Z\setminus\{n\}} \left| \int_{\gamma_k} \wt{\omega}_\ell \right|$\, does not depend
on how many times the curve \,$\gamma_k$\, winds around the pair of branch points \,$\vkap_{k,1}$\,, \,$\vkap_{k,2}$\,, and therefore 
Proposition~\ref{P:jacobi:asymp-omega}(3) shows that the sum \,$\sum_{k\in \Z\setminus\{n\}} \left| \int_{\gamma_k} \wt{\omega}_\ell \right|$\, converges.
Because \,$\omega_n$\, is a finite linear combination of the \,$\wt{\omega}_\ell$\,, it follows that the
infinite sum \eqref{eq:jacobi:jaccoord:infinite-sum} converges absolutely.

For proving the asymptotic estimate~\eqref{eq:jacobi:jaccoord:wtvi-asymp} we let \,$(\gamma_k) \in \mathfrak{C}_{D^o}$\, be given and consider \,$n\in \Z$\, with \,$|n|>N$\,.
The constants \,$C_k>0$\, occurring in the sequel depend on \,$(\gamma_k)$\, but not on \,$n$\,. 
We have
\begin{equation}
\label{eq:jacobi:jaccoord:wtvi-gamma-decomp}
\wt{\vi}_n((\gamma_k)) = \int_{\gamma_n}\omega_n + \sum_{|k|\leq N} \int_{\gamma_k} \omega_n + \sum_{\substack{|k|>N \\ k\neq n}} \int_{\gamma_k} \omega_n \; . 
\end{equation}
In the sequel we will estimate the three summands on the right hand side of this equation separately. For this purpose we will use the fact that by
Theorem~\ref{T:jacobi:canonical}(2) there exist constants \,$s_{n,n}\in\C$\, and \,$s_{n,\ell}\in\C$\, for \,$|\ell|\leq N$\, so that
\begin{equation}
\label{eq:jacobi:jaccoord:omega-n}
\omega_n = s_{n,n}\cdot \wt{\omega}_n + \sum_{|\ell|\leq N} s_{n,\ell}\cdot \wt{\omega}_\ell 
\end{equation}
holds; here the constants \,$s_{n,n}$\, and \,$s_{n,\ell}$\, are of order \,$n^{-1}$\,, hence there exists \,$C_1>0$\, so that we have
\begin{equation}
\label{eq:jacobi:jaccoord:snell}
|s_{n,n}|\,,\; |s_{n,\ell}| \leq C_1 \cdot \frac{1}{|n|} \;. 
\end{equation}

For the first summand in \eqref{eq:jacobi:jaccoord:wtvi-gamma-decomp} we note that 
by Theorem~\ref{T:jacobi:canonical}(2) there exists a sequence \,$(b_n^{[1]})\in \ell^2_{1,1}(n)$\, so that
\begin{equation*}
s_{n,n} = 
\begin{cases}
\frac{1}{4\pi} \cdot \lambda_{n,0}^{-1/2} + b_n^{[1]} & \text{for \,$n>N$\,} \\
\frac{1}{4\pi} \cdot \tau_{\xi_{n,k}}^{2}\cdot \lambda_{n,0}^{1/2} + b_n^{[1]} & \text{for \,$n<-N$\,}
\end{cases} 
\end{equation*}
holds, and by Proposition~\ref{P:jacobi:asymp-omega}(1) there exists a sequence \,$(b_n^{[2]}) \in \ell^2_{-1,-1}(n)$\, so that we have
$$ \int_{\gamma_n} \wt{\omega}_n =  
\begin{cases}
(-2i\,\lambda_{n,0}^{1/2}+b_n^{[2]})\cdot \ln\left( \frac{\lambda_n-\vkap_{n,*}+\Psi_n(\lambda_n,\mu_n)}{\lambda_n^o-\vkap_{n,*}+\Psi_n(\lambda_n^o,\mu_n^o)} \right) & \text{for \,$n>N$\,} \\
\left( \frac{2i}{\tau_{\xi_{n,k}}^2}\cdot \lambda_{n,0}^{-1/2}+ b_n^{[2]} \right) \cdot \ln\left( \frac{\lambda_n-\vkap_{n,*}+\Psi_n(\lambda_n,\mu_n)}{\lambda_n^o-\vkap_{n,*}+\Psi_n(\lambda_n^o,\mu_n^o)} \right) & \text{for \,$n<N$\,}
\end{cases} \;, $$
where the choice of the branch of the complex logarithm corresponds to the winding number \,$m_n$\, of \,$\gamma_n$\, around the pair of branch points \,$\vkap_{n,1}$\,, \,$\vkap_{n,2}$\, by \,$\ln(1)=2\pi i m_n$\,.
By multiplying the preceding two equations, we obtain that there exist sequences \,$(b_n^{[3]}),(b_n) \in \ell^2_{0,0}(n)$\, with
\begin{align}
s_{n,n}\cdot \int_{\gamma_n} \wt{\omega}_n 
& = \begin{cases} 
\left( \tfrac{1}{2\pi i} + b_n^{[3]} \right)\cdot \ln\left( \frac{\lambda_n-\vkap_{n,*}+\Psi_n(\lambda_n,\mu_n)}{\lambda_n^o-\vkap_{n,*}+\Psi_n(\lambda_n^o,\mu_n^o)} \right) & \text{for \,$n>N$\,} \\
\left( -\tfrac{1}{2\pi i} + b_n^{[3]} \right)\cdot \ln\left( \frac{\lambda_n-\vkap_{n,*}+\Psi_n(\lambda_n,\mu_n)}{\lambda_n^o-\vkap_{n,*}+\Psi_n(\lambda_n^o,\mu_n^o)} \right) & \text{for \,$n<-N$\,} 
\end{cases} \notag \\
\label{eq:jacobi:jaccoord:bn}
& = \frac{\sign(n)}{2\pi i}\cdot \ln\left( \frac{\lambda_n-\vkap_{n,*}+\Psi_n(\lambda_n,\mu_n)}{\lambda_n^o-\vkap_{n,*}+\Psi_n(\lambda_n^o,\mu_n^o)} \right) \cdot (1+b_n) \; .
\end{align}
Moreover because of \,$\int_{A_n} \wt{\omega}_\ell=0$\, for \,$|\ell|\leq N$\,, the value of \,$\left|\int_{\gamma_n} \wt{\omega}_\ell\right|$\, does not depend on the winding number of \,$\gamma_n$\,
around the pair of branch points \,$\vkap_{n,1}$\,, \,$\vkap_{n,2}$\,. Therefore it follows from Proposition~\ref{P:jacobi:asymp-omega}(2) that there exists a constant \,$C_2>0$\, with
\,$\left|\int_{\gamma_n} \wt{\omega}_\ell\right| \leq C_2 \cdot \tfrac{1}{n^2}$\, for \,$|\ell|\leq N$\,. Because of \eqref{eq:jacobi:jaccoord:snell} it follows that there exists a constant \,$C_3>0$\, with
$$ \sum_{|\ell|\leq N} |s_{n,\ell}| \cdot \left|\int_{\gamma_n} \wt{\omega}_\ell \right| \leq C_3 \cdot \frac{1}{n^3} \; . $$
We thus obtain by combining the preceding estimates with Equation~\eqref{eq:jacobi:jaccoord:omega-n} that there exists a sequence \,$(c_n^{[1]}) \in \ell^2_{0,0}(n)$\, with
\begin{equation}
\label{eq:jacobi:jaccoord:part1}
\int_{\gamma_n}\omega_n = \frac{\sign(n)}{2\pi i}\cdot \ln\left( \frac{\lambda_n-\vkap_{n,*}+\Psi_n(\lambda_n,\mu_n)}{\lambda_n^o-\vkap_{n,*}+\Psi_n(\lambda_n^o,\mu_n^o)} \right) \cdot (1+b_n) + c_n^{[1]} \; . 
\end{equation}

For the second summand in \eqref{eq:jacobi:jaccoord:wtvi-gamma-decomp} we note that by applying Proposition~\ref{P:jacobi:asymp-omega}(4) to all the finitely many paths of integration \,$\gamma_k$\, with \,$|k|\leq N$\,
it follows that there exists \,$C_4>0$\, so that we have for all \,$|k|\leq N$\, and \,$\ell\in \Z$\,
$$ \left| \int_{\gamma_k} \wt{\omega}_\ell \right| \leq C_4 \;; $$
via Equation~\eqref{eq:jacobi:jaccoord:omega-n} we obtain
\begin{align}
\left| \sum_{|k|\leq N} \int_{\gamma_k} \omega_n \right|
& = |s_{n,n}| \cdot \sum_{|k|\leq N} \left|\int_{\gamma_k} \wt{\omega}_n \right| + \sum_{|\ell|\leq N} |s_{n,\ell}| \cdot \left( \sum_{|k|\leq N} \left|\int_{\gamma_k} \wt{\omega}_\ell \right|\right) \notag \\
& \leq \left( |s_{n,n}| + \sum_{|\ell|\leq N} |s_{n,\ell}| \right) \cdot (2N+1) \cdot C_4 \notag \\
\label{eq:jacobi:jaccoord:part2}
& =: c_n^{[2]} \;, 
\end{align}
where the sequence \,$(c_n^{[2]})$\, is of order \,$n^{-1}$\, by \eqref{eq:jacobi:jaccoord:snell}, in particular \,$(c_n^{[2]}) \in \ell^2_{0,0}(n)$\, holds.

For the third summand in \eqref{eq:jacobi:jaccoord:wtvi-gamma-decomp} we note that for \,$|k|>N$\,, \,$k\neq n$\, we have \,$\int_{A_k} \wt{\omega}_\ell=0$\, for all \,$\ell\in \{-N,\dotsc,N\}\cup\{n\}$\,,
and therefore the value of \,$\left| \int_{\gamma_k} \wt{\omega}_\ell \right|$\, does not depend on the winding number of \,$\gamma_k$\, around the pair of branch points \,$\vkap_{k,1}$\,, \,$\vkap_{k,2}$\, for such \,$k,\ell$\,. 
It therefore follows from Proposition~\ref{P:jacobi:asymp-omega}(3) that for \,$\ell\in \{-N,\dotsc,N\}\cup\{n\}$\, we have
$$ \sum_{\substack{|k|>N \\ k\neq n}} \left| \int_{\gamma_k} \wt{\omega}_\ell \right| \;\in\; \ell^2_{-1,-1}(n) \;, $$
via Equations~\eqref{eq:jacobi:jaccoord:omega-n} and \eqref{eq:jacobi:jaccoord:snell} we conclude
\begin{equation}
\label{eq:jacobi:jaccoord:part3}
c_n^{[3]} := \left| \sum_{\substack{|k|>N \\ k\neq n}} \int_{\gamma_k} \omega_n \right| \;\in\; \ell^2_{0,0}(n) \;.
\end{equation}

By plugging the estimates \eqref{eq:jacobi:jaccoord:part1}, \eqref{eq:jacobi:jaccoord:part2} and \eqref{eq:jacobi:jaccoord:part3} into \eqref{eq:jacobi:jaccoord:wtvi-gamma-decomp}, we see
that \eqref{eq:jacobi:jaccoord:wtvi-asymp} holds with \,$b_n \in \ell^2(n)$\, as above and \,$c_n := c_n^{[1]} + c_n^{[2]} + c_n^{[3]} \in \ell^2(n)$\,.

\emph{For (2).}
It is clear that \,$\wt{\Jac}(\Sigma)$\, is a linear space, that \,$\|\,\cdot\,\|_{\wt{\Jac}(\Sigma)}$\, is a semi-norm on \,$\wt{\Jac}(\Sigma)$\, and that in the case \,$S=\varnothing$\,
(i.e.~\,$\vkap_{k,1}-\vkap_{k,2}\neq 0$\, for all \,$k$\,),
this semi-norm becomes a norm and \,$\wt{\Jac}(\Sigma)$\, a Banach space. 
We have \,$\vkap_{k,1}-\vkap_{k,2}\in \ell^2_{-1,3}(k)$\, and therefore \,$\ell^\infty(\Z) \subset \wt{\Jac}(\Sigma)$\,. To prove \,$(\wt{\vi}_n((\gamma_k))) \in \wt{\Jac}(\Sigma)$\,, we need to show that 
$$ \wt{\vi}_n((\gamma_k)) \cdot (\vkap_{n,1}-\vkap_{n,2}) \in \ell^2_{-1,3}(n) $$
holds. If we have \,$n\in S$\, and therefore \,$\vkap_{n,1}=\vkap_{n,2}$\,, there is nothing to show, so we consider \,$n\in \Z\setminus S$\, now. We may again suppose 
without loss of generality that \,$(\lambda_n^o,\mu_n^o) = \vkap_{n,1}$\, holds for \,$n \in \Z\setminus S$\,. Then we have
\begin{align}
\ln\left( \frac{\lambda_n-\vkap_{n,*}+\Psi_n(\lambda_n,\mu_n)}{\lambda_n^o-\vkap_{n,*}+\Psi_n(\lambda_n^o,\mu_n^o)} \right) 
& = \ln\left( \frac{\lambda_n-\vkap_{n,*}+\Psi_n(\lambda_n,\mu_n)}{\tfrac12(\vkap_{k,1}-\vkap_{k,2})} \right) \notag \\
\label{eq:jacobi:jaccoord:arcosh}
& = \arcosh\left( \frac{\lambda_n-\vkap_{n,*}}{\tfrac12(\vkap_{k,1}-\vkap_{k,2})} \right) + 2\pi i \, m_n \; ,
\end{align} 
where \,$m_n \in \Z$\, again is the winding number of \,$\gamma_n$\, around the pair of branch points \,$\vkap_{n,1},\vkap_{n,2}$\,. Here we denote by \,$\ln(z)$\, the branch of the complex logarithm function with 
\,$\ln(1)=2\pi i m_n$\,, and by \,$\arcosh(z)$\, the branch of the complex area cosine hyperbolicus function with \,$\arcosh(1)=0$\,. 

Because of \,$(\gamma_k) \in \mathfrak{C}_{D^o}$\,, the sequence \,$(m_n)_{n\in \Z}$\, is bounded, and therefore 
it follows from (1) and the estimate for \,$\arcosh$\, of Equation~\eqref{eq:jacobiprep:Psi-neu:arcosh-pre} 
that there exist constants \,$C_5,C_6>0$\, with 
\begin{equation}
\label{eq:jacobi:jaccoord:arcosh-estim}
|\wt{\vi}_n((\gamma_k))|
\leq C_5 \cdot \left| \arcosh\left( \frac{\lambda_n-\vkap_{n,*}}{\tfrac12(\vkap_{n,1}-\vkap_{n,2})} \right) \right| +C_6 \leq 2\,C_5 \cdot C_{\arcosh} \cdot \left| \frac{\lambda_n-\vkap_{n,1}}{\vkap_{n,1}-\vkap_{n,2}} \right| +C_6 \; . 
\end{equation}
Thus we have
$$ |\wt{\vi}_n((\gamma_k))| \cdot |\vkap_{n,1}-\vkap_{n,2}| \leq 2\,C_5\,C_{\arcosh}\cdot |\lambda_n-\vkap_{n,1}| + C_6\cdot |\vkap_{n,1}-\vkap_{n,2}| \in \ell^2_{-1,3}(n\in \Z\setminus S) $$
and therefore \,$(\wt{\vi}_n((\gamma_k))) \in \wt{\Jac}(\Sigma)$\,. 

\emph{For (3).}
For every fixed \,$k\in \Z$\,, \,$\int_{B_k} \wt{\omega}_n$\, is bounded with respect to \,$n$\, 
by Proposition~\ref{P:jacobi:asymp-omega}(4). Because both \,$s_{n,n}$\, and the \,$s_{n,\ell}$\, with \,$|\ell|\leq N$\, 
are of order \,$n^{-1}$\, by Theorem~\ref{T:jacobi:canonical}(2),
it follows that
$$ \int_{B_k} \omega_n = s_{n,n} \cdot \int_{B_k} \wt{\omega}_n + \sum_{\substack{|\ell|\leq N}} s_{n,\ell} \int_{B_k} \wt{\omega}_\ell $$
is of order \,$n^{-1}$\,. Thus we have \,$\left(\int_{B_k} \omega_n\right)_{n\in\Z} \in \ell^\infty_{1,1}(n) \subset \ell^\infty(n) \subset \wt{\Jac}(\Sigma)$\, by (2).
\end{proof}

\bigskip

With the following theorem, the construction of the Jacobi variety and the Abel map of \,$\Sigma$\, is completed:

\begin{thm}
\label{T:jacobi:jacobi}
We again fix an ``origin divisor'' \,$D^o \in \Div(\Sigma')$\, and use the notations of Proposition~\ref{P:jacobi:jaccoord}.

\begin{enumerate}
\item
For \,$k,n\in \Z$\, we let
$$ \alpha^{[k]}_n := \int_{A_k} \omega_n = \delta_{k,n} \qmq{and if \,$k\not\in S$\,} \beta^{[k]}_n := \int_{B_k} \omega_n \; . $$
Then we have \,$(\alpha^{[k]}_n)_n \in \wt{\Jac}(\Sigma)$\, for every \,$k\in \Z$\, and \,$(\beta^{[k]}_n)_n \in \wt{\Jac}(\Sigma)$\, for every \,$k\in \Z\setminus S$\,. 

We let \,$\Gamma$\, be the abelian group corresponding to the periods of all closed loops in \,$\mathfrak{C}_{D^o}$\,, i.e.~
$$ \Gamma := \left. \left\{ \sum_{k\in \Z} (a_k\,\alpha^{[k]}+b_k\,\beta^{[k]}) \right| \begin{matrix} a_k,b_k\in \Z \\ \exists N,m\in \N\;\forall k\in \Z,|k|>N\;:\; |a_k|\leq m,\;b_k=0 \end{matrix} \right\} \; . $$
Then \,$\Gamma$\, is an abelian subgroup of \,$\wt{\Jac}(\Sigma)$\,. 
We will call \,$\Gamma$\, the \emph{period lattice} of \,$\Sigma$\, (although it is not a discrete subset of \,$\wt{\Jac}(\Sigma)$\,).

We call the topological quotient space \,$\Jac(\Sigma) := \wt{\Jac}(\Sigma) / \Gamma$\, 
the \emph{Jacobi variety} of \,$\Sigma$\,, and we denote the canonical projection by \,$\pi: \wt{\Jac}(\Sigma) \to \Jac(\Sigma)$\,. 

\item
Let \,$\tau: \mathfrak{C}_{D^o} \to \Div(\Sigma')$\, be the map that associates to each \,$(\gamma_k)\in \mathfrak{C}_{D^o}$\, the divisor
\,$D\in \Div$\, induced by \,$(\gamma_k)$\,. 
Then there exists one and only one map \,$\vi: \Div(\Sigma') \to \Jac(\Sigma)$\, with \,$\vi\circ \tau = \pi \circ \wt{\vi}$\,,
i.e.~so that the following diagram commutes:
\begin{equation*}
\begin{minipage}{5cm}
\begin{xy}
\xymatrix{
\mathfrak{C}_{D^o} \ar[r]^{\wt{\vi}} \ar[d]_{\tau} & \wt{\Jac}(\Sigma) \ar[d]^{\pi} \\
\Div(\Sigma') \ar[r]_{\vi} & \Jac(\Sigma)
}
\end{xy}
\end{minipage}
\end{equation*}
We call \,$\vi$\, the \emph{Abel map} of \,$\Sigma$\,. 

\item
Change of the origin divisor \,$D^o$\, corresponds to a linear translation of \,$\vi$\,. More specifically, if we let
\,$D^{oo} \in \Div(\Sigma')$\, be another divisor, and denote the Abel maps with the origin divisor \,$D^o$\, resp.~\,$D^{oo}$\,
by \,$\vi_{D^o}$\, resp.~by \,$\vi_{D^{oo}}$\,, then we have for any \,$D\in \Div(\Sigma')$\,
$$ \vi_{D^{oo}}(D) = \vi_{D^o}(D) - \vi_{D^o}(D^{oo}) \; . $$
%
%
%
\end{enumerate}
\end{thm}

\begin{proof}[Proof of Theorem~\ref{T:jacobi:jacobi}.]
\emph{For (1).}
Because each sequence \,$(\alpha^{[k]}_n)_{n\in \Z}$\, has exactly one non-zero element, we have \,$(\alpha^{[k]}_n)_{n\in \Z} \in \ell^\infty(n)
\subset \wt{\Jac}(\Sigma)$\, by Proposition~\ref{P:jacobi:jaccoord}(2). Moreover, we have \,$(\beta^{[k]}_n)_{n\in \Z}
\in \wt{\Jac}(\Sigma)$\, by Proposition~\ref{P:jacobi:jaccoord}(3). 

Because of the definition of \,$\mathfrak{C}_{D^o}$\, (in Proposition~\ref{P:jacobi:jaccoord}) and the fact that the homology group of \,$\wh{U}_{k,\delta}$\, is generated by \,$A_k$\,
for \,$|k|$\, large, it is clear that \,$\Gamma$\, is the group of periods in \,$\mathfrak{C}_{D^o}$\,. Because \,$(\alpha_n^{[k]})_{n\in \Z}$\, is the \,$k$-th standard unit vector
of the space of sequences
\,$\C^{\Z}$\,, it also follows from the definition of \,$\Gamma$\, that \,$\Gamma \subset \ell^\infty(n) \subset \wt{\Jac}(\Sigma)$\,; for the inclusion \,$\ell^\infty \subset \wt{\Jac}(\Sigma)$\,
see Proposition~\ref{P:jacobi:jaccoord}(2).


\emph{For (2).}
Let \,$(\gamma_k),(\wt{\gamma}_k)\in \mathfrak{C}_{D^o}$\, be given such that \,$\tau((\gamma_k))=\tau((\wt{\gamma}_k))=:D\in\Div(\Sigma')$\, holds. 
We need to show \,$\wt{\vi}(({\gamma}_k))-\wt{\vi}((\wt{\gamma}_k))\in\Gamma$\,.
Without loss of generality we may suppose that the curves
\,$\gamma_k$\, and \,$\wt{\gamma}_k$\, are parameterized on the interval \,$[0,1]$\,. Because \,$(\gamma_k)$\, and \,$(\wt{\gamma}_k)$\, induce the same divisor \,$D$\,,
there exist permutations \,$h_0,h_1: \Z\to\Z$\, such that for every \,$k\in \Z$\, we have \,$\gamma_k(0)=\gamma_{h_0(k)}(0)=:(\lambda_k^o,\mu_k^o)$\, and
\,$\gamma_k(1)=\gamma_{h_1(k)}(1)=:(\lambda_k,\mu_k)$\,; here we have \,$D^o=\{(\lambda_k^o,\mu_k^o)\}$\, and \,$D=\{(\lambda_k,\mu_k)\}$\,. 
By the definition of \,$\mathfrak{C}_{D^o}$\, (in Proposition~\ref{P:jacobi:jaccoord}) there exists 
\,$N\in \N$\, so that \,$\gamma_k$\, and \,$\wt{\gamma}_k$\, run entirely in \,$\wh{U}_{k,\delta}$\, for \,$|k|>N$\,. 
Then \,$h_0(k)=h_1(k)=k$\, holds for all \,$k\in \Z$\, with \,$|k|>N$\,, and therefore \,$h_0|\{-N,\dotsc,N\}$\, and \,$h_1|\{-N,\dotsc,N\}$\, are  permutations
of \,$\{-N,\dotsc,N\}$\,. 

It follows that the curves \,$\gamma_k(t)$\, and the curves \,$\wt{\gamma}_k(1-t)$\,, where \,$k$\, runs through all of \,$\{-N,\dotsc,N\}$\,, together
define a cycle \,$Z_*$\, on \,$\Sigma$\,, hence \,$Z_*$\, is a \,$\Z$-linear combination of some \,$A_k$\, and \,$B_k$\,, say
$$ Z_* = \sum_{|j|\leq N'} (a_j\,A_j + b_j\,B_j) \qmq{with} N'\in \N, \; a_j,b_j \in \Z \;. $$
Moreover for each individual \,$k\in \Z$\, with \,$|k|>N$\,, the curves \,$\gamma_k(t)$\, and \,$\wt{\gamma}_k(1-t)$\,
form a cycle \,$Z_k$\, on \,$\Sigma$\, that runs entirely in \,$\wh{U}_{k,\delta}$\,, and is therefore a multiple of \,$A_k$\,, say \,$Z_k=a_k'\cdot A_k$\, with \,$a_k'\in \Z$\,; here 
there exists a constant \,$m>0$\, with \,$|a_k'|\leq m$\, for all \,$k$\, because \,$(\gamma_k),(\wt{\gamma}_k)\in \mathfrak{C}_{D^o}$\, holds. 
Therefore
\begin{align*}
\wt{\vi}((\gamma_k))-\wt{\vi}((\wt{\gamma_k})) & = \left( \int_{Z_*} \omega_n + \sum_{|k|>N} \int_{Z_k} \omega_n \right)_{n\in \Z} \\
& = \sum_{|j|\leq N'} (a_j\,\alpha^{[j]} + b_j\,\beta^{[j]}) + \sum_{|k|>N} a_k'\,\alpha^{[k]} \;\in\; \Gamma 
\end{align*}
holds.


\emph{For (3).} 
We write \,$D=\{(\lambda_k,\mu_k)\}$\,, \,$D^o=\{(\lambda_k^o,\mu_k^o)\}$\, and \,$D^{oo}=\{(\lambda_k^{oo},\mu_k^{oo})\}$\, with the usual asymptotic enumeration of asymptotic divisors.
For every \,$k\in \Z$\, we let \,$\gamma_k^{[1]}$\, be a path in \,$\Sigma$\, connecting \,$(\lambda_k^o,\mu_k^o)$\, to \,$(\lambda_k^{oo},\mu_k^{oo})$\,, and let \,$\gamma_k^{[2]}$\,
be a path connecting \,$(\lambda_k^{oo},\mu_k^{oo})$\, to \,$(\lambda_k,\mu_k)$\,. We choose these paths so that \,$(\gamma_k^{[1]})\in \mathfrak{C}_{D^o}$\, and 
\,$(\gamma_k^{[2]})\in \mathfrak{C}_{D^{oo}}$\, holds. For \,$k\in \Z$\,, let \,$\gamma_k$\, be the concatenation of the paths \,$\gamma_k^{[1]}$\, and \,$\gamma_k^{[2]}$\,,
then we also have \,$(\gamma_k) \in \mathfrak{C}_{D^o}$\,.

For every \,$k\in \Z$\, and \,$n\in \Z$\, we have
$$ \int_{\gamma_k} \omega_n = \int_{\gamma_k^{[1]}} \omega_n + \int_{\gamma_k^{[2]}} \omega_n $$
and therefore
$$ \wt{\vi}_{D^o}((\gamma_k)) = \wt{\vi}_{D^o}((\gamma_k^{[1]})) + \wt{\vi}_{D^{oo}}((\gamma_k^{[2]})) \; . $$
Because the divisor induced by \,$(\gamma_k)$\,, \,$(\gamma_k^{[1]})$\, resp.~\,$(\gamma_k^{[2]})$\, is \,$D$\,, \,$D^{oo}$\, resp.~\,$D$\,, we conclude
$$ \vi_{D^o}(D) = \vi_{D^o}(D^{oo}) + \vi_{D^{oo}}(D) \; . $$
%
%
\end{proof}

\enlargethispage{2em}

\begin{rem}
Note that the ``period lattice'' \,$\Gamma$\, from Theorem~\ref{T:jacobi:jacobi} is not discrete in \,$\wt{\Jac}(\Sigma)$\,, 
not even for \,$S=\varnothing$\,,
and is therefore in fact not a lattice. (We have \,$\|(\alpha^{[k]}_n)_n\|_{\wt{\Jac}(\Sigma)} = \tfrac{1}{k}\cdot |\vkap_{k,1}-\vkap_{k,2}| \in \ell^2(k)$\,
for \,$k>0$\,, and therefore the \,$(\alpha^{[k]}_n)_n \in \Gamma$\, accumulate near \,$0\in \wt{\Jac}(\Sigma)$\,.) For this reason we view the Jacobi variety \,$\Jac(\Sigma)=\wt{\Jac}(\Sigma)/\Gamma$\,
only as a quotient topological space, not as an (infinite-dimensional) manifold. The situation is similar to the one encountered by \textsc{McKean} and \textsc{Trubowitz} in \cite{McKean/Trubowitz:1976} concerning the Jacobi
variety for the integrable system associated to Hill's operator: 
There the period lattice is also not discrete in the respective Banach space, and the Jacobi variety is compact (topologically, it is a product of infinitely many circles)
and therefore does not carry the structure of an infinite dimensional manifold, see the discussion in \cite{McKean/Trubowitz:1976}, p.~154.

Another peculiarity concerning our period lattice \,$\Gamma$\, is that in the case \,$S\neq \varnothing$\, it no longer spans \,$\wt{\Jac}(\Sigma)$\, over \,$\R$\,. (So even if \,$\Gamma$\, were discrete,
it would not be a \emph{maximal} discrete abelian subgroup of \,$\wt{\Jac}(\Sigma)$\,.)
The reason is that for \,$k\in S$\,, there is no cycle \,$B_k$\,. Geometrically, one can think of a family of spectral curves with \,$k\not\in S$\, where the 
two branch points \,$\vkap_{k,1}$\, and \,$\vkap_{k,2}$\, ``close up'' to each other; in the limit we then have \,$\vkap_{k,1}=\vkap_{k,2}$\, and therefore \,$k\in S$\,. This process of ``closing up''
a pair of branch points causes the period corresponding to the cycle \,$B_k$\, to go to infinity, and this is the reason why there is no period in \,$\Gamma$\, corresponding to \,$B_k$\, in 
the limit spectral curve, where \,$k\in S$\, holds. 
As a consequence, \,$\Jac(\Sigma)$\, is not compact for \,$S\neq\varnothing$\,. 
%
\end{rem}

\begin{rem}
\label{R:jacobi:vidiffeo}
The Abel map \,$\vi: \Div(\Sigma') \to \Jac(\Sigma)$\, described in Theorem~\ref{T:jacobi:jacobi} is not a a local diffeomorphism in general (not even near non-special divisors), because
at least for a divisor \,$D=\{(\lambda_k,\mu_k)\} \in \Div(\Sigma')$\, that does \emph{not} satisfy the condition 
\begin{equation}
\label{eq:jacobi:vidiffeo:condition}
|\lambda_k-\vkap_{k,*}|\leq C_\xi\cdot |\vkap_{k,1}-\vkap_{k,2}|
\end{equation}
for any \,$C_\xi>0$\,, the Banach space structures of \,$\Div(\Sigma')$\, near \,$D$\, (locally diffeomorphic to \,$\ell^2_{-1,3}$\,) and of \,$\Jac(\Sigma)$\, near \,$\vi(D)$\,
(locally diffeomorphic to \,$\wt{\Jac}(\Sigma)$\,) do not match.

For some applications it would probably be preferable to have a Jacobi variety and associated Abel map for \,$\Sigma$\, that is a local diffeomorphism near every non-special, asymptotic divisor.
This is true in particular where one is interested in divisors with points close to (but not in) singularities of the spectral curve, as they occur for example as spectral divisors
of bubbletons (note that the condition \eqref{eq:jacobi:vidiffeo:condition} forces \,$\lambda_k=\vkap_{k,*}$\, for all \,$k\in S$\,). To obtain a Jacobi variety of this kind,
one would need to expand the space of 1-forms considered; more specifically, one might
work with 1-forms that are square-integrable only on \,$\wh{V}_\delta$\, for a fixed \,$\delta>0$\, and are regular on \,$\Sigma$\,. This would permit one to do away with the condition
\eqref{eq:jacobi:vidiffeo:condition}, also for the sequence \,$(\xi_{n,k})$\, of the 1-forms \,$\wt{\omega}_n$\, of Theorem~\ref{T:jacobi:canonical}(1). 
However the proof of the analogue of Theorem~\ref{T:jacobi:canonical} would become more difficult in this setting. 
\end{rem}

The space \,$\wt{\Jac}(\Sigma)$\, plays the role of a tangent space for the Jacobi variety \,$\Jac(\Sigma)$\,. In our setting where the period lattice \,$\Gamma$\, is not discrete,
the tangent space of \,$\Jac(\Sigma)$\, is not unique however, and similarly as it is the case for Hill's equation as studied by \textsc{McKean} and \textsc{Trubowitz} in \cite{McKean/Trubowitz:1976}, 
we need to pass to a larger tangent space so that the flow of translations of the potential (which we will study via the Jacobi variety in Section~\ref{Se:jacobitrans})
is tangential to \,$\Jac(\Sigma)$\,. This corresponds to a larger space of curve tuples \,$\mathfrak{C}_{D^o}$\,
and a larger Banach space \,$\wt{\Jac}(\Sigma)$\,, as we construct in the following proposition.

In fact McKean and Trubowitz construct in \cite{McKean/Trubowitz:1976} an entire ascending family of tangent spaces for their Jacobi variety, which correspond to the higher flows
of the integrable system associated with Hill's equation. In our setting we cannot define more than the first extension described in the following proposition, because
our potentials are only once weakly differentiable, in contrast to the infinitely differentiable potentials in \cite{McKean/Trubowitz:1976}.

\begin{prop}
\label{P:jacobi:larger}
We again fix an origin divisor \,$D^o \in \Div(\Sigma')$\,. We denote by
\,$\mathfrak{C}_{D^o,1}$\, the set of sequences \,$(\gamma_k)_{k\in \Z}$\, where each \,$\gamma_k$\, is a curve in \,$\Sigma'$\, running from
a point \,$(\lambda_k^o,\mu_k^o)\in \Sigma'$\, to another point \,$(\lambda_k,\mu_k) \in \Sigma'$\,, such that \,$(\lambda_k^o,\mu_k^o)_{k\in \Z}$\, equals the support of \,$D_0$\, and 
\,$D:=\Menge{(\lambda_k,\mu_k)}{k\in \Z}\in \Div(\Sigma')$\, holds; moreover for large \,$|k|$\, the curve \,$\gamma_k$\, runs entirely in \,$\wh{U}_{k,\delta}$\,,
and there is a number \,$m_\gamma\in \N$\, (depending on \,$\gamma$\, but not on \,$k$\,) so that the winding number of any \,$\gamma_k$\, around any branch point or puncture of \,$\Sigma'$\,
is at most \,$m_\gamma\cdot |k|$\,. 

For \,$n\in \Z$\, and \,$(\gamma_k)_{k\in \Z} \in \mathfrak{C}_{D^o}^{(1)}$\, the infinite sum
\begin{equation}
\label{eq:jacobi:better:infinite-sum}
\sum_{k\in \Z} \int_{\gamma_k} \omega_n
\end{equation}
converges absolutely in \,$\C$\,, and we define the map
\begin{equation}
\label{eq:jacobi:better:wtvi-def}
\wt{\vi}_{n}^{(1)} : \mathfrak{C}_{D^o}^{(1)} \to \C,\;(\gamma_k)_{k\in \Z} \mapsto \sum_{k\in \Z} \int_{\gamma_k} \omega_n \;. 
\end{equation}
For every \,$(\gamma_k) \in \mathfrak{C}_{D^o}^{(1)}$\, there exist sequences \,$(b_n),(c_n) \in \ell^2(n)$\, 
so that for \,$n\in \Z$\, with \,$|n|>N$\, (where \,$N\in \N$\, is as in Theorem~\ref{T:jacobi:canonical}(2)) we have
\begin{equation}
\label{eq:jacobi:better:wtvi-asymp}
\wt{\vi}_n((\gamma_k)) = \frac{\sign(n)}{2\pi i}\cdot \ln\left( \frac{\lambda_n-\vkap_{n,*}+\Psi_n(\lambda_n,\mu_n)}{\lambda_n^o-\vkap_{n,*}+\Psi_n(\lambda_n^o,\mu_n^o)} \right) \cdot (1+b_n) + c_n \; ,
\end{equation}
where \,$\Psi_n$\, is as in Lemma~\ref{L:jacobiprep:Psi-intro}, and
\,$\ln(z)$\, is the branch of the complex logarithm function with \,$\ln(1)=2\pi i m_n$\,
with \,$m_n\in \Z$\, being the winding number of \,$\gamma_n$\, around the pair of branch points \,$\vkap_{n,1}$\,, \,$\vkap_{n,2}$\,.

Moreover with
$$ \wt{\Jac}^{(1)}(\Sigma) := \Mengegr{(a_n)_{n\in \Z}}{a_n \cdot (\vkap_{n,1}-\vkap_{n,2}) \in \ell^2_{-2,2}(n)}\,, $$
the map \,$\wt{\vi}^{(1)} := (\wt{\vi}_n^{(1)})_{n\in \Z}$\, maps into \,$\wt{\Jac}^{(1)}(\Sigma)$\,. We have \,$\ell^\infty_{-1,-1}(n) \subset \wt{\Jac}^{(1)}(\Sigma)$\, and \,$\wt{\Jac}(\Sigma) \subset
\wt{\Jac}^{(1)}(\Sigma)$\,. 
\end{prop}

\begin{proof}
The proof is mostly analogous to that of Proposition~\ref{P:jacobi:jaccoord}(1),(2). The only substantial difference is in the treatment of \,$\wt{\vi}_n^{(1)}((\gamma_k))$\, in the proof
that \,$\wt{\vi}^{(1)}$\, maps into \,$\wt{\Jac}^{(1)}(\Sigma)$\, in Proposition~\ref{P:jacobi:jaccoord}(2). If we suppose that \,$(\gamma_k) \in \mathfrak{C}_{D^o}^{(1)}$\, is given, we still have
as in Equation~\eqref{eq:jacobi:jaccoord:arcosh}
$$ \ln\left( \frac{\lambda_n-\vkap_{n,*}+\Psi_n(\lambda_n,\mu_n)}{\lambda_n^o-\vkap_{n,*}+\Psi_n(\lambda_n^o,\mu_n^o)} \right) 
= \arcosh\left( \frac{\lambda_n-\vkap_{n,*}}{\tfrac12(\vkap_{k,1}-\vkap_{k,2})} \right) + 2\pi i \, m_n \; , $$
where \,$m_n$\, is the winding number of \,$\gamma_n$\, around the pair of branch points \,$\vkap_{n,1}$\,, \,$\vkap_{n,2}$\,. But in contrast to the situation of Proposition~\ref{P:jacobi:jaccoord},
the sequence \,$(m_n)_{n\in \Z}$\, is no longer bounded here, rather there exists a constant \,$m_\gamma\in \N$\, so that \,$|m_n| \leq m_\gamma \cdot |n|$\, holds for all \,$n$\,. Thus we obtain instead of 
\eqref{eq:jacobi:jaccoord:arcosh-estim} that there exist \,$C_5,C_6>0$\, so that
$$ |\wt{\vi}_n^{(1)}((\gamma_k))|
\leq 
C_5 \cdot \left| \arcosh\left( \frac{\lambda_n-\vkap_{n,*}}{\tfrac12(\vkap_{n,1}-\vkap_{n,2})} \right) \right| +C_6\cdot|n| \leq 2\,C_5 \cdot C_{\arcosh} \cdot \left| \frac{\lambda_n-\vkap_{n,1}}{\vkap_{n,1}-\vkap_{n,2}} \right| +C_6\cdot|n| $$
and therefore 
$$ |\wt{\vi}_n^{(1)}((\gamma_k))| \cdot |\vkap_{n,1}-\vkap_{n,2}| \leq 2\,C_5\,C_{\arcosh}\cdot |\lambda_n-\vkap_{n,1}| + C_6\cdot |\vkap_{n,1}-\vkap_{n,2}|\cdot|n| \in \ell^2_{-2,2}(n\in \Z\setminus S) \;,$$
whence \,$(\wt{\vi}_n^{(1)}((\gamma_k))) \in \wt{\Jac}^{(1)}(\Sigma)$\, follows.
\end{proof}

\section{The Jacobi variety and translations of the potential}
\label{Se:jacobitrans}

Let \,$(\Sigma,D)$\, be spectral data corresponding to a simply periodic solution \,$u:X \to \C$\, of the sinh-Gordon equation,
where \,$X\subset \C$\, is a horizontal strip with \,$0\in X$\,. 

Throughout this entire work, we have constructed the spectral data \,$(\Sigma,D)$\, corresponding to \,$u$\, with respect to the monodromy at the base point \,$z=0$\,. In the present section, we would now like to
describe the spectral data \,$(\Sigma_{z_0},D_{z_0})$\, corresponding to the monodromy based at some other point \,$z_0\in X$\,, or equivalently, the spectral data (in the previous sense, via the monodromy
at the base point \,$z=0$\,) of the translated potential \,$z\mapsto u(z+z_0)$\,. 

First we note that the spectral curve does not change at all under such a translation, i.e.~\,$\Sigma_{z_0}=\Sigma$\, holds. Indeed by Equation \eqref{eq:mono:monodromy-dgl}, the monodromy \,$M_{z_0}$\, satisfies
with respect to \,$z_0$\, the differential equation
\begin{equation*}
\mathrm{d}_{z_0} M_{z_0}(\lambda) = [\alpha_\lambda(z_0),M_{z_0}(\lambda)] \;,
\end{equation*}
where \,$\alpha_\lambda$\, is the connection form associated to the potential \,$u$\, as in Equation~\eqref{eq:mono:alphaxy}. If we denote by \,$F_\lambda(z)$\, the extended frame associated to \,$u$\,, 
then \,$\wt{M}_{z_0}(\lambda) := F_\lambda(z_0)\cdot M(\lambda)\cdot F_\lambda(z_0)^{-1}$\, satisfies the same differential equation with the same initial value \,$\wt{M}_{z_0=0}(\lambda) = M(\lambda)$\, at \,$z_0=0$\,. 
Hence we have \,$M_{z_0}(\lambda) = F_\lambda(z_0) \cdot M(\lambda) \cdot F_\lambda(z_0)^{-1}$\,. Therefore the eigenvalues of the monodromy, and thus the spectral curve, do not depend on \,$z_0$\,. 

But of course, the spectral divisor \,$D= \Menge{(\lambda_k,\mu_k)}{k\in\Z}$\, does change under translation. We will describe the motion of the divisor points on \,$\Sigma$\, in terms of the image of
the divisor under the Abel map, i.e.~in the Jacobi variety, via ordinary differential equations for the translation in \,$x$-direction and in \,$y$-direction. 
Similarly as is well-known for finite-type potentials, 
it will turn out that both translations correspond to linear motions, that is to constant vector fields, in the Jacobi variety.
In the case of the translation in \,$x$-direction, the direction of that linear motion is precisely a lattice direction of the Jacobi variety, corresponding to the fact that the potential \,$u$\, is 
periodic in the \,$x$-direction. These results are explicated in Theorem~\ref{T:jacobitrans:dgl} below.


In the sequel, we consider either a potential \,$(u,u_y)\in \Pot$\,, or a simply periodic solution of the sinh-Gordon equation \,$u: X \to \C$\, with period \,$1$\, defined on a horizontal strip \,$X\subset \C$\, with \,$0$\, in the 
interior of \,$X$\,.
In the former case we consider only the translation in \,$x$-direction (all references to translations in the direction of \,$y$\, and associated objects are then to be disregarded), and in the latter case,
we consider both translations in \,$x$-direction and in \,$y$-direction. We let \,$(\Sigma,D(x))$\, resp.~\,$(\Sigma,D(y))$\, be the spectral data of the potential translated in \,$x$-direction
resp.~\,$y$-direction, i.e.~of \,$u(z+x)$\, resp.~of \,$u(z+iy)$\,. 
Like in Section~\ref{Se:jacobi} we continue to require that the spectral curve \,$\Sigma$\, does not have any singularities besides ordinary double points, i.e.~that \,$\Delta^2-4$\, does not have
any zeros of order \,$\geq 3$\,, where \,$\Delta=\mu+\mu^{-1}$\, is the trace function of \,$\Sigma$\,. 
We then let
\,$(\omega_n)_{n\in \Z}$\, be the canonical basis of \,$\Omega(\Sigma)$\, (Theorem~\ref{T:jacobi:canonical}(2)), \,$\Jac(\Sigma)$\, be the Jacobi variety and \,$\vi: \Div(\Sigma')\to\Jac(\Sigma)$\,
the Abel map of \,$\Sigma$\, (Theorem~\ref{T:jacobi:jacobi}), where we choose \,$D^o := D(0)$\, as the origin divisor for the construction of the Abel map. 

In the sequel we will look at the derivatives \,$\tfrac{\partial \vi_n}{\partial x}$\, and \,$\tfrac{\partial \vi_n}{\partial y}$\, of the \,$n$-th Jacobi coordinate \,$\vi_n$\,. For these derivatives to make sense,
we need to define Jacobi coordinates \,$\vi_n$\, of \,$D(x)$\, resp.~\,$D(y)$\, at least for small \,$|x|$\, resp.~\,$|y|$\, for all \,$n\in \Z$\,. For this purpose
we write \,$D(x) = \{\lambda_k(x),\mu_k(x)\}$\, for small \,$|x|$\, and then consider for fixed \,$x$\, and all \,$k\in \Z$\, the curve 
\,$\gamma_{x,k}: [0,x] \to \Sigma,\; t\mapsto (\lambda_k(t),\mu_k(t))$\,. Because the map \,$\Pot\to\Div$\, is asymptotically close to the Fourier transform of the potential, \,$\gamma_{x=1,k}$\, winds \,$|k|$\, times around
the pair of branch points \,$\vkap_{k,1}$\,, \,$\vkap_{k,2}$\, for \,$|k|$\, large; it follows that we do not have \,$(\gamma_{x,k}) \in \mathfrak{C}_{D^o}$\, (Proposition~\ref{P:jacobi:jaccoord}),
but we do have \,$(\gamma_{x,k}) \in \mathfrak{C}_{D^o}^{(1)}$\, (Proposition~\ref{P:jacobi:larger}). Therefore we can define Jacobi coordinates for the translation in \,$x$-direction in the vicinity of \,$D(0)$\, by
$$ \vi_n(x) := \wt{\vi}_n^{(1)}(\gamma_{x,k}) \qmq{for \,$n\in \Z$\,,} $$
where \,$\wt{\vi}_n^{(1)}$\, is as in Proposition~\ref{P:jacobi:larger}. A similar construction applies for the translation in \,$y$-direction; here it is relevant that for large \,$|k|$\, the divisor
point \,$(\lambda_k(y),\mu_k(y))$\, remains in the excluded domain \,$\wh{U}_{k,\delta}$\, for sufficiently small \,$|y|$\, because the asymptotic estimates then apply to the translated potentials uniformly.

We denote by \,$\tfrac{\partial \vi_n}{\partial x}$\,
resp.~\,$\tfrac{\partial \vi_n}{\partial y}$\, the derivative of the Jacobi coordinate \,$\vi_n(x)$\, resp.~\,$\vi_n(y)$\, with respect to \,$x$\, resp.~\,$y$\,. 


\begin{thm}
\label{T:jacobitrans:dgl}
There exist sequences \,$a_n^x,a_n^y \in \ell^2_{-1,-1}(n)$\, (dependent only on the spectral curve \,$\Sigma$\,) so that 
under translation of the potential \,$u$\, in the direction of \,$x$\, resp.~\,$y$\,, the Jacobi coordinates \,$\vi_n$\, (\,$n\in \Z$\,) follow the differential equations
\begin{align*}
\frac{\partial \vi_n}{\partial x} & = n + a_n^x \;, \\
\frac{\partial \vi_n}{\partial y} & = -i|n| + a_n^y \; . 
\end{align*}
Moreover we have \,$a_n^x=0$\, for \,$|n|$\, large, and for every \,$n\in \Z$\,, \,$\tfrac{\partial \vi_n}{\partial x}$\, corresponds
to a member of the period lattice, i.e.~there exists a cycle \,$Z_n$\, of \,$\Sigma'$\, so that \,$\tfrac{\partial \vi_n}{\partial x}
= \int_{Z_n} \omega_n$\, holds.
(No similar statements apply to the sequence \,$(a_n^y)$\, in general.)
\end{thm}

The proof of this theorem is the objective of the remainder of the present section. 
%
%
%
%
At the heart of the proof of Theorem~\ref{T:jacobitrans:dgl} is a general construction of linear flows in the Picard variety of a Riemann surface \,$X$\,
(the space of isomorphy classes of line bundles on \,$X$\,) known as the \emph{Krichever construction}. In the sequel we will carry out the proof in our specific situation; we will then discuss
the relationship of the proof to the Krichever construction in Remark~\ref{R:jacobitrans:krichever}.

\begin{prop}
\label{P:jacobitrans:dgl-V1}
Suppose that the spectral divisor \,$D$\, is tame.
Under translation of the potential \,$u$\, in the direction of \,$x$\, resp.~\,$y$\,, the Jacobi coordinates \,$\vi_n$\, (\,$n\in \Z$\,) of the spectral divisor follow the differential equations
\begin{align*}
\frac{\partial \vi_n}{\partial x} & = -\sum_{k\in \Z} \alpha_{21}^x(\lambda_k)\cdot \Phi_n(\lambda_k) \cdot \frac{1}{c'(\lambda_k)} \\
\frac{\partial \vi_n}{\partial y} & = -\sum_{k\in \Z} \alpha_{21}^y(\lambda_k)\cdot \Phi_n(\lambda_k) \cdot \frac{1}{c'(\lambda_k)} 
\end{align*}
with
$$ \alpha_{21}^x(\lambda) = \frac{1}{4}\cdot (\lambda\,\tau+\tau^{-1}) \qmq{and} \alpha_{21}^y(\lambda) = \frac{i}{4}\cdot (-\lambda\,\tau+\tau^{-1}) \; , \qmq{where} \tau := e^{-u/2} \; . $$
\end{prop}

\newpage

\begin{proof}
%
We have
\begin{equation}
\label{eq:jacobitrans:dglpre:ansatz}
\frac{\partial \vi_n}{\partial x} = \sum_{k\in \Z} \frac{\partial \vi_n}{\partial \lambda_k} \cdot \frac{\partial \lambda_k}{\partial x} \; . 
\end{equation}
From the construction of the lift of the Abel map \,$\wt{\vi}^{(1)}$\, in Proposition~\ref{P:jacobi:larger} it follows that the ``canonical map'', i.e.~the derivative of the Jacobi coordinate \,$\vi_n$\,
with respect to the divisor coordinate \,$\lambda_k$\, is given by
\begin{equation}
\label{eq:jacobitrans:dglpre:canonical}
\frac{\partial \vi_n}{\partial \lambda_k} = \frac{\omega_n}{\mathrm{d}\lambda}(\lambda_k,\mu_k) = \frac{\Phi_n(\lambda_k)}{\mu_k-\mu_k^{-1}} \; , 
\end{equation}
where the holomorphic function \,$\Phi_n$\, is as in Theorem~\ref{T:jacobi:canonical}(2).
To calculate \,$\frac{\partial \lambda_k}{\partial x}$\,, we regard \,$c(\lambda)=c(\lambda,x)$\, and \,$\lambda_k=\lambda_k(x)$\, also as functions of \,$x$\,. 
By definition of the spectral divisor, we have \,$c(\lambda_n(x),x)=0$\,, and by differentiating this equation we obtain
\begin{equation}
\label{eq:jacobitrans:dglpre:dellambdax-pre}
\frac{\partial \lambda_k}{\partial x} = - \left( \left. \frac{\partial c}{\partial \lambda}\right|_{\lambda=\lambda_k} \right)^{-1} \cdot \left.\frac{\partial c}{\partial x}\right|_{\lambda=\lambda_k} \; ,
\end{equation}
where \,$\left. \frac{\partial c}{\partial \lambda}\right|_{\lambda=\lambda_k} \neq 0$\, holds because \,$D$\, is tame.
By Equation \eqref{eq:mono:monodromy-dgl}, the monodromy \,$M_{z_0}$\, satisfies
with respect to \,$z_0$\, the differential equation
\begin{equation*}
\mathrm{d}_{z_0} M_{z_0}(\lambda) = [\alpha_\lambda(z_0),M_{z_0}(\lambda)] \;,
\end{equation*}
where \,$\alpha_\lambda$\, is the connection form associated to the potential \,$u$\, as in Equation~\eqref{eq:mono:alphaxy}, and therefore we have 
$$ \frac{\partial c}{\partial x} = ([\alpha,M])_{21} = \alpha_{21}^x\,a + \alpha_{22}^x\,c - \alpha_{11}^x\,c - \alpha_{21}^x\,d \; . $$
Because of \,$c(\lambda_n(x),x)=0$\, it follows that 
$$ \left. \frac{\partial c}{\partial x} \right|_{\lambda=\lambda_n} = \alpha_{21}^x(\lambda_k)\,(a(\lambda_k)-d(\lambda_k)) = \alpha_{21}^x(\lambda_k)\,(\mu_k-\mu_k^{-1}) \; . $$
By plugging this equation into Equation~\eqref{eq:jacobitrans:dglpre:dellambdax-pre}, we obtain
\begin{equation}
\label{eq:jacobitrans:dglpre:dellambdax}
\frac{\partial \lambda_k}{\partial x} = - \frac{1}{c'(\lambda_k)} \cdot \alpha_{21}^x(\lambda_k) \cdot (\mu_k-\mu_k^{-1})\;, 
\end{equation}
and by plugging Equations~\eqref{eq:jacobitrans:dglpre:canonical} and \eqref{eq:jacobitrans:dglpre:dellambdax} into Equation~\eqref{eq:jacobitrans:dglpre:ansatz}, we obtain
\begin{align*}
\frac{\partial \vi_n}{\partial x} & = \sum_{k\in \Z} \frac{\Phi_n(\lambda_k)}{\mu_k-\mu_k^{-1}}\cdot \frac{-1}{c'(\lambda_k)} \cdot \alpha_{21}^x(\lambda_k) \cdot (\mu_k-\mu_k^{-1}) \\
& = - \sum_{k\in \Z} \alpha_{21}^x(\lambda_k) \cdot \Phi_n(\lambda_k) \cdot \frac{1}{c'(\lambda_k)} \; . 
\end{align*}
The differential equation for \,$\tfrac{\partial \vi_n}{\partial y}$\, is proven in literally the same way, by just replacing every \,$x$\, with \,$y$\,. 
\end{proof}

\newpage

We also note the following:

\begin{cor}
\label{C:jacobitrans:transsigma}
If \,$D=\{(\lambda_k,\mu_k)\} \in \Div$\, is given, and \,$\wt{D}=\{(\wt{\lambda}_k,\wt{\mu}_k)\}$\, is the image of \,$D$\, under the hyperelliptic involution \,$\sigma: \Sigma\to\Sigma$\,, then we have
\,$\tfrac{\partial \wt{\lambda}_n}{\partial x} = -\tfrac{\partial {\lambda}_n}{\partial x}$\, and \,$\tfrac{\partial \wt{\lambda}_n}{\partial y} = -\tfrac{\partial {\lambda}_n}{\partial y}$\,.
\end{cor}

\begin{proof}
This follows from Equation~\eqref{eq:jacobitrans:dglpre:dellambdax} resp.~the corresponding equation for the derivative in the direction of \,$y$\,:

We first note that \,$\tau$\, is determined by the \,$\lambda_k$\, up to sign by the trace formula in Theorem~\ref{T:asympfinal:monodromy}(3), and is therefore up to sign invariant under the 
hyperelliptic involution \,$\sigma$\,. Both \,$\alpha_{21}^x(\lambda_n)$\, and \,$c'(\lambda_n)$\, depend only on the \,$\lambda_k$\, and on \,$\tau$\,, and they both change their sign when \,$\tau$\,
changes sign. It follows that \,$-\tfrac{1}{c'(\lambda_k)}\cdot \alpha_{21}^x(\lambda_n)$\, is invariant under \,$\sigma$\,. Because \,$\mu_k-\mu_k^{-1}$\, changes sign under \,$\sigma$\,,
it follows from Equation~\eqref{eq:jacobitrans:dglpre:dellambdax} that \,$\tfrac{\partial \lambda_n}{\partial x}$\, changes sign under \,$\sigma$\,. 
\end{proof}

The presentation of the differential equations for translating the potential in Proposition~\ref{P:jacobitrans:dgl-V1} is not yet satisfactory, because their right hand side apparently varies with \,$\lambda_k$\,,
whereas we are expecting that the Jacobi coordinates change linearly under translation, as explained at the beginning of the section, that is, 
\,$\tfrac{\partial \vi_n}{\partial x}$\, and \,$\tfrac{\partial \vi_n}{\partial y}$\, should be constant. Moreover, the description in Proposition~\ref{P:jacobitrans:dgl-V1} is applicable only if the spectral
divisor \,$D$\, is tame. 
In the following Proposition~\ref{P:jacobitrans:dgl-V2} we will show that the translational flows
of the Jacobi coordinates are indeed linear, and also dispose of the restriction to tame spectral divisors \,$D$\,. 

In preparation for Proposition~\ref{P:jacobitrans:dgl-V2} we define
residues in \,$0$\, and in \,$\infty$\, for meromorphic \,$1$-forms on \,$\wh{V}_\delta$\,. 

\label{not:jacobitrans:Res}
For this purpose we let \,$\wh{K}_r(0)$\, be a cycle in \,$\Sigma$\, that is obtained by lifting a counter-clockwise parameterization of two traversals of a circle of radius \,$r>0$\, around \,$0$\, in \,$\C^*$\,
to \,$\Sigma$\,. 
We also put \,$\wh{K}_r(\infty) := -\wh{K}_{1/r}(0)$\,. Then we choose a sequence \,$(r_n)_{n\geq 1}$\, with \,$r_n>0$\, and \,$\lim_{n\to\infty} r_n=0$\,, such that \,$\wh{K}_{r_n}(0)$\, and \,$\wh{K}_{r_n}(\infty)$\,
are contained in \,$\wh{V}_\delta$\, for all \,$n\geq 1$\,. For any meromorphic 1-form \,$\eta$\, defined on \,$\widehat{V}_\delta$\, 
we then define
\begin{align*}
\Res_0(\eta) & := \frac{1}{2\pi i} \, \lim_{n\to\infty} \int_{\widehat{K}_{r_n}(0)} \eta \\
\qmq{and} \Res_\infty(\eta) & := \frac{1}{2\pi i} \, \lim_{n\to\infty} \int_{\widehat{K}_{r_n}(\infty)} \eta = -\frac{1}{2\pi i} \, \lim_{n\to \infty} \int_{\widehat{K}_{1/r_n}(0)} \eta \;, 
\end{align*}
provided that these limits exist. It is clear that this definition does not depend on the choice of the sequence \,$(r_n)$\,. 
It should be noted that because \,$\wh{K}_r(0)$\, resp.~\,$\wh{K}_{r}(\infty)$\, winds twice around \,$\lambda=0$\, resp.~\,$\lambda=\infty$\,, we have for the meromorphic \,$1$-form \,$\tfrac{\mathrm{d}\lambda}{\lambda}$\,
on~\,$\Sigma$\,
\begin{equation}
\label{eq:jacobitrans:resdlambda/lambda}
\Res_0\left( \frac{\mathrm{d}\lambda}{\lambda} \right) = 2 \qmq{and} \Res_\infty\left( \frac{\mathrm{d}\lambda}{\lambda} \right) = -2 \; .
\end{equation}

\enlargethispage{2em}
The orientation of the paths of integration was chosen such that if \,$\eta$\, is in fact a meromorphic \,$1$-form on \,$\Sigma$\,, then
\begin{equation}
\label{eq:jacobitrans:restheo}
\sum_{P} \Res_P(\eta) + \Res_0(\eta) + \Res_\infty(\eta) = 0 
\end{equation}
holds; here the sum runs over all \,$P\in \Sigma$\, where \,$\eta$\, has a pole.

\begin{lem}
\label{L:jacobitrans:reslemma}
\strut \\
\begin{enumerate}
\item Let \,$h$\, be a holomorphic function on \,$\widehat{V}_\delta$\,. 
\begin{enumerate}
\item If \,$h \in \As_0(\wh{V}_\delta,\ell^2_{-1},0)$\, holds, then we have \,$\Res_0(h\cdot \omega_n)=0$\,. 
\item If \,$h \in \As_\infty(\wh{V}_\delta,\ell^2_{-1},0)$\, holds, then we have \,$\Res_\infty(h\cdot \omega_n)=0$\,. 
\end{enumerate}
\item \begin{enumerate}
\item There exists a sequence \,$(a_n)_{n\in \Z}$\, (dependent only on the spectral curve \,$\Sigma$\,) with \,$(a_n)_{n>0} \in O(n^{-1})$\, and \,$(a_n)_{n<0} \in \ell^2_{-1}(n)$\, such that for every 
\,$h\in \As_0(\wh{V}_\delta,\ell^\infty_{-1},0)$\, for which \,$(h\cdot \sqrt{\lambda})(0) := \lim_{\lambda\to 0,\lambda\in V_\delta} \big( h(\lambda)\cdot \sqrt{\lambda})$\, exists in \,$\C$\, we have
\begin{equation}
\label{eq:jacobitrans:reslemma:Res0}
\Res_0(h\cdot \omega_n) = \begin{cases} \left( h\cdot \sqrt{\lambda} \right)(0) \cdot a_n & \text{for \,$n>0$\,} \\
\left( h\cdot \sqrt{\lambda} \right)(0) \cdot \left( -4in + a_n \right) & \text{for \,$n<0$\,} \end{cases} \; . 
\end{equation}
\item
There exists a sequence \,$(a_n)_{n\in \Z}$\, (dependent only on the spectral curve \,$\Sigma$\,) with \,$(a_n)_{n>0} \in \ell^2_{-1}(n)$\, and \,$(a_n)_{n<0} \in O(n^{-1})$\, such that for every 
\,$h\in \As_\infty(\wh{V}_\delta,\ell^\infty_{-1},0)$\, for which \,$(\tfrac{h}{\sqrt{\lambda}})(\infty) := \lim_{\lambda\to \infty,\lambda\in V_\delta} \big( \tfrac{h(\lambda)}{\sqrt{\lambda}}\big)$\, exists in \,$\C$\, we have
\begin{equation}
\label{eq:jacobitrans:reslemma:Resinf}
\Res_\infty(h\cdot \omega_n) = \begin{cases} \left( \tfrac{h}{\sqrt{\lambda}} \right)(\infty) \cdot \left( -4in + a_n \right) & \text{for \,$n>0$\,} \\
\left( \tfrac{h}{\sqrt{\lambda}} \right)(\infty) \cdot a_n & \text{for \,$n<0$\,} \end{cases} \; .
\end{equation}
\end{enumerate}
\end{enumerate}
\end{lem}

\begin{rem}
In Equations~\eqref{eq:jacobitrans:reslemma:Res0} and \eqref{eq:jacobitrans:reslemma:Resinf}, the factors \,$-4in+a_n$\, resp.~\,$a_n$\, are essentially the value of \,$2\,\omega_n$\, at \,$\lambda=0$\, resp.~\,$\lambda=\infty$\,
in a sense that is made precise by the proof of the lemma below. We write this factor differently for \,$n>0$\, resp.~\,$n<0$\, to obtain a sequence \,$a_n$\, that becomes small for \,$|n|\to\infty$\,. 
The different behavior of \,$a_n$\, for \,$n\to\infty$\, and for \,$n\to -\infty$\, (for example in Lemma~\ref{L:jacobitrans:reslemma}(2)(a) we have \,$a_n=O(n^{-1})$\, for \,$n>0$\, but only
\,$a_n\in \ell^2_{-1}(n)$\, for \,$n<0$\,) occurs because in one case the \,$\lambda_n$\, approach the point where the residue is taken, whereas in the other case the \,$\lambda_n$\, are far away from this point.
\end{rem}

\begin{proof}[Proof of Lemma~\ref{L:jacobitrans:reslemma}.]
At first, we consider in both (1) and (2) the case where the residue of \,$h\cdot \omega_n$\, is calculated at \,$\lambda=0$\,. For this purpose, we suppose that 
\,$h\in \As_0(\wh{V}_\delta,\ell^\infty_{-1},0)$\, holds and that for \,$\wt{h} := h\cdot \sqrt{\lambda} \in \As_0(\wh{V}_\delta,\ell^\infty_0,0)$\,, the limit
\,$\wt{h}(0) := \lim_{\lambda\to 0,\lambda\in V_\delta} \wt{h}(\lambda)$\, exists in \,$\C$\,. 

By Theorem~\ref{T:jacobi:canonical}, we have \,$\omega_n = \tfrac{\Phi_n(\lambda)}{\mu-\mu^{-1}}\,\mathrm{d}\lambda$\,, where \,$\Phi_n$\, is a linear combination of product functions of the kind investigated in 
Proposition~\ref{P:jacobi:asymp-Phi}. By Proposition~\ref{P:jacobi:asymp-f}(1) it therefore follows that we have \,$f_n := \tfrac{\Phi_n(\lambda)}{\mu-\mu^{-1}}|\wh{V}_\delta \in \As_0(\wh{V}_\delta,\ell^\infty_{-1},0)$\,. 
Thus \,$\wt{f}_n := \sqrt{\lambda}\cdot f_n$\, is holomorphic on \,$\wh{V}_\delta$\, with \,$\wt{f}_n \in \As(\wh{V}_\delta,\ell^\infty_0,0)$\,. The asymptotic estimate of Proposition~\ref{P:jacobi:asymp-f}(1) 
(where we have \,$\rho=-1$\,) moreover shows that the limit \,$\wt{f}_n(0) := \lim_{\lambda\to 0,\lambda\in V_\delta} \wt{f}_n(\lambda)$\, exists in \,$\C$\,. 

We now obtain
\begin{align}
\Res_0(h\cdot \omega_n) 
& = \Res_0\left( h\cdot \frac{\Phi_n(\lambda)}{\mu-\mu^{-1}}\,\mathrm{d}\lambda\right)
= \Res_0\left( \frac{1}{\sqrt{\lambda}}\,\wt{h}(\lambda)\cdot \frac{1}{\sqrt{\lambda}}\,\wt{f}_n(\lambda)\,\mathrm{d}\lambda \right) \notag \\
\label{eq:jacobitrans:reslemma:Res0homega}
& = \wt{h}(0)\cdot \wt{f}_n(0)\cdot \Res_0\left( \frac{\mathrm{d}\lambda}{\lambda} \right) \overset{\eqref{eq:jacobitrans:resdlambda/lambda}}{=} 2\cdot \wt{h}(0)\cdot \wt{f}_n(0) \; . 
\end{align}

From this calculation, (1)(a) follows: If we have \,$h\in \As_0(\wh{V}_\delta,\ell^2_{-1},0)$\,, then we have \,$\wt{h} \in \As_0(\wh{V}_\delta,\ell^2_0,0)$\, and therefore \,$\wt{h}(0)=0$\,. 
Thus we obtain \,$\Res_0(h\cdot \omega_n)=0$\, from Equation~\eqref{eq:jacobitrans:reslemma:Res0homega}.

For the proof of (2)(a) we put 
$$ a_n := \begin{cases} 2\,\wt{f}_n(0) & \text{for \,$n>0$\,} \\ 2\wt{f}_n(0)+4in & \text{for \,$n<0$\,} \end{cases} \;; $$
it follows from Equation~\eqref{eq:jacobitrans:reslemma:Res0homega} that with this choice of \,$(a_n)$\,, Equation \eqref{eq:jacobitrans:reslemma:Res0} holds.
Because \,$\wt{f}_n(0)$\, depends only on the spectral curve \,$\Sigma$\,, not on \,$h$\,, also the sequence \,$(a_n)$\, depends only on the spectral curve. 
It remains to show that \,$(a_n)$\, has the claimed asymptotic behavior, i.e.~that \,$(a_n)_{n>0} \in O(n^{-1})$\, and \,$(a_n)_{n<0} \in \ell^2_{-1}(n)$\, holds.

For this purpose we need to determine \,$\wt{f}_n(0)$\, more precisely; it suffices to consider \,$|n|>N$\,, where \,$N\in\N$\, is the constant from Theorem~\ref{T:jacobi:canonical}(2). 
From Theorem~\ref{T:jacobi:canonical}(2)  and Proposition~\ref{P:jacobi:asymp-f}(1) it follows that for \,$n>N$\, we have
$$ \wt{f}_n - (-2i)\cdot \left( s_{n,n}\,\tau_{\xi_n}^{-2}\,\frac{16\pi^2\,n^2}{\lambda_{n,0}-\lambda} + \sum_{\substack{|\ell|\leq N}} s_{n,\ell}\,\tau_{\xi_\ell}^{-2}\,\frac{16\pi^2\,\ell^2}{\lambda_{\ell,0}-\lambda} \right) \;\in\; \As_0(\wh{V}_\delta,\ell^2_{0},0) \;, $$
where \,$s_{n,n}$\, and \,$s_{n,\ell}$\, are the constants from Theorem~\ref{T:jacobi:canonical}(2), and therefore 
\begin{align*}
\wt{f}_n(0)
& = (-2i)\cdot \left( s_{n,n}\,\tau_{\xi_n}^{-2}\,\frac{16\pi^2\,n^2}{\lambda_{n,0}} + \sum_{\substack{|\ell|\leq N}} s_{n,\ell}\,\tau_{\xi_\ell}^{-2}\,\frac{16\pi^2\,\ell^2}{\lambda_{\ell,0}} \right) \;. 
\end{align*}
Because the \,$s_{n,n}$\, and the \,$s_{n,\ell}$\, are of order \,$n^{-1}$\, by Theorem~\ref{T:jacobi:canonical}(2), it follows that we have \,$\wt{f}_n(0)\in O(n^{-1})$\, and 
hence \,$a_n = 2\,\wt{f}_n(0) \in O(n^{-1})$\, for \,$n>0$\,. 

For \,$n<-N$\, we have again by Theorem~\ref{T:jacobi:canonical}(2)  and Proposition~\ref{P:jacobi:asymp-f}(1)
$$ \wt{f}_n(0) = (-2i)\cdot\left( s_{n,n}\,\tau_{\xi_n}^{-2}\,\frac{-1}{\lambda_{n,0}} + \sum_{\substack{|\ell|\leq N}} s_{n,\ell}\,\tau_{\xi_\ell}^{-2}\,\frac{-1}{\lambda_{\ell,0}} \right) \;. $$
By the asymptotic description of \,$s_{n,n}$\, in Equation~\eqref{eq:jacobi:canonical:snn} and again the fact that \,$s_{n,\ell}$\, is of order \,$|n|^{-1}$\, for \,$|\ell|\leq N$\,
we conclude
$$ \wt{f}_n(0) = -\frac{1}{2\pi i} \lambda_{n,0}^{-1/2} + \ell^2_{-1}(n) = -2in + \ell^2_{-1}(n) $$
and therefore \,$a_n = 2\,\wt{f}_n(0)+4in \in \ell^2_{-1}(n)$\, for \,$n<0$\,. 


To also obtain the residue in \,$\lambda=\infty$\, in (1)(b) and (2)(b), we could apply (1)(a) resp.~(2)(a) to the function
\,$\check{h} := h\circ (\lambda\mapsto \lambda^{-1})$\, on the surface \,$\check{\Sigma} := ((\lambda,\mu)\mapsto (\lambda^{-1},\mu))(\Sigma)$\,. But to avoid difficulties in phrasing this transformation, 
we rather carry out a calculation analogous to the one for \,$\lambda=0$\,. We suppose 
\,$h\in \As_\infty(\wh{V}_\delta,\ell^\infty_{-1},0)$\,, and that for \,$\wt{h} := \tfrac{h}{\sqrt{\lambda}} \in \As_\infty(\wh{V}_\delta,\ell^\infty_0,0)$\, the limit 
\,$\wt{h}(\infty) := \lim_{\lambda\to\infty, \lambda\in V_\delta} \wt{h}(\lambda)$\, exists in \,$\C$\,. 
By Proposition~\ref{P:jacobi:asymp-f}(1) we have \,$\tfrac{\Phi_n(\lambda)}{\mu-\mu^{-1}}|\wh{V}_\delta \in \As_\infty(\wh{V}_\delta,\ell^\infty_{3},0)$\,
and thus \,$\wt{f}_n := \lambda^{3/2}\cdot \tfrac{\Phi_n(\lambda)}{\mu-\mu^{-1}}|\wh{V}_\delta \in \As_\infty(\wh{V}_\delta,\ell^\infty_0,0)$\,;
it also follows from Proposition~\ref{P:jacobi:asymp-f}(1) that the limit \,$\wt{f}_n(\infty) := \lim_{\lambda\to\infty,\lambda\in V_\delta} \wt{f}_n(\lambda)$\, exists in \,$\C$\,. 
We thus obtain
\begin{align}
\Res_\infty(h\cdot \omega_n) 
& = \Res_\infty\left( h\cdot \frac{\Phi_n(\lambda)}{\mu-\mu^{-1}}\,\mathrm{d}\lambda\right)
= \Res_\infty\left( \sqrt{\lambda}\,\wt{h}(\lambda)\cdot \frac{1}{\lambda^{3/2}}\,\wt{f}_n(\lambda)\,\mathrm{d}\lambda \right) \notag \\
\label{eq:jacobitrans:reslemma:Resinftyhomega}
& = \wt{h}(\infty)\cdot \wt{f}_n(\infty)\cdot \Res_\infty\left( \frac{\mathrm{d}\lambda}{\lambda} \right) \overset{\eqref{eq:jacobitrans:resdlambda/lambda}}{=} -2\cdot \wt{h}(\infty)\cdot \wt{f}_n(\infty) \; . 
\end{align}
Similarly as before, (1)(b) immediately follows from Equation~\eqref{eq:jacobitrans:reslemma:Resinftyhomega}, and we calculate
$$ \wt{f}_n(\infty)
= \begin{cases} 2in + \ell^2_{-1}(n) & \text{for \,$n>0$\,} \\
O(n^{-1}) & \text{for \,$n<0$\,}
\end{cases} \;, $$
whence (2)(b) also follows.
\end{proof}

The following proposition is crucial. It shows in particular that the translations of the potential correspond to linear motions in the Jacobi variety. 

\begin{prop}
\label{P:jacobitrans:dgl-V2}
Let \,$M(\lambda)=\left( \begin{smallmatrix} a(\lambda) & b(\lambda) \\ c(\lambda) & d(\lambda) \end{smallmatrix} \right)$\, be the monodromy associated to the spectral divisor \,$D$\,,
and let \,$\alpha_{21}^x$\, and \,$\alpha_{21}^y$\, be as in Proposition~\ref{P:jacobitrans:dgl-V1}.
For \,$n \in \Z$\,, we then consider the meromorphic \,$1$-forms on \,$\Sigma$\,
$$ \chi_n^x := \alpha_{21}^x \cdot \frac{\mu-d}{c}\cdot \omega_n \qmq{and} \chi_n^y := \alpha_{21}^y \cdot \frac{\mu-d}{c}\cdot \omega_n \; . $$
Then the residues \,$\Res_0(\chi_n)$\, and \,$\Res_\infty(\chi_n)$\, depend only on the spectral curve \,$\Sigma$\,, not on the spectral divisor \,$D$\, under consideration.
Moreover, there exist sequences \,$(a_n^x)_{n\in \Z}, (a_n^y)_{n\in \Z} \in \ell^2_{-1,-1}(n)$\, (depending only on the spectral curve \,$\Sigma$\,), so that 
the Jacobi coordinates \,$\vi_n$\, (\,$n\in \Z$\,) satisfy the following differential equations
under translation of the potential \,$u$\, in the direction of \,$x$\, resp.~\,$y$\,:
\begin{align*}
\frac{\partial \vi_n}{\partial x} & = \Res_0(\chi_n^x) + \Res_\infty(\chi_n^x) = n + a_n^x \\
\frac{\partial \vi_n}{\partial y} & = \Res_0(\chi_n^y) + \Res_\infty(\chi_n^y) = -i|n| + a_n^y\; . 
\end{align*}
\end{prop}

\begin{proof}
Because the set of tame divisors on \,$\Sigma'$\, is dense in \,$\Div(\Sigma')$\, and the right hand side of the differential equations claimed does not depend on the spectral divisor \,$D$\,, it suffices to consider the case
where the spectral divisor \,$D$\, is tame. For the same reason, we may further suppose without loss of generality that no point in the support of \,$D$\, is a branch point of \,$\Sigma'$\,. 

We let \,$\alpha_{21}$\, be one of the functions \,$\alpha_{21}^x$\,, \,$\alpha_{21}^y$\,, and correspondingly let \,$\chi_n$\, be one of \,$\chi_n^x$\,, \,$\chi_n^y$\,. 
\,$\chi_n$\, is a meromorphic \,$1$-form on \,$\Sigma$\,; because we supposed that \,$D$\, is tame, \,$\chi_n$\,  
has simple poles in the points \,$(\lambda_k,\mu_k)$\,, \,$k\in \Z$\, comprising the support of \,$D$\,, and no other poles. We compute the residue of \,$\chi_n$\, in \,$(\lambda_k,\mu_k)$\,: We have
$$ \chi_n = \alpha_{21}\cdot \frac{\mu-d}{c}\cdot \omega_n = \frac{\alpha_{21}(\lambda)\cdot (\mu-d(\lambda)) \cdot \Phi_n(\lambda)}{\mu-\mu^{-1}}\cdot \frac{\mathrm{d}\lambda}{c(\lambda)} \;. $$
Because \,$(\lambda_k,\mu_k)$\, is not a branch point of \,$\Sigma'$\,, the function 
\,$\tfrac{\alpha_{21}(\lambda)\cdot (\mu-d(\lambda)) \cdot \Phi_n(\lambda)}{\mu-\mu^{-1}}$\, is holomorphic at \,$(\lambda_k,\mu_k)$\,. Moreover, \,$c(\lambda)$\, has a zero of order \,$1$\, and \,$d(\lambda_k)=\mu_k^{-1}$\, holds, 
and thus we obtain
\begin{align*}
\Res_{(\lambda_k,\mu_k)} (\chi_n) & = \frac{\alpha_{21}(\lambda_k) \cdot (\mu_k-d(\lambda_k)) \cdot \Phi_n(\lambda_k)}{\mu_k-\mu_k^{-1}} \cdot \frac{1}{c'(\lambda_k)} = \alpha_{21}(\lambda_k) \cdot \Phi_n(\lambda_k) \cdot \frac{1}{c'(\lambda_k)}  \; . 
\end{align*}
By Proposition~\ref{P:jacobitrans:dgl-V1} and Equation~\eqref{eq:jacobitrans:restheo} we therefore obtain
$$ \frac{\partial \vi_n}{\partial x} = - \sum_{k\in \Z} \Res_{(\lambda_k,\mu_k)}(\chi_n) = \Res_0(\chi_n) + \Res_\infty(\chi_n) \; . $$

We now proceed to calculate the residues of \,$\chi_n$\, in \,$\lambda=0$\, and in \,$\lambda=\infty$\,. We will see that these residues do not depend on the divisor \,$D$\,. For this, 
the trace formula of Theorem~\ref{T:asympfinal:monodromy}(3) turns out to be crucial: This formula shows 
that the constant \,$\tau=e^{-u(0)/2}$\, occurring in the component \,$\alpha_{21}$\, of the connection form \,$\alpha$\, associated to the potential \,$u$\, equals the constant 
\,$\tau=\left( \prod_{k\in \Z} \frac{\lambda_{k,0}}{\lambda_k} \right)^{1/2}$\, occurring in the description of the asymptotic behavior of \,$c$\, given in Proposition~\ref{P:interpolate:lambda}.

By Corollary~\ref{C:asympfinal:monodromy2}(6) we have
$$ \frac{\mu-d(\lambda)}{c(\lambda)}-\frac{i}{\tau\,\sqrt{\lambda}} \in \As_\infty(\wh{V}_\delta,\ell^2_1,0) \qmq{and} \frac{\mu-d(\lambda)}{c(\lambda)}-\frac{i}{\tau^{-1}\,\sqrt{\lambda}} \in \As_0(\wh{V}_\delta,\ell^2_{-1},0) \;, $$
and therefore
\begin{align*}
\alpha_{21}(\lambda) \cdot \left( \frac{\mu-d(\lambda)}{c(\lambda)}-\frac{i}{\tau\,\sqrt{\lambda}}\right) & \in \As_\infty(\wh{V}_\delta,\ell^2_{-1},0) \\
\qmq{and}
\alpha_{21}(\lambda) \cdot \left( \frac{\mu-d(\lambda)}{c(\lambda)}-\frac{i}{\tau^{-1}\,\sqrt{\lambda}}\right) & \in \As_0(\wh{V}_\delta,\ell^2_{-1},0) \; .
\end{align*}
We therefore obtain by Lemma~\ref{L:jacobitrans:reslemma}(1),(2)
\begin{align*}
\Res_0(\chi_n) & = \Res_0\left(\alpha_{21}\cdot \frac{\mu-d}{c}\cdot \omega_n \right) = \Res_0\left(\alpha_{21}\cdot \frac{i}{\tau^{-1}\,\sqrt{\lambda}} \cdot \omega_n \right) \\
& = \left( \frac{\alpha_{21}\cdot i}{\tau^{-1}} \right)(0) \cdot 
\begin{cases} O(n^{-1}) & \text{for \,$n>0$\,} \\ -4in + \ell^2_{-1}(n) & \text{for \,$n<0$\,} \end{cases} 
\end{align*}
and
\begin{align*}
\Res_\infty(\chi_n) & = \Res_\infty\left(\alpha_{21}\cdot \frac{\mu-d}{c}\cdot \omega_n \right) = \Res_\infty\left(\alpha_{21}\cdot \frac{i}{\tau\,\sqrt{\lambda}} \cdot \omega_n \right) \\
& = \left( \frac{\alpha_{21}\cdot i}{\tau\cdot \lambda} \right)(\infty) \cdot 
\begin{cases} -4in + \ell^2_{-1}(n) & \text{for \,$n>0$\,} \\ O(n^{-1}) & \text{for \,$n<0$\,} \end{cases} \; ,
\end{align*}
where the sequences in \,$O(n^{-1})$\, and \,$\ell^2_{-1}(n)$ occurring here depend only on the spectral curve \,$\Sigma$\,, not on the divisor \,$D$\,.

When the translation in \,$x$\, direction is concerned, we have \,$\alpha_{21}^x=\tfrac14(\lambda\tau + \tau^{-1})$\, and therefore
$$ \left( \frac{\alpha_{21}^x\cdot i}{\tau^{-1}} \right)(0) = \left( \frac{\alpha_{21}^x\cdot i}{\tau\cdot \lambda} \right)(\infty) = \frac{i}{4} \;, $$
whence
$$ \frac{\partial \vi_n}{\partial x} = \Res_0(\chi_n^x) + \Res_\infty(\chi_n^x) = \frac{i}{4}\cdot (-4in + \ell^2_{-1,-1}(n)) = n + \ell^2_{-1,-1}(n) $$
follows, again with an \,$\ell^2_{-1,-1}(n)$-sequence that depends only on \,$\Sigma$\,. 
On the other hand, for the translation in \,$y$\, direction, we have \,$\alpha_{21}^y = \tfrac{i}{4}(-\lambda \tau + \tau^{-1})$\, and therefore
$$ \left( \frac{\alpha_{21}^y\cdot i}{\tau^{-1}} \right)(0) = -\frac{1}{4} \qmq{and} \left( \frac{\alpha_{21}^y\cdot i}{\tau\cdot \lambda} \right)(\infty) = \frac{1}{4} \;, $$
whence
\begin{align*}
\frac{\partial \vi_n}{\partial y} 
& = \Res_0(\chi_n^y) + \Res_\infty(\chi_n^y) = \begin{cases} \tfrac14\cdot (-4in + \ell^2_{-1}(n)) & \text{for \,$n>0$\,} \\ -\tfrac14\cdot (-4in + \ell^2_{-1}(n)) & \text{for \,$n<0$\,} \end{cases} \\
& = -i|n| + \ell^2_{-1,-1}(n)
\end{align*}
follows, where the \,$\ell^2_{-1,-1}(n)$-sequence depends only on \,$\Sigma$\,. 
\end{proof}

For the translation in \,$x$-direction, we expect a better result still. Our potentials are periodic in that direction, so we expect that \,$\tfrac{\partial \vi_n}{\partial x}$\, is an \emph{exact}, not only
an asymptotically approximate multiple of a lattice vector of the Jacobi variety of \,$\Sigma$\,. (For the translation in \,$y$-direction we cannot expect a similar result, because \,$u$\, is not periodic
in the direction of \,$y$\,.)
Because the coordinates \,$(\lambda_n,\mu_n)$\, of the spectral divisor corresponding to some \,$(u,u_y)\in \Pot$\, are asymptotically close to the \,$n$-th Fourier coefficient of \,$u_z$\, resp.~\,$u_{\overline{z}}$\,,
we more specifically expect that \,$\tfrac{\partial \vi_n}{\partial x}=n$\, holds at least for \,$|n|$\, large.

The periodicity of the potential \,$u$\, corresponds to the existence of the global function \,$\mu$\, on the spectral curve \,$\Sigma$\,,
and we base our proof of the expected statement on the existence of a global logarithm \,$\ln(\mu)$\, of that function on a surface \,$\wt{\Sigma}$\, with boundary obtained from \,$\Sigma$\, by cutting along certain curves.
The following lemma serves to establish the topological prerequisites for the cutting process.

\begin{lem}
\label{L:jacobi:ordering-the-vkap}
The branch points resp.~singularities \,$\vkap_{k,\nu}$\, of the spectral curve \,$\Sigma$\, can be numbered in such a way that besides the previous requirements (\,$\vkap_{k,\nu}\in \wh{U}_{k,\delta}$\, for \,$|k|$\, large
and \,$\vkap$\, is a double point of \,$\Sigma$\, if and only if there exists \,$k\in S$\, with \,$\vkap_{k,1}=\vkap_{k,2}=\vkap$\,) the following holds:

For every \,$k\in \Z$\, there exists a curve \,$\vartheta_k$\, in the \,$\lambda$-plane, which for \,$k\not\in S$\, connects \,$\vkap_{k,1}$\, to \,$\vkap_{k,2}$\, whereas for \,$k\in S$\, is a small circle
around \,$\vkap_{k,*}$\, that does not include any other branch points of \,$\Sigma$\,, such that none of the \,$\vartheta_k$\, intersect, and such that for every \,$k\in \Z$\, we have
$$ \int_{\vartheta_k} \frac{\mathrm{d}\mu}{\mu} = 0 \;. $$
\end{lem}

\newpage

\begin{proof}
We consider the potential \,$(u,u_y)\in \Pot$\, corresponding to the spectral data \,$(\Sigma,D^o)$\,, and the family \,$u_t := (t\cdot u,t\cdot u_y) \in \Pot$\, with \,$t\in [0,1]$\,
describing the deformation of the vacuum (\,$u=0$\,) into that potential. For \,$t\in [0,1]$\, we let \,$\Sigma(t)$\, be the spectral curve corresponding to \,$u_t$\,. The set \,$\Menge{u_t}{t\in [0,1]}$\, is relatively 
compact in \,$\Pot$\,, therefore the asymptotic estimates apply uniformly to this set of potentials. Hence there exist continuous functions \,$\vkap_{k,\nu}(t): [0,1]\to \C^*$\, so that 
\,$\vkap_{k,\nu}(0)=\lambda_{k,0}$\, holds (these points are the double points of the spectral curve \,$\Sigma(0)=\Sigma_0$\, of the vacuum) and so that for every \,$t\in [0,1]$\,, the \,$\vkap_{k,\nu}(t)$\, are
the branch points resp.~singularities of \,$\Sigma(t)$\,; moreover there exists \,$N\in \N$\, so that \,$\vkap_{k,\nu}(t) \in U_{k,\delta}$\, holds for all \,$k\in \Z$\, with \,$|k|>N$\, and for all \,$t\in [0,1]$\,.
Note that it is possible to choose these functions as continuous even if some of the \,$\vkap_{k,\nu}(t)$\, coincide for some \,$t\in[0,1]$\,.

Further there exist for every \,$t\in [0,1]$\, continuous curves \,$\vartheta_{k,t}:[0,1]\to \C^*$\, which connect \,$\vkap_{k,1}(t)$\, to \,$\vkap_{k,2}(t)$\,. They can be chosen such that they depend
continuously also on \,$t\in [0,1]$\,, that they run entirely within \,$U_{k,\delta}$\, for \,$|k|>N$\, and that they are constant if \,$\vkap_{k,1}(t)=\vkap_{k,2}(t)$\, holds. Moreover after possibly
modifying the \,$\vkap_{k,\nu}(t)$\, past the times \,$t$\, where two \,$\vkap_{k,\nu}(t)$\, intersect, it is possible to choose the \,$\vartheta_{k,t}$\, so that for every \,$t\in [0,1]$\, no two \,$\vartheta_{k,t}$\, intersect. 

Note that for every \,$t\in [0,1]$\, the lift onto \,$\Sigma(t)$\, of \,$\vartheta_{k,t}$\, traversed once in the usual direction and then once in the opposite direction is a closed curve through \,$\vkap_{k,1}(t)$\, 
and \,$\vkap_{k,2}(t)$\, that meets no other branch points of \,$\Sigma(t)$\,, and therefore \,$\int_{\vartheta_{k,t}} \tfrac{\mathrm{d}\mu}{\mu}$\, is in any event an integer multiple of \,$i\pi$\, (here \,$\mu$\,
denotes the corresponding parameter of the hyperelliptic curve \,$\Sigma(t)$\,). Because this integral depends continuously on \,$t$\,, and we have \,$\int_{\vartheta_{k,0}} \tfrac{\mathrm{d}\mu}{\mu}=0$\,
because \,$\vartheta_{k,0}$\, is constant, it follows that 
$$ \int_{\vartheta_{k,t}} \frac{\mathrm{d}\mu}{\mu}=0 \qmq{holds for all \,$t\in [0,1]$\, and \,$k\in \Z$\,.} $$

Because of the hypothesis that \,$\Sigma=\Sigma(1)$\, does not have any singularities other than ordinary double points, no more than two of the \,$\vkap_{k,\nu}(t=1)$\, can coincide. However it is possible
for \,$\vkap_{k_1,\nu_1}(t=1)=\vkap_{k_2,\nu_2}(t=1)$\, to hold, where \,$k_1\neq k_2$\, and \,$\nu_1,\nu_2\in \{1,2\}$\,. In this case, we put
\begin{align*}
\vkap_{k_1,1} & := \vkap_{k_1,\nu_1}(t=1) & \vkap_{k_2,1} & := \vkap_{k_1,3-\nu_1}(t=1) \\
\vkap_{k_1,2} & := \vkap_{k_2,\nu_2}(t=1) & \vkap_{k_2,2} & := \vkap_{k_2,3-\nu_2}(t=1) \;,
\end{align*}
and we let \,$\vartheta_{k_1}$\, be a small circle around \,$\vkap_{k_1,1}=\vkap_{k_1,2}$\, that does not encircle any other branch points of \,$\Sigma$\,. Moreover if \,$\vkap_{k_2,1}\neq \vkap_{k_2,2}$\, we let 
\,$\vartheta_{k_2}$\, be a curve close and homologically equivalent to the concatenation of \,$\vartheta_{k_1,t=1}$\, and \,$\vartheta_{k_2,t=1}$\, in the appropriate direction to connect \,$\vkap_{k_2,1}$\, with \,$\vkap_{k_2,2}$\,,
but which avoids the circle \,$\vartheta_{k_1}$\, around \,$\vkap_{k_1,1}=\vkap_{k_1,2}$\,; if \,$\vkap_{k_2,1}= \vkap_{k_2,2}$\, we let \,$\vartheta_{k_2}$\, again 
be a small circle around \,$\vkap_{k_2,1}=\vkap_{k_2,2}$\, that does not encircle any other branch points of \,$\Sigma$\,.

In all other cases we let \,$\vkap_{k,\nu} := \vkap_{k,\nu}(t=1)$\,, and if \,$\vkap_{k,1}\neq \vkap_{k,2}$\, holds we let \,$\vartheta_k := \vartheta_{k,t=1}$\,, whereas for \,$\vkap_{k,1}=\vkap_{k,2}$\, we let 
\,$\vartheta_k$\, be a small circle around \,$\vkap_{k,1}=\vkap_{k,2}$\, that does not encircle any other branch points of \,$\Sigma$\,.
\end{proof}

\begin{prop}
\label{P:jacobi:lnmu}
Suppose that the branch points \,$\vkap_{k,\nu}$\, of \,$\Sigma$\, are numbered as in Lemma~\ref{L:jacobi:ordering-the-vkap}, and that the cycles \,$A_k$\, from the canonical basis \,$(A_k,B_k)$\, of the homology
of \,$\Sigma$\, are chosen as the lift onto \,$\Sigma$\, of two traversals of the curves \,$\vartheta_k$\, from Lemma~\ref{L:jacobi:ordering-the-vkap}.

Then we let \,$\wt{\Sigma}$\, be the surface with boundary that is obtained from \,$\Sigma$\, by cutting along all the cycles \,$A_n$\,.
\,$\wt{\Sigma}$\, contains as boundaries two copies of each \,$A_n$\,, which we denote by \,$A_n^+$\, and \,$A_n^-$\, (as cycles with the same orientation as \,$A_n$\,), 
where \,$B_n$\, runs from \,$A_n^+$\, to \,$A_n^-$\,. We view \,$\wh{V}_\delta$\, as a subset of \,$\wt{\Sigma}$\,. 

On \,$\wt{\Sigma}$\,, \,$\ln(\mu)$\, exists as a global, holomorphic function. For \,$|n|$\, large, and corresponding points \,$\mu^+$\, and \,$\mu^-$\, on \,$A_n^+$\, resp.~on \,$A_n^-$\,, we have
\,$\ln(\mu^+)-\ln(\mu^-)=-2\pi i n$\,. 
Moreover \,$\bigr(\ln(\mu)-i\,\zeta(\lambda)\bigr)|\wh{V}_\delta \in \As(\wh{V}_\delta,\ell^2_{0,0},0)$\, holds. 
\end{prop}

\begin{proof}
By Lemma~\ref{L:jacobi:ordering-the-vkap} we have for any \,$k\in \Z$\,
\begin{equation}
\label{eq:jacobi:lnmu:Ak}
\int_{A_k} \mathrm{d}(\ln(\mu)) = 2\int_{\vartheta_k} \frac{\mathrm{d}\mu}{\mu}=0 \; .
\end{equation}

For any cycle \,$Z$\, on \,$\Sigma$\,,
$$ \int_Z \mathrm{d}(\ln(\mu)) = \int_Z \frac{\mathrm{d}\mu}{\mu} $$
is \,$2\pi i$\, times the winding number of \,$\mu \circ Z$\, around \,$0$\,. 
If the cycle \,$Z$\, passes from an excluded domain \,$\wh{U}_{k,\delta}$\, to \,$\wh{U}_{k\pm 1,\delta}$\,, \,$\mu$\, changes from values near \,$\pm 1$\, to values near \,$\mp 1$\,, and therefore we have
for the cycle \,$Z=B_k$\, with \,$|k|$\, large (this cycle connects \,$\wh{U}_{k,\delta}$\, to \,$\wh{U}_{0,\delta}$\,)
\begin{equation}
\label{eq:jacobi:lnmu:Bk}
\int_{B_k} \mathrm{d}(\ln(\mu)) = 2\pi i k \;. 
\end{equation}

Because the homology group of \,$\wt{\Sigma}$\, is generated by the \,$A_k^\pm$\,, it follows from Equation~\eqref{eq:jacobi:lnmu:Ak} 
that \,$\ln(\mu)$\, is a global holomorphic function on \,$\wt{\Sigma}$\,.
Moreover, it follows from Equation~\eqref{eq:jacobi:lnmu:Bk} that the values of \,$\ln(\mu)$\, at corresponding points of \,$A_n^+$\, and \,$A_n^-$\, differ by \,$-2\pi i n$\,. 

Finally, concerning the asymptotic behavior of \,$\ln(\mu)$\, near \,$\lambda=0$\,, \,$\lambda=\infty$\,, we note that 
by  Corollary~\ref{C:asympfinal:monodromy2}(4), we have \,$(\mu-\mu_0)|\wh{V}_\delta \in \As(\wh{V}_\delta,\ell^2_{0,0},1)$\, 
and therefore \,$\left. \left( \tfrac{\mu-\mu_0}{\mu_0} \right) \right|\wh{V}_\delta \in \As(\wh{V}_\delta,\ell^2_{0,0},0)$\, because on \,$\wh{V}_\delta$\,,
\,$\mu_0$\, is comparable to \,$w(\lambda)$\,. Herefrom we obtain by the equality
$$ \ln(\mu)-\ln(\mu_0) = \ln\left( \frac{\mu}{\mu_0} \right) = \ln\left( 1 + \frac{\mu-\mu_0}{\mu_0} \right) $$
and the Taylor expansion \,$\ln(1+z)=z+O(z^2)$\, of the logarithm near \,$z=1$\, that 
\,$\ln(\mu)-\ln(\mu_0) \in \As(\wh{V}_\delta,\ell^2_{0,0},0)$\, holds. Because of \,$\ln(\mu_0)=i\zeta(\lambda)$\, (compare
Equation~\eqref{eq:vacuum:Sigma0}), the claimed asymptotic behavior follows.
\end{proof}

\begin{prop}
\label{P:jacobitrans:dgl-V3}
Under translation of the potential \,$u$\, in the direction of \,$x$\,, the Jacobi coordinates \,$\vi_n$\, with \,$|n|$\, large satisfy the differential equation
$$ \frac{\partial \vi_n}{\partial x} = n \;. $$
\end{prop}

\begin{proof}
For \,$n\in \Z$\, we consider the holomorphic 1-form \,$\eta_n := \ln(\mu) \cdot \omega_n$\, on the surface with boundary \,$\widetilde{\Sigma} \supset \wh{V}_\delta$\, from Proposition~\ref{P:jacobi:lnmu}.
We claim that
\begin{equation}
\label{eq:jacobitrans:dgl-V3:resetachi}
\Res_0(\eta_n) = \Res_0(\chi^x_n) \qmq{and} \Res_\infty(\eta_n) = \Res_\infty(\chi^x_n)
\end{equation}
holds. 
Indeed, by the calculations in the proof of Proposition~\ref{P:jacobitrans:dgl-V2} we have
\begin{align*}
\left. \left( \alpha_{21}^x\,\frac{\mu-d}{c} - \frac{i\,\sqrt{\lambda}}{4}\right)\right|\wh{V}_\delta & \in \As_\infty(\wh{V}_\delta,\ell^2_{-1},0) \\
\qmq{and}
\left. \left( \alpha_{21}^x\,\frac{\mu-d}{c} - \frac{i}{4\,\sqrt{\lambda}} \right)\right|\wh{V}_\delta & \in \As_0(\wh{V}_\delta,\ell^2_{-1},0) \;,
\end{align*}
and therefore
$$ \left. \left( \alpha_{21}^x\,\frac{\mu-d}{c} - i\,\zeta(\lambda) \right)\right|\wh{V}_\delta \in \As(\wh{V}_\delta,\ell^2_{-1,-1},0) $$
holds. By Proposition~\ref{P:jacobi:lnmu} we also have
$$ \left. \left( \ln(\mu) - i\,\zeta(\lambda) \right)\right|\wh{V}_\delta \in \As(\wh{V}_\delta,\ell^2_{0,0},0) \;, $$
and thus we obtain 
$$ \left. \left( \alpha_{21}^x\,\frac{\mu-d}{c} - \ln(\mu) \right)\right|\wh{V}_\delta \in \As(\wh{V}_\delta,\ell^2_{-1,-1},0) \; . $$
It therefore follows from Lemma~\ref{L:jacobitrans:reslemma}(1) that Equations~\eqref{eq:jacobitrans:dgl-V3:resetachi} hold.

Due to the topology of \,$\wt{\Sigma}$\,, we have
\begin{equation}
\label{eq:jacobitrans:dgl-V3:restheo-eta}
\sum_{k\in \Z} \frac{1}{2\pi i}\left( \int_{A_k^+} \eta_n - \int_{A_k^-} \eta_n\right)  + \Res_0(\eta_n) + \Res_\infty(\eta_n) = 0 \; .
\end{equation}
For \,$|k|$\, large, the function values of \,$\ln(\mu)$\, at corresponding points of \,$A_k^+$\,
and \,$A_k^-$\, differ by \,$-2\pi i k$\, (see Proposition~\ref{P:jacobi:lnmu}), and thus 
\begin{equation*}
\int_{A_k^+} \eta_n - \int_{A_k^-} \eta_n = -2\pi i k \cdot \int_{A_k} \omega_n = -2\pi i k \cdot \delta_{k,n} 
\end{equation*}
holds.
Therefore it follows from Equation~\eqref{eq:jacobitrans:dgl-V3:restheo-eta} for \,$|n|$\, large that 
\begin{equation}
\label{eq:jacobitrans:dgl-V3:reseta}
\Res_0(\eta_n) + \Res_\infty(\eta_n) = n
\end{equation}
holds.

From Proposition~\ref{P:jacobitrans:dgl-V2} we now obtain
$$ \frac{\partial \vi_n}{\partial x} 
= \Res_0(\chi_n^x) + \Res_\infty(\chi_n^x) \overset{\eqref{eq:jacobitrans:dgl-V3:resetachi}}{=} \Res_0(\eta_n) + \Res_\infty(\eta_n) \overset{\eqref{eq:jacobitrans:dgl-V3:reseta}}{=} n \; . $$
\end{proof}

\enlargethispage{2em}

\begin{proof}[Proof of Theorem~\ref{T:jacobitrans:dgl}.]
In view of Propositions~\ref{P:jacobitrans:dgl-V2} and \ref{P:jacobitrans:dgl-V3} it only remains to show that
for \,$n\in \Z$\, there exists a cycle \,$Z_n$\, on \,$\Sigma'$\, so that \,$\tfrac{\partial \vi_n}{\partial x}=\int_{Z_n}\omega_n$\,
holds. To show this we choose a tame divisor \,$D$\, on \,$\Sigma$\,;
by Theorem~\ref{T:diffeo:diffeo} there exists a periodic potential \,$(u,u_y)\in \Pot$\, with spectral data \,$(\Sigma,D)$\,. Because this potential is periodic in \,$x$-direction, 
the flow corresponding to \,$x$-translation in the Jacobi variety is also periodic, 
and therefore the Jacobi coordinate \,$\vi_n$\, changes under the flow of \,$x$-translation along a period of \,$(u,u_y)$\, by a lattice vector of the Jacobi variety.
Thus there exists a cycle \,$Z_n$\, on \,$\Sigma'$\, so that \,$\tfrac{\partial \vi_n}{\partial x}=\int_{Z_n}\omega_n$\, holds.
\end{proof}


\begin{rem}
\label{R:jacobitrans:krichever}
The preceding proof of Theorem~\ref{T:jacobitrans:dgl} bears a relationship to a construction principle for linear flows on the Picard variety of a compact Riemann surface \,$X$\, described by 
\textsc{Krichever}, see \cite{Krichever:1977}. 

We now describe this construction in general terms for a compact Riemann surface \,$X$\,. We mark points \,$x_1,\dotsc, x_n\in X$\, on \,$X$\, and fix local coordinates \,$z_1,\dotsc,z_n$\, of \,$X$\, around these points
with \,$z_k(x_k)=0$\,. Let \,$H$\, be the linear space of meromorphic function germs near \,$0\in \C$\, (i.e.~\,$H$\, is the space of Laurent series in one complex variable with positive convergence radius). For every
\,$(h_1,\dotsc,h_n)\in H^n$\, there exist neighborhoods \,$U_k$\, of \,$x_k$\, in \,$X$\, such that \,$z_k^* h_k$\, is a holomorphic function on \,$U_k \setminus \{x_k\}$\,. Together with
\,$U_0 := X \setminus \{x_1,\dotsc,x_n\}$\,, \,$\mathcal{U} := \{U_0,U_1,\dotsc,U_n\}$\, is an open covering of \,$X$\,, and \,$(z_k^*h_k)_{k=1,\dotsc,n}$\, defines a cocycle with respect to the covering \,$\mathcal{U}$\,
and thus induces an element in the cohomology group \,$H^1(X,\mathcal{O})$\, (where \,$\mathcal{O}$\, denotes the sheaf of holomorphic functions on \,$X$\,). Thereby we obtain a map
\,$\psi: H^n \to H^1(X,\mathcal{O})$\,.

The short exact sequence of sheaves
$$ 0 \longrightarrow \Z \longrightarrow \calO \longrightarrow \calO^* \longrightarrow 0 $$
(where \,$\calO^*$\, denotes the subsheaf of \,$\calO$\, of invertible function germs, and the map \,$\calO \longrightarrow \calO^*$\, is \,$f \mapsto \exp(f)$\,)
induces a long exact sequence of cohomology groups, see \cite{Forster:1981}, Theorem~15.12, p.~123
\begin{gather*}
0 \longrightarrow H^0(X,\Z) \longrightarrow H^0(X,\calO) \longrightarrow H^0(X,\calO^*) \longrightarrow H^1(X,\Z) \longrightarrow \\
\longrightarrow H^1(X,\calO) \longrightarrow H^1(X,\calO^*) \longrightarrow H^2(X,\Z) \longrightarrow H^2(X,\calO) \longrightarrow \dotsc \; . 
\end{gather*}
The sequence \,$0 \longrightarrow H^0(X,\Z) \longrightarrow H^0(X,\calO) \longrightarrow H^0(X,\calO^*) \longrightarrow 0$\, splits off, and moreover we have \,$H^2(X,\calO)=0$\,. Thus we obtain the exactness of the sequence
$$ 0 \longrightarrow H^1(X,\Z) \longrightarrow H^1(X,\calO) \longrightarrow H^1(X,\calO^*) \longrightarrow H^2(X,\Z) \longrightarrow 0 \; . $$
Here we note that the group \,$H^1(X,\calO^*)$\, is the Picard variety of \,$X$\,, which is isomorphic to the space of isomorphy classes of line bundles on \,$X$\,, and the map 
\,$H^1(X,\calO^*) \longrightarrow H^2(X,\Z)$\, in the above sequence is the degree map of line bundles. Therefore the Jacobi variety \,$H^1(X,\calO)$\, of \,$X$\, is the Lie algebra of the Picard variety
\,$H^1(X,\calO^*)$\,. 

For this reason, every \,$h=(h_1,\dotsc,h_n)\in H^n$\, defines a family \,$L_h(t)$\, of elements of \,$H^1(X,\calO^*)$\,, i.e.~of line bundles on \,$X$\,; the cocycle defining \,$L_h(t)$\, is given by
\,$(z_k^*\exp(h_k))_{k=1,\dotsc,n}$\,. This family is a one-parameter subgroup of \,$H^1(X,\calO^*)$\,, i.e.~we have \,$L_h(t+t') = L_h(t) \otimes L_h(t')$\,. Moreover it can be shown that 
the flow \,$L_h(t)$\, thus induced by some \,$h\in H^n$\, is periodic with period \,$1$\, (i.e.~\,$L_h(t+1)=L_h(t)$\, for all \,$t$\,) if and only if the Mittag-Leffler distribution given by \,$h$\, is 
solvable by means of a multi-valued, meromorphic function on \,$X$\, whose function values at a point differ by elements of \,$2\pi i \Z$\,. 

The Krichever construction applies to our situation of the present section in the following way: We carry out the analogous construction on the spectral curve \,$\Sigma$\,
or on the surface \,$\wt{\Sigma}$\, constructed in Proposition~\ref{P:jacobi:lnmu}, which unlike the previous surface \,$X$\, 
are non-compact and can have singularities. We mark the two points \,$0$\, and \,$\infty$\,, and consider the local coordinates \,$\sqrt{\lambda}$\, and \,$\tfrac{1}{\sqrt{\lambda}}$\, around these points.
Let \,$h_x=(h_{x,0},h_{x,\infty})$\, resp.~\,$h_y=(h_{y,0},h_{y,\infty})$\, be the pair of Laurent series near \,$0\in \C$\, which each have a pole of order \,$1$\, and whose principal part is characterized by
$$ \Res_0(h_{x,0}) = \Res_0\left( \alpha_{21}^x\cdot \frac{\mu-d}{c} \right) \qmq{and}\Res_0(h_{x,\infty}) = \Res_\infty\left( \alpha_{21}^x\cdot \frac{\mu-d}{c} \right) $$
resp.~the analogous equations for \,$h_y$\,; here Proposition~\ref{P:jacobitrans:dgl-V2} shows that these residues do not depend on the monodromy \,$M(\lambda)=\left(\begin{smallmatrix} a & b \\c & d \end{smallmatrix} \right)$\,
with spectral curve \,$\Sigma$\, used in the equations. Then the spectral divisor \,$D$\, of this monodromy corresponds to the line bundle \,$\Lambda$\, with the holomorphic section \,$(\tfrac{\mu-d}{c},1)$\,;
when we translate the divisor in \,$x$-direction resp.~in \,$y$-direction, the translated divisor corresponds to the line bundle \,$\Lambda \otimes L_{h_x}(t)$\, resp.~\,$\Lambda \otimes L_{h_y}(t)$\,.
This is the interpretation of Proposition~\ref{P:jacobitrans:dgl-V2} in the context of the Krichever construction. Concerning the periodicity of the \,$x$-translation, Proposition~\ref{P:jacobitrans:dgl-V3} shows
that the function \,$\ln(\mu)$\, constructed on \,$\wt{\Sigma}$\, in Proposition~\ref{P:jacobi:lnmu} solves the Mittag-Leffler distribution given by \,$h_x$\, as a multi-valued, meromorphic function,
whose function values at a point differ by a multiple of \,$2\pi i$\,. Thus it follows that the 1-parameter group \,$L_{h_x}(t)$\, induced by the \,$x$-translation is periodic, as expected.
\end{rem}

\section{Asymptotics of spectral data for potentials on a horizontal strip}
\label{Se:strip}

As a final result, we study the asymptotic behavior of the spectral data \,$(\Sigma,D)$\, corresponding to a simply periodic solution \,$u: X \to \C$\, of the sinh-Gordon equation defined on an entire 
horizontal strip \,$X \subset \C$\, with positive height. Because such a solution is real analytic on the interior of \,$X$\,, we expect a far better asymptotic for such spectral data than for the
spectral data of Cauchy data potentials \,$(u,u_y)$\, with only the weak requirements \,$u\in W^{1,2}([0,1])$\,, \,$u_y \in L^2([0,1])$\, we have been using throughout most of the paper.
More specifically, we expect both the distance of branch points \,$\vkap_{k,1}-\vkap_{k,2}$\, of the spectral curve \,$\Sigma$\, and the distance of the corresponding spectral divisor points
to the branch points to fall off exponentially for \,$k\to\pm \infty$\,. 

The following theorem shows that our expectation is correct:

\begin{thm}
\label{T:strip:falloff}
Let \,$X \subset \C$\, be a closed, horizontal strip in \,$\C$\, of height \,$2y_0$\, with \,$y_0>0$\,, i.e.~\,$X=\Mengegr{z\in \C}{|\IM(z)|\leq y_0}$\,, and let \,$u:X \to \C$\, be a simply periodic solution of the
sinh-Gordon equation \,$\Delta u + \sinh(u)=0$\,, such that for every \,$y \in [-y_0,y_0]$\, (even on the boundary!) we have \,$u(\,\cdot\, + iy)\in \Pot$\, and \,$\|u(\,\cdot\, + iy)\|_\Pot \leq C_1$\, with \,$C_1>0$\,.
We let \,$\Sigma$\, be the spectral curve corresponding to \,$u$\, (with branch points \,$\vkap_{n,\nu}$\,, and \,$\vkap_{n,*}:=\tfrac12(\vkap_{n,1}+\vkap_{n,2})$\,) and let \,$D := \{(\lambda_n,\mu_n)\}_{n\in \Z}$\,
be the spectral divisor of \,$u$\, along the real axis.

Then there exists a constant \,$C>0$\, and a sequence \,$(s_n)_{n\in \Z} \in \ell^2_{0,0}(n)$\, of real numbers so that 
\begin{align} 
\label{eq:strip:falloff:vkap}
|\vkap_{n,1}-\vkap_{n,2}| & \leq C\,e^{-2\pi\,(1-s_n)\,|n|\,y_0}  \;,\\
\label{eq:strip:falloff:lambda}
|\lambda_n-\vkap_{n,*}| & \leq C\,e^{-2\pi\,(1-s_n)\,|n|\,y_0} \;,  \\
\label{eq:strip:falloff:mu}
\qmq{and} |\mu_n - (-1)^n| & \leq C\,e^{-\pi\,(1-s_n)\,|n|\,y_0} \;. 
\end{align}
\end{thm}

\newpage

\begin{proof}
For \,$y\in [-y_0,y_0]$\, we let \,$D(y) := \{(\lambda_n(y),\mu_n(y))\}_{n\in \Z}$\, be the spectral divisor of \,$u(\,\cdot\,+iy) \in \Pot$\,, then \,$D(0)=D$\, holds. 
Note that all the divisors \,$D(y)$\, lie on the same spectral curve \,$\Sigma$\, as we noted at the beginning of Section~\ref{Se:jacobitrans}. 

If \,$\Sigma$\, has a singularity at some
\,$\vkap_{n,*}$\, with \,$n\in \Z$\,, and \,$\lambda_n(y)=\vkap_{n,*}$\, holds for some \,$y\in [-y_0,y_0]$\,, then we will also have \,$\lambda_n(y)=\vkap_{n,*}$\, for all \,$y\in [-y_0,y_0]$\,; conversely if there exists
\,$y\in [-y_0,y_0]$\, with \,$\lambda_n(y)\neq \vkap_{n,*}$\,, then \,$\lambda_n(y)$\, will avoid the singularity \,$\vkap_{n,*}$\, for all times \,$y\in [-y_0,y_0]$\,. For this reason we may suppose without loss of generality
that \,$D(y)\in \Div(\Sigma')$\, holds for all \,$y \in [-y_0,y_0]$\,. 


For given \,$n\in \Z$\,, we let \,$\Psi_n$\, be as in Lemma~\ref{L:jacobiprep:Psi-intro} and define the abbreviation
$$ \Theta_n^\pm(\lambda,\mu) := \lambda-\vkap_{n,*}\pm \Psi_n(\lambda,\mu) \;. $$
We clearly have
\begin{equation}
\label{eq:strip:falloff:Theta-sum}
\Theta_n^+(\lambda,\mu)+\Theta_n^-(\lambda,\mu) = 2\cdot (\lambda-\vkap_{n,*})
\end{equation}
and because we have \,$\Psi_n \circ \sigma = -\Psi_n$\, (where \,$\sigma:\Sigma\to\Sigma$\, denotes the hyperelliptic involution of \,$\Sigma$\,), we have
\begin{equation}
\label{eq:strip:falloff:Theta-pm}
\Theta_n^\pm \circ \sigma = \Theta_n^\mp \;. 
\end{equation}
Moreover a straight-forward calculation using Equation~\eqref{eq:jacobiprep:Psi:psidef} shows 
\begin{equation}
\label{eq:strip:falloff:Theta-prod}
\Theta_n^+(\lambda,\mu) \cdot \Theta_n^-(\lambda,\mu)= (\lambda-\vkap_{n,*})^2 - \Psi_n(\lambda)^2 = \frac14\,(\vkap_{n,1}-\vkap_{n,2})^2 \; . 
\end{equation}

We now let \,$y_1 := \sign(n)\cdot y_0$\,, and 
choose the divisor \,$D^o := D(y_1) \in \Div(\Sigma')$\, as origin divisor; with this divisor we consider Jacobi coordinates \,$\vi_n(D(y))$\, for \,$D(y)$\,
like in Section~\ref{Se:jacobitrans}. 
We then have on one hand by Proposition~\ref{P:jacobi:larger}
$$ \vi_n(D(0)) = \frac{\sign(n)}{2\pi i}\cdot \ln\left( \frac{\Theta_n^+(\lambda_n(0),\mu_n(0))}{\Theta_n^+(\lambda_n(y_1),\mu_n(y_1))} \right) \cdot (1+a_n) + b_n $$
with sequences \,$(a_n),(b_n) \in \ell^2_{0,0}(n)$\,, on the other hand by Theorem~\ref{T:jacobitrans:dgl} 
$$ \vi_n(D(0)) = \vi_n(D(0))-\vi_n(D(y_1)) = (-i|n|+c_n)\cdot (-y_1) $$
with a sequence \,$c_n \in \ell^2_{-1,-1}(n)$\,. By comparing these two equations, it follows that we have
\begin{align*}
\ln\left( \frac{\Theta_n^+(\lambda_n(0),\mu_n(0))}{\Theta_n^+(\lambda_n(y_1),\mu_n(y_1))} \right) & = \frac{2\pi i \,\sign(n)\cdot \bigr((i\,|n|-c_n)\,y_1-b_n\bigr)}{1+a_n} \\
& = -2\pi n y_1 + r_n = -2\pi |n|y_0 + r_n
\end{align*}
with
$$ r_n := -2\pi n y_1 \cdot \left( \frac{1}{1+a_n}-1 \right) - \frac{2\pi i \,\sign(n)\,(c_n\,y_0 + b_n)}{1+a_n} \;\in\; \ell^2_{-1,-1}(n) \; , $$
and therefore
\begin{equation}
\label{eq:strip:falloff:y}
\Theta_n^+(\lambda_n(0),\mu_n(0)) = e^{-2\pi |n|y_0+r_n} \cdot \Theta_n^+(\lambda_n(y_1),\mu_n(y_1)) \; . 
\end{equation}

We now also consider the divisor \,$\wt{D}(0) := \sigma(D(0))$\,.
Because of Corollary~\ref{C:jacobitrans:transsigma}, the \,$y$-translational flow \,$\wt{D}(y)$\, of \,$\wt{D}(0)$\, is defined for all times \,$y\in [-y_0,y_0]$\,,
and we have 
\begin{equation}
\label{eq:strip:falloff:wtDy}
\wt{D}(y) = \sigma(D(-y)) \; . 
\end{equation}
Repeating the calculation leading to Equation~\eqref{eq:strip:falloff:y} with \,$\wt{D}(y)$\, in the place of \,$D(y)$\,, we obtain
$$ \Theta_n^+(\wt{\lambda}_n(0),\wt{\mu}_n(0)) = e^{-2\pi |n|y_0+\wt{r}_n} \cdot \Theta_n^+(\wt{\lambda}_n(y_1),\wt{\mu}_n(y_1)) $$
with another sequence \,$\wt{r}_n \in \ell^2_{-1,-1}(n)$\,. 
Because of Equations~\eqref{eq:strip:falloff:wtDy} and \eqref{eq:strip:falloff:Theta-pm}, it follows that we have
\begin{equation}
\label{eq:strip:falloff:-y}
\Theta_n^-(\lambda_n(0),\mu_n(0)) = e^{-2\pi |n|y_0+\wt{r}_n} \cdot \Theta_n^-(\lambda_n(-y_1),\mu_n(-y_1)) \; .
\end{equation}

By taking the sum of Equations~\eqref{eq:strip:falloff:y} and \eqref{eq:strip:falloff:-y} and applying Equation~\eqref{eq:strip:falloff:Theta-sum}, we obtain
\begin{align}
\label{eq:strip:falloff:lambdapre}
2\cdot (\lambda_n(0)-\vkap_{n,*}) & = e^{-2\pi |n| y_0} \cdot \bigr(e^{r_n}\cdot \Theta_n^+(\lambda_n(y_1),\mu_n(y_1)) + e^{\wt{r}_n} \cdot \Theta_n^-(\lambda_n(-y_1),\mu_n(-y_1)) \bigr) \; . 
\end{align}
We have \,$\Theta_n^\pm(\lambda_n(\pm y_1),\mu_n(\pm y_1)) \in \ell^2_{-1,3}(n)$\, and therefore in any case
\begin{equation}
\label{eq:strip:falloff:exp-estimate}
|\Theta_n^\pm(\lambda_n(\pm y_1),\mu_n(\pm y_1))| \leq C_1 \cdot |n| = C_1 \cdot e^{\wh{r}_n} 
\end{equation}
with \,$\wh{r}_n := \ln|n| \in \ell^2_{-1,-1}(n)$\, and a constant \,$C_1>0$\,. By plugging this into Equation~\eqref{eq:strip:falloff:lambdapre} we see that
$$ |\lambda_n(0)-\vkap_{n,*}| \leq C_2\cdot e^{-2\pi |n| \,(1-s_n)\,y_0} $$
holds with another constant \,$C_2>0$\, and with 
$$ s_n := \frac{1}{2\pi |n|\,y_0} \cdot \left( \max\{ \RE(r_n),\RE(\wt{r}_n)\} + \wh{r}_n \right) \;\in\;\ell^2_{0,0}(n) \; ,$$
showing Equation~\eqref{eq:strip:falloff:lambda}. 

Moreover, by taking the product of  Equations~\eqref{eq:strip:falloff:y} and \eqref{eq:strip:falloff:-y}, we obtain
\begin{gather*}
\Theta_n^+(\lambda_n(0),\mu_n(0)) \cdot \Theta_n^-(\lambda_n(0),\mu_n(0)) \\
= e^{-4\pi |n| y_0+r_n+\wt{r}_n} \cdot \Theta_n^+(\lambda_n(y_1),\mu_n(y_1)) \cdot \Theta_n^-(\lambda_n(-y_1),\mu_n(-y_1)) \; ,
\end{gather*} 
whence we obtain by applying Equation~\eqref{eq:strip:falloff:Theta-prod} and taking the square root
$$ \frac{1}{2}(\vkap_{n,1}-\vkap_{n,2}) = e^{-2\pi |n| y_0+(r_n+\wt{r}_n)/2} \cdot \bigr( \Theta_n^+(\lambda_n(y_1),\mu_n(y_1)) \cdot \Theta_n^-(\lambda_n(-y_1),\mu_n(-y_1)) \bigr)^{1/2} \; . $$
Taking absolute values and again applying the estimate \eqref{eq:strip:falloff:exp-estimate}, we obtain
$$ |\vkap_{n,1}-\vkap_{n,2}| \leq C_3 \cdot e^{-2\pi\,(1-s_n')\,|n|\,y_0} $$
with a constant \,$C_3>0$\, and 
$$ s_n' := \frac{1}{2\pi |n| \, y_0} \cdot \left( \frac{\RE(r_n+\wt{r}_n)}{2} + \wh{r}_n \right) \;\in\; \ell^2_{0,0}(n) \; . $$
This shows Equation~\eqref{eq:strip:falloff:vkap}.

Finally, because of \,$\Delta' \in \As(\C^*,\ell^\infty_{1,-3},1)$\, we have
\begin{align*}
\left( \mu_n-(-1)^n \right)^2
& = \mu_n\cdot \bigr(\mu_n+\mu_n^{-1}-2(-1)^n\bigr) = \mu_n \cdot (\Delta(\lambda_n)-\Delta(\vkap_{n,1})) \\
& = \mu_n \cdot \int_{\vkap_{n,1}}^{\lambda_n} \Delta'(\lambda)\,\mathrm{d}\lambda = \ell^\infty_{1,-3}(n) \cdot (\lambda_n-\vkap_{n,1}) \\
& = \ell^\infty_{1,-3}(n) \cdot \bigr( \, (\lambda_n-\vkap_{n,*}) \, - \, \tfrac12(\vkap_{n,1}-\vkap_{n,2}) \,\bigr) 
\end{align*}
and therefore
$$ |\mu_n-(-1)^n| \leq C_4\cdot e^{-\pi\,(1-s_n'')\,|n|\,y_0} $$
with another constant \,$C_4>0$\, and a sequence \,$(s_n'')\in \ell^2_{0,0}(n)$\,, concluding the proof of Equation~\eqref{eq:strip:falloff:mu}.
\end{proof}

\part{Perspectives}

\section{Perspectives}
\label{Se:perspectives}

In this paper we have studied simply periodic solutions \,$u:X \to \C$\, of the sinh-Gordon equation via their spectral data \,$(\Sigma,D)$\, in the style
of Bobenko. In particular we gained insight into the asymptotic behavior of the spectral data.
As was described in Section~\ref{Se:minimal}, real-valued such solutions \,$u$\, correspond to minimal immersions \,$f:X \to S^3$\, without
umbilical points.

One possible direction for the extension of this research would be to investigate minimal immersions \,$f: X \to S^3$\, of \emph{compact Riemann
surfaces} into \,$S^3$\,. Here the case where \,$X$\, is of genus \,$g=0$\, is of no interest, as only the round sphere can occur.
The case \,$g=1$\, has been classified completely by the work of \textsc{Pinkall/Sterling} \cite{Pinkall/Sterling:1989} 
resp.~\textsc{Hitchin} \cite{Hitchin:1990}. Therefore one would be interested mostly in the case where \,$X$\, is of genus \,$g\geq 2$\,.

The main obstacle for applying the theory described in this work to such compact surfaces is that a minimal immersion \,$f: X \to S^3$\,
of a compact surface of genus \,$g\geq 2$\, necessarily has umbilical points. (In fact, the Codazzi equation for the minimal immersion \,$f$\, 
shows that the Hopf differential \,$E\,\mathrm{d}z^2$\, of \,$f$\, is holomorphic, and therefore this quadratic differential has 
\,$3g-3$\, zeros on \,$X$\,.) In the vicinity of an umbilical point, \,$f$\, cannot be reparameterized in such a way that
\,$E$\, is constant non-zero, which is the pre-requisite for the conformal factor \,$u$\, of the immersion to be a solution of the
sinh-Gordon equation. 

On the other hand, when we reparameterize \,$f$\, for \,$E$\, to be constant, the umbilical points of \,$f$\, correspond to isolated
singularities of the solution \,$u$\, of the sinh-Gordon equation.

Therefore I propose as a possible next step for research in the direction of understanding compact immersed minimal surfaces in \,$S^3$\, 
(and compact immersed constant mean curvature surfaces in \,$\R^3$\,): Take a periodic solution \,$u$\, of the sinh-Gordon equation,
which may now have isolated singularities. For every base point \,$z_0$\, for which the period does not run through a singularity,
define the monodromy \,$M(\lambda)$\, and the associated spectral data \,$(\Sigma,D)$\, as before. Now investigate the asymptotic
behavior of \,$(\Sigma,D)$\, as \,$z_0$\, approaches a singularity, and use this information to relate the monodromy resp.~the 
spectral data for a period ``above'' the isolated singularity to these data for a period ``below'' the singularity.

By pursuing this path of research, one might gain an overview of how the monodromy resp.~the spectral data of the compact minimal immersed
surface of genus \,$g\geq 2$\, are related for any two different base points \,$z_0,z_1$\,. In this way, one might gain some insight into the structure
of such surfaces. The distant goal would of course be a complete classification of these surfaces.

\appendix

\part{Appendices}


\section{Some infinite sums and products}
\label{Ap:inf}

In the present section we prove the convergence and estimates for some infinite sums and products. These results are required for the proofs in Section~\ref{Se:interpolate}. 
We shall work in a setting that is modelled after the one ``half'' of the 
setting of 
Sections~\ref{Se:vacuum} and \ref{Se:excl}, in that we will consider sequences \,$(\lambda_k)_{k\geq 1}$\, only for positive indices \,$k$\,, corresponding to the ``half'' of the sequence
that goes to \,$\infty$\,. We will work in a somewhat more general setting than in the main text in that we fix a sequence \,$(\lambda_{k,0})_{k\in \Z}$\, in \,$\C^*$\, with
\begin{equation}
\label{eq:inf:lambdak0}
\lambda_{k,0}=ck^2 + O(1) 
\end{equation}
for \,$k\geq 1$\, and a fixed constant \,$c\in \C^*$\,, and then define ``excluded domains'' for \,$\delta>0$\, and \,$k\geq 1$\, for the purposes of this appendix by
$$ U_{k,\delta} := 
\Mengegr{\lambda\in \C}{|\lambda-\lambda_{k,0}|<\delta\,k} 
\;; $$
we also put \,$U_\delta := \bigcup_{k\geq 1} U_{k,\delta}$\, and \,$V_\delta := \C\setminus U_\delta$\,. Proposition~\ref{P:vac2:excldom-new}(2) shows that this definition of excluded domains
is essentially equivalent to the one in Definition~\ref{D:vac2:excldom}.

In an analogous vein, we define the annuli \,$S_k$\, for \,$k\geq 1$\, for the purposes of this appendix by
$$ S_k := 
\Mengegr{\lambda\in \C}{|c|\,(k-\tfrac12)^2 \leq |\lambda| \leq |c|\,(k+\tfrac12)^2} 
\; . $$
Again, this definition is equivalent to the one in Equation~\eqref{eq:As:Sk}. The \,$S_k$\, intersect only on their boundary,
and we have \,$\bigcup_{k\geq 1} S_k = \Mengegr{\lambda\in \C}{|\lambda| \geq |c|/4}$\,. 

\begin{prop}[Holomorphic functions defined by infinite sums.]
\label{P:inf:sumholo}
Let \,$R_0>0$\, and sequences \,$(\lambda_k)_{k\geq 1},(\lambda_k^{[1]})_{k\geq 1},(\lambda_k^{[2]})_{k\geq 1}$\, with \,$\|\lambda_k-\lambda_{k,0}\|_{\ell^2_{-1}}, \|\lambda_k^{[\nu]}-\lambda_{k,0}\|_{\ell^2_{-1}} \leq R_0$\, be given, 
so that there exists \,$0<\delta_0\leq \tfrac{|c|}{4}$\, such that \,$\lambda_k,\lambda_k^{[1]},\lambda_k^{[2]} \in U_{k,\delta_0}$\, holds for every \,$k\geq 1$\, with at most finitely many exceptions.
Also let a further sequence \,$(a_k)_{k\geq 1}$\, be given. 
\begin{enumerate}
\item \,$\sum_{k=1}^\infty \frac{1}{\lambda-\lambda_k}$\, defines a holomorphic function on \,$V_0 := \C \setminus \Menge{\lambda_k}{k\in\N}$\,.
For every \,$\delta>\delta_0$\, there
exist \,$R,C>0$\, (depending only on \,$\delta$\,, \,$\delta_0$\, and \,$R_0$\,) such that for every \,$\lambda\in V_\delta$\, with \,$|\lambda|\geq R$\, we have
$$ \left| \sum_{k=1}^\infty \frac{1}{\lambda-\lambda_k} \right| \leq C\cdot|\lambda|^{-1/2} \; . $$ 
\item If \,$a_k \in \ell^2_{-1}(k)$\, and \,$\|a_k\|_{\ell^2_{-1}} \leq R_0$\, holds, then \,$\sum_{k=1}^\infty \frac{a_k}{\lambda-\lambda_k}$\, defines a holomorphic function on \,$V_0$\,.
For every \,$\delta>\delta_0$\, there exists a constant \,$C>0$\, (depending only on \,$\delta$\,, \,$\delta_0$\, and \,$R_0$\,), such that with the 
sequence \,$(r_n)_{n\geq 1} \in \ell^2(n)$\, given by
$$ r_n = C \cdot \left( \frac{|a_k|}{k} * \frac{1}{|k|} \right)_n $$
(where we put \,$\tfrac{1}{0}:=1$\,, and extend \,$\tfrac{|a_k|}{k}=0$\, for \,$k\leq 0$\,),
we have for every \,$n\geq 1$\, and \,$\lambda \in S_n \cap V_\delta$\, 
\begin{equation}
\label{eq:inf:sumholo:2:claim}
\sum_{k=1}^\infty \frac{|a_k|}{|\lambda-\lambda_k|} \leq r_n \; .
\end{equation}
\item 
If in fact \,$a_k \in \ell^2(k)$\, and \,$\|a_k\|_{\ell^2} \leq R_0$\, holds in the situation of (2), then
there exists a constant \,$C>0$\, (depending only on \,$\delta$\,, \,$\delta_0$\, and \,$R_0$\,), such that 
\eqref{eq:inf:sumholo:2:claim} holds with the sequence \,$(r_n)_{n\geq 1} \in \ell^2_1(n)$\, given by
\begin{equation}
\label{eq:inf:sumholo:3:rn}
r_n = \frac{C}{n} \cdot \left( |a_k| * \frac{1}{|k|} \right)_n \;.
\end{equation}
\item If \,$a_k \in \ell^2_{-1}(k)$\, and \,$\|a_k\|_{\ell^2_{-1}} \leq R_0$\, holds, 
then \,$\sum_{k=1}^\infty \frac{a_k}{(\lambda-\lambda_k^{[1]})\cdot (\lambda-\lambda_k^{[2]})}$\, defines a holomorphic function on \,$\C \setminus \Menge{\lambda_k^{[1]},\lambda_k^{[2]}}{k\in\N}$\,.
For every \,$\delta>\delta_0$\, there exists a constant \,$C>0$\, (depending only on \,$\delta$\,, \,$\delta_0$\, and \,$R_0$\,) such that with the sequence \,$(r_n)_{n\geq 1} \in \ell^2_1(n)$\, 
given by Equation~\eqref{eq:inf:sumholo:3:rn} we have for every \,$n\geq 1$\, and \,$\lambda \in S_n\cap V_\delta$\, 
$$ \sum_{k=1}^\infty \frac{|a_k|}{|\lambda-\lambda_k^{[1]}|\cdot |\lambda-\lambda_k^{[2]}|} \leq r_n \; . $$ 
\end{enumerate}
\end{prop}

\begin{proof}
We may suppose without loss of generality that \,$\lambda_k, \lambda_k^{[1]},\lambda_k^{[2]} \in U_{k,\delta_0}$\, holds for all \,$k\geq 1$\,. We then have
\begin{equation}
\label{eq:inf:sumholo:lambdak-ck2}
|\lambda_k -ck^2| \leq |\lambda_k-\lambda_{k,0}| + |\lambda_{k,0}-ck^2| \leq \delta_0\cdot k + C_1 \leq \frac{|c|}{4}\cdot k+C_1
\end{equation}
with a constant \,$C_1>0$\,.

To prepare the proof of the proposition, we show the following approximation:
There exists \,$N\in \N$\,, such that for every \,$n\geq N$\,, \,$\lambda \in S_n$\, and \,$k\in \N$\,, we have
\begin{equation}
\label{eq:inf:sumholo:lambdaest}
|\lambda-\lambda_k| \geq \frac{|c|}{4}\cdot |n^2-k^2| \; .
\end{equation}
Indeed, let \,$n,k\in \N$\, be given. 
In the case \,$k<n$\,, the \,$\lambda\in S_n$\, for which \,$|\lambda-ck^2|$\, is minimal is \,$\lambda=c(n-\tfrac12)^2$\,,
and therefore we have
\begin{align*}
|\lambda-\lambda_k| & \;\;\geq\;\; |\lambda-ck^2| - |\lambda_k-ck^2| \geq \left|c(n-\tfrac12)^2-ck^2 \right| - |ck^2-\lambda_k| \\
& \overset{\eqref{eq:inf:sumholo:lambdak-ck2}}{\geq} \left(|c|\,(n-\tfrac12)^2-|c|\,k^2 \right) - \frac{|c|}{4}\,k-C_1 
= |c|\cdot (n^2-n+\tfrac14-k^2-\tfrac14 k)-C_1 \\
& \;\;\geq\;\; \frac{|c|}{4}\cdot (n^2-k^2) + \frac{|c|}{4}\cdot (3n^2-4n+1-3k^2-k) - C_1 \\
& \;\;\overset{(*)}{\geq}\;\; \frac{|c|}{4}\cdot (n^2-k^2) + \frac{|c|}{4}\cdot (3n^2-4n+1-3(n-1)^2-(n-1)) - C_1 \\
& \;\;=\;\; \frac{|c|}{4}\cdot (n^2-k^2) + \frac{|c|}{4}\cdot (n-1) - C_1 
\;\;\overset{(\dagger)}{\geq}\;\; \frac{|c|}{4}\cdot (n^2-k^2)\;,
\end{align*}
where the \,$\geq$\, sign marked $(*)$ is true because of \,$k\leq n-1$\,, and the \,$\geq$\, sign marked $(\dagger)$ holds for 
all \,$n\geq \tfrac{4}{|c|}\,C_1+1$\,. Therefore \eqref{eq:inf:sumholo:lambdaest} holds for such \,$n$\, and all \,$k<n$\,. 

In the case \,$k=n$\, there is nothing to show. And in the case \,$k>n$\,, the \,$\lambda\in S_n$\, for which \,$|\lambda-ck^2|$\, 
is minimal is \,$\lambda=c(n+\tfrac12)^2$\,, and therefore we have
\begin{align*}
|\lambda-\lambda_k| & \;\;\geq\;\; |\lambda-ck^2| - |\lambda_k-ck^2| \geq \left|c(n+\tfrac12)^2-ck^2 \right| - |ck^2-\lambda_k| \\
& \overset{\eqref{eq:inf:sumholo:lambdak-ck2}}{\geq} \left(|c|\,k^2- |c|\,(n+\tfrac12)^2 \right) - \frac{|c|}{4}\cdot k-C_1 \\
& \;\;=\;\; |c|\,(k^2-n^2) - |c|\,(n+\tfrac14\,k)-(C_1+\tfrac14\,|c|) \\
& \;\;=\;\; \tfrac{|c|}{4}\,(k^2-n^2) + |c|\,\left(\tfrac12\,(k^2-n^2)-n\right) + \tfrac{|c|}{4}\,\left(\tfrac12\,(k^2-n^2)-k+\tfrac12\right) \\ 
& \qquad\qquad + \left(\tfrac18\,(k^2-n^2)-(C_1+\tfrac38\,|c|)\right) \; . 
\end{align*}
Because of \,$k\geq n+1$\, we have
$$ k^2-n^2 \geq (n+1)^2-n^2 = 2n+1 \geq 2n \qmq{and thus} \tfrac12(k^2-n^2)-n \geq 0 $$
and because the function \,$f(x)=\tfrac12\,x^2-x$\, is monotonously increasing for \,$x\geq 1$\,,
$$ \tfrac12\,(k^2-n^2)-k+\tfrac12 = \tfrac12\,k^2-k -\tfrac12\,n^2+\tfrac12 \geq \tfrac12\,(n+1)^2-(n+1) -\tfrac12\,n^2+\tfrac12 = 0 \; . $$
Thus we obtain
\begin{align*}
|\lambda-\lambda_k| & 
\geq \tfrac{|c|}{4}\,(k^2-n^2) + (\tfrac14\,n-(C_1+\tfrac38\,|c|)) \geq \tfrac{|c|}{4}\,(k^2-n^2) \;, 
\end{align*}
where the last \,$\geq$\, sign holds whenever \,$n \geq 4\,C_1+\tfrac32\,|c|$\,, and therefore \eqref{eq:inf:sumholo:lambdaest} holds
for these \,$n$\, and \,$k>n$\,. 

The preceding calculations show that \eqref{eq:inf:sumholo:lambdaest} holds for all \,$k\in \N$\, and \,$n\geq N$\,, if we choose \,$N$\, as the smallest integer larger than \,$\max\{\tfrac{4}{|c|}\,C_1+1, 4\,C_1+\tfrac32\,|c|\}$\,. 


We now turn to the proof of (1)--(4). For proving the estimates claimed in the respective parts of the proposition, it suffices to consider \,$\lambda \in S_n$\, with \,$n\geq N$\,; the estimates can then be extended
to \,$\lambda\in \C$\, by enlarging the constant \,$C$\, via an argument of compactness.

\emph{For (1).}
The series \,$\sum_{k=1}^\infty \tfrac{1}{\lambda-\lambda_k}$\, converges absolutely and locally uniformly for \,$\lambda\in V_0$\,, and thus defines a holomorphic function. Now
let \,$\delta>\delta_0$\, and \,$\lambda\in V_\delta$\, be given. We write 
\begin{equation}
\label{eq:inf:sumholo:1:master}
\left| \sum_{k=1}^\infty \frac{1}{\lambda-\lambda_k} \right| \leq \left| \sum_{k=1}^\infty \frac{1}{\lambda-ck^2} \right| + \sum_{k=1}^\infty \left| \frac{1}{\lambda-\lambda_k} - \frac{1}{\lambda-ck^2} \right| \; .
\end{equation}
To estimate \,$\sum_{k=1}^\infty \tfrac{1}{\lambda-ck^2}$\,, we use the partial fraction decomposition of the cotangent function
$$ \pi\,\cot(\pi z) = z^{-1} + \sum_{k=1}^\infty \frac{2z}{z^2-k^2} \qmq{for \,$z\in \C \setminus \Menge{k\pi}{k\in \Z}$\,}\;, $$
which implies
\begin{equation}
\label{eq:inf:sumholo:1:a}
\sum_{k=1}^\infty \frac{1}{\lambda-ck^2} = \frac{1}{2\sqrt{c\,\lambda}}\,\left( \pi\,\cot(\pi\,\sqrt{c^{-1}\,\lambda}) - \sqrt{c\,\lambda^{-1}} \right) = O\left( |\lambda|^{-1/2} \right) 
\end{equation}
for \,$\lambda\in V_\delta$\,; note that \,$\cot(\pi\sqrt{c^{-1}\,\lambda})$\, is bounded on \,$V_\delta$\,. 

To estimate the second summand in \eqref{eq:inf:sumholo:1:master}, we note that we have for \,$\lambda \in S_n \cap V_\delta$\, with \,$n\geq N$\,
if \,$k\neq n$\, by \eqref{eq:inf:sumholo:lambdak-ck2} and \eqref{eq:inf:sumholo:lambdaest}
$$ \left| \frac{1}{\lambda-\lambda_k} - \frac{1}{\lambda-ck^2} \right| = \left| \frac{\lambda_k-ck^2}{(\lambda-\lambda_k)(\lambda-ck^2)} \right| \leq \frac{16}{|c|^2}\frac{\tfrac{|c|}{4}k+C_1}{|n^2-k^2|^2} \leq C_2 \cdot \frac{k}{|n^2-k^2|^2} $$
with a constant \,$C_2>0$\,, and if \,$k=n$\,
$$ \left| \frac{1}{\lambda-\lambda_n} - \frac{1}{\lambda-cn^2} \right| = \left| \frac{\lambda_n-cn^2}{(\lambda-\lambda_n)(\lambda-cn^2)} \right| \leq \frac{n\,\delta_0}{(n\,(\delta-\delta_0)) \cdot (n\,\delta)} \leq \frac{1}{n\,(\delta-\delta_0)} \;. $$
Thus we obtain
\begin{align}
\left| \sum_{k=1}^\infty \left( \frac{1}{\lambda-\lambda_k} - \frac{1}{c\lambda-k^2} \right) \right|
& \leq \frac{1}{n\,(\delta-\delta_0)} + C_2\sum_{k\neq n}\frac{k}{|n^2-k^2|^2} \notag \\
\label{eq:inf:sumholo:1:b}
& \leq \frac{1}{n\,(\delta-\delta_0)} + \frac{C_2}{n}\sum_{k\neq n}\frac{1}{|n-k|^2} = O(n^{-1})\;.
\end{align}
Because \,$n^{-1}$\, is comparable to \,$|\lambda|^{-1/2}$\, for \,$\lambda\in S_n$\,, it follows by plugging \eqref{eq:inf:sumholo:1:a} and \eqref{eq:inf:sumholo:1:b} into \eqref{eq:inf:sumholo:1:master} that
$$ \left| \sum_{k=1}^\infty \frac{1}{\lambda-\lambda_k} \right| = O(|\lambda|^{-1/2}) $$
holds for \,$\lambda\in V_\delta$\,, and this is the claimed statement.

\emph{For (2).}
The series \,$\sum_{k=1}^\infty \tfrac{a_k}{\lambda-\lambda_k}$\, converges absolutely and locally uniformly for \,$\lambda\in V_0$\,, and thus defines a holomorphic function. 
For \,$n\geq N$\, and \,$\lambda\in S_n \cap V_\delta$\, we have
\begin{align}
\sum_{k=1}^\infty \frac{|a_k|}{|\lambda-\lambda_k|} 
& \overset{\eqref{eq:inf:sumholo:lambdaest}}{\leq} \frac{4}{|c|} \sum_{k\neq n} \frac{|a_k|}{|n^2-k^2|} + \frac{|a_n|}{|\lambda-\lambda_n|} 
\leq \frac{4}{|c|} \sum_{k\neq n} \frac{|a_k|}{k}\cdot \frac{1}{|n-k|} + \frac{|a_n|}{n\cdot (\delta-\delta_0)} \notag \\
\label{eq:inf:sumholo:2:est}
& \;\;=\;\; \frac{4}{|c|} \, \left( \frac{|a_k|}{k} * \frac{1}{|k|} \right)_{n} + \frac{1}{\delta-\delta_0} \cdot \frac{|a_n|}{n} 
\leq C \cdot \left( \frac{|a_k|}{k} * \frac{1}{|k|} \right)_n =: r_n \;, 
\end{align}
where we put \,$C := \max\{ \tfrac{4}{|c|}, \frac{1}{\delta-\delta_0}\}$\,, and we 
we extend \,$\tfrac{|a_k|}{k}$\, to \,$k\in \Z$\, by setting \,$\tfrac{|a_k|}{k}=0$\, for \,$k\leq 0$\,, and again put \,$\tfrac{1}{0} := 1$\,. 
The convolution of the \,$\ell^2$-sequence \,$\tfrac{|a_k|}{k}$\, with \,$\tfrac{1}{|k|}$\, is 
again an \,$\ell^2$-sequence by the version \eqref{eq:fasymp:fourier-small:weakyoung} of Young's inequality for the convolution with \,$\tfrac{1}{k}$\,. 
This shows that \,$r_n \in \ell^2(n)$\, holds. 


\emph{For (3).}
We proceed similarly as in (2): For \,$n\geq N$\, and \,$\lambda\in S_n \cap V_\delta$\, we now have
\begin{align}
\sum_{k=1}^\infty \frac{|a_k|}{|\lambda-\lambda_k|}
& \overset{\eqref{eq:inf:sumholo:lambdaest}}{\leq} \frac{4}{|c|} \sum_{k\neq n} \frac{|a_k|}{|n^2-k^2|} + \frac{|a_n|}{|\lambda-\lambda_n|} 
\leq \frac{4}{|c|} \sum_{k\neq n} \frac{|a_k|}{n}\cdot \frac{1}{|n-k|} + \frac{|a_n|}{n\cdot (\delta-\delta_0)} \notag \\
& \;\;=\;\; \frac{1}{n} \cdot \left( \frac{4}{|c|}\,\left( |a_k| * \frac{1}{|k|} \right)_{n} + \frac{1}{\delta-\delta_0}\cdot |a_n| \right) \leq \frac{C}{n}\cdot \left( |a_k| * \frac{1}{|k|} \right)_{n} =:r_n\; ,
\end{align}
where we again extend \,$|a_k|$\, to \,$k\in \Z$\, by zero, and where \,$C>0$\, is a constant.
The convolution of \,$|a_k|$\, with \,$\frac{1}{|k|}$\, is an \,$\ell^2$-sequence by \eqref{eq:fasymp:fourier-small:weakyoung}, hence \,$r_n \in \ell^2_1(n)$\, holds.


\emph{For (4).} For the proof of (4), we choose \,$N$\, so that \eqref{eq:inf:sumholo:lambdaest} applies for both \,$(\lambda_n^{[1]})$\, and \,$(\lambda_n^{[2]})$\, in the place of \,$(\lambda_n)$\,. 
Then we have
\begin{align*}
\sum_{k=1}^\infty \frac{a_k}{(\lambda-\lambda_k^{[1]}) \cdot (\lambda-\lambda_k^{[2]})}
& \overset{\eqref{eq:inf:sumholo:lambdaest}}{\leq} \frac{16}{|c|^2} \sum_{k\neq n} \frac{|a_k|}{|n^2-k^2|^2} + \frac{|a_n|}{|\lambda-\lambda_n^{[1]}| \cdot |\lambda-\lambda_n^{[2]}|} \\
& \;\;\leq\;\; \frac{16}{|c|^2}\cdot \frac{1}{n} \sum_{k\neq n} \frac{|a_k|}{k}\cdot \frac{1}{|n-k|^2} + \frac{|a_n|}{(n\,(\delta-\delta_0))^2} \\
& \;\;=\;\; \frac{1}{n} \cdot \left( \frac{16}{|c|^2} \left( \frac{|a_k|}{k} * \frac{1}{k^2} \right)_n + \frac{1}{(\delta-\delta_0)^2} \cdot \frac{|a_n|}{n} \right) \\
& \;\;\leq\;\; \frac{C}{n}\cdot \left( \frac{|a_k|}{k} * \frac{1}{k} \right)_n  =: r_n \; . 
\end{align*}
Here, \,$C>0$\, is again a constant, and similarly as before, we see that \,$r_n \in \ell^2_1(n)$\, holds. 
\end{proof}

The relationship between convergent infinite products \,$\prod_{k=1}^\infty (1+a_k)$\, and absolutely convergent series \,$\sum_{k=1}^\infty a_k$\, described in the following proposition is very well-known
(see, for example, \cite{Conway:1978}, Corollary~VII.5.6 and Lemma~VII.5.8, p.~166), as is the estimate of the product by the \,$\ell^1$-norm of \,$(a_k)$\,. We state and prove this result here
for the sake of completeness.


\begin{prop}[Estimating products by sums.]
\label{P:inf:prod2sum}
Let a sequence \,$(a_k) \in \ell^1(k)$\, be given. 
The infinite product \,$\prod_{k=1}^\infty (1+a_k)$\, converges in \,$\C^*$\,, and we have 
$$ \left| \prod_{k=1}^\infty (1+a_k) - 1 \right| \leq \exp\left( \|a_k\|_{\ell^1} \right) - 1 \; . $$
Moreover, if the \,$(a_k)$\, depend on a parameter so that \,$\sum_{k=1}^\infty |a_k|$\, converges uniformly with respect to that parameter,
then the convergence of \,$\prod_{k=1}^\infty (1+a_k)$\, is also uniform.
\end{prop}

\begin{proof}
For \,$N \in \N$\, we have
$$ \prod_{k=1}^N (1+a_k) - 1 = \sum_{k=1}^N \left( \sum_{1\leq j_1<\dotsc<j_k\leq N} a_{j_1}\cdot\dotsc\cdot a_{j_k} \right) \; . $$
For \,$1 \leq k \leq N$\, we have
\begin{align*}
\left| \sum_{1\leq j_1<\dotsc<j_k\leq N} a_{j_1}\cdot\dotsc\cdot a_{j_k} \right|
& \leq \frac{1}{k!} \sum_{1 \leq j_1,\dotsc,j_k \leq N} |a_{j_1}|\cdot\dotsc\cdot|a_{j_k}| = \frac{1}{k!} \left( \sum_{j=1}^N |a_j|\right)^k \leq \frac{1}{k!} \|a_j\|_{\ell^1(j)}^k
\end{align*}
and therefore
$$ \left| \prod_{k=1}^N (1+a_k) - 1 \right| \leq \sum_{k=1}^N \frac{1}{k!} \|a_j\|_{\ell^1(j)}^k \leq \sum_{k=1}^\infty \frac{1}{k!} \|a_j\|_{\ell^1(j)}^k = \exp(\|a_j\|_{\ell^1(j)})-1 \; . $$
By taking the limit \,$N\to\infty$\,, we obtain the claimed result.
\end{proof}

\begin{prop}[An infinite product.]
\label{P:inf:prod}
Let \,$R_0>0$\, and
sequences \,$(\lambda_k^{[1]})_{k\geq 1}$\,, \,$(\lambda_k^{[2]})_{k \geq 1}$\, with \,$\lambda_k^{[1]}-\lambda_k^{[2]} \in \ell^2_{-1}(k)$\, and \,$\|\lambda_k^{[1]}-\lambda_k^{[2]}\|_{\ell^2_{-1}} \leq R_0$\, be given.
Then the infinite product \,$\prod_{k=1}^\infty \frac{\lambda_k^{[1]}}{\lambda_k^{[2]}}$\, converges in \,$\C^*$\,, and there exists a constant \,$C>0$\,, dependent only on \,$R_0$\,, such that
$$ \left| \prod_{k=1}^\infty \frac{\lambda_k^{[1]}}{\lambda_k^{[2]}} - 1 \right| \leq C \cdot \|\lambda_k^{[1]}-\lambda_k^{[2]}\|_{\ell^2_{-1,3}} $$
holds.
\end{prop}

\begin{proof}
We have \,$\prod_{k=1}^\infty \frac{\lambda_k^{[1]}}{\lambda_k^{[2]}} = \prod_{k=1}^\infty(1+a_k)$\, with 
$$ a_k := \frac{\lambda_k^{[1]}-\lambda_k^{[2]}}{\lambda_k^{[2]}} \; . $$
There exists \,$C_1>0$\, with \,$|\lambda_k^{[\nu]}| \geq C_1\cdot k^2$\, and therefore we have by Cauchy-Schwarz's inequality:
\begin{align*}
\|a_k\|_{\ell^1} 
& = \left\| \frac{\lambda_k^{[1]}-\lambda_k^{[2]}}{\lambda_k^{[2]}} \right\|_{\ell^1} \leq \frac{1}{C_1} \cdot \left\|\frac{\lambda_k^{[1]}-\lambda_k^{[2]}}{k^2} \right\|_{\ell^1} \\
& \leq \frac{1}{C_1}\cdot \left\| \frac{1}{k} \right\|_{\ell^2} \cdot \left\| \frac{\lambda_k^{[1]}-\lambda_k^{[2]}}{k} \right\|_{\ell^2} 
\leq C_2 \cdot \left\| \lambda_k^{[1]}-\lambda_k^{[2]} \right\|_{\ell^2_{-1}} \; . 
\end{align*}
It follows by Proposition~\ref{P:inf:prod2sum}(1) that \,$\prod_{k=1}^\infty \frac{\lambda_k^{[1]}}{\lambda_k^{[2]}}$\, converges, and that
$$ \left| \prod_{k=1}^\infty \frac{\lambda_k^{[1]}}{\lambda_k^{[2]}} -1 \right| \leq \exp\left( \|a_k\|_{\ell^1} \right)-1 \leq \exp\left( C_2 \cdot \left\| \lambda_k^{[1]}-\lambda_k^{[2]} \right\|_{\ell^2_{-1}} \right)-1 $$
holds. Because \,$\exp$\, is Lipschitz continuous on \,$[0,C_2\cdot R_0]$\,, it follows that there exists a constant \,$C_3>0$\, so that 
$$ \left| \prod_{k=1}^\infty \frac{\lambda_k^{[1]}}{\lambda_k^{[2]}} -1 \right| \leq  C_3 \cdot \left\| \lambda_k^{[1]}-\lambda_k^{[2]} \right\|_{\ell^2_{-1}} $$
holds.
\end{proof}

\begin{prop}[Holomorphic functions defined by infinite products.]
\label{P:inf:prodholo}
\strut
\begin{enumerate}
\item 
Let \,$(\lambda_k)_{k\geq 1}$\, be given with \,$\lambda_k-\lambda_{k,0}\in \ell^2_{-1}(k)$\,. Then 
\,$f(\lambda) := \prod_{k=1}^\infty \left( 1-\frac{\lambda}{\lambda_k} \right)$\, defines a holomorphic function on \,$\C$\, with \,$f(0)=1$\,. 
\item 
Let \,$R_0>0$\, 
and sequences \,$(\lambda_k^{[1]})_{k\geq 1}, (\lambda_k^{[2]})_{k\geq 1}$\, with 
\,$\lambda_k^{[1]}-\lambda_k^{[2]} \in \ell^2_{-1}(k)$\, and \,$\|\lambda_k^{[1]}-\lambda_k^{[2]}\|_{\ell^2_{-1}} \leq R_0$\, be given,
so that there exists \,$0<\delta_0\leq \tfrac{|c|}{4}$\, such that \,$\lambda_k^{[\nu]} \in U_{k,\delta_0}$\, holds for every \,$k\geq 1$\, with at most
finitely many exceptions.

Then \,$g(\lambda) := \prod_{k=1}^\infty \frac{\lambda_k^{[1]}-\lambda}{\lambda_k^{[2]}-\lambda}$\, defines a holomorphic function on \,$\C \setminus \Menge{\lambda_k^{[2]}}{k \geq 1}$\,
with \,$g(0)=\prod_{k=1}^\infty \frac{\lambda_k^{[1]}}{\lambda_k^{[2]}}$\,, 
and \,$h(\lambda) := \prod_{k=1}^\infty \frac{(\lambda_k^{[1]})^{-1}-\lambda}{(\lambda_k^{[2]})^{-1}-\lambda}$\, defines a holomorphic function on \,$(\C^* \cup \{\infty\}) \setminus \Menge{(\lambda_k^{[2]})^{-1}}{k \geq 1}$\,
with \,$h(\infty) = 1$\,. 

For \,$\delta>\delta_0$\,, there exists a constant \,$C>0$\,, depending only on \,$R_0$\, and \,$\delta$\,, such that with the sequence \,$(r_n)_{n \in \Z} \in \ell^2(n)$\, given by
\begin{equation*}
r_n = C \cdot \left( \frac{|\lambda_k^{[1]}-\lambda_k^{[2]}|}{k} * \frac{1}{|k|} \right)_n
\end{equation*}
(where we put \,$\tfrac{1}{0}:=1$\,, and extend \,$\tfrac{|\lambda_k^{[1]}-\lambda_k^{[2]}|}{k}=0$\, for \,$k\leq 0$\,), we have for every \,$n\geq 1$\, and \,$\lambda \in S_n \cap V_\delta$\, 
\begin{equation}
\label{eq:inf:prodholo:2:claim-g}
|g(\lambda)-1| \leq r_n 
\end{equation}
and
\begin{equation}
\label{eq:inf:prodholo:2:claim-h}
| h(\lambda)-1| \leq r_{-n} \;.  
\end{equation}
\end{enumerate}
\end{prop}

\begin{rem}
In Proposition~\ref{P:inf:prodholo}(2), in fact a far better estimate than \eqref{eq:inf:prodholo:2:claim-h} is possible regarding the function \,$h$\,: in fact \,$|h(\lambda)-1|\leq \wh{r}_{-n}$\, holds with
\begin{equation*}
\wh{r}_n = C \cdot \left( \frac{|\lambda_k^{[1]}-\lambda_k^{[2]}|}{k^3} * \frac{1}{|k|} \right)_n \; . 
\end{equation*}
\end{rem}

\begin{proof}[Proof of Proposition~\ref{P:inf:prodholo}.]
\emph{For (1).}
We have \,$\prod_{k=1}^\infty \left( 1 -\tfrac{\lambda}{\lambda_k} \right) = \prod_{k=1}^\infty (1+a_k(\lambda))$\, with \,$a_k(\lambda) := -\tfrac{\lambda}{\lambda_k}$\,, and 
$$ \|a_k(\lambda)\|_{\ell^1(k)} = |\lambda| \cdot \sum_{k=1}^\infty \underbrace{\frac{1}{|\lambda_k|}}_{=O(k^{-2})} \leq C_1\cdot |\lambda| $$
holds with some constant \,$C_1>0$\,. Therefore it follows from Proposition~\ref{P:inf:prod2sum}(1) that the infinite product \,$\prod_{k=1}^\infty (1+a_k(\lambda))$\, converges locally uniformly in \,$\lambda$\,,
hence this product defines a holomorphic function \,$f$\, in \,$\lambda \in \C$\,. Moreover, we have \,$a_k(0)=0$\, for all \,$k$\,, and therefore \,$f(0)=1$\,.

\emph{For (2).}
We have \,$g(\lambda) = \prod_{k=1}^\infty (1+a_k(\lambda))$\, with \,$a_k(\lambda) := \tfrac{\lambda_k^{[1]}-\lambda_k^{[2]}}{\lambda_k^{[2]}-\lambda}$\,. Because we have \,$|\lambda_k^{[1]}-\lambda_k^{[2]}| \in \ell^2_{-1}(k)$\,,
Proposition~\ref{P:inf:sumholo}(2) is applicable to \,$\sum_{k=1}^\infty |a_k(\lambda)|$\,. By that Proposition, \,$\sum_{k=1}^\infty a_k(\lambda)$\, defines a holomorphic function
on \,$\lambda \in \C \setminus \Menge{\lambda_{k}^{[2]}}{k \geq 1}$\,. Moreover, if we suppose that \,$\delta>\delta_0$\, is given, there exists a constant \,$C_2>0$\, such that with 
the sequence \,$(\wt{r}_n)_{n\geq 1} \in \ell^2(n)$\, given by
$$ \wt{r}_n = C_2 \cdot \left( \frac{|\lambda_k^{[1]}-\lambda_k^{[2]}|}{k} * \frac{1}{|k|} \right)_n $$
we have for any \,$n\geq 1$\, and \,$\lambda \in S_n \cap V_\delta$\,
$$ \|a_k(\lambda)\|_{\ell^1} = \sum_{k=1}^\infty \frac{|\lambda_k^{[1]}-\lambda_k^{[2]}|}{|\lambda_k^{[2]}-\lambda|} \leq \wt{r}_n \; .  $$
By Proposition~\ref{P:inf:prod2sum}(1) it follows that
\begin{equation}
\label{eq:inf:prodholo:2:pre}
|g(\lambda)-1| \leq \exp(\wt{r}_n)-1 
\end{equation}
holds.

By the variant of Young's inequality for weak \,$\ell^1$-sequences (see \eqref{eq:fasymp:fourier-small:weakyoung}) there exists a constant \,$C_3>0$\, so that we have for every \,$n \in \N$\,
$$ |\wt{r}_n| \leq \|\wt{r}_k\|_{\ell^2(k)} \leq C_2 \cdot C_3 \cdot \|\lambda_k^{[1]}-\lambda_k^{[2]}\|_{\ell^2_{-1}(k)} \leq C_2 \cdot C_3 \cdot R_0 \;, $$
and therefore there exists a constant \,$C_4>0$\, so that 
$$ \exp(\wt{r}_n)-1 \leq C_4 \cdot \wt{r}_n $$
holds. From the estimate~\eqref{eq:inf:prodholo:2:pre} we thus obtain that for \,$n\geq 1$\, and \,$\lambda \in S_n \cap V_\delta$\,,
$$ |g(\lambda)-1| \leq C_4 \cdot \wt{r}_n $$
holds; therefore the claimed estimate \eqref{eq:inf:prodholo:2:claim-g} holds with \,$C := C_4 \cdot C_2 >0$\,. 

Similarly, we have \,$h(\lambda) = \prod_{k=1}^\infty (1+b_k(\lambda))$\, with 
$$ b_k(\lambda) := \tfrac{(\lambda_k^{[1]})^{-1}-(\lambda_k^{[2]})^{-1}}{(\lambda_k^{[2]})^{-1}-\lambda} = \frac{\lambda_k^{[2]}-\lambda_k^{[1]}}{\lambda_k^{[1]}\cdot(1-\lambda_k^{[2]}\cdot\lambda)} \; . $$
We have \,$|\lambda_k^{[\nu]}| \geq C_5\cdot k^2$\, and for \,$\lambda \in S_n$\,, \,$n\in \N$\,,
$$ |1-\lambda_k^{[2]}\cdot\lambda| \geq C_6 \cdot |1-k^2\cdot n^2| = C_6 \cdot |1-k\,n|\cdot |1+k\,n| \geq C_6 \cdot k \cdot |1-k\,n| \geq C_7 \cdot (k+n) \; , $$
and therefore
$$ |b_k(\lambda)| \leq C_8 \cdot \frac{|\lambda_k^{[1]}-\lambda_k^{[2]}|}{k} \cdot \frac{1}{k+n} \;, $$
whence
$$ \|b_k(\lambda)\|_{\ell^1} \leq C_8 \cdot \sum_{k=1}^\infty \frac{|\lambda_k^{[1]}-\lambda_k^{[2]}|}{k} \cdot \frac{1}{k+n} = C_8 \cdot \left( \frac{|\lambda_k^{[1]}-\lambda_k^{[2]}|}{k} * \frac{1}{|k|} \right)_{-n} = C_9 \cdot r_{-n} $$
follows. By an analogous application of Proposition~\ref{P:inf:prod2sum}(1) and Young's inequality for weak \,$\ell^1$-sequences as before, 
the estimate \eqref{eq:inf:prodholo:2:claim-h} follows (with an appropriately chosen \,$C$\,). 
\end{proof}

\begin{rem}
\label{R:inf:startingindex}
Mutatis mutandis, the results of the present appendix remain true if the starting index \,$k=1$\, of the infinite sums or products is replaced by \,$k=0$\,.
\end{rem}



\section{Index of Notations}

\begin{longtable}{lp{12cm}}
\multicolumn{2}{l}{\textbf{Latin Letters.}} \\
\,$A_n$\, & ``small'' cycle in the homology basis of \,$\Sigma$\,, p.~\pageref{not:jacobi:homology-basis}. \\
\,$A_n^+$\,, \,$A_n^-$\, & boundary cycles of the surface \,$\wt{\Sigma}$\, from Proposition~\ref{P:jacobi:lnmu} on p.~\pageref{P:jacobi:lnmu}. \\
\,$\As(G,\ell^p_{n,m},s)$\, & space of asymptotic functions on \,$G\subset \C^*$\, or \,$G\subset \Sigma$\,, Definition~\ref{D:As:As} on p.~\pageref{D:As:As}. \\
\,$\As_\infty(G_\infty,\ell^p_{n,m},s)$\, & space of asymptotic functions on \,$G_\infty\subset \C^*$\, or \,$G_\infty\subset \Sigma$\, near \,$\lambda=\infty$\,, Definition~\ref{D:As:As0inf} on p.~\pageref{D:As:As0inf}. \\
\,$\As_0(G_0,\ell^p_{n,m},s)$\, & space of asymptotic functions on \,$G_0\subset \C^*$\, or \,$G_0\subset \Sigma$\, near \,$\lambda=0$\,, Definition~\ref{D:As:As0inf} on p.~\pageref{D:As:As0inf}. \\
\,$a(\lambda)$\, & upper-left entry of the monodromy \,$M(\lambda)$\,, Equation~\eqref{eq:spectral:M-entries} on p.~\pageref{eq:spectral:M-entries}. \\
\,$a_0(\lambda)$\, & upper-left entry of the monodromy \,$M_0(\lambda)$\, of the vacuum, Equation~\eqref{eq:vacuum:M0} on p.~\pageref{eq:vacuum:M0}. \\
\,$B_n$\, & ``large'' cycle in the homology basis of \,$\Sigma$\,, p.~\pageref{not:jacobi:homology-basis}. \\
\,$b(\lambda)$\, & upper-right entry of the monodromy \,$M(\lambda)$\,, Equation~\eqref{eq:spectral:M-entries} on p.~\pageref{eq:spectral:M-entries}. \\
\,$b_0(\lambda)$\, & upper-right entry of the monodromy \,$M_0(\lambda)$\, of the vacuum, Equation~\eqref{eq:vacuum:M0} on p.~\pageref{eq:vacuum:M0}. \\
\,$C^+$\,, \,$C^-$\, & boundary cycles of the surface \,$\wt{\Sigma}$\, from Proposition~\ref{P:jacobi:lnmu} on p.~\pageref{P:jacobi:lnmu}. \\
\,$\mathfrak{C}_{D^o}$\, & a certain space of sequences of curves emanating from the support of the divisor \,$D^o \in \Div(\Sigma')$\,, Proposition~\ref{P:jacobi:jaccoord} on p.~\pageref{P:jacobi:jaccoord}. \\
\,$c(\lambda)$\, & lower-left entry of the monodromy \,$M(\lambda)$\,, Equation~\eqref{eq:spectral:M-entries} on p.~\pageref{eq:spectral:M-entries}. \\
\,$c_0(\lambda)$\, & lower-left entry of the monodromy \,$M_0(\lambda)$\, of the vacuum, Equation~\eqref{eq:vacuum:M0} on p.~\pageref{eq:vacuum:M0}. \\
\,$D$\, & classical spectral divisor of a potential \,$(u,u_y)$\, or of a monodromy \,$M(\lambda)$\,, p.~\pageref{not:spectral:D-classical}. \\
\,$D_0$\, & classical spectral divisor of the vacuum, Equation~\eqref{eq:vacuum:D0} on p.~\pageref{eq:vacuum:D0}. \\
\,$\calD$\, & generalized spectral divisor of a potential \,$(u,u_y)$\, or of a monodromy \,$M(\lambda)$\,, Definition~\ref{D:spectral:D-generalized} on p.~\pageref{D:spectral:D-generalized}. \\
\,$\calD_0$\, & generalized spectral divisor of the vacuum, p.~\pageref{not:vacuum:gen-D0}. \\
\,$\Div$\, & space of asymptotic divisors, Definition~\ref{D:asympdiv:asympdiv} on p.~\pageref{not:asympdiv:Div} and the following remarks. \\
\,$\Div(\Sigma)$\, & space of asymptotic divisors with support on \,$\Sigma$\,, Equation~\eqref{eq:jacobi:DivSigma} on p.~\pageref{eq:jacobi:DivSigma}. \\
\,$\Div(\Sigma')$\, & space of asymptotic divisors with support on \,$\Sigma'$\,, Equation~\eqref{eq:jacobi:DivSigma} on p.~\pageref{eq:jacobi:DivSigma}. \\
\,$\Div_{fin}$\, & set of asymptotic divisors of finite type, p.~\pageref{not:diffeo:Div-fin}. \\
\,$\Div_{tame}$\, & space of tame asymptotic divisors, Definition~\ref{D:special:tame}(1) on p.~\pageref{D:special:tame}. \\
\,$\Div_{tame}(\Sigma)$\, & space of tame asymptotic divisors with support on \,$\Sigma$\,, Equation~\eqref{eq:jacobi:DivSigma-tame} on p.~\pageref{eq:jacobi:DivSigma-tame}. \\
\,$\Div_{tame}(\Sigma')$\, & space of tame asymptotic divisors with support on \,$\Sigma'$\,, Equation~\eqref{eq:jacobi:DivSigma-tame} on p.~\pageref{eq:jacobi:DivSigma-tame}. \\
\,$d(\lambda)$\, & lower-right entry of the monodromy \,$M(\lambda)$\,, Equation~\eqref{eq:spectral:M-entries} on p.~\pageref{eq:spectral:M-entries}. \\
\,$d_0(\lambda)$\, & lower-right entry of the monodromy \,$M_0(\lambda)$\, of the vacuum, Equation~\eqref{eq:vacuum:M0} on p.~\pageref{eq:vacuum:M0}. \\
\,$F$\,, \,$F_\lambda$\, & extended frame of a solution \,$u$\, or of a potential \,$(u,u_y)$\,, p.~\pageref{not:spectrum:F}, p.~\pageref{eq:asymp:Fode}. \\
\,$F_0$\, & extended frame of the vacuum, Equation~\eqref{eq:vacuum:F0} on p.~\pageref{eq:vacuum:F0}. \\
\,$f_{n,\xi,\rho}$\, & a meromorphic function on \,$\C^*$\, defined in Proposition~\ref{P:jacobi:asymp-f} on p.~\pageref{P:jacobi:asymp-f}. \\
\,$\wh{f}(k)$\, & Fourier transform of a function \,$f$\,, Equation~\eqref{eq:fasymp:fourier} on p.~\pageref{eq:fasymp:fourier}. \\
\,$\Jac(\Sigma)$\, & Jacobi variety of \,$\Sigma$\,, Theorem~\ref{T:jacobi:jacobi}(1) on p.~\pageref{T:jacobi:jacobi}. \\
\,$\wt{\Jac}(\Sigma)$\, & covering space for the Jacobi variety of \,$\Sigma$\,, Proposition~\ref{P:jacobi:jaccoord}(2) on p.~\pageref{P:jacobi:jaccoord}. \\
\,$j_0$\, & a parameter for a generalized divisor \,$\calD$\,, Proposition~\ref{P:spectrum:locally-free} on p.~\pageref{P:spectrum:locally-free}. \\ 
\,$\wh{K}_r(0)$\,, \,$\wh{K}_r(\infty)$\, & two circles in \,$\Sigma$\, around \,$0$\, resp.~\,$\infty$\,, p.~\pageref{not:jacobitrans:Res}. \\
\,$L^p(M)$\, & Banach space of \,$p$-integrable functions on some domain \,$M$\,. \\
\,$L^2(\Sigma,T^*\Sigma)$\, & space of square-integrable 1-forms on \,$\Sigma$\,. \\
\,$\ell^p(M)$\, & space of \,$p$-summable sequences on an index set \,$M\subset \Z$\,. \\
\,$\ell^p_n(M)$\, & space of \,$p$-summable sequences shifted by \,$k^n$\, on an index set \,$M\subset \Z$\,, p.~\pageref{not:fasymp:ellp-n}. \\ 
\,$\ell^p_{n,m}(M)$\, & space of \,$p$-summable sequences shifted by \,$k^n$\, resp.~\,$k^m$\, on an index set \,$M\subset \Z$\,, p.~\pageref{not:asympdiv:ellp-nm}. \\ 
\,$M(\lambda)$\, & monodromy of a solution \,$u$\, or of a potential \,$(u,u_y)$\,, p.~\pageref{not:spectrum:M}, p.~\pageref{not:asymp:M}. \\ 
\,$M_0(\lambda)$\, & monodromy of the vacuum, Equation~\eqref{eq:vacuum:M0} on p.~\pageref{eq:vacuum:M0}. \\
\,$\calM$\, & sheaf of meromorphic functions on the spectral curve \,$\Sigma$\,, p.~\pageref{not:spectral:calMcalO}. \\
\,$\Mon$\, & space of asymptotic monodromies, p.~\pageref{not:asympfinal:Mon}. \\
\,$\Mon_{np}$\, & space of non-periodic asymptotic monodromies, p.~\pageref{not:asympfinal:Mon}. \\
\,$\Mon_\tau$\, & space of asymptotic monodromies with asymptotic parameter \,$\tau$\,, p.~\pageref{not:asympfinal:Mon}. \\
\,$\Mon_{\tau,\upsilon}$\, & space of non-periodic asymptotic monodromies with asymptotic parameters \,$\tau$\, and \,$\upsilon$\,, p.~\pageref{not:asympfinal:Mon}. \\
\,$\calO$\, & sheaf of holomorphic functions on the spectral curve \,$\Sigma$\,, p.~\pageref{not:spectral:calMcalO}. \\ 
\,$\calO_M$\, & sheaf of holomorphic \,$(2\times 2)$-matrices which commute with the monodromy \,$M(\lambda)$\,, p.~\pageref{not:spectral:eigensheaf}. \\
\,$\wh{\calO}$\, & direct image of the sheaf of holomorphic functions on the normalization \,$\wh{\Sigma}$\, of the spectral curve \,$\Sigma$\, onto \,$\Sigma$\,, p.~\pageref{not:spectral:whSigma}. \\
\,$\Pot$\, & space of potentials (Cauchy data for the sinh-Gordon equation), Equation~\eqref{eq:asymp:Pot} on p.~\pageref{eq:asymp:Pot}. \\
\,$\Pot_{np}$\, & space of non-periodic potentials, p.~\pageref{not:asymp:Potnp}. \\
\,$\Pot_{np}^1$\, & space of once more differentiable, non-periodic potentials, Equation~\eqref{eq:asymp:Potnp1} on p.~\pageref{eq:asymp:Potnp1}. \\
\,$\Pot_{tame}$\, & space of tame potentials, Definition~\ref{D:special:tame}(3) on p.~\pageref{D:special:tame}. \\
\,$\calR$\, & the subsheaf of \,$\wh{\calO}$\, over which a generalized divisor \,$\calD$\, on \,$\Sigma$\, is locally free, Proposition~\ref{P:spectrum:locally-free} on p.~\pageref{P:spectrum:locally-free}. \\
\,$\calR_M$\, & sheaf of eigenvalues of the monodromy \,$M(\lambda)$\,, p.~\pageref{not:spectral:eigensheaf}. \\
\,$\Res_0(\eta)$\,, \,$\Res_\infty(\eta)$\, & residue of a \,$1$-form \,$\eta$\, on \,$\Sigma$\, at \,$0$\, resp.~\,$\infty$\,, p.~\pageref{not:jacobitrans:Res}. \\
\,$S$\, & \,$=\Menge{k\in \Z}{\vkap_{k,1}=\vkap_{k,2}}$\,, Equation~\eqref{eq:holo1forms:S-def} on p.~\pageref{eq:holo1forms:S-def}. \\
\,$S^3$\,, \,$S^3(\vkap)$\, & the 3-sphere (of curvature \,$\vkap>0$\,), p.~\pageref{not:minimal:S3}. \\
\,$S_k$\, & annulus near \,$\lambda_{k,0}$\, in \,$\C^*$\,, Equation~\eqref{eq:As:Sk} on p.~\pageref{eq:As:Sk}. \\
\,$\wh{S}_k$\, & annulus near \,$\lambda_{k,0}$\, in \,$\Sigma$\,, Equation~\eqref{eq:As:whSk} on p.~\pageref{eq:As:whSk}. \\ 
\,$s_{n,\ell}$\,, \,$s_{n,n}$\, & numbers describing the member \,$\omega_n$\, of the canonical basis of \,$1$-forms on \,$\Sigma$\,, Theorem~\ref{T:jacobi:canonical}(2) on p.~\pageref{T:jacobi:canonical}. \\
\,$U_{k,\delta}$\, & \,$k$-th excluded domain in \,$\C^*$\,, Definition~\ref{D:vac2:excldom} on p.~\pageref{D:vac2:excldom}. \\ 
\,$U_{k,\delta}'$\, & possibly punctured \,$k$-th excluded domain in \,$\C^*$\,, Equation~\eqref{eq:jacobiprep:Ukd'} on p.~\pageref{eq:jacobiprep:Ukd'}. \\
\,$\wh{U}_{k,\delta}$\, & \,$k$-th excluded domain in \,$\Sigma$\,, p.~\pageref{not:excl:excldom-Sigma}. \\
\,$\wh{U}_{k,\delta}'$\, & possibly punctured \,$k$-th excluded domain in \,$\Sigma$\,, Equation~\eqref{eq:jacobiprep:Ukd'} on p.~\pageref{eq:jacobiprep:Ukd'}. \\
\,$U_{\delta}$\, & union of all excluded domains in \,$\C^*$\,, Definition~\ref{D:vac2:excldom} on p.~\pageref{D:vac2:excldom}. \\ 
\,$\wh{U}_{\delta}$\, & union of all excluded domains in \,$\Sigma$\,, p.~\pageref{not:excl:excldom-Sigma}. \\
\,$u$\, & solution of the sinh-Gordon equation \,$\Delta u+\sinh(u)=0$\,, p.~\pageref{not:minimal:u}. \\
\,$(u,u_y)$\, & a potential in \,$\Pot$\, or \,$\Pot_{np}$\,, i.e.~Cauchy data for the sinh-Gordon equation. \\
\,$V_{\delta}$\, & area outside of the excluded domains in \,$\C^*$\,, Definition~\ref{D:vac2:excldom} on p.~\pageref{D:vac2:excldom}. \\ 
\,$\wh{V}_{\delta}$\, & area outside of the excluded domains in \,$\Sigma$\,, p.~\pageref{not:excl:excldom-Sigma}. \\
\,$v_k$\, & symplectic basis vector of \,$T_{(u,u_y)}\Pot$\,, Equation~\eqref{eq:darboux:vkwk} on p.~\pageref{eq:darboux:vkwk} and Theorem~\ref{T:darboux:darboux}(1) on p.~\pageref{T:darboux:darboux}. \\
\,$W^{n,p}(M)$\, & Sobolev space of weakly \,$n$ times differentiable functions with \,$p$-integrable derivative on a domain \,$M$\,. \\
\,$w(\lambda)$\, & \,$= |\cos(\zeta(\lambda))|+|\sin(\zeta(\lambda))|$\,, Equation~\eqref{eq:asymp:w} on p.~\pageref{eq:asymp:w}. \\
\,$w_k$\, & symplectic basis vector of \,$T_{(u,u_y)}\Pot$\,, Equation~\eqref{eq:darboux:vkwk} on p.~\pageref{eq:darboux:vkwk} and Theorem~\ref{T:darboux:darboux}(1) on p.~\pageref{T:darboux:darboux}. \\
\\
\multicolumn{2}{l}{\textbf{Greek Letters.}} \\
\,$\alpha$\,, \,$\alpha_\lambda$\, & the \,$\mathfrak{sl}(2,\C)$-valued flat connection form of a solution \,$u$\, or of a potential \,$(u,u_y)$\,, p.~\pageref{not:spectrum:alpha}, p.~\pageref{eq:asymp:alpha}. \\
\,$\alpha_0$\, & the \,$\mathfrak{sl}(2,\C)$-valued flat connection form associated to the vacuum, Equation~\eqref{eq:vacuum:alpha0} on p.~\pageref{eq:vacuum:alpha0}. \\
\,$\Gamma$\, & period lattice of \,$\Sigma$\,, Theorem~\ref{T:jacobi:jacobi}(1) on p.~\pageref{T:jacobi:jacobi}. \\
\,$\Delta(\lambda)$\, & trace of the monodromy \,$M(\lambda)$\,, Equation~\eqref{eq:spectral:tr} on p.~\pageref{eq:spectral:tr}.\\
\,$\Delta_0(\lambda)$\, & trace of the monodromy \,$M_0(\lambda)$\, of the vacuum, Equation~\eqref{eq:vacuum:Delta0} on p.~\pageref{eq:vacuum:Delta0}.\\
\,$\delta f$\, & variation of the function \,$f$\, defined on \,$\Pot$\,, p.~\pageref{not:darboux:deltaf}. \\
\,$(\delta u,\delta u_y)$\, & a tangent vector in \,$T_{(u,u_y)}\Pot$\,, p.~\pageref{not:darboux:deltaf}. \\
\,$\delta_{k\ell}$\, & Kronecker delta, i.e.~\,$\delta_{k\ell}=1$\, for \,$k=\ell$\,, \,$\delta_{k\ell}=0$\, for \,$k \neq \ell$\,. \\
\,$\zeta(\lambda)$\, & \,$=\tfrac14(\lambda^{1/2}+\lambda^{-1/2})$\,, Equation~\eqref{eq:vacuum:zeta} on p.~\pageref{eq:vacuum:zeta}. \\
\,$\wt{\zeta}(\lambda)$\, & \,$=\tfrac14(\lambda^{1/2}-\lambda^{-1/2})$\,, Equation~\eqref{eq:vacuum:wtzeta} on p.~\pageref{eq:vacuum:wtzeta}. \\
\,$\eta_k$\, & zero of \,$\Delta'$\,, Lemma~\ref{L:finite:eta}(1) on p.~\pageref{L:finite:eta}. \\
\,$\vartheta_k$\, & factor for the symplectic basis \,$(v_k,w_k)$\, of \,$T_{(u,u_y)}\Pot$\,, Equation~\eqref{eq:darboux:vartheta} on p.~\pageref{eq:darboux:vartheta} and Theorem~\ref{T:darboux:darboux}(1) on p.~\pageref{T:darboux:darboux}. \\
\,$\vkap_{k,\nu}$\, & zero of \,$\Delta^2-4$\, (branch point or singularity of the spectral curve \,$\Sigma$\,), Proposition~\ref{P:excl:basic}(1) on p.~\pageref{P:excl:basic}. \\
\,$\vkap_{k,*}$\, & \,$=\tfrac12\,(\vkap_{k,1}+\vkap_{k,2})$\,, Equation~\eqref{eq:jacobiprep:vkap-k-star} on p.~\pageref{eq:jacobiprep:vkap-k-star}. \\
\,$\Lambda$\, & eigenline bundle of the monodromy \,$M(\lambda)$\,, p.~\pageref{not:spectral:Lambda}. \\
\,$\lambda$\, & spectral parameter, p.~\pageref{not:minimal:lambda}. \\
\,$\lambda_k$\, & \,$\lambda$-parameter of a spectral divisor point, Proposition~\ref{P:excl:basic}(2) on p.~\pageref{P:excl:basic}. \\
\,$\lambda_{k,0}$\, & \,$\lambda$-parameter of a divisor point of the vacuum, Equation~\eqref{eq:vacuum:lambdak0-def} on p.~\pageref{eq:vacuum:lambdak0-def}. \\
\,$\mu$\, & eigenvalue of the monodromy \,$M(\lambda)$\,, Equation~\eqref{eq:spectral:mu} on p.~\pageref{eq:spectral:mu} \\
\,$\mu_k$\, & \,$\mu$-parameter of a spectral divisor point, Proposition~\ref{P:excl:basic}(2) on p.~\pageref{P:excl:basic}. \\
\,$\mu_{k,0}$\, & \,$\mu$-parameter of a divisor point of the vacuum, Equation~\eqref{eq:vacuum:muk0-def} on p.~\pageref{eq:vacuum:muk0-def}. \\
\,$\Sigma$\, & spectral curve, Equation~\eqref{eq:spectral:Sigma} on p.~\pageref{eq:spectral:Sigma}. \\ 
\,$\Sigma'$\, & the spectral curve punctured at its singularities, Equation~\eqref{eq:holo1forms:Sigma'-def} on p.~\pageref{eq:holo1forms:Sigma'-def}. \\ 
\,$\wh{\Sigma}$\, & normalization of the spectral curve, p.~\pageref{not:spectral:whSigma}. \\
\,$\wt{\Sigma}$\, & surface obtained from \,$\Sigma$\, by cutting along the cycles \,$A_n$\,, Proposition~\ref{P:jacobi:lnmu} on p.~\pageref{P:jacobi:lnmu}. \\
\,$\Sigma_0$\, & spectral curve of the vacuum, Equation~\eqref{eq:vacuum:Sigma0} on p.~\pageref{eq:vacuum:Sigma0}.\\
\,$\sigma$\, & hyperelliptic involution of the spectral curve, Equation~\eqref{eq:spectral:sigma} on p.~\pageref{eq:spectral:sigma}. \\
\,$\tau$\, & asymptotic parameter for the component functions \,$b(\lambda)$\,, \,$c(\lambda)$\, of the monodromy \,$M(\lambda)$\,, Theorem~\ref{T:asymp:basic} on p.~\pageref{T:asymp:basic}. \\
\,$\upsilon$\, & asymptotic parameter for the component function \,$a(\lambda)$\,, \,$d(\lambda)$\, of the monodromy \,$M(\lambda)$\,, Theorem~\ref{T:asymp:basic} on p.~\pageref{T:asymp:basic}. \\
\,$\Phi_{n,\xi,\rho}$\, & holomorphic function on \,$\C^*$\, defined by an infinite product in Proposition~\ref{P:jacobi:asymp-Phi} on p.~\pageref{P:jacobi:asymp-Phi}. \\
\,$\vi$\, & Abel map of \,$\Sigma$\,, Theorem~\ref{T:jacobi:jacobi}(2) on p.~\pageref{T:jacobi:jacobi}. \\
\,$\vi_n$\, & \,$n$-th Jacobi coordinate for \,$\Sigma$\,, Theorem~\ref{T:jacobi:jacobi}(4) on p.~\pageref{T:jacobi:jacobi}. \\
\,$\wt{\vi}$\, & lift of the Abel map of \,$\Sigma$\,, Proposition~\ref{P:jacobi:jaccoord}(2) on p.~\pageref{P:jacobi:jaccoord}. \\
\,$\wt{\vi}_n$\, & lift of the \,$n$-th Jacobi coordinate for \,$\Sigma$\,, Proposition~\ref{P:jacobi:jaccoord}(1) on p.~\pageref{P:jacobi:jaccoord}. \\
\,$\Psi_k$\, & \,$=\sqrt{(\lambda-\vkap_{k,1})\cdot (\lambda-\vkap_{k,2})}$\,, Lemma~\ref{L:jacobiprep:Psi-intro} on p.~\pageref{L:jacobiprep:Psi-intro}. \\
\,$\Omega$\, & symplectic form on \,$\Pot$\,, Equation~\eqref{eq:darboux:Omega-defined} on p.~\pageref{eq:darboux:Omega-defined}. \\
\,$\Omega(\Sigma)$\, & space of holomorphic 1-forms on \,$\Sigma$\,, p.~\pageref{not:holo1forms:Omega-L2}. \\
\,$\omega_n$\, & canonical basis of holomorphic 1-forms on \,$\Sigma$\, resp.~\,$\Sigma'$\,, Theorem~\ref{T:jacobi:canonical}(2) on p.~\pageref{T:jacobi:canonical}. \\
\,$\wt{\omega}_n$\, & certain holomorphic 1-forms on \,$\Sigma$\, that are not quite as good as the member \,$\omega_n$\, of the canonical basis, Theorem~\ref{T:jacobi:canonical}(1) on p.~\pageref{T:jacobi:canonical}. \\
\\
\multicolumn{2}{l}{\textbf{Others.}} \\
\,$*$\, & convolution of two sequences, Equation~\eqref{eq:fasymp:convolution} on p.~\pageref{eq:fasymp:convolution}. 
\end{longtable}

\end{document}